\documentclass[11pt, twoside]{uwthesis}
\usepackage{amsmath,amssymb,amsthm,xspace,stmaryrd}

\theoremstyle{plain} 
\newtheorem{thm}{Theorem}[chapter]
\newtheorem{cor}[thm]{Corollary}
\newtheorem{lem}[thm]{Lemma}
\newtheorem{prop}[thm]{Proposition}
\newtheorem{conj}[thm]{Conjecture}

\theoremstyle{definition}
\newtheorem{definition}[thm]{Definition}
\newtheorem{example}[thm]{Example}

\newtheorem{claim}{Claim}[thm]
\newtheorem{thmremark}{Remark}[thm] 
\newtheorem{remark}[thm]{Remark} 

\newcommand{\extraline}{\vskip\baselineskip}

\newcommand{\proofpart}[1]{\noindent \underline{Part #1:}\xspace}

\newcommand{\textdef}{\textbf} 

\newcommand{\textstack}[2]{\genfrac{}{}{0pt}{}{\text{#1}}{\text{#2}}}
\newcommand{\mathstack}[2]{\genfrac{}{}{0pt}{}{#1}{#2}}

\renewcommand{\th}{^\text{th}} 

\newcommand{\vect}{\mathbf} 

\newcommand{\leftexp}[2]{{\vphantom{#2}}^{#1}{#2}} 

\newcommand{\N}{\mathbb{N}} 
\newcommand{\naturalzero}{\mathbb{N}} 
\newcommand{\naturalpositive}{\mathbb{N}_{\ne 0}}  
\newcommand{\R}{\mathbb{R}} 
\newcommand{\Rplus}{\R_+} 
\newcommand{\Z}{\mathbb{Z}} 
\newcommand{\halfspace}[1][]{\mathbb{H}^{#1}} 
\newcommand{\rationals}{\mathbb{Q}} 

\newcommand{\Zd}{{\Z^d}} 
\newcommand{\Rd}{{\R^d}} 

\newcommand{\delimit}[4][0]{
\ifcase#1
	#2 #4 #3				
	\or \left#2 #4 \right#3		
	\or \bigl#2 #4\bigr#3		
	\or \Bigl#2 #4\Bigr#3		
	\or \biggl#2 #4\biggr#3	
	\or \Biggl#2 #4\Biggr#3	
	\else #4				
\fi}

\newcommand{\tripledelimit}[6][0]{%
\ifcase#1
	#2 #5 \;#3\; #6 #4					
	\or \left#2 #5 \;\middle#3\; #6 \right#4		
	\or \bigl#2 #5 \bigm#3 #6 \bigr#4		
	\or \Bigl#2 #5 \Bigm#3 #6 \Bigr#4		
	\or \biggl#2 #5 \biggm#3 #6 \biggr#4		
	\else \Biggl#2 #5 \Biggm#3 #6 \Biggr#4	
\fi}

\newcommand{\rightdelimit}[3][0]{
\ifcase#1
	#3 #2			
	\or \left. #3 \right#2	
	\or #3 \bigr#2		
	\or #3 \Bigr#2		
	\or #3 \biggr#2		
	\or #3 \Biggr#2		
	\else #3			
\fi}

\newcommand{\norm}[2][0]{\delimit[#1]{\lVert}{\rVert}{#2}}
\newcommand{\abs}[2][0]{\delimit[#1]{\lvert}{\rvert}{#2}}
\newcommand{\braces}[2][0]{\delimit[#1]{\lbrace}{\rbrace}{#2}}
\newcommand{\brackets}[2][0]{\delimit[#1]{\lbrack}{\rbrack}{#2}}
\newcommand{\parens}[2][0]{\delimit[#1]{(}{)}{#2}}
\newcommand{\abrackets}[2][0]{\delimit[#1]{\langle}{\rangle}{#2}}
\newcommand{\bbrackets}[2][0]{\delimit[#1]{\llbracket}{\rrbracket}{#2}}
\newcommand{\floor}[2][0]{\delimit[#1]{\lfloor}{\rfloor}{#2}}
\newcommand{\ceil}[2][0]{\delimit[#1]{\lceil}{\rceil}{#2}}

\newcommand{\ifandonlyif}{\Longleftrightarrow}
\newcommand{\followsfrom}{\Longleftarrow}

\newcommand{\card}[2][0]{\abs[#1]{#2}} 
\newcommand{\powerset}[1]{2^{#1}}

\newcommand{\definedas}{:=} 
\newcommand{\defines}{=:} 

\newcommand{\restrict}[3][0]{\rightdelimit[#1]{\rvert}{#2}_{#3}}

\newcommand{\meet}{\wedge}
\newcommand{\join}{\vee}

\newcommand{\setof}[2][0]{\braces[#1]{#2}}
\newcommand{\pair}[3][0]{\parens[#1]{#2,#3}}

\newcommand{\setbuilder}[3][0]{\braces[#1]{#2 : #3}}

\newcommand{\setbuilderbar}[3][0]{\tripledelimit[#1]{\lbrace}{\vert}{\rbrace}{#2}{#3}}

\newcommand{\bd}{\partial} 
\newcommand{\interior}[1]{#1^{\circ}} 
\newcommand{\closure}[1]{\overline{#1}} 
\DeclareMathOperator{\ext}{ext} 
\newcommand{\exterior}[2][10]{\ext \parens[#1]{#2}} 

\DeclareMathOperator{\shell}{shell} 
\newcommand{\shellof}[2][10]{\shell \parens[#1]{#2}} 

\newcommand{\innerprod}[3][0]{\abrackets[#1]{#2,#3}}
\newcommand{\basisvec}[1]{\vect e_{#1}}

\DeclareMathOperator{\lengthOp}{Len}
\newcommand{\lengthfcn}[1]{\lengthOp_{#1}}
\newcommand{\length}[3][0]{\lengthfcn{#2} \parens[#1]{#3}}

\DeclareMathOperator{\dil}{dil} 
\newcommand{\dilnorm}[3][0]{\dil_{#2} \parens[#1]{#3}} 
\newcommand{\esupnorm}[2][0]{\norm[#1]{#2}_{L^\infty}} 
\newcommand{\esupnormon}[3][0]{\norm[#1]{#3}_{L^\infty(#2)}} 

\newcommand{\lebmeas}[1]{\mathfrak{m}_L^{#1}}
\newcommand{\lebmeasof}[3][0]{\lebmeas{#2} \parens[#1]{#3}}

\DeclareMathOperator{\supp}{supp} 

\newcommand{\lspace}[2]{\ell^{#1}_{#2}}
\newcommand{\linfty}[1]{\lspace{\infty}{#1}}
\newcommand{\ltwo}[1]{\lspace{2}{#1}}
\newcommand{\lone}[1]{\lspace{1}{#1}}

\DeclareMathOperator{\latOp}{lat} 

\newcommand{\lat}[2][0]{%
\ifcase#1
	\dddot{#2}				
	\or \ddddot{#2}			
	\or \overset{\dotfill}{#2}	
	\or \latOp \parens[0]{#2}	
	\or \latOp \parens[1]{#2} 	
	\or \latOp \parens[2]{#2} 	
	\else \latOp \parens[3]{#2} 	
\fi}

\newcommand{\cubify}[2][0]{\bbrackets[#1]{#2}} 
\newcommand{\intcubify}[2][0]{\interior{\cubify[#1]{#2}}}

\newcommand{\normshellbound}[1]{\varpi_{#1}}
\newcommand{\mushellbound}{\normshellbound{\mu}}

\newcommand{\normcylconst}[2]{\varpi_{#1,#2}}

\newcommand{\graphfont}{\mathsf}
\newcommand{\edgesymb}{\graphfont{E}}
\newcommand{\nbrsymb}{\graphfont{N}}

\newcommand{\edge}{\mathsf{e}} 

\newcommand{\nbrset}[2][0]{\nbrsymb \parens[#1]{#2}} 
\newcommand{\nbrhood}[2][0]{\nbrsymb \brackets[#1]{#2}} 

\newcommand{\edges}[2][0]{\edgesymb \parens[#1]{#2}} 

\newcommand{\jedges}[3][0]{\edgesymb \parens[#1]{#2\text{-}#3}}

\newcommand{\compedges}[2][0]{\edgesymb^\complement \parens[#1]{#2}} 
\newcommand{\bdedges}[2][0]{\edgesymb_\bd \parens[#1]{#2}} 
\newcommand{\staredges}[2][0]{\edgesymb^* \parens[#1]{#2}} 

\newcommand{\edgecomplement}[2][10]{%
\ifnum#1>5 #2^\complement
\else \edgesymb^\complement \parens[#1]{#2}
\fi}

\newcommand{\stargraph}[2][10]{\parens[#1]{#2}^*}

\newcommand{\subpath}[4][0]{#2\brackets[#1]{#3,#4}}

\newcommand{\pathclass}{\Pi}

\renewcommand{\Pr}{\operatorname{\mathbf{Pr}}} 
\newcommand{\E}{\operatorname{\mathbf{E}}} 

\newcommand{\eventcomplement}[2][10]{\parens[#1]{#2}^\complement}

\newcommand{\cprob}[3][0]{\Pr \tripledelimit[#1]{(}{\vert}{)}{#2}{#3}}
\newcommand{\cprobbrack}[3][0]{\Pr \tripledelimit[#1]{\lbrack}{\vert}{\rbrack}{#2}{#3}}

\newcommand{\expectation}[3][0]{
  \ifcase#1 
     E( #2 \mid #3 ) 
     \or E \bigl( #2 \bigm\vert #3 \bigr)
     \or E \Bigl( #2 \Bigm\vert #3 \Bigr) 
     \or E \biggl( #2 \biggm\vert #3 \biggr)
     \or E \Biggl( #2 \Biggm\vert #3 \Biggr) 
  \else 
     E \left( #2  \;\middle\vert\; #3 \right)
  \fi}

\newcommand{\ind}[1]{{\mathbf{1}_{#1}}} 

\newcommand{\law}{\mathcal{L}} 

\newcommand{\iid}{i.i.d.\xspace} 

\newcommand{\pcrit}{p_c(\Zd)} 

\newcommand{\eventthat}[2][0]{\setof[#1]{#2}}

\newcommand{\indist}[1]{\overset{d}{#1}}
\newcommand{\eqd}{\indist{=}}

\newcommand{\assubseteq}{\subseteq_{\rm a.s.}} 
\newcommand{\asequal}{=_{\rm a.s.}} 

\newcommand{\explaw}[2][0]{\operatorname{Exponential} \parens[#1]{#2}}

\newcommand{\samplespace}[1][0]{%
\ifcase#1
	\Omega
	\or \parens{\R_+\times\R_+}^{\edges{\Zd}}
	\else #1
\fi
}
\newcommand{\goodsamplespace}{\samplespace_{\rm good}}
\newcommand{\rsamplespace}[1]{\restrict{\samplespace}{#1}}
\newcommand{\outcome}{\omega}
\newcommand{\sigmafield}{\mathcal{A}}
\newcommand{\probabilityspace}{(\samplespace,\sigmafield,\Pr)}

\newcommand{\jointlaw}{\nu}

\newcommand{\marginals}[1]{\jointlaw^{(#1)}}

\newcommand{\statespacesymb}{\mathbf{S}}

\newcommand{\statespace}[1][0]{%
\ifcase #1
	\statespacesymb
	\or \setof{0,1,2}^{\Zd}
	\else #1
\fi}

\newcommand{\statesigmafield}{\mathcal{S}}
\newcommand{\statemeasurespace}{(\statespace,\statesigmafield)}

\newcommand{\statesubspace}[2][0]{%
\ifcase #1
	\statespacesymb_{\setof{#2}}
	\or \setof{#2}^{\Zd}
	\else #1
\fi}

\newcommand{\twotypeproduct}{\samplespace\times\statespace}

\newcommand{\processspace}{\statespace^{[0,\infty]}}
\newcommand{\processsigmafield}{\statesigmafield^{\otimes [0,\infty]}}

\newcommand{\processlaw}[2]{\Pr_{#1}^{#2}}

\newcommand{\onetypespace}{\Rplus^{\edges{\Zd}}} 
\newcommand{\onetypestates}{2^{\Zd}} 
\newcommand{\onetypeproduct}{\onetypespace\times \onetypestates} 

\newcommand{\bigcsymb}{C} 
\newcommand{\smallcsymb}{c} 
\newcommand{\konstsymb}{K} 
\newcommand{\smallkonstsymb}{\kappa} 

\newcommand{\bigc}[1]{\bigcsymb_{#1}} 
\newcommand{\smallc}[1]{\smallcsymb_{#1}} 
\newcommand{\konst}[1]{\konstsymb_{#1}} 

\newcommand{\bigcref}[1]{\bigc{\ref{#1}}} 
\newcommand{\smallcref}[1]{\smallc{\ref{#1}}} 
\newcommand{\konstref}[1]{\konst{\ref{#1}}} 

\newcommand{\symboleqref}[2]{#1_{\eqref{#2}}} 
\newcommand{\symbolref}[2]{#1_{\ref{#2}}} 

\newcommand{\eventsymb}{E} 
\newcommand{\eventref}[1]{\symbolref{\eventsymb}{#1}} 

\DeclareMathOperator{\conv}{conv} 
\DeclareMathOperator{\aff}{aff} 

\newcommand{\affrelbd}{\bd_{\aff}} 

\newcommand{\affspan}[2][0]{\aff \parens[#1]{#2}}
\newcommand{\convhull}[2][0]{\conv \parens[#1]{#2}}

\newcommand{\Hspace}{J} 
\newcommand{\Hplane}{H} 

\newcommand{\segment}[3][0]{\brackets[#1]{#2, #3}}
\newcommand{\ocsegment}[3][0]{\delimit[#1]{(}{]}{#2, #3}}
\newcommand{\cosegment}[3][0]{\delimit[#1]{[}{)}{#2, #3}}
\newcommand{\osegment}[3][0]{\delimit[#1]{(}{)}{#2, #3}}

\newcommand{\ray}[2]{\overrightarrow{#1 #2}}
\newcommand{\cray}[2]{\leftexp{\bullet}{\ray{#1}{#2}}}

\newcommand{\trianglev}[4][0]{\triangle \brackets[#1]{#2,#3,#4}}

\newcommand{\homothe}[2]{#1_{#2}}

\newcommand{\asiset}{W} 

\newcommand{\wedgesetsymb}{\mathcal{W}} 

\newcommand{\wedgesat}[1]{\wedgesetsymb_{#1}} 
\newcommand{\cwedgesat}[1]{\closure{\wedgesetsymb}_{#1}} 

\DeclareMathOperator{\dir}{dir} 
\newcommand{\dirset}[2][0]{\dir \parens[#1]{#2}}
\newcommand{\dirsetat}[3][0]{\dir_{#2} \parens[#1]{#3}}

\newcommand{\direquiv}{\sim}

\DeclareMathOperator{\apex}{apex} 
\newcommand{\apexset}[2][0]{\apex \parens[#1]{#2}}

\newcommand{\anapex}{\vect a} 

\newcommand{\siversion}[1]{#1_{\vect 0}}

\newcommand{\dirsymb}[1][0]{%
\ifcase#1
	\boldsymbol{\upsilon}	
	\or \boldsymbol{\varphi}	
\fi}

\newcommand{\direction}{\vec} 
\newcommand{\dirvec}[1][0]{\direction{\dirsymb[#1]}} 
\newcommand{\cdirvec}[1][\dirsymb]{\leftexp{{\scriptscriptstyle \bullet}}{\direction{#1}}} 
\newcommand{\unitvec}[1]{\widehat{#1}} 
\newcommand{\unitvecnorm}[2]{\widehat{#1}_{#2}} 
\newcommand{\munitvec}[1][]{\unitvecnorm{\dirsymb}{\mu_{#1}}}

\newcommand{\unitball}[1]{\mathcal B_{#1}} 
\newcommand{\sphere}[2]{\mathbb{S}^{#1}_{#2}} 

\newcommand{\unitcube}[2][1]{\brackets[#1]{-\frac{1}{2},\frac{1}{2}}^{#2}}

\newcommand{\conesymb}{\mathcal C} 
\newcommand{\cone}[4]{\conesymb_{#1,#2}^{#3,#4}} 
\newcommand{\conetip}[6][0]{\cone{#2}{#3}{#4}{#5} \parens[#1]{#6}} 
\newcommand{\conesegment}[7][0]{\cone{#2}{#3}{#4}{#5} \brackets[#1]{#6,#7}} 

\newcommand{\cones}{\mathcal C} 

\newcommand{\conetext}[3][0]{$\parens[#1]{#2,#3}$-cone\xspace} 
\newcommand{\conesegmenttext}[3][0]{\conetext[#1]{#2}{#3} segment\xspace} 

\newcommand{\bpinsymb}{\mathcal{P}}

\newcommand{\bpin}[7][0]{\bpinsymb_{#2,#3,#4}^{#5,#6} \parens[#1]{#7}}

\newcommand{\conefamily}[3][0]{\operatorname{Cones}_{#2} \parens[#1]{#3}}
\newcommand{\conesubfam}{\mathcal{J}}
\newcommand{\coneunion}[1]{\conesymb_{#1}}

\newcommand{\tubesymb}{\mathcal T}
\newcommand{\tube}[4]{\tubesymb_{#1,#2}^{\,#3,#4}}
\newcommand{\tubetip}[6][0]{\tube{#2}{#3}{#4}{#5} \parens[#1]{#6}}
\newcommand{\rfcn}{\rho}

\newcommand{\stars}[2][2]{\abrackets[#1]{\bigstar^{#2}}}
\newcommand{\starbodies}[2][2]{\abrackets[#1]{\bigstar_{+}^{#2}}}

\newcommand{\genstars}[3][2]{\stars[#1]{\exists \bullet \in #3}} 
\newcommand{\starscent}[3][2]{\stars[#1]{\forall \bullet \in #3}} 

\newcommand{\normstars}[3][2]{\abrackets[#1]{\bigstar_{#2}^{#3}}}
\newcommand{\fatstars}[4][2]{\abrackets[#1]{\bigstar_{#2,#3}^{#4}}}

\newcommand{\dist}[1]{\operatorname{dist}_{#1}} 

\DeclareMathOperator{\distOp}{dist}
\newcommand{\normmetric}[1]{\distOp_{#1}} 
\newcommand{\imetric}[2]{\normmetric{#1}^{#2}}

\newcommand{\distnew}[4][0]{\normmetric{#2} \parens[#1]{#3,#4}}
\newcommand{\idist}[5][0]{\imetric{#2}{#3} \parens[#1]{#4,#5}}

\newcommand{\starify}[2][10]{\parens[#1]{#2}^*}
\newcommand{\starmetric}[3][10]{\normmetric{#2}^{\starify[#1]{#3}}}
\newcommand{\stardist}[5][0]{\starmetric{#2}{#3} \parens[#1]{#4,#5}}

\newcommand{\speciesname}[1]{{\sc #1}}
\newcommand{\shapefunction}{\mu}
\newcommand{\limitshape}{\unitball{\shapefunction}}

\newcommand{\mtsymb}{Z}
\newcommand{\mttime}{\mtsymb}
\newcommand{\mttimes}[1]{\mtsymb^{(#1)}}

\newcommand{\cmttimes}[2]{%
\bar{\mtsymb}^{%
\def\testemptyA{#1}%
\ifx\testemptyA\empty \empty
\else (#1)%
\fi}_{#2}
}

\newcommand{\mttimenew}{\ttime(\edge)}

\newcommand{\cmttime}[1]{\ttime_{#1}(\edge)}

\newcommand{\labelededge}[2]{%
\edge^{%
\def\testemptyA{#1}
\ifx\testemptyA\empty \empty
\else \scriptscriptstyle (#1)
\fi}_{#2}%
}

\newcommand{\labeledmttime}[3][0]{\ttime_{#2} \parens[#1]{\labelededge{#2}{#3}}}

\newcommand{\ttimesymb}{\tau} 
\newcommand{\ttime}{\ttimesymb} 
\newcommand{\tmeasure}{\ttimesymb} 

\newcommand{\tmeasureof}[2][0]{\tmeasure \parens[#1]{#2}}

\newcommand{\ttimepair}{\bar{\ttime}} 
\newcommand{\tmeasurepair}{\bar{\tmeasure}}

\newcommand{\mintime}{\ttime_{\min}}

\newcommand{\ptimesymb}{T} 

\newcommand{\ptmetric}[1][]{\ptimesymb_{#1}}
\newcommand{\ptime}{T} 
\newcommand{\ptimenew}[4][0]{\rptime[#1]{#2}{}{#3}{#4}} 

\newcommand{\rptmetric}[2][]{\ptmetric[#1]^{#2}}
\newcommand{\rptime}[5][0]{\ptimesymb_{#2}^{#3} \parens[#1]{#4,#5}} 

\newcommand{\starptmetric}[2][]{\ptmetric[#1]^{\stargraph{#2}}}
\newcommand{\starptime}[1]{\ptime^{\stargraph{#1}}} 
\newcommand{\starptimenew}[5][0]{\rptime[#1]{#2}{\stargraph{#3}}{#4}{#5}}

\newcommand{\ctimesymb}{\widetilde{\ptime}} 
\newcommand{\ctmetric}[2]{\ctimesymb_{#1}^{#2}} 
\newcommand{\ctime}[5][0]{\ctmetric{#2}{#3} \parens[#1]{#4,#5}} 

\newcommand{\ictimesymb}{\ctimesymb^{\rm Int}}
\newcommand{\ictmetric}[1]{\ictimesymb_{#1}}
\newcommand{\ictime}[4][0]{\ictmetric{#2} \parens[#1]{#3,#4}}

\newcommand{\entangledmetric}[2]{\ptimesymb_{#1}^{#2}}
\newcommand{\entangledtime}[4][0]{\entangledmetric{#2}{#3} \parens[#1]{#4}}

\newcommand{\infectiontime}[3][0]{\entangledtime[#1]{1\cup 2}{#2}{#3}}

\newcommand{\entangledctmetric}[2]{\ctimesymb_{#1}^{#2}}
\newcommand{\entangledctime}[4][0]{\entangledctmetric{#2}{#3} \parens[#1]{#4}}

\newcommand{\ptimefamily}[1][]{\ptimesymb_{#1}^\cdot} 
\newcommand{\starptimefamily}[1][]{\ptimesymb_{#1}^{\cdot *}} 
\newcommand{\fppprocfamily}[1][]{\reachedsymb_{#1}^\cdot} 

\newcommand{\reachedsymb}{\eta} 
\newcommand{\reached}{\eta} 
\newcommand{\creachedsymb}{\widetilde{\reachedsymb}} 

\newcommand{\reachedx}[2][]{\reached_{#1}^{#2}}
\newcommand{\rreached}[3][]{\reachedx[#1]{#2;#3}}

\newcommand{\reachedfcn}[2]{\reachedsymb_{#1}^{#2}}
\newcommand{\rreachedfcn}[3]{\reachedfcn{#1}{#2;\smash[b]{#3}}}
\newcommand{\reachedset}[4][0]{\reachedfcn{#2}{#3} \parens[#1]{#4}}
\newcommand{\rreachedset}[5][0]{\rreachedfcn{#2}{#3}{#4} \parens[#1]{#5}}

\newcommand{\creachedfcn}[2]{\creachedsymb_{#1}^{\: #2}}
\newcommand{\crreachedfcn}[3]{\creachedfcn{#1}{#2;\smash[b]{#3}}}

\newcommand{\crreachedset}[5][0]{\crreachedfcn{#2}{#3}{#4} \parens[#1]{#5}}

\newcommand{\bothreachedset}[3][0]{\reachedset[#1]{1\cup 2}{#2}{#3}}

\newcommand{\normpair}{\bar\mu} 
\newcommand{\adjustednormpair}[2]{\normpair^{\scriptscriptstyle{(#1)}}_{#2}}

\newcommand{\reachedd}{\mathcal{B}} 
\newcommand{\reacheddsymb}{\mathcal{B}} 

\newcommand{\reacheddfcn}[2]{\reacheddsymb_{#1}^{#2}}
\newcommand{\rreacheddfcn}[3]{\reacheddfcn{#1}{#2;#3}}
\newcommand{\reacheddset}[4][0]{\reacheddfcn{#2}{#3} \parens[#1]{#4}}
\newcommand{\rreacheddset}[5][0]{\rreacheddfcn{#2}{#3}{#4} \parens[#1]{#5}}

\newcommand{\procsymb}{\xi} 

\newcommand{\processt}[2]{\procsymb^{#1}_{#2}}
\newcommand{\process}[1]{\processt{#1}{\cdot}}
\newcommand{\processtoutcome}[4][0]{\parens[#1]{\processt{#2}{#3}}_{#4}}
\newcommand{\processoutcome}[3][0]{\processtoutcome[#1]{#2}{\cdot}{#3}} 

\newcommand{\processtv}[4][0]{\processt{#2}{#3} \parens[#1]{#4}}

\newcommand{\domain}{D}
\newcommand{\startdomain}{\domain_0}
\newcommand{\domainspace}{\samplespace\times\statespace}

\newcommand{\stateset}[3][0]{\reachedset[#1]{#2}{}{#3}}
\newcommand{\nonzeroset}[2][0]{\stateset[#1]{1\cup 2}{#2}}

\newcommand{\shiftsymb}{\theta} 
\newcommand{\shiftop}[1]{\shiftsymb_{#1}} 

\newcommand{\pshiftsymb}{\bar{\theta}} 
\newcommand{\pshiftop}[1]{\pshiftsymb_{#1}} 

\newcommand{\apshiftsymb}{\tilde{\theta}} 
\newcommand{\apshiftop}[1]{\apshiftsymb_{#1}} 

\newcommand{\rpshiftop}[2]{\pshiftsymb_{#1}^{#2}}

\newcommand{\shiftstate}[4][0]{\pshiftop{#2}^{#3} \parens[#1]{#4}} 
\newcommand{\shifttimes}[4][0]{\pshiftop{#2}^{#3} \parens[#1]{#4}} 
\newcommand{\shiftpair}[4][0]{\pshiftop{#2} \parens[#1]{#3,#4}} 

\newcommand{\ashiftstate}[4][0]{\apshiftop{#2}^{#3} \parens[#1]{#4}} 
\newcommand{\ashifttimes}[4][0]{\apshiftop{#2}^{#3} \parens[#1]{#4}} 
\newcommand{\ashiftpair}[4][0]{\apshiftop{#2} \parens[#1]{#3,#4}} 

\newcommand{\rpshiftstate}[5][0]{\pshiftsymb_{#2}^{#4;#3} \parens[#1]{#5}}
\newcommand{\rpshifttimes}[5][0]{\pshiftsymb_{#2}^{#4;#3} \parens[#1]{#5}}
\newcommand{\rpshiftpair}[5][0]{\rpshiftop{#2}{#3} \pair[#1]{#4}{#5}}

\newcommand{\pbetter}{\gtrsim}
\newcommand{\pworse}{\lesssim}
\newcommand{\pequiv}{\sim}

\newcommand{\rpworse}[1]{\lesssim_{\,#1}}

\newcommand{\rple}[1]{\le_{\,#1}}

\newcommand{\rmeet}[2]{%
\ifcase#1
	BAD!							
	\or \varrho_{#2}						
	\or \varrho_{#2}^{\scriptscriptstyle (#1)}	
	\else BAD!!						
\fi}

\DeclareMathOperator{\Vor}{Vor}

\newcommand{\closerset}[6][0]{\Vor_{#2,#3}^{(#4)}\parens[#1]{#5,#6}}

\newcommand{\cstarset}[6][0]{ \leftexp{\bigstar}{\Vor}_{#2,#3}^{(#4)}\parens[#1]{#5,#6}}

\newcommand{\fatcstarset}[7][0]{\leftexp{\bigstar}{\Vor}_{#2,#3,#4}^{(#5)}\parens[#1]{#6,#7}}
\newcommand{\fatcstar}[7][0]{\bigstar_{{\scriptscriptstyle (#5)},#3,#4}^{\parens[#1]{#6,#7}}} 

\DeclareMathOperator{\adv}{adv} 
\DeclareMathOperator{\Adv}{Adv} 

\newcommand{\stateadvantage}[4][0]{\Adv_{#3} \parens[#1]{#4}}
\newcommand{\configadvantage}[5][0]{\Adv_{#3} \pair[#1]{#4}{#5}}
\newcommand{\coneadvantage}[4][0]{\Adv_{#3} \parens[#1]{#4}}
\newcommand{\diradvantage}[5][0]{\adv_{#3}^{#4} \parens[#1]{#5}}

\DeclareMathOperator{\fav}{fav} 

\newcommand{\favwedge}[4][0]{\fav_{#2}^+ \pair[#1]{#3}{#4}}
\newcommand{\cfavwedge}[4][0]{\fav_{#2} \pair[#1]{#3}{#4}}

\DeclareMathOperator{\conq}{conq}

\newcommand{\conqset}[5][0]{\conq_{#3} \pair[#1]{#4}{#5}}

\newcommand{\spp}{\epsilon} 
\newcommand{\thp}{\delta} 

\newcommand{\thpadv}[4][0]{\delta_{#2}^{#3} \parens[#1]{#4}}

\newcommand{\propertyfont}[1]{\mathbf{#1}}

\newcommand{\finitespeed}[1]{\ensuremath{(\propertyfont{FPS}_{#1})}} 

\newcommand{\notiesprob}{\ensuremath{(\propertyfont{NTP1})}} 

\newcommand{\expmoment}[1][]{\ensuremath{(\propertyfont{EM}_{#1})}} 
\newcommand{\Lspace}[1]{\ensuremath{(\propertyfont{L}^{#1})}} 
\newcommand{\minimalmoment}[1]{\ensuremath{\parens[2]{\propertyfont{MM}_{#1}}}}

\newcommand{\geodesicsexist}{\ensuremath{(\exists\propertyfont{Geo})}} 
\newcommand{\stargeodesicsexist}{\ensuremath{(\exists\propertyfont{Geo}^*)}} 

\newcommand{\boundedgrowth}{\ensuremath{(\propertyfont{FS})}} 
\newcommand{\normexists}{\ensuremath{(\sim\!\propertyfont{Norm})}} 

\newcommand{\notiesdet}{\ensuremath{(\propertyfont{NT})}} 
\newcommand{\mintimebddgrowth}{\ensuremath{(\propertyfont{MTFS})}} 

\newcommand{\initialsetsymb}{A}
\newcommand{\initialset}[1]{\initialsetsymb_{#1}}
\newcommand{\finalsetsymb}{\mathbf{C}}
\newcommand{\finalset}[1]{\finalsetsymb_{#1}}

\newcommand{\initconfig}{\pair{\initialset{1}}{\initialset{2}}}

\newcommand{\dfinalsetsymb}{\mathbf{D}}
\newcommand{\dfinalset}[1]{\dfinalsetsymb_{#1}}

\DeclareMathOperator{\CoexOp}{Coex}
\newcommand{\coexconfig}[4][0]{\CoexOp_{#2} \pair[#1]{#3}{#4}}
\newcommand{\coexstate}[3][0]{\CoexOp_{#2} \parens[#1]{#3}}

\newcommand{\expmean}[1][0]{%
\ifcase#1
	\beta				
	\or \exprate^{-1}		
\fi}
\newcommand{\exprate}[1][0]{%
\ifcase#1
	\lambda			
	\or \expmean^{-1}	
\fi}

\newcommand{\expnormpair}{\normpair_{\exprate}}

\newcommand{\statespacefont}{\mathsf}
\newcommand{\setaasymb}{\statespacefont{ADV}} 
\newcommand{\setbbsymb}{\statespacefont{CON}} 
\newcommand{\setccsymb}{\statespacefont{ESC}} 
\newcommand{\setddsymb}{\statespacefont{WIN!}} 

\newcommand{\pvbsymb}{\statespacefont{PvB}} 
\newcommand{\pvbset}[2]{\pvbsymb_{#1} \parens{#2}}

\newcommand{\advset}[2]{\setaasymb_{#1}^{#2}}
\newcommand{\advsetr}[3]{\advset{#1}{#2} \parens{#3}}
\newcommand{\coneset}[2][0]{\setbbsymb \parens[#1]{#2}}
\newcommand{\escset}[2][0]{\setccsymb \parens[#1]{#2}}
\newcommand{\winset}{\setddsymb}

\newcommand{\Haggstrom}{H{\"a}ggstr{\"o}m\xspace}

\newcommand{\CD}{Cox:1981aa}
\newcommand{\DHa}{Deijfen:2006aa}
\newcommand{\DHb}{Deijfen:2006ab}
\newcommand{\DHunbounded}{Deijfen:2007aa}

\newcommand{\GM}{Garet:2005aa}
\newcommand{\GMdensity}{Garet:2008aa}

\newcommand{\Hofa}{Hoffman:2005aa}
\newcommand{\Hofb}{Hoffman:2008aa}

\newcommand{\HPa}{Haggstrom:1998aa}
\newcommand{\HPb}{Haggstrom:2000aa}
\newcommand{\HW}{Hammersley:1965aa}

\newcommand{\KesAspects}{Kesten:1986aa}

\newcommand{\PimMultitype}{Pimentel:2007aa}
\newcommand{\PimAsymptotic}{Pimentel:2011aa}

\newcommand{\Ric}{Richardson:1973aa}

\begin{document}

\prelimpages


\Title{A Geometric Perspective on First-Passage Competition}
\Author{Nathaniel D. Blair-Stahn}
\Program{Mathematics}
\Year{2012}

\Chair{Christopher Hoffman}{Professor}{Mathematics}
\Signature{David Wilson}
\Signature{Krzysztof Burdzy}

\copyrightpage

\titlepage

\abstract{
We study the macroscopic geometry of first-passage competition on the integer lattice $\Zd$, with a particular interest in describing the behavior when one species initially occupies the exterior of a cone. First-passage competition is a stochastic process modeling two infections spreading outward from initially occupied disjoint subsets of $\Zd$. Each infecting species transmits its infection at random times from previously infected sites to neighboring uninfected sites. The infection times are governed by species-specific probability distributions, and every vertex of $\Zd$ remains permanently infected by whichever species infects it first.

We introduce a new, simple construction of first-passage competition that works for an arbitrary pair of disjoint starting sets in $\Zd$, and we analogously define a deterministic first-passage competition process in the Euclidean space $\Rd$, providing a formal definition for a model of crystal growth that has previously been studied computationally. We then prove large deviations estimates for the random $\Zd$-process, showing that on large scales it is well-approximated by the deterministic $\Rd$-process, with high probability.
Analyzing the geometry of the deterministic process allows us to identify critical phenomena in the random process when one of the two species initially occupies the entire exterior of a cone and the other species initially occupies a single interior site.
Our results generalize those in a 2007 paper of Deijfen and \Haggstrom, who considered the case where the cone is a half-space. 
Moreover, we use our results about competition in cones to strengthen a 2000 result of \Haggstrom and Pemantle about competition from finite starting configurations.
%
%
%
%
%
%
}
\tableofcontents


\acknowledgments{%
For all the personal support they have given me, I thank my advisor Chris Hoffman, my committee member and coauthor David Wilson, my committee member and professor Krzysztof Burdzy, my professors at the University of Washington and the University of Arizona, my mother Jennifer Blair, my girlfriend Caitlin Jones, my family, the Jones family, my friends, my housemates in Seattle, the cryptography group at Sandia National Laboratories, my fellow graduate students in the Math Department, the Swing Kids at UW, the Lindy Hop and blues dancing communities in Seattle, and the community of ultimate players in Seattle.  I thank Seattle's caf\'{e}s and coffee shops, particularly those in the U District and Capitol Hill, for providing the venues in which I wrote most of this thesis.  I thank the taxpayers of Washington State and the rest of the country for paying my living expenses during most of my time in graduate school, and I thank the voters and politicians who made this possible.  I thank the book \textit{Ishmael} and its sequels, by Daniel Quinn, for enabling me to make sense of the world around me. Finally, I thank the triangle inequality, without which almost none of the proofs in this thesis would have been possible.
}

\dedication{%
\begin{center}
To my father,\\
Ruggles Monroe Stahn,
for nurturing my innate sense of curiosity\\
And to my abuelo,\\
Aubrey Nathan Blair,
for being an inspiration to me and my family
\end{center}
}

\textpages

\part{A Unified Treatment of Random and Deterministic First-Passage Competition}


\chapter{Introduction}
\label{intro2_chap}

%
%
%

First-passage competition is a stochastic process modeling two infections spreading outward from initially occupied disjoint vertex sets in a graph.
This competition model was introduced by \Haggstrom and Pemantle \cite{\HPa} as a generalization of first-passage percolation, which models a single infection spreading throughout the graph.
First-passage percolation is described by the shortest path metric on a graph with random edge weights, and thus it is essentially a model of random geometry.
We will take the underlying graph to be the integer lattice $\Zd$ with $d\ge 2$, and we will assume that the edge weights are independent and identically distributed (\iid). 
We are mainly interested in the two-type competition model, but since this two-type process is constructed from two one-type first-passage percolation processes spreading simultaneously in $\Zd$, we begin by describing both the one-type and two-type processes in a bit more detail.


\section{Overview of First-Passage Percolation}
\label{intro:fpp_overview_sec}

\subsection{Brief Description of First-Passage Percolation (the One-Type Process)}
\label{intro:one_type_description_sec}

The one-type process (first-passage percolation) consists of a single species, say \speciesname{red}, spreading out from an initial set of vertices $A\subseteq \Zd$. To construct the one-type process, we start with a collection of \iid\ nonnegative random variables $\tmeasure = \setof{\ttime(\edge)}_{\edge\in\edges{\Zd}}$, indexed by the edges in the lattice $\Zd$, where two vertices in $\Zd$ are joined by an edge if the Euclidean distance between them is 1. The random variable $\ttime(\edge)$ is called the \textdef{traversal time}\footnote{%
In the literature on first-passage percolation, the random variable $\mttimenew$ is more commonly called the \emph{passage time} of the edge $\edge$. However, the term ``passage time" is also used to refer to the distance $\ptimenew{}{\vect u}{\vect v}$ between two vertices in the induced shortest path metric on the graph. To avoid confusion, I prefer to use different terms for these two concepts. I follow the convention of Kordzakhia and Lalley \cite{Kordzakhia:2005aa} and use the term ``passage time" to mean the {infimum} of ``traversal times" of paths between two sets, and I treat a single edge as a path of length 1; see Sections~\ref{fpc_basics:traversal_meas_sec} and \ref{fpc_basics:passage_time_sec} for the precise general definitions.
}
of the edge $\edge$, and is interpreted as the time it takes \speciesname{red} to cross the edge $\edge$ in either direction. \speciesname{Red} occupies the initial set $A$ at time $t=0$, and then spreads out from $A$ by crossing edges in the graph according to the traversal times. 
More precisely, this means:
\begin{quote}
\emph{All the vertices in $A$ are infected by \speciesname{red} at time $t=0$. If \speciesname{red} infects the vertex $\vect u$ at time $t\ge 0$, and $\edge = \setof{\vect u, \vect u'}$ is an edge of $\Zd$, then \speciesname{red} infects $\vect u'$ from $\vect u$ at time $t+\ttime(\edge)$.}
\end{quote}
If $\vect u'$ was previously infected from a vertex other than $\vect u$, then its state does not change --- $\vect u'$ simply remains infected by \speciesname{red}.
The growth of \speciesname{red} can then be described as the set $\reachedset{\tmeasure}{A}{t}$ of infected vertices reached by time $t$ from the set $A$, using the collection of traversal times $\tmeasure$.

The above construction of first-passage growth gives rise to a random pseudometric $\ptime$ on $\Zd$, defined by $\ptimenew{}{\vect u}{\vect v} \definedas \inf \setbuilder[2]{t\ge 0}{\vect v\in \reachedset{\tmeasure}{\vect u}{t}}$ for $\vect u,\vect v\in\Zd$. The random variable $\ptimenew{}{\vect u}{\vect v}$ is called the \textdef{first passage time} or simply the \textdef{passage time} from $\vect u$ to $\vect v$, and the random pseudometric $\ptmetric$ is a metric if and only if $\ttime(\edge)>0$ almost surely.
The pseudometric $\ptmetric$ coincides with the shortest path (pseudo-)metric induced by the ``edge weights" $\setof{\tmeasure(\edge)}_{\edge\in\edges{\Zd}}$; that is, $\ptimenew{}{\vect u}{\vect v}$ is the minimal (or infimal) weight of a path from $\vect u$ to $\vect v$, where the weight of a path is computed as the sum of the traversal times of its constituent edges.
For a fixed $\vect u\in\Zd$, the set $\reachedset{\tmeasure}{\vect u}{t}$ of vertices reached from $\vect u$ by time $t$ is just a deterministically expanding ball in the random pseudometric $\ptmetric$, and more generally, $\reachedset{\tmeasure}{A}{t}$ is a union of such balls for any $A\subseteq\Zd$.

In Section~\ref{fpc_basics:restricted_sec}, we will take a slightly different, more general approach to constructing the one-type first-passage percolation process, by using the traversal times $\ttime$ to construct an entire family of pseudometrics whose metric balls correspond to a family of ``restricted" one-type processes, growing in subsets of $\Zd$ rather than the whole space. We use this generalized construction to give a new   and simple formal definition of the two-type process in Section~\ref{fpc_basics:two_type_construction_sec}, and in Chapter~\ref{cone_growth_chap} we study the restricted first-passage percolation process  in its own right. 



\subsection{A Brief History of First-Passage Percolation}
\label{intro:fpp_history_sec}

First-passage percolation was first introduced by Hammersley and Welsh \cite{\HW} to model fluid flow through a porous medium. For recent surveys of the subject, see \cite{Howard:2004aa}, \cite{Blair-Stahn:2010aa}, \cite{LaGatta:2011aa}, \cite{Ahlberg:2011ac}, \cite{Grimmett:2012aa}, and for earlier accounts, see \cite{Kesten:1986aa}, \cite{Kesten:1987aa}. Here we highlight some of the main results in first-passage percolation so that the reader can gain a general sense of the subject and its relation to the first-passage competition process described in the next section.


The principal feature of the first-passage growth process $\reachedset{\tmeasure}{A}{t}$ is the \emph{Shape Theorem}, first proved by Richardson \cite{\Ric} in the sense of convergence in probability for exponentially distributed traversal times, and strengthened to the following almost sure version for general traversal time distributions by Cox and Durrett \cite{\CD} in dimension $d=2$ and by Kesten \cite{\KesAspects} for general $d\ge 2$. The Shape Theorem says that, under some mild assumptions on the traversal times $\tmeasure = \setof{\ttime(\edge)}_{\edge\in\edges{\Zd}}$, the the set $\reachedset{\tmeasure}{\vect 0}{t}$ of sites reached by time $t$ from the origin $\vect 0\in\Zd$ looks roughly like the ball of radius $t$ in some norm $\mu$ on $\Rd$ which depends only on the distribution of $\ttime(\edge)$. To state the result, we define for $V\subseteq\Zd$ the \textdef{cube expansion} of $V$ to be the set $\cubify{V} = V+\unitcube{d}\subseteq \Rd$ obtained by placing a unit cube around each point in $V$.

\begin{thm}[Shape Theorem \cite{\CD}, \cite{\KesAspects}]
\label{shape_thm}

Fix an integer $d\ge 2$, and let $\tmeasure = \setof{\ttime(\edge)}_{\edge\in\edges{\Zd}}$ be a collection of \iid\ nonnegative random variables satisfying the following two conditions:
\begin{enumerate}
\item $\Pr \eventthat[2]{\ttime(\edge) = 0} < \pcrit$, where $\pcrit$ is the critical probability for Bernoulli bond percolation on $\Zd$.

\item $\E \min \setof[2]{\labeledmttime{}{1}^d,\dotsc,\labeledmttime{}{2d}^d}<\infty$, where $\setof{\labelededge{}{1},\dotsc,\labelededge{}{2d}}$ is any collection of $2d$ distinct edges in $\edges{\Zd}$.
\end{enumerate}
Then there is a deterministic norm $\mu$ on $\Rd$ such that for any $\epsilon>0$,
\begin{equation}
\label{shape_thm_eqn}
\Pr \eventthat[3]{(1-\epsilon) t\unitball{\mu} \subseteq 
	\cubify[1]{\reachedset{\tmeasure}{\vect 0}{t}}
	\subseteq (1+\epsilon) t\unitball{\mu}
	\text{\ \ for all large } t } = 1,
\end{equation}
where $\unitball{\mu} = \setbuilder[2]{\vect x\in\Rd}{\mu(\vect x)\le 1}$ is the closed unit ball of the norm $\mu$, and $\cubify[2]{\reachedset{\tmeasure}{\vect 0}{t}} = \reachedset{\ttime}{\vect 0}{t} + \brackets[1]{-\frac{1}{2}, \frac{1}{2}}^d$ is the cube-expanded first-passage percolation process started from the origin $\vect 0\in \Zd$ and using the traversal times $\tmeasure$.
\end{thm}

Another way to interpret the Shape Theorem is that large distances in the random pseudometric $\ptmetric$ are well-approximated by distances in the deterministic norm metric corresponding to $\mu$, in a sense we formulate precisely in Appendix~\ref{asymptotic_chap}. The norm $\mu$ in Theorem~\ref{shape_thm} is called the \textdef{shape function} for $\ttime$, and the unit ball $\unitball{\mu}$ is called the \textdef{limit shape} for $\tau$, since the scaled process $\reachedset{\ttime}{\vect 0}{t}/t$ converges to $\unitball{\mu}$ in the sense of \eqref{shape_thm_eqn}.  In general, little is known about the limit shape $\unitball{\mu}$ for a given collection of traversal times $\ttime$, other than the obvious facts that it is a compact convex body in $\Rd$ which must satisfy all the symmetries of $\Zd$ (i.e.\ invariance under reflection or permutation of the coordinate axes); note that the convexity and the fact that $\unitball{\mu}$ is a body (i.e.\ has nonempty interior) follow from the fact that $\mu$ is a norm.
Kesten \cite[Section~8]{Kesten:1986aa} shows that for a large class of traversal time distributions (including exponential), the limit shape $\unitball{\mu}$ is not a Euclidean ball for large $d$, casting doubt on the early conjecture that $\unitball{\mu}$ might be a disc for $d=2$ when $\ttime(\edge)$ is exponentially distributed.
%
Durrett and Liggett \cite{Durrett:1981aa} show that for certain traversal time distributions, the limit shape in $\Z^2$ has flat edges in the diagonal directions but is not the full diamond $\unitball{\lone{2}} = \setbuilder[2]{(x_1,x_2)\in \R^2}{\abs{x_1}+\abs{x_2} \le 1}$. In particular, this occurs if $\ttime(\edge)$ is nondegenerate and attains some nonzero minimum value with probability greater than $p_c^{\text{dir}}(\Z^2)$, where $p_c^{\text{dir}}(\Z^2)$ is the critical value for \emph{directed} Bernoulli bond percolation on $\Z^2$.
Damron and Hochman \cite{Damron:2010aa} prove that there exist traversal time distributions for which the limit shape in $\Z^2$ is not a polygon; previously, even this basic result was not known rigorously for \iid\ traversal times.

The assumption that the collection $\tmeasure =   \setof{\ttime(\edge)}_{\edge\in\edges{\Zd}}$ is \iid\ can be relaxed. 
In particular, Boivin \cite{Boivin:1990aa} proves a Shape Theorem in $\Zd$ for stationary and ergodic traversal times with finite moment of order $d+\epsilon$. \Haggstrom and Meester \cite{Haggstrom:1995aa} prove that any compact convex body that is symmetric with respect to reflection through the origin can arise as the limit shape of some collection of stationary ergodic traversal times. Note that symmetry with respect to permutation of the axes doesn't necessarily hold in the stationary case because the traversal times of edges in different directions can have different distributions.


Some recent papers of interest about first-passage percolation include \cite{Chatterjee:2011aa}, \cite{Auffinger:2011aa}, \cite{Auffinger:2011ab}, \cite{Chatterjee:2009aa}, \cite{Ahlberg:2011aa}, \cite{Ahlberg:2011ab}, \cite{Ahlberg:2011ac}, \cite{Damron:2010aa}, \cite{LaGatta:2010aa}, \cite{LaGatta:2011aa}, \cite{Zhang:2010aa}. 



\section{Overview of First-Passage Competition}
\label{intro:fpc_overview_sec}

\subsection{Brief Description of First-Passage Competition (the Two-Type Process)}
\label{intro:two_type_description_sec}
The two-type first-passage competition process consists of two species, say \speciesname{species 1} and \speciesname{species 2}, competing for space in the lattice $\Zd$. To construct the process, each edge $\edge\in \edges{\Zd}$ gets a \textdef{pair of traversal times} $\ttimepair(\edge) = \parens[2]{\ttime_1(\edge),\ttime_2(\edge)}$, where for $i\in\setof{1,2}$, the nonnegative random variable $\ttime_i(\edge)$ is interpreted as the time it takes \speciesname{species $i$} to cross the edge $\edge$ in either direction. \speciesname{Species $i$} starts on some initial set of vertices $\initialset{i}$ (which we assume is disjoint from $\initialset{3-i}$) and, until it encounters the other species, grows as a first-passage percolation process using its own set of traversal times $\tmeasure_i = \setof{\ttime_i(\edge)}_{\edge\in\edges{\Zd}}$. When the two species come into contact, the interaction between them can be summarized succinctly as follows:
\begin{quote}
\emph{In the two-type model, each vertex in $\Zd$ is conquered by whichever species arrives there first, and it can never be re-conquered.} 
\end{quote}
Thus, if one of the species tries to infect a vertex that has already been conquered, then nothing happens, and once a vertex changes from an uninfected state to either state 1 or state 2 (corresponding to the two species), it remains in that state forever. 
The evolution of the two-type process started from the initial configuration $\initconfig$ and using the pair of traversal times $\tmeasurepair = \pair{\tmeasure_1}{\tmeasure_2}$ can be described by the pair $\reachedset{\tmeasurepair}{\initconfig}{t} = \pair[1]{\reachedset{1}{\initconfig}{t}}{\,\reachedset{2}{\initconfig}{t}}$, where $\reachedset{i}{\initconfig}{t}$ is the set of vertices occupied by \speciesname{species $i$} at time $t$.

Since each vertex will be infected exactly once in the two-type process, each edge can be crossed by at most one of the two species. This means that for each $\edge\in \edges{\Zd}$, at most one of the random variables $\ttime_1(\edge)$ and $\ttime_2(\edge)$ will be used in the construction of the process, and hence it does not matter how these two variables are coupled. However, independence between disjoint edges will be important, as will translation invariance of the entire model, so we assume that the sequence of pairs $\ttimepair = \setof{\ttimepair(\edge)}_{\edge\in\edges{\Zd}}$ is \iid\ (though this assumption can probably be weakened to stationary and ergodic without losing any essential features of the model). We will impose some additional mild restrictions on the distributions of $\ttime_1(\edge)$ and $\ttime_2(\edge)$ in Chapter~\ref{fpc_basics_chap} when we formally construct the two-type process.

When both of the starting sets $A_1$ and $A_2$ are finite, determining which species arrives at a vertex first can be done step-by-step, e.g.\ using a modified version of Dijkstra's algorithm to simultaneously find two shortest path trees in the graph $\Zd$, with edges weighted according to the $\ttime_i$'s (see e.g.\ \cite{\GMdensity}). However, if we allow the starting sets $A_i$ to be infinite, the determination of which vertices are conquered by which species constitutes a somewhat subtle problem, because infinitely many infection events can occur within a finite time interval. We will deal with this issue and provide an explicit construction of the two-type process in Section~\ref{fpc_basics:two_type_construction_sec}.

In \Haggstrom and Pemantle's \cite{\HPa} original definition of the two-type process, $\ttime_1(\edge)$ and $\ttime_2(\edge)$ are taken to be exponentially distributed, i.e.\ $\Pr \eventthat{\ttime_i(\edge)>x} = e^{-\lambda_i x}$ for some $\lambda_1,\lambda_2>0$.
The memoryless property of the exponential distribution implies that in this case, $\reachedset{\tmeasurepair}{\initconfig}{t}$ is a Markov process in which the uninfected vertices of $\Zd$ are infected by \speciesname{species $i$} at rate $\lambda_i$ times the number of their neighbors already in $\reachedset{i}{\initialset{i}}{t}$.
\Haggstrom and Pemantle named this version of first-passage competition the \textdef{two-type Richardson model}, because the Markov version of first-passage percolation is known as Richardson's \cite{\Ric} growth model.
We will study the Markov case in Chapter~\ref{coex_finite_chap}, but in the remainder of the paper, we assume no more than the minimal assumptions described in Chapter~\ref{fpc_basics_chap}, unless explicitly stated otherwise.

\subsection{Questions of Interest in First-Passage Competition}
\label{intro:fpc_questions_sec}

In this section we summarize several key results about first-passage competition and describe the motivation for the present work.
For more in-depth surveys, see \cite{Deijfen:2008aa}, \cite{Blair-Stahn:2010aa}. 
Most of the results in this section are about the two-type Richardson model, i.e.\ the Markov version of first-passage competition in which the traversal times $\ttime_1(\edge)$ and $\ttime_2(\edge)$ are exponentially distributed with rates $\lambda_1$ and $\lambda_2$, respectively.
Note that by time scaling, most questions about the two-type Richardson model depend only on the \emph{ratio} $\lambda = \lambda_2/\lambda_1$, so in the results below it suffices to assume $\lambda_1=1$ and $\lambda_2 = \lambda$; however, we will generally keep the subscripts to make it easier to remember which parameter belongs to which species.

Since its introduction in \cite{\HPa}, the majority of work on first-passage competition has focused on the question of \emph{coexistence} of the two species as time goes to $\infty$. That is, do both species continue growing indefinitely, or does one species
end up surrounded by the other so that it is only able to infect a
finite number of sites? If $\initialset{1}$ and $\initialset{2}$ are two disjoint subsets
of $\Zd$ and $\tmeasurepair = \pair{\tmeasure_1}{\tmeasure_2}$ is the collection of random traversal times used to run the two-type process, we denote by $\coexconfig{\tmeasurepair}{\initialset{1}}{\initialset{2}}$ the event that both species
eventually infect an infinite number of sites when \speciesname{species $i$}
initially occupies the sites in $\initialset{i}$, and we call this event
\textdef{coexistence} or \textdef{mutual unbounded growth} for the
initial configuration $\initconfig$. 
For two-type Richardson model with rates $\lambda_1$ and $\lambda_2$, we will use the notation $\coexconfig{\lambda_1,\lambda_2}{\initialset{1}}{\initialset{2}}$ for the event of coexistence.
We also say that \speciesname{species $i$} \textdef{survives} if it conquers an infinite number of sites, so that $\coexconfig{\tmeasurepair}{\initialset{1}}{\initialset{2}}$ is the event that both species survive from the initial configuration $\initconfig$.

It is easy to see that
$\Pr \parens[2]{\coexconfig{\tmeasurepair}{\initialset{1}}{\initialset{2}}}<1$ unless both of the sets $\initialset{1}$ and $\initialset{2}$ are
already infinite or the distributions of $\cmttime{1}$ and $\cmttime{2}$ are degenerate, so the first nontrivial question to ask is whether $\Pr \parens[2]{\coexconfig{\tmeasurepair}{\initialset{1}}{\initialset{2}}} >0$.
Clearly coexistence is impossible if one of the initial sets
$\initialset{i}$ surrounds the other set $\initialset{j}$, i.e.\ if every infinite simple path starting in $\initialset{j}$ intersects $\initialset{i}$. Following the terminology of Deijfen and \Haggstrom \cite{\DHa}, we say that
the initial configuration $\initconfig$ is \textdef{fertile} if neither initial set surrounds
the other. Deijfen and \Haggstrom show that, as long as the
initial configuration of two-type Richardson process is finite and fertile, the
choice of configuration is irrelevant to the question of whether
coexistence has positive probability:

\begin{thm}[Irrelevance of starting configuration for possibility of coexistence \cite{\DHa}]
\label{dh_irrelevance_thm}
\newcommand{\Aone}{\initialset{1}}
\newcommand{\Atwo}{\initialset{2}}

\newcommand{\Aonep}{\Aone'}
\newcommand{\Atwop}{\Aone'}
\newcommand{\initconfigp}{\pair{\Aonep}{\Atwop}}

If $\initconfig$ and $\initconfigp$ are two fertile pairs of disjoint
finite sets in $\Zd$, then for the two-type Richardson model with fixed growth rates $\lambda_1$
and $\lambda_2$,
\[
\Pr\parens[2]{\coexconfig{\lambda_1,\lambda_2}{\Aone}{\Atwo}} >0
\iff
\Pr\parens[2]{\coexconfig{\lambda_1,\lambda_2}{\Aonep}{\Atwop}}>0.
\]
\end{thm}

Garet and Marchand \cite[Lemma~5.1]{\GM} prove an analogue of Theorem~\ref{dh_irrelevance_thm} in the case where the two species share a single collection of stationary traversal times.
Both results are proved by modifying finitely many edge-traversal times in order to couple the two processes to have comparable growth outside some bounded region.
In the present paper, we consider the situation where $\initialset{1}$ is finite but
the set $\initialset{2}$ is already infinite, in which case the question of coexistence becomes the question of survival of \speciesname{species 1}. In Section~\ref{fpc_basics:unbounded_survival_sec} we use a similar finite modification argument to prove an analogue of Theorem~\ref{dh_irrelevance_thm} in this setting; see Theorem~\ref{intro_irrelevance_thm} in Section~\ref{intro:ch2_sec} below.


Intuitively, if the growth rates $\lambda_1$ and $\lambda_2$ in the two-type Richardson model are equal, we might expect that coexistence from finite fertile starting configurations occurs with positive probability since neither species has an inherent advantage over the other. On the other hand, if the growth rates are different, say $\lambda_1> \lambda_2$, then unless species 2 gets lucky and surrounds species 1 relatively quickly, species 1 is likely to overtake species 2 by virtue of its superior speed, making coexistence implausible. \Haggstrom and Pemantle conjectured that this is indeed how the process behaves:

\begin{conj}[Coexistence occurs precisely when speeds are equal \cite{\HPa}]
\label{hp_coexistence_conj}
In the two-type Richardson model on $\Zd$ with rates $\lambda_1$ and $\lambda_2$, if the initial configuration $\initconfig$ is finite and fertile, then coexistence has positive probability if and only if $\lambda_1=\lambda_2$.
\end{conj}

In fact, \Haggstrom and Pemantle \cite{\HPa} proved the ``if" direction of Conjecture~\ref{hp_coexistence_conj} in dimension $d=2$, showing that coexistence is possible when the two species have the same growth rate.
Subsequently, Garet and Marchand \cite{\GM} and Hoffman \cite{\Hofa} independently generalized this result to show that in any dimension $d\ge 2$, coexistence has positive probability for a wide class of ergodic stationary traversal times, when the collections of traversal times $\ttime_1$ and $\ttime_2$ for the two species have the same distribution. We state here the main coexistence result for the two-type Richardson model.


\begin{thm}[Coexistence is possible with equal speeds \cite{Haggstrom:1998aa}, \cite{Garet:2005aa}, \cite{Hoffman:2005aa}]
\label{richardson_coex_thm}
In the two-type Richardson model on $\Zd$ with growth rates $\lambda_1= \lambda_2 =\lambda$,
if $\initconfig$ is any finite, fertile initial configuration, then
$
\Pr \parens[2]{\coexconfig{\lambda,\lambda}{\initialset{1}}{\initialset{2}}} >0.
$
\end{thm}

Additionally, Hoffman \cite{\Hofb} extends the results of \cite{\Hofa} to show that coexistence of four species is possible for a large class of ergodic stationary traversal times, and moreover, that if the limit shape is not a polygon, then coexistence of arbitrarily many species is possible. Combined with the results of \cite{Damron:2010aa}, it follows that there exists an \iid\ family of traversal times for which coexistence of $k$ species has positive probability for every $k\ge 1$.

%
%
%

In contrast to the progress that has been made in studying equally powerful infections, the situation when the two species have different growth rates remains stubbornly unresolved. The best result to date is the following weakened version of the ``only if" direction of Conjecture~\ref{hp_coexistence_conj}, proved by \Haggstrom and Pemantle \cite{\HPb}.

\begin{thm}[Coexistence impossible for all but countably many speeds \cite{\HPb}]
\label{hp_countable_parameters_thm}

In the two-type Richardson model on $\Zd$ with a finite initial configuration $\initconfig$, if the rate $\lambda_1$ is fixed, then
$
\Pr \parens[2]{\coexconfig{\lambda_1,\lambda_2}{\initialset{1}}{\initialset{2}}} = 0
$
for all but countably many values of $\lambda_2$.
\end{thm}


Although Theorem~\ref{hp_countable_parameters_thm} shows that coexistence is impossible for almost all speed ratios $\lambda = \lambda_2/\lambda_1$, we are presently unable to identify a single ratio $\lambda$ where coexistence does not occur.
Theorem~\ref{hp_countable_parameters_thm} has been extended beyond the two-type Richardson model by Garet and Marchand \cite{\GMdensity} to general first-passage competition with stochastically comparable traversal time distributions when one species' times depend on a continuous parameter, and by Deijfen, \Haggstrom, and Bagley \cite{Deijfen:2004aa} to a two-type continuum ``outburst" process in $\Rd$ analogous to first-passage competition.


At first glance, it may seem strange that Theorem~\ref{hp_countable_parameters_thm}
has not been extended to include all values of $\lambda_2\ne
\lambda_1$. Intuitively, we expect that $\Pr \parens[2]{\coexconfig{\lambda_1,\lambda_2}{\initialset{1}}{\initialset{2}}}$ should decrease as $\lambda=\lambda_2/\lambda_1$ moves farther
away from 1. Since Theorem~\ref{hp_countable_parameters_thm} implies that we can
choose $\lambda$ arbitrarily close to 1 such that $\Pr \parens[2]{\coexconfig{\lambda_1,\lambda_2}{\initialset{1}}{\initialset{2}}} =0$, such monotonicity would imply that coexistence is
impossible for all $\lambda\ne 1$. However, it is not obvious how to
prove that the probability of coexistence is monotone in $\lambda$.
In fact, although this monotonicity property is certainly plausible
for the integer lattice $\Zd$, Deijfen and H\"aggstr\"om \cite{\DHb}
have shown that there are other (highly non-symmetric) graphs where monotonicity of coexistence probabilities does not hold.

\Haggstrom and Pemantle's proof of Theorem~\ref{hp_countable_parameters_thm} relies critically on the following result, which shows that if coexistence occurs, then the total region occupied by both species
satisfies the same Shape Theorem as the slow species.
To state the result, let
\begin{equation}
\label{intro_both_reached_def}
\bothreachedset{\initconfig}{t} = \reachedset{1}{\initconfig}{t} \cup \reachedset{2}{\initconfig}{t}
\end{equation}
denote the region in $\Zd$ occupied by both species at time $t$ in a two-type process started from $\initconfig$, and define the \textdef{shape function} $\mu_i$ for \speciesname{species~$i$} to be the norm on $\Rd$ given by Theorem~\ref{shape_thm} for the one-type process using \speciesname{species~$i$}'s traversal times $\tmeasure_i$.

\begin{thm}[Coexistence implies slow growth {\cite[Lemma~5.2]{\HPb}}]
\label{hp_slow_growth_thm}
%
Consider the two-type Richardson model on $\Zd$ with rates $\lambda_1$ and $\lambda_2$ and corresponding shape functions $\mu_1$ and $\mu_2$ for the two species. If $\lambda_1\ge \lambda_2$ and $\initconfig$ is a finite starting configuration in $\Zd$, then almost surely on the event  $\coexconfig{\lambda_1,\lambda_2}{\initialset{1}}{\initialset{2}}$, for any $\epsilon\in (0,1]$ there exists $t_\epsilon<\infty$ such that
\[
(1-\epsilon) t \unitball{\mu_2} \subseteq
\cubify[1]{\bothreachedset{\initconfig}{t}}
\subseteq (1+\epsilon) t \unitball{\mu_2}
\quad\text{for all $t\ge t_\epsilon$},
\]
where $\unitball{\mu_2}$ is the unit ball for species~2's shape function, $\bothreachedset{\initconfig}{t}$ is defined by \eqref{intro_both_reached_def}, and $\cubify{V} = V+\unitcube{d}$ for $V\subseteq \Zd$.
\end{thm}

To prove Theorem~\ref{hp_countable_parameters_thm} from Theorem~\ref{hp_slow_growth_thm} , \Haggstrom and Pemantle employ
a monotone coupling over all $\lambda_2\in (0,\lambda_1]$ of the two-type processes $\reachedset{\lambda_1,\lambda_2}{\initconfig}{t}$ with rates $(\lambda_1,\lambda_2)$, where the monotonicity of the coupling is with respect to $\lambda_2$ and the ordering on processes we describe in Section~\ref{fpc_basics:monotonicity_sec}; for more details, see \cite{Deijfen:2008aa}, \cite{Blair-Stahn:2010aa}. 
The lower bound in Theorem~\ref{hp_slow_growth_thm}, i.e.\ the fact that the overall growth is at least as fast as the slow species, is essentially trivial; see e.g.\ \cite[Lemma~3.2]{\HPb} or \cite[Lemma~1.1]{\GMdensity}. The hard part of Theorem~\ref{hp_slow_growth_thm} is proving the upper bound, i.e.\ showing that if coexistence occurs, then the fast species is actually constrained to grow at the same 
macroscopic 
rate as its slower competitor. The idea is to show that if the fast species gets a big enough head start, then with high probability it is able to conquer a family of expanding $\mu_1$-spherical shells that spiral around and cut off the slow species.
\Haggstrom and Pemantle prove this key result in  \cite[Proposition~2.2]{\HPb}; once this is established, it follows from the strong Markov property that if species~1 reaches outside the ball $(1+\epsilon) t \unitball{\mu_2}$ infinitely often, then species~2 will lose.
The key result \cite[Proposition~2.2]{\HPb}, which we will state more formally in Chapter~\ref{random_fpc_chap} as Proposition~\ref{hp_point_ball_prop} (p.~\pageref{hp_point_ball_prop}), is essentially a stochastic version of a purely geometric result describing the growth of an analogous deterministic competition process started with the slow species initially occupying a ball and the fast species initially occupying a single point on the boundary of the ball. We describe this deterministic process in Section~\ref{intro:determ_crystal_sec} below, and in Section~\ref{deterministic:teardrop_sec}, we will see that from this ``point vs.\ ball" starting configuration for the deterministic process in $\R^2$, the slow species is only able to conquer a bounded ``heart-shaped" region enclosed by two partial logarithmic spirals, which is the same geometry appearing in the proof of \cite[Proposition~2.2]{\HPb}.

One of the main motivations for the present work was to prove a strengthened version of Theorem~\ref{hp_slow_growth_thm}, putting additional restrictions on the growth of the fast species when coexistence occurs. Our main result in this direction is Theorem~\ref{intro_conv_hull_thm} in Section~\ref{intro:ch6_sec} below.
In the same vein as Theorem~\ref{hp_slow_growth_thm} and Theorem~\ref{intro_conv_hull_thm} are the results proved by Garet and Marchand in \cite{\GMdensity}.
The main result of \cite{Garet:2008aa} shows that almost surely on the event of coexistence, the fast species cannot occupy too large a fraction of the boundary of the growing region, and in dimension $d=2$, it is almost surely impossible for the both species to finally occupy a set of positive density in the plane.

In addition to strengthening Theorem~\ref{hp_slow_growth_thm}, the other main focus of this thesis is to consider the question of coexistence when one species initially occupies an unbounded set and the other species starts at a single point. This question was first addressed by Deijfen and \Haggstrom \cite{\DHunbounded}, who considered the case where $\initialset{1}= \setof{\vect 0}$ and $\initialset{2}$ is either the half-space $\Hspace = \setbuilder[2]{\vect v\in\Zd}{v_1<0}$ or the half-line $L = \setbuilder[2]{\vect v \in\Zd}{\text{$v_1<0$ and $v_j = 0$ for all $j\ne 1$}}$, where we write $\vect v = (v_1,\dotsc,v_d)$ for a point in $\Zd$. Their main result is the following.

\begin{thm}[Survival from ``point vs.\ half-space" or ``point vs.\ half-line" \cite{\DHunbounded}]
\label{dh_unbounded_survival_thm}
\newcommand{\zdhalfspace}{\Hspace}
\newcommand{\zdhalfline}{L}

Let $\zdhalfspace = \setbuilder[2]{\vect v\in\Zd}{v_1<0}$ and $\zdhalfline = \setbuilder[2]{\vect v \in\Zd}{\text{$v_1<0$ and $v_j = 0$ for all $j\ne 1$}}$. Then for the two-type Richardson model on $\Zd$ with rates $\lambda_1$ and $\lambda_2$,
\begin{enumerate}
\item \label{dh_unbounded_survival:Hpart}
$\Pr \parens[2]{\coexconfig{\lambda_1,\lambda_2}{\vect 0}{\zdhalfspace}}>0$ if and only if $\lambda_1>\lambda_2$.
\item \label{dh_unbounded_survival:Lpart}
$\Pr \parens[2]{\coexconfig{\lambda_1,\lambda_2}{\vect 0}{\zdhalfline}}>0$ if and only if $\lambda_1\ge\lambda_2$.
\end{enumerate}
\end{thm}



In short, Theorem~\ref{dh_unbounded_survival_thm} says that if species~1 is strictly faster, then it is able to survive against either the half-line $L$ or the larger half-space $\Hspace$, but if the two species have the same growth rate, then species~1 can only survive against the smaller set $L$.  Note that, in contrast to the situation for finite initial configurations, Theorem~\ref{dh_unbounded_survival_thm} identifies exactly which speed ratios allow coexistence from the infinite starting configurations $(\vect 0,\Hspace)$ and $(\vect 0,L)$. The main reason this is possible is that while the event of coexistence is non-monotone for finite initial configurations, for semi-infinite configurations such as $(\vect 0,\Hspace)$ and $(\vect 0,L)$, coexistence is reduced to the monotone event that species~1 survives.

We now briefly discuss the proof of Theorem~\ref{dh_unbounded_survival_thm}. First consider the subcritical case $\lambda_1<\lambda_2$, so species~1 is strictly slower than species~2.
The fact that coexistence is impossible
from either configuration in this case follows easily from the key result \cite[Proposition~2.2]{Haggstrom:2000aa} of \Haggstrom and Pemantle, mentioned above. In fact, \cite[Proposition~2.2]{Haggstrom:2000aa} (cf.\ Proposition~\ref{hp_point_ball_prop}) implies that if $\lambda_1<\lambda_2$ and $\initialset{1}$ is finite, then species~1 cannot survive against any infinite set $\initialset{2}$.

Now consider the critical case $\lambda_1=\lambda_2$. When the two species have equal growth rates, the two-type process can be reduced to a one-type process by using the same collection of traversal times $\ttime_1=\ttime_2=\ttime$ for the two species, simplifying the construction of the process and making it possible to exploit various symmetries.
The proof that coexistence is possible when $\lambda_1=\lambda_2$ from the configuration $(\vect 0,L)$ is similar in spirit to the coexistence proofs appearing in \cite{Haggstrom:1998aa} and \cite{Garet:2005aa};
in fact, Deijfen and \Haggstrom use their result to give a simple proof of coexistence for finite starting configurations in the case $\lambda_1=\lambda_2$ \cite[Theorem~6.1]{Deijfen:2007aa}. The most difficult part of Theorem~\ref{dh_unbounded_survival_thm} is showing that coexistence is impossible from the configuration $(\vect 0,\Hspace)$ when $\lambda_1=\lambda_2$. The basic idea is to consider the family of hyperplanes $\Hplane_n = \setbuilder[2]{\vect v\in \Zd}{v_1=n}$ and argue that the number of vertices in $\Hplane_n$ infected by species~1 converges almost surely to 0 as $n\to\infty$. This is accomplished by first using a symmetry argument to conclude that the expected number of vertices in $\Hplane_n$ infected from the origin equals 1 for any $n$ (and in particular is bounded as $n\to\infty$), then introducing a filtration and applying Levy's 0-1 law. We mention that Deijfen and \Haggstrom's result was in fact stated for the initial configuration $\pair[2]{\vect 0}{\Hplane_0 \setminus \setof{\vect 0}}$ rather than the configuration $(\vect 0,\Hspace)$ above. However, this is equivalent to the version stated in Theorem~\ref{dh_unbounded_survival_thm}  because of Theorem~\ref{intro_irrelevance_thm} below about the irrelevance of the initial location of species~1 when species~2 starts on an infinite set. In particular, note that the combination of Theorem~\ref{intro_irrelevance_thm} with Theorem~\ref{dh_unbounded_survival_thm} shows that in the critical case $\lambda_1=\lambda_2$, species~1 cannot survive when species~2 initially occupies a hyperplane (or half-space), no matter how far away from the hyperplane species~1 starts.

The remaining case in Theorem~\ref{dh_unbounded_survival_thm} is the supercritical case $\lambda_1>\lambda_2$, in which species~1 is faster. Note that once it is established that coexistence is possible from the initial configuration $(\vect 0,\Hspace)$ when $\lambda_1>\lambda_2$, the same result follows trivially for the configuration $(\vect 0,L)$ since $L\subseteq \Hspace$. This leaves the proof of the ``if" direction of Part~\ref{dh_unbounded_survival:Hpart} of Theorem~\ref{dh_unbounded_survival_thm}. Deijfen and \Haggstrom prove this result \cite[Proposition~3.1]{Deijfen:2007aa} by combining a Shape Theorem for a one-type process started from $\Hspace$ (which follows from large deviations estimates for first-passage growth) with a Shape Theorem for first-passage percolation restricted to a half-cylinder in the direction of the first coordinate axis.
The idea is to show that if $\lambda_1>\lambda_2$, then there is a positive probability that species~1 gets a big enough head start over species~2 that it is able to take over the entire half-cylinder without interference. The proof of this result was the main inspiration for several of the results in this thesis. In particular, in Chapter~\ref{random_fpc_chap} we use the same basic idea to generalize the results of Theorem~\ref{dh_unbounded_survival_thm} for non-critical speed ratios to the case where $\initialset{1}$ is a point and $\initialset{2}$ is the complement of a cone-shaped region in $\Zd$. See Section~\ref{intro:ch5_sec} below.

One challenge in writing this thesis was to formulate an explicit construction of first-passage competition that works for infinite starting configurations when the two species' collections of traversal times ($\tmeasure_1$ and $\tmeasure_2$) have different distributions.
Deijfen and \Haggstrom \cite{\DHunbounded} study the two-type Richardson model with unbounded starting configurations but do not explain how to construct the process for such configurations.
Garet and Marchand \cite{\GMdensity} give an explicit algorithmic construction for general traversal times, but it only works for finite starting configurations. The papers \cite{\GM} and \cite{\Hofa} give constructions that work for any starting configuration, but only when the two species share a common set of traversal times $\tmeasure$. 
In this case, the regions finally conquered by the two species are just the \textdef{Voronoi cells}  for the initial sets $\initialset{1}$ and $\initialset{2}$ with respect to the $\tmeasure$-induced shortest-path pseudometric $\ptmetric$ on $\Zd$. That is, the region conquered by species~1 is just the set of vertices in $\Zd$ that are $\ptmetric$-closer to the set $\initialset{1}$ than to the set $\initialset{2}$, and similarly for species~2. Once the finally conquered sets are identified, it is easy to define the time-evolution of each species as a restricted first-passage percolation process.

In Chapter~\ref{fpc_basics_chap} we develop a variant of this Voronoi-cell construction of the two-type process that works for infinite starting configurations and a general pair of collections of traversal times.
Moreover, our construction generalizes naturally to any pair or tuple of metrics defined on a common space.\footnote{%
More precisely, the two-type process is constructed using a tuple of \emph{length structures}, a concept from metric geometry (cf.\ Burago et al.\ \cite{Burago:2001aa} or Gromov \cite{Gromov:1999aa}) that we will discuss briefly in Chapter~\ref{fpc_basics_chap}. One way to obtain these length structures is to start with metrics, but for first-passage competition on a graph, the fundamental objects are the length structures, not the metrics, which contain less information.
} 
In particular, we are able to explicitly construct the limiting deterministic $\Rd$-process in Chapter~\ref{deterministic_chap}. The resulting deterministic process is a generalization of ``multiplicatively weighted crystal Voronoi diagrams," which have been studied in computational geometry in relation to crystal growth, ecology, path-finding algorithms, and other applications, but never formally defined except in an ``operational" or ``computational" sense. We discuss this deterministic process further in Section~\ref{intro:determ_crystal_sec} below.

First-passage competition and related models have been studied by a number of other authors, including \cite{\PimMultitype}, \cite{\PimAsymptotic}, \cite{Antunovic:2011aa}, \cite{Deijfen:2004aa}, \cite{Deijfen:2004ab}, \cite{Kordzakhia:2005aa}, \cite{Gouere:2007aa}, \cite{Hoffman:2008aa}, \cite{\GMdensity}.

\subsection{Deterministic First-Passage Competition in $\Rd$ and Connections with Crystal-Growth Voronoi Diagrams}
\label{intro:determ_crystal_sec}

%

In Chapter~\ref{deterministic_chap} we construct deterministic first-passage competition in $\Rd$ in a manner that parallels the construction of the two-type competition process in $\Zd$ from Chapter~\ref{fpc_basics_chap}. The purpose of the deterministic process is to provide an approximation to the random process on large scales, via large deviations estimates for the Shape Theorem.
Instead of using a collection of edge traversal times $\tmeasurepair = \pair{\tmeasure_1}{\tmeasure_2}$ to run the process, the growth of the deterministic process is governed by a pair of norms $\normpair = \pair{\mu_1}{\mu_2}$ on $\Rd$. If $\vect x,\vect y\in\Rd$, the distance $\mu_i(\vect y-\vect x)$ represents the time it takes \speciesname{species $i$} to travel in a straight line from $\vect x$ to $\vect y$.
Thus, in the absence of the other species, the set of points \speciesname{species $i$} reaches from the point $\vect x$ by time $t$
is just $\reacheddset{\mu_i}{\vect x}{t} = \setbuilder[2]{\vect y\in\Rd}{\mu_i(\vect y-\vect x)\le t}$, i.e.\ the closed $\mu_i$-ball of radius $t$ centered at $\vect x$.
When both species are present, then just like in the $\Zd$-process, the two species are not allowed to occupy the same space (except along the common boundary between the two species' conquered regions), and whichever species arrives at a point first permanently occupies that point thereafter. The time it takes \speciesname{species $i$} to reach a point within its conquered region is computed using a shortest path lying entirely within this region, where the path length is measured with respect to the norm $\mu_i$.

%



%

%

One can imagine this deterministic competition process as a model of crystal growth. Suppose two types of crystals start on disjoint ``seed sets" $\initialset{1}$ and $\initialset{2}$ in $\Rd$, and each type of crystal then begins growing from its seed $\initialset{i}$ with deterministic speed dictated by the norm $\mu_i$.
The two crystals will eventually cover the whole space $\Rd$, partitioning it into two disjoint connected regions, each containing one of the seed sets. If one type of crystal grows faster than the other in some direction, then as the faster crystal grows it will ``wrap around" the region already occupied by the slower crystal.

In fact, Schaudt and Drysdale \cite{Schaudt:1991aa} have already introduced this growth model for an arbitrary finite number of crystals in $\R^2$, in the special case where each crystal's norm is a multiple of the Euclidean norm and each seed $\initialset{i}$ is a single point. The authors call the resulting partition of $\R^2$ the \textdef{multiplicatively weighted crystal-growth Voronoi diagram}, and the analogous definition makes sense in any dimension $d$.
The descriptor ``multiplicatively weighted" refers to the fact that each crystal's norm is obtained by scaling a single norm by some multiplicative constant, meaning that the crystals are allowed to grow at different rates. 
If the multiplicative weights are all equal, so that all the crystals grow at the same rate, then the resulting partition of $\Rd$ is just the ordinary Voronoi diagram induced by the seed sets, in which the seed $\initialset{i}$ captures all the points in $\Rd$ that are closer to $\initialset{i}$ than to $\bigcup_{j\ne i}\initialset{j}$ in the crystals' common norm metric.
On the other hand, if one computes the Voronoi regions for seed sets with different multiplicative weights, without taking into account the ``wrap-around" effect exhibited by physical crystal growth, then one obtains the \emph{multiplicatively weighted Voronoi diagram}, in which some cells may be disconnected. There is also an \emph{additively weighted Voronoi diagram}, obtained by allowing identical crystals to start growing at different times. When we analyze the deterministic first-passage competition process in Chapter~\ref{deterministic_chap}, we will define (in Section~\ref{deterministic:closer_sets_sec}) a type of weighted Voronoi region that generalizes both of these weighted versions, as a first approximation to the crystal-growth Voronoi regions we are ultimately interested in. In Section~\ref{deterministic:euclid_comp_sec} we will explicitly describe the geometry of the crystal-growth Voronoi diagrams for certain starting configurations $\initconfig$, in the simplest case of two crystals with speeds given by multiples of the Euclidean norm.

Crystal-growth Voronoi diagrams were introduced by Schaudt and Drysdale \cite{Schaudt:1991aa} and have been studied further by Schaudt \cite{Schaudt:1992aa} and by Kobayashi and Sugihara \cite{Kobayashi:2002aa}. The papers \cite{Schaudt:1992aa} and \cite{Kobayashi:2002aa} introduce two different algorithms for approximating the diagram, which is a nontrivial problem. The three papers also suggest some possible applications of crystal-growth Voronoi diagrams. For example, Schaudt \cite{Schaudt:1992aa} mentions some applications previously modeled using ordinary Voronoi diagrams, including modeling territories dominated by animals (possibly with additional natural obstacles), or regions of ground eventually covered by, e.g.\ patches of clovers. Kobayashi and Sugihara \cite{Kobayashi:2002aa} apply their algorithm to compute a collision-free path for a robot moving among enemy robots. Kobayashi and Sugihara also introduced a ``generalized crystal Voronoi diagram" which generalizes both the multiplicatively weighted Voronoi diagram and the crystal-growth Voronoi diagram, and can be used to interpolate between the two.

However, all three papers \cite{Schaudt:1991aa}, \cite{Schaudt:1992aa}, \cite{Kobayashi:2002aa} treat the crystal Voronoi regions from a computational perspective, describing algorithms to approximate the regions without ever providing a formal definition of the region being approximated, instead relying on the informal description of crystal growth given above. My contribution to this problem in Chapter~\ref{deterministic_chap} is to provide a concrete definition of the crystal-growth Voronoi cells for an arbitrary pair of norms, demonstrating that the above crystal growth process is a continuum analogue of the first-passage competition process in $\Zd$. As far as I am aware, no such formal definition has been previously proposed.

\section{Main Results and Organization of This Thesis}
\label{intro:main_results_sec}

The preceding sections were intended to familiarize the reader with first-passage percolation and first-passage competition, providing background and motivation for the present work. We now provide a more detailed description of the content of this work.

\subsection{Existence of First-Passage Competition and Irrelevance of Starting Configuration (Chapter~\ref{fpc_basics_chap})}
\label{intro:ch2_sec}

In Chapter~\ref{fpc_basics_chap} we formally describe the probability space for the one-type and two-type processes and introduce all the definitions and elementary results needed to construct the two-type first-passage competition process in $\Zd$, in Definition~\ref{two_type_proc_def} (p.~\pageref{two_type_proc_def}). Our first main result is the following theorem, which is proved in Propositions~\ref{entangled_full_prop} and \ref{entangled_well_defined_prop} (pp.~\pageref{entangled_full_prop}--\pageref{entangled_well_defined_prop}).


\begin{thm}[Existence of a well-defined FPC process on $\Zd$]
\label{intro_existence_thm}
Suppose that for $i\in\setof{1,2}$, the distribution function of $\ttime_i(\edge)$ is continuous and $\E \brackets{\ttime_i(\edge)}^{1/2} <\infty$, and let $\tmeasurepair = \setof[2]{\pair[2]{\ttime_1(\edge)}{\ttime_2(\edge)}}_{\edge\in \edges{\Zd}}$ be i.i.d. Then there almost surely exists a well-defined two-type competition process $\reachedset{\tmeasurepair}{\initconfig}{t}$ that behaves in the manner described in Section~\ref{intro:two_type_description_sec} for all pairs of disjoint starting sets $\initconfig \in \powerset{\Zd}\times\powerset{\Zd}$. In particular, the construction of the two-type process in Definition~\ref{two_type_proc_def} works when one or both starting sets is infinite and when the two species use different collections of traversal times.
\end{thm}

Theorem~\ref{intro_existence_thm} is actually true for a wider class of traversal times. In particular, the assumption $\E \brackets{\ttime_i(\edge)}^{1/2} <\infty$ can be replaced by the weaker moment condition appearing in the Shape Theorem, and the continuity assumption can be relaxed somewhat; the precise assumptions on $\tmeasurepair$ will be given in Chapter~\ref{fpc_basics_chap}, namely Lemma~\ref{two_type_a_s_properties_lem} (p.~\pageref{two_type_a_s_properties_lem}).

The highlight of Chapter~\ref{fpc_basics_chap} is Theorem~\ref{irrelevance_thm} (p.~\pageref{irrelevance_thm}), which is an analogue for infinite starting configurations of the ``irrelevance of starting configuration for the possibility of coexistence" results proved in \cite{\DHa} and \cite{\GM} for finite starting configurations (Theorem~\ref{dh_irrelevance_thm} above). The following result is a corollary of Theorem~\ref{irrelevance_thm}.

\begin{thm}[Irrelevance of starting configuration for the possibility of survival]
\label{intro_irrelevance_thm}
Suppose $\tmeasurepair$ satisfies the hypotheses of Theorem~\ref{intro_existence_thm} and also $\Pr \eventthat{\ttime_1(\edge)<\epsilon}>0$ for all $\epsilon>0$. In a two-type process constructed from the traversal times $\tmeasurepair$,
if species 2 initially occupies a fixed infinite subset $\initialset{2}$ of $\Zd$, and species 1 can survive with positive probability from some finite initial set $\initialset{1}\subseteq \Zd\setminus \initialset{2}$, then species 1 can also survive with positive probability from any other finite initial set $\initialset{1}'$ in the same component of $\Zd\setminus \initialset{2}$.
\end{thm}

The proof of Theorem~\ref{irrelevance_thm} will require several elementary but technical properties of the two-type process (in particular Lemmas~\ref{one_type_two_type_lem}, \ref{monotonicity_lem}, and \ref{markovesque_lem}), which will be developed throughout Chapter~\ref{fpc_basics_chap}.

\subsection{Definition of First-Passage Competition in $\Rd$ and Deterministic Analogues of Stochastic Competition (Chapter~\ref{deterministic_chap})}
\label{intro:ch3_sec}


In Chapter~\ref{deterministic_chap} we will define deterministic first-passage percolation and first-passage competition in $\Rd$. We will construct these one-type and two-type deterministic processes in a manner that parallels the constructions in Chapter~\ref{fpc_basics_chap}, giving a formal definition of the crystal-growth Voronoi cells in Definition~\ref{determ_two_type_proc_def} (p.~\pageref{determ_two_type_proc_def}). The main purpose of the deterministic $\Rd$-processes is to provide approximations to the random $\Zd$-processes on large scales, via large deviations estimates for the Shape Theorem. The tools for making such approximations will be developed in Chapters~\ref{cone_growth_chap} and \ref{random_fpc_chap}.

In addition to providing the framework for defining deterministic first-passage competition in Section~\ref{deterministic:two_type_construction_sec}, Chapter~\ref{deterministic_chap} serves two other functions. In the first half of the chapter we introduce various elements of basic convex geometry that will be needed to analyze both the deterministic and random processes. The latter half of the chapter focuses on describing the behavior of the deterministic competition process for the same class of starting configurations we are interested in for the random two-type process, namely those in which \speciesname{species 1} starts at a single point and \speciesname{species 2} starts on the entire exterior of a cone. In particular, the behavior of the deterministic process described in Sections~\ref{deterministic:surrounded_sec} and \ref{deterministic:cone_competition_sec} serves as a model for the behavior of the stochastic first-passage competition process in Chapter~\ref{random_fpc_chap}. We will describe the main results from Section~\ref{deterministic:cone_competition_sec} in parallel with the main results for Chapter~\ref{random_fpc_chap} below.

\subsection{Large Deviations for Growth in Cones and Other Star Sets (Chapter~\ref{cone_growth_chap})}
\label{intro:ch4_sec}


In Chapter~\ref{cone_growth_chap} we study the (one-type) first-passage percolation process restricted to subgraphs of $\Zd$. We are in particular interested in obtaining lower bounds on the growth of the process in cone-shaped regions in $\Zd$, because our analysis of the two-type process in Chapter~\ref{random_fpc_chap} will focus on competition that takes place entirely within a cone. It will be convenient to describe these cone-shaped subgraphs in terms of subsets of $\Rd$. Namely, for $A\subseteq S\subseteq\Rd$, we will define an ``$S$-restricted" first-passage percolation process $\rreachedset{\tmeasure}{A}{S}{t}$ by restricting the process's growth to the induced subgraph of $\Zd$ containing the vertices closest to $S$, and then using cube expansion to treat this process as a subset of $S$ rather than a subset of the lattice. Our goal will be to obtain large-deviations estimates for the growth of such $S$-restricted processes, for various choices of the restricting set $S$.

Starting with a basic large deviations estimate (Lemma~\ref{GM_large_dev_lem}, p.~\pageref{GM_large_dev_lem}, taken from \cite[Proposition~2.1]{\GMdensity}) for the unrestricted one-type process, we use a sequence of bootstrap arguments to prove large deviations estimates for the one-type process restricted to increasingly general classes of subsets of $\Rd$.
Our main result applies to a certain class of ``thick" star-shaped sets which we call ``$\mu$-stars" because they are defined in terms of the norm $\mu$ corresponding to the one-type process. 
More precisely, for any $\vect z\in\Rd$ and $\thp>0$, we say that $S\subseteq\Rd$ is a \textdef{$(\mu,\thp)$-star at $\vect z$} if for every $\vect x\in S$, the line segment $\segment{\vect z}{\vect x}$ is contained in $S$ (i.e.\ $S$ is star-shaped at $\vect z$) and there is a $\mu$-ball $B$ of radius $\frac{\thp}{1+\thp}\cdot \mu(\vect x-\vect z)$ with $\vect x\in B\subseteq S$. A \textdef{$\mu$-star at $\vect z$} is then any set that is a $(\mu,\thp)$-star at $\vect z$ for some $\thp>0$.
The main result of Chapter~\ref{cone_growth_chap} is Theorem~\ref{cone_segment_growth_thm} (p.~\pageref{cone_segment_growth_thm}), which can be stated in simplified form as follows.


\begin{thm}[Large deviations estimate for growth in thick star-shaped sets]
\label{intro_star_growth_thm}
Let $\tmeasure = \setof{\ttime(\edge)}_{\edge\in\edges{\Zd}}$ be a collection of \iid\ nonnegative random variables satisfying $\Pr \eventthat[2]{\ttime(\edge) = 0} < \pcrit$ and $\E e^{b \ttime(\edge)}<\infty$ for some $b>0$, and let $\mu$ be the shape function for $\tmeasure$ from Theorem~\ref{shape_thm}.
Given $\spp\in (0,1]$, there exist positive constants $\bigcsymb$ and $\smallcsymb$ such that if $\thp\in (0,1]$ and $S$ is a $(\mu,\thp)$-star at $\vect z\in\Rd$, then for any $t_0\ge 0$,
\[
\Pr \eventthat[3]{ \rreachedset{\tmeasure}{\vect z}{S}{t} \supseteq 
\reacheddset[2]{\mu}{\vect z}{(1-\epsilon)t}\cap S \text{ for all } t\ge t_0}
\ge 1-  \tfrac{1}{\thp^d}  \bigcsymb e^{-\smallcsymb \thp t_0},
\]
where $\rreachedset{\tmeasure}{\vect z}{S}{t}$ is the one-type process started from $\vect z$ and restricted to $S$, and for any $r\ge 0$, $\reacheddset{\mu}{\vect z}{r} \definedas \setbuilder[1]{\vect x\in\Rd}{\mu(\vect x-\vect z)\le r}$ is the $\mu$-ball of radius $r$ centered at $\vect z$.
\end{thm}

%
%

For example, note that if $S$ is a $(\mu,\thp)$-star at the origin $\vect 0$, then so is $rS$ for any $r>0$. Thus, if we consider a one-type process started from $\vect 0$ and restricted to $rS$, then taking $t_0$ proportional to $r$ in Theorem~\ref{intro_star_growth_thm} shows that as we scale the picture up by $r$, the probability that the growth of the restricted random process looks like that of the corresponding restricted deterministic process (up to an error of $\spp t$ at time $t$) converges exponentially to 1. Note, however, that the large deviations estimate for growth in $(\mu,\thp)$-stars gets worse as the thickness $\thp$ approaches 0. This is to be expected since the asymptotic growth rate of the $S$-restricted process in Theorem~\ref{intro_star_growth_thm} is bounded below by the growth rate of the corresponding unrestricted process, whereas a first-passage percolation process restricted to, for example, a line (corresponding to $\thp=0$) grows strictly slower than an unrestricted process.


As a special case of Theorem~\ref{intro_star_growth_thm}, we can take the $(\mu,\thp)$-star $S$ to be an infinite cone generated by a $\mu$-ball.
Namely, for $\thp\le 1$, we define a \textdef{$(\mu,\thp)$-cone} at $\vect z= \vect 0$ to be a set of the form $\cones = \bigcup_{h\ge 0} \reacheddset{\mu}{h \unitvec{\dirsymb}}{\thp h}$, where $\unitvec{\dirsymb}\in\Rd$ is a $\mu$-unit vector (i.e.\ $\mu(\unitvec{\dirsymb})=1$) and $\reacheddset{\mu}{h \unitvec{\dirsymb}}{\thp h}$ is the $\mu$-ball of radius $\thp h$ centered at $h \unitvec{\dirsymb}$. Note that by the triangle inequality for $\mu$, the $(\mu,\thp)$-cone $\cones$ is in fact a $(\mu,\thp)$-star according to the definition above, and moreover, $\cones$ is \textdef{scale-invariant}, meaning that $r\cones = \cones$ for $r>0$. Taking $\vect z = \vect 0$ and $S = \cones$ in Theorem~\ref{intro_star_growth_thm} shows that a one-type process started from the origin and restricted to $\cones$ grows asymptotically as fast as an unrestricted process, and moreover we get a large deviations estimate for the corresponding shape theorem in $\cones$. More generally, the same result holds with $S$ equal to any union of $(\mu,\thp)$-cones at $\vect 0$. This large deviations estimate for growth in cones can be compared with the results of Ahlberg \cite{Ahlberg:2011aa}, who proved a shape theorem and large deviations estimates for first-passage growth in a large class of infinite ``tube-shaped" subgraphs of $\Zd$, which may have widths that increase much more slowly than the linearly expanding cones we consider.


\subsection{Competition in Cones (Chapters~\ref{deterministic_chap} and \ref{random_fpc_chap})}
\label{intro:ch5_sec}
{
\newcommand{\critthp}{\thp_c}
\newcommand{\thisthp}{\thp}
\newcommand{\speedparam}{\Lambda}
\newcommand{\critspeed}{\speedparam_c}


In Chapter~\ref{random_fpc_chap} we are interested in describing the two-type competition process in the case where one species starts on the exterior of a cone in $\Rd$, and the other species starts at a single point inside the cone. In particular, when is the species inside the cone able to survive, that is, conquer an unbounded set?

We answer this question in terms of the geometry of the cone by considering the behavior of the deterministic competition process in Chapter~\ref{deterministic_chap}. Our results generalize those of Deijfen and \Haggstrom \cite{Deijfen:2007aa} (i.e.\ Theorem~\ref{dh_unbounded_survival_thm} above) and we will use them to prove Theorem~\ref{intro_conv_hull_thm} below. The main idea is as follows. Suppose that in the deterministic process, the species starting at the single point is faster than the species starting on the infinite exterior of the cone. If we consider some family of cones with a common axis but varying thickness, then there is some critical thickness above which the fast species can survive inside the cone and below which it can't. Conversely, for a fixed cone there is some critical speed for the species starting inside the cone, above which it can survive and below which it can't. For the deterministic process, the critical thickness and critical speed can be computed from the geometry of the cone and the norms $\mu_1$ and $\mu_2$ of the two species. It turns out that the critical values for the random process are the same as those of the deterministic process with $\mu_1$ and $\mu_2$ chosen to be the shape functions for the two species in the random process, i.e.\ the norms appearing in each species' Shape Theorem.

For the present discussion, we can define a cone $\cones$ in $\Rd$ to be either a closed half-space or else the $d$-dimensional analogue of an infinite cone as defined in elementary Euclidean geometry. For example, an infinite pyramid in $\R^3$ is such a cone, and in this case we call the top of the pyramid the \textdef{apex} of the cone, and we call the ray pointing vertically downward from the apex the \textdef{axis} of the cone.
We will call $\cones$ a \textdef{round cone} if its cross-section is a Euclidean ball of dimension $d-1$, in which case we can define the \textdef{thickness} of the cone to be the angle between the axis and any boundary ray of the cone, measured at the apex.
%
%
%
%
%
%
%
Any cone $\cones$ according to the informal definition given here is homeomorphic to a closed half-space in $\Rd$. 
In the following theorems, we also assume that the cone is convex in order to simplify the statements, though we will relax this assumption when we analyze the competition processes in Chapters~\ref{deterministic_chap} and \ref{random_fpc_chap}.
We first state the main results for the deterministic competition process, found in Section~\ref{deterministic:cone_competition_sec}, and then we state the analogous results for the random competition process, found in Section~\ref{random_fpc:cone_competition_sec}.

\begin{thm}[Deterministic first-passage competition in cones]
\label{intro_determ_cone_comp_thm}

Consider a deterministic two-type process run using the pair of norms $\normpair = \pair{\mu_1}{\mu_2}$ on $\Rd$, in which species 1 starts at the origin $\vect 0$, and species 2 starts on the exterior of a convex cone $\cones$ with apex $\vect 0$. 
\begin{enumerate}
\item 
\label{intro_determ_cone:speeds_fixed_part}

Suppose the norms $\pair{\mu_1}{\mu_2}$ of the two species are fixed, and there is some direction in which species 1 is faster than species 2, meaning that there is some $\dirsymb\in\Rd\setminus \setof{\vect 0}$ such that $\mu_1(\dirsymb)<\mu_2(\dirsymb)$. If we vary the thickness $\thisthp$ of a family of round cones $\cones$ with apex $\vect 0$ and axis in the direction of $\dirsymb$, then there is some critical thickness $\critthp>0$ above which species 1 conquers a subcone of $\cones$ and below which it conquers nothing but $\vect 0$.

\item 
\label{intro_determ_cone:cone_fixed_part}

Suppose the cone $\cones$ is fixed, and we vary the speed of species 1 by scaling its norm while holding the norm of species 2 constant.  That is, the norm $\mu_2$ is fixed, and we set $\mu_1 = \speedparam^{-1} \mu$ for some norm $\mu$ and some $\speedparam>0$. Then there is some critical speed factor $\critspeed \in (0,\infty)$ above which species 1 conquers a subcone of $\cones$ and below which it conquers only $\vect 0$.
\end{enumerate}
\end{thm}

Theorem~\ref{intro_determ_cone_comp_thm} follows from Proposition~\ref{conquering_from_apex_prop} (p.~\pageref{conquering_from_apex_prop}). Here is the corresponding statement for the random process, which follows from Theorem~\ref{wide_narrow_survival_thm} (p.~\pageref{wide_narrow_survival_thm}).

\begin{thm}[Stochastic first-passage competition in cones]
\label{intro_random_cone_comp_thm}

Let $\tmeasurepair = \pair{\ttime_1}{\ttime_2}$ be a collection of \iid\ random traversal time pairs on $\edges{\Zd}$ satisfying the conditions of Theorem~\ref{intro_irrelevance_thm}
and also
$\E e^{b\ttime_i(\edge)} <\infty$ for some $b>0$ and $i\in \setof{1,2}$. Let $\mu_i$ be the shape function (norm) for $\ttime_i$ from the Shape Theorem~\ref{shape_thm},
and consider a random two-type process run with the traversal times $\tmeasurepair$, in which species 1 starts at the origin $\vect 0$, and species 2 starts on the exterior of a convex cone $\cones$ with apex $\vect 0$.
\begin{enumerate}
\item Let the distributions of $\ttime_1(\edge)$ and $\ttime_2(\edge)$ be fixed. Suppose there is some $\dirsymb\in\Rd\setminus \setof{\vect 0}$ such that $\mu_1(\dirsymb)<\mu_2(\dirsymb)$, and suppose we vary the thickness $\thisthp$ of a family of round cones $\cones$ with apex $\vect 0$ and axis in the direction of $\dirsymb$, as in Part~\ref{intro_determ_cone:speeds_fixed_part} of Theorem~\ref{intro_determ_cone_comp_thm}. If $\critthp$ is the corresponding critical thickness for the norm pair $\pair{\mu_1}{\mu_2}$, then for $\thisthp>\critthp$, species 1 conquers a subcone of $\cones$ with positive probability, and for $\thisthp<\critthp$, species 1 almost surely conquers only a bounded set.
%
%
%

\item Suppose the cone $\cones$ is fixed, and we vary the speed of species 1 by scaling its traversal times $\ttime_1(\edge)$ by some constant multiple $\speedparam^{-1}$ while fixing the distribution of species 2's traversal times $\ttime_2(\edge)$.  Then there is some critical speed factor $\critspeed \in (0,\infty)$ above which species 1 can survive in $\cones$ and below which it can't.
\end{enumerate}
\end{thm}

%
%
%
%
%

}

\subsection{Coexistence from Finite Starting Configurations (Chapter~\ref{coex_finite_chap})}
\label{intro:ch6_sec}

In Chapter~\ref{coex_finite_chap} we consider Markov first-passage competition (the two-type Richardson model) started from finite initial configurations. 
The main result is Theorem~\ref{conv_hull_thm} (p.~\pageref{conv_hull_thm}),
which has the same flavor as \Haggstrom and Pemantle's Theorem~\ref{hp_slow_growth_thm} above, putting restrictions on the growth of the fast species when coexistence occurs.
The following is a restatement of Theorem~\ref{conv_hull_thm}.

\begin{thm}[Coexistence implies fast species can't touch support hyperplanes]
\label{intro_conv_hull_thm}
\newcommand{\chevent}{H}
%
Consider the two-type Richardson model with rates $\lambda_1$ and $\lambda_2$, and suppose $\lambda_1>\lambda_2$. Let $\initconfig$ be a finite initial configuration, and let $\chevent$ be the event that in the process started from $\initconfig$, there exists a sequence of times $\setof{t_n}_{n\in\N}$ with $t_n\to\infty$ such that at each time $t_n$, \speciesname{species~1} occupies some vertex on the boundary of the convex hull of the set $\bothreachedset{\initconfig}{t_n}$. Then $\Pr \parens[2]{\coexconfig{\lambda_1,\lambda_2}{\initialset{1}}{\initialset{2}} \cap \chevent} = 0$.
\end{thm}


It is not hard to see that Theorem~\ref{intro_conv_hull_thm} actually implies Theorem~\ref{hp_slow_growth_thm}, as we will explain in Chapter~\ref{coex_finite_chap}. In fact, the proof of Theorem~\ref{intro_conv_hull_thm} is a bootstrap argument, employing Theorem~\ref{hp_slow_growth_thm} at a crucial step to obtain the stronger result. The idea of the proof is to combine the results about competition in cones from Section~\ref{intro:ch5_sec} and growth in cones from Section~\ref{intro:ch4_sec} with Theorem~\ref{hp_slow_growth_thm} and the strong Markov property. In particular, if species 1 is faster than species 2, then Theorem~\ref{intro_random_cone_comp_thm} implies that species 1 can conquer a cone with positive probability if it starts at a single point on the boundary hyperplane of any half-space initially occupied by species~2.
Thus, if the starting configuration is finite and species 1 ever reaches a support hyperplane of the total occupied region, then species~1 has a positive probability of conquering a cone, because at any such time, species~1 is on the boundary of a half-space containing species 2. Now, if species~1 conquers a cone, the results of Section~\ref{intro:ch4_sec} imply that the growth of species 1 in this conquered cone is asymptotically as fast as its unrestricted growth. On the other hand, since Theorem~\ref{hp_slow_growth_thm} implies that on the event of coexistence the overall growth rate must be the same as that of the slow species, we see that if species 1 manages to conquer some cone, then coexistence cannot occur because the overall growth would be too fast. To summarize, any time  species~1 reaches a support hyperplane of the occupied region, it has a positive probability of outrunning and, by Theorem~\ref{hp_slow_growth_thm}, surrounding species 2, preventing coexistence from occurring. Applying the strong Markov property at the times when species~1 reaches a support hyperplane will show that there can only be finitely many such times if coexistence occurs.

\section{Basic Definitions and Notation}

In this section we summarize the standard notation and basic definitions used throughout the paper for easy reference.
%

\subsection{Sets, Functions, and General Notation}

$Y^X= \setof{\text{all functions from $X$ to $Y$}}$; $2^X=$ power set of $X$; $\restrict{f}{A} =$ the restriction of the function $f$ to $A$; pairwise max $\join$ and min $\meet$; $\abs{\,\cdot\,}$ denotes absolute value or cardinality; $f(t+)$ and $f(t-)$ denote the right-hand and left-hand limits of the function $f$ at $t$; $A\definedas B$ or $B\defines A$ means the definition of $A$ is $B$.

\subsection{Probability}

$\Pr =$ probability measure; $\E=$ expectation operator; $\ind{A} =$ indicator function of $A$; $\eventthat{\ldots} = $ ``the event that\ldots"; $\law$ denotes the law (distribution) of a random element; \iid\ means ``independent and identically distributed"; $\eqd$ means equal in distribution; ``$X\sim\ldots$" indicates that the random variable $X$ is distributed according to a specified distribution; $E^\complement=$ the complement of the event $E$; product $\otimes$ for $\sigma$-algebras or measures, e.g.\ $\nu^{\otimes n}$; convolution $*$ of measures, e.g.\ $\nu^{*m}$; if $\setof{A_n}_{n\ge 0}$ is a sequence of events, ``i.o."\ and ``ev."\ stand for ``infinitely often" and ``eventually," i.e.\ $\eventthat{A_n \text{ i.o.}} \definedas \bigcap_{m\ge 0} \bigcup_{n\ge m} A_n$ and $\eventthat{A_n \text{ ev.}} \definedas \bigcup_{m\ge 0} \bigcap_{n\ge m} A_n$.

\subsection{Euclidean and Lattice Geometry}
\label{intro:euc_lat_geo_sec}
Let $d$ be a positive integer. The $d$-dimensional \textdef{Eucliedan space} $\Rd$ is the set of all $d$-tuples $\vect x= (x_1,\dotsc, x_d)$ with $x_j\in \R$ (the real numbers) for all $j$. We reserve bold type lower case letters (e.g.\ $\vect x,\vect y,\vect a,\vect u,\vect v$) for points in Euclidean space. We write $\vect 0 = (0,\dotsc, 0)$ for the \textdef{origin} in $\Rd$ and $\basisvec{1},\dotsc,\basisvec{d}$ for the \textdef{standard basis vectors} in $\Rd$, i.e.\ $\basisvec{1} = (1,0,0,\dotsc,0)$, $\basisvec{2} = (0,1,0,\dotsc,0)$, etc. The $d$-dimensional \textdef{integer lattice} $\Zd$ is the set of $\vect u\in \Rd$ with each $u_j$ an integer, treated as a graph in which two vertices are adjacent if the Euclidean distance between them is one. The $d$-dimensional \textdef{upper half-space} is $\halfspace[d] = \setbuilder{\vect x\in\Rd}{x_d>0}$. We write $\R_+$ or $[0,\infty)$ for the set of \textdef{nonnegative real numbers}, $\N$ for the set of \textdef{nonnegative integers}, and $\naturalpositive$ for the set of \textdef{positive integers}. If $A\subseteq \Rd$, we write $\interior{A}$, $\closure{A}$, and $\bd A$ for its \textdef{interior}, \textdef{closure}, and \textdef{boundary}, respectively. We write $\exterior[0]{A} = \Rd\setminus \closure{A}$ for the \textdef{exterior} of $A$, and we define the \textdef{shell} of $A$ by $\shellof[0]{A} \definedas \Rd\setminus \interior{A}$. We also occasionally use the term ``spherical shell" informally to mean the region lying between two concentric spheres. We denote the \textdef{Euclidean inner product} on $\Rd$ by $\innerprod{\cdot}{\cdot}$; that is, $\innerprod{\vect x}{\vect y} = \sum_{i=1}^d x_i y_i$ for $\vect x,\vect y\in\Rd$. We denote \textdef{Lebesgue measure} on $\Rd$ by $\lebmeas{d}$.

For two subsets $A,B\subseteq \Rd$, the \textbf{Minkowski sum} of $A$ and $B$ is $A+B \definedas \setbuilder{\vect x+\vect y}{\vect x\in A, \vect y\in B}$. Similarly, if $I\subseteq \R$ and $A\subseteq \Rd$, we define the \textbf{product} $I\cdot A = IA \definedas \setbuilder{r\vect x}{r\in I, \vect x\in A}$. If $A$, $B$, or $I$ is a singleton, we write its unique element instead of the set when using set operations such as these.
%
If $\vect x,\vect y\in\Rd$, we let $\segment{\vect x}{\vect y} \definedas \vect x + \segment{0}{1}\cdot (\vect y - \vect x)$ denote the closed \textdef{line segment} from $\vect x$ to $\vect y$, and we similarly define the ``open" line segment $\osegment{\vect x}{\vect y}$ and ``half-open" segments $\ocsegment{\vect x}{\vect y}$ and $\cosegment{\vect x}{\vect y}$. 
If additionally $\vect x\ne \vect y$, we let $\ray{\vect x}{\vect y} \definedas \vect x + (0,\infty)\cdot (\vect y-\vect x)$ denote the ``open" or ``blunt" \textdef{ray} originating at $\vect x$ and passing through $\vect y$, and we write $\cray{\vect x}{\vect y}$ for the corresponding \textdef{closed ray}.
In the case $\vect x=\vect 0$, we use the shorthand
$
\direction{\vect y} \definedas \ray{\vect 0}{\vect y}
=\setbuilder[1]{r\vect y}{r>0}
$
and $\cdirvec[\vect y] \definedas \cray{\vect 0}{\vect y} =\setbuilder[1]{r\vect y}{r\ge0}$.

For $r>0$, the map $\vect x\mapsto r\vect x$ is called a \emph{uniform scaling} or \emph{homothety} of $\Rd$ with center $\vect 0$. More generally, for $r\ge 0$ and $\vect a\in\Rd$, we let $\homothe{r}{\vect a}$ denote the \textdef{homothety} with scaling factor $r$ and center $\vect a$, defined by $\homothe{r}{\vect a} \vect x \definedas \vect a + r(\vect x-\vect a)$ for $\vect x\in\Rd$. The map $\homothe{r}{\vect a} \colon \Rd\to\Rd$ ``zooms in" on the point $\vect a$ if $r>1$ and ``zooms out" from $\vect a$ if $r<1$. That is, any closed ray $\cray{\vect a}{\vect y}$ is invariant under $\homothe{r}{\vect a}$, and each point $\vect x\in \cray{\vect a}{\vect y}$ mapped to the point $\vect x'\in \cray{\vect a}{\vect y}$ whose distance from $\vect a$ is $r$ times the distance from $\vect x$ to $\vect a$.
By analogy with the Minkowski product notation, if $I\subseteq \Rplus$ and $A\subseteq\Rd$, we define $\homothe{I}{\anapex} \definedas \setbuilder{\homothe{r}{\anapex}}{r\in I}$ and $\homothe{I}{\anapex} A \definedas \setbuilder{\homothe{r}{\anapex} \vect x}{ r\in I,\vect x\in A}$. Note that $\homothe{(\Rplus)}{\vect a}$ is a commutative semigroup with identity $1_{\vect a}$, and the subset $\homothe{(0,\infty)}{\vect a}$ is an abelian group.

Let $\norm{\cdot}$ stand for an arbitrary norm on $\Rd$.
We write  $\unitball{\norm{\cdot}}$ and $\sphere{d-1}{\norm{\cdot}}$ for the closed \textbf{unit ball} and \textbf{unit sphere} with respect to $\norm{\cdot}$, i.e.\ $\unitball{\norm{\cdot}} = \setbuilder{\vect x\in \Rd}{\norm{\vect x}\le 1}$, and $\sphere{d-1}{\norm{\cdot}} = \bd \unitball{\norm{\cdot}} =  \setbuilder{\vect x\in \Rd}{\norm{\vect x}= 1}$. Also, for $\vect x\in\Rd$ and $r\in\R$, we write $\reacheddset{\norm{\cdot}}{\vect x}{r} \definedas \vect x + r\unitball{\norm{\cdot}}$ for the closed ball with $\norm{\cdot}$-radius $r$ and center $\vect x$, and we refer to any such ball (or its interior) as a \textdef{$\norm{\cdot}$-ball}. If $\dirsymb\in\Rd\setminus \setof{\vect 0}$, we write $\unitvecnorm{\dirsymb}{\norm{\cdot}} \definedas \dirsymb /\norm{\dirsymb}$ for the \textdef{$\norm{\cdot}$-unit vector} in the direction of $\dirsymb$. We write $\dist{\norm{\cdot}}$ for the translation-invariant \textbf{norm metric} associated with $\norm{\cdot}$; thus, if $A,B\subseteq \Rd$, then
\[
	\dist{\norm{\cdot}}(A,B) \definedas \inf\setbuilder{\norm{\vect y-\vect x}}{\vect x\in A, \vect y\in B}.
\]
If $\vect x\in \Rd$ and $A\subseteq\Rd$, we say that $\vect y$ is a \textdef{$\norm{\cdot}$-closest point in $A$ to $\vect x$} if $\vect y\in A$ and $\distnew{\norm{\cdot}}{\vect y}{\vect x} = \distnew{\norm{\cdot}}{A}{\vect x}$. Such a $\vect y$ always exists if $A$ is nonempty and  closed, and any such $\vect y$ lies in $\bd A$ (cf.\ Lemma~\ref{closest_bd_point_lem}), but $\vect y$ need not be unique in general.
If $\mu_1$ and $\mu_2$ are norms on $\Rd$, the \textdef{dilatation}\footnote{%
\newcommand{\distX}{\distOp_X}%
\newcommand{\distY}{\distOp_Y}%
If $(X, \distX)$ and $(Y, \distY)$ are metric spaces, the \emph{dilatation} (cf.\ \cite{Gromov:1999aa} or \cite{Burago:2001aa}) of a function $f\colon X\to Y$ is defined to be the minimal Lipschitz constant of $f$ (or $\infty$ if $f$ is not Lipschitz), i.e.
\[
\dil f = \sup_{\substack{x,x'\in X \\ x\ne x'}} \frac{\distY\parens[2]{f(x),f(x')}}{\distX(x,x')}.
\]
Our definition of $\dilnorm{\mu_1}{\mu_2}$ above coincides with the dilatation of the norm $\mu_2$ viewed as a function between the metric spaces $(\Rd, \normmetric{\mu_1})$ and $(\Rplus, \normmetric{\abs{\,\cdot\,}})$.
}
of $\mu_2$ with respect to $\mu_1$ is
\[
\newcommand{\thispt}{\dirsymb}
\dilnorm{\mu_1}{\mu_2} \definedas \esupnorm[1]{\frac{\mu_2}{\mu_1}}
= \: \sup_{\dirsymb \ne \vect 0}\: \frac{\mu_2(\dirsymb)}{\mu_1(\dirsymb)},
\]
where $\esupnorm{\cdot}$ denotes the \emph{essential sup norm} for measurable functions on $\Rd$. Geometrically, $\dilnorm{\mu_1}{\mu_2}$ is the $\mu_2$-radius of the smallest $\mu_2$-ball that contains a unit $\mu_1$-ball, i.e.\ $\dilnorm{\mu_1}{\mu_2} = \inf \setbuilder[2]{r> 0}{\unitball{\mu_1} \subseteq r\unitball{\mu_2}}$. The fact that any two norms on $\Rd$ are equivalent means that $0<\dilnorm{\mu_1}{\mu_2}<\infty$ for any norms $\mu_1$ and $\mu_2$. We denote the $\lspace{p}{d}$ norm on $\Rd$ by $\norm{\cdot}_{\lspace{p}{d}}$ for $1\le p\le \infty$ (that is, $\norm{\vect x}_{\lspace{p}{d}} = (\sum_{j=1}^d \abs{x_j}^p)^{1/p}$ for $1\le p<\infty$, and $\norm{\vect x}_{\linfty{d}} = \max \{\abs{x_1},\ldots, \abs{x_d}\}$). For the $\lspace{p}{d}$ norms, we abbreviate the above notations to $\unitball{\lspace{p}{d}}$, $\sphere{d-1}{\lspace{p}{d}}$, $\dist{\lspace{p}{d}}$, and $\dil_{\lspace{p}{d}}$. Note that $\unitball{\ltwo{d}}$ and $\sphere{d-1}{\ltwo{d}}$ are the closed \textdef{Euclidean unit ball} and \textdef{Euclidean unit sphere} in $\Rd$, and that $\frac{1}{2}\unitball{\linfty{d}} = \segment[1]{-\frac{1}{2}}{\frac{1}{2}}^d$ is the closed \textdef{unit cube} centered at $\vect 0$.


For any $A\subseteq \Rd$, we define the \textdef{cube expansion} of $A$ to be $\cubify{A} \definedas A + \frac{1}{2}\unitball{\linfty{d}}$, and we define the \textdef{lattice approximation} of $A$ to be $\lat{A} \definedas \Zd\cap \cubify{A}$. To make formulas more readable, we will sometimes denote the lattice approximation using more than three dots, as in $\lat[1]{AA}$ or $\lat[2]{AAAAAAA}$, or, failing that, with the symbol ``$\latOp$", as in $\lat[5]{A^{A^{A^{A^A}}}}$. Cube expansion of sets $V\subseteq \Zd$ and lattice approximation of sets $A\subseteq \Rd$ exhibit a sort of ``duality" between the lattice and Euclidean space; for example, note that $\cubify{\Zd} = \Rd$ and $\lat[1]{\Rd} = \Zd$. Lemma~\ref{lattice_cube_properties_lem} enumerates several properties of these operations that we will use throughout the rest of the paper.

\subsection{The Graph $\Zd$}
\label{graph_theory_sec}


We consider the integer lattice $\Zd$ both as a graph and as a subset of $\Rd$. Two vertices $\vect u,\vect v\in \Zd$ are adjacent if $\dist{\ell^1}(\vect u,\vect v) = 1$, in which case we write $\edge = \setof{\vect u,\vect v}$ for the edge joining $\vect u$ to $\vect v$, and we say that $\vect u$ and $\vect v$ are the \textdef{endvertices} of $\edge$. We denote the edge set of $\Zd$ by $\edges{\Zd}$. Observe that $\dist{\ell^1}$ coincides with the graph metric on $\Zd$. If $V$ is any subset of $\Zd$, its \textbf{induced subgraph} is the graph with vertex set $V$ and edge set $\edges{V}$ consisting of all edges $\setof{\vect u,\vect v}\in \edges{\Zd}$ with $\vect u,\vect v\in V$. We will frequently identify a set of vertices in $\Zd$ with its induced subgraph, using the two notions interchangeably. The \textdef{boundary} of $V$ is the set $\bd V$ of vertices in $V$ that are adjacent to some vertex outside $V$, and the \textdef{neighbor set} of $V$ is the set $\nbrset{V}$ of vertices in $\Zd\setminus V$ that are adjacent to $V$; thus, $\nbrset{V} = \bd \parens{\Zd\setminus V}$. The \textdef{graph neighborhood} of $V$ is $\nbrhood{V} \definedas V\cup \nbrset{V}$, which we can view as a set of vertices or an induced subgraph.

If $\edge = \setof{\vect u,\vect v}$ with $\vect u\in U$ and $\vect v\in V$, we say that \textdef{$\edge$ joins $U$ to $V$} or that $\edge$ is a \textdef{$U$-$V$ edge}, and we write $\jedges{U}{V}$ for the set of all such edges. A \textdef{boundary edge} of $V$ is any edge $\setof{\vect u,\vect v}\in \edges{\Zd}$ with $\vect u\in V$ and $\vect v\not\in V$ (equivalently, $\vect u\in \bd V$ and $\vect v\in \nbrset{V}$), and we denote the set of boundary edges by $\bdedges{V}$. Observe that $\bdedges{V}= \bdedges{\Zd\setminus V}$. We further define the \textdef{complementary edge set} and \textdef{star edge set} of $V$, respectively, by
\[
\text{$\compedges{V} \definedas \edges{\Zd} \setminus \edges{V}$
and $\staredges{V} \definedas \edges{V}\cup \bdedges{V}$.}
\]
Thus, an edge is in $\compedges{V}$ if and only if \emph{at most} one of its endvertices is in $V$, and an edge is in $\staredges{V}$ if an only if \emph{at least} one of its endvertices is in $V$. Note that $\compedges{V} = \edges{\Zd\setminus V} \cup \bdedges{V} = \staredges{\Zd\setminus V}$.
The star edge set $\staredges{V}$ can be identified with the edge set of the ``star graph" $\stargraph{V}$ of $V$, defined in Section~\ref{fpc_basics:star_path_sec}.
Lemmas~\ref{edge_sets_lem} and \ref{bd_nbr_lem} enumerate several relationships between various sets of edges, boundary vertices, and neighbor vertices.

A \textdef{path} in $\Zd$ is a finite or infinite sequence of vertices in which consecutive vertices are adjacent, together with the sequence of edges joining the consecutive vertices. We view a path both as a pair of sequences, which give the path a well-defined direction or orientation, as well as the undirected connected subgraph of $\Zd$ whose vertex and edge sets are the images of these sequences. (Our definition of ``path" coincides with what is commonly called a ``walk" on the lattice.) We call a path \textdef{edge-distinct} if its edge sequence has no repeated elements, and \textdef{simple} if its vertex sequence has no repeated elements. A simple path must be edge-distinct, but not conversely. 

If $\gamma = \pair[2]{(\vect u_0,\vect u_1,\dotsc,\vect u_{n-1},\vect u_n)}{(\edge_1,\dotsc,\edge_n)}$ is a finite path, we call $\vect u_0$ and $\vect u_n$ the \textdef{initial vertex} and \textdef{final vertex} of $\gamma$, respectively.
If $U,V\subseteq \Zd$, we say that $\gamma$ is a \textdef{path from $U$ to $V$} if $\gamma$ is a finite path whose initial vertex lies in $U$ and final vertex lies in $V$. If we do not care about the orientation of the path, we may call any such $\gamma$ a path \textdef{between $U$ and $V$}. 
A \textdef{subpath} of $\gamma$ is a path whose vertex and edge sets are subsequences of those of $\gamma$. If $\gamma$ is a path in which the vertex $\vect u$ precedes the vertex $\vect v$, we write $\subpath{\gamma}{\vect u}{\vect v}$ for the subpath of $\gamma$ whose vertex sequence comprises all the vertices of $\gamma$ between $\vect u$ and $\vect v$, inclusive.

If $\pathclass$ is a collection of paths in $\Zd$, we call a path in $\pathclass$ \textdef{minimal} if it is a minimal element of the partially ordered set $\pair{\pathclass}{\subseteq}$, where $\subseteq$ is the inclusion partial order on subgraphs, inherited from $2^\Zd$. That is, $\gamma$ is a minimal path in $\pathclass$ if there is no proper subpath of $\gamma$ that is also an element of $\pathclass$. Lemma~\ref{minimal_path_lem} gives an elementary result about minimal paths that will be used in Chapter~\ref{fpc_basics_chap}.

\chapter{First-Passage Competition in a Random Environment}
\label{fpc_basics_chap}

First-passage percolation and the related first-passage competition models are perhaps best described as ``deterministic motion in a random environment."\footnote{%
I borrow this phrase from a talk given by Tom LaGatta at the 2009 Cornell Probability Summer School about his work on Riemannian first-passage percolation in $\Rd$.
}
For first-passage processes on $\Zd$ or other graphs, the random environment is provided by a collection of random edge traversal times which induce a random metric (or pseudometric if some edges can be crossed instantaneously) on the underlying graph. Once the random environment, i.e.\ randomly edge-weighted graph, is given,
the ``motion," i.e.\ the evolution of the process, is completely determined.
I feel that this distinction between the ``deterministic motion" and the ``random environment" is helpful for understanding the proofs of a number of results in first-passage percolation and first-passage competition, because many of the arguments are geometric rather than probabilistic in nature. This is certainly true of the present work, particularly for several proofs in Chapters~\ref{fpc_basics_chap}, \ref{cone_growth_chap}, and \ref{random_fpc_chap}. Thus, throughout the paper I make an effort to separate the probabilistic arguments, i.e.\ those pertaining to the properties of the random traversal times, from the deterministic arguments, which pertain to the evolution of the process once a suitable collection of traversal times is obtained via a realization of some probability measure. I now describe in a bit more detail one way to formally conceptualize the ``random environment" for the first-passage processes.

In the one-type first-passage percolation process, one can think of the environment as a random (pseudo)metric space, and the ``motion" or ``growth" consists of a deterministically expanding metric ball in this space, as described in Section~\ref{intro:one_type_description_sec}. However, the induced pseudometric on the graph carries strictly less information than the underlying collection of edge traversal times, and it is therefore better to think of the environment as a random \textdef{length structure} on the graph. A length structure consists of a class of admissible paths in some space, together with a function for measuring their lengths.\footnote{%
For a general length structure, both the class of admissible paths and the length function are required to satisfy certain natural properties. See Burago et al.\ \cite{Burago:2001aa} or Gromov \cite{Gromov:1999aa} for an axiomatic treatment of length structures in the setting of metric geometry.
} 
For our purposes, the admissible paths comprise all graph paths in $\Zd$ (as defined in Section~\ref{graph_theory_sec}), and the length of a path is the sum of the traversal times of its edges. Once we know how to measure the length of any path, we have all the information needed to construct the first-passage percolation process. Moreover, given the length structure, we can define \textdef{restricted first-passage percolation} in any subgraph $S\subseteq\Zd$ by considering paths that are confined to $S$ (see Section~\ref{fpc_basics:restricted_sec}). On the other hand, knowledge of the induced pseudometric in $\Zd$ on its own does not provide enough information to construct the restricted first-passage percolation process in an arbitrary subgraph $S$, which will be essential to our definition and analysis of the two-type first-passage competition process.
In the two-type (or $n$-type) first-passage competition model described in Section~\ref{intro:two_type_description_sec}, the random environment consists of a graph endowed with a pair (or $n$-tuple) of random length structures, and each species grows according to its corresponding first-passage percolation process until it is forced to take detours around the other species, perhaps being blocked entirely. Given a pair of length structures on $\Zd$, we will see (in Section~\ref{fpc_basics:two_type_construction_sec}) how to construct the two-type process in terms of two families of restricted first-passage percolation processes. Our construction of the two-type process generalizes to pairs of length structures in other spaces. In particular, in Chapter~\ref{deterministic_chap} we will define first-passage competition using the length structures induced by a pair of norm metrics on $\Rd$, providing a deterministic analogue and large-scale limit of the random process studied in the present chapter.

Chapter~\ref{fpc_basics_chap} is organized as follows. In Section~\ref{fpc_basics:setup_sec} we describe the probability space for the one-type and two-type processes, and we introduce the concept of \emph{traversal measure} as a convenient way to encode the random environment. In Section~\ref{fpc_basics:restricted_sec} we define restricted first-passage percolation and introduce a number of important basic concepts and notation that will be used throughout the paper. In Section~\ref{fpc_basics:almost_sure_environment_sec}, we show that on the probability space constructed in Section~\ref{fpc_basics:prob_space_sec}, the random environments for the one-type and two-type processes almost surely satisfy several deterministic properties that will be needed in our construction of the two-type process and in our subsequent analysis of the two processes. In Section~\ref{fpc_basics:two_type_construction_sec}, we construct two-type first-passage competition on $\Zd$ (Definition~\ref{two_type_proc_def}) using the restricted first-passage percolation process defined in Section~\ref{fpc_basics:restricted_sec}, and we prove several basic properties about the two-type process. In particular, Section~\ref{fpc_basics:two_type_construction_sec} contains the first non-trivial result of the paper, namely Proposition~\ref{entangled_well_defined_prop}, in which we show that the process constructed in Definition~\ref{two_type_proc_def} is almost surely ``well-defined" and behaves in accordance with the informal description given in Section~\ref{intro:two_type_description_sec}. In Section~\ref{fpc_basics:basic_properties_sec}, we prove several basic but technical properties of the one-type and two-type processes, which will be needed in Section~\ref{fpc_basics:unbounded_survival_sec} and in later chapters. The results in Section~\ref{fpc_basics:basic_properties_sec} are essentially deterministic, applying to any realization of the random environment that satisfies certain deterministic properties which almost surely hold in our probability space. In Section~\ref{fpc_basics:unbounded_survival_sec}, we prove one of our main results about the two-type process, Theorem~\ref{irrelevance_thm}, which shows that if one species starts on a fixed infinite set, the particular location of the other species' starting set is essentially irrelevant to the question of whether it has a positive probability of surviving.

\section{Setup and Properties of the Traversal Times}
\label{fpc_basics:setup_sec}
In this section we we describe the general setup for the first-passage processes we will be studying. In Section~\ref{fpc_basics:prob_space_sec} we specify conditions on the collections of random traversal times $\setof{\ttime(\edge)}_{\edge\in\edges{\Zd}}$ and $\setof[2]{\parens[2]{\ttime_1(\edge),\ttime_2(\edge)}}_{\edge\in\edges{\Zd}}$ that are sufficient to obtain random length structures suitable for constructing a one-type or two-type process, and  in Section~\ref{fpc_basics:traversal_meas_sec} we introduce the concept of ``traversal measure" as a convenient way of interpreting the random environment provided by the collection of edge traversal times on $\Zd$.

\subsection{Assumptions on the Probability Space}
\label{fpc_basics:prob_space_sec}


\subsubsection{Model Traversal Times}

When describing the one-type or two-type process, it will be convenient to speak about generic random variables --- say $\mttime_0$ for the one-type process, or $\mttimes{1}_0$ and $\mttimes{2}_0$ for the two-type process --- whose distributions match the common distribution of the edge-traversal times $\ttime(\edge)$ or $\ttime_i(\edge)$ for the corresponding species, but which are not tied to any particular edge $\edge$. In fact, it will be convenient to have available whole \iid\ sequences of such variables which are ``disinterested" in particular edges of the graph. We will call these generic random variables the \textdef{model traversal times} for the corresponding species, because we can use them to model the behavior of the traversal times. That is, we can think of drawing repeated samples from the distribution of $\mttime_0$ or $\mttimes{i}_0$ to obtain a prediction of how long the corresponding species will have to wait to cross each new edge it encounters, without having to specify any particular edges. To describe the traversal times, we take the approach of first describing the properties of the generic \emph{model} traversal times, and then requiring the traversal times to behave the same way. The prototypical example to keep in mind for the model traversal time is an exponential random variable, i.e.\ $\Pr \eventthat{\mttime_0 > t} = e^{-\lambda t}$ for some $\lambda >0$.

\subsubsection{Properties of the Model Traversal Times}
To set up our probability space, we first make several definitions for properties we may require of a given random variable. Suppose $Y$ is a nonnegative random variable on some probability space $\probabilityspace$ with expectation operator $\E$; we will take $Y$ to be the model traversal time $\mttime_0$ or $\mttimes{i}_0$ when describing the probability space for the one-type or two type process.

\subsubsection{Non-Percolation Assumption}

The first property places a restriction on the mass of $Y$ at 0:
\begin{description}
\item[\finitespeed{d}] $\Pr \eventthat[2]{Y = 0} < \pcrit$, where $d$ is a positive integer, and $\pcrit$ is the critical probability for Bernoulli bond percolation on $\Zd$.
\end{description}
The condition \finitespeed{d} stands for ``Finite Percolation Speed in $d$ dimensions." When $Y$ is a model traversal time, \finitespeed{d} guarantees that the infection cannot instantaneously infect an infinite cluster. In \cite{\KesAspects}, it is shown that this is sufficient for the first-passage percolation process to have a finite asymptotic speed in every direction (cf.\ Theorem~\ref{shape_thm}, the Shape Theorem). Observe that \finitespeed{d} holds for any $d$ if $Y$ does not have an atom at 0.

\subsubsection{Moment Conditions}

The next three properties are moment conditions on $Y$: 
\begin{description}
\item[\expmoment] There exists $b>0$ such that $\E e^{b Y} < \infty$.
\item[\Lspace{p}] $\E Y^p <\infty$, where $p>0$.
\item[\minimalmoment{d}] $\E \min \setof{Y_1^d,\dotsc,Y_{2d}^d} < \infty$, where $d$ is a positive integer, and $Y_1,\dotsc,Y_{2d}$ are \iid\ random variables with $Y_j \eqd Y$ for $1\le j\le 2d$.
\end{description}
Condition \expmoment\ stands for ``Exponential Moment," and condition \Lspace{p} stands for ``$Y\in L^p$."
The final condition, \minimalmoment{d}, stands for ``Minimal Moment condition in $d$ dimensions."
Clearly \expmoment\ implies that \Lspace{p} holds for all $p>0$. In Lemma~\ref{sufficient_moment_cond_lem}, we show that \Lspace{1/2} implies that \minimalmoment{d} holds for all $d\ge 1$.

When $Y$ is a model traversal time, the condition \minimalmoment{d} means that the ``escape time" from a vertex, i.e.\ the time it takes to move from a given vertex to a neighboring vertex,
has finite $d\th$ moment. In \cite{\KesAspects}, it is shown that  \minimalmoment{d} is the weakest moment condition necessary to get a Shape Theorem (Theorem~\ref{shape_thm}) for an \iid\ first-passage percolation process, in that if the condition fails, then the infected region will have persistent holes far from its outer boundary.

In the present chapter, the Shape Theorem is the only relevant technical result, so we will assume the weakest possible moment condition \minimalmoment{d} for the model traversal times. Starting in Chapter~\ref{cone_growth_chap}, we will assume that the model traversal times satisfy the stronger condition \expmoment\ in order to obtain large deviations estimates on the growth of the process.

\subsubsection{``No Ties" Assumption on the Traversal Times}

The last property we consider deals with a pair of \iid\ sequences of nonnegative random variables, say $\pair{\boldsymbol{X}}{\boldsymbol{Y}}$, where $\boldsymbol{X} = \setof{X_j}_{j\in\N}$ and $\boldsymbol{Y} = \setof{Y_j}$:
\begin{description}
\item[\notiesprob] For any $m,n\in \N$, $\Pr \eventthat[1]{\sum_{j=0}^m X_j = \sum_{j=0}^n Y_j} = 0$.
\end{description}
The condition \notiesprob\ stands for ``No Ties with Probability 1." If $\boldsymbol{X}$ and $\boldsymbol{Y}$ represent sequences of traversal times for the two species in the competition process, then \notiesprob\ guarantees that the two species cannot reach the same vertex at the same time. This property will be necessary to ensure that the two-type process is well-defined. Note that if the sequences $\boldsymbol{X}$ and $\boldsymbol{Y}$ are mutually independent, then \notiesprob\ is satisfied if the random variables $X_0$ and $Y_0$ are both continuous (i.e.\ nonatomic).

\subsubsection{The Probability Space for the One-Type Process}

Fix an integer $d\ge 2$, and let  $\ttime = \setof{\ttime(\edge)}_{\edge\in\edges{\Zd}}$ be an \iid\ collection of nonnegative random variables (the \textdef{traversal times}) in some probability space $\probabilityspace$ with expectation operator $\E$, and assume that for any $\edge\in\edges{\Zd}$, the random variable $\ttime(\edge)$ satisfies \finitespeed{d}\ and \minimalmoment{d}. Also let $\mttime = \setof[2]{\mttime_j}_{j\in\N}$ be an \iid\ sequence of random variables (the \textdef{model traversal times}) with $\mttime_j \eqd \ttime(\edge)$ for each $j\in\N$ and $\edge\in\edges{\Zd}$.

As described above, the non-percolation condition \finitespeed{d}\ guarantees that the growth of the one-type process is not too fast, whereas the moment condition \minimalmoment{d} guarantees that the growth is not too slow, so these can be seen as complementary restrictions on the growth of the process. As noted above, we will replace \minimalmoment{d} with the stronger condition \expmoment\ in later chapters.

%

\subsubsection{The Probability Space for the Two-Type Process}

%

Let $\mttimes{1}= \setof[2]{\mttimes{1}_j}_{j\in\N}$ and $\mttimes{2}= \setof[2]{\mttimes{2}_j}_{j\in\N}$ be independent sequences of \iid\ nonnegative random variables (the \textdef{model traversal times} for the two species) on some probability space $\probabilityspace$ with expectation operator $\E$, and let $\marginals{i}$ be the distribution of $\mttimes{i}_0$. Suppose the random variables $\mttimes{1}_0$ and $\mttimes{2}_0$ both satisfy \finitespeed{d}\ and \minimalmoment{d}, and suppose that the pair of sequences $\pair[2]{\mttimes{1}}{\mttimes{2}}$ satisfies \notiesprob.

Now let $\jointlaw$ be any probability measure on $\R_+\times \R_+$ with marginals $\marginals{1}$ and $\marginals{2}$ (i.e.\ $\jointlaw$ specifies some coupling of $\mttimes{1}_0$ and $\mttimes{2}_0$), and let $\cmttimes{}{} = \setof[2]{\cmttimes{}{j}}_{j\in\N} =  \setof[2]{\pair[2]{\cmttimes{1}{j}}{\cmttimes{2}{j}}}_{j\in\N}$ be an \iid\ sequence of random vectors in $\Rplus\times\Rplus$ with common distribution $\jointlaw$. We call $\cmttimes{}{}$ the sequence of \textdef{coupled model traversal times}. Finally, for a fixed integer $d\ge 2$, let $\ttimepair = \setof[2]{\ttimepair(\edge)}_{\edge\in\edges{\Zd}} = \setof[2]{\parens[2]{\ttime_1(\edge),\ttime_2(\edge)}}_{\edge\in\edges{\Zd}}$ be an \iid\ collection of random vectors with $\ttimepair(\edge) \eqd \cmttimes{}{j}$ for $\edge\in\edges{\Zd}$ and $j\in\N$. For $i\in\setof{1,2}$, the random variables $\setof[2]{\ttime_i(\edge)}_{\edge\in\edges{\Zd}}$ are the \textdef{traversal times} for species $i$.

Note that both \finitespeed{d}\ and \notiesprob\ hold if the random variables $\mttimes{1}_0$ and $\mttimes{2}_0$ are continuous (i.e.\ non-atomic), for example. Also note that if \notiesprob\ holds, then at most one of the random variables $\mttimes{1}_0$ and $\mttimes{2}_0$ can have an atom at zero, which renders one of the \finitespeed{d} assumptions redundant. Again, we will replace the moment condition \minimalmoment{d} with the stronger condition \expmoment\ in later chapters.

\subsubsection{The Canonical Sample Space}

In the preceding construction, we have followed the pedagogy of deliberately leaving the sample space ambiguous, focusing instead on objects and properties which are independent of the particular probability space. However, we can readily construct a canonical sample space for the traversal times of both a one-type and a two-type process as follows.

For the two-type process, let $\marginals{1}$ and $\marginals{2}$ be the distributions of two nonnegative random variables satisfying \finitespeed{d}, \minimalmoment{d}, and \notiesprob. That is, $\marginals{1}$ and $\marginals{2}$ are probability measures on $\R_+$ which satisfy
\begin{description}
\item[\finitespeed{d}] $\marginals{i}\setof{0} < \pcrit$.
\item[\minimalmoment{d}] $d \int_0^\infty \brackets[2]{\marginals{i}(x,\infty)}^{2d} x^{d-1} dx <\infty$. 
\item[\notiesprob] For any $m_1, m_2\in \naturalpositive$, $\parens{\marginals{1}}^{* m_1} \otimes \parens{\marginals{2}}^{* m_2} \parens[3]{\setbuilder[2]{(x,x)}{x\in \R_+}} = 0$.
\end{description}
Starting in Chapter~\ref{cone_growth_chap}, we will replace the moment condition \minimalmoment{d} with the stronger condition
\begin{description}
\item[\expmoment] There exists $b>0$ such that $\int_0^\infty e^{bx}\, d\marginals{i}(x) < \infty$.
\end{description}
Now let $\jointlaw$ be a coupling of $\marginals{1}$ and $\marginals{2}$. Fix an integer $d\ge 2$ and let $\samplespace = \parens{\R_+ \times \R_+}^{\edges{\Zd}}$, endowed with the product Borel $\sigma$-field $\sigmafield$, and let $\Pr = \jointlaw^{\otimes \edges{\Zd}}$. For a realization $\outcome\in \samplespace$ of $\Pr$, define the two-type traversal times by $\ttimepair(\edge)_\outcome \equiv \parens[2]{\ttime_1(\edge)_\outcome, \ttime_2(\edge)_\outcome} = \outcome(\edge)$ for $\edge\in\edges{\Zd}$. To recover the model traversal times, label the edges along the positive axis of the first (second) coordinate as $\labelededge{1}{0}, \labelededge{1}{1}, \labelededge{1}{2},\dotsc$ (resp.\ $\labelededge{2}{0}, \labelededge{2}{1}, \labelededge{2}{2},\dotsc$), then set $\mttimes{i}_j(\outcome) = \ttime_i \parens{\labelededge{i}{j}}_\outcome$ for $i\in\setof{1,2}$ and $j\in\N$.

For the one-type process, take $(\samplespace,\sigmafield,\Pr)$, $\ttimepair$, and $\mttimes{i}$ as defined above for the two-type process, then set $\ttime(\edge) = \ttime_1(\edge)$ for $\edge\in\edges{\Zd}$, and set $\mttime_j = \mttimes{1}_j$ for $j\in\N$.

%

\subsection{Traversal Measure = The Random Environment}
\label{fpc_basics:traversal_meas_sec}

For a fixed realization $\outcome$ of the one-type process, the collection $\setof{\ttime(\edge)}_{\edge\in\edges{\Zd}}$ defines a function $\ttime \colon \edges{\Zd}\to \R_+$. We can equivalently think of $\ttime$ as a \emph{measure} on $\edges{\Zd}$, with domain $2^{\edges{\Zd}}$, by defining\footnote{%
The construction of a purely atomic measure as in \eqref{traversal_measure_def} works for an arbitrary $[0,\infty]$-valued function in place of $\ttime$; see \cite[p.~25]{Folland:1999aa} for further examples. Note that our convention of using the symbol $\ttime$ to denote both a function on $\edges{\Zd}$ and the induced measure in \eqref{traversal_measure_def} conflicts with the usual notation for the image of $F\subseteq \edges{\Zd}$ under the function $\ttime$, i.e.\ $\ttime(F) = \setbuilder{\ttime(\edge)}{\edge\in F}$; however, we will never use the notation $\ttime(F)$ in this sense, so there should be no confusion.}%
\begin{equation}\label{traversal_measure_def}
\ttime(F) \definedas \sum_{\edge\in F} \ttime(\edge)
\quad\text{for } F\subseteq \edges{\Zd}.
\end{equation}
We will call the random measure $\ttime$ the \textdef{traversal measure for the one-type process}. Similarly, we can think of the collection of random vectors $\ttimepair = \setof[2]{\parens[2]{\ttime_1(\edge),\ttime_2(\edge)}}_{\edge\in\edges{\Zd}}$ as a pair of measures or a vector-valued measure on $\edges{\Zd}$, and we will refer to $\ttimepair$ as the \textdef{traversal measure for the two-type process}. More generally, we will use the term ``traversal measure" to refer to any $\sigma$-finite $[0,\infty]$-valued or $[0,\infty]^2$-valued measure on $\edges{\Zd}$, and we will use the terms ``traversal measure" and ``collection of traversal times" interchangeably when referring to the objects $\ttime$ and $\ttimepair$. 

Extending the definition \eqref{traversal_measure_def}, we further overload the use of the symbol $\tau$ by declaring that for any ``suitable" object $X$, the notation $\tmeasure(X)$ means the sum of the traversal times of ``edges in $X$," where the the final quoted phrase must be interpreted based on the nature of $X$. Our two primary uses of this convention will be: (1) If $V\subseteq \Zd$, then $\tmeasure(V)$ means the traversal measure of the edge set of $V$, viewed as an induced subgraph of $\Zd$; and (2) if $A\subseteq \Rd$, then $\tmeasure(A) \definedas \tmeasure(\lat A)$, where the right-hand appearance of $\tmeasure$ refers to the usage just defined in (1).


The introduction of traversal \emph{measure} has two main advantages over viewing $\ttime$ (or $\ttimepair$) merely as a random \emph{function} on edges. First, it simplifies notation, because our primary use of the function $\ttime$ will be to measure the time it takes to traverse a lattice path $\gamma \subseteq \parens[1]{\Zd,\edges{\Zd}}$, by adding up the traversal times $\ttime(\edge)$ for the edges in $\gamma$. By interpreting $\ttime$ as a measure, the traversal time of an edge-distinct path $\gamma$ is written simply as $\ttime(\gamma)$, and thus all the relevant information of the traversal times is conveniently encoded in the single high-level object $\ttime$. Second, and perhaps more importantly, identifying $\ttime$ as a measure makes it easier to draw direct analogies with the deterministic process in the next chapter, as well as with other continuum variants of first-passage percolation: The fundamental object that is needed to define any\footnote{%
One possible exception is\ Deijfen's ``outburst" model \cite{Deijfen:2003aa} --- it is not immediately clear whether it is possible to identify a length measure analogous to $\tmeasure$ in this case.
On the other hand, a continuum model with an obvious analogue of $\tmeasure$ is LaGatta and Wehr's Riemannian first-passage percolation \cite{LaGatta:2010aa}, which is constructed using the length structure of a random Riemannian metric.
}
type of first-passage percolation process is a (random or deterministic) \emph{measure}, or equivalently,\footnote{%
Any metric induces a length structure that coincides with the metric's one-dimensional Hausdorff measure on simple paths, and any lower semi-continuous length structure
can be induced by a metric (see \cite[ pp.~39, 53]{Burago:2001aa}), hence corresponds to some measure.
Conversely, given a suitable measure ($\tmeasure$ in our case), one can construct a corresponding length structure.%
}
\emph{length structure} that assigns lengths to paths in the underlying space. The first-passage percolation process is then defined using the \emph{induced intrinsic metric} arising from the length structure (cf.\ Burago et al.\ \cite{Burago:2001aa} or Gromov \cite{Gromov:1999aa} for more information about intrinsic metrics and length structures). We will return to these ideas and define intrinsic metrics when discussing the deterministic process in Chapter~\ref{deterministic_chap}.

The random measure $\ttime$ contains all the information needed to construct the first-passage percolation process, and moreover it contains all the randomness in the model; that is, the traversal measure essentially \emph{is} the random environment for first-passage percolation, and treating it as an object of interest makes it relatively easy to separate out the deterministic motion from the background randomness. Reversing our perspective from that of the previous section, we say that a random traversal measure $\tmeasure$ on $\edges{\Zd}$ is \iid\ with model traversal time $\mttime_0$ if the collection of random variables $\setof{\ttime(\edge)}_{\edge\in\edges{\Zd}}$ is \iid\ with $\ttime(\edge)\eqd \mttime_0$ for all $\edge$. In this case, we use the term ``model traversal times" to refer to any \iid\ sequence $\setof{\mttime_j}_{j\in\N}$ with $\mttime_j\eqd \ttime(\edge)$.

\section{The Restricted First-Passage Percolation Process}
\label{fpc_basics:restricted_sec}

%
%
%
%
%

In this section we define first-passage percolation restricted to a subset of $\Zd$, obtaining the ordinary unrestricted version as a special case. Restricting the growth of, say, the \speciesname{red} species to $S\subseteq \Zd$ means that \speciesname{red} starts on some set $A\subseteq S$ and is only allowed to travel along paths contained in $S$. Thus, \speciesname{red} is confined to $S$ and cannot infect any vertices in $\Zd\setminus S$. The restricted one-type process will play a central role both in the construction of the two-type process in Section~\ref{fpc_basics:two_type_construction_sec} and in the analysis of first-passage competition in Sections~\ref{fpc_basics:basic_properties_sec} and \ref{fpc_basics:unbounded_survival_sec} and in Chapters~\ref{random_fpc_chap} and \ref{coex_finite_chap}, because it allows us to separate the growth of the two species, treating them as two separate, simpler one-type processes. The idea of restricting the set of paths in first-passage percolation has been around as long as the subject itself
and has played an important role in the analysis of both the one-type and two-type models.
See, for example, \cite{Kesten:1986aa}, \cite{Chatterjee:2009aa}, \cite{Ahlberg:2011aa}, \cite{Haggstrom:2000aa}, \cite{Deijfen:2007aa}.



Although the restricted process itself is not new, we present a simple new idea here that will allow us to construct the two-type process in terms of two restricted one-type processes.
Namely, when \speciesname{red} is confined to the set $S$, every vertex \emph{neighboring} $S$ has a well-defined time at which it \emph{would} be reached by \speciesname{red}\ldots that is, \emph{if} \speciesname{red} were allowed to take one extra step from a boundary vertex of $S$ to a neighbor outside of $S$. This observation is the key to the new construction in Section~\ref{fpc_basics:two_type_construction_sec} because it allows us to detect the mutual interaction of two species growing in disjoint subsets of $\Zd$, by comparing the times the two species take to reach vertices along the common boundaries of the restricting sets. The precise formulation of this statement is Proposition~\ref{conquering_property_prop} in Section~\ref{fpc_basics:two_type_construction_sec} below, which motivates the definition of the two-type process in Definition~\ref{two_type_proc_def}. In order to formalize the notion of restricting the growth of \speciesname{red} to $S$ while allowing ``peeks" at neighboring vertices, we make the following definitions.

%
%

\subsection{Restricted and Star-restricted Paths}
\label{fpc_basics:star_path_sec}


For $S\subseteq \Zd$, we call a lattice path $\gamma$ an \textdef{$S$-path} if all of its edges are contained in $\edges{S}$, and we call $\gamma$ an \textdef{$\stargraph{S}$-path} if $\gamma$ is contained in $S$ except possibly for its two endvertices and their incident edges (thus every $S$-path is also an $\stargraph{S}$-path).
We can also think of an $\stargraph{S}$-path as the projection in $\Zd$ of a path in the graph $\stargraph{S}$ (the \textdef{star graph\footnote{%
We use the term ``star graph" because this definition generalizes the notion of the star centered at a vertex (see, e.g.\ Diestel \cite{Diestel:2010aa}) to a star ``centered at $S$," as $\stargraph{S}$ consists of the subgraph $S$ with a union of stars centered at its boundary vertices.
} of $S$}),
which we define by attaching to $S$ (viewed as an induced subgraph of $\Zd$) all the edges connecting $S$ to its complement, and adding a distinct vertex to the outer end of each of these boundary edges. That is, different boundary edges of $S$ lead to distinct vertices in $\stargraph{S}$, even if they are adjacent to the same vertex in $\Zd\setminus S$.
The graph $\stargraph{S}$ then projects naturally onto the graph $\pair[2]{\nbrhood{S}}{\staredges{S}}\subseteq \Zd$: Each vertex and edge of $S$ is mapped to itself; each of the added boundary edges maps onto the corresponding edge in $\Zd$; and each added endvertex is mapped to the endvertex of the corresponding boundary edge in $\Zd$.
This projection map is injective on $S$ and on the added boundary edges, but not necessarily on the added endvertices of these edges, since a single vertex in $\nbrset{S}$ may be adjacent to several vertices in $S$, each of which gives rise to a distinct new vertex in $\stargraph{S}$. With this definition, the $\stargraph{S}$-paths in $\Zd$ are precisely the images of paths in $\stargraph{S}$ under this map. Additionally, the connected components of $\stargraph{S}$ are in one-to-one correspondence with the connected components of $S$.

\subsection{Restricted and Unrestricted Passage Times}
\label{fpc_basics:passage_time_sec}

Now let $\ttime$ be a traversal measure (i.e.\ any $\sigma$-finite $[0,\infty]$-valued measure) on $\edges{\Zd}$.  For each subset $S\subseteq \Zd$, we define the $\ttime$-induced \textdef{$S$-restricted passage time} between two sets $U,V\subseteq \Zd$ as
\begin{equation}
\label{restricted_ptime_def_eqn}
\rptime{\tmeasure}{S}{U}{V} := \inf \setbuilder[3]{\ttime(\gamma)}
	{\text{$\gamma$ is an $S$-path from $U$ to $V$}},
\end{equation}
and we define the \textdef{$\stargraph{S}$-restricted passage time} between $U$ and $V$ as
\begin{equation}
\label{star_ptime_def_eqn}
\starptimenew{\tmeasure}{S}{U}{V} := \inf \setbuilder[3]{\ttime(\gamma)}
	{\text{$\gamma$ is an $\stargraph{S}$-path from $U$ to $V$}}.
\end{equation}
In both definitions we follow the convention that $\inf \emptyset = \infty$.
If the traversal measure is clear from context or is irrelevant, we may omit it from the notation and write $\rptmetric{S}$ or $\starptmetric{S}$ instead of $\rptmetric[\ttime]{S}$ or $\starptmetric[\ttime]{S}$.
If $S=\Zd$, we omit the superscript, and we call $T(U,V)\definedas T^\Zd(U,V)$ the \textdef{unrestricted passage time}, or simply the \textdef{passage time}, from $U$ to $V$. For operations such as $T^S(\cdot,\cdot)$ that take sets as arguments, we will follow the convention of dropping braces around singletons, e.g.\ $T^S(\vect u,V) \definedas T^S(\setof{\vect u},V)$. 

Taking both $U$ and $V$ to be singletons, $\rptmetric{S}$ and $\starptmetric{S}$ define functions on $\Zd\times\Zd$, which we can think of generically as distance functions on $\Zd$.
For each $S\subseteq\Zd$, the function $\rptmetric{S}\colon \Zd\times \Zd \to [0,\infty]$ is a pseudometric on $\Zd$ (i.e.\ $\rptmetric{S}$ is symmetric, satisfies the triangle inequality, and $\rptmetric{S}(\vect u,\vect u)=0$), and $\rptmetric{S}$ restricts to a finite pseudometric on each connected component of $S$. In fact, $\rptmetric{S}$ is precisely the intrinsic pseudometric on $S$ induced by the restriction to $S$ of the length structure corresponding to $\ttime$ (see \cite[pp.~31, 42]{Burago:2001aa}). Moreover, $\rptmetric[\ttime]{S}$ is a metric if and only if $\ttime(\edge) \ne 0$ for all $\edge\in \edges{S}$. Observe that $\rptmetric{S}(\vect u,\vect v)$ and $\starptime{S}(\vect u,\vect v)$ agree whenever neither $\vect u$ nor $\vect v$ is a neighbor of $S$; in particular $\rptmetric{S}$ and $\starptmetric{S}$ define the same pseudometric on $S$. The distance function $\starptmetric{S}$ is not in general a pseudometric on $\Zd$, as the triangle inequality may fail for points in $\nbrset{S}$. However, $\starptmetric{S}$ is a pseudometric on the star graph $\stargraph{S}$ which is finite on components.

\subsection{Geodesics and Star Geodesics}
\label{fpc_basics:geodesics_sec}

Any (necessarily simple and finite) path $\gamma$ that achieves the infimum in \eqref{restricted_ptime_def_eqn} or \eqref{star_ptime_def_eqn} for some $U,V\subseteq\Zd$ is called a \textdef{shortest path} or \textdef{geodesic}\footnote{%
In a broader metic geometry setting (see, e.g.\ \cite[p.~51]{Burago:2001aa}), the term ``geodesic" refers to any path that is \emph{locally} distance minimizing, not just the \emph{globally} shortest paths we refer to as geodesics. In general, shortest paths are more properly called \emph{minimizing geodesics}. However, the concept of ``locally distance minimizing" is not as universally useful on graphs as in continuum settings, since it is unclear what ``locally" should mean in general. Hence, in graph theoretic contexts, it is standard to use the term ``geodesic" to mean ``shortest path."%
}
for the respective distance function $\rptmetric{S}$ or $\starptmetric{S}$; we refer to such a path as a $\rptmetric{S}$-geodesic or $\starptmetric{S}$-geodesic, respectively, from $U$ to $V$. Equivalently, $\gamma$ is a $\rptmetric{S}$-geodesic if and only if $\rptmetric{S}(\vect u,\vect v) = \rptmetric{\gamma}(\vect u,\vect v)$ for all vertices $\vect u,\vect v\in\gamma$, and similarly for $\starptmetric{S}$ (observe that $\rptime{\tmeasure}{\gamma}{\vect u}{\vect v} = \ttime \parens[2]{\gamma[\vect u,\vect v]}$, where $\gamma[\vect u,\vect v]$ is the subpath of $\gamma$ from $\vect u$ to $\vect v$); this definition of geodesics applies equally well to infinite paths. 
We will consider the question of the existence of finite geodesics in the next section. 

\subsection{Restricted and Unrestricted Growth}
\label{fpc_basics:restricted_growth_sec}

If $A,S\subseteq\Zd$ (typically with $A\subseteq S$), we define the \textdef{$S$-restricted first-passage percolation process} started from $A$ as
\[
\rreachedset{\ttime}{A}{S}{t} \definedas
\setbuilder[2]{\vect v\in S}{\rptime{\ttime}{S}{A}{\vect v}\le t}
\quad\text{for } t\ge 0,
\]
and the \textdef{$\stargraph{S}$-restricted first-passage percolation process} started from $A$ as
\[
\rreachedset{\ttime}{A}{\stargraph{S}}{t} \definedas
\setbuilder[2]{\vect v\in \Zd}{\starptimenew{\ttime}{S}{A}{\vect v}\le t}
\quad\text{for } t\ge 0,
\]
and we set
\[
\rreachedset{\ttime}{A}{S}{\infty} \definedas \bigcup_{t\ge 0} \rreachedset{\ttime}{A}{S}{t}
\quad\text{and}\quad
\rreachedset{\ttime}{A}{\stargraph{S}}{\infty} \definedas \bigcup_{t\ge 0} \rreachedset{\ttime}{A}{\stargraph{S}}{t}.
\]
Note that $\rreachedset{\ttime}{A}{S}{t}\subseteq S$ by definition, and $\rreachedset{\ttime}{A}{\stargraph{S}}{t}\subseteq \nbrhood{S}$. If $S=\Zd$, we call $\rreached[\ttime]{A}{\Zd}(t)$ the \textdef{unrestricted first-passage percolation process}, and as with the pseudometrics $\ptime$, we drop the superscript $S$ from the notation: $\reachedx[\ttime]{A}(t) \definedas \rreached[\ttime]{A}{\Zd}(t)$ for $t\in [0,\infty]$.
Note that it follows from the above definition that for any $A,B,S\subseteq\Zd$,
\[
\rptime{\tmeasure}{S}{A}{B}
= \inf \setbuilder[1]{t}{\rreachedset{\ttime}{A}{S}{t} \cap B \ne \emptyset}.
\]
That is, $\rptime{\tmeasure}{S}{A}{B}$ is the \textdef{hitting time} of the set $B$ for the restricted growth process $\rreachedfcn{\ttime}{A}{S}$.

Finally, we write $\ptimefamily[\ttime]$ for the family of pseudometrics $\setbuilder[2]{\rptmetric[\ttime]{S}}{S\subseteq \Zd}$ and $\starptimefamily[\ttime]$ for the family $\setbuilder[2]{\starptime{S}_\ttime}{S\subseteq \Zd}$, and we write $\fppprocfamily[\ttime]$ for the family of first-passage percolation processes $\setbuilder[3]{\rreached[\ttime]{A}{S}}{A\subseteq S\subseteq\Zd}$. Thus, every traversal measure $\tau$ induces a family $\ptimefamily[\ttime]$ of pseudometrics on $\Zd$, which in turn induces a family $\fppprocfamily[\ttime]$ of first-passage percolation processes.


\subsection{Extensions to $\Rd$ via Lattice Approximation}
\label{fpc_basics:continuum_extension_sec}


It will often be convenient to view the first-passage percolation and competition processes as embedded in $\Rd$ rather than just $\Zd$, particularly in Chapters~\ref{cone_growth_chap} and \ref{random_fpc_chap}. To this end, we make the following definitions. Recall that for $A\subseteq \Rd$, the cube expansion of $A$ is $\cubify{A} = A+ \brackets[1]{-\frac{1}{2}, \frac{1}{2}}^d$, and the lattice approximation of $A$ is $\lat{A} = \Zd\cap \cubify{A}$. Using lattice approximation, we can extend the definitions of $\ptimefamily[\tmeasure]$ and $\fppprocfamily[\ttime]$ to subsets of $\Rd$.

If $A,B,S\subseteq\Rd$, we define the \textdef{$S$-restricted passage time between continuum sets $A$ and $B$} by
\begin{equation}\label{continuum_ptime_def_eqn}
\rptmetric{S}(A,B) \definedas \rptmetric{\lat{S}}\parens[2]{\lat{A},\lat{B}}.
\end{equation}
As before, if $\lat{S} = \Zd$ then we drop the superscript $S$ from the notation, so that $\ptime(A,B)\definedas \rptmetric{\Zd}(A,B)$ represents the unrestricted passage time between $A,B\subseteq \Rd$.

For $A,S\subseteq\Rd$ (typically with $\lat{A}\subseteq \lat{S}$), we use the generalized passage times in \eqref{continuum_ptime_def_eqn} to define the \textdef{$S$-restricted continuum process started from $A$}:
\begin{equation}\label{continuum_process_def_eqn}
\rreachedset{\tmeasure}{A}{S}{t} \definedas \setbuilder[1]{\vect x\in S}{\rptime{\tmeasure}{S}{A}{\vect x}\le t}
\quad \text{for } t\ge 0.
\end{equation}
It follows from the definitions that $\rreachedset{\tmeasure}{A}{S}{t} = S\cap \cubify[1]{\rreachedset{\ttime}{\lat{A}}{\lat{S}}{t}}$. That is, $\rreachedset{\tmeasure}{A}{S}{t}$ is the process obtained by placing a closed unit cube around each point of the lattice process $\rreachedset{\ttime}{\lat{A}}{\lat{S}}{t}$, and then intersecting this ``cubified" process with $S$. Again, if $S=\Zd$ we omit it from the notation, so that for $A\subseteq\Rd$, $\reachedset{\tmeasure}{A}{t} = \reachedset{\tmeasure}{\lat{A}}{t}$ is the lattice process started from $\lat A$.
Note that if we take $S=\Rd$, the process
\begin{equation}\label{unrestricted_continuum_process_eqn}
\rreachedset{\tmeasure}{A}{\Rd}{t}
= \setbuilder[1]{\vect x\in \Rd}{\ptimenew{\tmeasure}{A}{\vect x}\le t}
= \cubify[1]{\reachedset{\ttime}{A}{t}}
\end{equation}
gives a cube-expanded continuum version of the unrestricted lattice process $\reachedset{\tmeasure}{\lat A}{t}$.

Some caution is needed when using the continuum passage times defined in \eqref{continuum_ptime_def_eqn}: Unlike the lattice passage times defined in \eqref{restricted_ptime_def_eqn}, the continuum passage times in \eqref{continuum_ptime_def_eqn} fail to satisfy the triangle inequality for general points $\vect x,\vect y,\vect z\in\Rd$. In Chapter~\ref{cone_growth_chap}, we will get around this problem by defining \emph{covering times} $\ctime{}{S}{A}{B}$ (cf.\ \eqref{continuum_covering_time_def}, p.~\pageref{continuum_covering_time_def}), which lack the symmetry between $A$ and $B$ enjoyed by passage times, but satisfy the triangle inequality for arbitrary $A,B,C\subseteq \Rd$.

\subsection{Elementary Properties of Restricted First-Passage Percolation}
\label{fpc_basics:elementary_fpp_properties_sec}

Here we enumerate some elementary properties of the first-passage percolation process.
For $U,V,S\subseteq \Zd$, define a \textdef{minimal $S$-path from $U$ to $V$} to be an $S$-path from $U$ to $V$ which contains no proper subpath that is also an $S$-path from $U$ to $V$, and similarly for $\stargraph{S}$-paths. By Lemma~\ref{minimal_path_lem}, any $S$-path from $U$ to $V$ contains a minimal $S$-path from $U$ to $V$, and similarly for $\stargraph{S}$-paths.

\begin{lem}[{Elementary properties of $\ptimefamily[\tmeasure]$}]
\label{elementary_fpp_properties_lem}
Let $\tmeasure$ be a one-type traversal measure on $\edges{\Zd}$ with corresponding family of pseudometrics $\ptimefamily[\tmeasure]$, and let $U,V,S\subseteq\Zd$ and $A,B\subseteq\Rd$.
\begin{enumerate}
\item \label{elementary_fpp_properties:monotonicity_part}
The passage times $\rptime{\tmeasure}{S}{A}{B}$ and $\starptimenew{\tmeasure}{S}{A}{B}$ are monotone with respect to $S$, $A$, $B$, and $\tmeasure$.

\item \label{elementary_fpp_properties:connectivity_part}
If $S$ is a connected subgraph of $\Zd$, then for any and $U\subseteq \Zd$ and $\vect v\in\Zd$,
\begin{enumerate}
\item $\rptime{\tmeasure}{S}{U}{\vect v} < \infty \ifandonlyif \vect v\in U$ or ($\vect v\in S$ and $U\cap S\ne \emptyset$).
\item $\starptimenew{\tmeasure}{S}{U}{\vect v} < \infty \ifandonlyif \vect v\in U$ or ($\vect v\in \nbrhood{S}$ and $U\cap \nbrhood{S} \ne\emptyset$).
\end{enumerate}

\item \label{elementary_fpp_properties:minimal_paths_part}
The infima defining $\rptmetric[\tmeasure]{S}$ and $\starptmetric[\tmeasure]{S}$ in \eqref{restricted_ptime_def_eqn} and \eqref{star_ptime_def_eqn} are unchanged if $S$-paths from $U$ to $V$ are replaced with \emph{minimal} $S$-paths from $U$ to $V$ (and similarly for $\stargraph{S}$).

\end{enumerate}
\end{lem}

The following lemma gives a stronger version of the monotonicity of $\ptmetric[\tmeasure]$ with respect to $\tmeasure$, which will be needed in some later proofs.


\begin{lem}[Monotonicity of passage times with traversal measure]
\label{ptime_tmeasure_monotone_lem}
Let $U,V,S\subseteq \Zd$, and let $\tmeasure$ and $\tmeasure'$ be traversal measures on $\edges{\Zd}$.
\begin{enumerate}
\item Suppose\ \
$\displaystyle
\restrict{\tmeasure}{\edges{S}\setminus \brackets[2]{\edges{U}\cup \edges{V}}}
\le \restrict{\tmeasure'}{\edges{S}\setminus \brackets[2]{\edges{U}\cup \edges{V}}}.
$

Then any minimal $S$-path from $U$ to $V$ satisfies $\tmeasure(\gamma)\le \tmeasure'(\gamma)$, and $\rptime{\tmeasure}{S}{U}{V} \le \rptime{\tmeasure'}{S}{U}{V}$. 

\item Suppose\ \ $\displaystyle \restrict{\tmeasure}{\staredges{S}\setminus \brackets[2]{\edges{U}\cup \edges{V}}} \le \restrict{\tmeasure'}{\staredges{S}\setminus \brackets[2]{\edges{U}\cup \edges{V}}}$.

Then any minimal $\stargraph{S}$-path from $U$ to $V$ satisfies $\tmeasure(\gamma)\le \tmeasure'(\gamma)$, and $\starptimenew{\tmeasure}{S}{U}{V} \le \starptimenew{\tmeasure'}{S}{U}{V}$.

%
\end{enumerate}
\end{lem}

\section{Almost Sure Properties of the Random Environment}
\label{fpc_basics:almost_sure_environment_sec}
In this section we show that certain nice properties are satisfied almost surely by the traversal measures $\tmeasure$ and $\tmeasurepair$ constructed in Section~\ref{fpc_basics:setup_sec}.

\subsection{The Shape Theorem and the One-Type Process}
\label{fpc_basics:one_type_a_s_sec}

Our main tool for analyzing both the one-type and two-type processes in the present chapter will be the Shape Theorem (Theorem~\ref{shape_thm}), which provides a bridge between the random environment and the deterministic motion of the processes. We restate the Shape Theorem here for easy reference:

\begin{thm}[Shape Theorem, \cite{Cox:1981aa}, \cite{Kesten:1986aa}]
\label{shape_thm1}
If $\tmeasure$ is an \iid\ traversal measure on $\edges{\Zd}$ with model traversal times satisfying \finitespeed{d}\ and \minimalmoment{d}, then there is a norm $\mu$ on $\Rd$ (the shape function) such that for any $\epsilon>0$,
\begin{equation}\label{shape_thm_eqn1}
\Pr \eventthat[3]{(1-\epsilon) t\unitball{\mu} \subseteq 
	\rreachedset{\tmeasure}{\vect 0}{\Rd}{t} \subseteq (1+\epsilon) t\unitball{\mu}
	\text{\ \ for all large } t } = 1,
\end{equation}
where $\unitball{\mu}$ is the closed unit ball of the norm $\mu$, and $\rreachedset{\tmeasure}{\vect 0}{\Rd}{t}$ is defined by \eqref{unrestricted_continuum_process_eqn}.
\end{thm}

By considering a countable sequence $\epsilon_n\to 0$ in Theorem~\ref{shape_thm1}, it follows that under the same hypotheses,
\begin{equation}\label{shape_thm_eqn2}
\Pr \eventthat[3]{\forall \epsilon>0,\; \exists t_\epsilon<\infty
	\text{ such that }
	(1-\epsilon) t\unitball{\mu} \subseteq 
	\rreachedset{\tmeasure}{\vect 0}{\Rd}{t} \subseteq (1+\epsilon) t\unitball{\mu}
	\text{ for all $t\ge t_\epsilon$}} = 1. 
\end{equation}
The event in \eqref{shape_thm_eqn2} says precisely that the random pseudometric $\ptmetric[\tmeasure]$ is asymptotic to the norm $\mu$, in the sense defined in Appendix~\ref{asymptotic_chap} (see Proposition~\ref{asymptotic_equivalence_prop}).
Once we know that $\ptmetric[\tmeasure]$ is asymptotic to a norm, there are several other useful properties of the traversal measure $\ttime$ that follow deterministically from this result.
We list these deterministic properties in the following definition. The first property, \normexists, is a restatement of the event in \eqref{shape_thm_eqn2}. The remaining properties will be useful in constructing and analyzing the one-type and two-type processes.

\begin{definition}[Deterministic properties of traversal measures]
\label{traversal_meas_properties_def}
Let $\ttime$ be a traversal measure on $\edges{\Zd}$ as defined in Section~\ref{fpc_basics:traversal_meas_sec}, with $\ptimefamily[\ttime]$ and $\fppprocfamily[\ttime]$ defined as in Section~\ref{fpc_basics:restricted_sec}. We define the following properties, which may or may not be satisfied by a given measure $\ttime$:
\begin{description}
\item[\normexists] There exists a norm $\mu$ on $\Rd$ such that  the following equivalent conditions hold:
\begin{enumerate}

\item For every $\epsilon>0$ there exists $t_0 = t_0(\epsilon)<\infty$ such that if $t\ge t_0$, then
\[
(1-\epsilon) t\unitball{\mu}
	\subseteq \rreachedset{\tmeasure}{\vect 0}{\Rd}{t}
	\subseteq (1+\epsilon) t\unitball{\mu},
\]
where $\unitball{\mu}$ is the unit $\mu$-ball centered at $\vect 0$.

\item For every $\vect u\in\Zd$ and $\epsilon>0$, there exists $t_0 = t_0(\vect u,\epsilon) <\infty$ such that if $t\ge t_0$, then
\[
\reacheddset[2]{\mu}{\vect u}{(1-\epsilon) t}
	\subseteq\rreachedset{\tmeasure}{\vect u}{\Rd}{t}
	\subseteq \reacheddset[2]{\mu}{\vect u}{(1+\epsilon) t},
\]
where $\reachedd_\mu^{\vect u}(r)$ is the $\mu$-ball of radius $r$ centered at $\vect u$.

\item For every $\vect u\in\Zd$,
$\displaystyle
\lim_{\substack{\norm{\vect v}\to\infty\\ \vect v\in\Zd}}
	\frac{1}{\norm{\vect v}}\cdot
	\abs[2]{\ptimenew{\tmeasure}{\vect u}{\vect v} - \distnew{\mu}{\vect u}{\vect v}} = 0.
$
\end{enumerate}

\item[\boundedgrowth] The following equivalent conditions hold:
\begin{enumerate}

\item The set $\reached_\ttime^{\vect 0}(t)$ is finite for all $t<\infty$.

\item The set $\reached_\ttime^{\vect u}(t)$ is finite for all $\vect u\in\Zd$ and $t<\infty$.

\item For every $\vect u\in\Zd$, 
$\displaystyle
\lim_{\substack{\norm{\vect v}\to\infty\\ \vect v\in\Zd}}
	\ptimenew{\tmeasure}{\vect u}{\vect v} = \infty.
$
\end{enumerate}

\item[\stargeodesicsexist] For any $\vect u\in \Zd$ and $A\subseteq \Zd$, the following holds: If $S$ is any connected subset of $\Zd$ such that $\nbrhood{S}$ contains both $\vect u$ and $A$, then there exists a $\starptmetric[\tmeasure]{S}$-geodesic from $\vect u$ to $A$.

\item[\geodesicsexist] For any $\vect u\in \Zd$ and $A\subseteq \Zd$, the following holds: If $S$ is any connected subset of $\Zd$ containing both $\vect u$ and $A$, then there exists a $\rptmetric[\tmeasure]{S}$-geodesic from $\vect u$ to $A$.
%
\end{description}
\end{definition}

The property \mbox{\normexists}\ stands for ``Asymptotic to a Norm," \boundedgrowth\ stands for ``Finite Speed," \stargeodesicsexist\ stands for ``Star Geodesics Exist," and \geodesicsexist\ stands for ``Geodesics Exist." Proving the equivalence of the three conditions listed under either \normexists\ or \boundedgrowth\ is straightforward, if somewhat tedious; the ``inversion formulas" in Lemma~\ref{inversion_formula_lem} can be used to relate the second and third conditions in \normexists.


The next lemma shows that the properties in Definition~\ref{traversal_meas_properties_def} are in fact listed from strongest to weakest; that is, any traversal measure satisfying one of the properties also satisfies all the subsequent properties.

\begin{lem}[Relations between deterministic properties of traversal measures]
\label{traversal_meas_properties_lem}
For any ($\sigma$-finite) traversal measure $\ttime$ on $\edges{\Zd}$, the following implications hold:
\[\normexists\implies \boundedgrowth\implies \stargeodesicsexist\implies \geodesicsexist.
\]
\end{lem}

\begin{proof}
The implication $\normexists\implies \boundedgrowth$ is obvious since the ball $(1+\epsilon) t\unitball{\mu}$ is bounded with respect to {any} norm ($\norm{\cdot}_{\ell^1}$, for example), and hence contains only a finite number of lattice points. It is also obvious that $\stargeodesicsexist\implies \geodesicsexist$, since if $\vect u$ and $A$ are contained in $S$, then a $\starptmetric{S}$-geodesic between them cannot use any vertices in $\nbrhood{S}$, and hence is a $\rptmetric{S}$-geodesic.
Thus it remains to prove $\boundedgrowth\implies \stargeodesicsexist$.

Let $S$ be a connected subset of $\Zd$ such that $\nbrhood{S}$ contains both $\vect u$ and $A$.
Since $S$ is connected, there exists some finite path $\gamma\subseteq S$ from $\vect u$ to $A$, and since $\mttime <\infty$ a.s., we have $b\definedas \ttime (\gamma)<\infty$ a.s. Then the infimum of the traversal times of all paths in $S$ from $\vect u$ to $A$ is at most $b$; call this infimum $a\le b<\infty$. Now suppose the shape theorem holds for all $\vect u\in \Zd$, and fix some $\epsilon>0$, say $\epsilon=1/2$. Then there exists $t_{\vect u}<\infty$ such that $\reached_\ttime^{\vect u}(t)\subseteq \unitball{\mu}^{\vect u} \parens[2]{(1+\epsilon)t}$ for all $t\ge t_{\vect u}$. Let $c = t_{\vect u}\vee b$. Then we have $\reached_\ttime^{\vect u}(b)\subseteq \reached_\ttime^{\vect u}(c)\subseteq \unitball{\mu}^{\vect u} \parens[2]{(1+\epsilon)c}$. Since $a$ is the infimum of traversal times from $\vect u$ to $A$, there is some sequence of paths $\gamma_n\subseteq S$ from $\vect u$ to $A$ such that $\ttime(\gamma_n)$ decreases to $a$, and we can assume that $\gamma_0 = \gamma$. Thus, all the paths $\gamma_n$ have $\ttime(\gamma_n)\le b$, and hence $\gamma_n\subseteq \reached_\ttime^{\vect u}(b)\subseteq \unitball{\mu}^{\vect u} \parens[2]{(1+\epsilon)c}$. Since this ball is bounded, it contains only a finite number of paths originating from $\vect u$.
Thus, there are only a finite number of $\gamma_n$'s, so one of them must achieve the infimum. That is, $\ttime(\gamma_n)= a$ for some $n\in\naturalzero$, in which case $\gamma_n$ is a $\starptmetric[\tmeasure]{S}$-geodesic from $\vect u$ to $A$.
\end{proof}

As noted in \eqref{shape_thm_eqn2} above, the property \normexists\ holds almost surely for any random traversal measure $\tmeasure$ satisfying the Shape Theorem. Thus, combining Theorem~\ref{shape_thm1} with Lemma~\ref{traversal_meas_properties_lem} immediately yields the following result, which will be used in the proof of Lemma~\ref{two_type_a_s_properties_lem} below.

\begin{lem}[Almost sure properties of $\ttime$]
\label{one_type_a_s_properties_lem}
Suppose $\ttime$ is a random traversal measure on $\edges{\Zd}$ which is \iid\ with model traversal time $\mttime_0$ satisfying \finitespeed{d}\ and \minimalmoment{d}. Then $\Pr$-almost surely, $\ttime$ satisfies \normexists\ and hence also \boundedgrowth, \stargeodesicsexist, and \geodesicsexist.
\end{lem}



\subsection{The Two-Type Process and Minimum Traversal Measure}
\label{fpc_basics:two_type_a_s_sec}



Now we turn to the properties of two-type traversal measures that will be desirable when constructing the two-type process in the next section.
First note that the \iid\ random two-type traversal measure $\ttimepair=(\ttime_1, \ttime_2)$  constructed in Section~\ref{fpc_basics:prob_space_sec} has a certain almost sure behavior built into it by assumption, namely the ``No Ties with Probability 1" property \notiesprob\ for the model traversal times.
Observe that 
\notiesprob\ is equivalent to the statement that if $F_1$ and $F_2$ are nonempty, disjoint, finite subsets of $\edges{\Zd}$, then $\Pr \eventthat[2]{\ttime_1(F_1) = \ttime_2(F_2)} = 0$. In particular, this holds when $F_1$ and $F_2$ are finite edge-disjoint paths of positive length. Since the number of finite paths is countable, this implies that the following deterministic ``No Ties" property holds $\Pr$-almost surely for the random traversal measure $\ttimepair$.
This incarnation of the ``No Ties" property is the one that will be relevant for constructing the competition process.

\begin{definition}[Deterministic ``No Ties" property]
\label{det_no_ties_def}
\mbox{}
\begin{description}
\item[\notiesdet] If $\gamma_1$ and $\gamma_2$ are finite, edge-disjoint paths in $\Zd$ of positive length, then $\ttime_1(\gamma_1) \ne \ttime_2(\gamma_2)$.
\end{description}
\end{definition}

We now introduce one more idea that will be important for the construction of the two-type process. Given a two-type traversal measure $\ttimepair = (\ttime_1,\ttime_2)$ on $\edges{\Zd}$, we define its \textdef{minimum traversal measure} to be $\mintime = \ttime_1\wedge \ttime_2$; that is, the $\mintime$-measure of an edge is
\begin{equation}
\label{mintime_def_eqn}
\mintime(\edge) \definedas \ttime_1(\edge)\wedge \ttime_2(\edge)
\quad\text{for } \edge\in \edges{\Zd}.
\end{equation}
The relevance of the measure $\mintime$ will become clear in the proof of Proposition~\ref{entangled_well_defined_prop} in the next section.
Namely, in order to get a well-defined two-type process, we will assume that $\mintime$ satisfies the property \boundedgrowth\ in order to guarantee that no vertex can be affected by other vertices that are infinitely far away.

The following lemma explicitly identifies the relevant deterministic properties that are satisfied almost surely by a random traversal measure $\tmeasurepair$ constructed as in Section~\ref{fpc_basics:setup_sec}. We will appeal to Lemma~\ref{two_type_a_s_properties_lem} in Section~\ref{fpc_basics:two_type_construction_sec} below to conclude that our construction of the two-type process is well-behaved almost surely.


\begin{lem}[Almost sure properties of $\ttimepair$]
\label{two_type_a_s_properties_lem}
Suppose $\ttimepair = (\ttime_1,\ttime_2)$ is a random two-type traversal measure on $\edges{\Zd}$ which is \iid\ with model traversal times $\mttimes{1}$ and $\mttimes{2}$ satisfying \finitespeed{d}, \minimalmoment{d}, and \notiesprob. Then $\Pr$-almost surely, $\ttimepair$ satisfies \notiesdet, and the three one-type traversal measures $\ttime_1$, $\ttime_2$, and $\mintime = \ttime_1\wedge \ttime_2$ all satisfy \normexists, hence also \boundedgrowth, \stargeodesicsexist, and \geodesicsexist.
\end{lem}

\begin{proof}
It was already noted above that \notiesdet\ follows $\Pr$-a.s.\ from \notiesprob. Moreover, because of the property \notiesprob, at most one of the two species can have a positive probability of crossing an edge instantaneously; that is, it is impossible to have both $\Pr \eventthat[2]{\mttimes{1}_0 = 0}>0$ and $\Pr \eventthat[2]{\mttimes{2}_0 = 0}>0$, because then we would have
\[
\Pr \eventthat[1]{\mttimes{1}_0 = \mttimes{2}_0 = 0} 
= \Pr \eventthat[1]{\mttimes{1}_0 = 0} \Pr \eventthat[1]{\mttimes{2}_0 = 0}
>0,
\]
violating \notiesprob. Therefore, if $\mttimes{\min}_0$ denotes the model traversal time for $\mintime$, i.e.\ $\mttimes{\min}_0 = \cmttimes{1}{0}\meet \cmttimes{2}{0}$, then no matter what coupling $\jointlaw$ we use for $\cmttimes{1}{0}$ and $\cmttimes{2}{0}$, it follows that
\begin{align*}
 \Pr \eventthat[2]{\mttimes{\min}_0=0}
&\le \Pr \eventthat[2]{\cmttimes{1}{0} = 0}+\Pr\eventthat[2]{\cmttimes{2}{0} = 0}\\
&= \max \setof[1]{\Pr \eventthat[2]{\mttimes{1}_0 = 0}, \Pr\eventthat[2]{\mttimes{2}_0 = 0}}
	<\pcrit.
\end{align*}
This shows that $\mttimes{\min}_0$ satisfies \finitespeed{d}, regardless of how the traversal times for an edge are coupled. Moreover, $\mttimes{\min}_0$ obviously satisfies \minimalmoment{d}  whenever either $\mttimes{1}_0$ or $\mttimes{2}_0$ does, so all three model traversal times $\mttimes{\min}_0$, $\mttimes{1}_0$, and $\mttimes{2}_0$ satisfy both \finitespeed{d}\ and \minimalmoment{d}. Thus, the traversal measures $\ttime_1$, $\ttime_2$, and $\mintime$ all satisfy \normexists\ $\Pr$-a.s.\ by Lemma~\ref{one_type_a_s_properties_lem}.
\end{proof}



\section{A New Construction of First-Passage Competition with Different Speeds}
\label{fpc_basics:two_type_construction_sec}

In this section we introduce a new geometric construction of two-type first-passage competition that works for species with different traversal measures $\tmeasure_1$ and $\tmeasure_2$, even when the starting configuration is infinite. We first briefly discuss finite starting configurations in Section~\ref{fpc_basics:finite_configs_sec}, and then move to the general case in Section~\ref{fpc_basics:infinite_configs_sec}. The new construction appears in Definition~\ref{two_type_proc_def}, which is followed by several lemmas identifying some elementary properties of the competition process. In Propositions~\ref{entangled_full_prop} and \ref{entangled_well_defined_prop} we show that the two-type process almost surely behaves in accordance with the informal description given in Section~\ref{intro:two_type_description_sec}.
In Section~\ref{fpc_basics:ips_def_sec} we describe how to interpret the two-type competition process as an interacting particle system, and in Section~\ref{fpc_basics:additional_notation_sec} we introduce some additional notation that will be needed later.

\subsection{Finite Starting Configurations}
\label{fpc_basics:finite_configs_sec}

Recall from Section~\ref{intro:two_type_description_sec} that the two-type first-passage competition process consists of two species, species 1 and species 2, racing to capture vertices in $\Zd$. The two species start from disjoint sets $\initialset{1}\subseteq\Zd$ and $\initialset{2}\subseteq\Zd$, respectively, and move according to their respective traversal times, $\ttime_1$ and $\ttime_2$, which we can view in conjunction as a two-type traversal measure $\tmeasurepair$ on $\edges{\Zd}$. We call any pair $\initconfig$ of disjoint subsets of $\Zd$ an \textdef{initial configuration} or \textdef{starting configuration} for the two-type process. We say the initial configuration $\initconfig$ is \emph{finite} if both $\initialset{1}$ and $\initialset{2}$ are finite sets, and we say that it is \emph{infinite} otherwise. By analogy with the one-type process, we write $\reachedset{\tmeasurepair}{\initconfig}{t}$ for the two-type process started from the initial configuration $\initconfig$ and using the two-type traversal measure $\tmeasurepair$. The object $\reachedfcn{\tmeasurepair}{\initconfig}$ is a pair of $2^{\Zd}$-valued functions on $[0,\infty]$, which we write as
\[
\reachedset{\tmeasurepair}{\initconfig}{t}
=\pair[1]{\reachedset{1}{\initconfig}{t}}{\,\reachedset{2}{\initconfig}{t}}_{\tmeasurepair},
\]
where for each $t\in[0,\infty]$ and $i\in\setof{1,2}$, the component $\reachedset{i}{\initconfig}{t} = \reachedset{i}{\initconfig}{t}_{\tmeasurepair}$ is the set of vertices species $i$ has reached by time $t$.
%
%
%

As mentioned previously, for finite starting configurations, determining which species arrives first at a given vertex is a ``simple" matter of constructing the process step-by-step, as done, for example, in \cite{\HPb}, \cite{Deijfen:2006aa}, or \cite{\GMdensity}. Each step consists of one infection event. At each step, we search all vertices neighboring the infected region, find the one (or more) which is infected next, and assign it the appropriate ``color" 1 or 2, depending on which species infects it. In order for this algorithmic construction to work, in addition to the finiteness assumption on $\initconfig$, we also require that $\tmeasurepair = \pair{\tmeasure_1}{\tmeasure_2}$ satisfies \notiesdet\ (so that no vertex is simultaneously infected by both species) and that each of the component measures $\tmeasure_1$ and $\tmeasure_2$ satisfies \boundedgrowth\ (so that only finitely many infections can take place in any finite time interval).

\subsection{Infinite Starting Configurations and the New Construction}
\label{fpc_basics:infinite_configs_sec}

If one or both of the starting sets $A_i$ is infinite, the above algorithmic construction breaks down in general. For example, if the traversal measure $\tmeasurepair$ is \iid\ and 0 is in the support of $\ttime_i(\edge)$, then almost surely there will be infinitely many infections in any finite time interval, and hence the notion of which vertex is infected ``next" does not make sense.
In the special case that $\tmeasure_1 = \tmeasure_2 = \tmeasure$ (see e.g.\ \cite{\GM} or \cite{Hoffman:2005aa}), this problem is avoided because the regions finally conquered by the two species are precisely the Voronoi cells for the configuration $\initconfig$ in the pseudometric $\ptmetric[\tmeasure]$, and the occupied regions at time $t$ can be identified explicitly as
\begin{equation}
\label{two_type_same_times_def}
\begin{split}
\reachedset{1}{\initconfig}{t}_{\pair{\tmeasure}{\tmeasure}}
&= A_1\cup \setbuilder[1]{\vect v\in \Zd}
	{\ptimenew{\tmeasure}{A_1}{\vect v}< \ptimenew{\tmeasure}{A_2}{\vect v}
	\text{ and } \ptimenew{\tmeasure}{A_1}{\vect v}\le t},
	\\
\reachedset{2}{\initconfig}{t}_{\pair{\tmeasure}{\tmeasure}}
&= A_2\cup \setbuilder[1]{\vect v\in \Zd}
	{\ptimenew{\tmeasure}{A_2}{\vect v}< \ptimenew{\tmeasure}{A_1}{\vect v}
	\text{ and } \ptimenew{\tmeasure}{A_2}{\vect v}\le t}.
\end{split}
\end{equation}
The description in \eqref{two_type_same_times_def} works only in the case where the two species use the same set of traversal times $\tmeasure$, because this guarantees that the geodesics for the two species cannot ``interfere" with each other, i.e.\ the two conquered sets as defined in \eqref{two_type_same_times_def} are disjoint, whereas there would be no such guarantee if two different collections of traversal times were used. However, we now make a simple observation about the two-type process that is reminiscent of \eqref{two_type_same_times_def} and will provide a method of defining the process more generally.

\begin{prop}[Conquering property of the two-type process]
\label{conquering_property_prop}
Consider a two-type process $\reachedset{\tmeasurepair}{\initconfig}{t}$, started from a nonempty initial configuration $\initconfig$ and using the traversal measure $\tmeasurepair = \pair{\tmeasure_1}{\tmeasure_2}$.
If $S$ is a subset of $\Zd$ such that 
\begin{equation}
\label{conquering_property_eqn}
S\cap A_2=\emptyset \quad \text{and}\quad
\rptime{\tmeasure_1}{S}{A_1}{\vect v}<\starptimenew{\tmeasure_2}{(\Zd\setminus S)}{A_2}{\vect v}
\text{ for all } \vect v\in \bd S\setminus A_1,
\end{equation}
then {species 1} conquers $S$ in the two-type process, i.e.\ $S\subseteq \reachedset{1}{\initconfig}{\infty}_{\tmeasurepair}$. A symmetric statement holds for species 2.
\end{prop}


We refer to any set $S$
satisfying the conditions in \eqref{conquering_property_eqn} as a \textdef{conquering set for species~1}, and we make the analogous definition for species~2. Implicit in the statement of Proposition~\ref{conquering_property_prop} is that the process $\reachedset{\tmeasurepair}{\initconfig}{t}$ is well-defined, i.e.\ we can unambiguously determine the conquered sets $\reachedset{i}{\initconfig}{t}_{\tmeasurepair}$ at each $t\ge 0$; for example, this is the case when the initial configuration $\initconfig$ is finite, $\tmeasurepair$ satisfies \notiesdet, and each $\tmeasure_i$ satisfies \boundedgrowth, as noted above.
On the other hand, Proposition~\ref{conquering_property_prop} provides a condition that we should expect the two-type process $\reachedset{\tmeasurepair}{\initconfig}{t}$ to possess whenever it does happen to be defined, and it therefore gives us a way to \emph{define} the process in general. Namely, the region finally conquered by species~1 should be the largest $S\subseteq\Zd$ satisfying
\eqref{conquering_property_eqn}, and similarly for species~2. Once the finally conquered regions are identified, the evolution of the process can be defined by restricting the growth of each species to its finally conquered set. We formalize this in the following definition.

\begin{definition}[Two-type first-passage competition]
\label{two_type_proc_def}
Let $\initialset{1}$ and $\initialset{2}$ be disjoint subsets of $\Zd$, let $\tmeasurepair = \pair{\tmeasure_1}{\tmeasure_2}$ be a two-type traversal measure on $\edges{\Zd}$, and set
\begin{align*}
\finalset{1} &\definedas \bigcup \setbuilder[1]{S\subseteq \Zd\setminus \initialset{2}}
	{
	\rptime{\tmeasure_1}{S}{\initialset{1}}{\vect v}<
	\starptimenew{\tmeasure_2}{(\Zd\setminus S)}{\initialset{2}}{\vect v}
	\text{ for all } \vect v\in \bd S\setminus \initialset{1}},\\
\finalset{2} &\definedas \bigcup \setbuilder[1]{S\subseteq \Zd\setminus \initialset{1}}
	{
	\rptime{\tmeasure_2}{S}{\initialset{2}}{\vect v}<
	\starptimenew{\tmeasure_1}{(\Zd\setminus S)}{\initialset{1}}{\vect v}
	\text{ for all } \vect v\in \bd S\setminus \initialset{2}},
\end{align*}
if at least one of $\initialset{1}$, $\initialset{2}$ is nonempty, and $\finalset{1}=\finalset{2} \definedas\emptyset$ otherwise.
We define the \textdef{two-type competition process} or \textdef{entangled process} with starting configuration $\initconfig$ to be the pair of $2^\Zd$-valued functions
\[
\reachedset{\tmeasurepair}{\initconfig}{t}
=\pair[1]{\reachedset{1}{\initconfig}{t}}{\,\reachedset{2}{\initconfig}{t}}_{\tmeasurepair},
\]
defined for $t\in [0,\infty]$ by
\[
\reachedset{1}{\initconfig}{t} \definedas
	\rreachedset{\tmeasure_1}{\initialset{1}}{\finalset{1}}{t}
\quad\text{and}\quad
\reachedset{2}{\initconfig}{t} \definedas
	\rreachedset{\tmeasure_2}{\initialset{2}}{\finalset{2}}{t}
\quad\text{for } t\ge 0,
\]
and
\[
\reachedset{i}{\initconfig}{\infty} := \bigcup_{t\ge 0} \reachedset{i}{\initconfig}{t}
\quad\text{for $i\in\setof{1,2}$}.
\]
We call the sets $\finalset{1}$ and $\finalset{2}$ the \textdef{finally conquered sets} for species 1 and 2, respectively, and we refer to $\reachedset{i}{\initconfig}{t}$ as the set of vertices \textdef{occupied} by species $i$ at time $t$. We say that the entangled process is \textdef{well-defined} if $\finalset{1}\cap \finalset{2} = \emptyset$ and \textdef{full} if $\finalset{1}\cup \finalset{2} = \Zd$.
\end{definition}

Note that $\finalset{i}$ is precisely the union of all the conquering sets for species $i$ as defined by Proposition~\ref{conquering_property_prop}; we will sometimes refer to species $i$'s conquering sets as \textdef{$\finalset{i}$-sets} for short. Lemma~\ref{final_sets_properties_lem} below shows that $\finalset{i} = \reachedset{i}{\initconfig}{\infty}$, so that the name ``finally conquered set" is appropriate; moreover, it then follows immediately that Proposition~\ref{conquering_property_prop} does in fact hold for the two-type process defined in Definition~\ref{two_type_proc_def}. Note that Definition~\ref{two_type_proc_def} makes sense for \emph{every} two-type traversal measure $\tmeasurepair$ and initial configuration $\initconfig$, because the restricted one-type process is defined for any one-type traversal measure. However, if we allow arbitrary initial configurations and traversal measures, it is possible that $\finalset{1}\cap \finalset{2} \ne \emptyset$ and/or $\finalset{1}\cup \finalset{2} \ne \Zd$, so the two-type process is neither well-defined nor full in general, meaning that it may not behave in accordance with the informal description given in Section~\ref{intro:two_type_description_sec}. We will later investigate (in Lemma~\ref{trivial_init_configs_lem} and Propositions~\ref{entangled_full_prop} and \ref{entangled_well_defined_prop}) sufficient conditions on $\initconfig$ and $\tmeasurepair$ for the entangled process to be well-defined and full. First we enumerate some elementary properties of the finally conquered sets $\finalset{1}$ and $\finalset{2}$ in Lemma~\ref{final_sets_properties_lem}, which will be used in the proofs of several results below; Lemma~\ref{final_sets_properties_lem} is proved in Section~\ref{leftover:fpc_basics_sec} of Appendix~\ref{leftover_chap}.

\begin{lem}[Properties of the finally conquered sets]
\label{final_sets_properties_lem}
For any two-type traversal measure $\tmeasurepair$ and initial configuration $\initconfig$, the finally conquered sets $\finalset{1}$ and $\finalset{2}$ in Definition~\ref{two_type_proc_def} satisfy the following for $i\in\setof{1,2}$.
\begin{enumerate}
\item \label{final_sets:initial_sets_part}
$\finalset{i}\supseteq \initialset{i}$ and $ \finalset{i} \cap \initialset{3-i} = \emptyset$. Moreover, $\finalset{i}$ is empty if and only if $\initialset{i}$ is empty.
\item \label{final_sets:characterization_part}
$\vect v\in \finalset{i}$ if and only if $\vect v\in A_i\text{ or } \rptime{i}{\finalset{i}}{A_i}{\vect v} < \starptimenew{3-i}{(\Zd\setminus \finalset{i})}{A_{3-i}}{\vect v}$.
\item \label{final_sets:components_part}
Every component of $\finalset{i}$ contains a component of $\initialset{i}$.
\item \label{final_sets:are_union_part}
$\finalset{i} = \reachedset{i}{\initconfig}{\infty}$.
\end{enumerate}
\end{lem}

We refer to the two-type process as ``entangled" to emphasize that its evolution depends simultaneously on both components of the traversal measure, $\tmeasure_1$ and $\tmeasure_2$. By contrast, we will refer to the passage time families $\ptimefamily[\tmeasure_i]$ and $\starptimefamily[\tmeasure_i]$ ($i\in\setof{1,2}$) as \textdef{disentangled passage times} for $\tmeasurepair = \pair{\tmeasure_1}{\tmeasure_2}$, and to the family of one-type processes $\fppprocfamily[\tmeasure_i]$ as the \textdef{disentangled processes} corresponding to the entangled process $\reachedfcn{\tmeasurepair}{\initconfig}$, to emphasize that we are forgetting about one of the components of $\tmeasurepair$. When the traversal measure $\tmeasurepair$ is clear from context, we will frequently omit it from the notation for the disentangled passage times, instead writing
\begin{equation}
\rptmetric[i]{S} \definedas \rptmetric[\tmeasure_i]{S}
\quad\text{and}\quad
\starptmetric[i]{S} \definedas \starptmetric[\tmeasure_i]{S}
\quad\text{for } S\subseteq\Zd.
\end{equation}
From a probabilistic standpoint, if $\tmeasurepair$ is a random two-type traversal measure, the disentangled process $\rreachedfcn{\tmeasure_i}{A}{S}$ is a one-type process coupled to use the same traversal times as one of the species in $\reachedfcn{\tmeasurepair}{\initconfig}$.
In general we have the following relationships between each species and its disentangled version.

\begin{lem}[Comparison of entangled and disentangled processes]
\label{entangled_disentangled_lem}
\label{entangled_disentangled_comparison_lem}
\label{conquered_set_growth_lem}
Let $\tmeasurepair = \pair{\tmeasure_1}{\tmeasure_2}$ be any two-type traversal measure, let $\initconfig$ be any initial configuration, and let $i\in\setof{1,2}$.
\begin{enumerate}
\item \label{entangled_disentangled:unrestricted_part}
In the two-type process $\reachedfcn{\tmeasurepair}{\initconfig}$, each species is dominated by its unrestricted disentangled version. That is, $\reachedset{i}{\initconfig}{t} \subseteq \reachedset{\tmeasure_i}{\initialset{i}}{t}$ for all $t\in [0,\infty]$.

\item \label{entangled_disentangled:restricted_part}
If species $i$ conquers some set $S\subseteq \Zd$ in the two-type process $\reachedfcn{\tmeasurepair}{\initconfig}$, then
the entangled growth of species~$i$ dominates the growth of its disentangled version restricted to $S$. That is, $\reachedset{i}{\initconfig}{t} \supseteq \rreachedset{\tmeasure_i}{\initialset{i}}{S}{t}$ for all $t\in [0,\infty]$.
\end{enumerate}
\end{lem}

\begin{proof}
Both parts follow directly from Definition~\ref{two_type_proc_def} and the monotonicity of the restricted first-passage percolation process with respect to the restricting set. More explicitly, for Part 1, we have
\[
\reachedset{i}{\initconfig}{t} 
=\rreachedset{\tmeasure_i}{\initialset{i}}{\finalset{i}}{t}
\subseteq \rreachedset{\tmeasure_i}{\initialset{i}}{\Zd}{t}
= \reachedset{\tmeasure_i}{\initialset{i}}{t}
\quad\text{for all $t\in [0,\infty]$,}
\]
where $\finalset{i}$ is species $i$'s finally conquered set from Definition~\ref{two_type_proc_def}.
For Part 2, the hypothesis is that $S\subseteq \reachedset{i}{\initconfig}{\infty}$, and Lemma~\ref{final_sets_properties_lem} shows that $\reachedset{i}{\initconfig}{\infty} = \finalset{i}$, so $S\subseteq \finalset{i}$. By Definition~\ref{two_type_proc_def}, $\reachedset{i}{\initconfig}{t} = \rreachedset{\tmeasure_i}{\initialset{i}}{\finalset{i}}{t}$ for all $t\in [0,\infty]$, and by monotonicity, we have $\rreachedset{\tmeasure_i}{\initialset{i}}{\finalset{i}}{t}\supseteq \rreachedset{\tmeasure_i}{\initialset{i}}{S}{t}$ since $\finalset{i}\supseteq S$.
\end{proof}

In fact, we show in Lemma~\ref{one_type_two_type_lem} below that under some additional assumptions on $\tmeasurepair$, a stronger version of Lemma~\ref{entangled_disentangled_lem} holds, in which each species' entangled growth is in fact \emph{equal to} any member of an appropriate family of disentangled processes up to a given time.
The next lemma gives two trivial situations where the entangled process is easy to describe explicitly.

\begin{lem}[Trivial starting configurations]
\label{trivial_init_configs_lem}
\label{absent_species_lem} 
Let $\tmeasurepair = \pair{\tmeasure_1}{\tmeasure_2}$ be any two-type traversal measure, and let $\initconfig$ be an initial configuration.
\begin{enumerate}
\item If $\initialset{i} = \emptyset$ for some $i\in\setof{1,2}$, and $j=3-i$, then
\[
\reachedset{i}{\initconfig}{t} = \emptyset
\quad\text{and}\quad
\reachedset{j}{\initconfig}{t} = \reachedset{\tmeasure_j}{\initialset{j}}{t}
\quad\text{ for all $t\ge 0$}.
\]

\item If $\initialset{1}\cup \initialset{2} = \Zd$, then $\reachedset{\tmeasurepair}{\initconfig}{t} = \initconfig$ for all $t\ge 0$.
\end{enumerate}
Therefore,  if at least one of the sets $\initialset{1}$, $\initialset{2}$, or $\Zd\setminus \parens{\initialset{1}\cup \initialset{2}}$ is empty, then the two-type process with initial configuration $\initconfig$ is well-defined and full for all two-type traversal measures $\tmeasurepair$.
\end{lem}

\begin{proof}
Both parts follow immediately from Definition~\ref{two_type_proc_def}. For Part 1, note that if $\initialset{i} =\emptyset$ and $\initialset{3-i}\ne\emptyset$, then $\finalset{i} = \emptyset$ and $\finalset{3-i} = \Zd$, and if $\initialset{1}=\initialset{2}=\emptyset$, then $\finalset{1} = \finalset{2} = \emptyset$.
For Part 2, note that $\initialset{i} \subseteq \finalset{i} \subseteq \Zd\setminus \initialset{3-i}$, so if $\initialset{1}\cup \initialset{2} = \Zd$, then we must have $\finalset{i} = \initialset{i}$.
\end{proof}

The following proposition gives sufficient conditions for the entangled process to eventually cover the whole lattice. The proof uses
an elementary result stated in Section~\ref{leftover:fpc_basics_sec}.

\begin{prop}[Sufficient conditions for a full entangled process]
\label{entangled_full_prop}
If $\initconfig$ is any nonempty initial configuration and $\tmeasurepair$ is a two-type traversal measure on $\edges{\Zd}$ satisfying \notiesdet\ and whose component measures satisfy \stargeodesicsexist, then $\finalset{1}\cup \finalset{2} = \Zd$, where $\finalset{1}$ and $\finalset{2}$ are the finally conquered sets for the entangled process $\reachedfcn{\tmeasurepair}{\initconfig}$.
\end{prop}

\begin{proof}
\newcommand{\compset}{S}
\newcommand{\badcomponent}{S_0}
\newcommand{\badvertex}{\vect v_0}
\newcommand{\genvertex}{\vect v}

First note that the assumption that the initial configuration $\initconfig$ is nonempty implies that any $\finalset{i}$-set is contained in $\finalset{i}$; we will use this fact repeatedly throughout the proof.
To prove that $\finalset{1}\cup \finalset{2} = \Zd$, we let $\compset \definedas \Zd\setminus \finalset{1}$, and then we show that $\compset \subseteq \finalset{2}$ by proving that $\compset$ is a $\finalset{2}$-set.

First we show by contradiction that
\begin{equation}
\label{entangled_full:finite_eqn}
\rptime{2}{S}{\initialset{2}}{\vect v}<\infty
\quad\text{for all $\vect v\in S$.}
\end{equation}
To see that \eqref{entangled_full:finite_eqn} holds, we make the following claim.

\begin{claim}
\label{entangled_full:finite_claim}

If $\rptime{2}{S}{\initialset{2}}{\badvertex}=\infty$ for some $\badvertex\in\compset$, and $\badcomponent$ is the component of $\compset$ containing $\badvertex$, then $\finalset{1}\cup \badcomponent$ is a $\finalset{1}$-set.
\end{claim}

\begin{proof}[Proof of Claim~\ref{entangled_full:finite_claim}]
\newcommand{\newfset}{\widehat{\finalsetsymb}_{1}}

Let $\newfset \definedas \finalset{1}\cup \badcomponent$. First, we claim that $\bd \newfset \subseteq \bd \finalset{1}$. To see this, note that by Lemma~\ref{bd_nbr_lem}, we have
\begin{equation}
\label{entangled_full:bd_eqn}
\bd \newfset = \bd (\finalset{1}\cup \badcomponent) =
\setbuilder[2]{\genvertex\in \bd\finalset{1}\cup \bd\badcomponent}
	{\nbrset{\genvertex}\not\subseteq \finalset{1}\cup\badcomponent}.
\end{equation}
Thus, to show that $\bd \newfset \subseteq \bd \finalset{1}$, it will suffice to show that for any $\genvertex\in \badcomponent$ we have $\nbrset{\genvertex}\subseteq \finalset{1}\cup\badcomponent$. Let $\genvertex\in \badcomponent$. Then $\nbrset{\genvertex}\setminus \badcomponent \subseteq \Zd\setminus \badcomponent = \parens{\compset\setminus \badcomponent} \cup \finalset{1}$. Now if $\nbrset{\genvertex} \cap \parens{\compset\setminus \badcomponent} \ne \emptyset$, then $\badcomponent$ would be connected to some point in $\compset\setminus \badcomponent$, contradicting the maximality of the component $\badcomponent$, so we must have $\nbrset{\genvertex} \cap \parens{\compset\setminus \badcomponent} = \emptyset$. Hence $\nbrset{\genvertex}\setminus \badcomponent \subseteq \finalset{1}$, or $\nbrset{\genvertex}\subseteq \finalset{1}\cup\badcomponent$ for every $\genvertex\in \badcomponent$. Therefore, \eqref{entangled_full:bd_eqn} implies that $\bd \newfset\cap \badcomponent = \emptyset$ and hence $\bd \newfset\subseteq \bd \finalset{1}$ as claimed.

Now, since $\bd \newfset\subseteq \bd \finalset{1}$ and $\newfset\supseteq \finalset{1}$, it follows from Part~\ref{final_sets:characterization_part} of Lemma~\ref{final_sets_properties_lem} that
\begin{multline}
\label{entangled_full:bdvertex_ineq_eqn}
\forall \genvertex\in \bd \newfset\setminus 
	\parens[10]{\initialset{1}}, \quad\\
\rptime{1}{\newfset}{\initialset{1}}{\genvertex}
\le \rptime{1}{\finalset{1}}{\initialset{1}}{\genvertex}
<  \starptimenew{2}{(\Zd\setminus \finalset{1})}{\initialset{2}}{\genvertex}
\le \starptimenew{2}{(\Zd\setminus \newfset)}{\initialset{2}}{\genvertex}.
\end{multline}
Finally, note that since $\badcomponent$ is connected and contains a vertex $\badvertex$ with $\rptime{2}{\badcomponent}{\initialset{2}}{\badvertex}=\infty$, 
Part~\ref{elementary_fpp_properties:connectivity_part} of Lemma~\ref{elementary_fpp_properties_lem}
implies that $\initialset{2}\cap \badcomponent =\emptyset$.
Since $\initialset{2}\cap \finalset{1}=\emptyset$ by definition, we therefore have $\initialset{2}\cap \newfset = \emptyset$, and so \eqref{entangled_full:bdvertex_ineq_eqn} shows that $\newfset= \finalset{1}\cup \badcomponent$ is a $\finalset{1}$-set.
\end{proof}

By Claim~\ref{entangled_full:finite_claim}, if $\rptime{2}{S}{\initialset{2}}{\badvertex}=\infty$ for some $\badvertex\in\compset$, then we get $\badvertex\in \finalset{1}$ (since $\badvertex$ is contained in some $\badcomponent\subseteq \compset$ such that $\finalset{1}\cup \badcomponent$ is a $\finalset{1}$-set), but this is a contradiction since $\compset \cap \finalset{1} = \emptyset$ by definition. Thus we conclude that \eqref{entangled_full:finite_eqn} holds.

Next, we again use contradiction to show that
\begin{equation}
\label{entangled_full:ineq_eqn}
\rptime{2}{\compset}{\initialset{2}}{\genvertex}
\le \starptimenew{1}{(\Zd\setminus \compset)}{\initialset{1}}{\genvertex}
\quad\text{for all $\genvertex\in \compset$.}
\end{equation}
To see that \eqref{entangled_full:ineq_eqn} holds, we make the following claim.

\begin{claim}
\label{entangled_full:ineq_claim}
If $\starptimenew{1}{(\Zd\setminus \compset)}{\initialset{1}}{\badvertex}< \rptime{2}{\compset}{\initialset{2}}{\badvertex}$ for some $\badvertex\in \compset$, then $\finalset{1}\cup \setof{\badvertex}$ is a $\finalset{1}$-set.
\end{claim}

\begin{proof}[Proof of Claim~\ref{entangled_full:ineq_claim}]
\newcommand{\newfset}{\widehat{\finalsetsymb}_{1}}
\newcommand{\newcompset}{(\Zd\setminus \newfset)}

First note that $\Zd\setminus \compset = \finalset{1}$, so the hypothesis is that there exists $\badvertex\in \compset$ with $\starptimenew{1}{\finalset{1}}{\initialset{1}}{\badvertex}< \rptime{2}{\compset}{\initialset{2}}{\badvertex}$. This implies that $\rptime{2}{\compset}{\initialset{2}}{\badvertex}>0$, so $\badvertex\not\in \initialset{2}$, and $\starptimenew{1}{\finalset{1}}{\initialset{1}}{\badvertex}<\infty$, so $\badvertex\in \nbrset{\finalset{1}} = \bd \compset$.
Let $\newfset \definedas \finalset{1}\cup \setof{\badvertex}$.
We claim that
\begin{equation}
\label{entangled_full:equal_times_eqn}
\rptime{1}{\newfset}{\initialset{1}}{\badvertex}
= \starptimenew{1}{\finalset{1}}{\initialset{1}}{\badvertex}
\quad\text{and}\quad
\starptimenew{2}{\newcompset}{\initialset{2}}{\badvertex} = \rptime{2}{\compset}{\initialset{2}}{\badvertex}.
\end{equation}
To see that \eqref{entangled_full:equal_times_eqn} holds, first observe that any  $\stargraph{\finalset{1}}$-path from $\initialset{1}$ to $\badvertex$ is a $\newfset$-path from $\initialset{1}$ to $\badvertex$, and any minimal $\newfset$-path from $\initialset{1}$ to $\badvertex$ is a $\stargraph{\finalset{1}}$-path from $\initialset{1}$ to $\badvertex$. Similarly, noting that $\Zd\setminus \newfset = \compset \setminus \setof{\badvertex}$, observe that any $\stargraph[0]{\Zd\setminus \newfset}$-path from $\initialset{2}$ to $\badvertex$ is an $\compset$-path from $\initialset{2}$ to $\badvertex$, and any minimal $\compset$-path from $\initialset{2}$ to $\badvertex$ is a $\stargraph[0]{\Zd\setminus \newfset}$-path from $\initialset{2}$ to $\badvertex$. Thus, Part~\ref{elementary_fpp_properties:minimal_paths_part} of Lemma~\ref{elementary_fpp_properties_lem} implies that \eqref{entangled_full:equal_times_eqn} holds.
%
%
%
Combining \eqref{entangled_full:equal_times_eqn} with the hypothesis of Claim~\ref{entangled_full:ineq_claim}, we have
\begin{equation}
\label{entangled_full:badvertex_ineq_eqn}
\rptime{1}{\newfset}{\initialset{1}}{\badvertex}
= \starptimenew{1}{\finalset{1}}{\initialset{1}}{\badvertex}
< \rptime{2}{\compset}{\initialset{2}}{\badvertex}
= \starptimenew{2}{(\Zd\setminus \newfset)}{\initialset{2}}{\badvertex}.
\end{equation}
On the other hand, 
combining the fact that $\newfset\supseteq \finalset{1}$ with Part~\ref{final_sets:characterization_part} of Lemma~\ref{final_sets_properties_lem}, we have
\begin{multline}
\label{entangled_full:genvertex_ineq_eqn}
\forall \genvertex\in \newfset\setminus 
	\parens[2]{\initialset{1}\cup \setof{\badvertex}}, \quad\\
\rptime{1}{\newfset}{\initialset{1}}{\genvertex}
\le \rptime{1}{\finalset{1}}{\initialset{1}}{\genvertex}
<  \starptimenew{2}{(\Zd\setminus \finalset{1})}{\initialset{2}}{\genvertex}
\le \starptimenew{2}{\newcompset}{\initialset{2}}{\genvertex}.
\end{multline}
Combining \eqref{entangled_full:badvertex_ineq_eqn} and \eqref{entangled_full:genvertex_ineq_eqn}, we have
\[
\rptime{1}{\newfset}{\initialset{1}}{\genvertex}
< \starptimenew{2}{\newcompset}{\initialset{2}}{\genvertex}
\quad\text{for all $\genvertex\in \newfset\setminus \initialset{1}$}.
\]
Finally, since $\badvertex\notin \initialset{2}$ and $\finalset{1}\cap\initialset{2}=\emptyset$, we also have $\newfset\cap \initialset{2} = \emptyset$, so this shows that $\newfset = \finalset{1}\cup \setof{\badvertex}$ is a $\finalset{1}$-set.
\end{proof}

Claim~\ref{entangled_full:ineq_claim} immediately implies that if $\starptimenew{1}{(\Zd\setminus \compset)}{\initialset{1}}{\badvertex}< \rptime{2}{\compset}{\initialset{2}}{\badvertex}$ for some $\badvertex\in \compset$, then $\badvertex\in \finalset{1}$, which is a contradiction since $\compset \cap \finalset{1} = \emptyset$ by definition. Therefore we conclude that \eqref{entangled_full:ineq_eqn} holds.

Now, since $\tmeasurepair$ satisfies \notiesdet\ and \stargeodesicsexist, Lemma~\ref{conq_set_equiv_conditions_lem} implies that the conditions \eqref{entangled_full:finite_eqn} and \eqref{entangled_full:ineq_eqn} together are equivalent to
\[
\rptime{2}{\compset}{\initialset{2}}{\genvertex}< 
\starptimenew{1}{(\Zd\setminus \compset)}{\initialset{1}}{\genvertex}
\quad \forall \genvertex\in \compset\setminus \initialset{2}.
\]
Therefore, $\compset$ is a $\finalset{2}$-set and hence $\compset \subseteq \finalset{2}$. Since $\compset = \Zd\setminus \finalset{1}$, we thus have $\finalset{1}\cup \finalset{2} = \Zd$.
\end{proof}

We now turn to the question of when the entangled process is well-defined, i.e.\ when the finally conquered sets $\finalset{1}$ and $\finalset{2}$ are disjoint. It was noted above that in order to construct the two-type process $\reachedfcn{\tmeasurepair}{\initconfig}$ algorithmically from a finite starting configuration $\initconfig$, we had to assume that the two-type traversal measure $\tmeasurepair = \pair{\tmeasure_1}{\tmeasure_2}$ satisfied the ``No Ties" property \notiesdet, and that its two component measures $\tmeasure_1$ and $\tmeasure_2$ satisfied the ``Finite Speed" property \boundedgrowth. It turns out that this second condition is not quite sufficient to get a well-defined process for arbitrary infinite starting configurations.
Instead, we will impose the following stronger condition on $\tmeasurepair$:

\begin{definition}[``Finite Speed" property for minimum traversal measure]
\label{mintimebbdgrowth_def}
\mbox{}
\begin{description}
\item[\mintimebddgrowth] The minimum traversal measure $\mintime = \ttime_1\wedge \ttime_2$ satisfies \boundedgrowth.
\end{description}
\end{definition}

The abbreviation \mintimebddgrowth\ stands for ``Minimum Traversal (Measure) has Finite Speed." Observe that if $\tmeasurepair$ satisfies \mintimebddgrowth, then the two traversal measures $\tmeasure_1$ and $\tmeasure_2$ both satisfy \boundedgrowth, since obviously $\reachedset{\tmeasure_i}{\vect 0}{t}\subseteq \reachedset{\mintime}{\vect 0}{t}$ for $i\in\setof{1,2}$ and $t\ge 0$. 
By Lemma~\ref{two_type_a_s_properties_lem}, any \iid\ random traversal measure $\ttimepair$ constructed as in Section~\ref{fpc_basics:prob_space_sec} almost surely satisfies both \notiesdet\ and \mintimebddgrowth, and we prove in the following lemma that these two properties together are sufficient to get a well-defined process.


\begin{prop}[Sufficient conditions for a well-defined entangled process]
\label{entangled_well_defined_prop}
If $\tmeasurepair$ is a two-type traversal measure on $\edges{\Zd}$ which satisfies \notiesdet\ and \mintimebddgrowth, then the two-type process $\reachedfcn{\tmeasurepair}{\initconfig}$ is well-defined for all initial configurations $\initconfig$. In particular, this holds almost surely for any random \iid\ $\tmeasurepair$ satisfying the hypotheses in Lemma~\ref{two_type_a_s_properties_lem}.
\end{prop}

\begin{proof}
Lemma~\ref{absent_species_lem} showed that the process is well-defined if either $\initialset{1}$ or $\initialset{2}$ is empty, so assume that $\initialset{1}$ and $\initialset{2}$ are both nonempty. We will prove that $\finalset{1}\cap \finalset{2} = \emptyset$ by contradiction. The strategy will be to assume that the intersection is nonempty, and then use this assumption to construct an infinite path $\gamma$ with $\mintime(\gamma)<\infty$, contradicting \mintimebddgrowth.

Suppose there exists some vertex $\vect w\in \finalset{1}\cap \finalset{2}$. Let $P_0$ be a $\rptmetric[1]{\finalset{1}}$-geodesic from $\initialset{1}$ to $\vect w$, and let $\vect u_0$ be the first point in $P_0$ that's also in $\finalset{2}$ (where ``first" refers to the order in which the vertices in $P_0$ are traversed, starting in $\initialset{1}$ and ending at $\vect w$; this ordering is well-defined since $P_0$ is finite and simple).
Note that $\vect u_0$ exists because $\vect w\in \finalset{2}$, and we must have $\vect u_0\in \finalset{1}\cap \finalset{2}$ since $P_0\subseteq \finalset{1}$. Now let $Q_0$ be a $\rptmetric[2]{\finalset{2}}$-geodesic from $\initialset{2}$ to $\vect u_0$, and let $\vect v_0$ be the first point in $Q_0$ that's also in $\finalset{1}$ (again, ``first" refers to the ordering of the vertices in $Q_0$ induced by its direction, and $\vect v_0\in \finalset{1}\cap\finalset{2}$ exists because the end vertex $\vect u_0$ is in $\finalset{1}$).

We continue this procedure inductively for all $n\in\N$, constructing two sequences of geodesics $P_n$, $Q_n$ with intermediate points $\vect u_n$, $\vect v_n$. For $n\ge 1$,
\begin{itemize}
\item Let $P_n$ be a $\rptmetric[1]{\finalset{1}}$-geodesic from $\initialset{1}$ to $\vect v_{n-1}$, and let $\vect u_n$ be the first vertex of $P_n$ lying in $\finalset{2}$. The geodesic $P_n$ exists because (by induction) $\vect v_{n-1}\in \finalset{1}$, and every component of $\finalset{1}$ contains a component of $\initialset{1}$; the vertex $\vect u_n$ exists because (by induction) $\vect v_{n-1}\in \finalset{2}$.
\item Let $Q_n$ be a $\rptmetric[2]{\finalset{2}}$-geodesic from $\initialset{2}$ to $\vect u_n$, and let $\vect v_n$ be the first vertex of $Q_n$ lying in $\finalset{1}$. The geodesic $Q_n$ exists because $\vect u_n\in \finalset{2}$, and every component of $\finalset{2}$ contains a component of $\initialset{2}$; the vertex $\vect v_n$ exists because $\vect u_n\in \finalset{1}$.
\end{itemize}
Our goal will be to construct an infinite path $\gamma$ with $\mintime(\gamma)<\infty$ by gluing together subpaths of $P_n$ and $Q_n$. In order to know that this works, we need to know several things about $P_n$ and $Q_n$, which we enumerate in the following claims.

Let $\subpath{P_n}{\initialset{1}}{\vect u_n}$ and $\subpath{P_n}{\vect u_n}{\vect v_{n-1}}$ denote the subpaths of $P_n$ from $\initialset{1}$ to $\vect u_n$, and from $\vect u_n$ to $\vect v_{n-1}$, respectively. Similarly, let $\subpath{Q_n}{\initialset{2}}{\vect v_n}$ and $\subpath{Q_n}{\vect v_n}{\vect u_n}$ denote the appropriate subpaths of $Q_n$. 

\begin{claim}
\label{entangled_well_defined:star_path_claim}
For $n\ge 0$, $\subpath{P_n}{\initialset{1}}{\vect u_n}\cap \finalset{2} = \setof{\vect u_n}$, and $\subpath{Q_n}{\initialset{2}}{\vect v_n}\cap \finalset{1} = \setof{\vect v_n}$.
\end{claim}

This follows because $\vect u_n$ is the \emph{first} vertex of $P_n$ in $\finalset{2}$, and $\vect v_n$ is the \emph{first} vertex of $Q_n$ in $\finalset{1}$. In particular, Claim~\ref{entangled_well_defined:star_path_claim} implies that $\subpath{P_n}{\initialset{1}}{\vect u_n}$ is a $\stargraph[0]{\Zd\setminus \finalset{2}}$-path, hence is edge-disjoint from $\finalset{2}$, and that $\subpath{Q_n}{\initialset{2}}{\vect v_n}$ is a $\stargraph[0]{\Zd\setminus \finalset{1}}$-path, hence is edge-disjoint from $\finalset{1}$.

Using Claim~\ref{entangled_well_defined:star_path_claim}, we prove the following inequalities for the traversal times of $P_n$ and $Q_n$, which are the crux of the argument:

\begin{claim}
\label{entangled_well_defined:time_inequality_claim}
For $n\ge 0$, $\ttime_1\parens[2]{\subpath{P_n}{\initialset{1}}{\vect u_n}} > \ttime_2(Q_n)$, and $\ttime_2\parens[2]{\subpath{Q_n}{\initialset{2}}{\vect v_n}} > \ttime_1(P_{n+1})$.
\end{claim}

Essentially, Claim~\ref{entangled_well_defined:time_inequality_claim} holds because in order for the vertex $\vect u_n$ to be in $\finalset{2}$, species 2 has to reach $\vect u_n$ from within $\finalset{2}$ before species 1 can reach it from outside of $\finalset{2}$; similarly, in order for the vertex $\vect v_n$ to be in $\finalset{1}$, species 1 has to reach $\vect v_n$ from within $\finalset{1}$ before species 2 can reach it from outside of $\finalset{1}$.

\begin{proof}[Proof of Claim~\ref{entangled_well_defined:time_inequality_claim}]
Suppose the first inequality fails for some $n\ge 0$, i.e.\ $\ttime_1\parens[2]{\subpath{P_n}{\initialset{1}}{\vect u_n}} \le \ttime_2(Q_n)$. Then since \notiesdet\ holds, we must in fact have $\ttime_1\parens[2]{\subpath{P_n}{\initialset{1}}{\vect u_n}} < \ttime_2(Q_n)$, because $\subpath{P_n}{\initialset{1}}{\vect u_n}$ is edge-disjoint from $\finalset{2}$ by Claim~\ref{entangled_well_defined:star_path_claim}, hence edge-disjoint from $Q_n\subseteq \finalset{2}$. It then follows that
\begin{equation}\label{point_not_in_C2_eqn}
\starptime{(\Zd\setminus\finalset{2})}_1(\initialset{1},\vect u_n)
\le \ttime_1\parens[2]{\subpath{P_n}{\initialset{1}}{\vect u_n}} 
< \ttime_2(Q_n) = \rptmetric[2]{\finalset{2}}(\initialset{2}, \vect u_n).
\end{equation}
The first inequality in \eqref{point_not_in_C2_eqn} holds because $\subpath{P_n}{\initialset{2}}{\vect u_n}$ is a $\stargraph{(\Zd\setminus\finalset{2})}$-path from $\initialset{1}$ to $\vect u_n$ by Claim~\ref{entangled_well_defined:star_path_claim}, and the final equality holds because $Q_n$ is a $\finalset{2}$-geodesic from $\initialset{2}$ to $\vect u_n$. But then it follows from \eqref{point_not_in_C2_eqn} and Part~\ref{final_sets:characterization_part} of Lemma~\ref{final_sets_properties_lem} that $\vect u_n\not\in \finalset{2}$, contradicting the choice of $\vect u_n\in \finalset{1}\cap\finalset{2}$. Thus we conclude that $\ttime_1\parens[2]{\subpath{P_n}{\initialset{1}}{\vect u_n}} > \ttime_2(Q_n)$ for all $n\ge 0$.

Similarly, if $\ttime_2\parens[2]{\subpath{Q_n}{\initialset{2}}{\vect v_n}} \le \ttime_1(P_{n+1})$ for some $n\ge 0$, it follows that
\[
\starptime{(\Zd\setminus\finalset{1})}_2(\initialset{2},\vect v_n)
\le \ttime_2\parens[2]{\subpath{Q_n}{\initialset{2}}{\vect v_n}} 
< \ttime_1(P_{n+1}) = \rptmetric[1]{\finalset{1}}(\initialset{1}, \vect v_n),
\]
contradicting the fact that $\vect v_n\in \finalset{1}$, so we conclude that $\ttime_2\parens[2]{\subpath{Q_n}{\initialset{2}}{\vect v_n}} > \ttime_1(P_{n+1})$ for all $n\ge 0$.
\end{proof}

The following relations are trivial, but worth highlighting:

\begin{claim}
\label{entangled_well_defined:additivity_claim}
For $n\ge 0$,
\begin{align*}
\ttime_1(P_n)
	&= \ttime_1\parens[2]{\subpath{P_n}{\initialset{1}}{\vect u_n}}
	+ \ttime_1 \parens[2]{\subpath{P_n}{\vect u_n}{\vect v_{n-1}}}, \text{ and}\\
\ttime_2(Q_n)
	&= \ttime_2\parens[2]{\subpath{Q_n}{\initialset{2}}{\vect v_n}}
	+ \ttime_2 \parens[2]{\subpath{Q_n}{\vect v_n}{\vect u_n}}.
\end{align*}
\end{claim}

Claim~\ref{entangled_well_defined:additivity_claim} follows from the additivity of the measures $\ttime_1$ and $\ttime_2$ and the fact that $P_n$ and $Q_n$ are simple paths (because they are geodesics), so they are the (edge-)disjoint union of the subpaths on the right.

The next claim shows that our inductive construction of the paths $P_n$ and $Q_n$ does not ``break down" after a finite number of steps, and actually yields an infinite sequence of distinct vertices $\vect u_n$ and $\vect v_n$:

\begin{claim}
\label{entangled_well_defined:distinct_claim}
The vertices $\vect u_n$ and $\vect v_n$ ($n\in\N$) are all distinct. That is, for any $m,n\in\N$, $\vect u_m\ne \vect v_n$, and if $m\ne n$, then $\vect u_m\ne \vect u_n$ and $\vect v_m\ne \vect v_n$. Moreover, the sequences $\setof[2]{\rptmetric[1]{\finalset{1}}(\initialset{1},\vect u_n)}_{n\in\N}$ and $\setof[2]{\rptmetric[2]{\finalset{2}}(\initialset{2},\vect v_n)}_{n\in\N}$ are strictly decreasing.
\end{claim}

\begin{proof}[Proof of Claim~\ref{entangled_well_defined:distinct_claim}]
First note that since subpaths of geodesics are also geodesics, Claim~\ref{entangled_well_defined:time_inequality_claim} implies that for all $m,n\in\N$ we have
\begin{align*}
\rptmetric[1]{\finalset{1}}(\initialset{1},\vect u_m)
	&=\ttime_1\parens[2]{\subpath{P_m}{\initialset{1}}{\vect u_m}}
	> \ttime_2(Q_m)
	=\rptmetric[2]{\finalset{2}}(\initialset{2},\vect u_m),\text{ and }\\
\rptmetric[2]{\finalset{2}}(\initialset{2},\vect v_n)
	&=\ttime_2\parens[2]{\subpath{Q_n}{\initialset{2}}{\vect v_n}}
	> \ttime_1(P_{n+1})
	=\rptmetric[_1]{\finalset{1}}(\initialset{1},\vect v_n).
\end{align*}
If we had $\vect u_m = \vect v_n$ for some $m,n\in\N$, this would imply that
\[
\rptmetric[1]{\finalset{1}}(\initialset{1},\vect u_m)
> \rptmetric[2]{\finalset{2}}(\initialset{2},\vect u_m)
=\rptmetric[2]{\finalset{2}}(\initialset{2},\vect v_n)
>\rptmetric[1]{\finalset{1}}(\initialset{1},\vect v_n)
= \rptmetric[1]{\finalset{1}}(\initialset{1},\vect u_m),
\]
which is a contradiction, so we conclude that $\vect u_m \ne \vect v_n$ for all $m,n\in\N$. Next, by combining Claim~\ref{entangled_well_defined:time_inequality_claim} and Claim~\ref{entangled_well_defined:additivity_claim}, we get the following chain of inequalities for each $n\in\N$:
\begin{align*}
\ttime_1 \parens[2]{\subpath{P_n}{\initialset{1}}{\vect u_n}}
>\ttime_2(Q_n)
\ge \ttime_2 \parens[2]{\subpath{Q_n}{\initialset{2}}{\vect v_n}}
>\ttime_1(P_{n+1})
&\ge \ttime_1 \parens[2]{\subpath{P_{n+1}}{\initialset{1}}{\vect u_{n+1}}}\\
&>\ttime_2(Q_{n+1})
\ge \ttime_2 \parens[2]{\subpath{Q_{n+1}}{\initialset{2}}{\vect v_{n+1}}}.
\end{align*}
In particular, for all $n\in\N$ we have
\begin{align*}
\rptmetric[1]{\finalset{1}}(\initialset{1},\vect u_n)
	&=\ttime_1 \parens[2]{\subpath{P_n}{\initialset{1}}{\vect u_n}}
	> \ttime_1 \parens[2]{\subpath{P_{n+1}}{\initialset{1}}{\vect u_{n+1}}}
	=\rptmetric[1]{\finalset{1}}(\initialset{1},\vect u_{n+1}),
	\text{ and}\\
\rptmetric[2]{\finalset{2}}(\initialset{2},\vect v_n)
	&=\ttime_2 \parens[2]{\subpath{Q_n}{\initialset{2}}{\vect v_n}}
	> \ttime_2 \parens[2]{\subpath{Q_{n+1}}{\initialset{2}}{\vect v_{n+1}}}
	=\rptmetric[2]{\finalset{2}}(\initialset{2},\vect v_{n+1}).
\end{align*}
Thus, the sequences $\setof[2]{\rptmetric[1]{\finalset{1}}(\initialset{1},\vect u_n)}_{n\in\N}$ and $\setof[2]{\rptmetric[2]{\finalset{2}}(\initialset{2},\vect v_n)}_{n\in\N}$ are strictly decreasing, so the maps $n\mapsto \vect u_n$ and $n\mapsto \vect v_n$ must be injective.
\end{proof}

Finally, we are ready to construct the path $\gamma$ and derive a contradiction. For $j\in\N$, let $\gamma_{2j} = \subpath{Q_j}{\vect u_j}{\vect v_j} \definedas -\subpath{Q_j}{\vect v_j}{\vect u_j}$, and let $\gamma_{2j+1} = \subpath{P_{j+1}}{\vect v_j}{\vect u_{j+1}} \definedas -\subpath{P_{j+1}}{\vect u_{j+1}}{\vect v_j}$, where the negative signs indicate that the subpaths are traversed in the reverse direction from the paths $Q_j$ and $P_{j+1}$. Then the concatenation $\gamma \definedas \gamma_0 \gamma_1\gamma_2\gamma_3\dotso$ is a path starting at $\vect u_0$ and traversing the $\vect u_n$'s and $\vect v_n$'s in order:
\[
\gamma: \vect u_0 \xrightarrow{\gamma_0} \vect v_0 \xrightarrow{\gamma_1} \vect u_1
\xrightarrow{\gamma_2} \vect v_1\xrightarrow{\gamma_3} \vect u_2 
\xrightarrow{\gamma_4} \vect v_2\xrightarrow{\gamma_5} \dotsb.
\]
By Claim~\ref{entangled_well_defined:distinct_claim}, the vertices in the above sequence are all distinct, so the path $\gamma$ is infinite. Note that although the vertices $\vect u_n$ and $\vect v_n$ are all distinct, we have no guarantee that the subpaths $\gamma_n$ don't intersect at interior vertices, so the path $\gamma$ may not be simple. This is not a problem, however, since we only need an upper bound on the traversal time of $\gamma$, which we now derive.
By Claim~\ref{entangled_well_defined:time_inequality_claim} and Claim~\ref{entangled_well_defined:additivity_claim}, for $j\in\N$ we have
\begin{align*}
\ttime_2(\gamma_{2j}) 
	+ \ttime_2 \parens[2]{\subpath{Q_j}{\initialset{2}}{\vect v_j}} 
	&= \ttime_2(Q_j) 
	< \ttime_1 \parens[2]{\subpath{P_j}{\initialset{1}}{\vect  u_j}},
	\text{ and}\\
\ttime_1(\gamma_{2j+1}) 
	+ \ttime_1 \parens[2]{\subpath{P_{j+1}}{\initialset{1}}{\vect u_{j+1}} }
 	&= \ttime_1(P_{j+1})
	<\ttime_2 \parens[2]{\subpath{Q_j}{\initialset{2}}{\vect v_j}},
\end{align*}
and hence
\begin{align*}
\ttime_2(\gamma_{2j}) 
	&< \ttime_1 \parens[2]{\subpath{P_j}{\initialset{1}}{\vect  u_j}}
	- \ttime_2 \parens[2]{\subpath{Q_j}{\initialset{2}}{\vect v_j}},
	\text{ and}\\
\ttime_1(\gamma_{2j+1}) 
	&<\ttime_2 \parens[2]{\subpath{Q_j}{\initialset{2}}{\vect v_j}}
	- \ttime_1 \parens[2]{\subpath{P_{j+1}}{\initialset{1}}{\vect u_{j+1}} }.
\end{align*}
Thus, we get a telescoping sum bounding $\mintime(\gamma)$:
\begin{align*}
\mintime(\gamma)
&\le \ttime_2(\gamma_0) + \ttime_1(\gamma_1) 
	+ \ttime_2(\gamma_2) + \ttime_1(\gamma_3) + \dotsb\\
&< \sum_{j=0}^\infty \brackets[3]{%
	\ttime_1 \parens[2]{\subpath{P_j}{\initialset{1}}{\vect  u_j}}
	- \ttime_2 \parens[2]{\subpath{Q_j}{\initialset{2}}{\vect v_j}}
	+ \ttime_2 \parens[2]{\subpath{Q_j}{\initialset{2}}{\vect v_j}}
	- \ttime_1 \parens[2]{\subpath{P_{j+1}}{\initialset{1}}{\vect u_{j+1}} }
	}\\
&= \sum_{j=0}^\infty \brackets[1]{%
	\ttime_1 \parens[2]{\subpath{P_j}{\initialset{1}}{\vect  u_j}}
	- \ttime_1 \parens[2]{\subpath{P_{j+1}}{\initialset{1}}{\vect u_{j+1}} }
	}\\
&= \ttime_1 \parens[2]{\subpath{P_0}{\initialset{1}}{\vect  u_0}}
	-\lim_{j\to\infty} \rptmetric[1]{\finalset{1}}(\initialset{1},\vect u_j).
\end{align*}
Therefore, $\gamma$ is an infinite path with $\mintime(\gamma)<  \rptmetric[1]{\finalset{1}}(\initialset{1},\vect u_0) <\infty$. Since this contradicts \mintimebddgrowth, we conclude that $\finalset{1}\cap\finalset{2}=\emptyset$.
\end{proof}

\begin{thmremark}
If $A_1$ and $A_2$ are not both infinite, the assumption \mintimebddgrowth\ in Proposition~\ref{entangled_well_defined_prop} can be weakened. For example, if $A_1$ is finite, then instead of \mintimebddgrowth, we can assume merely that $\ttime_1$ satisfies \boundedgrowth\ and $\ttime_2$ satisfies \stargeodesicsexist. Note that this assumption is weaker than the one needed for the algorithmic construction, and that it is not necessary for $A_2$ to be finite. The proof proceeds exactly as above up through Claim~\ref{entangled_well_defined:distinct_claim}. Then, rather than constructing the infinite path $\gamma$, we simply observe that if $\initialset{1}$ is finite, then the fact that $\setof[2]{\rptmetric[1]{\finalset{1}}(\initialset{1},\vect u_n)}_{n\in\N}$ is strictly decreasing implies (via a compactness argument) that there is an infinite $\rptmetric[1]{\finalset{1}}$-geodesic with finite traversal time, contradicting \boundedgrowth\ for $\ttime_1$.
\end{thmremark}

\subsection{First-Passage Competition as an Interacting Particle System}
\label{fpc_basics:ips_def_sec}

%
%

In order to more easily make formal probabilistic arguments, we can view the first-passage competition process as a (possibly non-Markovian) interacting particle system (cf.\ Liggett \cite{Liggett:1985aa}) with state space $\statemeasurespace$, where $\statespace = \setof{0,1,2}^{\Zd}$, and $\statesigmafield$ is the product $\sigma$-field (generated by cylinder sets). The state of a vertex $\vect v\in \Zd$ at time $t\ge 0$ is 0 if $\vect v$ has not yet been infected, and is 1 or 2 if $\vect v$ has been infected by species 1 or 2, respectively.
If the process is well-defined and full, the state of every vertex eventually changes from 0 to either 1 or 2, and then remains in that state forever.

To describe the interacting particle system perspective
more formally, first note that the state space $\statespace$ is in one-to-one correspondence with pairs of disjoint subsets of $\Zd$:
\begin{quote}
For $X\in\statespace$ and $A_1,A_2\subseteq\Zd$ with $A_1\cap A_2=\emptyset$, we write
\begin{equation}
\label{state_config_equiv_eqn}
\begin{split}
X\equiv (A_1,A_2) &\iff X = 1\cdot \ind{A_1} + 2\cdot \ind{A_2}\\
&\iff A_1 = X^{-1}\setof{1} \text{ and } A_2 = X^{-1}\setof{2}.
\end{split}
\end{equation}
\end{quote}
Now let $\tmeasurepair$ be a random two-type traversal measure in some probability space $\probabilityspace$, and suppose that $\tmeasurepair$ satisfies \notiesdet\ and \mintimebddgrowth\ $\Pr$-a.s., so that the entangled process $\reachedfcn{\tmeasurepair}{X} = \reachedfcn{\tmeasurepair}{\initconfig}$ is almost surely well-defined for all starting configurations $X\equiv \initconfig$ by Proposition~\ref{entangled_well_defined_prop} (and also full by Proposition~\ref{entangled_full_prop}). Then, for any $X\in \statespace$, we let $\process{X} =\processoutcome[2]{X}{\tmeasurepair}$ denote the $\statespace$-valued entangled process on $[0,\infty]$ with initial configuration $X$ and run with the random traversal measure $\tmeasurepair$. That is, using the correspondence in \eqref{state_config_equiv_eqn}, we have by definition
\begin{equation}\label{ips_def_eqn}
\processt{X}{t} \equiv \reachedset{\tmeasurepair}{X}{t}
\quad\text{for all } t\in [0,\infty],
\end{equation}
so that for $\vect v\in\Zd$ and $ i\in \setof{1,2}$, $\processtv{X}{t}{\vect v} = i$ iff $\vect v\in \reachedset{i}{X}{t}$. Note that \eqref{ips_def_eqn} is a \emph{pointwise} definition, i.e.\ $\processtoutcome[2]{X}{t}{\tmeasurepair_{\outcome}}\equiv \reachedset{\tmeasurepair_{\outcome}}{X}{t}$ for all $\outcome$ in the $\Pr$-a.s.\ subset of $\samplespace$ on which $\reachedfcn{\tmeasurepair}{X}$ is well-defined. We leave it to the reader to verify that for all $X\in \statespace$ and $t\ge 0$, the map $\processtoutcome[2]{X}{t}{\tmeasurepair}\colon \samplespace \to\statespace$ is $\sigmafield/\statesigmafield$-measurable whenever the traversal measure $\tmeasurepair \colon \samplespace \to \samplespace[1]$ is Borel measurable.



%
%

The two-type competition process $\process{X}$ was originally defined in \cite{\HPa} and \cite{\HPb} as a Markov process on $\statespace$. This is equivalent (by the memoryless property of the exponential distribution) to specifying in \eqref{ips_def_eqn} that $\tmeasurepair$ is \iid\ with exponentially distributed model traversal times, say $\mttimes{i}_0 \sim \operatorname{Exponential}(\lambda_i)$ for some $\lambda_1,\lambda_2>0$.
In the Markov case, the process can be completely described by its infinitesimal dynamics (cf.\ \cite{Liggett:2010aa}): A vertex in state 1 or 2 never changes state, while a vertex in state 0 changes to state $i\in \setof{1,2}$ at infinitesimal rate $\lambda_i$ times the number of neighbors it has in state $i$.

For traversal measures $\tmeasurepair$ that are not \iid\ exponential, the process $\process{X}$ as defined in \eqref{ips_def_eqn} is not Markov. However, $\process{X}$ does always satisfy a sort of ``spatio-temporal Markov property" rather than a purely temporal one, because if the ``state" of a vertex is enlarged to carry extra information taking into account the time at which it was infected, then the future evolution of the process $\process{X}$ only depends on the states of vertices on the \emph{boundary} of the conquered region. We will use this idea to prove Lemma~\ref{markovesque_lem} in Section~\ref{fpc_basics:shift_ops_sec}, which will substitute for the Markov property in the proof of Theorem~\ref{irrelevance_thm} below.

%
%
%
%
%
%
%
%
%

\subsection{Additional Notation for the Two-Type Process}
\label{fpc_basics:additional_notation_sec}

In this subsection we introduce some additional definitions and notation for the two-type process that will be needed in Sections~\ref{fpc_basics:monotonicity_sec}, \ref{fpc_basics:shift_ops_sec}, and \ref{fpc_basics:unbounded_survival_sec} below, as well as in Chapters~\ref{random_fpc_chap} and \ref{coex_finite_chap}.

\subsubsection{Occupied Sets and Entangled Passage Times}

Let $\initconfig$ be an initial configuration and let $\ttimepair$ be a two-type traversal measure on $\Zd$. For $t\ge 0$, we define the \textdef{total occupied region} at time $t$ by
\begin{equation}
\label{total_occupied_def_eqn}
\bothreachedset{\initconfig}{t} \definedas \reachedset{1}{\initconfig}{t}\cup \reachedset{2}{\initconfig}{t},
\end{equation}
i.e.\ the set of sites conquered by either species by time $t$. If $\vect v\in\Zd$ and $\vect v\notin \bothreachedset{\initconfig}{t}$, we say that $\vect v$ is \textdef{unoccupied} at time $t$. For $\vect v\in\Zd$ and $i\in\setof{1,2}$, we define \textdef{species $i$'s $\initconfig$-entangled passage time} to $\vect v$ by
\begin{equation}
\label{entangled_time_def_eqn}
\entangledtime{i}{\initconfig}{\vect v} \definedas
\inf \setbuilder[1]{t}{\vect v\in \reachedset{i}{\initconfig}{t}}
=\rptime{i}{\finalset{i}}{A_i}{\vect v},
\end{equation}
i.e.\ the time $\vect v$ is infected by species $i$ in the entangled process started from $\initconfig$. We also define the \textdef{$\initconfig$-infection time} of $\vect v$ (by either species) as
\begin{equation}
\label{infection_time_def_eqn}
\infectiontime{\initconfig}{\vect v} \definedas
\inf \setbuilder[1]{t}{\vect v \in \bothreachedset{\initconfig}{t}}
=\entangledtime{1}{\initconfig}{\vect v} \meet \entangledtime{2}{\initconfig}{\vect v}.
\end{equation}
More generally, we can similarly define the \textdef{entangled hitting time} of any subset $B\subseteq\Zd$ by one or either species in the two-type process.
Note that if $\reachedfcn{\tmeasurepair}{\initconfig}$ is well-defined and full, then exactly one of $\entangledtime{1}{\initconfig}{\vect v}$ and $\entangledtime{2}{\initconfig}{\vect v}$ is finite. With the above definitions, observe that
\begin{equation*}
\reachedset{i}{\initconfig}{t} =
\setbuilder[3]{\vect v\in\Zd}{\entangledtime{i}{\initconfig}{\vect v} \le t},
\text{ and\ \ }
\bothreachedset{\initconfig}{t} =
\setbuilder[3]{\vect v\in\Zd}{\infectiontime{\initconfig}{\vect v} \le t}.
\end{equation*}
If $\initconfig\equiv X\in \statespace$, the notation for the total occupied region, entangled passage time, and infection time becomes $\bothreachedset{X}{t}$, $\entangledtime{i}{X}{\vect v}$, and $\infectiontime{X}{\vect v}$, respectively.

\subsubsection{Notation for Elements of the State Space $\statespace$}

Here we introduce some convenient notation for elements of the state space $\statespace = \statespace[1]$. The following definitions are intended to be intuitive shorthand for converting between states in $\statespace$ and various related subsets of $\Zd$, following the general convention of identifying a state with its set of occupied sites. If $X\in\statespace$, we define the \textdef{occupied regions for $X$} by
\begin{equation}
\label{state_occupied_sets_def_eqn}
\stateset{i}{X} \definedas X^{-1}\setof{i}
\text{ for $i\in \setof{1,2}$,\quad and}\quad
\nonzeroset{X} \definedas X^{-1} \setof{1,2}.
\end{equation}
That is, if $X\equiv \initconfig$, then $\stateset{i}{X} = \initialset{i}$ is the region occupied by species $i$ in $X$, and $\nonzeroset{X}=\initialset{1}\cup \initialset{2}$ is the region occupied by either species. Note that it follows from \eqref{state_occupied_sets_def_eqn} and the definition \eqref{ips_def_eqn} of $\process{X}$ that 
\[
\stateset[2]{i}{\processt{X}{t}} = \reachedset{i}{X}{t}
\quad\text{and}\quad
\nonzeroset[2]{\processt{X}{t}} = \bothreachedset{X}{t}
\quad\text{for all $t\in [0,\infty]$.}
\]
Just as we did for subsets of $\Zd$, we can identify a state $X\in\statespace$ with an induced subgraph of $\Zd$, namely the subgraph induced by its total occupied region $\nonzeroset{X}$. In this way, we extend all of our \textdef{graph theoretic notation} to any $X\in \statespace$, e.g.
\[
\text{$\nbrhood{X}$, $\nbrset{X}$, $\bd X$, $\edges{X}$, $\compedges{X}$, $\bdedges{X}$, $\staredges{X}$,}
\]
by replacing each instance of $X$ in this notation with the subgraph of $\Zd$ induced by $\nonzeroset{X}$. Furthermore, for $X,X'\in\statespace$, we define the \textdef{union of $X$ and $X'$} to be the set of vertices
\begin{equation}
\label{state_union_def_eqn}
X\cup X'\definedas \nonzeroset{X} \cup \nonzeroset{X'} \in \powerset{\Zd}.
\end{equation}
That is, $X\cup X' \subseteq\Zd$ is the set of all occupied vertices in either $X$ or $X'$. Following our convention for initial configurations $\initconfig$, we call a state $X\in\statespace$ \textdef{finite} if $\card{\nonzeroset{X}}<\infty$ and \textdef{infinite} if $\card{\nonzeroset{X}}=\infty$. We also call $X$ \textdef{final} if $\nonzeroset{X}=\Zd$; thus any full process ends in a final state by definition.
More generally, we define the following three \textdef{subspaces of $\statespace$}:
\begin{equation}
\label{state_subspace_def_eqn}
\statesubspace{0,1} \definedas \statesubspace[1]{0,1},\qquad
\statesubspace{0,2} \definedas \statesubspace[1]{0,2},
\qquad\text{and}\qquad
\statesubspace{1,2} \definedas \statesubspace[1]{1,2}.
\end{equation}
That is, $\statesubspace{0,i}$ is the set of states corresponding to a one-type process of type $i$, and $\statesubspace{1,2}$ is the set of final states. Note that the three subspaces in \eqref{state_subspace_def_eqn} are precisely the collections of trivial initial configurations to which Lemma~\ref{trivial_init_configs_lem} applies.

\subsubsection{Notation for the Canonical Sample Space $\samplespace$ and Subspaces $\rsamplespace{\Lambda}$}

Let $\samplespace = \samplespace[1]$ be the canonical sample space for the two-type process, as defined in Section~\ref{fpc_basics:prob_space_sec}. We write $\outcome = (\outcome_1,\outcome_2)$ for an element of $\samplespace$, where $\outcome_1,\outcome_2\in (\Rplus)^{\edges{\Zd}}$ are the component functions of $\outcome$. Recall that in the canonical case, if $\ttimepair = \pair{\ttime_1}{\ttime_2}$ is the collection of traversal time pairs corresponding to the outcome $\outcome\in \samplespace$, then $\ttimepair = \outcome$ and $\ttime_i = \omega_i$. If we want to view $\tmeasurepair$ as a measure rather than a function (as in definition \eqref{traversal_measure_def}, Section~\ref{fpc_basics:traversal_meas_sec}), we will write $\tmeasurepair_\outcome$ instead of $\outcome$. If $\Lambda \subseteq \edges{\Zd}$, we write $\restrict{\outcome}{\Lambda}$ for the function $\outcome$ restricted to the edge set $\Lambda$, and $\restrict{\outcome_i}{\Lambda}$ ($i\in\setof{1,2}$) for its component functions. We define $\rsamplespace{\Lambda} \definedas (\R_+\times\R_+)^\Lambda$, so that $\restrict{\outcome}{\Lambda}\in \rsamplespace{\Lambda}$ for any $\outcome\in \samplespace$.

\subsubsection{The Canonical Process Space $\processspace$}

The first-passage competition process $\process{X}$ is defined in Section~\ref{fpc_basics:ips_def_sec} above as an $\statespace$-valued stochastic process on $[0,\infty]$, where the state space $\statespace = \setof{0,1,2}^{\Zd}$ is endowed with the product $\sigma$-field $\statesigmafield$. For each $X\in\statespace$ and $t\ge 0$, the state of the process started from $X$ at time $t$ is $\processt{X}{t}\in\statespace$, which is a random element of $\statemeasurespace$ in some unspecified probability space $\probabilityspace$. Equivalently, for a given initial state $X$, the process $\process{X}$ is a random element of $\pair[2]{\processspace}{\processsigmafield}$, which we call the \textdef{canonical process space}.  We write $\processlaw{\tmeasurepair}{X}$ for the \textdef{law of the two-type process} started from $X$ and using the random traversal measure $\tmeasurepair$; that is, $\processlaw{\tmeasurepair}{X} \definedas \Pr \circ \parens[2]{\processoutcome{X}{\tmeasurepair}}^{-1}$ is the pushforward measure on the canonical process space.

\subsubsection{The Two-Type Process as a Map Between Canonical Spaces}

In Sections~\ref{fpc_basics:monotonicity_sec}, \ref{fpc_basics:shift_ops_sec}, and \ref{fpc_basics:unbounded_survival_sec} below, we will assume that $\probabilityspace$ is the canonical probability space defined in Section~\ref{fpc_basics:prob_space_sec}, i.e.\ $\samplespace = \samplespace[1]$, endowed with its Borel $\sigma$-field $\sigmafield$ and a product measure $\Pr$ with marginals satisfying $\finitespeed{d}$, $\minimalmoment{d}$, and $\notiesprob$. In order to take the starting configuration $X$ into account, it will be convenient to treat $\process{}$ as a function $\domainspace \to \processspace$, which maps the pair $\pair{\outcome}{X}$ to the process $\processoutcome[2]{X}{\outcome}$ started from $X$ on the outcome $\outcome$. This is not strictly correct, however, since for a given starting state $X$, the process $\process{X}$ may only be defined on a $\Pr$-a.s.\ subset of $\samplespace$, and so $\process{}$ is not well-defined for all pairs $\pair{\outcome}{X}\in \domainspace$. Thus, we refer to any $\domain\subseteq \samplespace\times\statespace$ on which $\process{}$ is defined as a \textdef{domain of definition} for the two-type process, and we write $\process{}\colon \domain \to \processspace$.

We explicitly identify one particular domain of definition $\domain = \startdomain$ as follows. First define the set of \textdef{good traversal measures} by
\begin{equation}
\label{goodsamplespace_def_eqn}
\goodsamplespace \definedas
\setbuilder[1]{\outcome \in\samplespace[1]}
	{\outcome \text{ satisfies \notiesdet\ and \mintimebddgrowth}}.
\end{equation}
In the canonical probability space we have $\goodsamplespace\in\sigmafield$, and $\Pr(\goodsamplespace) = 1$ by Lemma~\ref{two_type_a_s_properties_lem}. Now define the \textdef{initial domain} of the two-type process to be
\begin{equation}
\label{startdomain_def_eqn}
\startdomain \definedas \samplespace \times
	\parens[1]{\statesubspace{0,1} \cup \statesubspace{0,2} \cup \statesubspace{1,2}}
\cup \parens[2]{\goodsamplespace \times \statespace},
\end{equation}
where $\statesubspace{0,1}$, $\statesubspace{0,2}$, and $\statesubspace{1,2}$ are the subspaces of $\statespace$ defined in \eqref{state_subspace_def_eqn}. Then by Lemma~\ref{trivial_init_configs_lem} and Propositions~\ref{entangled_full_prop} and \ref{entangled_well_defined_prop}, the two-type process $\reachedfcn{\outcome}{X}$ is well-defined and full for all $\pair{\outcome}{X}\in \startdomain$.
In particular, the fact that $\reachedfcn{\outcome}{X}$ is well-defined on $\startdomain$ means that the occupied regions $\reachedset{1}{X}{t}_\outcome$ and $\reachedset{2}{X}{t}_\outcome$ are disjoint for each $t\in [0,\infty]$, so the map $\processt{}{t}$ is a legitimate function from $\startdomain$ into the subspace $\statespace = 3^{\Zd}$ of the larger \textit{a priori} codomain $\powerset{\Zd}\times \powerset{\Zd} \cong 4^{\Zd}$, given by $\pair{\outcome}{X} \mapsto \pair[2]{\reachedset{1}{X}{t}}{\,\reachedset{2}{X}{t}}_\outcome \equiv \processtoutcome{X}{t}{\outcome}$. 
Thus, $\startdomain$ is a domain of definition for the two-type process, meaning that $\processspace$ is a
valid codomain for the function $\process{}\colon \startdomain \to \processspace$ defined on $\startdomain$. We will see below (in Corollary~\ref{domain_extension_cor} and Lemma~\ref{markovesque_lem}) that we can in fact extend $\domain_0$ to larger domains, which will be necessary to state some of the properties of $\process{}$ that will be used later.

\section{Basic Properties of First-Passage Percolation \& First-Passage Competition}
\label{fpc_basics:basic_properties_sec}

\subsection{Locality and Equivalent Processes}
\label{fpc_basics:locality_sec}

In this section we prove various properties about the restricted one-type process and the two-type process which show, essentially, that if an $R$-restricted process is actually contained in some set $S\subseteq R$ at a given time $t$, then knowledge about the process growing in $S$ is equivalent to knowledge about the process growing in $R$, up to time $t$.
We refer to this general property as ``locality," and it implies (Lemma~\ref{one_type_two_type_lem}) that the growth of species $i$ in the entangled process $\reachedfcn{\tmeasurepair}{\initconfig}$ up to time $t$ is equivalent to the growth of species $i$'s disentangled process restricted to any set containing $\reachedset{i}{\initconfig}{t}$ and avoiding $\reachedset{3-i}{\initconfig}{t}$. This reduction of a two-type process to a one-type process is similar to \Haggstrom and Pemantle's ``separator lemma" \cite[p.~5]{Haggstrom:2000aa}, and will be a key step in analyzing the entangled process in the next section.

\begin{lem}[Locality of restricted passage times] 
\label{restricted_times_equiv_lem}
Let $\tmeasure$ be a traversal measure on $\edges{\Zd}$ with $\ptimefamily = \ptimefamily[\tmeasure]$. If $A\subseteq S\subseteq R\subseteq \Zd$, then
\begin{enumerate}
\item \label{restricted_times_equiv:set_times_part}
$\rptime{}{R}{A}{\,R\setminus S} = \starptimenew{}{S}{A}{\,R\setminus S}$.
\item \label{restricted_times_equiv:vertex_times_part}
If $\vect v\in R$ and $\rptime{}{R}{A}{\vect v} < \starptimenew{}{S}{A}{\,R\setminus S}$, then $\rptime{}{R}{A}{\vect v} = \rptime{}{S}{A}{\vect v}$.
\end{enumerate}
\end{lem}


\begin{proof}
\newcommand{\starpath}[1]{#1^*}


\proofpart{\ref{restricted_times_equiv:set_times_part}} Since $A\subseteq R$, $R\setminus S\subseteq R$, and $S\subseteq R$, we trivially have 
\[
\rptime{}{R}{A}{R\setminus S} = \starptimenew{}{R}{A}{R\setminus S}
\le \starptimenew{}{S}{A}{R\setminus S},
\]
so we need to prove the reverse inequality. Given any $R$-path $\gamma$ from $A$ to $R\setminus S$, let $\starpath{\gamma}$ be the subpath of $\gamma$ ending at the \emph{first} vertex of $\gamma$ lying outside of $S$.
Then every vertex of $\starpath{\gamma}$ is contained in $S$ except for the final vertex, so $\starpath{\gamma}$ is an $\stargraph{S}$-path from $A$ to $R\setminus S$, and $\tmeasure(\starpath{\gamma})\le \tmeasure (\gamma)$ since $\starpath{\gamma}$ is a subpath of $\gamma$. Therefore,
\begin{align*}
\rptime{}{R}{A}{R\setminus S}
&= \inf \setbuilder[1]{\tmeasure (\gamma)}
	{\text{$\gamma$ is an $R$-path from $A$ to $R\setminus S$}}\\
&\ge \inf \setbuilder[1]{\tmeasure (\starpath{\gamma})}
	{\text{$\gamma$ is an $R$-path from $A$ to $R\setminus S$}}\\
&\ge \inf \setbuilder[1]{\tmeasure (\lambda)}
	{\text{$\lambda$ is an $\stargraph{S}$-path from $A$ to $R\setminus S$}}\\
&= \starptimenew{}{S}{A}{R\setminus S}.
\end{align*}

\proofpart{\ref{restricted_times_equiv:vertex_times_part}} 
Let $\vect v\in R$, and suppose $\rptime{}{R}{A}{\vect v}< \starptimenew{}{S}{A}{R\setminus S}$.
Since $S\subseteq R$, the inequality $\rptime{}{R}{A}{\vect v}\le \rptime{}{S}{A}{\vect v}$ is trivial, so we need to prove the reverse inequality. Since $\rptime{}{R}{A}{\vect v}< \starptimenew{}{S}{A}{R\setminus S}$, there is some $R$-path $\gamma$ from $A$ to $\vect v$ such that $\tmeasure(\gamma) < \starptimenew{}{S}{A}{R\setminus S} = \rptime{}{R}{A}{R\setminus S}$, where the final equality is from Part~\ref{restricted_times_equiv:set_times_part}. We claim that any such path satisfies $\gamma \subseteq S$. Indeed, suppose $\gamma$ contained some vertex $\vect u\in R\setminus S$. Then the subpath $\subpath{\gamma}{A}{\vect u}$ would be an $R$-path from $A$ to $R\setminus S$, and we would then have
\[
\rptime{}{R}{A}{R\setminus S} \le
\tmeasure \parens[2]{\subpath{\gamma}{A}{\vect u}}
\le \tmeasure (\gamma) < \rptime{}{R}{A}{R\setminus S},
\]
which is a contradiction. Therefore, since $\rptime{}{R}{A}{\vect v}< \rptime{}{R}{A}{R\setminus S}$, we have
\begin{align*}
\rptime{}{R}{A}{\vect v}
&= \inf \setbuilder[2]{\tmeasure(\gamma)}
	{\text{$\gamma$ is an $R$-path from $A$ to $\vect v$,
	and $\tmeasure(\gamma) <  \rptime{}{R}{A}{R\setminus S}$}}\\
&\ge \inf \setbuilder[2]{\tmeasure(\gamma)}
	{\text{$\gamma$ is an $R$-path from $A$ to $\vect v$,
	and $\gamma\subseteq S$}}\\
&= \rptime{}{S}{A}{\vect v}.\qedhere
\end{align*}
\end{proof}

The following lemma gives several equivalent conditions for checking that the left-continuous $R$-restricted process is contained in some set $S\subseteq R$ at time $t$.

\begin{lem}[Restricting the process to a smaller set] 
\label{restricted_proc_equiv_lem}
Let $\tmeasure$ be a traversal measure on $\edges{\Zd}$ with $\ptimefamily = \ptimefamily[\tmeasure]$, and let $A\subseteq S\subseteq R\subseteq \Zd$. Then for any $t\ge 0$, the following are equivalent.
\begin{enumerate}
\item\label{restricted_proc_equiv:hitting_time_part}
$\starptimenew{}{S}{A}{\,R\setminus S}\ge t$.
\item \label{restricted_proc_equiv:vertex_times_part}
If  $\vect v\in R$ and $\rptime{}{R}{A}{\vect v} < t$, then $\rptime{}{R}{A}{\vect v} = \rptime{}{S}{A}{\vect v}$.
\item \label{restricted_proc_equiv:process_times_part}
$\rreachedset{\tmeasure}{A}{R}{t'} = \rreachedset{\tmeasure}{A}{S}{t'}$ for all $t'< t$.
\item \label{restricted_proc_equiv:final_time_part}
$\rreachedset{\tmeasure}{A}{R}{t-} = \rreachedset{\tmeasure}{A}{S}{t-}$.
\item \label{restricted_proc_equiv:containment_part}
$\rreachedset{\tmeasure}{A}{R}{t-}\subseteq S$.
\item \label{restricted_proc_equiv:star_containment_part}
$R\cap \rreachedset{\tmeasure}{A}{\stargraph{S}}{t-}\subseteq S$.
\end{enumerate}
\end{lem}


\begin{proof}

($\ref{restricted_proc_equiv:hitting_time_part} \implies \ref{restricted_proc_equiv:vertex_times_part}$) Suppose $\starptimenew{}{S}{A}{\,R\setminus S}\ge t$. If  $\vect v\in R$ and $\rptime{}{R}{A}{\vect v} < t$, then $\rptime{}{R}{A}{\vect v}< \starptimenew{}{S}{A}{\,R\setminus S}$, so $\rptime{}{R}{A}{\vect v} = \rptime{}{S}{A}{\vect v}$ by Part~\ref{restricted_times_equiv:vertex_times_part} of Lemma~\ref{restricted_times_equiv_lem}.

($\ref{restricted_proc_equiv:vertex_times_part} \implies \ref{restricted_proc_equiv:process_times_part}$) Suppose $\rptime{}{R}{A}{\vect v} = \rptime{}{S}{A}{\vect v}$ for all $\vect v\in R$ with $\rptime{}{R}{A}{\vect v} < t$. Then for any $t'<t$ and $\vect v\in R$, we have
\[
\rptime{}{R}{A}{\vect v}\le t' \ifandonlyif \rptime{}{S}{A}{\vect v}\le t'
\ifandonlyif \rptime{}{S}{A}{\vect v}\le t' \text{ and } \vect v\in S
\]
since $\rptime{}{S}{A}{\vect v}=\infty$ for $\vect v\not\in S$. Therefore, for any $t'<t$,
\begin{align*}
\rreachedset{\tmeasure}{A}{R}{t'}
= \setbuilder[1]{\vect v\in R}{\rptime{}{R}{A}{\vect v}\le t'}
= \setbuilder[1]{\vect v\in S}{\rptime{}{S}{A}{\vect v}\le t'}
= \rreachedset{\tmeasure}{A}{S}{t'}.
\end{align*}

($\ref{restricted_proc_equiv:process_times_part}\implies \ref{restricted_proc_equiv:final_time_part}$) If Statement~\ref{restricted_proc_equiv:process_times_part} holds, then by definition we have
\[
\rreachedset{\tmeasure}{A}{R}{t-} = \bigcup_{t'<t} \rreachedset{\tmeasure}{A}{R}{t'}
= \bigcup_{t'<t} \rreachedset{\tmeasure}{A}{S}{t'} = \rreachedset{\tmeasure}{A}{S}{t-}.
\]

($\ref{restricted_proc_equiv:final_time_part}\implies \ref{restricted_proc_equiv:containment_part}$) This is trivial since $\rreachedset{\tmeasure}{A}{S}{t-}\subseteq S$ by definition.


($\ref{restricted_proc_equiv:containment_part}\implies \ref{restricted_proc_equiv:star_containment_part}$) Suppose $\rreachedset{\tmeasure}{A}{R}{t-} \subseteq S$. Since $A\subseteq R$ and $S\subseteq R$, we have
\[
\rptime{}{R}{A}{\vect v} = \starptimenew{}{R}{A}{\vect v}
\le \starptimenew{}{S}{A}{\vect v}
\quad\text{for all } \vect v\in R.
\]
Therefore,
\[
R\cap \rreachedset{\tmeasure}{A}{\stargraph{S}}{t-}
= \setbuilder[2]{\vect v\in R}{ \starptimenew{}{S}{A}{\vect v}< t}
\subseteq  \setbuilder[2]{\vect v\in R}{ \rptime{}{R}{A}{\vect v}< t}
= \rreachedset{\tmeasure}{A}{R}{t-} \subseteq S.
\]

($\ref{restricted_proc_equiv:star_containment_part}\implies \ref{restricted_proc_equiv:hitting_time_part}$) Suppose $R\cap \rreachedset{\tmeasure}{A}{\stargraph{S}}{t-}\subseteq S$, i.e.\ if $\vect v\in R$ and $\starptimenew{}{S}{A}{\vect v} < t$, then $\vect v\in S$. Then for any $\vect v\in R\setminus S$ we have $\starptimenew{}{S}{A}{\vect v}\ge t$, so $\starptimenew{}{S}{A}{R\setminus S}\ge t$.
\end{proof}


The corresponding statement for the right-continuous process is slightly more complicated.

\begin{lem}[Re-restricting the right-continuous process]
\label{strong_restricted_proc_equiv_lem}
Let $\tmeasure$ be a traversal measure on $\edges{\Zd}$ with $\ptimefamily = \ptimefamily[\tmeasure]$, and let $A\subseteq S\subseteq R\subseteq \Zd$. Then
\begin{enumerate}

\item \label{strong_restricted_proc:equiv_part}
For any $t\ge 0$, the following are equivalent.
\begin{enumerate}
\item \label{strong_restricted_proc:vertex_times_part}
If  $\vect v\in R$ and $\rptime{}{R}{A}{\vect v} \le t$, then $\rptime{}{R}{A}{\vect v} = \rptime{}{S}{A}{\vect v}$.
\item \label{strong_restricted_proc:process_times_part}
$\rreachedset{\tmeasure}{A}{R}{t'} = \rreachedset{\tmeasure}{A}{S}{t'}$ for all $t'\le t$.
\item \label{strong_restricted_proc:final_time_part}
$\rreachedset{\tmeasure}{A}{R}{t} = \rreachedset{\tmeasure}{A}{S}{t}$.
\item \label{strong_restricted_proc:star_process_part}
$\rreachedset{\tmeasure}{A}{\stargraph{R}}{t} = \rreachedset{\tmeasure}{A}{\stargraph{S}}{t}$ and $R\cap \rreachedset{\tmeasure}{A}{\stargraph{S}}{t}\subseteq S$.
\end{enumerate}

\item \label{strong_restricted_proc:containment_part}
The equivalent conditions in Part~\ref{strong_restricted_proc:equiv_part} imply that $\rreachedset{\tmeasure}{A}{R}{t}\subseteq S$, which in turn implies that $R\cap \rreachedset{\tmeasure}{A}{\stargraph{S}}{t}\subseteq S$. If $\tmeasure$ satisfies \geodesicsexist, then $R\cap \rreachedset{\tmeasure}{A}{\stargraph{S}}{t}\subseteq S$ implies the statements in Part~\ref{strong_restricted_proc:equiv_part}, so all six statements are equivalent in this case.

\item \label{strong_restricted_proc:hitting_time_part}
The equivalent statements in Part~\ref{strong_restricted_proc:equiv_part} are all satisfied if $\starptimenew{}{S}{A}{R\setminus S}> t$. If $\bd S\setminus \bd R$ is finite, then $R\cap \rreachedset{\tmeasure}{A}{\stargraph{S}}{t}\subseteq S$ implies that $\starptimenew{}{S}{A}{R\setminus S}> t$, so all seven conditions from Parts 1--3 are equivalent in this case.

%
%

\end{enumerate}
\end{lem}

\begin{proof}

\proofpart{\ref{strong_restricted_proc:equiv_part}} The implication $\ref{strong_restricted_proc:vertex_times_part} \implies \ref{strong_restricted_proc:process_times_part}$ is proved the same way as the corresponding statement from Lemma~\ref{restricted_proc_equiv_lem} by changing the appropriate strict inequality into a weak one, and the implication $\ref{strong_restricted_proc:process_times_part} \implies \ref{strong_restricted_proc:final_time_part}$ is trivial. We start by proving $\ref{strong_restricted_proc:final_time_part} \implies \ref{strong_restricted_proc:vertex_times_part}$.

Fix $t\ge 0$, and suppose $\rreachedset{\tmeasure}{A}{R}{t} = \rreachedset{\tmeasure}{A}{S}{t}$. Then for any $\vect v\in R$,
\begin{equation}
\label{strong_restricted_proc:R_implies_S_time_eqn}
\rptime{}{R}{A}{\vect v}\le t \implies \rptime{}{S}{A}{\vect v}\le t.
\end{equation}
If $\rptime{}{R}{A}{\vect v}= t$, then \eqref{strong_restricted_proc:R_implies_S_time_eqn} implies that $\rptime{}{S}{A}{\vect v} = t$ since we always have $\rptime{}{S}{A}{\vect v}\ge \rptime{}{R}{A}{\vect v}$, so the conclusion of \ref{strong_restricted_proc:vertex_times_part} holds in this case. Thus, suppose $\rptime{}{R}{A}{\vect v}< t$.
Note that we must have $\rptime{}{R}{A}{R\setminus S}\ge t$ since \eqref{strong_restricted_proc:R_implies_S_time_eqn} implies that $\rptime{}{R}{A}{\vect u}> t$ for all $\vect u\in R\setminus S$. Thus, $\rptime{}{R}{A}{\vect v}< t\le \rptime{}{R}{A}{R\setminus S}$, so Part~\ref{restricted_times_equiv:vertex_times_part} of Lemma~\ref{restricted_times_equiv_lem} implies that $\rptime{}{R}{A}{\vect v} = \rptime{}{S}{A}{\vect v}$ in this case as well.


Next we show that the first statement in \eqref{strong_restricted_proc:star_process_part} follows from \eqref{strong_restricted_proc:vertex_times_part}, and the second statement follows from \eqref{strong_restricted_proc:final_time_part}. Suppose \eqref{strong_restricted_proc:vertex_times_part} holds, and let $\vect v\in\Zd$ with $\starptimenew{}{R}{A}{\vect v} \le t$. Then there exists a sequence $\setof{\gamma_n}_{n\in\N}$ of minimal $\stargraph{R}$-paths from $A$ to $\vect v$ with $\inf_n \tmeasure(\gamma_n) = \starptimenew{}{R}{A}{\vect v}$. Since there are only finitely many edges at $\vect v$, there is some edge $\edge = \setof{\vect v,\vect u}$ contained in infinitely many of the $\gamma_n$'s; without loss of generality, assume $\edge\in \gamma_n$ for all $n$. Since each $\gamma_n$ is a minimal $\stargraph{R}$-path and $A\subseteq R$, each subpath $\subpath{\gamma_n}{A}{\vect u}$ must be an $R$-path, so
\begin{equation}
\label{strong_restricted_proc:TR<TRstar_eqn}
\rptime{}{R}{A}{\vect u} \le \inf_n \tmeasure \parens[2]{\subpath{\gamma_n}{A}{\vect u}}
= \inf_n \brackets[2]{\tmeasure(\gamma_n) -\ttime(\edge)}
=  \starptimenew{}{R}{A}{\vect v} - \ttime(\edge).
\end{equation}
Since $\starptimenew{}{R}{A}{\vect v}\le t$ by assumption, \eqref{strong_restricted_proc:TR<TRstar_eqn} implies that $\rptime{}{R}{A}{\vect u}\le t$, so $\rptime{}{R}{A}{\vect u} = \rptime{}{S}{A}{\vect u}$ by \eqref{strong_restricted_proc:vertex_times_part}. 
%
Now let $\setof{\lambda_n}_{n\in\N}$ be a sequence of $S$-paths from $A$ to $\vect u$ with $\inf_n \tmeasure(\lambda_n) = \rptime{}{S}{A}{\vect u} = \rptime{}{R}{A}{\vect u}$, and let $\lambda_n\edge$ denote the concatenation of $\lambda_n$ with the edge $\edge$ (plus its endvertices $\vect u$ and $\vect v$). Then each $\lambda_n\edge$ is an $\stargraph{S}$-path from $A$ to $\vect v$, so we have
\[
\starptimenew{}{S}{A}{\vect v} \le \inf_n \tmeasure(\lambda_n \edge)
= \inf_n \brackets[2]{\tmeasure(\lambda_n) + \ttime(\edge)}
= \rptime{}{R}{A}{\vect u} + \ttime(\edge)
\le \starptimenew{}{R}{A}{\vect v},
\]
where the final inequality follows from \eqref{strong_restricted_proc:TR<TRstar_eqn}. Since $S\subseteq R$, the reverse inequality $\starptimenew{}{S}{A}{\vect v}\ge \starptimenew{}{R}{A}{\vect v}$ is trivial, so we must have $\starptimenew{}{S}{A}{\vect v} = \starptimenew{}{R}{A}{\vect v}$. The above argument shows that if \eqref{strong_restricted_proc:vertex_times_part} holds, then
\[
\forall \vect v\in\Zd,\quad
\starptimenew{}{R}{A}{\vect v} \le t \implies
\starptimenew{}{R}{A}{\vect v} = \starptimenew{}{S}{A}{\vect v},
\]
which implies that $\rreachedset{\tmeasure}{A}{\stargraph{R}}{t} = \rreachedset{\tmeasure}{A}{\stargraph{S}}{t}$. Now, clearly \eqref{strong_restricted_proc:final_time_part} implies that $\rreachedset{\tmeasure}{A}{R}{t}\subseteq S$, and the implication $\rreachedset{\tmeasure}{A}{R}{t}\subseteq S \implies R\cap \rreachedset{\tmeasure}{A}{\stargraph{S}}{t}\subseteq S$ is proved the same way as the corresponding statement from Lemma~\ref{restricted_proc_equiv_lem} by changing the appropriate strict inequality into a weak one. Thus, \eqref{strong_restricted_proc:star_process_part} holds.

Finally we show that  $\ref{strong_restricted_proc:star_process_part} \implies \ref{strong_restricted_proc:final_time_part}$. Suppose that \eqref{strong_restricted_proc:star_process_part} holds, and let $\vect v\in R$ with $\rptime{}{R}{A}{\vect v}\le t$. Then since $\vect v\in R$, we have $\starptimenew{}{R}{A}{\vect v} = \rptime{}{R}{A}{\vect v}\le t$, so $\starptimenew{}{S}{A}{\vect v}\le t$ by the first statement in \eqref{strong_restricted_proc:star_process_part}. But then $\vect v\in R\cap \rreachedset{\tmeasure}{A}{\stargraph{S}}{t}$, so the second statement in \eqref{strong_restricted_proc:star_process_part} implies that $\vect v\in S$ and hence $\rptime{}{S}{A}{\vect v} = \starptimenew{}{S}{A}{\vect v}\le t$. Therefore, $\rptime{}{R}{A}{\vect v}\le t \implies \rptime{}{S}{A}{\vect v}\le t$, and since the reverse implication is trivial, \eqref{strong_restricted_proc:final_time_part} holds.

\proofpart{\ref{strong_restricted_proc:containment_part}} Clearly \ref{strong_restricted_proc:final_time_part} implies that $\rreachedset{\tmeasure}{A}{R}{t}\subseteq S$, and the implication $\rreachedset{\tmeasure}{A}{R}{t}\subseteq S \implies R\cap \rreachedset{\tmeasure}{A}{\stargraph{S}}{t}\subseteq S$ is proved the same way as the orresponding statement from Lemma~\ref{restricted_proc_equiv_lem} by changing the appropriate strict inequality into a weak one. We will prove that if $\tmeasure$ satisfies \geodesicsexist, then $R\cap \rreachedset{\tmeasure}{A}{\stargraph{S}}{t}\subseteq S \implies \ref{strong_restricted_proc:vertex_times_part}$.

Suppose that $R\cap \rreachedset{\tmeasure}{A}{\stargraph{S}}{t}\subseteq S$, i.e.\ if $\vect u\in R$ and $\starptimenew{}{S}{A}{\vect u}\le t$, then $\vect u\in S$. Let $\vect v\in \Zd$ with $\rptime{}{R}{A}{\vect v}\le t$. Since $\tmeasure$ satisfies \geodesicsexist, there is some $\rptmetric{R}$-geodesic $\gamma$ from $A$ to $\vect v$, so $\tmeasure(\gamma) = \rptime{}{R}{A}{\vect v}\le t$. We will show that $\gamma\subseteq S$ by contradiction. If $\gamma$ contains some point outside $S$, let $\vect u$ be the \emph{first} such point. Then $\vect u\in \nbrset{S}$, and the subpath $\subpath{\gamma}{A}{\vect u}$ is an $\stargraph{S}$-path. Therefore,
\[
\starptimenew{}{S}{A}{\vect u} \le \tmeasure \parens[2]{\subpath{\gamma}{A}{\vect u}}
\le \tmeasure(\gamma) \le t,
\]
which implies that $\vect u\in S$ by hypothesis. This contradicts our choice of $\vect u$, so we conclude that $\gamma\subseteq S$. Therefore, $\gamma$ is an $S$-path, so
\[
\rptime{}{S}{A}{\vect v} \le \tmeasure(\gamma) = \rptime{}{R}{A}{\vect v}.
\]
Since $S\subseteq R$, the reverse inequality is trivial, so we must have $\rptime{}{S}{A}{\vect v}= \rptime{}{R}{A}{\vect v}$.

\proofpart{\ref{strong_restricted_proc:hitting_time_part}} The first part follows from Lemma~\ref{restricted_times_equiv_lem} the same way as in the proof of Lemma~\ref{restricted_proc_equiv_lem}, and the second part follows because if $\bd S\setminus \bd R$ is finite, then $\bd (R\setminus S) \subseteq \nbrset{\bd S\setminus \bd R}$ is also finite, so the infimum defining $\starptimenew{}{S}{A}{R\setminus S}$ becomes a minimum. That is, $R\cap \rreachedset{\tmeasure}{A}{\stargraph{S}}{t}\subseteq S$ means that $\starptimenew{}{S}{A}{\vect v} > t$ for all $\vect v\in R\setminus S$, and since any $\stargraph{S}$-path from $A$ to $R\setminus S$ must end in $\bd (R\setminus S)$, the assumption that $\bd (R\setminus S)$ is finite implies that
\[
\starptimenew{}{S}{A}{R\setminus S}
= \inf_{\vect v\in R\setminus S} \starptimenew{}{S}{A}{\vect v}
= \inf_{\vect v\in \bd(R\setminus S)} \starptimenew{}{S}{A}{\vect v}
= \min_{\vect v\in \bd(R\setminus S)} \starptimenew{}{S}{A}{\vect v}>t.
\qedhere
\]
%
%
\end{proof}

Lemma~\ref{strong_restricted_proc_equiv_lem} implies the following relationship between the one-type and two-type processes.

\begin{lem}[Equivalence of one-type and two-type processes]
\label{one_type_two_type_lem}
Let $\tmeasurepair$ be a two-type traversal measure on $\edges{\Zd}$ satisfying \notiesdet\ and \mintimebddgrowth (i.e.\  $\tmeasurepair\in \goodsamplespace$), and let $\reachedfcn{\tmeasurepair}{\initconfig}$ be the corresponding two-type process started from $\initconfig$, as defined in Definition~\ref{two_type_proc_def}. For $i\in\setof{1,2}$ and $t\in [0,\infty]$, if $S$ is any subset of $\Zd$ such that
\[
\reachedset{i}{\initconfig}{t}\subseteq S \subseteq
\Zd\setminus \reachedset{3-i}{\initconfig}{t-},
\]
then
\[
\reachedset{i}{\initconfig}{t'} = \rreachedset{\tmeasure_i}{\initialset{i}}{S}{t'}
\quad\text{for all } t'\le t,
\]
where for $t=0$, we use the convention $\reachedset{3-i}{\initconfig}{0-} = \initialset{2}$.
\end{lem}

\begin{proof}
\newcommand{\smallset}{S_t}
\newcommand{\bigset}{R_t}


Let $\finalset{1}$ and $\finalset{2}$ be the finally conquered sets from Definition~\ref{two_type_proc_def}. By Propositions~\ref{entangled_full_prop} and \ref{entangled_well_defined_prop}, the assumptions on $\tmeasurepair$ imply that $\finalset{1}\cup \finalset{2} = \Zd$ and $\finalset{1}\cap \finalset{2} = \emptyset$, and also that each $\tmeasure_i$ satisfies \geodesicsexist. These three properties will be needed for the proof.

For concreteness, take $i=1$. First note that since $\finalset{1} = \Zd\setminus \finalset{2}$, if $t=\infty$, then we must have $S=\finalset{1}$ by Part~\ref{final_sets:are_union_part} of Lemma~\ref{final_sets_properties_lem}, so the result is trivial in this case. Thus, assume $t<\infty$. Let $\smallset = \reachedset{1}{\initconfig}{t}$ and $\bigset = \Zd\setminus \reachedset{2}{\initconfig}{t-}$. Then $\smallset\subseteq S\subseteq \bigset$ by assumption, so
\begin{equation}
\label{one_type_two_type:S_trapped_eqn}
\rreachedset{\tmeasure_1}{\initialset{1}}{\smallset}{t}
\subseteq \rreachedset{\tmeasure_1}{\initialset{1}}{S}{t}
\subseteq \rreachedset{\tmeasure_1}{\initialset{1}}{\bigset}{t}.
\end{equation}
Now, $\smallset = \reachedset{1}{\initconfig}{t} \subseteq \finalset{1}$ by definition, and since $\finalset{1}\cap \finalset{2} = \emptyset$ and $\reachedset{2}{\initconfig}{t-} \subseteq \finalset{2}$, we have $\finalset{1} \cap \reachedset{2}{\initconfig}{t-} = \emptyset$, so $\finalset{1}\subseteq \bigset$.  Therefore, $\smallset \subseteq \finalset{1} \subseteq \bigset$, and since $\reachedset{1}{\initconfig}{t} = \rreachedset{\tmeasure_1}{\initialset{1}}{\finalset{1}}{t}$, we have
\begin{equation}
\label{one_type_two_type:sp1_trapped_eqn}
\rreachedset{\tmeasure_1}{\initialset{1}}{\smallset}{t}
\subseteq \reachedset{1}{\initconfig}{t}
\subseteq \rreachedset{\tmeasure_1}{\initialset{1}}{\bigset}{t}.
\end{equation}
Our strategy will be to show that $\bigset \cap \rreachedset{\tmeasure_1}{\initialset{1}}{\stargraph{\smallset}}{t} \subseteq \smallset$ by contradiction, whence Lemma~\ref{strong_restricted_proc_equiv_lem} implies that $\rreachedset{\tmeasure_1}{\initialset{1}}{\bigset}{t} = \rreachedset{\tmeasure_1}{\initialset{1}}{\smallset}{t}$, and the desired result then follows from \eqref{one_type_two_type:S_trapped_eqn} and \eqref{one_type_two_type:sp1_trapped_eqn}.



Suppose that  $\bigset \cap \rreachedset{\tmeasure_1}{\initialset{1}}{\stargraph{\smallset}}{t} \not\subseteq \smallset$, i.e.\ $\starptimenew{\tmeasure_1}{\smallset}{\initialset{1}}{\vect v} \le t$ for some $\vect v\in \bigset\setminus \smallset$. Then $\starptimenew{\tmeasure_1}{\finalset{1}}{\initialset{1}}{\vect v} \le t$ since $\smallset\subseteq \finalset{1}$. If $\vect v\in \finalset{1}$, this implies $\rptime{\tmeasure_1}{\finalset{1}}{\initialset{1}}{\vect v}\le t$, so $\vect v\in \rreachedset{\tmeasure_1}{\initialset{1}}{\finalset{1}}{t} = \smallset$, a contradiction. Thus we have $\vect v\not\in \finalset{1}$, so $\vect v\in \finalset{2}$ since $\finalset{1}\cup \finalset{2} = \Zd$. Then since $\starptimenew{\tmeasure_1}{\finalset{1}}{\initialset{1}}{\vect v}\le t$, we must have $\vect v\in \initialset{2}$ or $\rptime{\tmeasure_2}{\finalset{2}}{\initialset{2}}{\vect v}< t$ by Part~\ref{final_sets:characterization_part} of Lemma~\ref{final_sets_properties_lem}. But then $\vect v\in \reachedset{2}{\initconfig}{t-} = \Zd\setminus \bigset$, also a contradiction. Thus we conclude that $\starptimenew{\tmeasure_1}{\smallset}{\initialset{1}}{\vect v} > t$ for all $\vect v\in \bigset\setminus \smallset$, which means $\bigset \cap \rreachedset{\tmeasure_1}{\initialset{1}}{\stargraph{\smallset}}{t} \subseteq \smallset$. Since $\tmeasure$ satisfies \geodesicsexist, Part~\ref{strong_restricted_proc:containment_part} of Lemma~\ref{strong_restricted_proc_equiv_lem} then implies that
\begin{equation}
\label{one_type_two_type:big_small_equal_eqn}
\rreachedset{\tmeasure_1}{\initialset{1}}{\bigset}{t}
= \rreachedset{\tmeasure_1}{\initialset{1}}{\smallset}{t}.
\end{equation}
Combining \eqref{one_type_two_type:big_small_equal_eqn} with \eqref{one_type_two_type:S_trapped_eqn} and \eqref{one_type_two_type:sp1_trapped_eqn}, we get
\[
\rreachedset{\tmeasure_1}{\initialset{1}}{S}{t}
= \reachedset{1}{\initconfig}{t}
= \rreachedset{\tmeasure_1}{\initialset{1}}{\finalset{1}}{t},
\]
and the conclusion of the lemma then follows from Part~\ref{strong_restricted_proc:equiv_part} of Lemma~\ref{strong_restricted_proc_equiv_lem}.
%
\end{proof}

\begin{thmremark}
The proof of Lemma~\ref{one_type_two_type_lem} shows that we can replace the assumption that $\tmeasurepair$ satisfies \notiesdet\ and \mintimebddgrowth\ with the weaker assumption that each $\tmeasure_i$ satisfies \geodesicsexist\ and that the finally conquered sets $\finalset{1}$ and $\finalset{2}$ partition $\Zd$. Moreover, if we drop the assumption \geodesicsexist\ and assume only that $\finalset{1}$ and $\finalset{2}$ partition $\Zd$, then using Lemma~\ref{restricted_proc_equiv_lem} instead of Lemma~\ref{strong_restricted_proc_equiv_lem}, the same proof shows that if $\reachedset{i}{\initconfig}{t-}\subseteq S \subseteq \Zd\setminus \reachedset{3-i}{\initconfig}{t-}$, then $\reachedset{i}{\initconfig}{t-} = \rreachedset{\tmeasure_i}{\initialset{i}}{S}{t-}$.
\end{thmremark}

The next result, which generalizes Lemma~\ref{entangled_disentangled_lem} for $\tmeasurepair\in\goodsamplespace$, follows immediately from Lemma~\ref{one_type_two_type_lem} and the monotonicity of the one-type process with respect to the restricting set.

\begin{cor}[Comparison with disentangled processes]
\label{entangled_disentangled_cor}
Let $\tmeasurepair$ be a two-type traversal measure on $\edges{\Zd}$ satisfying \notiesdet\ and \mintimebddgrowth, and let $\reachedfcn{\tmeasurepair}{\initconfig}$ be the corresponding two-type process started from $\initconfig$. If $C\subseteq \Zd\setminus \reachedset{3-i}{\initconfig}{t-}\subseteq D$, then
\[
\rreachedset{\tmeasure_i}{\initialset{i}}{C}{t'}
\subseteq \reachedset{i}{\initconfig}{t'} \subseteq
\rreachedset{\tmeasure_i}{\initialset{i}}{D}{t'}
\quad\text{for all } t'\le t.
\]
\end{cor}


\subsection{Orderings and Monotonicity}
\label{fpc_basics:monotonicity_sec}

It is intuitively clear that the two-type process $\reachedset{\tmeasurepair}{\initconfig}{t}$ should satisfy a sort of monotonicity property with respect to the initial configuration $\initconfig$ and the collection of traversal times $\tmeasurepair$ used to run the process. For example, suppose we run the process twice with the same set of traversal times $\tmeasurepair$, but with different starting configurations. If we enlarge species 1's starting set from $\initialset{1}$ in the first process to some set $\initialset{1}'\supseteq \initialset{1}$ in the second process, while keeping species 2's starting set $\initialset{2}$ fixed, then species 1 should retain its advantage throughout the entire run of the second process. That is, at any time $t\ge 0$, the set of sites occupied by species 1 in the second process should be at least as large as the corresponding set in the first process, and the set of 2's in the second process should be no larger than the set of 2's in the first process. Similarly, if we keep the starting configuration fixed but decrease some of species 1's traversal times while leaving species 2's traversal times untouched, then again, species 1 will have an advantage throughout the entire run of the second process compared with the first.

Our goal in this section will be to prove this intuitive monotonicity property, in Lemma~\ref{monotonicity_lem} below. In order to state the result, we will work on the canonical sample space $\samplespace=\samplespace[1]$, and as described in Section~\ref{fpc_basics:additional_notation_sec}, we consider the two-type process as a map $\process{}\colon \domain\to \processspace$, where $\statespace=\statespace[1]$ is the state space, $\domain\subseteq\domainspace$ is some domain of definition for the process, and $\processspace$ is the canonical process space. Lemma~\ref{monotonicity_lem} is stated in terms of appropriate orderings on $\domainspace$ and $\processspace$, which we now define by starting with orderings on the lower-level spaces $\setof{0,1,2}$, $\statespace$, and $\samplespace$.

\subsubsection{Ordering the Canonical Spaces by Favorability to Species 1}

By convention, we will choose our orderings to favor species 1. That is, for comparable objects $\mathsf{x}$ and $\mathsf{y}$, we will say that $\mathsf{x}\ge \mathsf{y}$ if $\mathsf{x}$ is better for species 1 (or equivalently, worse for species 2) than $\mathsf{y}$. For example, a vertex in state 1 is strictly better than a vertex in state 0, which is strictly better than a vertex in state 2. This defines a linear ordering on the set of vertex states $\setof{0,1,2}$, which we can easily remember by making the identification $\setof{0,1,2}\equiv \setof{0,1,-1} \pmod 3$:
%
%
\begin{equation}\label{basic_order_eqn}
1>0>-1\equiv 2.
\end{equation}
The ordering \eqref{basic_order_eqn} induces the (componentwise) product partial order on $\statespace = \statespace[1] \equiv \setof{0,1,-1}^{\Zd} \pmod 3$, which agrees with the usual partial order on functions: If $X,Y\in\statespace$, then $X\ge Y$ iff $X(\vect v)\ge Y(\vect v)$ for all $\vect v\in\Zd$, where the latter instance of $\ge$ refers to \eqref{basic_order_eqn}. Equivalently,
\begin{equation}\label{state_order_eqn}
\forall X,Y\in\statespace,\quad
\text{$X\ge Y$ iff $\stateset{1}{X}\supseteq \stateset{1}{Y}$ and $\stateset{2}{X}\subseteq \stateset{2}{Y}$,}
\end{equation}
where the occupied sets $\reachedfcn{i}{}$ are defined in \eqref{state_occupied_sets_def_eqn}. The ordering \eqref{state_order_eqn} in turn induces the product partial order on the space $\statespace^{[0,\infty]}$ where the two-type process lives:
\begin{equation}\label{process_order_eqn}
\forall \zeta_\cdot,\vartheta_\cdot \in \statespace^{[0,\infty]},\quad
\text{$\zeta_\cdot\ge \vartheta_\cdot$ iff $\zeta_t\ge \vartheta_t$ for all $t\in [0,\infty]$,}
\end{equation}
where the latter instance of $\ge$ refers to \eqref{state_order_eqn}.

For two collections of traversal time pairs $\outcome,\outcome'\in \samplespace = \samplespace[1]$, the configuration $\outcome$ is better for species 1 than the configuration $\outcome'$ if the traversal times in $\outcome$ allow species 1 to move faster and species 2 to move slower than those in $\outcome'$. More generally, we define the following partial order on $\rsamplespace{\Lambda} = (\R_+\times\R_+)^\Lambda$, for $\Lambda\subseteq \edges{\Zd}$:
\begin{equation}\label{restricted_ttime_order_eqn}
\forall f,g \in \rsamplespace{\Lambda},\quad
\text{$f \ge g$ iff $f_1(\edge) \le g_1(\edge)$ and $f_2(\edge)\ge g_2(\edge)$ for all $\edge\in\Lambda$,}
\end{equation}
where $f=(f_1,f_2)$ and $g=(g_1,g_2)$ for some $f_i,g_i\colon \Lambda\to\R_+$. Note that \eqref{restricted_ttime_order_eqn} is just the product partial order induced by the partial order on $\R_+\times \R_+$ given by $(x_1,x_2)\ge (y_1,y_2)$ iff $x_1\le y_1$ and $x_2\ge y_2$.

The orderings \eqref{state_order_eqn} and \eqref{restricted_ttime_order_eqn} induce the product partial order on $\samplespace\times\statespace$, namely, for two pairs $\pair{\outcome}{X}$ and $\pair{\outcome'}{X'}$ in $\samplespace\times\statespace$,
\begin{equation}\label{pair_product_order_eqn}
\text{$\pair{\outcome}{X}\ge \pair{\outcome'}{X'}$ iff $\outcome\ge \outcome'$ and $X\ge X'$.}
\end{equation}
The partial order \eqref{pair_product_order_eqn} is actually stronger than necessary to get monotonicity for the two-type process. In particular, in order to compare two processes started from initial states $X$ and $X'$ with $X\ge X'$, we do not need any information about the traversal times within the initially occupied regions $\nonzeroset{X}$ and $\nonzeroset{X'}$.
Therefore, we define the following \emph{preorder} (i.e.\ partial order without antisymmetry) $\pbetter$ on $\samplespace\times\statespace$, which treats the traversal times of two input pairs as equivalent if they agree outside the total occupied region.
That is, for $\pair{\outcome}{X}$ and $\pair{\outcome'}{X'}$ in $\samplespace\times\statespace$, we define
\begin{equation}\label{pair_preorder_eqn}
\pair{\outcome}{X}\pbetter \pair{\outcome'}{X'}
\text{ iff $X\ge X'$ and $\restrict{\outcome}{\compedges{X\cup X'}}\ge \restrict{\outcome'}{\compedges{X\cup X'}}$},
\end{equation}
where the partial orders on the right refer to \eqref{state_order_eqn} and \eqref{restricted_ttime_order_eqn}, and $X\cup X' = \nonzeroset{X} \cup \nonzeroset{X'}$. More explicitly, \eqref{pair_preorder_eqn} says that $\pair{\outcome}{X}\pbetter \pair{\outcome'}{X'}$ if
\[
\stateset{1}{X}\supseteq \stateset{1}{X'}
\text{ and }
\stateset{2}{X}\subseteq \stateset{2}{X'},
\]
and for all edges $\edge\notin \edges[2]{\stateset{1}{X}\cup \stateset{2}{X'}}$,
\[
\outcome_1(\edge)\le \outcome_1'(\edge)
\text{ and }
\outcome_2(\edge)\ge \outcome_2'(\edge).
\]
As with any preorder, \eqref{pair_preorder_eqn} induces an equivalence relation defined by $\pair{\outcome}{X}\pequiv \pair{\outcome'}{X'}$ iff $\pair{\outcome}{X}\pbetter \pair{\outcome'}{X'}$ and $\pair{\outcome}{X}\pworse \pair{\outcome'}{X'}$. That is,
\begin{equation}\label{pair_equiv_eqn}
\text{$\pair{\outcome}{X}\pequiv \pair{\outcome'}{X'}$ iff $X=X'$ and $\outcome(\edge) = \outcome'(\edge)$ for all $\edge\notin \edges{X}$.}
\end{equation}
The preorder $\pbetter$ corresponds to a partial order on the set of equivalence classes $\samplespace\times \statespace / \pequiv$ in the obvious way.

\subsubsection{Monotonicity of $\process{}$ and Extension of Domain via Equivalence Classes}

We are now ready to state the main monotonicity result for the first-passage competition process.
\begin{lem}[Monotonicity of infection]
\label{monotonicity_lem}
Let $\samplespace = \samplespace[1]$ and $\statespace = \statespace[1]$. If $\domain \subseteq \samplespace\times\statespace$ is any domain of definition for the two-type process $\process{}$, then $\process{}: \domain \to \statespace^{[0,\infty]}$ is increasing with respect to the preorder \eqref{pair_preorder_eqn} and the partial order \eqref{process_order_eqn}. That is, if $\pair{\outcome}{X}$ and $\pair{\outcome'}{X'}$ are two pairs in $\samplespace\times\statespace$ on which $\process{}$ is defined, and $\pair{\outcome}{X}\pbetter \pair{\outcome'}{X'}$, then $\processtoutcome[2]{X}{t}{\outcome} \ge \processtoutcome[2]{X'}{t}{\outcome'}$ for all $t\ge 0$.
\end{lem}

\begin{proof}
Let $X\equiv \pair{A_1}{A_2}$ and $X'\equiv \pair{A_1'}{A_2'}$. Then $\pair{\outcome}{X}\pbetter \pair{\outcome'}{X'}$ means that
\[
A_1\supseteq A_1' \text{ and } A_2\subseteq A_2',
\]
and for all edges $\edge\notin \edges{A_1\cup A_2'}$,
\[
\outcome_1(\edge)\le \outcome_1'(\edge)
\text{ and }
\outcome_2(\edge)\ge \outcome_2'(\edge).
\]
For $i\in\setof{1,2}$, let $\finalset{i} = \reachedset{i}{X}{\infty}_\outcome$ and $\finalset{i}' = \reachedset{i}{X'}{\infty}_{\outcome'}$ be the final conquered sets for the entangled processes run with inputs $\pair{\outcome}{X}$ and $\pair{\outcome'}{X'}$, respectively, as constructed in Definition~\ref{two_type_proc_def}.
Our first goal is to show that $\finalset{1}\supseteq \finalset{1}'$ and $\finalset{2}\subseteq \finalset{2}'$, i.e.\ $\processtoutcome{X}{\infty}{\outcome}\ge \processtoutcome{X'}{\infty}{\outcome'}$.

Let $S'$ be a subset of $\Zd$ satisfying the criterion defining $\finalset{1}'$, i.e.\ $A_1'\subseteq S'$, $A_2'\cap S' = \emptyset$, and
\begin{equation}
\label{monotonicity:final_prime_ineq}
\rptime{\outcome_1'}{S'}{A_1'}{\vect v}< \starptimenew{\outcome_2'}{\parens{\Zd\setminus S'}}{A_2'}{\vect v}
\text{ for all } \vect v\in\bd S'\setminus A_1'.
\end{equation}
We claim that the set $S\definedas S'\cup A_1$ satisfies the criterion defining $\finalset{1}$, which then implies that $\finalset{1}'\subseteq \finalset{1}$ since $S'\subseteq S$. We have $A_1\subseteq S$ by definition, and since $A_1\cap A_2=\emptyset$ and $A_2\subseteq A_2'$, we have $A_2\cap S \subseteq A_2'\cap S' = \emptyset$. Furthermore, since $A_1'\subseteq A_1$, if $\vect v\in \bd S\setminus A_1$, then $\vect v\in \bd S'\setminus A_1'$.
Now let $\vect v\in \bd S\setminus A_1$, let $\gamma_1$ be any $S'$-path from $A_1$ to $\vect v$, and let $\gamma_2$ be any $\stargraph[0]{\Zd\setminus S}$-path from $A_2'$ to $\vect v$. Then since $S'\cap A_2'=\emptyset$, we have $\edges{\gamma_1}\subseteq \edgecomplement{\edges{A_1\cup A_2'}}$, and since $A_1\subseteq S$, we have $\edges{\gamma_2}\subseteq \edgecomplement{\edges{A_1\cup A_2'}}$, and hence
\begin{equation}
\label{monotonicity:path_comparison_eqn}
\tmeasure_{\outcome_1}(\gamma_1) \le \tmeasure_{\outcome_1'}(\gamma_1)
\text{ and }
\tmeasure_{\outcome_2}(\gamma_2) \ge \tmeasure_{\outcome_2'}(\gamma_2)
\end{equation}
since $\restrict{\outcome}{\edgecomplement{\edges{A_1\cup A_2'}}}\ge \restrict{\outcome'}{\edgecomplement{\edges{A_1\cup A_2'}}}$. Since $\gamma_1$ was arbitrary, and since $S\supseteq S'$ and $A_1\supseteq A_1'$, \eqref{monotonicity:path_comparison_eqn} implies
\[
\rptime{\outcome_1}{S}{A_1}{\vect v}
\le \rptime{\outcome_1}{S'}{A_1}{\vect v}
\le \rptime{\outcome_1'}{S'}{A_1}{\vect v}
\le\rptime{\outcome_1'}{S'}{A_1'}{\vect v}.
\]
Since $\gamma_2$ was arbitrary, and since $A_2\subseteq A_2'$ and  $S\supseteq S'$, \eqref{monotonicity:path_comparison_eqn} implies
\[
\starptimenew{\outcome_2}{\parens{\Zd\setminus S}}{A_2}{\vect v}
\ge \starptimenew{\outcome_2}{\parens{\Zd\setminus S}}{A_2'}{\vect v}
\ge \starptimenew{\outcome_2'}{\parens{\Zd\setminus S}}{A_2'}{\vect v}
\ge \starptimenew{\outcome_2'}{\parens{\Zd\setminus S'}}{A_2'}{\vect v}.
\]
Combining these inequalities with \eqref{monotonicity:final_prime_ineq}, for all $\vect v\in\bd S\setminus A_1$ we have
\begin{equation}
\rptime{\outcome_1}{S}{A_1}{\vect v}\le \rptime{\outcome_1'}{S}{A_1}{\vect v}
\le \rptime{\outcome_1'}{S'}{A_1'}{\vect v}
< \starptimenew{\outcome_2'}{\parens{\Zd\setminus S'}}{A_2'}{\vect v}
\le \starptimenew{\outcome_2}{\parens{\Zd\setminus S}}{A_2}{\vect v},
\end{equation}
and so $S$ satisfies the criterion defining $\finalset{1}$ as claimed. Therefore, $\finalset{1}'\subseteq \finalset{1}$, and the same argument with 1's, 2's, and primes switched shows that $\finalset{2}\subseteq \finalset{2}'$.

Now, $\finalset{1}'$ satisfies all the properties assumed of $S'$, so the above argument shows that
$
\rptime{\outcome_1}{\finalset{1}'}{A_1}{\vect v}
\le \rptime{\outcome_1'}{\finalset{1}'}{A_1'}{\vect v}
$
for all $\vect v\in \finalset{1}'$, and thus, since $\finalset{1}\supseteq \finalset{1}'$,
\[
\reachedset{1}{X}{t}_\outcome = \rreachedset{\outcome_1}{A_1}{\finalset{1}}{t}
\supseteq \rreachedset{\outcome_1}{A_1}{\finalset{1}'}{t}
\supseteq \rreachedset{\outcome_1'}{A_1'}{\finalset{1}'}{t}
= \reachedset{1}{X'}{t}_{\outcome'}
\quad\text{for all } t\ge 0.
\]
Similarly, $\rptime{\outcome_2}{\finalset{2}}{A_2}{\vect v} \ge \rptime{\outcome_2'}{\finalset{2}}{A_2'}{\vect v}$ for all $\vect v\in \finalset{2}$, and thus, since $\finalset{2}\subseteq \finalset{2}'$,
\[
\reachedset{2}{X}{t}_\outcome = \rreachedset{\outcome_2}{A_2}{\finalset{2}}{t}
\subseteq \rreachedset{\outcome_2'}{A_2'}{\finalset{2}}{t}
\subseteq \rreachedset{\outcome_2'}{A_2'}{\finalset{2}'}{t}
= \reachedset{2}{X'}{t}_{\outcome'}
\quad\text{for all } t\ge 0.
\]
Therefore, $\processtoutcome{X}{t}{\outcome}\ge \processtoutcome{X'}{t}{\outcome'}$ for all $t\in [0,\infty]$.
\end{proof}

The following result is an immediate consequence of Lemma~\ref{monotonicity_lem}.

\begin{cor}[Irrelevance of internal traversal times; extension of domain]
\label{domain_extension_cor}
If $\pair{\outcome}{X}\pequiv \pair{\outcome'}{X'}$ as defined in \eqref{pair_equiv_eqn}, then
\begin{equation}\label{quotient_space_eqn}
\processtoutcome[2]{X}{t}{\outcome} = \processtoutcome[2]{X'}{t}{\outcome'}
\text{\ \ for all } t\ge 0.
\end{equation}
That is, if $\process{}$ is already defined on the equivalent pairs $\pair{\outcome}{X}$ and $\pair{\outcome'}{X'}$, then \eqref{quotient_space_eqn} holds, and if $\process{}$ is defined on $\pair{\outcome}{X}$ but not $\pair{\outcome'}{X'}$, then we take \eqref{quotient_space_eqn} as the definition of $\processtoutcome[2]{X'}{t}{\outcome'}$.
In this way, if $\domain \subseteq \samplespace\times\statespace$ is any domain of definition for $\process{}$, then $\process{}$ extends consistently via \eqref{quotient_space_eqn} to $\widetilde{\domain} \definedas \setbuilder[2]{\pair{\outcome'}{X'}\in \samplespace\times\statespace}{\exists \pair{\outcome}{X}\in \domain \text{ with } \pair{\outcome'}{X'}\pequiv \pair{\outcome}{X}}$.
\end{cor}

\subsubsection{Monotone Subsets of Ordered Spaces}

If $\pair{P}{\lesssim}$ is a preordered set, a subset $S\subseteq P$ is called \textdef{increasing} (\textdef{decreasing}) if its indicator function $\ind{S}$ is increasing (resp.\ decreasing) with respect to $\lesssim$ and the usual order on $\R$. Equivalently, $S$ is increasing (decreasing) if whenever $x\in S$ and $x\lesssim y$ (resp.\ $x\gtrsim y$), it follows that $y\in S$.

%

Note that Lemma~\ref{monotonicity_lem} implies that if $H$ is an increasing subset of the state space $\statespace$, and $\domain$ is any domain of definition for $\process{}$, then for any $t\in [0,\infty]$, the set
$
\setbuilder[1]{\pair{\outcome}{X}\in D}{\processtoutcome{X}{t}{\outcome}\in H}
$
is an increasing subset of $\samplespace\times\statespace$ with respect to the preorder \eqref{pair_preorder_eqn}, since its indicator function is the composition of monotone maps $\ind{H}\circ \processt{}{t}$.

For our purposes, we consider increasing subsets of the canonical sample space $\samplespace$ (i.e.\ events) and of the state space $\statespace$. The following lemma will be used in the proof of Theorem~\ref{irrelevance_thm} below and in Chapter~\ref{coex_finite_chap}.

\begin{lem}[Increasing sets]
\label{increasing_sets_lem}
Let $\probabilityspace$ and $\statemeasurespace$ be the canonical probability space and state space defined in Sections~\ref{fpc_basics:prob_space_sec} and \ref{fpc_basics:ips_def_sec}. If $H\in \statesigmafield$ is increasing, then for any $X\in \statespace$ and $t\in [0,\infty]$, the event $\eventthat[2]{\processt{X}{t} \in H}$ is increasing. If additionally $X'\in S$ and $X\le X'$, then
\[
\Pr \eventthat[3]{\processt{X}{t}\in H} \le \Pr \eventthat[3]{\processt{X'}{t}\in H}.
\]
\end{lem}

\begin{proof}
Suppose $\processtoutcome[2]{X}{t}{\outcome}\in H$ and $\outcome\le \outcome'$. Then $\pair{\outcome}{X} \le \pair{\outcome'}{X}$, and since the partial order \eqref{pair_product_order_eqn} is stronger than the preorder \eqref{pair_preorder_eqn}, Lemma~\ref{monotonicity_lem} implies that $\processtoutcome[2]{X}{t}{\outcome}\le \processtoutcome[2]{X}{t}{\outcome'}$, so $\processtoutcome[2]{X}{t}{\outcome'}\in H$ since $H$ is increasing. Thus, $\eventthat[2]{\processt{X}{t} \in H}$ is an increasing event. If $X\le X'$, then we also have $\processtoutcome[2]{X}{t}{\outcome}\le \processtoutcome[2]{X'}{t}{\outcome}$ by Lemma~\ref{monotonicity_lem}, so $\processtoutcome[2]{X'}{t}{\outcome}\in H$ since $H$ is increasing. Thus,
\[
\setbuilder[3]{\outcome}{\processtoutcome[2]{X}{t}{\outcome}\in H}
\subseteq
\setbuilder[3]{\outcome}{\processtoutcome[2]{X'}{t}{\outcome}\in H}.\qedhere
\]
\end{proof}

\subsection{Shift Operators on $\samplespace\times \statespace$, and a ``Markov-esque" Property}
\label{fpc_basics:shift_ops_sec}

%
%
%
%
%
%
%
%

Our goal in this section will be to prove a property similar to the Markov property for the first-passage competition process, even when the process is not Markov. However, rather than being a probabilistic result, our ``Markov-esque" property will hold pointwise on the initial domain $\startdomain \subset \domainspace$ defined in \eqref{startdomain_def_eqn}, Section~\ref{fpc_basics:additional_notation_sec}. Essentially, this Markovesque property says that if we keep track of the infection times of the vertices on the boundary of the total occupied region, then we can pause the process at time $t$ and then restart it, feeding in the current configuration and the infection times of the boundary vertices as the new initial input, and the process will continue to evolve as it would have if we had allowed it to evolve from time 0 without pausing.

We first prove the Markovesque property for the one-type process in Lemma~\ref{one_type_markovesque_lem} below, and then we use this result to show that the analogous property holds for the two-type process in Lemma~\ref{markovesque_lem}.
The first step will be to define ``shift operators" that compute adjusted traversal times for edges on the boundary of the infected region, using the current state of the process and the infection times of the boundary vertices.
These shift operators, defined in Definitions~\ref{one_type_shift_def} and \ref{two_type_shift_def} below, operate on the set $\onetypeproduct$ for the one-type process, and on domains of definition $\domain\subseteq \twotypeproduct$ for the two-type process, mapping a given input pair to a new input pair that we will use to run the restarted process. By comparison, for $t\in [0,\infty]$, let $\shiftop{t}\colon \processspace \to\processspace$ be the \textdef{natural shift operator} on the canonical process space for the two-type process, defined by $\shiftop{t} \zeta_\cdot \definedas \zeta_{t+\cdot}$, i.e.
\begin{equation}
\label{natural_shiftop_def_eqn}
\parens[2]{\shiftop{t} \zeta_\cdot}_s = \zeta_{t+s}\in \statespace
\quad\text{for all $s,t\in [0,\infty]$ and $\zeta_\cdot\in \processspace$.}
\end{equation}
The shift operators $\shiftop{t}$ can be used to state the Markov Property in the case when the two-type process is Markov, and we will come back to the definition \eqref{natural_shiftop_def_eqn} in Chapter~\ref{coex_finite_chap} when we study the Markov version of the process. One way to interpret the Markovesque property (Lemma~\ref{markovesque_lem} below) is to view the shift operator $\pshiftop{t}\colon \domain\to \twotypeproduct$ in Definition~\ref{two_type_shift_def} below as a sort of ``pullback" of the shift $\shiftop{t}$ through the process $\process{}\colon \domain\to \processspace$.

%
%
%
%

\subsubsection{Pausing and Restarting the One-Type Process}

We start by defining, for each $S\subseteq\Zd$ and $t\ge 0$, a shift operator $\rpshiftop{t}{S}\colon \onetypeproduct\to \onetypeproduct$ for the restricted one-type process $\rreachedfcn{\tmeasure}{A}{S}$, where $A\subseteq\Zd$ and $\tmeasure\in \onetypespace$.

\begin{definition}[Shift operators for the one-type process]
\label{one_type_shift_def}
For any restricting set $S\subseteq\Zd$ and $t\ge 0$, we define a shift operator $\rpshiftop{t}{S}$ on $\onetypeproduct$ as follows. Given any initial set $A\in \onetypestates$ and collection of traversal times $\ttime\in \onetypespace$, define $\rpshiftstate{t}{S}{\ttime}{A} \in \onetypestates$ and $\rpshifttimes{t}{S}{A}{\ttime}\in \onetypespace$ by
\[
\rpshiftstate{t}{S}{\ttime}{A} \definedas \rreachedset{\ttime}{A}{S}{t},
\]
and for $\edge\in\edges{\Zd}$,
\[
\rpshifttimes{t}{S}{A}{\ttime}(\edge) \definedas
\begin{cases}
0 & \text{if 
	$\edge$ joins $\vect u\in \rreachedset{\ttime}{A}{S}{t}$
	to $\vect v\in \rreachedset{\ttime}{A}{\stargraph{S}}{t}$},\\
\ttime(\edge) - \parens[2]{t - \rptime{\ttime}{S}{A}{\vect u}} &
	\text{if $\edge$ joins $\vect u\in \rreachedset{\ttime}{A}{S}{t}$
	to $\vect v\not\in \rreachedset{\ttime}{A}{\stargraph{S}}{t}$},\\
\ttime(\edge) & \text{otherwise}.
\end{cases}
\]
Then set $\rpshiftpair{t}{S}{\ttime}{A} \definedas \pair[1]{\rpshifttimes{t}{S}{A}{\ttime}}{\,\rpshiftstate{t}{S}{\ttime}{A}} \in \onetypeproduct$.
\end{definition}

\begin{thmremark}
To be sure that Definition~\ref{one_type_shift_def} makes sense, we need to verify that $\rpshifttimes{t}{S}{A}{\ttime}\in \onetypespace$, i.e.\ that the shifted traversal times are nonnegative. Clearly $\rpshifttimes{t}{S}{A}{\ttime}(\edge)\ge 0$ for any $\edge$ falling into the first case or third case of the definition, so we only need to worry about the second case, when $\edge$ joins $\vect u\in \rreachedset{\ttime}{A}{S}{t}$ to $\vect v\not\in \rreachedset{\ttime}{A}{\stargraph{S}}{t}$. Suppose $\edge = \setof{\vect u,\vect v}$ is such an edge. Since $\vect u\in S$ and $\edge = \setof{\vect u,\vect v}$, we have $\starptimenew{\ttime}{S}{A}{\vect v} \le \rptime{\ttime}{S}{A}{\vect u} + \ttime(e)$, because if $\gamma$ is any $S$-path from $A$ to $\vect u$, then the concatenation $\gamma \edge$ is an $\stargraph{S}$-path from $A$ to $\vect v$. On the other hand, since $\vect v\not\in \rreachedset{\ttime}{A}{\stargraph{S}}{t}$, we have $\starptimenew{\ttime}{S}{A}{\vect v}>t$. Therefore,
\[
\ttime(\edge)+ \rptime{\ttime}{S}{A}{\vect u} 
\ge\starptimenew{\ttime}{S}{A}{\vect v}>t,
\]
so we have
\[
\rpshifttimes{t}{S}{A}{\ttime}(\edge)
=\ttime(\edge) + \rptime{\ttime}{S}{A}{\vect u}-t
>0
\]
for all $\edge = \setof{\vect u,\vect v}$ with $\vect u\in \rreachedset{\ttime}{A}{S}{t}$ and $\vect v\not\in \rreachedset{\ttime}{A}{\stargraph{S}}{t}$.
This shows that the shifted traversal times are all nonnegative and hence define a valid traversal measure.
\end{thmremark}


\begin{lem}[Markovesque property for the one-type process]
\label{one_type_markovesque_lem}
\newcommand{\shiftedtimes}{\ttime_t}
\newcommand{\shiftedset}{A_t}
\newcommand{\shiftedpair}{\pair{\shiftedtimes}{\shiftedset}}

Let $\setbuilder[2]{\rpshiftop{t}{S}}{S\subseteq\Zd, t\ge 0}$ be the family of shift operators from Definition~\ref{one_type_shift_def}. Suppose $A\subseteq S\subseteq \Zd$, and for any $\ttime\in\onetypespace$ and $t\ge 0$, let $\shiftedpair = \rpshiftpair{t}{S}{\ttime}{A}$. Then
\begin{align}
\rptime{\ttime}{S}{A}{\vect v} &= \rptime{\shiftedtimes}{S}{\shiftedset}{\vect v} +t
\quad\text{whenever\quad $\rptime{\ttime}{S}{A}{\vect v}> t$ or $\rptime{\shiftedtimes}{S}{\shiftedset}{\vect v}>0$, and}
\label{one_type_markovesque:rptime_shift_eqn}\\
\starptimenew{\ttime}{S}{A}{\vect v} &= \starptimenew{\shiftedtimes}{S}{\shiftedset}{\vect v} +t
\quad\text{whenever\quad $\starptimenew{\ttime}{S}{A}{\vect v}> t$ or $\starptimenew{\shiftedtimes}{S}{\shiftedset}{\vect v}>0$.}
\label{one_type_markovesque:starptime_shift_eqn}
\end{align}
Moreover,
\begin{equation}
\label{one_type_markovesque:process_shift_eqn}
\rreachedset{\ttime}{A}{S}{t+s}
	= \rreachedset{\shiftedtimes}{\shiftedset}{S}{s}
\quad\text{and}\quad
\rreachedset{\ttime}{A}{\stargraph{S}}{t+s}
	= \rreachedset{\shiftedtimes}{\shiftedset}{\stargraph{S}}{s}
\quad \text{ for all } s> 0.
\end{equation}
\end{lem}

\begin{proof}
\newcommand{\shiftedtimes}{\ttime_t}
\newcommand{\shiftedset}{A_t}
\newcommand{\shiftedpair}{\pair{\shiftedtimes}{\shiftedset}}

\newcommand{\ivertex}[1]{\vect u_{}} 
\newcommand{\lastvertex}[1]{\vect u_{#1}} 


Definition~\ref{one_type_shift_def} was formulated specifically so that the shifted pair $\shiftedpair = \rpshiftpair{t}{S}{\ttime}{A}$ satisfies the following property, which will be the basis of the proof.

\begin{claim}
\label{one_type_markovesque:path_measure_claim}
If $\gamma$ is any $\stargraph{S}$-path
from $\vect u\in \bd\shiftedset$ to $\vect v\in \nbrhood{S} \setminus \rreachedset{\ttime}{A}{\stargraph{S}}{t}$ 
such that $\gamma\cap \shiftedset = \setof{\vect u}$, then
\begin{equation}
\label{one_type_markovesque:path_measure_eqn}
\ttime(\gamma) = \shiftedtimes(\gamma) + t - \rptime{\ttime}{S}{A}{\vect u}.
\end{equation}
In particular, if $\gamma$ is any minimal $\stargraph{S}$-path from $\shiftedset$ to $\vect v\in \nbrhood{S}\setminus \rreachedset{\ttime}{A}{\stargraph{S}}{t}$, then \eqref{one_type_markovesque:path_measure_eqn} holds for $\gamma$'s initial vertex $\vect u\in \bd \shiftedset$.
\end{claim}

\begin{proof}[Proof of Claim~\ref{one_type_markovesque:path_measure_claim}]
\eqref{one_type_markovesque:path_measure_eqn} holds because $\ttime$ agrees with $\shiftedtimes$ on all edges of $\gamma$ except the initial edge $\edge$ originating at $\vect u$, and for this edge, the difference is $\ttime(\edge) - \shiftedtimes(\edge) = t - \rptime{\ttime}{S}{A}{\vect u}$ (the assumption that $\gamma$ is an $\stargraph{S}$-path that ends outside of $\rreachedset{\ttime}{A}{\stargraph{S}}{t}$ is used to deduce this equation for the initial edge).
\end{proof}

Our first goal will be to prove the formula \eqref{one_type_markovesque:rptime_shift_eqn} relating the $\stargraph{S}$-restricted passage times. Once this is done, the formulas \eqref{one_type_markovesque:rptime_shift_eqn} and \eqref{one_type_markovesque:process_shift_eqn} for the $S$-restricted passage times and the growth processes will follow easily. We divide the proof of \eqref{one_type_markovesque:rptime_shift_eqn} into two claims.

\begin{claim}
\label{one_type_markovesque:shifted_times_bigger_claim}
 $\starptimenew{\ttime}{S}{A}{\vect v} \le \starptimenew{\shiftedtimes}{S}{\shiftedset}{\vect v} +t$ for all $\vect v\not\in \rreachedset{\ttime}{A}{\stargraph{S}}{t}$.
\end{claim}

\begin{proof}[Proof of Claim~\ref{one_type_markovesque:shifted_times_bigger_claim}]
Fix $\vect v\in \Zd \setminus \rreachedset{\ttime}{A}{\stargraph{S}}{t}$. If the right-hand side is infinite, then the inequality trivially holds, so we can assume that $\starptimenew{\shiftedtimes}{S}{\shiftedset}{\vect v}<\infty$, which means there is some $\stargraph{S}$-path connecting $\shiftedset$ to $\vect v$ (and hence $\vect v\in\nbrhood{S}$ since $\shiftedset\subseteq S$).
Then there exists a minimal $\stargraph{S}$-path from $\shiftedset$ to $\vect v$ by Lemma~\ref{minimal_path_lem}. Fix any such minimal path $\gamma$, and let $\ivertex{\gamma}\in \bd \shiftedset$ be the initial vertex of $\gamma$. Then $\gamma\cap \shiftedset = \setof{\ivertex{\gamma}}$, and hence 
$\ttime(\gamma) = \shiftedtimes(\gamma) + t - \rptime{\ttime}{S}{A}{\ivertex{\gamma}}$ by Claim~\ref{one_type_markovesque:path_measure_claim} since  $\vect v\in \nbrhood{S}\setminus \rreachedset{\ttime}{A}{\stargraph{S}}{t}$.
If $\lambda$ is any path ending at $\vect u$, let $\lambda \gamma$ denote the concatenation of $\lambda$ with $\gamma$, and note that $\ttime(\lambda\gamma) \le \ttime(\lambda)+\ttime(\gamma)$.
Then since the concatenation of any $S$-path with $\gamma$ is an $\stargraph{S}$-path,
\begin{align*}
\starptimenew{\ttime}{S}{A}{\vect v}
&= \inf \setbuilder[1]{\ttime (P)}{\text{$P$ is an $\stargraph{S}$-path from $A$ to $\vect v$}}\\
&\le \inf \setbuilder[1]{\ttime(\lambda \gamma)}{\text{$\lambda$ is an ${S}$-path from $A$ to $\vect u$}}\\
&\le \inf \setbuilder[1]{\ttime(\lambda)}{\text{$\lambda$ is an ${S}$-path from $A$ to $\vect u$}} + \ttime(\gamma)\\
&= \rptime{\ttime}{S}{A}{\ivertex{\gamma}} + \shiftedtimes(\gamma) + t - \rptime{\ttime}{S}{A}{\ivertex{\gamma}}\\
&= \shiftedtimes(\gamma) + t.
\end{align*}
Now, since $\gamma$ was an arbitrary minimal $\stargraph{S}$-path from $\shiftedset$ to $\vect v$, we have
\begin{align*}
\starptimenew{\ttime}{S}{A}{\vect v}
&\le \inf \setbuilder[2]{\shiftedtimes(\gamma) + t}{\text{$\gamma$ is a minimal $\stargraph{S}$-path from $\shiftedset$ to $\vect v$}}\\
&= \starptimenew{\shiftedtimes}{S}{\shiftedset}{\vect v} +t.
\qedhere
\end{align*}
\end{proof}

\begin{claim}
\label{one_type_markovesque:shifted_times_smaller_claim}
$\starptimenew{\ttime}{S}{A}{\vect v} \ge \starptimenew{\shiftedtimes}{S}{\shiftedset}{\vect v} +t$ for all $\vect v\not\in \rreachedset{\ttime}{A}{\stargraph{S}}{t}$.
\end{claim}

\begin{proof}[Proof of Claim~\ref{one_type_markovesque:shifted_times_smaller_claim}]
Fix $\vect v\in \Zd \setminus \rreachedset{\ttime}{A}{\stargraph{S}}{t}$. If the left-hand side of the inequality is infinite, then the statement is trivially true, so assume that $\starptimenew{\ttime}{S}{A}{\vect v}<\infty$. Then there is some $\stargraph{S}$-path from $A$ to $\vect v$ (so $\vect v\in\nbrhood{S}$ since $A\subseteq S$), and hence there exists some minimal such path by Lemma~\ref{minimal_path_lem}. 
If $\gamma$ is any minimal $\stargraph{S}$-path from $A$ to $\vect v$, let $\lastvertex{\gamma}$ be the last vertex of $\gamma$ to lie in $\shiftedset$; note that $\lastvertex{\gamma}\in \bd \shiftedset$ exists because the initial vertex of $\gamma$ is in $A\subseteq \shiftedset$ while the final vertex of $\gamma$ is $\vect v\notin \shiftedset$ by assumption. Then since $\gamma$ is minimal and therefore simple, the subpaths $\subpath{\gamma}{A}{\lastvertex{\gamma}}$ and $\subpath{\gamma}{\lastvertex{\gamma}}{\vect v}$ are edge-disjoint, so
$
\ttime(\gamma)
= \ttime \parens[2]{\subpath{\gamma}{A}{\lastvertex{\gamma}}}
	+ \ttime \parens[2]{\subpath{\gamma}{\lastvertex{\gamma}}{\vect v}}.
$
Now note that since $A\subseteq S$ and $\lastvertex{\gamma}\in \shiftedset\subseteq S$, the initial subpath $\subpath{\gamma}{A}{\lastvertex{\gamma}}$ must be an $S$-path, hence must have $\ttime$-measure at least $\rptime{\ttime}{S}{A}{\lastvertex{\gamma}}$. Moreover, the final subpath $\subpath{\gamma}{\lastvertex{\gamma}}{\vect v}$ is an $\stargraph{S}$-path with $\gamma\cap \shiftedset = \setof{\lastvertex{\gamma}}$, and we have $\vect v\in \nbrhood{S} \setminus \rreachedset{\ttime}{A}{\stargraph{S}}{t}$ by assumption, so Claim~\ref{one_type_markovesque:path_measure_claim} implies that
\begin{align*}
\ttime(\gamma)
&= \ttime \parens[2]{\subpath{\gamma}{A}{\lastvertex{\gamma}}}
	+ \ttime \parens[2]{\subpath{\gamma}{\lastvertex{\gamma}}{\vect v}}\\
&\ge \rptime{\ttime}{S}{A}{\lastvertex{\gamma}}
	+ \shiftedtimes \parens[2]{\subpath{\gamma}{\lastvertex{\gamma}}{\vect v}}
	+ t - \rptime{\ttime}{S}{A}{\lastvertex{\gamma}}\\
&= \shiftedtimes \parens[2]{\subpath{\gamma}{\lastvertex{\gamma}}{\vect v}} + t.
\end{align*}
Therefore, since $\gamma$ was an arbitrary minimal $\stargraph{S}$-path from $A$ to $\vect v$,  and since $\lastvertex{\gamma}\in \bd \shiftedset$ for any such $\gamma$ by definition, we have
\begin{align*}
\starptimenew{\ttime}{S}{A}{\vect v}
&= \inf \setbuilder[2]{\ttime(\gamma)}
	{\text{$\gamma$ is a minimal $\stargraph{S}$-path from $A$ to $\vect v$}}\\
&\ge \inf \setbuilder[3]{%
	\shiftedtimes \parens[2]{\subpath{\gamma}{\lastvertex{\gamma}}{\vect v}} + t}
	{\text{$\gamma$ is a minimal $\stargraph{S}$-path from $A$ to $\vect v$}}\\
&\ge \inf_{\vect u\in \bd \shiftedset} \setbuilder[2]{%
	\shiftedtimes \parens{P}+t}
	{\text{$P$ is an $\stargraph{S}$-path from $\vect u$ to $\vect v$}}\\
&\ge \inf \setbuilder[2]{%
	\shiftedtimes \parens{P}}
	{\text{$P$ is an $\stargraph{S}$-path from $\shiftedset$ to $\vect v$}} + t\\
&= \starptimenew{\shiftedtimes}{S}{\shiftedset}{\vect v} +t.
\qedhere
\end{align*}
\end{proof}

Claims~\ref{one_type_markovesque:shifted_times_bigger_claim} and \ref{one_type_markovesque:shifted_times_smaller_claim} 
together imply the formula \eqref{one_type_markovesque:starptime_shift_eqn} for the $\stargraph{S}$-restricted times. To obtain the $S$-restricted version \eqref{one_type_markovesque:rptime_shift_eqn}, let $\vect v\in \Zd\setminus \rreachedset{\ttime}{A}{S}{t}$. If $\vect v\in S$, then the result reduces to \eqref{one_type_markovesque:starptime_shift_eqn}. Otherwise, both sides of \eqref{one_type_markovesque:rptime_shift_eqn} are infinite, so the result is trivial. Finally, to obtain \eqref{one_type_markovesque:process_shift_eqn}, for $s>0$ we have
\begin{align*}
\rreachedset{\ttime}{A}{\stargraph{S}}{t+s}
&= \setbuilder[2]{\vect v\in \Zd}{\starptimenew{\ttime}{S}{A}{\vect v} \le t+s}\\
&= \setbuilder[2]{\vect v\in \Zd}{\starptimenew{\shiftedtimes}{S}{\shiftedset}{\vect v} \le s}
	&& \text{by \eqref{one_type_markovesque:starptime_shift_eqn}}\\
&= \rreachedset{\shiftedtimes}{\shiftedset}{\stargraph{S}}{s},
\end{align*}
so the proof is complete.
\end{proof}

In the next lemma
we introduce various partial orders and preorders on the space $\onetypeproduct$, similar to the orderings for the two-type process, and we show that the shift operator $\rpshiftop{t}{S}$ is increasing in all of its arguments. Our convention for the one-type process is that ``$\mathsf{x}\ge \mathsf{y}$" if ``$\mathsf{x}$ is better for \speciesname{red} than $\mathsf{y}$." The proof of Lemma~\ref{one_type_shift_monotone_lem} consists of checking definitions and is left to the reader; the only part that will be needed later is Part~\ref{one_type_shift_monotone:t_part}.

\begin{lem}[Monotonicity of the one-type shift]
\label{one_type_shift_monotone_lem}
Let $A,A',S,S'\subseteq \Zd$ and $\ttime,\ttime'\in \onetypespace$, and fix $t\ge 0$.
\begin{enumerate}
\item \label{one_type_shift_monotone:t_part}
The map $t\mapsto \rpshiftpair{t}{S}{\ttime}{A}$ is increasing: If $t \le t'$, then $\rpshiftpair{t}{S}{\ttime}{A} \le \rpshiftpair{t'}{S}{\ttime}{A}$, i.e.
\[
\rpshiftstate{t}{S}{\ttime}{A} \subseteq \rpshiftstate{t'}{S}{\ttime}{A}
\quad\text{and}\quad
\rpshifttimes{t}{S}{A}{\ttime} \ge \rpshifttimes{t'}{S}{A}{\ttime}.
\]

\item \label{one_type_shift_monotone:S_part}
The map $S\mapsto \rpshiftpair{t}{S}{\ttime}{A}$ is increasing: If $S\subseteq S'$, then $\rpshiftpair{t}{S}{\ttime}{A} \le \rpshiftpair{t}{S'}{\ttime}{A}$, i.e.
\[
\rpshiftstate{t}{S}{\ttime}{A} \subseteq \rpshiftstate{t}{S'}{\ttime}{A}
\quad\text{and}\quad
\rpshifttimes{t}{S}{A}{\ttime} \ge \rpshifttimes{t}{S'}{A}{\ttime}.
\]

\item \label{one_type_shift_monotone:pair_part}
The map $\pair{\ttime}{A} \mapsto \rpshiftpair{t}{S}{\ttime}{A}$ is increasing within $S$:
\[
\pair{\ttime}{A} \rpworse{S} \pair{\ttime'}{A'} \implies 
\rpshiftpair{t}{S}{\ttime}{A} \rple{S} \rpshiftpair{t}{S}{\ttime'}{A'}.
\]
This means that if
\[
\text{$A\cap S\subseteq A'\cap S$
\quad and\quad
$\restrict{\ttime}{\edges{S}\setminus \edges{A'}} \ge \restrict{\ttime'}{\edges{S}\setminus \edges{A'}}$,}
\]
then
\[
\rpshiftstate{t}{S}{\ttime}{A} \subseteq \rpshiftstate{t}{S}{\ttime'}{A'}
\quad\text{and}\quad
\restrict{\rpshifttimes{t}{S}{A}{\ttime}}{\edges{S}}
	\ge \restrict{\rpshifttimes{t}{S}{A'}{\ttime'}}{\edges{S}}.
\]
\end{enumerate}
\end{lem}

The next lemma shows that the locality properties from Section~\ref{fpc_basics:locality_sec} extend to the one-type shift operator.

\begin{lem}[Locality of one-type shifts]
\label{equiv_one_type_shifts_lem}
If $\ttime$ is any traversal measure and $A\subseteq S\subseteq R\subseteq\Zd$, the following statement can be added to the list of equivalent statements in Part~\ref{strong_restricted_proc:equiv_part} of Lemma~\ref{strong_restricted_proc_equiv_lem}:
\begin{equation}
\label{equiv_one_type_shifts_eqn}
\rpshiftpair{t'}{R}{\ttime}{A} = \rpshiftpair{t'}{S}{\ttime}{A}
\quad\text{for all $t'\le t$.}
\end{equation}
\end{lem}

\begin{proof}
Clearly \eqref{equiv_one_type_shifts_eqn} implies \eqref{strong_restricted_proc:process_times_part} from Lemma~\ref{strong_restricted_proc_equiv_lem}. On the other hand, if \eqref{strong_restricted_proc:process_times_part} holds, then combined with \eqref{strong_restricted_proc:star_process_part} we get
\[
\rreachedset{\ttime}{A}{R}{t'} = \rreachedset{\ttime}{A}{S}{t'}
\quad\text{and}\quad
\rreachedset{\ttime}{A}{\stargraph{R}}{t'} = \rreachedset{\ttime}{A}{\stargraph{S}}{t'}
\quad\text{for all $t'\le t$,}
\]
and combined with \eqref{strong_restricted_proc:vertex_times_part}, we get
\[
t'\le t\implies
\rptime{\ttime}{R}{A}{\vect u} = \rptime{\ttime}{S}{A}{\vect u}
\text{ for all $\vect v\in \rreachedset{\ttime}{A}{R}{t'} = \rreachedset{\ttime}{A}{S}{t'}$.}
\]
Thus, the formulas defining $\rpshiftpair{t'}{S}{\ttime}{A}$ in Definition~\ref{one_type_shift_def} remain unchanged if we replace $S$ with $R$, so \eqref{equiv_one_type_shifts_eqn} holds if the equivalent statements \eqref{strong_restricted_proc:process_times_part}, \eqref{strong_restricted_proc:star_process_part}, and \eqref{strong_restricted_proc:vertex_times_part} hold.
\end{proof}

\subsubsection{Pausing and Restarting the Two-Type Process}

\begin{definition}[Shifts on the domain of $\process{}$]
\label{two_type_shift_def}
Let $\domain\subseteq\domainspace$ be any domain of definition for the two-type process, and suppose $\pair{\outcome}{X}\in \domain$, so $\process{}$ is well-defined on the pair $\pair{\outcome}{X}$.
For $t\ge 0$, define $\shiftstate{t}{\outcome}{X}\in \statespace$ and $\shifttimes{t}{X}{\outcome} = \pair[2]{\shifttimes{t}{X}{\outcome}_1}{\shifttimes{t}{X}{\outcome}_2} \in \samplespace$ by
\[
\shiftstate{t}{\outcome}{X} \definedas \processtoutcome[2]{X}{t}{\outcome}
\in\statespace
\quad\text{and}\quad
\shifttimes{t}{X}{\outcome}_i \definedas
\rpshifttimes{t}{\finalset{i}}{\initialset{i}}{\outcome_i}\in \onetypespace,
\]
where $\initialset{i} = \stateset{i}{X}$ and $\finalset{i} = \reachedset{i}{X}{\infty}$ are species $i$'s initial and final sets, respectively, in the process started from $X$, and
$\rpshifttimes[10]{t}{\finalset{i}}{\initialset{i}}{}$ 
is the one-type traversal-shift operator from Definition~\ref{one_type_shift_def}. Then set
\[
\shiftpair{t}{\outcome}{X} \definedas \pair[2]{\shifttimes{t}{X}{\outcome}}{\shiftstate{t}{\outcome}{X}} \in \samplespace\times \statespace.
\]
%
\end{definition}

That is, $\shiftstate{t}{\outcome}{X}$ is just the state of the process started from $\pair{\outcome}{X}$ at time $t$, and $\shifttimes{t}{X}{\outcome}$ is the traversal time configuration obtained from $\outcome$ by zeroing out the traversal times within the conquered region at time $t$, and for each boundary edge of this conquered region, subtracting from its traversal time the length of time the infected endvertex of the edge has been in the infected state. Then $\shiftpair{t}{\outcome}{X}$ is simply this ``time-shifted" traversal/state pair. Note that the zeroing out of the traversal times within the conquered region is unimportant; this was done merely to emphasize the fact that the evolution of the process after time $t$ no longer depends on these times. The point of these definitions is Lemma~\ref{markovesque_lem} below.

\begin{lem}[Equivalent disentangled shifts]
\label{equiv_disentangled_shifts_lem}
\newcommand{\thispair}{\pair{\outcome}{X}}

Let $\thispair\in \goodsamplespace\times \statespace$, so $\outcome\in\samplespace$ satisfies \notiesdet\ and \mintimebddgrowth. For $i\in\setof{1,2}$ and $t\in [0,\infty]$, if $S$ is any subset of $\Zd$ satisfying
\[
\reachedset{i}{X}{t}_\outcome \subseteq S\subseteq
\Zd\setminus \reachedset{3-i}{X}{t-}_\outcome,
\]
then
\[
\shiftpair{t'}{\outcome}{X}_i = \rpshiftpair[2]{t'}{S}{\outcome_i}{\stateset{i}{X}}
\quad\text{for all $t'\le t$,}
\]
where $\shiftpair{t'}{\outcome}{X}_i \definedas
\pair[1]{\shifttimes{t'}{X}{\outcome}_i}{\,\stateset[2]{i}{\shiftstate{t'}{\outcome}{X}}}$.
\end{lem}

\begin{proof}
\newcommand{\fset}[1]{\finalset{#1}}
\newcommand{\fsett}[1]{\widehat{\finalsetsymb}_{#1}}
\newcommand{\iset}[1]{\initialset{#1}}
\newcommand{\isett}[1]{\initialset{#1}^t}
\newcommand{\Xt}{X_t}
\newcommand{\om}{\outcome}
\newcommand{\omi}[1]{\om_{#1}}
\newcommand{\omt}{\outcome_t}
\newcommand{\omti}[1]{\outcome_{#1}^t}
\newcommand{\thispair}{\pair{\om}{X}}
\newcommand{\thispairt}{\pair{\omt}{\Xt}}
\newcommand{\thispairi}[1]{\pair{\om_{#1}}{\iset{#1}}}
\newcommand{\thispairti}[1]{\pair[2]{\omti{#1}}{\isett{#1}}}


Unravelling the definitions, we have $\shiftpair{t'}{\outcome}{X}_i = \rpshiftpair[2]{t'}{\fset{i}}{\omi{i}}{\stateset{i}{X}}$, where the expression on the right is the one-type shift from Definition~\ref{one_type_shift_def}. Under the stated hypotheses, it follows directly from Lemma~\ref{one_type_two_type_lem} that
\[
\reachedset{i}{X}{t} = \rreachedset{i}{\iset{i}}{\fset{i}}{t}
= \rreachedset{i}{\iset{i}}{S}{t},
\]
and Lemma~\ref{equiv_one_type_shifts_lem} then implies that
\[
\rpshiftpair[2]{t'}{\fset{i}}{\omi{i}}{\stateset{i}{X}} = \rpshiftpair[2]{t'}{S}{\omi{i}}{\stateset{i}{X}}
\quad\text{for all $t'\le t$,}
\]
which proves the lemma.
\end{proof}

\begin{lem}[Markovesque property]
\label{markovesque_lem}
Let $\shiftop{\cdot}$ be the family of natural shift operators on the process space $\processspace$, defined in \eqref{natural_shiftop_def_eqn}, and let $\pshiftop{\cdot}$ be the family of shift operators on domains $\domain\subseteq\domainspace$ from Definition~\ref{two_type_shift_def}. Then on the domain $\startdomain\subset \domainspace$ defined in Section~\ref{fpc_basics:monotonicity_sec}, the first-passage competition process $\process{}$ satisfies
\[
\shiftop{t}\circ \process{} = \process{}\circ \pshiftop{t}
\quad\text{for all $t\ge 0$}.
\]
That is, if $\pair{\outcome}{X}\in \startdomain$, then $\process{}$ is well-defined on the pair $\pair{\outcome_t}{X_t} \definedas \shiftpair{t}{\outcome}{X}\in \domainspace$, and
\begin{equation}
\label{markovesque_eqn}
\processtoutcome[2]{X}{t+s}{\outcome}
= \processt{}{s} \circ \shiftpair{t}{\outcome}{X}
= \processtoutcome[2]{X_t}{s}{\outcome_t}
\quad\text{for all $t,s\ge 0$}.
\end{equation}
\end{lem}

\begin{proof}
\newcommand{\fset}[1]{\finalset{#1}}
\newcommand{\fsett}[1]{\widehat{\finalsetsymb}_{#1}}
\newcommand{\iset}[1]{\initialset{#1}}
\newcommand{\isett}[1]{\initialset{#1}^t}
\newcommand{\Xt}{X_t}
\newcommand{\om}{\outcome}
\newcommand{\omi}[1]{\om_{#1}}
\newcommand{\omt}{\outcome_t}
\newcommand{\omti}[1]{\outcome_{#1}^t}
\newcommand{\thispair}{\pair{\om}{X}}
\newcommand{\thispairt}{\pair{\omt}{\Xt}}
\newcommand{\thispairi}[1]{\pair{\om_{#1}}{\iset{#1}}}
\newcommand{\thispairti}[1]{\pair[2]{\omti{#1}}{\isett{#1}}}

Assume that $\thispair\in \goodsamplespace\times \statespace$; the proof for remaining elements of $\startdomain$ follows trivially from Lemmas~\ref{trivial_init_configs_lem} and \ref{one_type_markovesque_lem}. Fix $t\ge 0$ and let $\pair{\omt}{\Xt} = \shiftpair{t}{\om}{X}$. Our strategy will be to first show that $\reachedset{i}{X}{\infty}_\outcome = \reachedset{i}{\Xt}{\infty}_{\omt}$ for $i\in\setof{1,2}$ using Lemmas~\ref{one_type_markovesque_lem} and \ref{equiv_disentangled_shifts_lem}. It will then follow from equation \eqref{one_type_markovesque:process_shift_eqn} in Lemma~\ref{one_type_markovesque_lem} that the process started from $\thispairt$ is a time-shifted version of the process started from $\thispair$.

Let $\om = \pair{\omi{1}}{\omi{2}}$ and $\omt = \pair[2]{\omti{1}}{\omti{2}}$, let $X \equiv \pair{\iset{1}}{\iset{2}}$ and $\Xt \equiv \pair{\isett{1}}{\isett{2}}$, and let $\fset{i} = \reachedset{i}{X}{\infty}_\outcome$ and $\fsett{i} = \reachedset{i}{\Xt}{\infty}_{\omt}$ for $i\in\setof{1,2}$. Note that by definition we have $\thispairti{i} = \rpshiftpair{t}{\fset{i}}{\om_i}{\iset{i}}$ and $\isett{i} = \reachedset{i}{X}{t}_\om = \rreachedset{\om_i}{\iset{i}}{\fset{i}}{t}$. Also note that since $\om\in \goodsamplespace$, we have $\fset{1}\cap \fset{2} = \emptyset$ and $\fset{1} \cup \fset{2} = \Zd$ by Propositions~\ref{entangled_well_defined_prop} and \ref{entangled_full_prop}. We break the proof that $\fset{i}=\fsett{i}$ into two claims.

\begin{claim}
\label{markovesque:C1_smaller_claim}
$\fset{i}\subseteq\fsett{i}$ for $i\in\setof{1,2}$.
\end{claim}

\begin{proof}[Proof of Claim~\ref{markovesque:C1_smaller_claim}]
Take $i=1$ for concreteness.
We will show that $\fset{1}$ is a $\fsett{1}$-set, i.e.\ that $\fset{1}$ is one of the sets $S$ in the union defining $\fsett{1}$:
\[
\fsett{1} = \bigcup \setbuilder[1]{S\subseteq\Zd\setminus \isett{2}}
	{\rptime{\omti{1}}{S}{\isett{1}}{\vect v}
	<\starptimenew{\omti{2}}{(\Zd\setminus S)}{\isett{2}}{\vect v}
	\text{ for all } \vect v\in \bd S\setminus \isett{1}}.
\]
Since $\isett{2}\subseteq \fset{2}$, and $\fset{1}\cap\fset{2} = \emptyset$, we have $\fset{1}\cap \isett{2} = \emptyset$. Now let $\vect v\in \bd \fset{1}\setminus \isett{1}$. Then since $\thispairti{1} = \rpshiftpair{t}{\fset{1}}{\om_1}{\iset{1}}$ and $\vect v\not\in \isett{1} = \rreachedset{\om_1}{\iset{1}}{\fset{1}}{t}$, the time-shift equation \eqref{one_type_markovesque:rptime_shift_eqn} from Lemma~\ref{one_type_markovesque_lem} implies that
\begin{equation}
\label{markovesque:TC1_eqn}
\rptime{\omti{1}}{\fset{1}}{\isett{1}}{\vect v}
= \rptime{\om_1}{\fset{1}}{\iset{1}}{\vect v} - t.
\end{equation}
We need to show that the expression in \eqref{markovesque:TC1_eqn} is less than $\starptimenew{\omti{2}}{(\Zd\setminus \fset{1})}{\isett{2}}{\vect v}$. First note that since $\vect v\in \fset{1}\setminus \iset{1}$, Part~\ref{final_sets:characterization_part} of Lemma~\ref{final_sets_properties_lem} implies that
\begin{equation}
\label{markovesque:v_in_C1_eqn}
\rptime{\om_1}{\fset{1}}{\iset{1}}{\vect v} < \starptimenew{\om_2}{(\Zd\setminus \fset{1})}{\iset{2}}{\vect v}.
\end{equation}
Combined with the assumption that $\vect v\not\in \isett{1}$, \eqref{markovesque:v_in_C1_eqn} implies that $\vect v\not\in \rreachedset{\om_2}{\iset{2}}{\stargraph[0]{\Zd\setminus \fset{1}}}{t}$, because otherwise we'd have $\starptimenew{\om_2}{(\Zd\setminus \fset{1})}{\iset{2}}{\vect v}\le t < \rptime{\om_1}{\fset{1}}{\iset{1}}{\vect v}$, contradicting \eqref{markovesque:v_in_C1_eqn}. Now, since $\fset{2} = \Zd\setminus \fset{1}$, we have $\thispairti{2} =\rpshiftpair{t}{\fset{2}}{\om_2}{\iset{2}} = \rpshiftpair{t}{\Zd\setminus \fset{1}}{\om_2}{\iset{2}}$, and since $\vect v\not\in \rreachedset{\om_2}{\iset{2}}{\stargraph[0]{\Zd\setminus \fset{1}}}{t}$, the time-shift equation \eqref{one_type_markovesque:starptime_shift_eqn} from Lemma~\ref{one_type_markovesque_lem} then implies that
\begin{equation}
\label{markovesque:TC1_comp_eqn}
\starptimenew{\omti{2}}{(\Zd\setminus \fset{1})}{\isett{2}}{\vect v}
= \starptimenew{\om_2}{(\Zd\setminus \fset{1})}{\iset{2}}{\vect v}-t.
\end{equation}
Combining \eqref{markovesque:TC1_eqn}, \eqref{markovesque:v_in_C1_eqn}, and \eqref{markovesque:TC1_comp_eqn}, we get that for all $\vect v\in \bd \fset{1} \setminus \isett{1}$,
\begin{align*}
\rptime{\omti{1}}{\fset{1}}{\isett{1}}{\vect v}
\;= \;\rptime{\om_1}{\fset{1}}{\iset{1}}{\vect v} - t
\;<\; \starptimenew{\om_2}{(\Zd\setminus \fset{1})}{\iset{2}}{\vect v}-t
\;=\; \starptimenew{\omti{2}}{(\Zd\setminus \fset{1})}{\isett{2}}{\vect v},
\end{align*}
which shows that $\fset{1}$ is a $\fsett{1}$-set and hence $\fset{1}\subseteq \fsett{1}$. Switching all the 1's and 2's in the above argument shows that $\fset{2}$ is a $\fsett{2}$-set and hence $\fset{2}\subseteq \fsett{2}$.
\end{proof}

\begin{claim}
\label{markovesque:C1_bigger_claim}
$\fsett{i}\subseteq\fset{i}$ for $i\in\setof{1,2}$.
\end{claim}

\begin{proof}[Proof of Claim~\ref{markovesque:C1_bigger_claim}]
Take $i=1$ for concreteness. We will show that $\fsett{1}$ is a $\fset{1}$-set, i.e.\ that $\fsett{1}$ is one of the sets $S$ in the union defining $\fset{1}$:
\[
\fset{1} = \bigcup \setbuilder[1]{S\subseteq\Zd\setminus \iset{2}}
	{\rptime{\omi{1}}{S}{\iset{1}}{\vect v}
	<\starptimenew{\omi{2}}{(\Zd\setminus S)}{\iset{2}}{\vect v}
	\text{ for all } \vect v\in \bd S\setminus \iset{1}}.
\]
Since $\iset{2} \subseteq \isett{2}$ and $\isett{2}\cap \fsett{1} = \emptyset$ by definition, we have $\iset{2}\cap \fsett{1} = \emptyset$. Now let $\vect v\in \bd \fsett{1}\setminus \iset{1}$. First suppose $\vect v\in \isett{1}$. Since $\isett{1} = \rreachedset{\om_1}{\iset{1}}{\fset{1}}{t} \subseteq \fset{1}$, and $\fset{1}\subseteq \fsett{1}$ by Claim~\ref{markovesque:C1_smaller_claim}, Part~\ref{final_sets:characterization_part} of Lemma~\ref{final_sets_properties_lem} plus the monotonicity of $\ptimefamily[\om_1]$ and $\starptimefamily[\om_2]$ imply that
\begin{equation}
\label{markovesque:v_in_At_eqn}
\rptime{\omi{1}}{\fsett{1}}{\iset{1}}{\vect v}
\le \rptime{\omi{1}}{\fset{1}}{\iset{1}}{\vect v}
<\starptimenew{\omi{2}}{(\Zd\setminus \fset{1})}{\iset{2}}{\vect v}
\le \starptimenew{\omi{2}}{(\Zd\setminus \fsett{1})}{\iset{2}}{\vect v}
\quad \forall  \vect v\in  \isett{1}\setminus \iset{1}.
\end{equation}
Now suppose that $\vect v\not\in \isett{1}$. First note that by definition we have
\begin{equation*}
\reachedset{1}{X}{t}_\om = \isett{1}
\subseteq \fsett{1} \subseteq \Zd\setminus \isett{2}
= \Zd\setminus \reachedset{2}{X}{t}_\om,
\end{equation*}
and taking complements,
\begin{equation*}
\reachedset{2}{X}{t}_\om \subseteq \Zd\setminus \fsett{1}
\subseteq \Zd\setminus \reachedset{1}{X}{t}_\om.
\end{equation*}
Thus, since $\om\in \goodsamplespace$, Lemma~\ref{equiv_disentangled_shifts_lem} implies that
\begin{align}
\rpshiftpair{t}{\fsett{1}}{\om_1}{\iset{1}}
	&= \rpshiftpair{t}{\fset{1}}{\om_1}{\iset{1}} = \thispairti{1},
\quad\text{and}
\label{markovesque:C1_shift_eqn}\\
\rpshiftpair{t}{\Zd\setminus \fsett{1}}{\om_2}{\iset{2}}
	&=  \rpshiftpair{t}{\fset{2}}{\om_2}{\iset{2}} = \thispairti{2}.
\label{markovesque:C2_shift_eqn}
\end{align}
Now, \eqref{markovesque:C1_shift_eqn} implies that $\isett{1} = \rreachedset{\om_1}{\iset{1}}{\fsett{1}}{t}$, so since $\vect v\not\in \isett{1}$, the time-shift equation  \eqref{one_type_markovesque:rptime_shift_eqn} from Lemma~\ref{one_type_markovesque_lem} implies that
\begin{equation}
\label{markovesque:TC1t_eqn}
\rptime{\om_1}{\fsett{1}}{\iset{1}}{\vect v}
= \rptime{\omti{1}}{\fsett{1}}{\isett{1}}{\vect v} + t.
\end{equation}
Next note that since $\vect v\in \fsett{1}\setminus \isett{1}$, by the definition of $\fsett{1}$ and Part~\ref{final_sets:characterization_part} of Lemma~\ref{final_sets_properties_lem} we have
\begin{equation}
\label{markovesque:v_in_C1t_eqn}
\rptime{\omti{1}}{\fsett{1}}{\isett{1}}{\vect v} <
\starptimenew{\omti{2}}{(\Zd\setminus \fsett{1})}{\isett{2}}{\vect v}.
\end{equation}
We now claim that $\starptimenew{\om_2}{(\Zd\setminus \fsett{1})}{\iset{2}}{\vect v}>t$. For a contradiction, suppose  $\starptimenew{\om_2}{(\Zd\setminus \fsett{1})}{\iset{2}}{\vect v}\le t$, i.e.\ $\vect v\in \rreachedset{\om_2}{\iset{1}}{\stargraph[0]{\Zd\setminus \fsett{1}}}{t}$. Then it follows directly from \eqref{markovesque:C2_shift_eqn} and the construction of $\rpshiftop{t}{\Zd\setminus \fsett{1}}$ in Definition~\ref{one_type_shift_def} that $\starptimenew{\omti{2}}{(\Zd\setminus \fsett{1})}{\isett{2}}{\vect v}=0$, because either $\vect v\in \isett{2} = \rreachedset{\om_2}{\iset{1}}{\Zd\setminus \fsett{1}}{t}$, or $\vect v$ is connected to $\isett{2}$ by a single edge $\edge$ with $\omti{2}(\edge)=0$. On the other hand, since $\rptime{\om_1}{\fsett{1}}{\iset{1}}{\vect v}>t$ by assumption, \eqref{markovesque:TC1t_eqn} implies that $\rptime{\omti{1}}{\fsett{1}}{\isett{1}}{\vect v} =\rptime{\om_1}{\fsett{1}}{\iset{1}}{\vect v} - t >0$. This contradicts \eqref{markovesque:v_in_C1t_eqn}, so we conclude that we must have $\starptimenew{\om_2}{(\Zd\setminus \fsett{1})}{\iset{2}}{\vect v}>t$. Therefore, combined with \eqref{markovesque:C2_shift_eqn}, the time-shift equation \eqref{one_type_markovesque:starptime_shift_eqn} from Lemma~\ref{one_type_markovesque_lem} implies that
\begin{equation}
\label{markovesque:TC1t_comp_eqn}
\starptimenew{\om_2}{(\Zd\setminus \fsett{1})}{\iset{2}}{\vect v}
= \starptimenew{\omti{2}}{(\Zd\setminus \fsett{1})}{\isett{2}}{\vect v} + t.
\end{equation}
Combining \eqref{markovesque:TC1t_eqn}, \eqref{markovesque:v_in_C1t_eqn}, and \eqref{markovesque:TC1t_comp_eqn}, we get that for all $\vect v\in \bd \fsett{1} \setminus \isett{1}$,
\begin{align*}
\rptime{\om_1}{\fsett{1}}{\iset{1}}{\vect v}
\;= \;\rptime{\omti{1}}{\fsett{1}}{\isett{1}}{\vect v} + t
\;<\; \starptimenew{\omti{2}}{(\Zd\setminus \fsett{1})}{\isett{2}}{\vect v}+t
\;=\; \starptimenew{\om_2}{(\Zd\setminus \fsett{1})}{\iset{2}}{\vect v},
\end{align*}
and \eqref{markovesque:v_in_At_eqn} showed the same inequality for $\vect v\in \isett{1} \setminus \iset{1}$, so it holds for all $\vect v\in \bd \fsett{1} \setminus \iset{1}$. This shows that $\fsett{1}$ is a $\fset{1}$-set and hence $\fsett{1}\subseteq \fset{1}$. Switching all the 1's and 2's in the above argument shows that $\fsett{2}$ is a $\fset{2}$-set and hence $\fsett{2}\subseteq \fset{2}$.
\end{proof}

Combining Claims~\ref{markovesque:C1_smaller_claim} and \ref{markovesque:C1_bigger_claim}, we have $\fset{1}=\fsett{1}$ and $\fset{2}=\fsett{2}$. Therefore, by \eqref{one_type_markovesque:process_shift_eqn} in Lemma~\ref{one_type_markovesque_lem}, for $i\in\setof{1,2}$ we have
\begin{equation}
\label{markovesque:final_eqn}
\reachedset{i}{X}{t+s}_\om
=\rreachedset{\om_i}{\iset{i}}{\fset{i}}{t+s}
=\rreachedset{\omti{i}}{\isett{i}}{\fset{i}}{s}
=\rreachedset{\omti{i}}{\isett{i}}{\fsett{i}}{s}
=\reachedset{i}{\Xt}{s}_{\omt}
\quad \forall s\ge 0.
\end{equation}
Since $\process{}$ is well-defined on the pair $\thispair\in \startdomain$, the fact that \eqref{markovesque:final_eqn} holds implies that the sets $\reachedset{i}{\Xt}{s}_{\omt}$ are disjoint for all $s\ge 0$, which shows that $\process{}$ is also well-defined on $\thispairt$. Therefore, the expression $\processtoutcome[2]{\Xt}{s}{\omt}$ makes sense for all $t,s\ge 0$, and \eqref{markovesque_eqn} is equivalent to \eqref{markovesque:final_eqn}.
\end{proof}

The next result follows immediately by combining Lemmas~\ref{monotonicity_lem} and \ref{markovesque_lem}, and will be used in the proof of Theorem~\ref{irrelevance_thm} below.

\begin{cor}[Delayed domination]
\label{delayed_dom_cor}
If $\process{}$ is defined on the pairs $\pair{\outcome}{X}$ and $\pair{\outcome'}{X'}$, and if $\shiftpair{t_0}{\outcome}{X} \pbetter \shiftpair{t_0}{\outcome'}{X'}$ for some $t_0\ge 0$, then
$\processtoutcome[2]{X}{t}{\outcome} \ge \processtoutcome[2]{X'}{t}{\outcome'}$ for all $t\ge t_0$.
\end{cor}



The following is an alternate definition of a shift operator for the two-type process, defined directly in terms of the process rather than in terms of the one-type shifts.

\begin{definition}[Alternate two-type shift]
\label{direct_shift_def}
Suppose $\process{}$ is defined on the pair $\pair{\outcome}{X}\in \samplespace\times \statespace$. For $t\ge 0$, define $\ashiftstate{t}{\outcome}{X}\in \statespace$ and $\ashifttimes{t}{X}{\outcome} = \pair[2]{\ashifttimes{t}{X}{\outcome}_1}{\ashifttimes{t}{X}{\outcome}_2} \in \samplespace$ by
\[
\ashiftstate{t}{\outcome}{X} \definedas \processtoutcome[2]{X}{t}{\outcome},
\]
and for $\edge\in \edges{\Zd}$ and $i\in \setof{1,2}$,
\[
\ashifttimes{t}{X}{\outcome}_i(\edge) \definedas
\begin{cases}
0 & \text{if } e\in \edges[2]{\reachedset{i}{X}{t}_\outcome},\\
\outcome_i(\edge) - \parens[2]{t-\entangledtime{i}{X}{\vect u}_\outcome} &
	\text{if $\edge$ joins $\vect u\in \reachedset{i}{X}{t}_\outcome$
	to $\vect v\notin \bothreachedset{X}{t}_\outcome$},\\
\outcome_i(\edge) & \text{otherwise}.
\end{cases}
\]
Then set $\ashiftpair{t}{\outcome}{X} \definedas \pair[2]{\ashifttimes{t}{X}{\outcome}}{\ashiftstate{t}{\outcome}{X}} \in \samplespace\times \statespace$.
\end{definition}


The following lemma shows that for nice input pairs, the shift in Definition~\ref{direct_shift_def} is equivalent to the one in Definition~\ref{two_type_shift_def}.

\begin{lem}
\label{alternate_shift_lem}
If $\pair{\outcome}{X}\in \startdomain$, then $\shiftpair{t}{\outcome}{X} \pequiv \ashiftpair{t}{\outcome}{X}$ as defined in \eqref{pair_equiv_eqn}. Therefore, $\pshiftop{t}$ can be replaced with $\apshiftop{t}$ in Lemma~\ref{markovesque_lem} and Corollary~\ref{delayed_dom_cor}.
\end{lem}

\begin{proof}
The proof consists of checking definitions and is left to the reader.
\end{proof}

%

\section{Survival When One Initial Set Is Unbounded}
\label{fpc_basics:unbounded_survival_sec}

%
%

In this section we are interested in the case where one species starts on an infinite set while the other species starts on a finite set. The main result is Theorem~\ref{irrelevance_thm} below, which says that if $A_1$ is finite and $A_2$ is infinite, the possibility of survival for species 1 depends only on the component of $\Zd\setminus A_2$ in which species 1 starts, and not on the shape or location of its starting set $A_1$.
This is analogous to the result of \cite{\DHa} (and also \cite[Lemma~5.1]{\GM}), which says that when both $A_1$ and $A_2$ are finite, and neither species is initially surrounded by the other, the possibility of coexistence doesn't depend on the initial configuration, but only on the distribution of the traversal measure $\tmeasurepair$. Note that the present result, with $A_2$ infinite, is in one respect easier to prove than the corresponding theorems in \cite{\DHa} and \cite{\GM}, because the nonmonotone event that both species survive from finite initial sets is reduced to the monotone event that one species survives when the initial set of the other species is already infinite.
On the other hand, allowing the set $A_2$ to be infinite introduces the additional complication that species~2 can conquer infinitely many vertices in an arbitrarily small time interval. We will use a conditioning argument to control this potentially infinitely fast growth of species~2.

In fact, in Theorem~\ref{irrelevance_thm} we prove the stronger statement that for any fixed finite or infinite starting set $A_2$, if species 1 is able to conquer a given set $S\subseteq\Zd$ with positive probability from $A_1$, then it is also able to conquer $S$ with positive probability from $A_1'$ when $A_1'$ is in the same component of $\Zd\setminus A_2$ as $A_1$. This version of the statement will be needed for some of our results in Chapters~\ref{random_fpc_chap} and \ref{coex_finite_chap}.
The conditioning argument we use in the proof of Theorem~\ref{irrelevance_thm} is adapted from a similar argument in the proof of Proposition~3.1 in \cite{Deijfen:2007aa}. Recall from Section~\ref{fpc_basics:additional_notation_sec} that $\processlaw{\tmeasurepair}{\initconfig}$ denotes the law of the two-type process started from $\initconfig$ and using the random traversal measure $\tmeasurepair$.


\begin{thm}[Irrelevance of starting configuration]
\label{irrelevance_thm}
\newcommand{\startset}{\initialset{1}}
\newcommand{\newstartset}{\startset'}
\newcommand{\twoset}{\initialset{2}}
\newcommand{\startpair}{\pair{\startset}{\twoset}}
\newcommand{\newstartpair}{\pair{\newstartset}{\twoset}}
\newcommand{\targetset}{S}

Suppose $\tmeasurepair = \pair{\tmeasure_1}{\tmeasure_2}$ is an \iid\ random two-type traversal measure on $\edges{\Zd}$ which satisfies \minimalmoment{d} and \notiesprob, and also
\begin{enumerate}
\item $0\in \supp \cmttime{1}$.
\item Neither $\cmttime{1}$ nor $\cmttime{2}$ has an atom at 0.
\end{enumerate}
Let $\twoset \subset \Zd$, and let $\startset$ and $\newstartset$ be two finite connected subsets of $\Zd$ which lie in the same component of $\Zd\setminus \twoset$. Then if $\targetset$ is any subset of $\Zd$,
\[
\processlaw{\tmeasurepair}{\startpair}
\eventthat[3]{\text{\rm species 1 conquers $\targetset$}}>0
\iff \processlaw{\tmeasurepair}{\newstartpair}
\eventthat[3]{\text{\rm species 1 conquers $\targetset$}}>0.
\]
In particular, if $\twoset$ is a fixed infinite set, then species 1 has a positive probability of surviving from the initial configuration $\newstartpair$ if and only if it has a positive probability of surviving from $\startpair$.
\end{thm}

\begin{proof}
\newcommand{\targetset}{S}

\newcommand{\startset}{A_1}
\newcommand{\startsetp}{A_1'}
\newcommand{\startsett}{\widetilde{A}_1}

\newcommand{\startpair}{\pair[2]{\startset}{A_2}}
\newcommand{\startpairp}{\pair[2]{\startsetp}{A_2}}
\newcommand{\startpairt}{\pair[2]{\startsett}{A_2}}

\newcommand{\startstate}{X}
\newcommand{\startstatep}{X'}
\newcommand{\startstatet}{\widetilde{X}}

\newcommand{\neighbors}{\widetilde N_1}

\newcommand{\survives}{G_1}
\newcommand{\survivesp}{G_1'}
\newcommand{\survivest}{\widetilde{G}_1}

\newcommand{\holdingevent}[2]{\widetilde H_{#1} \parens{#2}}
\newcommand{\fastevent}[1]{F_1 \parens{#1}}

\newcommand{\thiseventsymb}{E}
\newcommand{\intersectionevent}{\symboleqref{\widetilde{\thiseventsymb}}{irrelevance:main_event_eqn}}
\newcommand{\conditionalevent}{\thiseventsymb_{\sigfield}^+}
\newcommand{\finalevent}{\symboleqref{\thiseventsymb'}{irrelevance:pos_prob_eqn}}
\newcommand{\ttimesevent}[1]{\widehat{\thiseventsymb}_{#1}}

\newcommand{\sigfield}{\mathcal{F}}
\newcommand{\edgeset}{\Lambda} 
\newcommand{\compedgeset}{\edgeset^\complement}

\newcommand{\outcomeh}{\widehat{\outcome}}
\newcommand{\outcomep}{\outcome'}
\newcommand{\outcomet}{\widetilde{\outcome}}

Since $\startset$ and $\startsetp$ are both finite and connected and lie in the same component of $\Zd\setminus A_2$, there is some finite path $\gamma\subseteq \Z_d\setminus A_2$ connecting $\startset$ and $\startsetp$, and hence the subgraph $\startsett \definedas \startset\cup \gamma\cup \startsetp$ is finite and connected. We will consider the process started from the following three initial states:
\[
\startstate\equiv\startpair,
\quad
\startstatep\equiv\startpairp,
\quad\text{and}\quad
\startstatet\equiv\startpairt.
\]
Let $\survives$, $\survivesp$, and $\survivest$ denote the event that species 1 conquers $S$ from the initial state $\startstate$, $\startstatep$, or $\startstatet$, respectively. More explicitly,
\[
\survives = \eventthat[3]{ \reachedset{1}{\startstate}{\infty} \supseteq \targetset},
\quad
\survivesp = \eventthat[3]{ \reachedset{1}{\startstatep}{\infty} \supseteq \targetset}
\quad\text{and}\quad
\survivest = \eventthat[3]{ \reachedset{1}{\startstatet}{\infty} \supseteq \targetset}.
\]
Note that the above definitions only make sense for outcomes where the process is well-defined for the given starting configuration; by convention, we assume that all events are subsets of the $\Pr$-a.s.\ subset of $\samplespace$ on which \notiesdet\ and \mintimebddgrowth\ hold, so that the entangled process is defined for all starting configurations.
We will show that $\Pr(\survives)>0 \implies \Pr (\survivesp)>0$, whence the reverse implication follows by symmetry.

Suppose $\Pr(\survives)>0$. Then since $\startstatet \ge \startstate$, Lemma~\ref{increasing_sets_lem} implies that $\Pr\parens[2]{\survivest}>0$, because the set of final states in which species 1 has conquered $\targetset$ is an increasing subset of $\statespace$.
%
%
Now let $\neighbors \definedas \nbrhood[2]{\startsett} \setminus A_2$, where $\nbrhood[2]{\startsett}$ denotes the graph neighborhood of $\startsett$ in $\Zd$, and for $t\ge 0$, define events
\[
\holdingevent{1}{t} \definedas \eventthat[1]{\ptimenew[2]{1}{\startsett}{\,\Zd\setminus \startsett}>t}
\quad\text{and}\quad
\holdingevent{2}{t} \definedas \eventthat[1]{\ptimenew[2]{2}{A_2}{\neighbors}>t},
\]
where $\ptmetric[1]$ and $\ptmetric[2]$ are the disentangled passage times for the two-type process corresponding to $\tmeasurepair$. Since neither $\cmttime{1}$ nor $\cmttime{2}$ has an atom at 0, and $\startsett$ and $\neighbors$ are finite, we have $\Pr \parens[2]{\holdingevent{1}{0} \cap \holdingevent{2}{0}} = 1$.
Thus, since $\Pr\parens[2]{\survivest}>0$,  we can choose some sufficiently small (deterministic) $t_0>0$ such that
\begin{equation}
\label{irrelevance:main_event_eqn}
\Pr \parens[1]{\survivest\cap \holdingevent{1}{t_0} \cap \holdingevent{2}{t_0}} >0.
\end{equation}
Let $\intersectionevent$ be the event in \eqref{irrelevance:main_event_eqn}, and let $\edgeset \definedas \edges{\neighbors}$, the induced edge set of $\neighbors$. Consider the $\sigma$-field
$
\sigfield = \sigma \setof[2]{\ttimepair(\edge)}_{\edge\notin \edgeset},
$
generated by all the traversal times outside the edge set of $\neighbors$, and define the event
\[
\conditionalevent \definedas \eventthat[3]{ \cprob[2]{\intersectionevent}{\sigfield} >0 }.
\]
Note that $\Pr \parens[2]{\conditionalevent}>0$ since $\Pr \parens[2]{\intersectionevent}>0$ (cf.\ Lemma~\ref{zero_cond_prob_lem}), and clearly $\conditionalevent\in \sigfield$. We pause briefly to describe in detail what the event $\conditionalevent$ means.

Intuitively, $\conditionalevent$ is the event that given the traversal times of edges in $\compedgeset$, there is a positive probability that when we choose the remaining traversal times at random according to the marginal distribution of $\Pr$ on $\edgeset$, the event $\intersectionevent$ occurs. To describe this more explicitly, assume we are working on the canonical sample space $\samplespace = \samplespace[1]$, so that $\ttimepair_\outcome = \outcome$ for any outcome $\outcome\in\samplespace$. It can be verified using Fubini's theorem and the definition of conditional expectation that for $\Pr$-almost every $\outcome\in\samplespace$,
\begin{equation}\label{irrelevance:cond_prob_eqn}
\cprob[2]{\intersectionevent}{\sigfield}(\outcome)=
\Pr \setbuilder[3]{\outcomeh \in \samplespace\:}{\:\restrict{\outcomeh}{\edgeset} \uplus \restrict{\outcome}{\compedgeset}\in \intersectionevent},
\end{equation}
where
for $f\in \rsamplespace{\edgeset}$
and
$g\in \rsamplespace{\compedgeset}$, we write $f\uplus g$ for the ``concatenated" function in $ \samplespace[0] = \parens{\R_+\times\R_+}^{\edgeset \sqcup \compedgeset}$ defined by
\[
f\uplus g(\edge) \definedas
\begin{cases}
f(\edge) & \text{if } \edge\in \edgeset,\\
g(\edge) & \text{if } \edge\in \compedgeset.
\end{cases}
\]
That is, \eqref{irrelevance:cond_prob_eqn} says that for a given traversal time configuration $\outcome\in\samplespace$, the conditional probability $\cprob[2]{\intersectionevent}{\sigfield}(\outcome)$ is the probability that if we choose another $\outcomeh\in\samplespace$ at random and replace $\outcome$'s traversal times on $\edgeset$ with those of $\outcomeh$, but keep $\outcome$'s traversal times for edges outside $\edgeset$, the resulting new configuration is in the event $\intersectionevent$. The event $\conditionalevent$ is then the set of $\outcome\in \samplespace$ such that the conditional probability in \eqref{irrelevance:cond_prob_eqn} is positive.

The point of conditioning the event in $\intersectionevent$ on the $\sigma$-field $\sigfield$ is that we can then independently control the traversal times of species~1 in the set $\edgeset$, as follows. For $t> 0$, define the event
\[
\fastevent{t} \definedas \eventthat[1]{\ttime_1(\edge) \le \frac{t}{\abs{\edgeset}}
\text{ for all } \edge\in \edgeset}.
\]
Then $\Pr \parens[2]{\fastevent{t}} >0$ for all $t>0$ since $\tmeasurepair$ is \iid, $0\in \supp \cmttime{1}$, and $\abs{\edgeset}<\infty$. Moreover, $\fastevent{t}$ is independent of $\sigfield$ since $\fastevent{t}$ only involves traversal times of edges in $\edgeset$. Therefore, since $\conditionalevent\in \sigfield$ we have
\begin{equation}
\label{irrelevance:pos_prob_eqn}
\Pr \parens[1]{\fastevent{t_0} \cap \conditionalevent}
= \Pr \parens[2]{\fastevent{t_0}} \Pr \parens[2]{\conditionalevent}
>0.
\end{equation}
Let $\finalevent \definedas \fastevent{t_0} \cap \conditionalevent$. We will show that $\survivesp$ occurs almost surely on the event $\finalevent$, which will finish the proof since $\Pr \parens[2]{\finalevent}>0$ by \eqref{irrelevance:pos_prob_eqn}. To show this, we make the following claim.

\begin{claim}
\label{irrelevance:domination_claim}
For $\Pr$-a.e.\ $\outcomep\in \finalevent$, there exists $\outcomet\in \intersectionevent$ such that
$
\ashiftpair[2]{t_0}{\outcomep}{\startstatep}\pbetter \ashiftpair[2]{t_0}{\outcomet}{\startstatet},
$
where $\apshiftop{t_0}$ is the shift operator in Definition~\ref{direct_shift_def} and the preorder $\pbetter$ is defined in \eqref{pair_preorder_eqn}.
\end{claim}

Assuming Claim~\ref{irrelevance:domination_claim} is true, given $\outcomep$ in the full-measure subset of $\finalevent$ where the claim holds, choose some $\outcomet\in \intersectionevent$ satisfying the inequality in the claim. Then Corollary~\ref{delayed_dom_cor} and Lemma~\ref{alternate_shift_lem} imply that $\processtoutcome[2]{\startstatep}{t}{\outcomep}\ge \processtoutcome[2]{\startstatet}{t}{\outcomet}$ for all $t\ge t_0$. In particular, $\processtoutcome[2]{\startstatep}{\infty}{\outcomep}\ge \processtoutcome[2]{\startstatet}{\infty}{\outcomet}$, which implies that $\reachedset{1}{\startstatep}{\infty}_{\outcomep} \supseteq \reachedset{1}{\startstatet}{\infty}_{\outcomet} \supseteq \targetset$ since $\outcomet\in \intersectionevent \subseteq \survivest$. Therefore $\outcomep\in \survivesp$, and since $\outcomep\in \finalevent$ was ``almost arbitrary," we have $\finalevent\assubseteq \survivesp$ as claimed.  Hence, $\Pr \parens[2]{\survivesp}\ge \Pr \parens[2]{\finalevent} >0$ by \eqref{irrelevance:pos_prob_eqn}.


To complete the proof, it remains only to prove Claim~\ref{irrelevance:domination_claim}.

\begin{proof}[Proof of Claim~\ref{irrelevance:domination_claim}]
\newcommand{\sstatep}{\startstatep_{t_0}}
\newcommand{\sstatet}{\startstatet_{t_0}}

\newcommand{\somp}{\outcomep_{t_0}}
\newcommand{\somt}{\outcomet_{t_0}}

\newcommand{\ssetp}[1]{A_{#1}^{\prime t_0}}
\newcommand{\ssett}[1]{\widetilde{A}_{#1}^{t_0}}

\newcommand{\ompi}[1]{\outcomep_{#1}}
\newcommand{\omti}[1]{\outcomet_{#1}}

\newcommand{\sompi}[1]{\outcomep^{t_0}_{#1}}
\newcommand{\somti}[1]{\outcomet^{t_0}_{#1}}

By \eqref{irrelevance:cond_prob_eqn} and the definition of $\conditionalevent$, for $\Pr$-a.e.\ $\outcomep\in \conditionalevent$, the event
\[
\ttimesevent{\outcomep}\definedas \setbuilder[3]{\outcomeh \in \samplespace\:}{\:\restrict{\outcomeh}{\edgeset} \uplus \restrict{\outcomep}{\compedgeset}\in \intersectionevent}
\]
has positive probability, and in particular is nonempty. Therefore, for $\Pr$-a.e.\ $\outcomep\in \finalevent \subseteq \conditionalevent$, we can define $\outcomet = \outcomet \parens[2]{\outcomep}$ by choosing any $\outcomeh\in \ttimesevent{\outcomep}$ and setting $\outcomet \definedas \restrict{\outcomeh}{\edgeset} \uplus \restrict{\outcomep}{\compedgeset}$. Then $\outcomet \in \intersectionevent$ by construction, and it remains to show that $\ashiftpair[2]{t_0}{\outcomep}{\startstatep}\pbetter \ashiftpair[2]{t_0}{\outcomet}{\startstatet}$.



Let $\pair[2]{\somp}{\sstatep} \definedas \ashiftpair[2]{t_0}{\outcomep}{\startstatep}$ and  $\pair[2]{\somt}{\sstatet} \definedas \ashiftpair[2]{t_0}{\outcomet}{\startstatet}$, and let
\[
\somp = \pair[2]{\sompi{1}}{\sompi{2}}
\quad\text{and}\quad
\somt = \pair[2]{\somti{1}}{\somti{2}},
\]
and
\[
\sstatep \equiv \pair[2]{\ssetp{1}}{\ssetp{2}}
\quad\text{and}\quad
\sstatet \equiv \pair[2]{\ssett{1}}{\ssett{2}}.
\]
Note that $\ssetp{i} = \reachedset{i}{\startstatep}{t_0}_{\outcomep}$ and $\ssett{i} = \reachedset{i}{\startstatet}{t_0}_{\outcomet}$ for $i\in \setof{1,2}$. Also note that by construction we have
\begin{equation}
\label{irrelevance:outcomes_agree_eqn}
\restrict{\outcomep}{\compedgeset} = \restrict{\outcomet}{\compedgeset}.
\end{equation}
This fact will be the basis of several of the arguments below. We break the proof of the desired inequality into two claims, one for the states, and one for the traversal times.

\begin{claim}
\label{irrelevance:state_comparison_claim}
$\sstatep\ge \sstatet$. That is, $\ssetp{1} \supseteq \ssett{1}$ and $\ssetp{2} \subseteq \ssett{2}$.
\end{claim}

\begin{proof}[Proof of Claim~\ref{irrelevance:state_comparison_claim}]
Our goal will be to prove the following pair of equations, from which Claim~\ref{irrelevance:state_comparison_claim} follows directly:
\begin{gather}
\label{irrelevance:sp1_comparison_eqn} 
\reachedset{1}{\startstatep}{t_0}_{\outcomep}
\overset{\eqref{irrelevance:nbrs_covered_eqn}}{\supseteq}
	\neighbors \supseteq \startsett
\overset{\eqref{irrelevance:sp1_unmoved_eqn}}{=}
	\reachedset{1}{\startstatet}{t_0}_{\outcomet},\\
\label{irrelevance:sp2_comparison_eqn}
\reachedset{2}{\startstatep}{t}_{\outcomep}
\overset{\eqref{irrelevance:sp2_contained_eqn}}{\subseteq}
	\rreachedset{\ompi{2}}{A_2}{\Zd\setminus \neighbors}{t}
\overset{\eqref{irrelevance:sp2_growth_same_eqn}}{=}
	\rreachedset{\omti{2}}{A_2}{\Zd\setminus \neighbors}{t}
\overset{\eqref{irrelevance:sp2_equiv_eqn}}{=}
	\reachedset{2}{\startstatet}{t}_{\outcomet}
\quad\forall t\le t_0.
\end{gather}
The labels above the relational symbols indicate the equation in which the relation is proved below. For the present proof, we only need \eqref{irrelevance:sp2_comparison_eqn} to hold with $t=t_0$, but the more general version will be needed in the subsequent proof of Claim~\ref{irrelevance:ttimes_comparison_claim}. We now proceed to prove the relations in \eqref{irrelevance:sp1_comparison_eqn} and \eqref{irrelevance:sp2_comparison_eqn}.

First, since $\startsetp \subseteq \neighbors$ and $\outcomep\in \fastevent{t_0}$, we have
\begin{equation}
\label{irrelevance:nbrs_reached_fast_eqn}
\sup_{\vect v\in \neighbors} \rptime{\ompi{1}}{\neighbors}{\startsetp}{\vect v}
\le \sum_{\edge\in \edges{\neighbors}} \ompi{1}(\edge)
\le \sum_{\edge\in \edgeset} \frac{t_0}{\abs{\edgeset}}
= t_0,
\end{equation}
or equivalently, $\rreachedset{\ompi{1}}{\startsetp}{\neighbors}{t_0} = \neighbors$ (this is a special case of the ``total traversal measure bound for covering times" in Lemma~\ref{covering_time_properties_lem}). Next, since $\compedgeset = \staredges{\Zd\setminus \neighbors}$, we have
$
\restrict{\ompi{2}}{\staredges{\Zd\setminus \neighbors}}
= \restrict{\omti{2}}{\staredges{\Zd\setminus \neighbors}}
$
by \eqref{irrelevance:outcomes_agree_eqn}, and Lemma~\ref{ptime_tmeasure_monotone_lem} then implies that
\begin{equation}
\label{irrelevance:star_metrics_eqn}
\starptmetric[\ompi{2}]{(\Zd\setminus \neighbors)}
= \starptmetric[\omti{2}]{(\Zd\setminus \neighbors)}.
\end{equation}
In particular, since $\outcomet\in \holdingevent{2}{t_0}$, \eqref{irrelevance:star_metrics_eqn} implies that
\begin{equation}
\label{irrelevance:nbrs_not_reached_eqn}
t_0 < \ptimenew{\omti{2}}{A_2}{\neighbors}
\le \starptimenew{\omti{2}}{(\Zd\setminus \neighbors)}{A_2}{\neighbors}
=\starptimenew{\ompi{2}}{(\Zd\setminus \neighbors)}{A_2}{\neighbors}.
\end{equation}
Combining \eqref{irrelevance:nbrs_not_reached_eqn} with \eqref{irrelevance:nbrs_reached_fast_eqn} we have
\[
\rptime{\ompi{1}}{\neighbors}{\startsetp}{\vect v}\le t_0
< \starptimenew{\ompi{2}}{(\Zd\setminus \neighbors)}{A_2}{\vect v}
\quad\text{for all } \vect v\in \neighbors,
\]
and thus species 1 conquers $\neighbors$ in the process $\processoutcome[2]{\startstatep}{\outcomep}$, by Proposition~\ref{conquering_property_prop}. Therefore, by Lemma~\ref{conquered_set_growth_lem} we have
\begin{equation}
\label{irrelevance:nbrs_covered_eqn}
\reachedset{1}{\startstatep}{t_0}_{\outcomep}
\supseteq \rreachedset{\ompi{1}}{\startsetp}{\neighbors}{t_0} = \neighbors.
\end{equation}
Taking complements in \eqref{irrelevance:nbrs_covered_eqn}, we have $\Zd\setminus \reachedset{1}{\startstatep}{t_0}_{\outcomep} \subseteq \Zd\setminus \neighbors$, and since $\outcomep\in \goodsamplespace$, Corollary~\ref{entangled_disentangled_cor} implies that
\begin{equation}
\label{irrelevance:sp2_contained_eqn}
\reachedset{2}{\startstatep}{t}_{\outcomep}
\subseteq \rreachedset{\ompi{2}}{A_2}{\Zd\setminus \neighbors}{t}
\quad\forall t\le t_0.
\end{equation}
Next, since $\outcomet\in \holdingevent{1}{t_0}$, we have $\ptimenew{\omti{1}}{\startsett}{\Zd\setminus \startsett} >t_0$. Lemmas~\ref{entangled_disentangled_comparison_lem} and \ref{strong_restricted_proc_equiv_lem} then imply that
$
\reachedset{1}{\startstatet}{t_0}_{\outcomet}
\subseteq \reachedset{\omti{1}}{\startsett}{t_0} \subseteq \startsett
$,
and since the reverse inclusion is trivial,
\begin{equation}
\label{irrelevance:sp1_unmoved_eqn}
\reachedset{1}{\startstatet}{t_0}_{\outcomet} = \startsett.
\end{equation}
Similarly, since $\outcomet\in \holdingevent{2}{t_0}$, we have $\ptimenew{\omti{2}}{A_2}{\neighbors} >t_0$, and then by Lemmas~\ref{entangled_disentangled_comparison_lem} and \ref{strong_restricted_proc_equiv_lem} we have
\begin{equation}
\label{irrelevance:sp2_separated_eqn}
\reachedset{2}{\startstatet}{t_0}_{\outcomet}
\subseteq \reachedset{\omti{2}}{A_2}{t_0}
\subseteq \Zd\setminus \neighbors.
\end{equation}
Since $\startsett\subseteq \neighbors$, combining \eqref{irrelevance:sp1_unmoved_eqn} and \eqref{irrelevance:sp2_separated_eqn} we get
$
\reachedset{2}{\startstatet}{t_0}_{\outcomet}
\subseteq \Zd\setminus \neighbors \subseteq
\Zd\setminus \reachedset{1}{\startstatet}{t_0}_{\outcomet},
$
and thus, since $\outcomet\in \goodsamplespace$, Lemma~\ref{one_type_two_type_lem} then implies that
\begin{equation}
\label{irrelevance:sp2_equiv_eqn}
\reachedset{2}{\startstatet}{t}_{\outcomet}
=\rreachedset{\omti{2}}{A_2}{\Zd\setminus \neighbors}{t}
\quad\forall t\le t_0.
\end{equation}
Now, observe that by \eqref{irrelevance:star_metrics_eqn} we have $\rptime{\ompi{2}}{\Zd\setminus \neighbors}{A_2}{\vect v} = \rptime{\omti{2}}{\Zd\setminus \neighbors}{A_2}{\vect v}$ for all $\vect v\in \Zd\setminus \neighbors$, which implies that
\begin{equation}
\label{irrelevance:sp2_growth_same_eqn}
\rreachedset{\ompi{2}}{A_2}{\Zd\setminus \neighbors}{t}
=\rreachedset{\omti{2}}{A_2}{\Zd\setminus \neighbors}{t}
\quad\forall t\in [0,\infty].
\end{equation}
Finally, combining \eqref{irrelevance:nbrs_covered_eqn} and \eqref{irrelevance:sp1_unmoved_eqn}, we get \eqref{irrelevance:sp1_comparison_eqn}, and combining \eqref{irrelevance:sp2_contained_eqn}, \eqref{irrelevance:sp2_growth_same_eqn} and \eqref{irrelevance:sp2_equiv_eqn}, we get \eqref{irrelevance:sp2_comparison_eqn}. This completes the proof of Claim~\ref{irrelevance:state_comparison_claim}.
\end{proof}


\begin{claim}
\label{irrelevance:ttimes_comparison_claim}
$\displaystyle \restrict{\somp}{\compedges[0]{\ssetp{1}\cup \ssett{2}}} \ge \restrict{\somt}{\compedges[0]{\ssetp{1}\cup \ssett{2}}}$. That is, for all $\edge\not\in \edges[2]{\ssetp{1}\cup \ssett{2}}$,
\begin{equation}
\label{irrelevance:ttimes_comparison_eqn}
\sompi{1}(\edge)\le \somti{1}(\edge)
\quad\text{and}\quad
\sompi{2}(\edge)\ge \somti{2}(\edge).
\end{equation}
\end{claim}

\begin{proof}[Proof of Claim~\ref{irrelevance:ttimes_comparison_claim}]
The idea behind the proof is that the relations in Claim~\ref{irrelevance:state_comparison_claim} imply that for most of the relevant edges, at most one of the traversal times $\ompi{i}$ and $\omti{i}$ will be affected by the shift $\apshiftop{t_0}$ (for fixed $i\in\setof{1,2}$), and because we chose the two outcomes so that \eqref{irrelevance:outcomes_agree_eqn} holds, the inequalities in \eqref{irrelevance:ttimes_comparison_eqn} will follow trivially. The only edges on which $\apshiftop{t_0}$ affects both $\outcomep$ and $\outcomet$ for the same species are those in $\bdedges[2]{\ssetp{2}}\cap \bdedges[2]{\ssett{2}}$; for these edges we will use \eqref{irrelevance:sp2_comparison_eqn} to obtain the desired inequality for species 2's traversal times. We now proceed with the proof.

First note that by \eqref{irrelevance:nbrs_covered_eqn}, we have $\neighbors \subseteq \ssetp{1}$, and hence
\[
\edgeset = \edges{\neighbors} \subseteq \edges{\ssetp{1}}
\subseteq \edges[2]{\ssetp{1}\cup \ssett{2}}.
\]
Taking complements, $\compedges[2]{\ssetp{1}\cup \ssett{2}} \subseteq \compedgeset$, and therefore \eqref{irrelevance:outcomes_agree_eqn} implies that
\begin{equation}
\label{irrelevance:unconquered_edge_times_agree_eqn}
\outcomep(\edge) = \outcomet(\edge)
\quad\text{for all } e\in \compedges[2]{\ssetp{1}\cup \ssett{2}}.
\end{equation}
We will use \eqref{irrelevance:unconquered_edge_times_agree_eqn} several times throughout the rest of the proof.

Our goal is to verify \eqref{irrelevance:ttimes_comparison_eqn} for all $\edge\in \compedges[2]{\ssetp{1}\cup \ssett{2}}$, so choose some such edge $\edge$.
Then either $\edge\in \bdedges[2]{\ssetp{1}\cup \ssett{2}}$ or $\edge\not\in \staredges[2]{\ssetp{1}\cup \ssett{2}}$. First suppose $\edge\not\in \staredges[2]{\ssetp{1}\cup \ssett{2}}$. Then, using Claim~\ref{irrelevance:state_comparison_claim} and Lemma~\ref{edge_sets_lem}, we have $\edge\not\in \staredges[2]{\reachedset{i}{\startstatep}{t_0}_{\outcomep}}$ and $\edge\not\in \staredges[2]{\reachedset{i}{\startstatet}{t_0}_{\outcomet}}$ for $i\in\setof{1,2}$,
so it follows directly from Definition~\ref{direct_shift_def} that $\sompi{i}(\edge) = \ompi{i}(\edge)$ and $\somti{i}(\edge) = \omti{i}(\edge)$ for $i\in\setof{1,2}$. Therefore, by \eqref{irrelevance:unconquered_edge_times_agree_eqn},
\begin{equation}
\label{irrelevance:non_star_times_agree_eqn}
\somp(\edge) = \outcomep(\edge) = \outcomet(\edge) = \somt(\edge)
\quad \text{for all } \edge\not\in \staredges[2]{\ssetp{1}\cup \ssett{2}},
\end{equation}
and thus \eqref{irrelevance:ttimes_comparison_eqn} holds for all such edges.

Now suppose $\edge\in \bdedges[2]{\ssetp{1}\cup \ssett{2}}$. Then $e = \setof{\vect u,\vect v}$ for some $\vect u\in \ssetp{1}\cup \ssett{2}$ and some $\vect v\not\in \ssetp{1}\cup \ssett{2}$. First note that, using Lemma~\ref{edge_sets_lem}, the definition of $\neighbors$, and \eqref{irrelevance:nbrs_covered_eqn}, we have
\[
\staredges[2]{\startsett}\subseteq \edges[2]{\nbrhood[0]{\startsett}}
\subseteq \edges[2]{\neighbors \cup A_2}
\subseteq \edges[2]{\ssetp{1} \cup \ssett{2}},
\]
so $\edge\not\in \staredges[2]{\startsett}$.
Since $\reachedset{1}{\startstatet}{t_0}_{\outcomet} = \startsett$ by \eqref{irrelevance:sp1_unmoved_eqn}, we have $\edge\not\in \staredges[2]{\reachedset{1}{\startstatet}{t_0}_{\outcomet}}$ and hence $\somti{1}(\edge) = \omti{1}(\edge)$ by Definition~\ref{direct_shift_def}. On the other hand, we trivially have $\sompi{1}(\edge) \le \ompi{1}(\edge)$ by Part~\ref{one_type_shift_monotone:t_part} of Lemma~\ref{one_type_shift_monotone_lem}, so using  \eqref{irrelevance:unconquered_edge_times_agree_eqn}, we get
\begin{equation}
\label{irrelevance:ttime1_comparison_eqn}
\sompi{1}(\edge) \le \ompi{1}(\edge)
=\omti{1}(\edge) = \somti{1}(\edge)
\quad\text{for all } e\in \bdedges[2]{\ssetp{1}\cup \ssett{2}}.
\end{equation}
To get the desired inequality for species 2's traversal times, we divide the argument into the two cases $\edge\in \staredges[2]{\ssetp{2}}$ and $\edge\not\in \staredges[2]{\ssetp{2}}$. Suppose first that $\edge\not\in \staredges[2]{\ssetp{2}} = \staredges[2]{\reachedset{2}{\startstatep}{t_0}_{\outcomep}}$. Then by Definition~\ref{direct_shift_def} we have
$\sompi{2}(\edge) = \ompi{2}(\edge)$, and combining this with \eqref{irrelevance:unconquered_edge_times_agree_eqn} and the trivial inequality for $\omti{2}$ from Part~\ref{one_type_shift_monotone:t_part} of Lemma~\ref{one_type_shift_monotone_lem}, we get
\begin{equation}
\label{irrelevance:simple_ttime2_eqn}
\sompi{2}(\edge) = \ompi{2}(\edge) = \omti{2}(\edge) \ge \somti{2}(\edge)
\quad\text{for all } e\in \bdedges[2]{\ssetp{1}\cup \ssett{2}} \setminus \staredges[2]{\ssetp{2}}.
\end{equation}
Now suppose $\edge\in \staredges[2]{\ssetp{2}}$. Then either $\vect u$ or $\vect v$ must be an element of $\ssetp{2}$. Since $\ssetp{2} \subseteq \ssett{2}$ by Claim~\ref{irrelevance:state_comparison_claim}, and $\vect v\not\in \ssetp{1}\cup \ssett{2}$ by definition, it must be that $\vect u\in \ssetp{2}$. It then follows from Claim~\ref{irrelevance:state_comparison_claim} that
\[
\vect u\in \ssetp{2} \text{ and } \vect v\not\in \ssetp{1}\cup \ssetp{2}
\quad\&\quad
\vect u\in \ssett{2} \text{ and } \vect v\not\in \ssett{1}\cup \ssett{2}.
\]
%
%
Therefore, by Definition~\ref{direct_shift_def} we must have
\begin{equation}
\label{irrelevance:shifted_times_eqn}
\sompi{2}(\edge) = \ompi{2}(\edge) - t_0 + \entangledtime{2}{\startstatep}{\vect u}_{\outcomep}
\quad\text{and}\quad
\somti{2}(\edge) = \omti{2}(\edge) - t_0 + \entangledtime{2}{\startstatet}{\vect u}_{\outcomet}.
\end{equation}
We now claim that 
\begin{equation}
\label{irrelevance:infection_times_eqn}
\entangledtime{2}{\startstatep}{\vect u}_{\outcomep}
\ge \entangledtime{2}{\startstatet}{\vect u}_{\outcomet}.
\end{equation}
To see this, observe that by \eqref{irrelevance:sp2_comparison_eqn} we have
\[
\setbuilder[3]{\vect w\in \Zd}{\entangledtime{2}{\startstatep}{\vect w}_{\outcomep} \le t}
=\reachedset{2}{\startstatep}{t}_{\outcomep}
\subseteq \reachedset{2}{\startstatet}{t}_{\outcomet}
= \setbuilder[3]{\vect w\in \Zd}{\entangledtime{2}{\startstatet}{\vect w}_{\outcomet} \le t}
\quad \forall t\le t_0,
\]
so \eqref{irrelevance:infection_times_eqn} follows by taking $t = \entangledtime{2}{\startstatep}{\vect u}_{\outcomep}$ and noting that $t\le t_0$ since $\vect u\in \ssetp{2} = \reachedset{2}{\startstatep}{t_0}_{\outcomep}$. Now, \eqref{irrelevance:unconquered_edge_times_agree_eqn} implies that $\ompi{2}(\edge) = \omti{2}(\edge)$, and combining this with \eqref{irrelevance:shifted_times_eqn} and \eqref{irrelevance:infection_times_eqn} we have
\begin{align}
\label{irrelevance:longer_ttime2_eqn}
\sompi{2}(\edge)
&= \ompi{2}(\edge) - t_0 + \entangledtime{2}{\startstatep}{\vect u}_{\outcomep}\notag \\
&\ge \omti{2}(\edge) - t_0 + \entangledtime{2}{\startstatet}{\vect u}_{\outcomet}
=\somti{2}(\edge)
\quad\text{for all } e\in \bdedges[2]{\ssetp{1}\cup \ssett{2}} \cap \staredges[2]{\ssetp{2}}.
\end{align}
Finally, \eqref{irrelevance:ttime1_comparison_eqn}, \eqref{irrelevance:simple_ttime2_eqn}, and \eqref{irrelevance:longer_ttime2_eqn} together imply that \eqref{irrelevance:ttimes_comparison_eqn} holds for all $\edge\in \bdedges[2]{\ssetp{1}\cup \ssett{2}}$, and hence for all $\edge\in \compedges[2]{\ssetp{1}\cup \ssett{2}}$ when combined with \eqref{irrelevance:non_star_times_agree_eqn}. This completes the proof of Claim~\ref{irrelevance:ttimes_comparison_claim}.
\end{proof}

Together, Claims~\ref{irrelevance:state_comparison_claim} and \ref{irrelevance:ttimes_comparison_claim} show that $\pair[2]{\somp}{\sstatep} \pbetter \pair[2]{\somt}{\sstatet}$ as defined in \eqref{pair_preorder_eqn}, which completes the proof of Claim~\ref{irrelevance:domination_claim}.
\end{proof}

This completes the proof of Theorem~\ref{irrelevance_thm}.
\end{proof}



\chapter{Euclidean Geometry \&\\ Deterministic First-Passage Competition}
\label{deterministic_chap}

In this chapter we introduce one-type and two-type deterministic processes in $\Rd$ which are analogues of the random one-type and two-type first-passage competition processes in $\Zd$ constructed in Chapter~\ref{fpc_basics_chap}. We call these processes simply \textit{deterministic first-passage percolation} and \textit{deterministic first-passage competition} in $\Rd$. The constructions will be essentially identical to the constructions of the lattice processes in Chapter~\ref{fpc_basics_chap}, except that we replace the length structure(s) arising from a traversal measure $\tmeasure$ or $\tmeasurepair$ on $\edges{\Zd}$ with the length structure(s) induced by a norm $\mu$ or pair of norms $(\mu_1,\mu_2)$ on $\Rd$, and we have to take into account some topological and geometric differences between $\Zd$ and $\Rd$.


The significance of the one-type deterministic process is that when we choose the norm $\mu$ to be the shape function corresponding to a random traversal measure $\tmeasure$, the Shape Theorem (Theorem~\ref{shape_thm}) implies that on large scales, with high probability the deterministic process induced by $\mu$ provides good approximations to the random process induced by $\tmeasure$. Similarly, the two-type deterministic process will provide an approximation to the random two-type competition process. Thus, we can hope to gain a better understanding of the random processes by making a careful analysis of their limiting deterministic processes, which are simpler.

Determining the evolution of the deterministic processes is a purely geometric problem, which has been studied computationally for Euclidean and polygonal norms in \cite{Schaudt:1991aa}, \cite{Schaudt:1992aa}, and \cite{Kobayashi:2002aa}. Our goal will be to reduce various probabilistic questions about the random processes to geometric questions about the deterministic processes, by showing that with high probability the random process behaves like the corresponding deterministic process in some sense. The tools for making such reductions, from the random processes to the deterministic, will be developed in Chapters~\ref{cone_growth_chap} and \ref{random_fpc_chap}. In the present chapter, we make a detailed analysis of the relevant geometric aspects of the deterministic growth and competition processes. In order to do so, we will first introduce various elementary concepts from Euclidean convex and metric geometry.

Chapter~\ref{deterministic_chap} is organized as follows. In Section~\ref{deterministic:convex_geom_sec} we introduce various elementary concepts from Euclidean convex geometry, and we define cones, $\mu$-cones, $\mu$-stars, and other related geometric objects that will appear repeatedly throughout the remaining chapters. In Section~\ref{deterministic:one_type_sec} we define induced intrinsic metrics in $\Rd$ and use them to define the deterministic restricted one-type process,
and then we explore some basic properties of this process; 
we will study the stochastic analogue of this one-type process in Chapter~\ref{cone_growth_chap}. In Section~\ref{deterministic:basic_properties_sec} we define the deterministic two-type process and prove some of its basic properties; many of the proofs involve a type of generalized weighted Voronoi cell, which will be one of our main tools for analyzing the geometry of the process.
In Section~\ref{deterministic:euclid_comp_sec} we describe the behavior of the two-type competition process from certain starting configurations when both species' norms are multiples of the Euclidean norm; this simplest case serves as a prototype for understanding the behavior of the more general process. In Section~\ref{deterministic:surrounded_sec} we describe the behavior of the competition process with general norms from certain bounded starting configurations that will be important for our analysis of the random process in Chapter~\ref{random_fpc_chap}. Finally, in Section~\ref{deterministic:cone_competition_sec} we analyze the deterministic process when one species starts on the exterior of a cone and the other species starts at a single interior site, providing the model and necessary geometric analysis needed for describing the analogous behavior of the random two-type process in Chapter~\ref{random_fpc_chap}. 

\section{Elements of Convex Geometry in Euclidean Space}
\label{deterministic:convex_geom_sec}



\subsection{Basic Definitions from Convex Geometry}
\label{deterministic:convex_geom_defs_sec}


This section contains several definitions and basic results from convex geometry that will be needed in our analysis of the deterministic process. Some of these concepts will only be used in proofs located in an appendix, but we collect all the definitions here for easy reference. Additional results are located in Section~\ref{convex_geometry_sec}.

An \textdef{affine combination} of the points $\vect x_1,\dotsc,\vect x_k\in\Rd$ is a linear combination $\sum_{j=1}^k \alpha_j \vect x_j$ with coefficients $\alpha_j\in\R$ satisfying $\sum_{j=1}^k \alpha_j =1$. A \textdef{convex combination} is an affine combination with all the coefficients $\alpha_j\in [0,1]$ (equivalently $\alpha_j\ge 0$). If $B$ is any subset of $\Rd$, the \textdef{affine hull} or \textdef{affine span} of $B$ is the set $\aff B$ of all affine combinations of points in $B$, and coincides with the smallest affine linear subspace containing $B$ (cf.\ Moszynska \cite[p.~25]{Moszynska:2006aa}). For any $B\subseteq\Rd$, we denote by $\affrelbd B$ the \textdef{relative boundary} of $B$ in its affine hull $\aff B$.
%
For a nonnegative integer $k\le d$, a (relatively) open or closed \textdef{$k$-dimensional affine half-space} in $\Rd$ is any set of the form $\Hspace = \setbuilder[2]{\vect x\in A}{\innerprod{\vect x-\anapex}{\dirsymb}>0}$ or $\Hspace = \setbuilder[2]{\vect x\in A}{\innerprod{\vect x-\anapex}{\dirsymb} \ge 0}$, respectively, where $A\subseteq \Rd$ is an affine linear subspace of dimension $k$, $\anapex\in A$, and $\dirsymb\in\Rd\setminus \setof{\vect 0}$ is any vector that is not orthogonal to $A$ (i.e.\ $\innerprod{\vect x -\anapex}{\dirsymb} \ne 0$ for some $\vect x\in A$).
Note that then we have $\aff \Hspace = A$, and $\affrelbd \Hspace = \setbuilder[2]{\vect x\in A}{\innerprod{\vect x-\anapex}{\dirsymb} = 0}$, which is a $(k-1)$-dimensional affine subspace of $A$.
In what follows, the term ``half-space in $\Rd$" used without additional qualifiers will mean a $d$-dimensional open or closed affine half-space in $\Rd$.
If $\Hspace$ is a closed half-space in $\Rd$, then $\Hspace$ is called a \textdef{support half-space} of the set $B\subseteq\Rd$ if $B\subseteq \Hspace$ and $B\cap \bd \Hspace \ne\emptyset$. In this case we say that any $\vect x\in B\cap \bd \Hspace$ is a \textdef{support point} of $B$ and that $\Hplane = \bd \Hspace$ is a \textdef{support hyperplane} of $B$ at $\vect x$. We call a closed half-space $J$ an \textdef{opposing half-space} of $B$ at $\vect x$ if $\shellof[0]{\Hspace} = \Rd\setminus \interior{\Hspace}$ is a support half-space of $B$ at $\vect x$, i.e.\ $B\cap \interior{\Hspace} =\emptyset$ and $\vect x\in B\cap \bd \Hspace$.
For every nonempty compact $B\subseteq\Rd$ and every $\dirsymb\in \Rd\setminus \setof{\vect 0}$, there is a unique support hyperplane $\Hplane$ of $B$ with outer normal vector $\dirsymb$, meaning that $\Hplane = \setbuilder[2]{\vect x\in\Rd}{\innerprod{\vect x - \anapex}{\dirsymb} = 0}$ for some support point $\anapex\in B$ such that $\ray{\anapex}{\dirsymb}\cap B = \emptyset$ \cite[Theorem~2.2.2, p.~14]{Moszynska:2006aa}.

A subset $C$ of $\Rd$ is \textdef{convex} if $\vect x,\vect y\in C \implies \segment{\vect x}{\vect y} \subseteq C$. It follows directly from this definition that an arbitrary intersection of convex sets is convex. For any $A\subseteq\Rd$, the \textdef{convex hull} of $A$, denoted by $\conv A$, is the smallest convex set containing $A$, i.e.\ the intersection of all convex sets containing $A$.
It can be verified (see e.g., \cite[Theorem 3.2.4, p.\ 29]{Moszynska:2006aa}) that the convex hull coincides with the set of convex combinations of points in $A$.
If $C$ is convex, then $\alpha C+ \beta C = (\alpha+\beta)C$ for any $\alpha,\beta \ge 0$ \cite[Proposition 2.3.6, p.~18]{Moszynska:2006aa}.
Every point in the boundary of a closed convex subset of $\Rd$ is a support point \cite[Corollary 3.3.6, p.~34]{Moszynska:2006aa}, and a nonempty closed subset of $\Rd$ is convex if and only if it equals the intersection of all its support half-spaces \cite[Theorem 3.3.7, p.~34]{Moszynska:2006aa}. It is easy to check that the interior and closure of any convex set are convex \cite[Exercise 2.3.8, p.~19]{Moszynska:2006aa}. Furthermore, Lemma~\ref{convex_euclidean_lem} shows that convex sets in $\Rd$ are topologically friendly, in that if $C$ is any closed convex subset of $\Rd$, then there is some $k\le d$ such that $C$ is homeomorphic to either the closed Euclidean ball $\unitball{\ltwo{k}}$, the closed half-space $\closure{\halfspace[k]}$, or the whole Euclidean space $\R^k$.

If $\vect z\in\Rd$, a subset $S$ of $\Rd$ is \textdef{star-shaped at $\vect z$} if $\vect y\in S \implies \segment{\vect z}{\vect y} \subseteq S$, and in this case we call $\vect z$ a \textdef{center} of $S$. Clearly, any convex set $C$ is star-shaped at every $\vect z\in C$. A set that is star-shaped at some point in $\Rd$ is called a \textdef{star set}.
We write $\stars{\vect z}$ for the collection of all sets that are star-shaped at $\vect z$, and $\stars{}$ for the collection of all star sets. More generally, for any $A\subseteq\Rd$, we define the collection of \textdef{star sets at $A$} by
\begin{equation}
\label{stars_at_A_def_eqn}
\starscent{\bullet}{A} \definedas \setbuilder[1]{S\subseteq \Rd}{\forall \vect y\in S,\, \forall \vect z\in A,\; \segment{\vect z}{\vect y}\subseteq S},
\end{equation}
and the collection of \textdef{generalized star sets at $A$} by
\begin{equation}
\label{genstars_at_A_def_eqn}
\genstars{\bullet}{A} \definedas \setbuilder[1]{S\subseteq \Rd}{\forall \vect y\in S,\, \exists \vect z\in A \text{ such that } \segment{\vect z}{\vect y}\subseteq S}.
\end{equation}
Thus, $\starscent{\vect z}{A} = \bigcap_{\vect z\in A} \stars{\vect z} \subseteq \stars{}$, and an element $S$ of $\genstars{\vect z}{A}$ is a union of star sets with centers in $A$, but is not necessarily star-shaped unless $A$ is a singleton.

We define a \textdef{body}\footnote{%
My definition of \emph{body} seems to differ from standard usage, in that bodies in $\Rd$ are usually required additionally to be compact (see, e.g.\ \cite[Definition 14.2.1, p.~177]{Moszynska:2006aa}).
}
\textdef{in $\Rd$} to be a nonempty set $B\subseteq \Rd$ such that $\closure{\interior{B}} = \closure{B}$ and $\interior{\parens[1]{\,\closure{B}\,}} = \interior{B}$, or equivalently $\bd \interior{B} = \bd \closure{B} = \bd B$.
%
%
%
More generally, for an integer $k\le d$, we define an \textdef{embedded $k$-body in $\Rd$} to be a nonempty set $B\subseteq \Rd$ for which there is a topological embedding $\phi\colon \closure{B} \hookrightarrow \R^k$ such that the image $\phi(B)$ is a body in $\R^k$ (it follows that a body in $\Rd$ is the same as an embedded $d$-body).
Any convex subset of $\Rd$ with nonempty interior is a body in $\Rd$ according to our definition (see \cite{Moszynska:2006aa}, Section 2.4, p.~19, and Exercise 14.2.2, p.~177), and moreover, Lemma~\ref{convex_euclidean_lem} implies that any convex set $B\subseteq\Rd$ with $\dim (\aff B) = k$ is an embedded $k$-body. 
A body in $\Rd$ that is star-shaped is called a \textdef{star body}, and we denote the set of star bodies with center $\vect z$ by $\starbodies{\vect z}$.

\subsection{Affine Scale-Invariance: Wedges, Directions, and Cones}
\label{deterministic:scale_inv_sec}

\subsubsection{Wedges (a.k.a.\ Affine Scale-Invariant Sets)}
We call a subset $\asiset$ of $\Rd$ a \textdef{wedge}\footnote{%
My definition of \emph{wedge} coincides with what may be generically
called an \emph{affine cone} in linear algebra, although the restrictiveness of the definitions of such objects varies depending on author and context.
For example, a ``wedge" as defined by Aliprantis and Tourky \cite{Aliprantis:2007aa} would be a convex pointed wedge at $\vect 0$ in my terminology, and then a ``cone" would be a such a wedge with unique apex. I present my own definition of \emph{cone} in the final subsection of Section~\ref{deterministic:scale_inv_sec}.
}
if there exists some $\anapex\in\Rd$ such that $\homothe{r}{\anapex} \asiset = \asiset$ for all $r>0$, where $\homothe{r}{\anapex}\colon \Rd\to\Rd$ is the homothety defined by $\homothe{r}{\anapex}\vect x = \anapex+ r(\vect x-\anapex)$. In this case we call the point $\anapex$ an \textdef{apex} of $\asiset$. For an arbitrary $\asiset \subseteq \Rd$, we define the \textdef{apex set} of $\asiset$ by $\apexset{\asiset} \definedas \setbuilder[2]{\anapex\in\Rd}{\homothe{r}{\anapex} \asiset = \asiset \; \forall r>0}$, so that $\asiset$ is a wedge if and only if $\apexset{\asiset}\ne \emptyset$.
We write $\wedgesat{\anapex} \definedas \setbuilder{W\subseteq\Rd}{\anapex\in\apexset{W}}$ for the set of all \textdef{wedges at $\anapex$}.
If $\vect 0\in \apexset{\asiset}$, we say that $\asiset$ is \textdef{scale-invariant}. More generally, for any $\anapex\in \apexset{\asiset}$, we say that $\asiset$ is \textdef{scale-invariant at $\anapex$} or that $\asiset$ is \textdef{affine scale-invariant} with apex $\anapex$. It is easily seen that $\apexset{\asiset} = \apexset[2]{\Rd\setminus \asiset}$, so $\asiset\in \wedgesat{\anapex}$ if and only if $\Rd\setminus \asiset\in \wedgesat{\anapex}$.
We call a wedge $\asiset$ \textdef{trivial} if $\asiset\in \setof[2]{\emptyset,\Rd}$.

\begin{example}[Prototypical wedges]
\label{wedge_prototype_ex}
Two prototypical examples of wedges are:
\begin{itemize}
\item 
An infinite cone as defined in elementary Euclidean geometry. Then the apex (a.k.a.\ vertex) of the cone is the unique element of the apex set.

\item 
The wedge-shaped region $\asiset$ defined as follows: Start with a plane in $\R^3$, then ``fold" or ``crease" the plane along some line contained in it, so that the intersection of the folded plane with any plane orthogonal to the fold-line comprises two rays with a common origin, and then take $\asiset$ to be all the points on one side of the folded plane. In this case the set of apexes is the line defining the fold (as long as $\asiset$ is assumed to be either open or closed).

\end{itemize}
In both of these examples, the wedge is a body whose closure is homeomorphic to a closed half-space, so the wedges are in fact \emph{cones} according to our definition at the end of Section~\ref{deterministic:scale_inv_sec} below.
\end{example}

A wedge is \textdef{pointed} if it contains at least one of its apexes, and \textdef{blunt} otherwise. It is easy to see that a set $\asiset\subseteq \Rd$ is a wedge at $\anapex$ if and only if either $\asiset = \setof{\anapex}$ or $\asiset$ is a union of pointed or blunt rays originating at $\anapex$. This implies that an arbitrary union or intersection of wedges with apex $\anapex$ is also a wedge at $\anapex$.
If $\asiset$ is a wedge with apex $\anapex$, we call the parallel scale-invariant set $\siversion{\asiset} \definedas \asiset-\anapex$ the \textdef{scale-invariant version} of $\asiset$. The following lemma justifies the use of the article ``the" in this definition and enumerates several other useful facts about wedges.


\begin{lem}[Properties of wedges]
\label{wedge_properties_lem}
\newcommand{\mysiset}{\siversion{\asiset}}
\newcommand{\myapexset}{\apexset{\asiset}}

Let $\asiset$ be a wedge in $\Rd$.
\begin{enumerate}
\item \label{wedge_properties:siversion_unique_part}
The scale-invariant version $\mysiset$ of $\asiset$ is unique.

\item \label{wedge_properties:apex_subspace_part}
$\myapexset$ is an affine subspace of $\Rd$ which can be identified with the vector space of translations that fix $\mysiset$. 

\item \label{wedge_properties:apex_bd_part}
If $\asiset$ is nontrivial (i.e.\ $\asiset\notin \setof[2]{\emptyset,\Rd}$),
then $\myapexset\subseteq \bd \asiset$. 

\item \label{wedge_properties:more_wedges_part}
The sets $\interior{\asiset}$, $\closure{\asiset}$, $\bd \asiset$, $\exterior{\asiset}$, and $\shellof{\asiset}$ and  are all wedges whose apex sets contain $\apexset{\asiset}$. 

\item \label{wedge_properties:pointed_part}
If $\asiset$ is pointed, then it contains all of its apexes, and $\Rd\setminus \asiset$ is blunt.

\item \label{wedge_properties:apex_support_part}
If $\asiset$ is pointed and contained in some half-space, then every apex of $\asiset$ is a support point of $\asiset$. 

\item \label{wedge_properties:convexity_part}
$\asiset \subseteq \asiset + \mysiset$, with equality if and only if $\asiset$ is convex.
\end{enumerate}
\end{lem}

We prove Lemma~\ref{wedge_properties_lem} in Appendix~\ref{leftover_chap} (Section~\ref{leftover:deterministic_sec}); many of the parts rely on Lemma~\ref{affine_si_span_lem}, which requires a bit of geometric intuition.

We now introduce a few more definitions for a wedge $\asiset$. We call $\asiset$ \textdef{sharp} if $\apexset{\asiset}$ contains exactly one element and \textdef{flat} otherwise. By Part~\ref{wedge_properties:apex_subspace_part} of Lemma~\ref{wedge_properties_lem}, $\asiset$ is flat if and only if its apex set contains a nontrivial affine subspace. We call $\asiset$ \textdef{small} if it is contained in a half-space and \textdef{large} otherwise. Moreover, we call $\asiset$ \textdef{extra-small} if $\asiset$ is small and the supremum of angles between rays in $\asiset$ is strictly smaller than $\pi$, and we call $\asiset$ \textdef{extra-large} if the complementary wedge $\Rd\setminus \asiset$ is extra-small. An extra-small or extra-large wedge must be sharp, but not conversely; for example, the ``pointed upper half-pace" $\setof{\vect 0}\cup \halfspace[d]$ and its (blunt) complement are both small sharp wedges in $\Rd$, but are not extra-small.

%

If $\asiset$ is a wedge and $\siversion{\asiset}$ is its scale-invariant version, we define the \textdef{spherical section} of $\asiset$ with respect to a norm $\norm{\cdot}$ to be the intersection $\siversion{\asiset}\cap \sphere{d-1}{\norm{\cdot}}$. Note that a pointed (or blunt) wedge is uniquely determined by its apex and its spherical section. For $\anapex\in\Rd$, we denote the set of \textdef{closed wedges at $\anapex$} by $\cwedgesat{\anapex} \definedas \setbuilder[1]{\closure{\asiset}}{\asiset\in \wedgesat{\anapex}}$, and we topologize $\cwedgesat{\anapex}$ using the Hausdorff metric on spherical sections of its wedges.

If $I\subseteq\R$, $B\subseteq \Rd$, and $\anapex\in\Rd$, we define the \textdef{$I$-partial wedge generated by $B$} at $\anapex$ to be $\homothe{I}{\anapex} B = \setbuilder{\homothe{r}{\anapex} \vect x}{r\in I,\vect x\in B}$; note that $\homothe{I}{\anapex} B$ is a pointed wedge if $I=\Rplus$ and a blunt wedge if $I=(0,\infty)$. It is easily seen that if $B$ is convex and $I$ is an interval, then $\homothe{I}{\anapex} B$ is convex; in fact, if $B$ is convex and $I=[\alpha,\beta]$ for some $0\le\alpha\le\beta< \infty$, then $\homothe{I}{\anapex} B = \convhull[2]{\homothe{\alpha}{\anapex} B\cup \homothe{\beta}{\anapex} B}$.

\subsubsection{Directions in Wedges; Visible Directions in Subets of $\Rd$}

If $\dirsymb\in \Rd\setminus \setof{\vect 0}$, we define the \textdef{direction} of $\dirsymb$ to be the blunt ray $\dirvec = \setbuilder{r\dirsymb}{r>0}$. If $\asiset$ is a wedge, we define the \textdef{direction space} of $\asiset$ to be $\dirset{\asiset} \definedas \setbuilder[2]{\dirvec}{\dirsymb\in \siversion{\asiset}\setminus \setof{\vect 0}}$, where $\siversion{\asiset}$ is the scale-invariant version of $\asiset$ (recall that $\siversion{\asiset}$ is unique by Part~\ref{wedge_properties:siversion_unique_part} of Lemma~\ref{wedge_properties_lem}). Then $\dirset{\asiset}$ is the quotient space of $\siversion{\asiset}\setminus \setof{\vect 0}$ under the equivalence relation ``$\dirsymb\direquiv r\dirsymb$ if $r>0$," and we endow $\dirset{\asiset}$ with the quotient topology.\footnote{%
Observe that the construction of $\dirset{\asiset}$ from $\asiset$ mimics the construction of the projective space $P(V)$ from a vector space $V$, except that our equivalence classes are rays instead of lines.
}
If we fix some norm $\norm{\cdot}$, then each direction $\dirvec$ in $\dirset{\asiset}$ corresponds to a unique $\norm{\cdot}$-unit vector $\unitvecnorm{\dirsymb}{\norm{\cdot}}$ in $\siversion{\asiset}$, and it is easily verified (e.g.\ using the Closed Map Lemma \cite[Lemma~4.25]{Lee:2000aa}) that this correspondence is a homeomorphism between $\dirset{\asiset}$ and the spherical section $\siversion{\asiset} \cap \sphere{d-1}{\norm{\cdot}}$.

\newcommand{\viswedge}[3][0]{\vec \asiset_{#2} \parens[#1]{#3}}

For any $A\subseteq \Rd$ and $\vect z\in \closure{A}$, let $\viswedge{\vect z}{A} \definedas \bigcup \setbuilder[2]{\vect z+\dirvec}{\vect z+\dirvec \subseteq A}$. Then $\viswedge{\vect z}{A}$ is a (possibly empty) blunt wedge at $\vect z$ and contained in $A$, which we call the \textdef{visible wedge in $A$ at $\vect z$}. Now define the \textdef{visible directions in $A$ at $\vect z$} by $\dirsetat{\vect z}{A} \definedas \dirset[2]{\viswedge{\vect z}{A}}$. If we think of $A$ as ``open space," with points in the complement $\Rd\setminus A$ being ``obstacles," then $\dirsetat{\vect z}{A}$ is the set of ``unblocked" directions in $A$ seen from $\vect z$, meaning that if a beam of light is aimed in the direction $\dirvec\in \dirsetat{\vect z}{A}$ from the point $\vect z$,
then it will travel through $A$ forever, avoiding the obstacles in the complement.
It follows from the definitions that if $\asiset$ is a wedge, then for any $\anapex\in \apexset{\asiset}$, $\viswedge{\anapex}{\asiset} = \anapex + \bigcup \dirset{\asiset} = \asiset \setminus \setof{\anapex}$, and $\dirsetat{\anapex}{\asiset} = \dirset{\asiset}$. The following lemma shows that for any $\vect z\in\asiset$, we have the containment $\dirsetat{\vect z}{\asiset} \subseteq \dirset{\closure{\asiset}}$, since $\viswedge{\vect z}{\asiset} \subseteq \asiset$ by definition.
\begin{lem}[Monotonicity of direction sets]
\label{dirset_monotone_lem}
If $\asiset$ and $\asiset'$ are both affine scale-invariant (possibly with different apex sets), and $\asiset'\subseteq \asiset$, then $\dirset{\asiset'} \subseteq \dirset{\closure{\asiset}}$, and the quotient topology on $\dirset{\asiset'}$ agrees with the subspace topology inherited from the quotient topology on $\dirset{\closure{\asiset}}$.
\end{lem}

\begin{proof}

\newcommand{\origin}{\vect 0}
\newcommand{\thisa}{\vect a}
\newcommand{\thisx}{\vect x}
\newcommand{\axray}{\ray{\thisa}{\thisx}}
\newcommand{\scaledpt}[1]{\homothe{#1}{\thisa} \thisx} 

Without loss of generality, assume $\asiset$ has apex $\origin$. Let $\thisa$ be an apex of $\asiset'$, and let $\dirvec \in\dirset{\asiset'}$. Then $\dirvec = -\thisa + \axray$ for some $\thisx\in \asiset'$, and the ray $\axray$
is contained in $\asiset'$ since $\asiset'$ is scale-invariant at $\thisa$. For $r>0$, let $\scaledpt{r} \definedas \thisa + r(\thisx-\thisa) \in \axray$. Since $\asiset'\subseteq \asiset$ we have $\scaledpt{r}\in \asiset$ for all $r>0$, and since $\asiset$ is scale-invariant, we then have $s\cdot  \scaledpt{r}\in \asiset$ for all $r,s>0$. Taking $s=\frac{1}{r}$, we have $\frac{1}{r}\cdot  \scaledpt{r} = \frac{1}{r} \thisa + (\thisx- \thisa)\in \asiset$ for all $r>0$. Letting $r\to\infty$, we see that $\thisx-\thisa\in \closure{\asiset}$. Therefore, $\dirset{\closure{\asiset}}$ contains the ray $\ray{\origin}{(\thisx - \thisa)} = -\thisa + \axray = \dirvec$. The statement about topologies then follows from abstract nonsense. 
\end{proof}


\subsubsection{Scale-Equivariance and Scale-Invariance of Functions}


If $A\subseteq \Rd$ is a scale-invariant set, a function $\mu : A\to\R_+$ is called \textdef{scale-equivariant} or \textbf{positive homogeneous of degree 1} if $\mu(r\vect x) = r\mu(\vect x)$ for all $\vect x\in A$ and $r> 0$. Note that a norm $\norm{\cdot}$ on $\Rd$ is a scale-equivariant function which is also subadditive, positive definite, and even. If $A$ is scale-invariant and $\norm{\cdot}$ is any norm on $\Rd$, it is easy to see that the function $f(\vect x) = \dist{\norm{\cdot}}(A,\vect x)$ is scale-equivariant and continuous on $\Rd$.

A function $\psi: A\setminus \setof{\vect 0} \to Y$, where $Y$ is any nonempty set, is called \textdef{scale-invariant} if $\psi(r\vect x) = \psi(\vect x)$ for all $\vect x\in A$ and $r>0$. Clearly, a ratio of two scale-equivariant functions is a scale-invariant function into $\Rplus$.
Moreover, for any set (or topological space) $Y\ne \emptyset$, there is a natural correspondence between (continuous) scale-invariant functions $\psi : A\setminus \setof{\vect 0}\to Y$ and arbitrary (continuous) functions $\tilde \psi: \dirset{A}\to Y$, or equivalently, (continuous) functions from any spherical section of $A$ into $Y$.

More generally, if $A$ is affine scale-invariant with apex $\vect a$, we say that a function $\mu\colon A\to\Rplus$ is \textdef{scale-equivariant at $\vect a$} if $\mu\parens[2]{\homothe{r}{\vect a} \vect x} = r\mu\parens[0]{\vect x}$ for all $\vect x\in A$ and $r>0$, where $\homothe{r}{\vect a} \definedas \vect a + r(\vect x -\vect a)$ is the homothety with scale factor $r$ and center $\vect a$, and we say that a function $\psi\colon A\setminus \setof{\vect a}\to Y$ is \textdef{scale-invariant at $\vect a$} if $\psi \parens[2]{\homothe{r}{\vect a} \vect x} = \psi(\vect x)$ for all $\vect x\in A$ and $r>0$.

\subsubsection{Cones in $\Rd$}

For an integer $k\le d$, we define a \textdef{cone of dimension $k$} or \textdef{$k$-cone} in $\Rd$ to be a wedge that is is also an embedded $k$-body in $\Rd$ whose closure is homeomorphic to the closed upper half-space $\closure{\halfspace[k]} = \setbuilder{\vect x \in \R^k}{ x_k \ge 0}$. Equivalently, $\cones\subseteq\Rd$ is a $k$-cone if $\cones$ is a
wedge whose spherical section $\siversion{\cones}\cap \sphere{d-1}{\ltwo{d}}$ is an embedded $(k-1)$-body with closure homeomorphic to the closed Euclidean ball $\unitball{\ltwo{k-1}}$. 
We say simply that $\cones\subseteq \Rd$ is a \textdef{cone} in $\Rd$ if $\cones$ is a $k$-cone for some $k\le d$, in which case we say that $\cones$ is \textdef{degenerate} if $k<d$ and \textdef{nondegenerate} if $k=d$.

A cone in $\Rd$ is degenerate if and only if it has empty interior. Any open wedge is a union of open cones.
Lemma~\ref{convex_euclidean_lem} implies that any convex wedge is either a cone or an affine subspace of $\Rd$.
Our prototypical example of a cone is the wedge generated by a $\norm{\cdot}$-ball for some norm $\norm{\cdot}$.
Specifically, if $\vect x\in\Rd$ and $\norm{\vect x}> r>0$, then the pointed scale-invariant set generated by the closed ball $\vect x+r\unitball{\norm{\cdot}}$ is a nondegenerate extra-small convex closed cone at $\vect 0$. We expand on this construction in the next section.


\begin{remark}
The primary reason for allowing nonconvex cones is that we want the flexibility to consider some family of extra-large cones interpolating between the upper half-space $\halfspace[d]$ and the ``improper cone" $\Rd\setminus \setof{\direction{\basisvec{k}}}$, as these are the two configurations considered by Deijfen and \Haggstrom \cite{\DHunbounded}. In the present work the only interesting result we are able to prove about such cones is the large-deviations estimate in Theorem~\ref{wide_cone_survival_thm}, but the general framework we develop should be useful for future projects. One such project would be to study competition in critical cones (defined in Section~\ref{deterministic:cone_competition_sec}), including extra-large ones. At present we are unable to prove any results regarding critical cones in the random process, but the consideration of critical cones is also one reason for allowing degenerate cones in our general definition, as we conjecture that species 1 will ``typically" conquer a degenerate subcone of a critical cone in the deterministic process (though we don't attempt to prove this statement).
\end{remark}

\subsection{$\mu$-Cones and Other Subsets of $\Rd$ Defined by a Norm}
\label{deterministic:norm_geometry_sec}


\subsubsection{$\mu$-Cones and $\mu$-Cone Segments}
\label{deterministic:mu_cones_sec}

Let $\mu$ be a norm on $\Rd$. We define a \textdef{$\mu$-cone} to be any wedge generated by a $\mu$-ball. A $\mu$-cone is a cone as defined in Section~\ref{deterministic:scale_inv_sec} as long as the chosen apex is not in the interior of the ball; otherwise, the $\mu$-``cone" is all of $\Rd$. More generally, we define a \textdef{$\mu$-cone segment} to be an $I$-partial wedge generated by a $\mu$-ball for any interval $I\subseteq \Rplus$, and we call this partial wedge an \textdef{initial $\mu$-cone segment} if $0\in I$.

More explicitly, for any $\vect z\in \Rd$, $\dirvec\in \dirset{\Rd}$, $\thp\ge 0$, and $0\le h_1\le h_2\le\infty$,\footnote{%
For $A\subseteq\Rd$ we interpret the expression $\infty\cdot A$ as a limit, i.e.\ $\infty \cdot A \definedas \lim_{r\to\infty} rA$, which can be defined rigorously using various notions of the limit of a sequence of sets. Then, extending the notation for homotheties, we define for $\anapex \in\Rd$ a map $\homothe{\infty}{\anapex} \colon \powerset{\Rd}\to \powerset{\Rd}$ by $\homothe{\infty}{\anapex} \cdot A \definedas \lim_{r\to\infty} \homothe{r}{\anapex}A \definedas \anapex + \lim_{r\to\infty} r(A-\anapex)$. Using this convention allows us to write various formulas without treating $\infty$ as a special case.
} we define the \textdef{$\mu$-cone segment} $\conesegment{\mu}{\thp}{\vect z}{\dirvec}{h_1}{h_2}$ to be the  $\segment{h_1}{h_2}$-partial wedge generated by the $\mu$-ball $\reacheddset{\mu}{\vect z+\unitvecnorm{\dirsymb}{\mu}}{\thp}$ at the apex $\vect z$, i.e.
\begin{equation}
\label{mu_cone_seg_def_eqn}
\conesegment{\mu}{\thp}{\vect z}{\dirvec}{h_1}{h_2}
\definedas \homothe{[h_1,h_2]}{\vect z} \cdot
\reacheddset{\mu}{\vect z+\unitvecnorm{\dirsymb}{\mu}}{\thp}.
\end{equation}
(Recall that $\unitvecnorm{\dirsymb}{\mu}$ is the $\mu$-unit vector in the direction $\dirvec$.) We call the ray $\vect z+\dirvec$ the \textdef{axis}, the radius $\thp$ the \textdef{thickness},
and the lengths $h_1$ and $h_2$ the \textdef{initial and final heights} of the cone segment. We also refer to $\vect z$ as the \textdef{apex} even though the cone segment may not be a full wedge.
If $h_1=0$, we omit it from the notation, writing $\conetip{\mu}{\thp}{\vect z}{\dirvec}{h} \definedas \conesegment{\mu}{\thp}{\vect z}{\dirvec}{0}{h}$ for the \textdef{initial $\mu$-cone segment} of height $h\in[0,\infty]$. Finally, we define the \textdef{$\mu$-cone} with apex $\vect z$, axis $\vect z+\dirvec$, and thickness $\thp$ to be
\begin{equation}
\label{mu_cone_def_eqn}
\cone{\mu}{\thp}{\vect z}{\dirvec}
\definedas \conetip{\mu}{\thp}{\vect z}{\dirvec}{\infty}
= \vect z + \Rplus \cdot \reacheddset{\mu}{\unitvecnorm{\dirsymb}{\mu}}{\thp},
\end{equation}
which is the wedge at $\vect z$ generated by the $\mu$-ball $\reacheddset{\mu}{\vect z+\unitvecnorm{\dirsymb}{\mu}}{\thp}$. If $\thp>1$, this $\mu$-ball contains the point $\vect z$ in its interior, and hence the $\mu$-``cone" $\cone{\mu}{\thp}{\vect z}{\dirvec}$ is the trivial wedge $\Rd$ rather than an actual cone, but we nevertheless use the term $\mu$-cone in this case.

It follows from the definition \eqref{mu_cone_seg_def_eqn} and the convexity of $\mu$-balls that if $h_2<\infty$, then
\begin{equation}
\label{mu_cone_seg_equiv_eqn}
\newcommand{\zballat}[2][0]{\reacheddset[#1]{\mu}{\vect z+#2 \unitvecnorm{\dirsymb}{\mu}}{\thp #2}}
\newcommand{\uballat}[2][0]{#2 \reacheddset[#1]{\mu}{\unitvecnorm{\dirsymb}{\mu}}{\thp}}
\conesegment{\mu}{\thp}{\vect z}{\dirvec}{h_1}{h_2}
=\bigcup_{h\in [h_1,h_2]} \zballat{h}
= \vect z+ \convhull[1]{\uballat{h_1}\cup \uballat{h_2}}.
\end{equation}
Evidently, $\thp\le \thp' \implies \conesegment{\mu}{\thp}{\vect z}{\dirvec}{h_1}{h_2} \subseteq \conesegment{\mu}{\thp'}{\vect z}{\dirvec}{h_1}{h_2}$, so it makes sense to define
\[
\conesegment{\mu}{\thp^-}{\vect z}{\dirvec}{h_1}{h_2}
\definedas \bigcup_{\thp'<\thp} \conesegment{\mu}{\thp}{\vect z}{\dirvec}{h_1}{h_2}
\quad\text{and}\quad
\conesegment{\mu}{\thp^+}{\vect z}{\dirvec}{h_1}{h_2}
\definedas \bigcap_{\thp'>\thp} \conesegment{\mu}{\thp}{\vect z}{\dirvec}{h_1}{h_2}.
\]
We will be primarily interested in full $\mu$-cones $\cone{\mu}{\thp}{\vect z}{\dirvec}$ and initial $\mu$-cone segments $\conetip{\mu}{\thp}{\vect z}{\dirvec}{h}$. The following lemma describes these sets in more detail for different values of $\thp$ and $h$. The proof consists of checking definitions and is left to the reader.

\begin{lem}[Description of initial $\mu$-cone segments]
\label{mu_cone_description_lem}

\newcommand{\thisz}{\vect z}
\newcommand{\coneth}[1]{\cone{\mu}{#1}{\thisz}{\dirvec}}
\newcommand{\conetipth}[2]{\conetip{\mu}{#1}{\thisz}{\dirvec}{#2}}
\newcommand{\thisunitvec}{\unitvec{\dirsymb}_\mu}

Let $\mu$ be a norm, let $\thisz\in\Rd$, $\dirsymb\in\Rd\setminus \setof{\vect 0}$, $\thp\ge 0$, and $h\in [0,\infty]$.

\begin{enumerate}

\item \label{mu_cone_description:th=0_part}
If $\thp = 0$, then $\conetipth{0}{h} = \thisz+[0,h]\cdot \thisunitvec$ is a line segment if $h<\infty$ and a ray if $h=\infty$. 
If $h=0$, then $\conetipth{\thp}{0} = \setof{\vect z}$ for any $\thp\ge 0$.

\item \label{mu_cone_description:th<1_part}
If $0<\thp<1$, then $\conetipth{\thp}{h}$ is a compact convex body if $0<h<\infty$, and if $h=\infty$, then $\coneth{\thp}$ is an extra-small, nondegenerate, convex closed cone with unique apex $\thisz$.

\item \label{mu_cone_description:th>1_part}
If $\thp\ge1$ and $h<\infty$, then $\conetipth{\thp}{h} = \reacheddset{\mu}{\thisz+h \unitvec{\dirsymb}_\mu}{\thp h}$. If $\thp>1$ and $h=\infty$, then $\coneth{\thp} = \Rd$.

\item \label{mu_cone_description:th=1_part}
If $\thp=1$, then $\coneth{1}$ is a small, pointed, nondegenerate convex cone with apex $\thisz$. If there is a unique support hyperplane $\Hplane$ at $\thisz$, then $\closure{\coneth{1}}$ is a closed half-space with boundary $\Hplane$.


\item \label{mu_cone_description:th-_part}
For any $\thp>0$, $\conetipth{\thp^-}{h} = \setof{\thisz} \cup \interior{\parens[1]{\conetipth{\thp}{h}}}$.


\item \label{mu_cone_description:th+_part}
If $\thp\ne 1$, then $\conetipth{\thp^+}{h} = \conetipth{\thp}{h}$. If $\thp = 1$, then 
$\conetipth{1^+}{h} = \reacheddset{\mu}{\thisz+h \unitvec{\dirsymb}_\mu}{h}$ if $h<\infty$,
and $\coneth{1^+} = \Rd$.

\end{enumerate}

\end{lem}


For any $\vect z\in\Rd$ and $\thp\ge 0$, we will use the notation $\conefamily{\mu}{\vect z,\thp}$ for the \textdef{set of all initial $\mu$-cone segments} of thickness $\thp$ with apex $\vect z$, i.e.
\begin{equation}
\label{cone_fam_def_eqn}
\conefamily{\mu}{\vect z,\thp} \definedas
\setbuilderbar[1]{\conetip{\mu}{\thp}{\vect z}{\dirvec}{h}}
	{\dirvec\in \dirset{\Rd},\, h\in [0,\infty]}.
\end{equation}
We call any member of $\conefamily{\mu}{\vect z,\thp}$ a \textdef{$\thp$-thick $\mu$-cone segment at $\vect z$} or simply a \textdef{$(\mu,\thp)$-cone segment}. We make the following observations about the collection $\conefamily{\mu}{\vect z,\thp}$.

\begin{lem}[Decomposing $\mu$-balls and $\mu$-cone segments]
\label{mu_cone_decomp_lem}
\newcommand{\smallthp}{\thp}
\newcommand{\bigthp}{\thp'}
\newcommand{\thisapex}{\vect z}
\newcommand{\ballcenter}{\vect z}
\newcommand{\thisball}{\reacheddset{\mu}{\ballcenter}{r}}
\newcommand{\smallconefam}{\conefamily{\mu}{\thisapex,\smallthp}}

Let $\thisapex\in\Rd$ and let $\smallthp\ge 0$.
\begin{enumerate}
\item \label{mu_cone_decomp:ball_part}
For any $r\ge 0$, the $\mu$-ball $\thisball$ is a union of members of $\smallconefam$.

\item \label{mu_cone_decomp:cone_part}
If $\bigthp\ge \smallthp$, then every member of $\conefamily{\mu}{\thisapex,\bigthp}$ is a union of members of $\smallconefam$.
\end{enumerate}
\end{lem}

\begin{proof}
Both parts follow from the triangle inequality for $\mu$.
\end{proof}

Here are some further results regarding $\mu$-cones that will be needed in Section~\ref{deterministic:cone_competition_sec} and in later chapters.
The proof of Lemma~\ref{mu_cone_bd_lem} is an elementary exercise, and Lemma~\ref{special_mu_properties_lem} follows from the convexity of $\mu$-balls and $\mu$-cones.

%
%
%

\begin{lem}[Boundaries of $\mu$-cones]
\label{mu_cone_bd_lem}
\newcommand{\thisaxis}{{\anapex+\cdirvec}}
\newcommand{\thiscone}{\cone{\mu}{\thp}{\anapex}{\dirvec}}
\newcommand{\thisx}{\vect x}
\newcommand{\thisy}{\vect y}

Let $\anapex\in\Rd$, $\dirvec\in \dirset{\Rd}$, and $\thp\in[0,1]$, and let $\cones = \thiscone$.
\begin{enumerate}

\item \label{mu_cone_bd:axis_part}
$\thisaxis = \setbuilder[2]{\thisx\in \cones}{\distnew{\mu}{\bd \cones}{\thisx} = \thp \distnew{\mu}{\anapex}{\thisx}}$.

\item \label{mu_cone_bd:bd_part}
$\bd \cones = \setbuilder[2]{\thisy\in\Rd}{\exists \thisx\in \thisaxis \text{ with } \distnew{\mu}{\thisy}{\thisx} = \thp \distnew{\mu}{\anapex}{\thisx}}$.

\end{enumerate}

\end{lem}


\begin{lem}[Special geometry of $\mu$-balls and $\mu$-cones]
\label{special_mu_properties_lem}
Let $\mu$ be a norm on $\Rd$. Then for any $\alpha\ge 0$, $r\ge 0$, and $\thp>0$,
\begin{enumerate}
\item \label{special_mu:ball_separation_part}
$\displaystyle
\setbuilder[1]{\vect x\in \reacheddset{\mu}{\vect z}{r}}{\distnew[2]{\mu}{\bd \reacheddset{\mu}{\vect z}{r}}{\,\vect x} \ge \alpha}
= \reacheddset{\mu}{\vect z}{r-\alpha}
$ 

\item \label{special_mu:cone_separation_part}
$\displaystyle
\setbuilder[1]{\vect x\in \cone{\mu}{\thp}{\vect z}{\dirvec}}
	{\distnew[2]{\mu}{\bd \cone{\mu}{\thp}{\vect z}{\dirvec}}{\,\vect x} \ge \alpha}
= \frac{\alpha }{\thp}\cdot \munitvec
	+ \cone{\mu}{\thp}{\vect z}{\dirvec}.
$

\item \label{special_mu:cone_intersection_part}
$\displaystyle
\bigcap_{\anapex\in \reacheddset{\mu}{\vect z}{\alpha}}
	\cone{\mu}{\thp}{\anapex}{\dirvec}
=\frac{\alpha}{\thp}\cdot \munitvec
	+\cone{\mu}{\thp}{\vect z}{\dirvec}.
$
\end{enumerate}
\end{lem}


%

\subsubsection{$\mu$-Bowling  Pins and $\mu$-Tubes}

An important construct in our analysis of the competition process will be to take the union of an initial $\mu$-cone segment with a $\mu$-ball centered at its tip (i.e.\ apex). Due to its shape, we call such a subset of $\Rd$ a \textdef{$\mu$-bowling pin}.

More explicitly, for $\vect z\in \Rd$, $\dirvec\in\dirset{\Rd}$, $r\ge 0$, $\thp \ge 0$, and $h\in [0,\infty]$, we define the \textdef{$\mu$-bowling pin}
\begin{equation}
\label{bpin_def_eqn}
\bpin{\mu}{r}{\thp}{\vect z}{\dirvec}{h} \definedas
\reacheddset{\mu}{\vect z}{r} \cup \conetip{\mu}{\thp}{\vect z}{\dirvec}{h}.
\end{equation}
Thus, $\bpin{\mu}{r}{\thp}{\vect z}{\dirvec}{h}$ is a ``bowling pin"-shaped region of height $h$ along the axis $\vect z+\dirvec$, with thickness $\thp$ and a ``head" of radius $r$ centered at $\vect z$. We call the point $\vect z$ the \textdef{origin} of the bowling pin; note that any bowling pin is star-shaped at its origin since the sets $\reacheddset{\mu}{\vect z}{r}$ and $\conetip{\mu}{\thp}{\vect z}{\dirvec}{h}$ are convex and contain $\vect z$. We call $\bpin{\mu}{r}{\thp}{\vect z}{\dirvec}{h}$ \textdef{nondegenerate} if $r$, $\thp$, and $h$ are all strictly positive; note that a nondegenerate bowling pin is a body in $\Rd$ which is compact if $h<\infty$. More generally, we may call any body in $\Rd$ a $\mu$-bowling pin if its interior is equal to the interior of some $\bpin{\mu}{r}{\thp}{\vect z}{\dirvec}{h}$ as defined above.

{
\newcommand{\hmax}{h}
\newcommand{\genh}{s}

All the ``$\mu$-objects" defined so far, i.e.\ $\mu$-balls, $\mu$-cones, $\mu$-cone segments, and $\mu$-bowling pins, are special cases of the following more general construction of ``$\mu$-tubes." For any $\vect z\in\Rd$, $\dirvec\in \dirset{\Rd}$, $\hmax \in [0,\infty)$, and nonnegative function $\rfcn\colon \Rplus\to\Rplus$, define the \textdef{initial $\mu$-tube segment} of height $\hmax$ with radius function $\rfcn$, origin $\vect z$, and axis $\vect z+\dirvec$, by
\begin{equation}
\label{mu_tube_segment_def_eqn}
\tubetip{\mu}{\rfcn}{\vect z}{\dirvec}{\hmax}
	\definedas \bigcup_{\genh\in [0,\hmax]}
	\reacheddset[2]{\mu}{\vect z+\genh\munitvec}{\rfcn(\genh)}.
\end{equation}
Then define the infinite \textdef{$\mu$-tube} by
\begin{equation}
\label{mu_tube_def_eqn}
\tube{\mu}{\rfcn}{\vect z}{\dirvec}
	\definedas \tubetip{\mu}{\rfcn}{\vect z}{\dirvec}{\infty}
	\definedas \bigcup_{\hmax\ge 0} \tubetip{\mu}{\rfcn}{\vect z}{\dirvec}{\hmax}.
\end{equation}
For example, $\conetip{\mu}{\thp}{\vect z}{\dirvec}{\hmax} = \tubetip{\mu}{\rfcn}{\vect z}{\dirvec}{\hmax}$ with $\rfcn(\genh) = \thp \genh$, and $\bpin{\mu}{r}{\thp}{\vect z}{\dirvec}{\hmax} = \tubetip{\mu}{\rfcn}{\vect z}{\dirvec}{\hmax}$ with $\rfcn(0) = r$ and $\rfcn(\genh) = \thp \genh$ for $\genh>0$. In Section~\ref{cone_growth:logarithmic_tube_sec} we will consider $\mu$-tubes with a radius function $\rfcn$ that grows faster than logarithmic but slower than linear.
}

\subsubsection{$\mu$-Stars}

Motivated by Lemma~\ref{mu_cone_decomp_lem} above, we make the following definitions.
For $\vect z\in\Rd$ and $\thp\ge 0$, we define the set of
\textdef{$(\mu,\thp)$-stars at $\vect z$} to be
\begin{equation}
\label{fat_stars_def_eqn}
\begin{split}
\fatstars{\mu}{\thp}{\vect z}
&\definedas \setof[2]{\text{Unions of elements in $\conefamily{\mu}{\vect z,\thp}$}}
= \setbuilderbar[2]{{\textstyle \bigcup}\, {\conesubfam}}
	{\conesubfam \subseteq \conefamily{\mu}{\vect z,\thp}}.
\end{split}
\end{equation}
That is, $\fatstars{\mu}{\thp}{\vect z}$ is the collection of star sets at $\vect z$ that are ``$(\mu,\thp)$-thick" in the sense that every point is contained not just in an infinitely thin line segment originating at $\vect z$, but in a $\thp$-thick $\mu$-cone segment at $\vect z$.
Using the triangle inequality, one can verify that the definition \eqref{fat_stars_def_eqn} coincides with the definition of $(\mu,\thp)$-stars at $\vect z$ given in Section~\ref{intro:ch4_sec}.
We call a subset of $\Rd$ simply a \textdef{$\mu$-star} if it is a $(\mu,\thp)$-star at $\vect z$ for some $\thp>0$ and $\vect z\in\Rd$, and we
denote the \textdef{set of all $\mu$-stars at $\vect z$} and the \textdef{set of all $\mu$-stars}, respectively, by
\begin{equation}
\label{more_fat_stars_defs_eqn}
\normstars{\mu}{\vect z} \definedas \bigcup_{\thp>0} \fatstars{\mu}{\thp}{\vect z},
\quad\text{and}\quad
\normstars{\mu}{} \definedas \bigcup_{\vect z\in\Rd} \normstars{\mu}{\vect z}.
\end{equation}
Part~\ref{mu_cone_decomp:ball_part} of Lemma~\ref{mu_cone_decomp_lem} shows that any $\mu$-ball with center $\vect z$ is a $(\mu,\thp)$-star at $\vect z$ for every $\thp\ge 0$. Moreover, Part~\ref{mu_cone_decomp:cone_part} of Lemma~\ref{mu_cone_decomp_lem} shows that if $\thp'\ge \thp$, then every element of  $\conefamily{\mu}{\vect z,\thp'}$ is a $(\mu,\thp)$-star at $\vect z$. It then follows from the definition \eqref{fat_stars_def_eqn} that
\begin{equation}
\label{fatstars_monotone_eqn}
\thp'\ge \thp \implies
\fatstars{\mu}{\thp'}{\vect z} \subseteq \fatstars{\mu}{\thp}{\vect z},
\end{equation}
so the collections $\fatstars{\mu}{\thp}{\vect z}$ form an increasing family as $\thp\searrow 0$, and it is easily seen that this family is strictly increasing. Moreover, Lemma~\ref{mu_cone_decomp_lem} implies that $\mu$-bowling pins are $\mu$-stars, namely
\begin{equation}
\label{bpins_are_fatstars_eqn}
\bpin{\mu}{r}{\thp}{\vect z}{\dirvec}{h} \in \fatstars{\mu}{\thp}{\vect z}.
\end{equation}
This fact will be important in Chapter~\ref{random_fpc_chap}.
Finally, note that for any $\thp>0$,
\[
\fatstars{\mu}{\thp}{\vect z}
\subsetneq 
\normstars{\mu}{\vect z}
\subsetneq 
\fatstars{\mu}{0}{\vect z} =\stars{\vect z}.
\]
Star sets will play an important role in our study of the deterministic process, as we will see later in Chapter~\ref{deterministic_chap}. In Chapters~\ref{cone_growth_chap} and \ref{random_fpc_chap}, we will see that $\mu$-stars play an analogous role for the random process.

\section{Induced Intrinsic Metrics and the One-Type Deterministic Process}
\label{deterministic:one_type_sec}

\subsection{Induced Intrinsic Metrics in Subsets of $\Rd$}
\label{deterministic:intrinsic_metrics_sec}

A \textdef{path} in $\Rd$ is any continuous function $\gamma\colon I\to\Rd$, where $I\subseteq\R$ is an interval. A path is called \textdef{simple} if it is injective. We will sometime abuse terminology and identify a path with its image in $\Rd$. If $A,B\subseteq\Rd$, we say that a path $\gamma\colon [a,b]\to \Rd$ is a \textdef{path from $A$ to $B$} if $\gamma(a)\in A$ and $\gamma(b)\in B$. If $S$ is any subset of $\Rd$ and $I$ is an interval, we say that a path $\gamma\colon I\to\Rd$ is an \textdef{$S$-path} if either $\gamma$ is constant or $\gamma(I)\subseteq S$.
Given a norm $\mu$ on $\Rd$, let $\lengthfcn{\mu}$ denote the length operator on paths in $\Rd$ induced by the norm metric $\normmetric{\mu}$, as defined in \eqref{metric_length_def_eqn} in Section~\ref{metric_geom_sec}. Then for $A,B,S\subseteq\Rd$ (typically with $A,B\subseteq S$) we define the \textdef{$S$-restricted $\mu$-distance} from $A$ to $B$ as
\begin{equation}
\label{S_dist_def_eqn}
\idist{\mu}{S}{A}{B} \definedas
\inf \setbuilder[2]{\length{\mu}{\gamma}}{\gamma \text{ is an $S$-path from $A$ to $B$}}.
\end{equation}
If we take $A$ and $B$ to be singletons, then $\imetric{\mu}{S}$ defines a metric on $\Rd$, which is (by definition) finite on each ``accessibility component" of $S$, and is called the \textdef{intrinsic metric} on $S$ induced by $\normmetric{\mu}$ (cf.\ \cite[pp.~28--29]{Burago:2001aa}). We further define, for $S\subseteq \Rd$, an \textdef{$\starify{S}$-path} to be a path which is contained in $S$ except perhaps for its endpoints, i.e.\ any path $\gamma\colon I\to\Rd$ such that $\gamma(\interior{I})\subseteq S$.
We then define the \textdef{$\starify{S}$-restricted $\mu$-distance} from $A$ to $B$ as
\begin{equation}
\label{star_dist_def_eqn}
\stardist{\mu}{S}{A}{B} \definedas
\inf \setbuilder[2]{\length{\mu}{\gamma}}{\gamma \text{ is an $\starify{S}$-path from $A$ to $B$}}.
\end{equation}
The $\starify{S}$-restricted distance will be used to define the two-type deterministic process on $\Rd$ in the same way we used $\stargraph{S}$-restricted passage times to define the random two-type competition process on the lattice. A basic property of the intrinsic distance defined in \eqref{S_dist_def_eqn} is the following.

\begin{lem}[Intrinsic distance along a line segment]
\label{segment_idist_lem}
If $\vect x,\vect y\in\Rd$ and $\segment{\vect x}{\vect y}\subseteq {S}$, then $\idist{\mu}{S}{\vect x}{\vect y} = \distnew{\mu}{\vect x}{\vect y}$.
\end{lem}

\begin{proof}
The generalized triangle inequality for metric-induced length structures (\eqref{gen_triangle_ineq_eqn} in Section~\ref{metric_geom_sec}) implies that $\idist{\mu}{S}{\vect x}{\vect y}\ge \distnew{\mu}{\vect x}{\vect y}$ for any $\vect x,\vect y\in {S}$.
On the other hand, the simple path $\gamma$ defined on $[0,1]$ by $\gamma(t) = (1-t)\vect x + t\vect y$ has image $\segment{\vect x}{\vect y}\subseteq {S}$, so
\[
\idist{\mu}{S}{\vect x}{\vect y} \le \length{\mu}{\gamma} =  \distnew{\mu}{\vect x}{\vect y},
\]
where the final equality follows from the fact that straight lines are distance minimizing paths with respect to any norm metric on $\Rd$ (Lemma~\ref{norm_strictly_intrinsic_lem}).
\end{proof}

Note that Lemma~\ref{segment_idist_lem} implies that $\idist{\mu}{\Rd}{A}{B} = \distnew{\mu}{A}{B}$ for any $A,B\subseteq \Rd$. More generally, if $S\subseteq \Rd$ is convex, then $\imetric{\mu}{S} = \restrict{\normmetric{\mu}}{S\times S}$.

\subsection{The Restricted and Unrestricted One-Type Deterministic Processes}
\label{deterministic:one_type_def_sec}

Fix a norm $\mu$ on $\Rd$. For $A,S\subseteq \Rd$ (typically with $A\subseteq S$), we define the \textdef{$S$-restricted (deterministic) $\mu$-first-passage percolation process} started from $A$ by
\begin{equation}
\label{rreachedd_def_eqn}
\rreacheddset{\mu}{A}{S}{t} \definedas
\setbuilder[1]{\vect x\in S}{\idist{\mu}{S}{A}{\vect x} \le t}
\quad\text{for } t\ge 0,
\end{equation}
and we set $\rreacheddset{\mu}{A}{S}{\infty} \definedas \bigcup_{t\ge 0} \rreacheddset{\mu}{A}{S}{t}$. If $S=\Rd$, we drop the superscript $S$ from the notation, and we call $\reacheddset{\mu}{A}{t} \definedas \rreacheddset{\mu}{A}{\Rd}{t}$, $t\in [0,\infty]$, the \textdef{unrestricted $\mu$-first-passage percolation process}. Note that it follows from the definition that for any $A,B,S\subseteq \Rd$,
\[
\idist{\mu}{S}{A}{B} = \inf \setbuilder[1]{t\ge 0}{\rreacheddset{\mu}{A}{S}{t} \cap B \ne \emptyset},
\]
so $\idist{\mu}{S}{A}{B}$ is the \textdef{hitting time} of the set $B$ for the restricted growth process $\rreacheddset{\mu}{A}{S}{t}$.

With this interpretation of $\imetric{\mu}{S}$ as a hitting time, note that we can think of $\mu(\vect x) = \distnew{\mu}{\vect 0}{\vect x}$ as the time it takes an unrestricted process started from $\vect 0$ to reach the point $\vect x\in\Rd$. For this reason, and in order to draw an analogy with the traversal measure $\tmeasure$ from the random process, we will
refer to the norm $\mu$ as the \textdef{traversal norm} for the collection of deterministic processes $\setbuilder[2]{\rreacheddfcn{\mu}{A}{S}}{A,S\subseteq \Rd}$, since $\mu$ measures the time it takes the process to traverse any fixed distance in a given direction. Since $\mu$ is interpreted as measuring time, we can fix any reference norm $\norm{\cdot}$ on $\Rd$ to measure distances, and then the speed of the process in the direction $\dirvec$ is given by $\norm{\dirsymb}/\mu(\dirsymb)$.
%
%

\subsection{Unrestricted Growth and Scale-Invariant Starting Sets}
\label{deterministic:unrestricted_growth_sec}

The unrestricted deterministic growth process $\reacheddset{\mu}{A}{t}$ is rather simple.
Since $\imetric{\mu}{\Rd} = \normmetric{\mu}$, it follows from the definitions that for any $A\subseteq \Rd$ and $t\ge 0$,
\begin{equation}\label{deterministic:unrestricted_eqn}
\reacheddset{\mu}{A}{t} = \closure{A} + t \unitball{\mu},
\end{equation}
where $\unitball{\mu}$ is the closed unit ball for the norm $\mu$. That is, at time $t$, the set of conquered sites consists of the union of all $\mu$-balls of radius $t$ centered at some point in $\closure{A}$. Note that \eqref{deterministic:unrestricted_eqn} is consistent with our definition in Chapter~\ref{intro2_chap} of $\reacheddset{\mu}{\vect x}{r}$ as a $\mu$-ball of radius $r\ge 0$ centered at $\vect x\in \Rd$.

It follows from \eqref{deterministic:unrestricted_eqn} that if we take $A$ to be a scale-invariant subset of $\Rd$, then $\reacheddset{\mu}{A}{t}/t = \closure{A} + \unitball{\mu}$ for all $t>0$. Thus, for scale-invariant starting sets, as soon as some nonzero amount of time has elapsed, the shape of the conquered region remains the same after rescaling by $t$. This invariance property is one reason why scale-invariant starting sets will play a central role in our analysis of both the deterministic and random first-passage processes.

\subsection{Restricted Growth in Star-Shaped Regions}
\label{deterministic:star_set_growth_sec}

The evolution of the restricted deterministic process depends heavily on the choice of the restricting set $S$, and in general there is no simple formula such as \eqref{deterministic:unrestricted_eqn} describing the conquered region $\rreacheddset{\mu}{A}{S}{t}$ at time $t$. As an extreme example of how the process depends on $S$, take $S$ to be a totally disconnected set, e.g.\ $S = \rationals^d$. Then there are no nonconstant paths in $S$, and it follows that for any $A\subseteq \Rd$, we have $\rreacheddset{\mu}{A}{S}{t} = A\cap S$ for all $t$. On the other hand, even if we choose a ``nice" restricting set $S$, e.g.\ an open subset of $\Rd$ with smooth boundary, then unless $S$ is convex, it can be difficult to describe the evolution of the process explicitly because shortest paths in $S$ may deviate from straight lines.
%
%
However, for star-shaped regions as defined in \eqref{stars_at_A_def_eqn}, we have the following result, which generalizes \eqref{deterministic:unrestricted_eqn}:

\begin{prop}[Deterministic growth in star-shaped sets]
\label{determ_star_growth_prop}
Let $A\subseteq S\subseteq \Rd$, and suppose that $S\in \starscent{\bullet}{A}$, i.e.\ $S$ is star-shaped at each point of $A$.
Then
$
\rreacheddset{\mu}{A}{S}{t} = \parens{\closure{A} + t\unitball{\mu}} \cap S
$
for all $t\ge 0$.
\end{prop}

\begin{proof}

First note that the inclusion $\subseteq$ is trivial since by definition the process is monotone with respect to the restricting set. That is, for any $t\ge 0$,
\[
\rreacheddset{\mu}{A}{S}{t}
\subseteq \rreacheddset{\mu}{A}{\Rd}{t} \cap S
= \parens{\closure{A} + t\unitball{\mu}} \cap S,
\]
where the inclusion follows from definitions \eqref{S_dist_def_eqn} and \eqref{rreachedd_def_eqn}, and the final equality is from \eqref{deterministic:unrestricted_eqn}. To get the reverse inclusion, suppose $\vect x\in \parens{\closure{A} + t\unitball{\mu}} \cap S$. Then there is some $\vect z\in \closure{A}$ with $\distnew{\mu}{\vect z}{\vect x}\le t$. Let $\setof{\vect z_n}_{n\in\N}$ be a sequence of points in $A$ with $\vect z_n\to \vect z$. Since $S$ is star-shaped at each point of $A$, we have $\segment{\vect z_n}{\vect x} \subseteq S$ for all $n\in\N$. Then definition \eqref{S_dist_def_eqn}, Lemma~\ref{segment_idist_lem}, and the triangle inequality for $\normmetric{\mu}$ imply that
\begin{align*}
\idist{\mu}{S}{A}{\vect x} \le \inf_{n\in\N} \idist{\mu}{S}{\vect z_n}{\vect x}
&= \inf_{n\in\N} \distnew{\mu}{\vect z_n}{\vect x}
\le \distnew{\mu}{\vect z}{\vect x} \le t,
\end{align*}
so we have $\vect x\in \rreacheddset{\mu}{A}{S}{t}$ by definition \eqref{rreachedd_def_eqn}.
\end{proof}

One of the main results of Chapter~\ref{cone_growth_chap}, Theorem~\ref{cone_segment_growth_thm}, is an analogue of Proposition~\ref{determ_star_growth_prop} for the random one-type process.

\section{Definition and Basic Properties of the Two-Type Deterministic Process}
\label{deterministic:basic_properties_sec}

%

\subsection{Construction of the Process and Monotonicity}
\label{deterministic:two_type_construction_sec}

We now define the deterministic two-type competition process in a way that parallels the definition of the random two-type process in Section~\ref{fpc_basics:two_type_construction_sec}. 
%
%

\begin{definition}[Deterministic two-type first-passage competition]
\label{determ_two_type_proc_def}
Let $\initialset{1}$ and $\initialset{2}$ be disjoint subsets of $\Rd$, let $\normpair = \pair{\mu_1}{\mu_2}$ be a pair of norms on $\Rd$, called the \textdef{traversal norm pair} for the two-type process, and set
\begin{align*}
\dfinalset{1} &\definedas \bigcup \setbuilder[1]{S\subseteq \Rd\setminus \initialset{2}}
	{\forall  \vect x\in \bd S,\ \infty\ne \idist{\mu_1}{S}{\initialset{1}}{\vect x}
		\le \stardist{\mu_2}{(\Rd\setminus S)}{\initialset{2}}{\vect x}},\\
\dfinalset{2} &\definedas \bigcup \setbuilder[1]{S\subseteq \Rd\setminus \initialset{1}}
	{\forall  \vect x\in \bd S,\ \infty\ne \idist{\mu_2}{S}{\initialset{2}}{\vect x}
		\le \stardist{\mu_1}{(\Rd\setminus S)}{\initialset{1}}{\vect x}
	},
\end{align*}
if at least one of $\initialset{1}$, $\initialset{2}$ is nonempty, and $\dfinalset{1} = \dfinalset{2} \definedas \emptyset$ otherwise. We define the \textdef{two-type deterministic $\normpair$-first-passage competition process} with \textdef{starting configuration} $\initconfig$ to be the pair of $2^{\Rd}$-valued functions
\[
\reacheddset{\normpair}{\initconfig}{t}
= \pair[1]{\reacheddset{1}{\initconfig}{t}}{\reacheddset{2}{\initconfig}{t}}_{\normpair},
\]
defined for $t\in [0,\infty]$ by
\[
\reacheddset{1}{\initconfig}{t} \definedas \rreacheddset{\mu_1}{\initialset{1}}{\dfinalset{1}}{t}
\quad\text{and}\quad
\reacheddset{2}{\initconfig}{t} \definedas \rreacheddset{\mu_2}{\initialset{2}}{\dfinalset{2}}{t}
\quad\text{for } t\ge 0,
\]
and
\[
\reacheddset{i}{\initconfig}{\infty} \definedas \bigcup_{t\ge 0} \reacheddset{i}{\initconfig}{t}
\quad\text{for } i\in\setof{1,2}.
\]
We call the sets $\dfinalset{1}$ and $\dfinalset{2}$ the \textdef{finally conquered sets} for species 1 and species 2, respectively, and we call the process \textdef{full} if $\dfinalset{1}\cup \dfinalset{2} =\Rd$.
\end{definition}

\begin{thmremark}

In contrast to Definition~\ref{two_type_proc_def} of the random process, the finally conquered sets $\dfinalset{1}$ and $\dfinalset{2}$ in Definition~\ref{determ_two_type_proc_def} will typically have nonempty intersection because there will usually be some set of
\textdef{tie points} that both species reach at the same time, defined by
\[
\setof{\text{tie points}} \definedas
\setbuilder[1]{\vect x\in \dfinalset{1}\cup \dfinalset{2} }
{\idist{\mu_1}{\dfinalset{1}}{\initialset{1}}{\vect x} = \idist{\mu_2}{\dfinalset{2}}{\initialset{2}}{\vect x}}.
\]
Rather than mimicking our definition for $\Zd$-process and calling the $\Rd$-process ``well-defined" when $\dfinalset{1}\cap \dfinalset{2}= \emptyset$, we propose defining the $\Rd$-process to be \textdef{well-defined} if $\dfinalset{1}\cap \dfinalset{2} = \setof{\text{tie points}}$, though we will not use this definition anywhere. We do note however, that at least for some choices of the norms $\mu_1,\mu_2$, it is possible to construct examples where the process is not well-defined in this sense. It is also possible to construct deterministic processes that are not full, i.e.\ $\dfinalset{1}\cup \dfinalset{2} \ne\Rd$. This is accomplished by interlacing the initial sets $\initialset{1}$ and $\initialset{2}$ in such a way that prevents either species from ever reaching certain regions in $\Rd\setminus (\initialset{1}\cup\initialset{2})$; two such examples in $\R^2$ are an ``infinite yin-yang" spiral, in which neither species can escape a disc, or an infinitely stretched ``topologist's sine curve" $y=\parens[1]{\frac{1}{x}} \sin \parens[1]{\frac{1}{x}}$ $(x>0)$ interlaced with a ``topological comb" with increasingly long teeth suspended from the graph of $y=1+\frac{1}{x}$ $(x>0)$, in which neither species can escape the right half-plane.
\end{thmremark}

\begin{thmremark}
The finally conquered sets $\dfinalset{1}$ and $\dfinalset{2}$ in  Definition~\ref{determ_two_type_proc_def} are precisely the \textdef{crystal-growth Voronoi cells} for crystals growing deterministically according to the norms $\mu_1$, $\mu_2$ from the ``seed sets" $\initialset{1}$ and $\initialset{2}$. This generalizes the definition of ``multiplicatively weighted crystal-growth Voronoi diagrams" introduced by Schaudt and Drysdale \cite{Schaudt:1991aa} to the case where the crystals' norms are not necessarily multiples of each other. The finally conquered sets for first-passage competition between more than two species would coincide with the crystal-growth Voronoi diagram for a larger number of crystals.
\end{thmremark}

The following three lemmas are the $\Rd$-process analogues of Lemmas~\ref{final_sets_properties_lem}, \ref{entangled_disentangled_lem}, and \ref{trivial_init_configs_lem} for the $\Zd$-process; the proofs are left to the reader as they are very similar to the corresponding proofs in Chapter~\ref{fpc_basics_chap}.
Note, however, that Part~\ref{determ_final_sets:components_part} of Lemma~\ref{determ_final_sets_properties_lem} demonstrates the presence of topological nuances that arise in the $\Rd$-process but were irrelevant for the $\Zd$-process.

\begin{lem}[Properties of the finally conquered sets in the deterministic process]
\label{determ_final_sets_properties_lem}
For any traversal norm pair $\normpair$ and initial configuration $\initconfig$, the finally conquered sets $\dfinalset{1}$ and $\dfinalset{2}$ in Definition~\ref{determ_two_type_proc_def} satisfy the following for $i\in\setof{1,2}$.
\begin{enumerate}
\item \label{determ_final_sets:initial_sets_part}
$\dfinalset{i}\supseteq \initialset{i}$ and $ \dfinalset{i} \cap \initialset{3-i} = \emptyset$. Moreover, $\dfinalset{i}$ is empty if and only if $\initialset{i}$ is empty.

\item \label{determ_final_sets:characterization_part}
The following are equivalent.
\begin{enumerate}
\item $\vect x\in \dfinalset{i}$
\item $\idist{\mu_i}{\Rd\setminus \initialset{3-i}}{\initialset{i}}{\vect x} <\infty$ and $\idist{\mu_i}{\dfinalset{i}}{\initialset{i}}{\vect x} \le \stardist{\mu_{3-i}}{(\Rd\setminus \dfinalset{i})}{\initialset{3-i}}{\vect x}$.
\item $\infty \ne \idist{\mu_i}{\dfinalset{i}}{\initialset{i}}{\vect x} \le \stardist{\mu_{3-i}}{(\Rd\setminus \dfinalset{i})}{\initialset{3-i}}{\vect x}$.
\end{enumerate}

\item \label{determ_final_sets:components_part}
Every path component of $\dfinalset{i}$ in the intrinsic metric topology induced by $\imetric{\mu}{\dfinalset{i}}$ contains a path component of $\initialset{i}$ in the subspace topology inherited from $\Rd$.

\item \label{determ_final_sets:are_union_part}
$\dfinalset{i} = \reacheddset{i}{\initconfig}{\infty}_{\normpair}$.
\end{enumerate}
\end{lem}


\begin{lem}[Comparison of entangled and disentangled deterministic processes]
\label{determ_entangled_disentangled_lem}

Let $\normpair = \pair{\mu_1}{\mu_2}$ be a pair of norms on $\Rd$, let $\initconfig$ be an initial configuration in $\Rd$, and let $i\in\setof{1,2}$.
\begin{enumerate}
\item \label{determ_entangled_disentangled:unrestricted_part}
In the two-type process $\reacheddfcn{\normpair}{\initconfig}$, each species is dominated by its unrestricted disentangled version. That is, $\reacheddset{i}{\initconfig}{t} \subseteq \reacheddset{\mu_i}{\initialset{i}}{t}$ for all $t\in [0,\infty]$.

\item \label{determ_entangled_disentangled:restricted_part}
If species $i$ conquers some set $S\subseteq \Rd$ in the two-type process $\reacheddfcn{\normpair}{\initconfig}$, then
the entangled growth of species~$i$ dominates the growth of its disentangled version restricted to $S$. That is, $\reacheddset{i}{\initconfig}{t} \supseteq \rreacheddset{\mu_i}{\initialset{i}}{S}{t}$ for all $t\in [0,\infty]$.
\end{enumerate}
\end{lem}

\begin{lem}[Trivial starting configurations in the deterministic process]
\label{determ_trivial_init_configs_lem}

Let $\normpair = \pair{\mu_1}{\mu_2}$ be a pair of norms on $\Rd$, and let $\initconfig$ be an initial configuration in $\Rd$.
\begin{enumerate}
\item If $\initialset{i} = \emptyset$ for some $i\in\setof{1,2}$ and $j=3-i$, then
\[
\reacheddset{i}{\initconfig}{t} = \emptyset
\quad\text{and}\quad
\reacheddset{j}{\initconfig}{t} = \reacheddset{\mu_j}{\initialset{j}}{t}
\quad\text{ for all $t\ge 0$}.
\]

\item If $\initialset{1}\cup \initialset{2} = \Rd$, then $\reacheddset{\normpair}{\initconfig}{t} = \initconfig$ for all $t\ge 0$.
\end{enumerate}
Therefore,  if at least one of the sets $\initialset{1}$, $\initialset{2}$, or $\Rd\setminus \parens{\initialset{1}\cup \initialset{2}}$ is empty, then the two-type deterministic process with initial configuration $\initconfig$ is well-defined and full for all traversal norm pairs $\normpair$.
\end{lem}

The following proposition is the deterministic analogue of Proposition~\ref{conquering_property_prop}.

\begin{prop}[Conquering property for the two-type deterministic process]
\label{determ_conquering_property_prop}
Consider a deterministic two-type competition process $\reacheddfcn{\normpair}{\initconfig}$, with nonempty starting configuration $\initconfig$ and traversal norm pair $\normpair = \pair{\mu_1}{\mu_2}$. If $S$ is a subset of $\Rd$ such that 
\begin{equation}
\label{determ_conquering_property_eqn}
S\cap \initialset{2} = \emptyset
\quad\text{and}\quad
\infty \ne \idist{\mu_1}{S}{\initialset{1}}{\vect x} \le
\stardist{\mu_2}{(\Rd\setminus S)}{\initialset{2}}{\vect x}
\quad\text{for all } \vect x\in \bd S,
\end{equation}
then species 1 conquers $S$ in the $\normpair$-process, i.e.\ $\reacheddset{1}{\initconfig}{\infty} _{\normpair} \supseteq S$. A symmetric statement holds with the roles of species 1 and 2 interchanged.
\end{prop}

\begin{proof}
If $S\subseteq \Rd$ satisfies the stated hypotheses, then it follows directly from Definition~\ref{determ_two_type_proc_def} that $S\subseteq \dfinalset{1}$, and by Lemma~\ref{determ_final_sets_properties_lem} we have $\dfinalset{1} = \reacheddset{1}{\initconfig}{\infty}_{\normpair}$.
\end{proof}

Analogously to the lattice process, we call a set $S$ that satisfies \eqref{determ_conquering_property_eqn} a \textdef{conquering set for species~1}, or a \textdef{$\dfinalset{1}$-set} for the $\normpair$-process, and similarly for species~2.
The following lemma about conquering sets is another example of the topological differences between the competition process in $\Rd$ versus $\Zd$.

\begin{lem}[Without loss of generality, conquering sets are closed]
\label{conq_sets_closed_lem}
If $S\subseteq \Rd$ is a conquering set for species~$i$ in the $\normpair$-process, then $\bd S\subseteq \initialset{i}\cup S$, and $\closure{S}$ is also a conquering set for species~$i$.
\end{lem}

\begin{proof}
Suppose $S$ is a conquering set for species $i$. Then by definition, for all $\vect x\in \bd S$ we have $\idist{\mu_i}{S}{\initialset{i}}{\vect x} <\infty$. By the definition of the intrinsic metric $\imetric{\mu_i}{S}$, this means that there is an $S$-path $\gamma$ from $\initialset{i}$ to $\vect x$, and the definition of an $S$-path then implies that either $\vect x\in \initialset{i}$ or $\vect x\in S$. Thus, $\bd S\subseteq \initialset{i}\cup S$ and hence $\closure{S}\subseteq \initialset{i}\cup S$. Since $\initialset{1} \cap \initialset{2} = \emptyset$ and $S\cap \initialset{3-i} = \emptyset$ by assumption, we thus have $\closure{S} \cap \initialset{3-i} = \emptyset$, and since $\bd \closure{S} \subseteq \bd S$ for any $S\subseteq\Rd$, we have
\[
\forall \vect x\in \bd \closure{S}, \quad
\infty \ne \idist{\mu_i}{\closure{S}}{\initialset{1}}{\vect x}
\le \idist{\mu_i}{S}{\initialset{1}}{\vect x}
\le \stardist{\mu_{3-i}}{(\Rd\setminus S)}{\initialset{2}}{\vect x}
\le \stardist{\mu_{3-i}}{(\Rd\setminus \closure{S})}{\initialset{2}}{\vect x}
\]
because $S$ is a conquering set, so $\closure{S}$ is also a conquering set.
\end{proof}

Note that although Lemma~\ref{conq_sets_closed_lem} shows that conquering sets must be closed, the finally conquered sets $\dfinalset{1}$ and $\dfinalset{2}$ need not be (for example if $\closure{\initialset{i}}\cap \initialset{3-i} \ne \emptyset$).

The following monotonicity property is a deterministic analogue of Lemma~\ref{monotonicity_lem}.

\begin{lem}[Monotonicity of infection for the deterministic process]
\label{determ_monotonicity_lem}
Let $\normpair = \pair{\mu_1}{\mu_2}$ and $\normpair' = \pair{\mu_1'}{\mu_2'}$ be traversal norm pairs, and let $\pair{A_1}{A_2}$ and $\pair{A_1'}{A_2'}$ be starting configurations. Suppose that $\mu_1\le \mu_1'$ and $\mu_2 \ge \mu_2'$, and that $A_1\supseteq A_1'$ and $A_2\subseteq A_2'$. Then
\[
\reacheddset{1}{\pair{A_1}{A_2}}{t} \supseteq \reacheddset{1}{\pair{A_1'}{A_2'}}{t}
\quad\text{and}\quad
\reacheddset{2}{\pair{A_1}{A_2}}{t} \subseteq \reacheddset{2}{\pair{A_1'}{A_2'}}{t}
\quad\text{for all } t\ge 0.
\]
More succinctly, this says that the deterministic process is monotone: If $\normpair \ge \normpair'$ and $\pair{A_1}{A_2}\ge \pair{A_1'}{A_2'}$, then $\reacheddfcn{\normpair}{\pair{A_1}{A_2}}\ge \reacheddfcn{\normpair'}{\pair{A_1'}{A_2'}}$, where the orderings are chosen to favor species 1.
\end{lem}

\subsection{Weighted Voronoi Cells and Star-Cells}
\label{deterministic:closer_sets_sec}

The finally conquered sets (a.k.a.\ crystal-growth Voronoi cells) $\dfinalset{1}$ and $\dfinalset{2}$ from Definition~\ref{determ_two_type_proc_def} are in general not easy to describe explicitly, due to the ``wrap-around" effect exhibited by the growth of the two species (or crystals) in the competition process. As a simpler, first approximation to the crystal-growth Voronoi cells, we can start with the ordinary Voronoi cells, which consist of the points in $\Rd$ that are closer to one species' starting set than the other. Here the word ``closer" means that we are comparing the (unrestricted) traversal times of the two species to a given point, i.e.\ we measure distances from $\initialset{1}$ using $\normmetric{\mu_1}$ and distances from $\initialset{2}$ using $\normmetric{\mu_2}$. This interpretation is more general than the usual definition of Voronoi cells using a single distance function. In fact, we will consider the even more general \emph{multiplicatively weighted Voronoi cells}, in which one species' distance function is adjusted by a multiplicative constant. This flexibility in perturbing distances by a multiplicative factor will be only marginally important in our analysis of the deterministic process, but will be essential for adapting the proofs in the present chapter to the random process in Chapter~\ref{random_fpc_chap}. We now proceed with the relevant definitions.

Let $\normpair = \pair{\mu_1}{\mu_2}$ be a pair of norms on $\Rd$, and let $\initialset{1}$ and $\initialset{2}$ be disjoint subsets of $\Rd$. For $\beta\ge 0$ and $i\in\setof{1,2}$, we define the \textdef{$\beta$-weighted Voronoi cell} (or simply \textdef{$\beta$-Voronoi cell}) for species~$i$ in the $\normpair$-process with starting configuration $\initconfig$ to be
\begin{equation}
\label{closerset_def_eqn}
\closerset{\normpair}{\beta}{i}{\initialset{1}}{\initialset{2}} \definedas
\setbuilder[1]{\vect x\in\Rd\setminus \initialset{3-i}}
	{\beta\distnew{\mu_i}{\initialset{i}}{\vect x} \le
	\distnew{\mu_{3-i}}{\initialset{3-i}}{\vect x}}.
\end{equation}
For example, species~1's $\beta$-Voronoi cell is the set of $\vect x\in\Rd\setminus \initialset{2}$ which are closer to $\initialset{1}$ than to $\initialset{2}$ by a factor of at least $\beta$, meaning that either the initial distance from species~$1$ to $\vect x$ is 0 or else the ratio $\frac{\distnew{\mu_{2}}{\initialset{2}}{\vect x}}{\distnew{\mu_1}{\initialset{1}}{\vect x}}$ is at least $\beta$.
When we study the random process in Chapter~\ref{random_fpc_chap}, we will typically want $\beta>1$ so that the points in species $i$'s $\beta$-Voronoi cell are strictly closer to species $i$ by some fixed ratio, giving us some wiggle room to account for random fluctuations.
By contrast, if $\beta<1$, then some points in species $i$'s $\beta$-Voronoi cell may be farther from $\initialset{i}$ than from $\initialset{3-i}$, but only by a factor of at most $\beta^{-1}$.

The case $\beta=1$ in \eqref{closerset_def_eqn} corresponds to the ordinary unweighted Voronoi cells induced by the pair of distance functions $\normmetric{\mu_1}$ and $\normmetric{\mu_2}$ for the seed sets $\initialset{1}$ and $\initialset{2}$. When $\mu_1=\mu_2$, the 1-Voronoi cells coincide precisely with the crystal-growth Voronoi cells  $\dfinalset{1}$ and $\dfinalset{2}$ from Definition~\ref{determ_two_type_proc_def}. On the other hand, for a general pair of norms $\normpair$, even norms that are multiples of one another, $\closerset{\normpair}{1}{i}{\initialset{1}}{\initialset{2}}$ neither contains nor is contained in the finally conquered set $\dfinalset{i}$, because definition \eqref{closerset_def_eqn} does not take into account the effect of the two species blocking each other's paths.
However, as long as the species are able to travel in straight lines without interference, this blocking effect does not come into play. For example, if we know that every point in a line segment originating from some point $\vect z\in \initialset{1}$ is closer to the point $\vect z$ than it is to any point in $\initialset{2}$, then it is intuitively clear that this line segment will be conquered by species~1. We formalize this idea in the following definition of ``Voronoi star-cells," which provide a second-step approximation of the sets $\dfinalset{1}$ and $\dfinalset{2}$. Namely, we will show in Proposition~\ref{determ_star_conquering_prop} below that species~$1$'s 1-Voronoi star-cell is contained in the set $\dfinalset{1}$, and for any $\beta>1$, species~2's $\beta$-Voronoi star-cell is contained in the complement of $\dfinalset{1}$.



With $\normpair$ and $\initconfig$ as above, we define the \textdef{$\beta$-weighted Voronoi star-cells} (or \textdef{$\beta$-Voronoi star-cells}) for the $\normpair$-competition process started from $\initconfig$ to be
\begin{equation}
\label{cstar_def_eqn}
\begin{split}
\cstarset{\normpair}{\beta}{1}{\initialset{1}}{\initialset{2}}
&\definedas \setbuilder[1]{\vect x\in \Rd}{\exists \vect z\in \initialset{1}
	\text{ such that } \segment{\vect z}{\vect x} \subseteq 
	\closerset{\normpair}{\beta}{1}{\vect z}{\initialset{2}} },\\
\cstarset{\normpair}{\beta}{2}{\initialset{1}}{\initialset{2}}
&\definedas \setbuilder[1]{\vect x\in \Rd}{\exists \vect z\in \initialset{2}
\text{ such that } \segment{\vect z}{\vect x} \subseteq 
\closerset{\normpair}{\beta}{2}{\initialset{1}}{\vect z} }.
\end{split}
\end{equation}
More explicitly,
\begin{equation*}
\cstarset{\normpair}{\beta}{i}{\initialset{1}}{\initialset{2}}=
\begin{cases}
\setbuilderbar[3]{\vect x\in\Rd}
	{\exists \vect z\in\initialset{i} \text{ such that }
	\segment{\vect z}{\vect x}\cap \initialset{3-i} = \emptyset }
	&\text{if } \beta=0,\\[4mm]
\setbuilderbar[1]{\vect x\in\Rd}{ \displaystyle
	\mathstack{ \exists \vect z\in\initialset{i} \text{ such that }
		\forall \vect y\in\segment{\vect z}{\vect x},}
	{\beta \distnew{\mu_i}{\vect z}{\vect y}
		\le \distnew{\mu_{3-i}}{\initialset{3-i}}{\vect y}}}
	&\text{if } \beta>0.
\end{cases}
\end{equation*}
It is clear from definition \eqref{cstar_def_eqn} and the monotonicity of Voronoi cells with respect to their seed sets that
\[
\cstarset{\normpair}{\beta}{i}{\initialset{1}}{\initialset{2}} \subseteq
\closerset{\normpair}{\beta}{i}{\initialset{1}}{\initialset{2}}
\quad\text{for $i\in\setof{1,2}$.}
\]
The name ``star cell" refers to the fact that $\cstarset{\normpair}{\beta}{i}{\initialset{1}}{\initialset{2}}$ is a union of star-shaped sets centered in $\initialset{i}$. That is, $\cstarset{\normpair}{\beta}{i}{\initialset{1}}{\initialset{2}}\in \genstars{}{\initialset{i}}$, where $\genstars{}{\initialset{i}}$ is the collection of ``generalized star sets at $\initialset{i}$" defined in \eqref{genstars_at_A_def_eqn}. 
Note, however, that $\cstarset{\normpair}{\beta}{i}{\initialset{1}}{\initialset{2}}$ is not necessarily star-shaped unless $\initialset{i}$ consists of a single point. Our main interest in the Voronoi star-cells is the above-mentioned intuition that the set $\cstarset{\normpair}{1}{i}{\initialset{1}}{\initialset{2}}$ should get conquered by species~$i$ in the $\normpair$-process started from $\initconfig$. We will prove this property in Proposition~\ref{determ_star_conquering_prop} below. First we make one more observation.

Note that the $\beta$-Voronoi cells and star-cells are decreasing in $\beta$, so we can define the right-hand and left-hand limits of the cells as $\beta'\searrow \beta$ or $\beta'\nearrow \beta$, by taking a union or intersection respectively. With these definitions, one can check that for all $\beta\ge 0$,
\begin{align}
\label{star_cell_plus_eqn}
\closerset{\normpair}{\beta^+}{i}{\initialset{1}}{\initialset{2}}
&= \setbuilderbar[1]{\vect x\in\Rd\setminus \initialset{3-i}}
	{\vect x\in \closure{\initialset{i}} \text{ or }
	\beta\distnew{\mu_i}{\initialset{i}}{\vect x} <
	\distnew{\mu_{3-i}}{\initialset{3-i}}{\vect x}}, \notag\\[2mm]
\cstarset{\normpair}{\beta^+}{i}{\initialset{1}}{\initialset{2}}
&= \setbuilderbar[1]{\vect x\in\Rd}{
	\mathstack{\vect x\in \initialset{i} \text{ or }
		\exists \vect z\in\initialset{i} \text{ such that }
		\forall \vect y\in\segment{\vect z}{\vect x},}
	{\beta \distnew{\mu_i}{\vect z}{\vect y}
		< \distnew{\mu_{3-i}}{\initialset{3-i}}{\vect y}}},
\end{align}
and
$
\closerset{\normpair}{\beta^-}{i}{\initialset{1}}{\initialset{2}} 
= \closerset{\normpair}{\beta}{i}{\initialset{1}}{\initialset{2}}
$
for all $\beta>0$. The corresponding description of $\cstarset{\normpair}{\beta^-}{i}{\initialset{1}}{\initialset{2}}$ is more complicated and is therefore omitted.

\begin{prop}[Conquering a union of star-shaped sets in the deterministic process]
\label{determ_star_conquering_prop}
In the deterministic $\normpair$-competition process started from $\initconfig$, species $i$ conquers every point in $\cstarset{\normpair}{1}{i}{\initialset{1}}{\initialset{2}}$ and no point in $\cstarset{\normpair}{1^+}{3-i}{\initialset{1}}{\initialset{2}}$.
\end{prop}

\begin{proof}
Take $i=1$ for concreteness. For the first statement, let $\vect x\in \cstarset{\normpair}{1}{1}{\initialset{1}}{\initialset{2}}$. Then there is some $\vect z\in \initialset{1}$ such that $\vect x \in \cstarset{\normpair}{1}{1}{\vect z}{\initialset{2}}$ and hence $\segment{\vect z}{\vect x}\subseteq \closerset{\normpair}{1}{1}{\vect z}{\initialset{2}}$, by definition \eqref{cstar_def_eqn}. We will show that the set $S= \segment{\vect z}{\vect x}$ satisfies the hypotheses of Proposition~\ref{determ_conquering_property_prop} and hence species 1 conquers $S$.

Let $\vect y\in \bd S = S$. Then since $\segment{\vect z}{\vect y}\subseteq S$, Lemma~\ref{segment_idist_lem} implies that
\begin{equation}
\label{determ_star_conq:dist_A1_y_eqn}
\idist{\mu_1}{S}{\initialset{1}}{\vect y} \le \idist{\mu_1}{S}{\vect z}{\vect y}
=\distnew{\mu_1}{\vect z}{\vect y}.
\end{equation}
In particular, \eqref{determ_star_conq:dist_A1_y_eqn} implies that $\idist{\mu_1}{S}{\initialset{1}}{\vect y}<\infty$. Also, since $S \subseteq \closerset{\normpair}{1}{1}{\vect z}{\initialset{2}}$, we have $S\cap \initialset{2} = \emptyset$ by the definition \eqref{closerset_def_eqn}. Finally, \eqref{determ_star_conq:dist_A1_y_eqn} and \eqref{closerset_def_eqn}, together with the monotonicity of the star intrinsic metric with respect to the restricting set, imply that
\[
\idist{\mu_1}{S}{\initialset{1}}{\vect y} \le \distnew{\mu_1}{\vect z}{\vect y} 
\le \distnew{\mu_2}{\initialset{2}}{\vect y}
\le \stardist{\mu_2}{(\Rd\setminus S)}{\initialset{2}}{\vect y}.
\]
Therefore, the set $S$ satisfies the three hypotheses in Proposition~\ref{determ_conquering_property_prop}, so species 1 conquers the segment $S= \segment{\vect z}{\vect x}$.
Since $\vect x$ was arbitrary, species 1 conquers all of $\cstarset{\normpair}{1}{1}{\initialset{1}}{\initialset{2}}$.

\extraline

\newcommand{\sset}{S_{\vect x}}
\newcommand{\bdpt}{\vect y_0}


For the second statement, let $\vect x\in \cstarset{\normpair}{1^+}{2}{\initialset{1}}{\initialset{2}}$.
Then by \eqref{star_cell_plus_eqn} there is some $\vect z\in \initialset{2}$ such that 
\begin{equation}
\label{determ_star_conq:dist_z_y_eqn}
\distnew{\mu_2}{\vect z}{\vect y} < \distnew{\mu_1}{\initialset{1}}{\vect y}
\quad\text{for all $\vect y\in \segment{\vect z}{\vect x}$.}
\end{equation}
For a contradiction, suppose $\vect x\in \reacheddset{1}{\initconfig}{\infty} = \dfinalset{1}$. Then by Lemma~\ref{conq_sets_closed_lem}, $\vect x$ is contained in some closed conquering set $\sset$ for species 1. Let $\bdpt$ be a $\mu_2$-closest point in $\segment{\vect z}{\vect x}\cap \sset$ to $\vect z$.
Then by Lemma~\ref{closest_bd_point_lem}, we have $\bdpt\in \bd \sset$ and $\cosegment{\vect z}{\bdpt} \cap \sset = \emptyset$. Thus, $\segment{\vect z}{\bdpt}$ is a $\starify[1]{\Rd\setminus \sset}$-path from $\initialset{2}$ to $\bdpt$, so (using Lemma~\ref{segment_idist_lem})
\begin{equation}
\label{determ_star_conq:star_path_eqn}
\stardist{\mu_2}{(\Rd\setminus \sset)}{\initialset{2}}{\bdpt}
\le \idist{\mu_2}{\segment{\vect z}{\bdpt}}{\vect z}{\bdpt}
= \distnew{\mu_2}{\vect z}{\bdpt}.
\end{equation}
Now, since $\sset$ is a $\dfinalset{1}$-set and $\bdpt\in \bd \sset$, we have
\begin{equation}
\label{determ_star_conq:conq_set_eqn}
\infty \ne \idist{\mu_1}{\sset}{\initialset{1}}{\bdpt}
\le \stardist{\mu_2}{(\Rd\setminus \sset)}{\initialset{2}}{\bdpt}.
\end{equation}
Combining \eqref{determ_star_conq:star_path_eqn}, \eqref{determ_star_conq:conq_set_eqn}, and the fact that the intrinsic metric in $\sset$ is at least as large as the unrestricted norm metric, we have
\[
\distnew{\mu_1}{\initialset{1}}{\bdpt}\le 
\idist{\mu_1}{\sset}{\initialset{1}}{\bdpt} 
\le \distnew{\mu_2}{\vect z}{\bdpt}.
\]
But since $\bdpt\in \segment{\vect z}{\vect x}$, this inequality contradicts \eqref{determ_star_conq:dist_z_y_eqn}. Thus we conclude that $\vect x\notin \reacheddset{1}{\initconfig}{\infty}$, so species~1 conquers no point in $\cstarset{\normpair}{1^+}{2}{\initialset{1}}{\initialset{2}}$.
\end{proof}

Proposition~\ref{determ_star_conquering_prop} will be our primary tool for describing the finally conquered sets $\dfinalset{1}$ and $\dfinalset{2}$ when we study competition in cones in Section~\ref{deterministic:cone_competition_sec}. One of our main goals in Chapter~\ref{random_fpc_chap} will be to prove an analogue of Proposition~\ref{determ_star_conquering_prop} for the random two-type process. The main general result of this sort is Theorem~\ref{cone_segment_conquering_thm} (and Corollary~\ref{star_point_conquering_cor}), which corresponds to the case $\initialset{1} = \setof{\vect z}$ and will be our primary tool for analyzing random first-passage competition in cones in Section~\ref{random_fpc:cone_competition_sec}. Corollary~\ref{star_conquering_cor} of Theorem~\ref{cone_segment_conquering_thm} is a general analogue of Proposition~\ref{determ_star_conquering_prop} in the case where $\initialset{1}$ is bounded.
Note that the converse of Proposition~\ref{determ_star_conquering_prop} is false in general, i.e.\ species~$i$ may conquer points outside of $\cstarset{\normpair}{1}{i}{\initialset{1}}{\initialset{2}}$ by following a curved path instead of a straight line, or by blocking the other species' shortest path to these points.

The following lemma allows us to convert between statements about the $\beta$-Voronoi cells or star-cells for different $\beta$'s, by scaling one or both of the traversal norms by an appropriate factor.

\begin{lem}[Transformation of the $\beta$-Voronoi cells under speed scaling]
\label{closer_set_scaling_lem}
Let $\mu_1$ and $\mu_2$ be two norms on $\Rd$, and consider the traversal norm pairs $\normpair = \pair{\mu_1}{\mu_2}$ and $\normpair' = \pair{\alpha_1 \mu_1}{\alpha_2 \mu_2}$, where $\alpha_1,\alpha_2>0$. Then for any starting configuration $\initconfig$, $i\in\setof{1,2}$, and $\beta,\beta'>0$, the following are equivalent.
\begin{enumerate}
\item $\displaystyle  \frac{\alpha_1}{\alpha_2} = \frac{\beta}{\beta'}$.

\item $\closerset{\normpair}{\beta}{i}{\initialset{1}}{\initialset{2}} =
\closerset{\normpair'}{\beta'}{i}{\initialset{1}}{\initialset{2}}$.

\item $\cstarset{\normpair}{\beta}{i}{\initialset{1}}{\initialset{2}} =
\cstarset{\normpair'}{\beta'}{i}{\initialset{1}}{\initialset{2}}$.
\end{enumerate}
\end{lem}

\begin{proof}
Since the distance functions for the respective norms get scaled by $\alpha_1$ and $\alpha_2$, it is evident from the definitions \eqref{closerset_def_eqn} and \eqref{cstar_def_eqn} that the $\beta$-Voronoi cells and $\beta$-Voronoi star-cells get scaled accordingly.
\end{proof}

Lemma~\ref{closer_set_scaling_lem} implies that the $\beta$-Voronoi (star-)cells can be interpreted as the $1$-Voronoi (star-)cells in a process in which one species' speed has been scaled by $\beta$, as illustrated by the following.

\begin{cor}[$\beta$-Voronoi star-cells and $\beta$-adjusted processes]
\label{beta_adjusted_conq_cor}

\newcommand{\pairone}{\pair{\beta \mu_1}{\mu_2}}
\newcommand{\pairtwo}{\pair{\mu_1}{\beta^{-1}\mu_2}}

Let $\normpair = \pair{\mu_1}{\mu_2}$ be a traversal norm pair, let $\initconfig$ be an initial configuration in $\Rd$, and let $\beta>0$. Then in either of the ``$\beta$-adjusted" deterministic processes $\reacheddset{\pairone}{\initconfig}{t}$ or $\reacheddset{\pairtwo}{\initconfig}{t}$, species 1 conquers every point in $\cstarset{\normpair}{\beta}{1}{\initialset{1}}{\initialset{2}}$ and no point in $\cstarset{\normpair}{\beta^+}{2}{\initialset{1}}{\initialset{2}}$.
\end{cor}

\begin{proof}
Combine Proposition~\ref{determ_star_conquering_prop}  and Lemma~\ref{closer_set_scaling_lem}.
\end{proof}

The following lemma gives a rather obvious but useful characterization of the $\beta$-Voronoi cells. We call this the ``dual characterization" due to its similarity to Lemma~\ref{dual_times_lem} in Section~\ref{random_fpc:conq_upper_bds_sec}; the ``duality" refers to a change in perspective from a $\mu_2$-process originating in some arbitrary set $\initialset{2}$ to a ``dual" $\mu_2$-process originating at a single point $\vect x$.

\begin{lem}[Dual characterization of the $\beta$-Voronoi cells]
\label{closer_char_lem}
Let $\initconfig$ be a starting configuration and $\normpair = \pair{\mu_1}{\mu_2}$ a traversal norm pair. Then for any $\beta\ge 0$ and $\vect x\in\Rd$,
\[
\vect x\in \closerset{\normpair}{\beta}{1}{\initialset{1}}{\initialset{2}}
\iff
\exists r\ge \beta \distnew{\mu_1}{\initialset{1}}{\vect x}
\text{ such that }
\reacheddset{\mu_2}{\vect x}{r-} \cap \initialset{2} = \emptyset,
\]
where $\reacheddset{\mu_2}{\vect x}{0-}$ is interpreted as $\setof{\vect x}$. A symmetric statement holds with species 1 and 2 switched.
\end{lem}

\begin{proof}
Let $\vect x\in\Rd$. First suppose $\reacheddset{\mu_2}{\vect x}{r-} \cap \initialset{2} = \emptyset$ for some $r\ge \beta \distnew{\mu_1}{\initialset{1}}{\vect x}$. Then for every $\vect y\in \initialset{2}$ we have $\distnew{\mu_2}{\vect x}{\vect y} \ge r \ge \beta \distnew{\mu_1}{\initialset{1}}{\vect x}$, so $\vect x\in \closerset{\normpair}{\beta}{1}{\initialset{1}}{\initialset{2}}$. On the other hand, suppose  $\reacheddset{\mu_2}{\vect x}{r-} \cap \initialset{2} \ne \emptyset$, where $r=  \beta \distnew{\mu_1}{\initialset{1}}{\vect x}$. Then there exists $\vect y\in \initialset{2}$ such that $\distnew{\mu_2}{\vect x}{\vect y} < r=  \beta \distnew{\mu_1}{\initialset{1}}{\vect x}$. Thus, $\distnew{\mu_2}{\vect x}{\initialset{2}} < \beta \distnew{\mu_1}{\initialset{1}}{\vect x}$, so $\vect x \not\in \closerset{\normpair}{\beta}{1}{\initialset{1}}{\initialset{2}}$. (Note: The special cases $\beta=0$ and $r=0$ should be handled separately.)
\end{proof}

Using Lemma~\ref{closer_char_lem}, we can prove the following similar ``dual" characterization of the $\beta$-Voronoi star-cells, which shows that the line segment $\segment{\vect z}{\vect x}$ is contained in $\closerset{\normpair}{\beta}{1}{\vect z}{A_2}$ for some $\beta>0$ if and only if $\initialset{2}$ does not intersect a particular nondegenerate $\mu_2$-cone segment with apex $\vect z$ and axis $\segment{\vect z}{\vect x}$. To state the result,
define
\begin{equation}
\label{thpadv_def_eqn}
\thpadv{2}{\dirvec}{\beta} \definedas \beta \frac{\mu_1(\dirsymb)}{\mu_2(\dirsymb)}
\quad\text{and}\quad
\thpadv{1}{\dirvec}{\beta} \definedas \beta \frac{\mu_2(\dirsymb)}{\mu_1(\dirsymb)}
\quad\text{for $\beta>0$ and $\dirsymb\in\Rd\setminus \setof{\vect 0}$.}
\end{equation}
We call the constants $\thpadv{2}{\dirvec}{\beta}$ and $\thpadv{1}{\dirvec}{\beta}$ defined in \eqref{thpadv_def_eqn} the \textdef{$\beta$-separation thickness} for a $\mu_2$-cone or $\mu_1$-cone, respectively, in the direction $\dirvec$. These constants will appear in numerous results about the two-type process throughout the remainder of Chapter~\ref{deterministic_chap} and Chapter~\ref{random_fpc_chap}.

\begin{lem}[Dual characterization of the $\beta$-Voronoi star-cells]
\label{star_closer_char_lem}
\newcommand{\thisx}{\vect x}
\newcommand{\thisdir}{\dirsymb}
\newcommand{\thisdirvec}{\dirvec}
\newcommand{\thisthp}{\thpadv{2}{\thisdirvec}{\beta}}

For any $\beta>0$, any $\vect x,\vect z\in\Rd$ with $\vect x\ne\vect z$, and any $\initialset{2}\subseteq \Rd\setminus \setof{\vect z}$, we have
\[
\thisx\in \cstarset{\normpair}{\beta}{1}{\vect z}{\initialset{2}}
\iff
\interior{\conetip{\mu_2}{\thisthp}{\vect z}{\thisdirvec}{h}} \cap \initialset{2} = \emptyset,
\text{ where $\thisdir = \thisx-\vect z$ and $h=\mu_2(\thisdir)$,} 
\]
and $\thisthp$ is defined by \eqref{thpadv_def_eqn}. A symmetric statement holds with species~1 and 2 switched.
\end{lem}

\begin{proof}
\newcommand{\segpt}{\vect y}
\newcommand{\xpt}{\vect x}
\newcommand{\zpt}{\vect z}
\newcommand{\thisthp}{\thpadv{2}{\dirvec}{\beta}}


First note that for any $\segpt\in \ocsegment{\zpt}{\xpt}$,
\begin{align}
\label{star_closer_char:ratio_eqn}
\beta \distnew{\mu_1}{\zpt}{\segpt}
&=\beta \distnew{\mu_2}{\zpt}{\segpt}
	\cdot \frac{\distnew{\mu_1}{\zpt}{\segpt}}{\distnew{\mu_2}{\zpt}{\segpt}}
	\notag\\
&=\beta \distnew{\mu_2}{\zpt}{\segpt}
	\cdot \frac{\distnew{\mu_1}{\zpt}{\xpt}}{\distnew{\mu_2}{\zpt}{\xpt}}
=\thisthp \distnew{\mu_2}{\zpt}{\segpt}.
\end{align}
Now, since $\zpt\not\in \initialset{2}$ we trivially have $\zpt\in \closerset{\normpair}{\beta}{1}{\zpt}{\initialset{2}}$, so
\begin{align*}
\xpt \in \cstarset{\normpair}{\beta}{1}{\zpt}{\initialset{2}}
&\ifandonlyif \forall \segpt\in \ocsegment{\zpt}{\xpt},\
	\segpt\in \closerset{\normpair}{\beta}{1}{\zpt}{\initialset{2}}
	&&\text{(by definition)}\\
&\ifandonlyif \forall \segpt\in \ocsegment{\zpt}{\xpt},\
	\reacheddset[2]{\mu_2}{\segpt}{\beta \distnew{\mu_1}{\zpt}{\segpt}-}
	\cap \initialset{2} = \emptyset
	&&\text{(by Lemma~\ref{closer_char_lem})}\\
&\ifandonlyif \bigcup_{\segpt\in \ocsegment{\zpt}{\xpt}}
	\reacheddset[2]{\mu_2}{\segpt}{\thisthp \distnew{\mu_2}{\zpt}{\segpt}-}
	\cap \initialset{2} = \emptyset
	&&\text{(by \eqref{star_closer_char:ratio_eqn})}\\
&\ifandonlyif \interior{\conetip{\mu_2}{\thisthp}{\vect z}{\dirvec}{h}}
	\cap \initialset{2} = \emptyset.
	&&\text{(by definition)}
	\qedhere
\end{align*}
\end{proof}

\subsection{Scale-Invariant Starting Configurations}
\label{deterministic:scale_inv_configs_sec}


The following lemma shows that scale-invariance is preserved by the two-type process.

\begin{lem}[Scale-invariant starting configurations]
\label{scale_inv_competition_lem}
If the initial configuration $\initconfig$ is scale-invariant (i.e.\ both $\initialset{1}$ and $\initialset{2}$ are scale-invariant), then the finally conquered regions are both scale-invariant, and $\frac{1}{t} \reacheddset{\normpair}{\initconfig}{t}$ is constant for $t>0$.
\end{lem}

\begin{proof}
This follows from the definition of the finally conquered sets and scale-equivariance of the norm metrics.
\end{proof}

For a scale-invariant starting configuration with $\initialset{i} = \setof{\vect 0}$, it is particularly easy to compute species~i's Voronoi star-cell.

\begin{lem}[Scale-invariant starting configurations with $\initialset{1} = \setof{\vect 0}$]
\label{zero_scale_inv_conquering_lem}
If $\initialset{1} = \setof{\vect 0}$ and $\initialset{2}$ is a blunt scale-invariant set, then $\cstarset{\normpair}{1}{1}{\vect 0}{\initialset{2}} = \closerset{\normpair}{1}{1}{\vect 0}{\initialset{2}}$, and hence species~1 conquers all of $\closerset{\normpair}{1}{1}{\vect 0}{\initialset{2}}$ in the process $\reacheddfcn{\normpair}{\pair{\vect 0}{\,\initialset{2}}}$.
\end{lem}


\begin{proof}
By definition \eqref{closerset_def_eqn},
\begin{align*}
\closerset{\normpair}{1}{1}{\vect 0}{\initialset{2}}
&= \setbuilder[1]{\vect x\in\Rd}{\mu_1(\vect x) \le \distnew{\mu_2}{\initialset{2}}{\vect x}}.
\end{align*}
Since $\initialset{2}$ is scale-invariant, the set on the right is evidently scale-invariant and contains $\vect 0$. Therefore, if the point $\vect x$ is contained in $\closerset{\normpair}{1}{1}{\vect 0}{\initialset{2}}$, so is the entire ray $[0,\infty)\cdot \vect x$. In particular, this implies that the segment $[0,1]\cdot \vect x  = \segment{\vect 0}{\vect x}$ is contained in $\closerset{\normpair}{1}{1}{\vect 0}{\initialset{2}}$, so $\vect x\in \cstarset{\normpair}{1}{1}{\vect 0}{\initialset{2}}$ by definition \eqref{cstar_def_eqn}. Since $\vect x$ was arbitrary, we have $\closerset{\normpair}{1}{1}{\vect 0}{\initialset{2}} = \cstarset{\normpair}{1}{1}{\vect 0}{\initialset{2}}$, and species 1 conquers this set by Proposition~\ref{determ_star_conquering_prop}.
\end{proof}

%
%

\section{Competition with Euclidean Norms}
\label{deterministic:euclid_comp_sec}

%
%
%
%

Deterministic competition in which both species' traversal norms are Euclidean (i.e.\ $\mu_1$ and $\mu_2$ are multiples of the $\lspace{2}{d}$ norm) can be considered as a prototypical example which exhibits many of the essential features arising from competition with arbitrary norms. Since explicit calculations are possible when the norms are Euclidean, at least for nice enough starting configurations, it is worth investigating the geometry of the competition process in this case. We will focus on several simple starting configurations in dimension $d=2$; higher dimensional analogues of these configurations can be obtained by rotation.

\subsection{Cartesian Ovals and Logarithmic Spirals (Point vs.\ Disc)}
\label{deterministic:teardrop_sec}

%
%

In this section we consider Euclidean competition from initial configurations in $\R^2$ in which $\initialset{1}$ and $\initialset{2}$ are either points or discs. We first consider the case of equal speeds, then we will consider different speeds. In both cases, the first step in describing the crystal-growth Voronoi diagram is to describe the ordinary Voronoi diagram.

\begin{lem}[Voronoi boundary for equal speeds]
\newcommand{\focus}[1]{\vect z_{#1}}
\newcommand{\thiscurve}{C}

Suppose species~1 and species~2 have equal speeds. For $i\in\setof{1,2}$, let $A_i$ be a nonempty open or closed disc of radius $r_i\ge 0$ centered at $\focus{i}\in\R^2$, and suppose that $\initialset{1}\cap\initialset{2} = \emptyset$. If $R= r_1-r_2$, 
then the boundary between the Voronoi cells for the two species is
\[
\thiscurve =
\setbuilder[1]{\vect x\in\R^2}{\distnew{\ltwo{2}}{\vect x}{\focus{1}} - \distnew{\ltwo{2}}{\vect x}{\focus{2}} = R}.
\]
If $\abs{R}=\distnew{\ltwo{2}}{\focus{1}}{\focus{2}}$, then $\thiscurve$ is a ray originating from $\focus{1}$ or $\focus{2}$. If $R=0$, then $\thiscurve$ is the perpendicular bisector of the segment $\segment{\focus{1}}{\focus{2}}$. If $0<\abs{R}<\distnew{\ltwo{2}}{\focus{1}}{\focus{2}}$, then $\thiscurve$ is (by definition) one branch of a hyperbola with foci $\focus{1}$ and $\focus{2}$.
\end{lem}

Since the Voronoi cells and crystal-growth Voronoi cells agree in the case of equal speeds, the preceding lemma implies the following result for the competition process.

%
%
%
%

\begin{prop}[Equal speeds: Half-planes and hyperbolas]
\label{euclidean_point_disc_equal_prop}
\newcommand{\Aone}{\initialset{1}}
\newcommand{\Atwo}{\initialset{2}}
\newcommand{\tworad}{r} 

Suppose $\mu_1=\mu_2$, and suppose $\initialset{1}$ is single point and $\initialset{2}$ is a disc of radius $r\ge 0$.
\begin{enumerate}
\item If $\tworad>0$ and $\Aone$ lies on the boundary of $\Atwo$, then species~1 conquers a ray perpendicular to the boundary of $\Atwo$.

\item If $\tworad=0$ (i.e.\ the two species initially occupy two distinct points), then each species conquers a half-space.

\item If $\tworad>0$, and the point $\Aone$ is separated from the disc $\Atwo$, then species~1 conquers the region bounded by one branch of a hyperbola.

\end{enumerate}
The third case interpolates between the first two as the point moves closer to the disc  or the radius of the disc increases from 0 to $\distnew{}{\Aone}{\Atwo}$.
\end{prop}

Next we describe the Voronoi diagram in the case of different speeds.

\begin{lem}[Voronoi boundary for different speeds]
\label{euclidean_point_disc_different_vor_lem}
\newcommand{\focus}[1]{\vect z_{#1}}
\newcommand{\thiscurve}{C}
\newcommand{\speedratio}{\Lambda}
\newcommand{\thisdiff}{R}

Suppose species~1 and species~2 have different speeds $\lambda_1$ and $\lambda_2$, respectively, and let $\speedratio = \lambda_1/\lambda_2$. For $i\in\setof{1,2}$, let $A_i$ be a nonempty open or closed disc of radius $r_i\ge 0$ centered at $\focus{i}\in\R^2$, and suppose that $\initialset{1}\cap\initialset{2} = \emptyset$. If $R= r_1-\speedratio r_2$, 
then the boundary between the Voronoi cells for the two species is
\[
\thiscurve =
\setbuilder[1]{\vect x\in\R^2}{\distnew{\ltwo{2}}{\vect x}{\focus{1}} 
	- \speedratio \distnew{\ltwo{2}}{\vect x}{\focus{2}} = R}.
\]
If ${R}=\distnew{\ltwo{2}}{\focus{1}}{\focus{2}}$, then $\thiscurve$ is the outer boundary of a looped lima\c{c}on of Pascal. If $R=0$, then $\thiscurve$ is (by definition) the Apollonius circle with distance ratio $\speedratio$ for the points $\focus{1}$ and $\focus{2}$. If $0<\abs{R}<\distnew{\ltwo{2}}{\focus{1}}{\focus{2}}$, then $\thiscurve$ is (by definition) one branch of the Cartesian ovals with foci $\focus{1}$ and $\focus{2}$ and distance ratio $\speedratio$.
\end{lem}

The Cartesian ovals are a family of fourth degree algebraic plane curves; each member of the family consists of two branches, some of which are ``oval"-shaped. For a description of these curves, see \cite[pp.~295--299]{Rice:1888aa}. The following proposition describes the crystal-growth Voronoi diagram for competing species with different speeds when the faster species starts at a point and the slower species starts on a disc.

%

\begin{prop}[Different speeds: Logarithmic spirals and Cartesian ovals]
\label{euclidean_point_disc_different_prop}
\newcommand{\Aone}{\initialset{1}}
\newcommand{\Atwo}{\initialset{2}}
\newcommand{\tworad}{r} 

Suppose species~1 is strictly faster than species~2, and species~1 starts at a single point $\Aone$ while species~2 starts on a disc $\Atwo$ of radius $\tworad\ge 0$.
\begin{enumerate}
\item If $\tworad>0$ and $\Aone$ is a single point on the boundary of the disc $\Atwo$, then species~2 conquers a ``heart-shaped" region bounded by two symmetric, partial logarithmic spirals.

\item If $\tworad=0$ (i.e.\ species~1 and species~2 start at two distinct points), then species~2 conquers a ``teardrop" shaped region, with the boundary of the teardrop's ``head" consisting of an arc of the Apollonius circle for the points $\Aone$ and $\Atwo$, and the boundary of the ``tail" of the teardrop consisting of two partial logarithmic spirals. 

\item If $\tworad>0$ and the point $\Aone$ does not lie on the boundary of $\Atwo$, then species~2 conquers a region shaped like either a ``heart" or ``teardrop" bounded by three curves, one being an arc of a Cartesian oval, the other two being arcs of symmetric logarithmic spirals.

\end{enumerate}
The third case interpolates between the first two as the point moves closer to the disc or the radius of the disc increases from 0 to $\distnew{}{\Aone}{\Atwo}$.
\end{prop}

\begin{proof}
The shape in Part~2 is described in \cite{Schaudt:1992aa} and \cite{Kobayashi:2002aa}. The other two parts are proved similarly, using Lemma~\ref{euclidean_point_disc_different_vor_lem} for Part~3. In Parts~2 and 3, the transition from the ordinary Voronoi boundary to the logarithmic spiral occurs at the points where a ray originating at $\initialset{1}$ is tangent to the Voronoi boundary, as these are the farthest points on the boundary that can be reached by a shortest path (i.e.\ straight line) contained completely within species~1's Voronoi cell.
\end{proof}

For a given point-disc starting configuration $\initconfig$, the finally conquered regions in each part of Proposition~\ref{euclidean_point_disc_different_prop} will approach the configuration in the corresponding part of Proposition~\ref{euclidean_point_disc_equal_prop} as the relative speed of the two species approaches~1.

\subsection{Conic Sections (Point vs.\ Hyperplane)}
\label{deterministic:conics_sec}

In this section we consider Euclidean competition in the case where $A_1$ is a point and $A_2$ is a hyperplane. Note that in $d=2$, the boundary between the (ordinary) Voronoi cells for this configuration is precisely a conic section with focus $A_1$, directrix $A_2$, and eccentricity equal to the ratio of species~1's speed to species~2's speed. It turns out that the crystal-growth Voronoi diagram is equal to the ordinary Voronoi diagram in the hyperbolic and parabolic regimes, but is slightly different in the elliptic regime.

\begin{prop}[Conic sections]
Suppose $A_1$ is a point and $A_2$ is a line. There are three cases, depending on the relative speed of the two species.
\begin{enumerate}
\item If species~1 is faster, then species~1 conquers a region bounded by one branch of a hyperbola with focus $A_1$ and directrix $A_2$.
\item If the speeds are the same, species~1 conquers a region bounded by a parabola with focus $A_1$ and directrix $A_2$.
\item If species~1 is slower, then species~1 conquers an ``elliptic teardrop," with the boundary of the teardrop's ``head" equal to an arc the ellipse with focus $A_1$ and directrix $A_2$, and the boundary of the ``tail" consisting of two partial logarithmic spirals.
\end{enumerate}
\end{prop}

\begin{proof}
The first two parts follow because in the hyperbolic and parabolic case, any point in species~2's Voronoi region is connected to $\initialset{2}$ via a shortest path contained completely within the Voronoi region.
In the elliptic case, the shortest path from $\initialset{2}$ to a point on the ``far" side of the ellipse passes through the ellipse's boundary, so species~2 must follow a curved path to reach these points. Thus there are two transition points on the ellipse where the crystal-growth Voronoi boundary changes to a logarithmic spiral as in Proposition~\ref{euclidean_point_disc_different_prop}; the transition occurs at the points where shortest paths from the hyperplane are tangent to the ellipse.
\end{proof}

Also, note what species 1's conquered set looks like in the limit as the point approaches the hyperplane: Cone for supercritical speeds; ray for critical (equal) speeds; point for subcritical speeds.

\subsection{Critical Speeds Greater than 1 (Point vs.\ Conical Shell)}
\label{deterministic:euclid_cone_sec}

Consider what happens when $A_1$ is a point and $A_2$ is a union of hyperplanes, e.g.\ tangent to some small convex cone. Then the region conquered by 1 should be the intersection of the corresponding regions from Section~\ref{deterministic:conics_sec}. For extra-small cones, the critical speed always occurs in the hyperbolic regime. In an extra-small circular cone, criticality occurs when the asymptotes of all the conquered hyperbolas are parallel.

%

\section{Competition When Species 1 Is Initially Surrounded}
\label{deterministic:surrounded_sec}

In this section we describe the deterministic growth from two starting configurations which will be important when analyzing the random process in Chapter~\ref{random_fpc_chap}. In both cases, species~1 starts at a single point $\vect z$, and species~2 initially occupies everything outside some bounded set containing $\vect z$, namely $\initialset{2} = \shellof{B} = \Rd\setminus \interior{B}$ for some bounded $B\subseteq\Rd$.

\subsection{Conquering a $\mu_1$-Ball}
\label{deterministic:conq_mu1_ball_sec}

For norms $\mu_1$ and $\mu_2$ on $\Rd$ and $\beta>0$, define
\begin{equation}
\label{rmeet_def_eqn}
\rmeet{1}{\beta} \definedas \frac{1}{\beta + \esupnorm{\mu_2/\mu_1}}
\quad\text{and}\quad
\rmeet{2}{\beta} \definedas \frac{1}{\beta + \esupnorm{\mu_1/\mu_2}}.
\end{equation}
We call $\rmeet{1}{\beta}$ the \textdef{$\beta$-meeting radius for species 1} and $\rmeet{2}{\beta}$ the \textdef{$\beta$-meeting radius for species 2}; this terminology is explained by the following lemma, which shows that in a $\beta$-adjusted process in which species 1 starts at some point $\vect z\in\Rd$ and species 2 starts outside a unit $\mu_2$-ball centered at $\vect z$, $\rmeet{1}{\beta}$ is the radius of the $\mu_1$-ball species 1 has conquered when the two species meet. The factors $\rmeet{1}{\beta}$ and $\rmeet{2}{\beta}$ will show up in the proofs of several later lemmas concerning both the deterministic and random processes in Chapters~\ref{deterministic_chap} and \ref{random_fpc_chap}.

\begin{lem}[Competition inside a $\mu_2$-spherical shell]
\label{determ_bullseye_lem}

\newcommand{\origin}{\vect z}
\newcommand{\rad}{r}
\newcommand{\thisrmeet}{ \rmeet{1}{\beta} \rad}

\newcommand{\tworad}{\esupnorm[2]{\tfrac{\mu_2}{\mu_1}} \thisrmeet}

\newcommand{\bigball}{\reacheddset{\mu_2}{\origin}{\rad}}
\newcommand{\smallball}{\reacheddset{\mu_1}{\origin}{\thisrmeet}}
\newcommand{\smalltwoball}{\reacheddset[1]{\mu_2}{\origin}{\tworad}}

\newcommand{\bigshell}{\shellof{\bigball}}
\newcommand{\smallshell}{\shellof{\smalltwoball}}

\newcommand{\Aone}{\origin}
\newcommand{\Atwo}{\initialset{2}}
\newcommand{\thisconfig}{\pair[1]{\origin}{\,\Atwo}}

\newcommand{\onecstarset}{\cstarset[1]{\normpair}{\beta}{1}{\Aone}{\Atwo}}
\newcommand{\twocstarset}{\cstarset[1]{\normpair}{\beta^+}{2}{\Aone}{\Atwo}}

\newcommand{\pairone}{\pair{\beta \mu_1}{\mu_2}}
\newcommand{\pairtwo}{\pair{\mu_1}{\beta^{-1}\mu_2}}

Fix $\rad >0$ and $\origin\in\Rd$, and let $\Atwo= \bigshell$.
Then for any $\beta>0$,
\begin{enumerate}

\item $\onecstarset\supseteq \smallball$, and if $R> \thisrmeet$, then $\twocstarset \cap \reacheddset{\mu_1}{\origin}{R} \ne \emptyset$.

\item In either of the $\beta$-adjusted processes $\reacheddset{\pairone}{\thisconfig}{t}$ or $\reacheddset{\pairtwo}{\thisconfig}{t}$, the largest $\mu_1$-ball conquered by species 1 is $\smallball$. In fact,
\begin{enumerate}
\item In the process $\reacheddset{\pairone}{\thisconfig}{t}$, the meeting time of the two species is $\beta\thisrmeet$, and
\[
\reacheddset{1}{\thisconfig}{\beta \thisrmeet} = \smallball.
\]
\item In the process $\reacheddset{\pairtwo}{\thisconfig}{t}$, the meeting time of the two species is $\thisrmeet$, and
\[
\reacheddset{1}{\thisconfig}{\thisrmeet} = \smallball.
\]
\end{enumerate}


\end{enumerate}
If the roles of species 1 and 2 are reversed, symmetric statements hold, with $\rmeet{1}{\beta}$ replaced by $\rmeet{2}{\beta}$.
\end{lem}

\begin{proof}
\newcommand{\origin}{\vect z}
\newcommand{\rad}{r}
\newcommand{\thisrmeet}{ \rmeet{1}{\beta} \rad}

\newcommand{\tworad}{\esupnorm[2]{\tfrac{\mu_2}{\mu_1}} \thisrmeet}

\newcommand{\bigball}{\reacheddset{\mu_2}{\origin}{\rad}}
\newcommand{\smallball}{\reacheddset{\mu_1}{\origin}{\thisrmeet}}
\newcommand{\smalltwoball}{\reacheddset[1]{\mu_2}{\origin}{\tworad}}

\newcommand{\bigshell}{\shellof{\bigball}}
\newcommand{\smallshell}{\shellof{\smalltwoball}}

\newcommand{\Aone}{\origin}
\newcommand{\Atwo}{\initialset{2}}
\newcommand{\thisconfig}{\pair[1]{\origin}{\,\Atwo}}

\newcommand{\onecstarset}{\cstarset[1]{\normpair}{\beta}{1}{\Aone}{\Atwo}}
\newcommand{\twocstarset}{\cstarset[1]{\normpair}{\beta^+}{2}{\Aone}{\Atwo}}

\newcommand{\pairone}{\pair{\beta \mu_1}{\mu_2}}
\newcommand{\pairtwo}{\pair{\mu_1}{\beta^{-1}\mu_2}}

Note that from this starting configuration, each species' $\beta$-Voronoi cell coincides with its $\beta$-Voronoi star-cell, so Corollary~\ref{beta_adjusted_conq_cor} implies that in either of the $\beta$-adjusted processes, both conquered sets coincide with the corresponding $\beta$-Voronoi cell.
In the case $\beta=1$, the meeting time of the two species is
\begin{align*}
t_\text{meet}
&= \sup\setbuilder[1]{t\ge 0}
	{\reacheddset{\mu_1}{\origin}{t} \subseteq \reacheddset{\mu_2}{\origin}{\rad -t}}\\
&= \sup \setbuilder[1]{t\ge 0}
	{\reacheddset[1]{\mu_2}{\origin}{\esupnorm[2]{\tfrac{\mu_2}{\mu_1}} t} 
	\subseteq \reacheddset{\mu_2}{\origin}{\rad -t}}\\
&= \sup \setbuilder[1]{t\ge 0}{\esupnorm[2]{\tfrac{\mu_2}{\mu_1}} t\le \rad - t}\\
&= \frac{\rad}{1+ \esupnorm[2]{\tfrac{\mu_2}{\mu_1}}}\\
&=\rmeet{1}{1} \rad.
\end{align*}
Thus, any $\mu_1$-ball of radius $\rmeet{1}{1} \rad$ or smaller is contained in its 1-Voronoi cell, and any $\mu_1$-ball of larger radius intersects the interior of species 2's conquered region, which coincides with species 2's $1^+$-Voronoi cell. The corresponding statements for general $\beta>0$ follow by time scaling and Corollary~\ref{beta_adjusted_conq_cor}.
\end{proof}

%
%

\subsection{Conquering a $\mu_1$-Cone Segment}
\label{deterministic:conq_mu1_cone_sec}

The next lemma computes the largest $\mu_1$-cone segment contained in species 1's $\beta$-Voronoi cell when species~1 starts at the point $\vect z$ and species~2 starts on the complement of a $\mu_2$-cone segment with apex $\vect z$. Recall the definitions
\[
\thpadv{2}{\dirvec}{\beta} \definedas \beta \frac{\mu_1(\dirsymb)}{\mu_2(\dirsymb)}
\quad\text{and}\quad
\rmeet{1}{\beta} \definedas \frac{1}{\beta + \esupnorm{\mu_2/\mu_1}}
\]
from \eqref{thpadv_def_eqn} and \eqref{rmeet_def_eqn}.


\begin{lem}[Conquering a $\mu_1$-cone segment inside the shell of a $\mu_2$-cone segment]
\label{cone_closer_set_lem}

\newcommand{\bigbeta}{\beta_0}
\newcommand{\smallbeta}{\beta}
\newcommand{\bigthp}{\thpadv{2}{\dirvec}{\bigbeta}}
\newcommand{\smallthp}{\thp}
\newcommand{\smallthpexp}{(\bigbeta-\smallbeta)\rmeet{1}{\smallbeta}}

\newcommand{\bigconetip}{\conetip[2]{\mu_2}{\bigthp}{\thisz}{\dirvec}{\mu_2 (\dirsymb)}}
\newcommand{\smallconetip}{\conetip[2]{\mu_1}{\smallthp}{\thisz}{\dirvec}{\mu_1(\dirsymb)}}
\newcommand{\bigcone}{\cone{\mu_2}{\bigthp}{\thisz}{\dirvec}}
\newcommand{\smallcone}{\cone{\mu_1}{\smallthp}{\thisz}{\dirvec}}

\newcommand{\thisconfig}{\pair{\thisz}{\,\initialset{2}}}
\newcommand{\betanorms}{\pair{\mu_1}{\beta^{-1}\mu_2}}
\newcommand{\betaprocess}{\reacheddset{\betanorms}{\thisconfig}{t}}

\newcommand{\thisz}{\vect z}

Fix a traversal norm pair $\normpair = \pair{\mu_1}{\mu_2}$ on $\Rd$, $\thisz\in\Rd$, $\dirsymb\in \Rd\setminus \setof{\vect 0}$, and $\bigbeta\ge 0$. Suppose $0\le \beta\le \beta_0$, and let $\smallthp \definedas \smallthpexp$.
Then
\begin{enumerate}
\item \label{cone_closer_set:finite_part}
If $\initialset{2} \subseteq \shell \brackets[10]{\bigconetip}$ and $\thisz\not\in \initialset{2}$, then  $\cstarset{\normpair}{\smallbeta}{1}{\thisz}{\initialset{2}} \supseteq \smallconetip,$ and hence species 1 conquers $\smallconetip$ in the $\beta$-adjusted process $\betaprocess$.

\item \label{cone_closer_set:infinite_part}
If $\initialset{2} \subseteq \shell \brackets[10]{\bigcone}$ and $\thisz\not\in \initialset{2}$, then  $\cstarset{\normpair}{\smallbeta}{1}{\thisz}{\initialset{2}} \supseteq \smallcone$, and hence species 1 conquers $\smallcone$ in the $\beta$-adjusted process $\betaprocess$.
\end{enumerate}
A symmetric statement holds with species 1 and 2 switched.
\end{lem}

\begin{thmremark}
In Lemma~\ref{cone_closer_set_lem}, the space inside the $\mu_2$-conical shell initially occupied by species 2 becomes wider as $\beta_0$ becomes larger. For large enough $\beta_0$, we have $\thpadv{2}{\dirvec}{\beta_0}>1$, and hence the  $\mu_2$-cone segment will actually be a $\mu_2$-ball if finite, or all of $\Rd$ if infinite; thus Lemma~\ref{cone_closer_set_lem} becomes trivial for large $\beta_0$ in the infinite case, but is always nontrivial in the finite case. The cone segment conquered by species 1 becomes narrower as $\beta\nearrow \beta_0$, approaching a line segment or ray in the limit. On the other hand, as $\beta$ becomes smaller, the cone segment conquered by species 1 in the $\beta$-adjusted process becomes wider, but the ``margin of error" with which this set gets conquered becomes increasingly bad (the ``margin of error" interpretation of $\beta$ makes sense for $\beta>1$; for $\beta<1$, we have a ``negative margin of error," meaning that species 1 needs to \emph{speed up} by a factor of $\beta^{-1}$ in order to conquer the cone segment).
\end{thmremark}

\begin{proof}[Proof of Lemma~\ref{cone_closer_set_lem}]

\newcommand{\bigbeta}{\beta_0}
\newcommand{\smallbeta}{\beta}
\newcommand{\bigthp}{\thpadv{2}{\dirvec}{\bigbeta}}
\newcommand{\smallthp}{\thp}
\newcommand{\smallthpexp}{(\bigbeta-\smallbeta)\rmeet{1}{\smallbeta}}

\newcommand{\thisz}{\vect z}
\newcommand{\thisy}{\vect y}
\newcommand{\thisx}{\vect x}

\newcommand{\thisapex}{\vect 0}

\newcommand{\twoht}{h_2}
\newcommand{\oneht}{h_1}

\newcommand{\bigconetip}{\conetip[2]{\mu_2}{\bigthp}{\thisapex}{\dirvec}{\twoht}}
\newcommand{\smallconetip}{\conetip[2]{\mu_1}{\smallthp}{\thisapex}{\dirvec}{\oneht}}

\newcommand{\bigcone}{\cone{\mu_2}{\bigthp}{\thisapex}{\dirvec}}
\newcommand{\smallcone}{\cone{\mu_1}{\smallthp}{\thisapex}{\dirvec}}

\newcommand{\thisconfig}{\pair{\thisz}{\,\initialset{2}}}
\newcommand{\betanorms}{\pair{\mu_1}{\beta^{-1}\mu_2}}
\newcommand{\betaprocess}{\reacheddset{\betanorms}{\thisconfig}{t}}

\newcommand{\ballcenter}{\thisx}
\newcommand{\ballpoint}{\thisy}
\newcommand{\coneball}{\reacheddset[2]{\mu_1}{\ballcenter}{\thp \mu_1(\ballcenter)}}
\newcommand{\lineseg}{\segment{\thisapex}{\oneht \dirsymb}}

\newcommand{\thisstarvor}{\cstarset{\normpair}{\smallbeta}{1}{\thisapex}{\initialset{2}}}
\newcommand{\thisvor}{\closerset{\normpair}{\smallbeta}{1}{\thisapex}{\initialset{2}}}

\newcommand{\thisdil}{\esupnorm{\mu_2/\mu_1}}

Assume without loss of generality that $\thisz=\thisapex$. For Part~\ref{cone_closer_set:finite_part}, set $\oneht = \mu_1(\dirsymb)$ and $\twoht = \mu_2(\dirsymb)$, and for Part~\ref{cone_closer_set:infinite_part}, set $\oneht = \twoht = \infty$. 
Note that since $\smallconetip$ is star-shaped at $\thisapex$, we have
\[
\smallconetip \subseteq \thisstarvor \iff \smallconetip \subseteq \thisvor.
\]
Thus, since $\smallconetip = \bigcup_{\ballcenter\in \lineseg} \coneball$, we seek a $\thp\ge 0$ so that $\coneball$ is guaranteed to be contained in $\thisvor$ for all $\ballcenter\in \lineseg$. That is, we need to find some $\thp\ge 0$ small enough that for any $\ballcenter\in \lineseg$ we have
\begin{equation}
\label{cone_closer_set:ball_in_vor_eqn}
\smallbeta \distnew{\mu_1}{\thisapex}{\ballpoint}
\le \distnew{\mu_2}{\initialset{2}}{\ballpoint}
\text{ for all }
\ballpoint\in \coneball.
\end{equation}
Fix some $\ballcenter\in \lineseg$, and let $\thp$ be a nonnegative number (yet to be determined). By the triangle inequality for $\mu_1$ we have
\[
\distnew{\mu_1}{\thisapex}{\ballpoint}\le (1+\thp) \mu_1(\ballcenter)
\text{ for all }
\ballpoint\in \coneball.
\]
Since $\initialset{2} \subseteq \shell \bigconetip$, we have $\distnew{\mu_2}{\initialset{2}}{\ballcenter} \ge \bigbeta \mu_1(\ballcenter)$, and if $\mu_1(\ballpoint-\ballcenter)\le \thp\mu_1(\ballcenter)$, then $\mu_2(\ballpoint - \ballcenter)\le \thp \thisdil \mu_1(\ballcenter)$. Thus, by the triangle inequality for $\mu_2$ we have
\[
\distnew{\mu_2}{\initialset{2}}{\ballpoint} \ge
	\parens[2]{\bigbeta + \thisdil} \mu_1(\ballcenter)
\text{ for all }
\ballcenter\in \coneball.
\]
It follows from the above two inequalities that \eqref{cone_closer_set:ball_in_vor_eqn} is guaranteed to hold if we choose $\thp$ small enough that
\[
\smallbeta (1+\thp) \mu_1(\ballcenter) \le
\parens[2]{\bigbeta - \thp \thisdil} \mu_1(\ballcenter),
\]
or
\[
\thp\le \frac{\bigbeta-\smallbeta}{\smallbeta+\thisdil} = \smallthpexp.
\]
Thus, if we take $\thp\definedas \smallthpexp$, then $\coneball \subseteq \thisvor$, and since this choice of $\thp$ is independent of $\ballcenter$, the entire union $\smallconetip$ is contained in $\thisvor$ as claimed.
\end{proof}

\section{Competition When Species 2 Starts on a Conical Shell}
\label{deterministic:cone_competition_sec}

In this section we consider a process in which species 2 starts on the boundary of a cone $\cones$ in $\Rd$, and species 1 starts at some point $\vect z$ inside the cone. Equivalently,
we can consider the starting configuration $\pair[1]{\vect z}{\,\Rd\setminus \cones}$, where $\vect z\in \cones$. With this initial configuration, species 1 is confined to the interior of the cone $\cones$, surrounded by an infinite expanse of species 2 on $\bd \cones$, and the only chance it has of escaping to infinity is to travel (asymptotically) in the direction of some ray in $\dirset{\cones}$, hoping to outrun species 2 as it encroaches from the boundary. The main point of interest is whether species 1 is able to survive indefinitely within the cone, conquering an unbounded set, or whether it eventually ends up completely surrounded by species~2.


In order to analyze the process from this starting configuration, it will be convenient to introduce some notation. For any set $\cones\subseteq\Rd$ and any $\vect z\in \cones$, we define $\conqset{1}{\normpair}{\vect z}{\cones}$ to be the region conquered by species 1 from the starting configuration $\pair[2]{\vect z}{\,\Rd\setminus \cones}$ in a process with traversal norm pair $\normpair$. That is,
\begin{equation}\label{conq_set_def_eqn}
\conqset{1}{\normpair}{\vect z}{\cones} \definedas 
\reacheddset{1}{\pair{\vect z}{\,\Rd\setminus \cones}}{\infty}_{\normpair}.
\end{equation}
Although this definition makes sense for an arbitrary subset $\cones$ of $\Rd$, we will always take $\cones$ to be a cone. We will start with the simplest case where $\vect z$ is an apex of the cone, and our first goal will be to classify cones according to whether species 1 can escape to infinity from this starting configuration. This classification is accomplished in Definition~\ref{cone_criticality_def} and Proposition~\ref{conquering_from_apex_prop} in Section~\ref{deterministic:cone_advantage_sec}. We then prove some technical results in Sections~\ref{deterministic:advantage_properties_sec} and \ref{deterministic:gen_advantage_sec} before considering the case of general $\vect z\in\cones$ in Section~\ref{deterministic:cones:conquered_regions_sec}. In Section~\ref{deterministic:cones:additional_sec} we describe some additional results about competition in cones that follow from the results in Section~\ref{deterministic:cones:conquered_regions_sec}.

\subsection{Speed Advantage; Favorable Directions; Wide and Narrow Cones}
\label{deterministic:cone_advantage_sec}

Let $\cones\subseteq\Rd$ be a pointed cone with apex $\vect a$, and let $\siversion{\cones}$ be the scale-invariant version of $\cones$, with apex $\vect 0$. We first consider the case when species 1 starts at the apex, i.e.\ the initial configuration is $\pair[1]{\vect a}{\Rd\setminus \cones}$. Intuitively, since this starting configuration is scale-invariant at $\vect a$, species 1 gains nothing by following a curved path --- the only chance it has of escaping is to
move along each ray in the cone from $\vect a$ out toward infinity, and hope that it can stay ahead of species 2 in some direction.
By the scale-invariance of the starting configuration, if at any time $t>0$ species 1 is ahead of species 2 in some direction, then it will stay ahead of species 2 for all $t>0$. That is, species 1 can survive indefinitely within the cone if and only if it can survive for some positive amount of time within the cone. We will prove these statements formally below (Proposition~\ref{conquering_from_apex_prop}), but first we use this intuition to classify directions within the cone according to how favorable or unfavorable they are to species 1, by quantifying
the rate at which species 1 outruns or falls behind species 2 when traveling in that direction.

\subsubsection{Favorability of Directions in a Cone}

Let $\dirsymb\in \siversion{\cones}\setminus \setof{\vect 0}$, and fix a traversal norm pair $\normpair$. We define the \textdef{speed advantage}, or simply the \textdef{advantage}, \textdef{of the direction $\dirvec$ in $\cones$} for species 1 as
\begin{equation}\label{dir_advantage_def_eqn}
\diradvantage{1}{\normpair}{\cones}{\dirvec} \definedas
\frac{\distnew{\mu_2}{\bd \siversion{\cones}}{\dirsymb}}{\mu_1(\dirsymb)}.
\end{equation}
Note that the advantage is well-defined and continuous on $\dirset{\cones}$ since it is defined as a scale-invariant continuous function of points in $\siversion{\cones}\setminus \setof{\vect 0}$.
We call the direction $\dirvec\in \dirset{\cones}$ \textdef{favorable}, \textdef{unfavorable}, or \textdef{critical} for species 1 in $\cones$ according to whether $\diradvantage{1}{\normpair}{\cones}{\dirvec}$ is greater than, less than, or equal to 1, respectively.
If $\dirvec$ is favorable in $\cones$, then species 1 outruns species 2 from the starting configuration $\pair[1]{\vect a}{\Rd\setminus \cones}$ by traveling along the ray $\vect a+\dirvec$. If $\dirvec$ is unfavorable in $\cones$,
then (at least intuitively; see Remark~\ref{unfav_conquerable_rem} below) species 2 beats species 1 to every point on the ray $\vect a+\dirvec$, and if the direction is critical, then the two species tie along this ray.

We also make the following definitions for $\vect a\in \apexset{\cones}$:
\begin{equation}
\label{favwedge_def_eqn}
\begin{split}
\cfavwedge{\normpair}{\vect a}{\cones}
&\definedas \setbuilder[2]{\vect x\in \cones}{\distnew{\mu_1}{\vect a}{\vect x} \le \distnew{\mu_2}{\bd \cones}{\vect x}},\\
\favwedge{\normpair}{\vect a}{\cones}
&\definedas \setbuilder[2]{\vect x\in \cones}{\distnew{\mu_1}{\vect a}{\vect x} < \distnew{\mu_2}{\bd \cones}{\vect x}}.
\end{split}
\end{equation}
We call the sets $\cfavwedge{\normpair}{\vect a}{\cones}$ and $\favwedge{\normpair}{\vect a}{\cones}$ the \textdef{favorable wedge} at $\vect a$ and the \textdef{strictly favorable wedge} at $\vect a$, respectively, because the ray $\vect a+\dirvec$ is contained in $\cfavwedge{\normpair}{\vect a}{\cones}$ (resp.\ $\favwedge{\normpair}{\vect a}{\cones}$) if and only if $\diradvantage{1}{\normpair}{\cones}{\dirvec} \ge 1$ (resp.\ $\diradvantage{1}{\normpair}{\cones}{\dirvec}>1$). Our first result is the following.

\begin{lem}[Favorable directions are conquered]
\label{favorable_conquered_lem}
Let $\cones\subseteq\Rd$ be a pointed cone with apex $\vect a$, and let $\normpair$ be a traversal norm pair. In the deterministic two-type process $\reacheddfcn{\normpair}{\pair{\vect a}{\,\Rd\setminus \cones}}$, species 1 conquers the favorable wedge at $\vect a$; that is,
$
\conqset{1}{\normpair}{\vect a}{\cones} \supseteq
\cfavwedge{\normpair}{\vect a}{\cones}.
$
\end{lem}

\begin{proof}
By translation-invariance, it suffices to assume $\vect a=\vect 0$. Note that
\[
\cfavwedge{\normpair}{\vect 0}{\cones}
= \setbuilder[2]{\vect x\in \cones}{\mu_1(\vect x) \le \distnew{\mu_2}{\bd \cones}{\vect x}}
= \closerset[2]{\normpair}{1}{1}{\vect 0}{\Rd\setminus \cones}.
\]
Since $\Rd\setminus \cones$ is scale-invariant, Lemma~\ref{zero_scale_inv_conquering_lem} implies that
\[
\closerset{\normpair}{1}{1}{\vect 0}{\Rd\setminus \cones} = \cstarset{\normpair}{1}{1}{\vect 0}{\Rd\setminus \cones},
\]
and this latter set is conquered by species 1 according to Proposition~\ref{determ_star_conquering_prop}. Thus, species 1 conquers the entire set $\cfavwedge{\normpair}{\vect 0}{\cones}$.
\end{proof}

\begin{thmremark}
\label{unfav_conquerable_rem}
The containment in Lemma~\ref{favorable_conquered_lem} can be strict in general, at least for large cones (i.e.\ those not contained in a half-space). That is, for some large cones $\cones$, one can find traversal norm pairs such that species 1 conquers some unfavorable directions in $\cones$ from the starting configuration $\pair[2]{\vect a}{\Rd\setminus \cones}$. This is accomplished by ``blocking off" a set of unfavorable directions from species 2 by surrounding the unfavorable region with favorable or critical directions.
\end{thmremark}

\subsubsection{The Advantage in a Cone}

By Lemma~\ref{favorable_conquered_lem}, as long as there is some favorable direction in $\cones$, then species 1 can survive from the apex by traveling in that direction.
This motivates the following definition.
For the traversal norm pair $\normpair$, we define the \textdef{advantage in $\cones$} of species 1 as
\begin{equation}
\label{cone_advantage_def_eqn}
\coneadvantage{1}{\normpair}{\cones} \definedas 
\sup_{\dirvec \in \dirset{\cones}} \diradvantage{1}{\normpair}{\cones}{\dirvec}
= \sup_{\vect x\in \cones\setminus \setof{\vect a}}
	\frac{\distnew{\mu_2}{\bd \cones}{\vect x}}{\distnew{\mu_1}{\vect a}{\vect x}}.
\end{equation}
It is clear by the translation-invariance of the norm metrics that $\coneadvantage{1}{\normpair}{\cones} = \coneadvantage{1}{\normpair}{\siversion{\cones}}$, so the advantage in the cone is independent of the choice of apex since the scale-invariant version of $\cones$ is unique. Similar to our above classification of directions according to \eqref{dir_advantage_def_eqn}, we use \eqref{cone_advantage_def_eqn} to classify the entire cone in terms of its favorability to species 1, as follows.

\begin{definition}[Super- and sub- criticality regimes for cones]
\label{cone_criticality_def}
Let $\cones$ be a cone in $\Rd$, and let $\alpha = \coneadvantage{1}{\normpair}{\cones}$, defined in \eqref{cone_advantage_def_eqn}. Then we say that species 1 is \textdef{supercritical}, \textdef{critical}, or \textdef{subcritical} for $\cones$, or that $\cones$ is \textdef{wide}, \textdef{critical}, or \textdef{narrow} for species 1, according to whether $\alpha>1$, $\alpha=1$, or $\alpha<1$, respectively.
\end{definition}

Before proceeding, we prove the following basic result about the advantage in a cone.

\begin{lem}[Existence of a direction of maximal advantage]
\label{max_adv_dir_lem}
Let $\cones$ be a cone in $\Rd$ with apex $\vect a$ and scale-invariant version $\siversion{\cones}$. If $\cones$ is degenerate, then $\coneadvantage{1}{\normpair}{\cones} = 0$. If $\cones$ is nondegenerate, then
$
0<\coneadvantage{1}{\normpair}{\cones} \le \esupnorm{\mu_2/\mu_1} <\infty,
$
and the first supremum in \eqref{cone_advantage_def_eqn} is achieved for some $\dirvec\in \dirset{\interior{\cones}}$; equivalently, there is some $\dirsymb\in \interior{\siversion{\cones}}$ such that the second supremum in \eqref{cone_advantage_def_eqn} is achieved by every $\vect x\in \interior{\cones}$ of the form $\vect x = \vect a + r\dirsymb$ with $r>0$.
\end{lem}

\begin{proof}
If $\cones$ is degenerate, then $\cones = \bd \cones$, so $\distnew{\mu_2}{\bd \cones}{\vect x} = 0$ for all $\vect x\in \cones$, and so $\coneadvantage{1}{\normpair}{\cones} = 0$ by the second formula in \eqref{cone_advantage_def_eqn}. For any cone $\cones$, we have $\vect 0\in \bd \siversion{\cones}$, so for any $\dirvec \in \dirset{\cones}$,
\[
\diradvantage{1}{\normpair}{\cones}{\dirvec} \le
\frac{\distnew{\mu_2}{\vect 0}{\dirsymb}}{\mu_1(\dirsymb)}
= \frac{\mu_2(\dirsymb)}{\mu_1(\dirsymb)}.
\]
Therefore, since $\dirvec$ was arbitrary,
\[
\coneadvantage{1}{\normpair}{\cones} \le \sup_{\dirvec \in \dirset{\cones}}
\frac{\mu_2(\dirsymb)}{\mu_1(\dirsymb)}
\le \esupnorm[1]{\frac{\mu_2}{\mu_1}} <\infty.
\]
If $\cones$ is nondegenerate, then $\interior{\siversion{\cones}}$ is nonempty. Choosing any $\dirsymb\in \interior{\siversion{\cones}}$, we have $\dirvec\in \dirset{\interior{\cones}}$ and
\[
\coneadvantage{1}{\normpair}{\cones} \ge
\frac{\distnew{\mu_2}{\bd \siversion{\cones}}{\dirsymb}}{\mu_1(\dirsymb)} >0.
\]
To see that the supremum in \eqref{cone_advantage_def_eqn} is achieved for some $\dirvec\in \dirset{\interior{\cones}}$ when $\cones$ is nondegenerate, first note that since every nondegenerate cone is a body by definition, we have $\bd \siversion{\closure{\cones}} = \bd \siversion{\cones}$ and hence $\diradvantage{1}{\normpair}{\cones}{\dirvec} = \diradvantage{1}{\normpair}{\closure{\cones}}{\dirvec}$ for any $\dirvec\in\dirset{\cones}$, and also $\interior{\parens[2]{\,\closure{\cones}\,}} = \interior{\cones}$. Therefore it suffices to assume $\cones$ is closed.

If $\cones$ is closed, then $\dirset{\cones}$ is compact (being homeomorphic to $\sphere{d-1}{\norm{\cdot}}\cap \siversion{\cones}$), so the continuous function $\diradvantage{1}{\normpair}{\cones}{\cdot}$ achieves a maximum on $\dirset{\cones}$.
Since $\diradvantage{1}{\normpair}{\cones}{\cdot}$ is zero on $\dirset{\bd \cones}$ and strictly positive on $\dirset{\interior{\cones}}$, the maximum must occur at some direction  $\dirvec \in \dirset{\interior{\cones}}$, which corresponds to some point $\dirsymb\in \interior{\siversion{\cones}}$. Finally, if $r>0$ and $\vect x = \vect a + r\dirsymb$, then $\vect x\in \interior{\cones}$, and by the translation-invariance and scale-equivariance of the norm metrics,
\[
\frac{\distnew{\mu_2}{\bd \cones}{\vect x}}{\distnew{\mu_1}{\vect a}{\vect x}}
= \frac{\distnew{\mu_2}{\bd \siversion{\cones}}{r\dirsymb}}{\mu_1(r\dirsymb)}
= \diradvantage{1}{\normpair}{\cones}{\dirsymb},
\]
so the supremum in the final expression in \eqref{cone_advantage_def_eqn} is achieved for any such $\vect x$.
\end{proof}

\subsubsection{The Size of the Conquered Region in Wide, Critical, and Narrow Cones}


The classification of cones in Definition~\ref{cone_criticality_def} leads us to the following proposition, which is one of the central ideas underlying this thesis.
Namely, in Proposition~\ref{conquering_from_apex_prop}, we show that the classification of cones as wide, critical, or narrow corresponds generally to the ``size" of the set species 1 conquers when it starts at an apex of the cone.
The proofs for wide and critical cones follow easily from Lemmas~\ref{favorable_conquered_lem} and \ref{max_adv_dir_lem} above, and the proof for narrow cones follows from Lemma~\ref{unobtrusive_star_set_lem} below. In turn, Lemma~\ref{unobtrusive_star_set_lem} relies on Lemmas~\ref{blocking_segments_lem} and \ref{strong_blocking_segments_lem}, whose proofs are rather intricate and have therefore been placed in an appendix to avoid interrupting the flow of the section. Proposition~\ref{conquering_from_point_prop} in Section~\ref{deterministic:cones:conquered_regions_sec} below generalizes Proposition~\ref{conquering_from_apex_prop} to the case where species 1 starts at a general point in $\cones$, not necessarily an apex.

\begin{prop}[The conquered region in a cone when species 1 starts at an apex]
\label{conquering_from_apex_prop}
Let $\cones\subseteq\Rd$ be a pointed cone with apex $\vect a$, and let $\conqset{1}{\normpair}{\vect a}{\cones}$ be species 1's conquered region from the starting configuration $\pair{\vect a}{\,\Rd\setminus \cones}$, as defined in \eqref{conq_set_def_eqn}.
\begin{enumerate}
\item \label{conquering_from_apex:wide_part}
If $\cones$ is wide for species 1 then $\conqset{1}{\normpair}{\vect a}{\cones}$ contains a nondegenerate cone at $\vect a$.

\item \label{conquering_from_apex:crit_part}
If $\cones$ is critical for species 1, then $\conqset{1}{\normpair}{\vect a}{\cones}$ contains at least a ray from $\vect a$.

\item \label{conquering_from_apex:narr_part}
If $\cones$ is narrow for species 1 and additionally is contained in a half-space, then
\[
\conqset{1}{\normpair}{\vect a}{\cones} = \setof{\vect a}.
\]
\end{enumerate}
\end{prop}

\begin{proof}
It suffices to assume $\vect a=\vect 0$. Recall from Lemma~\ref{favorable_conquered_lem} above that $\conqset{1}{\normpair}{\vect 0}{\cones} \supseteq \cfavwedge{\normpair}{\vect 0}{\cones}$, where
$
\cfavwedge{\normpair}{\vect 0}{\cones}=
\setbuilder[2]{\vect x\in \cones}{\mu_1(\vect x) \le \distnew{\mu_2}{\bd \cones}{\vect x}}
$
is the favorable wedge at $\vect 0$.

First suppose that $\cones$ is wide for species 1. Then the strictly favorable wedge,
$
\favwedge{\normpair}{\vect 0}{\cones} =
\setbuilder[2]{\vect x\in \cones}{\mu_1(\vect x) < \distnew{\mu_2}{\bd \cones}{\vect x}},
$
is nonempty, and it is open by the continuity of the norm metrics. Moreover, $\favwedge{\normpair}{\vect 0}{\cones}$ is scale-invariant and hence contains a nondegenerate cone since any open scale-invariant set contains an open cone (simply take the blunt scale-invariant set generated by an open ball). Thus, since $\favwedge{\normpair}{\vect 0}{\cones}\subseteq \cfavwedge{\normpair}{\vect 0}{\cones}$, species 1 conquers a nondegenerate subcone of $\cones$ when $\cones$ is wide.

Now, if $\cones$ is critical for species 1, then $\cfavwedge{\normpair}{\vect 0}{\cones}$
contains some $\vect x\ne \vect 0$ by Lemma~\ref{max_adv_dir_lem}, in which case it contains the entire ray $\direction{\vect x}$ since the set $\cfavwedge{\normpair}{\vect 0}{\cones}$ is scale-invariant, so species 1 conquers at least an entire ray from $\vect 0$ when $\cones$ is critical.

Finally, suppose $\cones$ is narrow for species 1 and additionally is contained in a half-space. Since $\cones$ is narrow, we can choose some $\beta$ with $1<\beta < \coneadvantage{1}{\normpair}{\cones}^{-1}$. Since $\cones$ is contained in a half-space and $\beta < \coneadvantage{1}{\normpair}{\cones}^{-1}$, Lemma~\ref{unobtrusive_star_set_lem} below shows that $\cstarset[2]{\normpair}{\beta}{2}{\vect 0}{\Rd\setminus \cones} = \Rd\setminus \setof{\vect 0}$, and since $\beta>1$, Proposition~\ref{determ_star_conquering_prop} then implies that species 1 conquers no point in $\Rd\setminus \setof{\vect 0}$. Thus we must have $\conqset{1}{\normpair}{\vect 0}{\cones} = \setof{\vect 0}$.
%
%
%
\end{proof}

\begin{thmremark}
The proof of Proposition~\ref{conquering_from_apex_prop} shows that Part~\ref{conquering_from_apex:crit_part} can be improved to say that species 1 conquers every critical direction in an arbitrary cone (though it ties with species 2 in these directions).
In a small critical cone, the set of critical directions should typically be a degenerate cone whose dimension depends on the specific geometry (though we do not attempt to prove this statement). The general statement in Part~\ref{conquering_from_apex:crit_part} merely states that the dimension of the conquered subcone is at least 1 in any critical cone, which may be suboptimal. In fact, for certain choices of $\normpair$, it is possible to find large critical cones in which species 1 conquers an entire nondegenerate subcone.
\end{thmremark}

\begin{thmremark}
The conclusion of Part~\ref{conquering_from_apex:narr_part} of Proposition~\ref{conquering_from_apex:narr_part} should also be true for cones that are not contained in a half-space. However, a different proof is needed because Lemma~\ref{unobtrusive_star_set_lem} below does not necessarily hold for large cones---if the cone is not contained in a half-space, species 2 may need to follow a curved path from the boundary (similar to the logarithmic spiral for Euclidean norms) in order to block species 1 at the apex.
\end{thmremark}


\begin{lem}[Species 2's star-closer set in a small cone]
\label{unobtrusive_star_set_lem}
Let $\cones$ be a cone at $\vect 0$ which is contained in some half-space. Then for any $\beta<\coneadvantage{1}{\normpair}{\cones}^{-1}$,
\[
\cstarset[2]{\normpair}{\beta}{2}{\vect 0}{\Rd\setminus \cones}
= \Rd\setminus \setof{\vect 0}.
\]
\end{lem}

\begin{proof}
\newcommand{\bluntcone}{\cones\setminus \setof{\vect 0}}
\newcommand{\conept}{\vect x} 
\newcommand{\bdpt}{\vect y_{\conept}} 

\newcommand{\ballpt}{\bdpt'} 
\newcommand{\segpt}{\vect y} 

\newcommand{\radius}{\epsilon_{\conept}} 
\newcommand{\thisball}{\reacheddset{\mu_2}{\bdpt}{\radius}}

If $\beta = 0$ the statement is trivial, so assume $\beta>0$. Then $\beta^{-1} > \coneadvantage{1}{\normpair}{\cones}$, so by Lemma~\ref{strong_blocking_segments_lem}, for any $\conept\in \bluntcone$ there exists a $\bdpt\in \bd \cones$ and $\radius>0$ such that for all $\ballpt\in \thisball$,
\begin{equation}
\label{unobtrusive_star_set:sup_eqn}
\sup_{\segpt\in \segment{\ballpt}{\conept}}
\frac{\distnew{\mu_2}{\ballpt}{\segpt}}{\distnew{\mu_1}{\vect 0}{\segpt}}
< \beta^{-1}.
\end{equation}
Since $\bdpt \in \bd \cones$, the ball $\reacheddset{\mu_2}{\bdpt}{\radius}$ intersects $\Rd\setminus \cones$, so choose any $\ballpt$ in this intersection. Then $\ballpt\in \Rd\setminus \cones$, and $\ballpt$  satisfies \eqref{unobtrusive_star_set:sup_eqn}, or equivalently,
\[
\forall \segpt\in \segment{\ballpt}{\conept},\qquad
\beta \distnew{\mu_2}{\ballpt}{\segpt} <  \distnew{\mu_1}{\vect 0}{\segpt}.
\]
This shows that $\segment{\ballpt}{\conept} \subseteq \closerset[2]{\normpair}{\beta}{2}{\vect 0}{\ballpt}$ and hence 
\[
\vect x\in \cstarset[2]{\normpair}{\beta}{2}{\vect 0}{\ballpt} \subseteq \cstarset[2]{\normpair}{\beta}{2}{\vect 0}{\Rd\setminus \cones}.\qedhere
\]
\end{proof}

Before proving a more general version of Proposition~\ref{conquering_from_apex_prop} in Section~\ref{deterministic:cones:conquered_regions_sec}, we first prove some additional properties of the speed advantage in Section~\ref{deterministic:advantage_properties_sec} and then generalize the notion of advantage in Section~\ref{deterministic:gen_advantage_sec}.

\subsection{Properties of the Speed Advantage}
\label{deterministic:advantage_properties_sec}

\begin{lem}[Monotonicity of the speed advantage]
\label{advantage_monotone_lem}
\newcommand{\bigcone}{\cones}
\newcommand{\smallcone}{\cones'}

\newcommand{\bigadv}{\diradvantage{1}{\normpair}{\bigcone}{\dirvec}}
\newcommand{\smalladv}{\diradvantage{1}{\normpair}{\smallcone}{\dirvec}}

\newcommand{\bigAdv}{\coneadvantage{1}{\normpair}{\bigcone}}
\newcommand{\smallAdv}{\coneadvantage{1}{\normpair}{\smallcone}}

Let $\bigcone$ be a closed cone in $\Rd$ and let $\smallcone$ be a subcone of $\bigcone$ (possibly with different apex set). Then $\bigadv \ge \smalladv$ for any $\dirvec\in\dirset{\smallcone}$, and $\bigAdv\ge \smallAdv$.
\end{lem}

\begin{proof}

\newcommand{\bigcone}{\cones}
\newcommand{\smallcone}{\cones'}

\newcommand{\thisx}{\vect x}
\newcommand{\thisa}{\vect a}
\newcommand{\origin}{\vect 0}

\newcommand{\thishom}{\homothe{r}{\thisa}}
\newcommand{\rpoint}{\thishom \thisx}

\newcommand{\rayrep}{\dirsymb}
\newcommand{\rrayrep}{\dirsymb_r}
\newcommand{\raydir}{\dirvec}
\newcommand{\rraydir}{\dirvec_r}

\newcommand{\bigadv}{\diradvantage{1}{\normpair}{\bigcone}{\raydir}}
\newcommand{\smalladv}{\diradvantage{1}{\normpair}{\smallcone}{\raydir}}
\newcommand{\rbigadv}{\diradvantage{1}{\normpair}{\bigcone}{\rraydir}}

Note that since $\bigcone$ is closed, $\dirset{\bigcone}\supseteq \dirset{\smallcone}$ by Lemma~\ref{dirset_monotone_lem}, so $\bigadv$ is well-defined.
Without loss of generality, assume $\bigcone$ has apex $\origin$. Let $\thisa$ be an apex of $\smallcone$, and let $\raydir\in \dirset{\smallcone}$. Then if $\rayrep$ is any representative of $\raydir$, there is some $\thisx\in \smallcone\setminus \setof{\thisa}$ such that $\rayrep = \thisx-\thisa$. For $r>0$, let $\thishom$ be the homothety with scale $r$ and center $\thisa$, and let $\rrayrep \definedas \rpoint = \thisa + r(\thisx-\thisa) \in \smallcone$. Then it follows from the definition \eqref{dir_advantage_def_eqn} plus translation-invariance and scale-equivariance of the norm metrics that
\begin{equation}
\label{advantage_monotone:smalladv_eqn}
\smalladv
= \frac{\distnew{\mu_2}{\siversion{\smallcone}}{\rayrep}}
	{\distnew{\mu_1}{\origin}{\rayrep}}
= \frac{\distnew{\mu_2}{\smallcone}{\rrayrep}}
	{\distnew{\mu_1}{\thisa}{\rrayrep}}
\quad\text{ for any $r>0$}.
\end{equation}
Now observe that $\frac{1}{r} \rrayrep = \frac{1}{r}\thisa + (\thisx-\thisa) \to \thisx - \thisa = \rayrep$ as $r\to\infty$, and since $\frac{1}{r} \rrayrep \direquiv \rrayrep$ for any $r>0$, we have $\rraydir \to \raydir$ in $\dirset{\bigcone}$. By continuity, we thus have 
\begin{equation}
\label{advantage_monotone:rbigadv_conv_eqn}
\lim_{r\to\infty} \rbigadv = \bigadv.
\end{equation}
Next, since $\bigcone\supseteq\smallcone$ and $\rrayrep\in\smallcone$, we have $\distnew{\mu_2}{\bd \bigcone}{\rrayrep} \ge \distnew{\mu_2}{\bd \smallcone}{\rrayrep}$.
Using this inequality and the triangle inequality for $\mu_1$, and then using \eqref{advantage_monotone:smalladv_eqn}, we have
\begin{align}
\label{advantage_monotone:rbigadv_larger_eqn}
\rbigadv
= \frac{\distnew{\mu_2}{\bd \bigcone}{\rrayrep}}
	{\distnew{\mu_1}{\origin}{\rrayrep}}
\ge \frac{\distnew{\mu_2}{\bd \smallcone}{\rrayrep}}
	{\distnew{\mu_1}{\thisa}{\rrayrep} + \mu_1(\thisa)} 
&= \frac{\distnew{\mu_2}{\bd \smallcone}{\rrayrep}}
		{\distnew{\mu_1}{\thisa}{\rrayrep} }
	\cdot \frac{1}{1+\frac{\mu_1(\thisa)}{\distnew{\mu_1}{\thisa}{\rrayrep}}} \notag\\[5mm]
&=\smalladv \cdot \parens[1]{1+\tfrac{\mu_1(\thisa)}{r\mu_1(\thisx -\thisa)}}^{-1}.
\end{align}
Combining \eqref{advantage_monotone:rbigadv_conv_eqn} and \eqref{advantage_monotone:rbigadv_larger_eqn} we get
\[
\bigadv = \lim_{r\to\infty} \rbigadv
\ge \smalladv \cdot \lim_{r\to\infty} 
	\parens[1]{1+\tfrac{1}{r}\cdot \tfrac{\mu_1(\thisa)}{\mu_1(\thisx -\thisa)}}^{-1}
=\; \smalladv.
\]
This proves the first statement, and the second statement then follows trivially from the definition \eqref{cone_advantage_def_eqn}.
\end{proof}

Note that Lemma~\ref{max_adv_dir_lem} implies that $\coneadvantage{1}{\normpair}{\cones} = \coneadvantage{1}{\normpair}{\closure{\cones}}$ for any cone $\cones$, so the conclusion of Lemma~\ref{advantage_monotone_lem} is also true for non-closed cones $\cones$ as long as $\dirvec$ is assumed to be in $\dirset{\cones}$. 

The following lemma gives explicit bounds on the advantage in a $\mu_2$-cone, improving the trivial bounds given in Lemma~\ref{max_adv_dir_lem} above. This result has several useful consequences.

\begin{lem}[The advantage in a small $\mu_2$-cone]
\label{mu2_cone_adv_lem}
\newcommand{\thiscone}{\cone{\mu_2}{\thp}{\vect a}{\dirvec}}
For any $\vect a\in\Rd$, $\dirsymb\in\Rd\setminus \setof{\vect 0}$, and $\thp\in [0,1]$,
\[
\thp\, \frac{\mu_2(\dirsymb)}{\mu_1(\dirsymb)}
\le \coneadvantage[1]{1}{\normpair}{\thiscone}
\le \thp\, \esupnormon[1]{\thiscone}{\frac{\mu_2}{\mu_1}}.
\]
The three quantities are equal if $\dirvec$ is any direction that maximizes $\mu_2/\mu_1$.
\end{lem}

\begin{proof}

\newcommand{\thisapex}{\vect 0}
\newcommand{\thiscone}{\cone{\mu_2}{\thp}{\thisapex}{\dirvec}}
\newcommand{\thisconeadv}{\coneadvantage{1}{\normpair}{\cones}}
\newcommand{\maxpt}{\dirsymb[1]}
\newcommand{\maxdir}{\dirvec[1]}
\newcommand{\genpt}{\vect x}

If $\thp=0$, the cone is degenerate and the result follows from Lemma~\ref{max_adv_dir_lem}, so suppose $\thp>0$. Assume $\anapex = \thisapex$, and let $\cones \definedas \thiscone$. Then $\dirsymb\in \cones\setminus \setof{\thisapex}$, so we have
\[
\thisconeadv
\ge \frac{\distnew{\mu_2}{\bd \cones}{\dirsymb}}{\distnew{\mu_1}{\thisapex}{\dirsymb}}
= \frac{\thp\, \mu_2(\dirsymb)}{\mu_1(\dirsymb)},
\]
where the equality $\distnew{\mu_2}{\bd \cones}{\dirsymb} = \thp \mu_2(\dirsymb)$ holds by Lemma~\ref{mu_cone_bd_lem}. This proves the lower bound. The upper bound relies upon the following geometric claim.

\begin{claim}
\label{mu2_cone_adv:closest_pt_claim}
\newcommand{\thispt}{\dirsymb[1]}
If $\thispt\in\cones = \thiscone$ and $\distnew{\mu_2}{\bd \cones}{\thispt}\ge \distnew{\mu_2}{\bd \cones}{\dirsymb}$, then $\mu_2(\thispt) \ge \mu_2(\dirsymb)$.
\end{claim}

\begin{proof}[Proof of Claim~\ref{mu2_cone_adv:closest_pt_claim}]
\newcommand{\thisdist}{\thp}
\newcommand{\thisset}{\cones_\thp}
\newcommand{\thispt}{\dirsymb[1]}
\newcommand{\newcone}{\cone{\mu_2}{\thp}{\dirsymb}{\dirvec}}
\newcommand{\bigcone}{\cone{\mu_2}{1}{\dirsymb}{\dirvec}}
\newcommand{\thisball}{\reacheddset{\mu_2}{\dirsymb +\dirsymb}{1}}
\newcommand{\flippedball}{\reacheddset{\mu_2}{\dirsymb - \dirsymb}{1}}
\newcommand{\flippedballo}{\reacheddset{\mu_2}{\vect 0}{1}}

By the scale-equivariance of $\mu_2$ and $\distnew{\mu_2}{\bd \cones}{\cdot}$ it suffices to assume $\mu_2(\dirsymb)=1$.
Then $\distnew{\mu_2}{\bd \cones}{\dirsymb} = \thisdist$ by Lemma~\ref{mu_cone_bd_lem}, so the claim is equivalent to the statement that $\dirsymb$ is a $\mu_2$-closest point to $\thisapex$ in the set $\thisset \definedas \setbuilder[2]{\thispt\in \cones}{\distnew{\mu_2}{\bd\cones}{\thispt} \ge \thisdist}$. By Lemma~\ref{special_mu_properties_lem} we have $\thisset =  \frac{\thp}{\thp} \cdot \dirsymb+ \thiscone = \newcone$.
Since $\thp\le 1$ we have $\newcone \subseteq \bigcone = \homothe{(\Rplus)}{\dirsymb} \cdot \thisball \subseteq \Hspace$, where $\Hspace$ is any support half-space of $\thisball$ at $\dirsymb$. By the symmetry of $\mu_2$ (i.e.\ evenness), $\Hspace$ must be an opposing half-space of $\flippedball = \flippedballo$ at $\dirsymb$, since $\flippedball$ is the reflection of $\thisball$ through the boundary point $\dirsymb$. Lemma~\ref{support_point_closest_lem} then implies that $\dirsymb$ is a $\mu_2$-closest point to $\vect 0$ in $\bd \Hspace$, hence in $\Hspace$ as well by Lemma~\ref{closest_bd_point_lem}. Since $\dirsymb\in \thisset\subseteq \Hspace$, $\dirsymb$ must be a $\mu_2$-closest point to $\vect 0$ in $\thisset$.
\end{proof}

Now, by Lemma~\ref{max_adv_dir_lem} there exists $\maxdir\in \dirset{\interior{\cones}}$ with $\thisconeadv = \diradvantage{1}{\normpair}{\cones}{\maxdir}$. Choose a representative $\maxpt\in\maxdir$ with $\distnew{\mu_2}{\bd\cones}{\maxpt} = \distnew{\mu_2}{\bd\cones}{\dirsymb}$; note that this is possible because for any representative $\maxpt_0\in\maxdir$, we can take $\maxpt \definedas \maxpt_0 \cdot \frac{\distnew{\mu_2}{\bd\cones}{\dirsymb}}{\distnew{\mu_2}{\bd\cones}{\maxpt_0}}\in \maxdir$. Then, using Claim~\ref{mu2_cone_adv:closest_pt_claim} in the second-to-last step, we have
\[
\thisconeadv
= \frac{\distnew{\mu_2}{\bd\cones}{\maxpt}}{\mu_1(\maxpt)}
= \frac{\distnew{\mu_2}{\bd\cones}{\dirsymb}}{\mu_1(\maxpt)}
= \frac{\thp \mu_2(\dirsymb)}{\mu_1(\maxpt)}
\le \frac{\thp \mu_2(\maxpt)}{\mu_1(\maxpt)}
\le \sup_{\genpt\in \cones\setminus \setof{\vect 0}}
\frac{\thp \mu_2(\genpt)}{\mu_1(\genpt)}.
\]
This proves the upper bound. The final statement then follows trivially since if $\dirvec$ maximizes $\mu_2/\mu_1$, then the upper and lower bounds must be equal.
\end{proof}

The first corollary of Lemma~\ref{mu2_cone_adv_lem} gives us another way of interpreting the separation thickness $\thpadv{2}{\dirvec}{\alpha}$ defined in \eqref{thpadv_def_eqn}.

\begin{cor}[Characterization of $\thpadv{2}{\dirvec}{\alpha}$ in terms of advantage]
\label{thp_adv_cor}
\newcommand{\thisthp}{\thpadv{2}{\dirvec}{\alpha}}
For $\alpha>0$ and $\dirsymb\in\Rd\setminus \setof{\vect 0}$, let $\thisthp =  \alpha\frac{\mu_1(\dirsymb)}{\mu_2(\dirsymb)}$ as defined in \eqref{thpadv_def_eqn}. Then
\[
\coneadvantage[1]{1}{\normpair}{\cone{\mu_2}{\thisthp}{\vect a}{\dirvec}} \ge \alpha.
\]
\end{cor}

\begin{proof}
This follows immediately from Lemma~\ref{mu2_cone_adv_lem} since $\thpadv{2}{\dirvec}{\alpha} \cdot \frac{\mu_2(\dirsymb)}{\mu_1(\dirsymb)} = \alpha$. 
\end{proof}

The next corollary guarantees that if $\esupnorm{\mu_2/\mu_1}>1$, then there exist wide convex  closed cones that are extra-small (i.e.\ interior angles at the apex are strictly less than $\pi$).

\begin{cor}[Existence of extra-small wide cones]
\label{acute_wide_cone_cor}
Let $\dirsymb\in \Rd\setminus \setof{\vect 0}$. If $\mu_1(\dirsymb)<\mu_2(\dirsymb)$, then there exists $\thp<1$ such that the $\mu_2$-cone $\cone{\mu_2}{\thp}{\vect 0}{\dirvec}$ is wide for species 1, and this cone is extra-small by Lemma~\ref{mu_cone_description_lem}.
\end{cor}

\begin{proof}
Choose any $\thp$ with $\frac{\mu_1(\dirsymb)}{\mu_2(\dirsymb)} <\thp <1$.
Then $\thp = \thpadv{2}{\dirvec}{\alpha}$, where $\alpha \definedas \thp \frac{\mu_2(\dirsymb)}{\mu_1(\dirsymb)} >1$, so the result follows from Corollary~\ref{thp_adv_cor}.
\end{proof}

The following result will be relevant in Chapter~\ref{coex_finite_chap}.

\begin{cor}[Advantage for norms that are scalar multiples of one another]
\label{exact_adv_cor}
Suppose $\mu_1 = \mu$ and $\mu_2 = \lambda^{-1}\mu$ for some norm $\mu$ and some $\lambda>0$ (meaning that species 2 is $\lambda$ times as fast as species 1 in every direction). Let $\normpair = \pair{\mu_1}{\mu_2} = \pair[2]{\mu}{\,\lambda^{-1}\mu}$.
Then for any $\anapex\in\Rd$, $\dirvec\in\dirset{\Rd}$, and $\thp\in [0,1]$, we have $\cone{\mu_1}{\thp}{\anapex}{\dirvec} = \cone{\mu_2}{\thp}{\anapex}{\dirvec} = \cone{\mu}{\thp}{\anapex}{\dirvec}$, and
\[
\coneadvantage[1]{1}{\normpair}{\cone{\mu}{\thp}{\anapex}{\dirvec}}
= \frac{\thp}{\lambda}.
\]
\end{cor}

\begin{proof}
In this case every direction $\dirvec$ maximizes $\mu_2/\mu_1 \equiv \lambda^{-1}$, so the three terms in Lemma~\ref{mu2_cone_adv_lem} are equal for any $\dirvec$.
\end{proof}

Combining Lemmas~\ref{max_adv_dir_lem} and \ref{advantage_monotone_lem}, we get the following useful result, which shows that the lower bound in Lemma~\ref{mu2_cone_adv_lem} is always achieved in some direction.

\begin{lem}[Existence of a $\mu_2$-subcone with equal advantage]
\label{mu2_subcone_lem}
%
%
\newcommand{\maxthp}{\thp}
\newcommand{\genthp}{\thp'}
\newcommand{\subcone}{\cones'}
\newcommand{\mutwocone}[1]{\cone{\mu_2}{#1}{\anapex}{\dirvec}}

Let $\cones\subseteq \Rd$ be a closed cone with apex $\anapex\in\Rd$, and let $\dirvec\in \dirset{\cones}$.
\begin{enumerate}
\item \label{mu2_subcone:dir_part}
Let $\maxthp =\distnew[2]{\mu_2}{\bd \siversion{\cones}}{\munitvec[2]}$, and let $\subcone = \mutwocone{\maxthp}$. Then $\subcone\subseteq \cones$, and $\diradvantage{1}{\normpair}{\subcone}{\dirvec} = \diradvantage{1}{\normpair}{\cones}{\dirvec}$.




\item \label{mu2_subcone:cone_part}
Suppose $\coneadvantage{1}{\normpair}{\cones} = \alpha$, and let $\thpadv{2}{\dirvec}{\alpha} = \alpha\frac{\mu_1(\dirsymb)}{\mu_2(\dirsymb)}$ as in \eqref{thpadv_def_eqn}. If $\dirvec$ achieves the supremum in \eqref{cone_advantage_def_eqn}, then the $\mu_2$-cone $\cone{\mu_2}{\thpadv{2}{\dirvec}{\alpha}}{\anapex}{\dirvec}$ is contained in $\cones$ and has advantage $\alpha$.
\end{enumerate}
\end{lem}

\begin{proof}
\newcommand{\maxthp}{\thp}
\newcommand{\mmaxthp}{\thpadv{2}{\dirvec}{\alpha}}
\newcommand{\subcone}{\cones'}
\newcommand{\thisapex}{\vect 0}
\newcommand{\mutwocone}[1]{\cone{\mu_2}{#1}{\thisapex}{\dirvec}}
\newcommand{\thisball}{\reacheddset{\mu_2}{\munitvec[2]}{\maxthp}}

Assume without loss of generality that $\anapex =\vect 0$ and hence $\siversion{\cones} = \cones$. For Part~\ref{mu2_subcone:dir_part}, suppose $\distnew[2]{\mu_2}{\bd \cones}{\munitvec[2]} = \maxthp$. Then since $\cones$ is closed, we must have $\thisball\subseteq \cones$, because otherwise there would be some $\vect y\in \bd \cones$ with $\distnew[2]{\mu_2}{\vect y}{\munitvec[2]} < \maxthp$. Since $\cones$ is scale-invariant and closed, this implies that $\cones$ contains the entire scale-invariant set $\Rplus \cdot \thisball = \mutwocone{\maxthp} \defines \subcone$. Moreover, by definition \eqref{dir_advantage_def_eqn} and Lemma~\ref{mu_cone_bd_lem}, we have
\[
\diradvantage{1}{\normpair}{\subcone}{\dirvec}
= \frac{\distnew{\mu_2}{\bd \subcone}{\munitvec[2]}}{\mu_1(\munitvec[2])}
=\frac{\maxthp}{\mu_1(\munitvec[2])}
=\frac{\distnew{\mu_2}{\bd \cones}{\munitvec[2]}}{\mu_1(\munitvec[2])}
=\diradvantage{1}{\normpair}{\cones}{\dirvec}.
\]
For Part~\ref{mu2_subcone:cone_part}, suppose $\diradvantage{1}{\normpair}{\cones}{\dirvec} = \coneadvantage{1}{\normpair}{\cones} = \alpha$ (note that there exists some such $\dirvec\in \dirset{\cones}$ by Lemma~\ref{max_adv_dir_lem}). Then $\distnew{\mu_2}{\bd \cones}{\dirsymb} = \alpha \mu_1(\dirsymb)$ for any $\dirsymb\in\dirvec$, so $\distnew{\mu_2}{\bd \cones}{\munitvec[2]} = \alpha \frac{\mu_1(\dirsymb)}{\mu_2(\dirsymb)} = \mmaxthp$. By Part~\ref{mu2_subcone:dir_part} we then have  $\subcone \definedas \mutwocone{\mmaxthp} \subseteq \cones$ and $\diradvantage{1}{\normpair}{\subcone}{\dirvec} = \diradvantage{1}{\normpair}{\cones}{\dirvec} = \alpha$. Therefore, $\coneadvantage{1}{\normpair}{\subcone} \ge \alpha$, and since $\subcone\subseteq \cones$, we also have $\coneadvantage{1}{\normpair}{\subcone} \le \alpha$ by Lemma~\ref{advantage_monotone_lem}.
\end{proof}


The following lemma gives some other useful ways of thinking about the advantage in a cone; for example, the last two properties can be visualized as packing a very solid, $\mu_2$-ball-shaped scoop of ice cream as tightly as possible into the ``ice cream" cone $\cones$. The proof is left to the reader.

\begin{lem}[Equivalent characterizations of the speed advantage]
\label{advantage_characterizations_lem}
Let $\cones$ be a nondegenerate cone in $\Rd$ with apex $\vect 0$, and let $\normpair = \pair{\mu_1}{\mu_2}$ be a traversal norm pair. Then
\begin{enumerate}

\item $\coneadvantage{1}{\normpair}{\cones}\ge \beta$ if and only if $\interior{\cones} \cap \closerset[1]{\normpair}{\beta}{1}{\vect 0}{\,\Rd\setminus \cones}$ is nonempty.

\item $\coneadvantage{1}{\normpair}{\cones} = \sup \setbuilder[2]{\alpha}
	{\exists \vect x\in \unitball{\mu_1}\cap \cones \text{\rm\ with } \distnew{\mu_2}{\bd \cones}{\vect x}>\alpha}$.
\item $\coneadvantage{1}{\normpair}{\cones} = \sup \setbuilder[2]{\alpha}{\reacheddset{\mu_2}{\vect x}{\alpha} \subseteq \cones  \text{ \rm for some } \vect x\in \unitball{\mu_1} }.$

\item $\coneadvantage{1}{\normpair}{\cones}^{-1} = \inf \setbuilder[2]{\mu_1(\vect x)}{\vect x +\unitball{\mu_2} \subseteq \cones}$.
\end{enumerate}
\end{lem}

\subsection{The Advantage of Species 1 in an Arbitrary Starting Configuration}
\label{deterministic:gen_advantage_sec}

We now generalize the idea of the advantage of species 1 in a cone to the advantage of species 1 in an arbitrary initial configuration $\initconfig$ of the deterministic process.
We will use this definition in the next section to analyze the process when species 1 starts at different locations within a given cone, and we will use a ``lattice-ized" version of this general definition in Chapter~\ref{coex_finite_chap} when analyzing the random process.

If $\normpair$ is a traversal norm pair, we define the \textdef{advantage of species 1} in an arbitrary initial configuration $\initconfig$ to be
\begin{equation}
\label{config_advantage_def_eqn}
\configadvantage{1}{\normpair}{\initialset{1}}{\initialset{2}}
\definedas \sup \setbuilderbar[1]{\coneadvantage{1}{\normpair}{\cones}}
{\textstack{$\cones$ is a cone with $\apexset{\cones}\cap \initialset{1} \ne \emptyset$}{and
	$\cones\cap \initialset{2} = \emptyset$}}.
\end{equation}
With this definition, it follows from Lemma~\ref{advantage_monotone_lem} (monotonicity) that if $\cones$ is a cone with apex $\vect a$, then
\begin{equation}
\label{cone_advantage_def2_eqn}
\coneadvantage{1}{\normpair}{\cones} =
\configadvantage[2]{1}{\normpair}{\vect a}{\Rd\setminus \cones},
\end{equation}
so we can view \eqref{config_advantage_def_eqn} as a generalization of the definition of the advantage in a cone. Moreover, it follows from Lemma~\ref{mu2_subcone_lem} that for any initial configuration $\initconfig$,
\begin{equation}
\label{config_advantage_def2_eqn}
\configadvantage{1}{\normpair}{\initialset{1}}{\initialset{2}}
= \sup \setbuilderbar[1]{\alpha \ge 0}
{\textstack{$\exists \dirsymb\in \Rd\setminus \setof{\vect 0}$ and $\vect a\in \initialset{1}$ with}
	 {$\cone{\mu_2}{\thpadv{2}{\dirvec}{\alpha}}{\vect a}{\dirvec} \cap \initialset{2} = \emptyset$,
	 where $\thpadv{2}{\dirvec}{\alpha} = \alpha \frac{\mu_1(\dirsymb)}{\mu_2(\dirsymb)}$}}.
\end{equation}
Thus, we can take \eqref{config_advantage_def2_eqn} as the most general definition of advantage, and treat the advantage in a cone as a special case defined by \eqref{cone_advantage_def2_eqn}. Just as we classified cones in Definition~\ref{cone_criticality_def}, we can classify an arbitrary starting configuration $\initconfig$ as \textdef{wide}, \textdef{narrow}, or \textdef{critical} for species 1 according to whether $\configadvantage{1}{\normpair}{\initialset{1}}{\initialset{2}}$ is greater than, less than, or equal to 1, respectively.

\begin{remark}
We mention that if one wants to analyze the geometry of the (random or deterministic) first-passage competition process in more depth, it might be more appropriate to redefine the advantage of species 1 using ``near-isometric" embeddings of cones into $\Rd$, e.g.\ ``bent" or ``twisted" cones, rather than using only ``straight" cones, because in general species 1 can escape by following curved paths rather than straight lines. It may also be fruitful to define the ``local advantage" of species 1 at an arbitrary location $\vect a\in \initialset{1}$, perhaps by homothetically zooming in on the configuration at $\vect a$, and taking the local advantage at $\vect a$ to be the advantage of a maximal cone in the scale-invariant limit. However, the definition \eqref{config_advantage_def_eqn} or \eqref{config_advantage_def2_eqn} will be sufficient for our purposes.
\end{remark}

\begin{lem}[The configuration advantage in cones]
\label{convex_advantage_lem}
\newcommand{\conept}{\vect z}
\newcommand{\thisconfigadv}{\configadvantage[2]{1}{\normpair}{\conept}{\,\Rd\setminus \cones}}
\newcommand{\thisconeadv}{\coneadvantage{1}{\normpair}{\cones}}

If $\cones\subseteq\Rd$ is a cone and $\conept\in\cones$, then $\thisconfigadv\le \thisconeadv$. If $\cones$ is convex, then $\thisconfigadv = \thisconeadv$ for all $\conept\in\cones$.
\end{lem}

\begin{proof}
\newcommand{\conept}{\vect z}
\newcommand{\thisadv}{\alpha}
\newcommand{\thiscone}{\cones}
\newcommand{\closedcone}{\closure{\cones}}
\newcommand{\opencone}{\interior{\cones}}
\newcommand{\origin}{\vect 0}
\newcommand{\thisapex}{\origin}

\newcommand{\thptwo}{\thpadv{2}{\dirvec}{\thisadv}}
\newcommand{\mutwocone}[1]{\cone{\mu_2}{\thptwo}{#1}{\dirvec}}
\newcommand{\omutwocone}[1]{\cone{\mu_2}{\thptwo^-}{#1}{\dirvec}}

\newcommand{\thisconfigadv}{\configadvantage[2]{1}{\normpair}{\conept}{\,\Rd\setminus \cones}}
\newcommand{\thisconeadv}{\coneadvantage{1}{\normpair}{\cones}}

\newcommand{\subcone}{\cones'}

Since cones are bodies and $\coneadvantage{1}{\normpair}{\thiscone} = \coneadvantage{1}{\normpair}{\closure{\thiscone}}$ by Lemma~\ref{max_adv_dir_lem}, it suffices to assume $\thiscone$ is closed.  First note that it follows directly from \eqref{config_advantage_def_eqn} and Lemma~\ref{advantage_monotone_lem} that
\begin{equation*}
\thisconfigadv
= \sup \setbuilderbar[1]{\coneadvantage{1}{\normpair}{\subcone}}
{\textstack{$\cones$ is a cone with $\conept\in \apexset{\subcone}$}{and
	$\subcone\subseteq \thiscone$}}
\le \thisconeadv,
\end{equation*}
which proves the first statement. Now suppose $\thiscone$ is convex.
Assume $\thisapex\in \apexset{\thiscone}$, and let $\thisadv = \coneadvantage{1}{\normpair}{\thiscone}$.
If $\thisadv = 0$, the statement is trivial, so assume $\thisadv>0$. 
Then by Lemmas~\ref{max_adv_dir_lem} and \ref{mu2_subcone_lem}, there exists $\dirvec\in \dirset{\opencone}$ such that the $\mu_2$-cone $\mutwocone{\thisapex}$ is contained in $\thiscone$ and has advantage $\thisadv$.
Since $\thiscone$ is convex with apex $\thisapex$, we have $\conept+\thiscone \subseteq \thiscone$ for any $\conept\in\thiscone$ by Part~\ref{wedge_properties:convexity_part} of Lemma~\ref{wedge_properties_lem}.  Therefore,
\[
\mutwocone{\conept} = \conept+\mutwocone{\thisapex}
\subseteq \conept + \thiscone \subseteq \thiscone,
\]
so we have $\thisconfigadv \ge \thisadv$ by \eqref{config_advantage_def2_eqn}. The first statement then implies that $\thisconfigadv = \thisadv$, which proves the second statement.
%
%
%
\end{proof}

\subsection{The Conquered Regions in Wide, Critical, and Narrow Cones}
\label{deterministic:cones:conquered_regions_sec}

Our main goal in this section will be to prove the following proposition, which generalizes Proposition~\ref{conquering_from_apex_prop} to the case where species~1 starts at an arbitrary location within the cone, not necessarily the apex. We break the proof into several smaller results which will be treated individually after the main statement. When we study the random two-type process in Chapter~\ref{random_fpc_chap}, the results of Section~\ref{random_fpc:cone_competition_sec} will parallel those in the present section, providing stochastic analogues for the deterministic growth described here. The main result of that section, Theorem~\ref{wide_narrow_survival_thm}, is the analogue of Proposition~\ref{conquering_from_point_prop}.

\begin{prop}[The conquered region in a cone]
\label{conquering_from_point_prop}

\newcommand{\conept}{\vect z}
\newcommand{\Atwo}{\Rd\setminus\cones}
\newcommand{\thisconfig}{\pair[2]{\conept}{\Atwo}}
\newcommand{\thisadv}{\configadvantage[1]{1}{\normpair}{\conept}{\Atwo}}

Let $\cones$ be a cone in $\Rd$, let $\conept\in\cones$, and  consider a $\normpair$-process with starting configuration $\thisconfig$.
\begin{enumerate}
\item \label{conquering_from_point:wide_part}
If $\cones$ is wide for species 1, and $\thisadv > 1$, then species 1 conquers at least a nondegenerate subcone of $\cones$ with apex $\conept$. Moreover, the thickness of the conquered cone is bounded below by a positive constant that depends only on $\thisadv$.

\item \label{conquering_from_point:crit_part}
If $\cones$ is wide or critical for species 1, and $\thisadv = 1$, then species 1 conquers at least a ray from $\conept$.

\item \label{conquering_from_point:narr_part}
If $\cones$ is narrow for species 1 and additionally is contained in a half-space, then species 1 conquers only a bounded set, no matter which $\conept\in \cones$ is chosen.
\end{enumerate}
\end{prop}

\begin{proof}
Parts~\ref{conquering_from_point:wide_part} and \ref{conquering_from_point:crit_part} follow from Proposition~\ref{wide_crit_survival_prop} below, and Part~\ref{conquering_from_point:narr_part} follows from Proposition~\ref{narrow_extinction_prop} below.
\end{proof}

\begin{thmremark}
\newcommand{\conept}{\vect z}
\newcommand{\Atwo}{\Rd\setminus\cones}
\newcommand{\thisconfig}{\pair[2]{\conept}{\Atwo}}
\newcommand{\thisadv}{\configadvantage[1]{1}{\normpair}{\conept}{\Atwo}}

If $\cones$ is a convex cone, then $\thisadv = \coneadvantage{1}{\normpair}{\cones}$ for all $\conept\in\cones$ by Lemma~\ref{convex_advantage_lem}, so
Part \ref{conquering_from_point:wide_part} of Proposition~\ref{conquering_from_point_prop} applies to all points in a wide convex cone, and Part~\ref{conquering_from_point:crit_part} applies to all points in a critical convex cone. However, if $\cones$ is nonconvex, there may be $\conept\in\cones$ with $\thisadv < \coneadvantage{1}{\normpair}{\cones}$, so Parts~\ref{conquering_from_point:wide_part} and \ref{conquering_from_point:crit_part} of Proposition~\ref{conquering_from_point_prop} may not apply to every point in a nonconvex cone. In fact, it is fairly easy to construct wide nonconvex cones in which species 1 cannot survive from some starting locations.
\end{thmremark}

\begin{thmremark}
\newcommand{\conept}{\vect z}
\newcommand{\Atwo}{\Rd\setminus\cones}
\newcommand{\thisconfig}{\pair[2]{\conept}{\Atwo}}
\newcommand{\thisadv}{\configadvantage[1]{1}{\normpair}{\conept}{\Atwo}}

It should be relatively straightforward to improve Part~\ref{conquering_from_point:crit_part} of Proposition~\ref{conquering_from_point_prop} to conclude that species 1 in fact conquers at least an entire half-cylinder of positive radius if $\thisadv = 1$ and $\conept\notin \bd\cones$.
\end{thmremark}

The following result, which follows from Part~\ref{cone_closer_set:infinite_part} of Lemma~\ref{cone_closer_set_lem}, implies Parts~\ref{conquering_from_point:wide_part} and \ref{conquering_from_point:crit_part} of Proposition~\ref{conquering_from_point_prop}.

\begin{prop}[Survival from the wide and critical regions in a cone]
\label{wide_crit_survival_prop}

\newcommand{\conept}{\vect z}
\newcommand{\thisadv}{\alpha}
\newcommand{\thisconfigadv}{\configadvantage[2]{1}{\normpair}{\conept}{\Rd\setminus\cones}}
\newcommand{\thisconfig}{\pair[2]{\conept}{\Rd\setminus \cones}}
\newcommand{\thisthp}{\symbolref{\thp}{wide_crit_survival_prop} (\thisadv)}

Let $\cones\subseteq\Rd$ be a cone that is wide or critical for species 1, and let $\conept\in\cones$ with $\thisconfigadv = \thisadv\ge 1$. Then there exists $\thisthp\ge 0$, with $\thisthp>0$ if $\thisadv>1$, such that species 1 conquers  a $\mu_1$-cone of thickness $\thisthp$ at $\conept$ in the $\normpair$-process started from $\thisconfig$. If $\cones$ is convex and has advantage $\thisadv\ge 1$, then this statement holds for all $\conept\in\cones$.

\end{prop}


\begin{proof}
\newcommand{\conept}{\vect z}
\newcommand{\thisadv}{\alpha}
\newcommand{\thisconfigadv}{\configadvantage[2]{1}{\normpair}{\conept}{\Rd\setminus\cones}}
\newcommand{\thisconfig}{\pair[2]{\conept}{\Rd\setminus \cones}}
\newcommand{\thisthp}{\symbolref{\thp}{wide_crit_survival_prop} (\thisadv)}

\newcommand{\newadv}{\alpha'}
\newcommand{\oldthp}{\symbolref{\thp}{determ_conq_mu1_subcone_lem} (\newadv)}

\newcommand{\bigthp}{\thpadv{2}{\dirvec}{\newadv}}
\newcommand{\mutwocone}{\cone{\mu_2}{\bigthp}{\conept}{\dirvec}}
\newcommand{\muonecone}{\cone{\mu_1}{\oldthp}{\conept}{\dirvec}}

Let $\newadv\definedas \frac{1}{2} (1+\thisadv)$. Since $\newadv<\thisadv = \thisconfigadv$, there is by definition some $\dirvec\in\dirset{\Rd}$ such that $\mutwocone\subseteq \cones$ (note that then $\mutwocone \ne \Rd$, which implies that $\bigthp\le 1$ and hence $\frac{\mu_2(\dirsymb)}{\mu_1(\dirsymb)} \ge \newadv$). Since $\newadv\ge 1$, by \eqref{conq_sets_monotone_eqn} and Lemma~\ref{determ_conq_mu1_subcone_lem}, we then have
\[
\conqset{1}{\normpair}{\conept}{\cones}
\supseteq \conqset[1]{1}{\normpair}{\conept}{\,\mutwocone}
\supseteq \muonecone.
\]
Thus we can take $\thisthp \definedas \oldthp\ge 0$, and $\thisthp>0$ for $\thisadv>1$ since $\newadv>1$ in this case.
The statement about convex cones follows from Lemma~\ref{convex_advantage_lem}.
\end{proof}

Our next goal will be to prove Part~\ref{conquering_from_point:narr_part} of Proposition~\ref{conquering_from_point_prop}, showing that species~1 can't survive in small narrow cones. The main result is restated in Proposition~\ref{narrow_extinction_prop} below; its proof will require two preliminary lemmas.
Recall from \eqref{bpin_def_eqn} that for a norm $\mu$, the \textdef{$\mu$-bowling pin} with parameters $r\ge 0$, $\thp \ge 0$, $\vect y\in \Rd$, $\dirvec \in \dirset{\Rd}$, and $h\in [0,\infty]$ is defined by
\[
\bpin{\mu}{r}{\thp}{\vect y}{\dirvec}{h} \definedas
\reacheddset{\mu}{\vect y}{r} \cup \conetip{\mu}{\thp}{\vect y}{\dirvec}{h}.
\]
The point $\vect y$ is the \textdef{origin} of $\bpin{\mu}{r}{\thp}{\vect y}{\dirvec}{h}$, and $\bpin{\mu}{r}{\thp}{\vect y}{\dirvec}{h}$ is \textdef{nondegenerate} if the parameters $r$, $\thp$, and $h$ are all strictly positive.

\begin{lem}[Fattening a $\beta_0$-closer segment into a $\beta$-closer bowling pin]
\label{star_closer_pin_lem}
\newcommand{\xpt}{\vect x}
\newcommand{\ypt}{\vect y}
\newcommand{\xdirsymb}{\dirsymb_{\xpt}}
\newcommand{\xdir}{\dirvec_{\xpt}}
\newcommand{\xdist}{\mathsf{d}_{\xpt}}
\newcommand{\xheight}{h_{\xpt}}


Let $\initialset{1}\subseteq \Rd$, and let $\vect y\in \Rd \setminus \initialset{1}$ with $\distnew{\mu_1}{\initialset{1}}{\vect y} = r_0$. Fix $\beta_0\ge 0$ and $\vect x\in \cstarset{\normpair}{\beta_0}{2}{\initialset{1}}{\vect y} \setminus \setof{\vect y}$, and let $\xdirsymb = \vect x -\vect y$. Then
\begin{enumerate}
\item $\initialset{1}\cap \interior{\bpin[2]{\mu_1}{r_0}{\thpadv{1}{\xdir}{\beta_0}}{\vect y}{\xdir}{\mu_1(\xdirsymb)}} = \emptyset$, where $\thpadv{1}{\xdir}{\beta_0} = \beta_0 \frac{\mu_2(\xdirsymb)}{\mu_1(\xdirsymb)}$ (as in \eqref{thpadv_def_eqn}).

\item For any $\beta \in \segment{0}{\beta_0}$, we have
\[
\cstarset{\normpair}{\beta}{2}{\initialset{1}}{\vect y}
\supseteq \bpin[2]{\mu_2}{r}{\thp}{\vect y}{\xdir}{\mu_2(\xdirsymb)},
\text{ where $r= r_0\rmeet{2}{\beta}$ and $\thp = (\beta_0-\beta)\rmeet{2}{\beta}$.}
\]
The $\mu_2$-bowling pin is nondegenerate if $r_0>0$ and $\beta<\beta_0$.
\end{enumerate}
%
%
%
\end{lem}

\begin{proof}
Combine Lemmas~\ref{star_closer_char_lem}, \ref{determ_bullseye_lem}, and \ref{cone_closer_set_lem}, with the roles of species~1 and species~2 switched.
\end{proof}

\begin{lem}[Covering a spherical section with bowling pins]
\label{cross_section_covering_lem}
\newcommand{\thisapex}{\vect a}
\newcommand{\crosssection}{K}
\newcommand{\pincollection}{\boldsymbol{P}}
\newcommand{\interiorcollection}{\interior{\pincollection}}

\newcommand{\pin}{P}
\newcommand{\pininterior}{\interior{P}}
\newcommand{\pinorigin}{\vect y_{\pin}}

\newcommand{\nbhd}{U_{\thisapex}}

Let $\cones$ be a closed cone at $\thisapex\in\Rd$ that is contained in some half-space, and let $\crosssection$ be a spherical section of $\cones$. Then for any $\beta< \coneadvantage{1}{\normpair}{\cones}^{-1}$, there exists a finite collection $\pincollection$ of nondegenerate $\mu_2$-bowling pins such that
\begin{enumerate}
\item \label{cross_section_covering:origin_part}
The origin $\pinorigin$ of each $\pin\in \pincollection$ lies in $\Rd\setminus \cones$.

\item \label{cross_section_covering:cover_part}
The collection $\interiorcollection \definedas \setbuilder{\pininterior}{\pin \in \pincollection}$ forms an open cover of $\crosssection$.

\item \label{cross_section_covering:closer_part}
There is some open neighborhood $\nbhd$ of $\thisapex$ such that for every $\pin\in \pincollection$,
\[
\pin \subseteq \cstarset{\normpair}{\beta}{2}{\nbhd}{\pinorigin}.
\]

\item \label{cross_section_covering:thickness_part}
Each $\pin\in \pincollection$ has the same thickness $\thp$, which depends only on  $\beta$ and $\coneadvantage{1}{\normpair}{\cones}$.

\end{enumerate}
\end{lem}

\begin{proof}
\newcommand{\thisapex}{\vect a}
\newcommand{\crosssection}{K}
\newcommand{\pincollection}{\boldsymbol{P}}
\newcommand{\interiorcollection}{\interior{\pincollection}}

\newcommand{\pinsymb}{P} 

\newcommand{\pin}{\pinsymb}
\newcommand{\pininterior}{\interior{P}}
\newcommand{\pinorigin}{\vect y_{\pin}}

\newcommand{\nbhd}{U_{\thisapex}}

\newcommand{\bigbeta}{\beta_1}
\newcommand{\biggerbeta}{\beta_2}

\newcommand{\conecomp}{\Rd\setminus \cones}
\newcommand{\bluntcone}{\cones\setminus \setof{\thisapex}}

\newcommand{\conept}[1][]{\vect x_{#1}}
\newcommand{\sourcept}[1][]{\vect y_{\conept[#1]}}

\newcommand{\xdirsymb}{\dirsymb_{\conept}}
\newcommand{\xdir}[1][]{\dirvec_{\conept[#1]}}

\newcommand{\thicksymb}{\thp} 
\newcommand{\radsymb}{r} 

\newcommand{\bigpin}[1][]{\widehat{\pinsymb}_{\conept[#1]}}
\newcommand{\bigthick}[1][]{\widehat{\thicksymb}_{\conept[#1]}}
\newcommand{\bigrad}[1][]{\widehat{\radsymb}_{\conept[#1]}}

\newcommand{\biggerpin}[1][]{\widetilde{\pinsymb}_{\conept[#1]}}
\newcommand{\biggerthick}[1][]{\widetilde{\thicksymb}_{\conept[#1]}}
\newcommand{\biggerrad}[1][]{\widetilde{\radsymb}_{\conept[#1]}}

\newcommand{\thispin}[1][]{\pinsymb_{\conept[#1]}}
\newcommand{\thisthick}{\thicksymb}
\newcommand{\thisrad}[1][]{\radsymb_{\conept[#1]}}


Let $\bigbeta = \frac{2}{3} \beta + \frac{1}{3} \coneadvantage{1}{\normpair}{\cones}^{-1}$ and $\biggerbeta = \frac{1}{3} \beta + \frac{2}{3} \coneadvantage{1}{\normpair}{\cones}^{-1}$, so
$
\beta<\bigbeta<\biggerbeta < \coneadvantage{1}{\normpair}{\cones}^{-1}.
$
Then by Lemma~\ref{unobtrusive_star_set_lem}, since $\biggerbeta < \coneadvantage{1}{\normpair}{\cones}^{-1}$ we have
\[
\cstarset[2]{\normpair}{\biggerbeta}{2}{\thisapex}{\conecomp}
= \Rd\setminus \setof{\thisapex}.
\]
Thus, for every $\conept \in \bluntcone$, there is some $\sourcept\in \conecomp$ such that $\conept \in \cstarset{\normpair}{\biggerbeta}{2}{\thisapex}{\sourcept}$, which then implies $\conept \in \cstarset{\normpair}{\bigbeta}{2}{\thisapex}{\sourcept}$ since $\bigbeta<\biggerbeta$. For each $\conept\in \bluntcone$, let $\xdirsymb = \conept - \sourcept \ne \vect 0$, $\biggerrad = \distnew{\mu_1}{\thisapex}{\sourcept}>0$, and $\biggerthick = \thpadv{1}{\xdir}{\biggerbeta} = \biggerbeta \frac{\mu_2(\xdirsymb)}{\mu_1(\xdirsymb)}>0$, and define
\[
\biggerpin \definedas \bpin[2]{\mu_1}{\biggerrad}{\biggerthick}{\sourcept}{\xdir}{\mu_1(\xdirsymb)}.
\]
Also let $\bigrad = \frac{\bigbeta}{\biggerbeta} \biggerrad$ and $\bigthick = \frac{\bigbeta}{\biggerbeta} \biggerthick$, and define
\[
\bigpin \definedas \bpin[2]{\mu_1}{\bigrad}{\bigthick}{\sourcept}{\xdir}{\mu_1(\xdirsymb)}.
\]
Since $\conept \in \cstarset{\normpair}{\biggerbeta}{2}{\thisapex}{\sourcept}$, Part 1 of Lemma~\ref{star_closer_pin_lem} implies that $\thisapex \notin \interior{\biggerpin}$. Since $\bigrad < \biggerrad$ and $\bigthick <\biggerthick$, we have $\bigpin \subset \interior{\biggerpin}$, so
\begin{equation}
\label{cross_section_covering:apex_exterior_eqn}
\thisapex \notin \bigpin \quad \forall\conept\in \bluntcone.
\end{equation}
Now let $\thisrad = \rmeet{2}{\beta} \bigrad$ and $\thisthick = (\bigbeta - \beta) \rmeet{2}{\beta}$, and define
\[
\thispin \definedas \bpin[2]{\mu_2}{\thisrad}{\thisthick}{\sourcept}{\xdir}{\mu_2(\xdirsymb)}.
\]
Then each of the $\mu_2$-bowling pins $\thispin$ is nondegenerate, and $\conept\in \interior{\thispin}$. Thus, the collection $\setof[2]{\interior{\thispin}}_{\conept\in \crosssection}$ is an open cover of the compact set $\crosssection$, so there is some finite subcover $\setof[2]{\interior{\thispin[1]},\dotsc,\interior{\thispin[k]}}$; let $\pincollection = \setof[2]{\thispin[1],\dotsc,\thispin[k]}$.
By \eqref{cross_section_covering:apex_exterior_eqn} we have $\thisapex\notin \bigcup_{j=1}^k \bigpin[j]$. Therefore, since $\bigcup_{j=1}^k \bigpin[j]$ is closed, there is some open neighborhood $\nbhd$ of $\thisapex$ such that
\[
\nbhd \cap \bigpin[j] = \emptyset
\quad \forall j\in\setof{1,\dotsc,k}.
\]
Since $\conetip[2]{\mu_1}{\bigthick}{\sourcept}{\xdir}{\mu_1(\xdirsymb)} \subseteq \bigpin$ for all $\conept$, Lemma~\ref{star_closer_char_lem} now implies that $\conept[j]\in \cstarset[2]{\normpair}{\bigbeta}{2}{\nbhd}{\sourcept[j]}$ for all $j$, and then Part 2 of Lemma~\ref{star_closer_pin_lem} implies that
\begin{equation}
\label{cross_section_covering:closer_eqn}
\thispin[j] \subseteq \cstarset[2]{\normpair}{\beta}{2}{\nbhd}{\sourcept[j]}
\quad \forall j\in\setof{1,\dotsc,k}.
\end{equation}
Thus, \eqref{cross_section_covering:closer_eqn} proves Part~\ref{cross_section_covering:closer_part} of Lemma~\ref{cross_section_covering_lem}, and Parts~\ref{cross_section_covering:origin_part} and \ref{cross_section_covering:cover_part} follow by construction, since $\sourcept[j]\in \conecomp$, and the collection $\interiorcollection = \setof[2]{\interior{\thispin[1]},\dotsc,\interior{\thispin[k]}}$ covers $\crosssection$.
For Part~\ref{cross_section_covering:thickness_part}, simply note that the thickness of each $\thispin[j]$ is $\thisthick = (\bigbeta - \beta) \rmeet{2}{\beta} = \frac{1}{3} \parens[2]{\coneadvantage{1}{\normpair}{\cones}^{-1}-\beta} \rmeet{2}{\beta}$.
\end{proof}

\begin{prop}[Extinction in small narrow cones]
\label{narrow_extinction_prop}
Let $\cones\subseteq\Rd$ be a closed cone that is narrow for species 1 and is contained in some half-space. If $\initialset{1}$ is any bounded subset of $\cones$, then $\conqset{1}{\normpair}{\initialset{1}}{\cones}$ is also bounded.
\end{prop}

\begin{proof}
\newcommand{\thisapex}{\vect 0}
\newcommand{\crosssection}{K}
\newcommand{\pincollection}{\boldsymbol{P}}
\newcommand{\interiorcollection}{\interior{\pincollection}}

\newcommand{\pinsymb}{P} 

\newcommand{\pin}{\pinsymb}
\newcommand{\pininterior}{\interior{P}}
\newcommand{\pinorigin}{\vect y_{\pin}}

\newcommand{\nbhd}{U_{\thisapex}}


\newcommand{\conecomp}{\Rd\setminus \cones}
\newcommand{\bluntcone}{\cones\setminus \setof{\thisapex}}
\newcommand{\conept}{\vect x}

\newcommand{\smallrad}{\alpha}
\newcommand{\bigrad}{R}

\newcommand{\smallball}{\smallrad \unitball{\mu_1}} 
\newcommand{\bigball}{\bigrad\smallrad \unitball{\mu_1}}
\newcommand{\biggerball}{\bigrad \unitball{\mu_1}}


By translation invariance we can assume that $\cones$ has apex $\vect 0$. Since $\cones$ is narrow for species 1, we can choose some $\beta$ with $1<\beta < \coneadvantage{1}{\normpair}{\cones}^{-1}$. Let $\crosssection = \sphere{d-1}{\mu_1}\cap \cones$, and let $\pincollection = \pincollection(\crosssection,\beta)$ be the family of $\mu_2$-bowling pins from Lemma~\ref{cross_section_covering_lem}. Then $\pincollection$ covers $\crosssection$, and there is some open neighborhood $\nbhd$ of $\thisapex$ such that 
\begin{equation}
\label{narrow_extinction:pin_closer_eqn}
\forall \pin\in \pincollection, \quad
\pin \subseteq \cstarset{\normpair}{\beta}{2}{\nbhd}{\pinorigin},
\end{equation}
where $\pinorigin \in \conecomp$ is the origin of $\pin$. Let $\smallrad>0$ be small enough that $\smallball \subset \nbhd$. We claim that $\smallrad<1$. To see this, first note that \eqref{narrow_extinction:pin_closer_eqn} implies that $\pin\cap \nbhd = \emptyset$ for all $\pin\in \pincollection$ since $\closerset{\normpair}{\beta}{2}{\nbhd}{\pinorigin}\subseteq \Rd\setminus \nbhd$ by definition. Since $\pincollection$ covers $\crosssection$, this implies that $\crosssection\cap \nbhd = \emptyset$. Now, since $\crosssection \subseteq \sphere{d-1}{\mu_1} \subset \unitball{\mu_1}$, if we had $\smallrad\ge 1$, then we would get $\crosssection \subset  \smallball \subset \nbhd$, which is a contradiction, so we must have $\smallrad<1$.

Now let $\initialset{1}\subset \cones$ be bounded, and define
\[
\bigrad \definedas \sup_{\vect z\in \initialset{1}} \frac{\mu_1(\vect z)}{\smallrad}<\infty.
\]
Then for all $\vect z\in \initialset{1}$ we have $\mu_1(\vect z) \le \bigrad\smallrad$, so
\begin{equation}
\label{narrow_extinction:A1_containment_eqn}
\initialset{1} \subseteq \bigball \subset \bigrad \nbhd.
\end{equation}
Thus, scaling the sets in \eqref{narrow_extinction:pin_closer_eqn} by $\bigrad$ and using \eqref{narrow_extinction:A1_containment_eqn}, we have
\begin{equation}
\label{narrow_extinction:scaled_pin_eqn}
\forall \pin\in \pincollection, \quad
\bigrad \pin
\subseteq \cstarset[2]{\normpair}{\beta}{2}{\bigrad \nbhd}{\bigrad \pinorigin}
\subseteq \cstarset[2]{\normpair}{\beta}{2}{\initialset{1}}{\bigrad \pinorigin}.
\end{equation}
Now, since $\conecomp$ is scale-invariant and $\pinorigin\in \conecomp$, we have $\bigrad \pinorigin\in \conecomp$ for all $\pin\in \pincollection$, and since $\pincollection$ covers $\crosssection$, the collection $\bigrad \pincollection \definedas \setof[2]{\bigrad \pin}_{\pin\in \pincollection}$ covers $\bigrad \crosssection$. Thus, since $\beta>1$, \eqref{narrow_extinction:scaled_pin_eqn} and Proposition~\ref{determ_star_conquering_prop} imply that species 1 does not conquer any of $\bigcup_{\pin\in \pincollection} \bigrad \pin \supseteq \bigrad \crosssection$ from the initial configuration $\pair{\initialset{1}}{\conecomp}$. Thus we have $\dfinalset{1}\cap \bigrad \crosssection = \emptyset$, where $\dfinalset{1} = \conqset[2]{1}{\normpair}{\initialset{1}}{\cones}$.

Finally, observe that since $\smallrad<1$, by \eqref{narrow_extinction:A1_containment_eqn} we have $\initialset{1}\subseteq \bigball \subset \biggerball$, and hence every path in $\cones$ from $\initialset{1}$ to $\cones\setminus \biggerball$ must pass through $\cones\cap \bd (\biggerball) = \bigrad \crosssection$. Since $\dfinalset{1}\cap \bigrad \crosssection = \emptyset$, this implies that there is no $\dfinalset{1}$-path from $\initialset{1}$ to $\cones\setminus \biggerball$. Therefore, species 1 cannot conquer any point in $\cones\setminus \biggerball$, because every point in $\dfinalset{1}$ must be connected to $\initialset{1}$ via a $\dfinalset{1}$-path by Part 3 of Lemma~\ref{determ_final_sets_properties_lem}.
%
%
Thus we have
$
\conqset[2]{1}{\normpair}{\initialset{1}}{\cones} \subseteq \biggerball,
$
which proves Proposition~\ref{narrow_extinction_prop}.
\end{proof}

\subsection{Additional Results for Competition in Cones}
\label{deterministic:cones:additional_sec}


Note that it follows immediately from definition \eqref{conq_set_def_eqn} and the monotonicity of the deterministic process (Lemma~\ref{determ_monotonicity_lem}) that for $B,B'\subseteq\Rd$,
\begin{equation}
\label{conq_sets_monotone_eqn}
B\subseteq B' \implies
\conqset{1}{\normpair}{\vect z}{B} \subseteq \conqset{1}{\normpair}{\vect z}{B'}.
\end{equation}
This fact implies the following result.

\begin{lem}[Conquered sets and convex subcones]
\label{conqset_relation_lem}
\newcommand{\thisset}{B}
\newcommand{\thispt}{\vect z}
\newcommand{\thisconfig}{\pair[1]{\thispt}{\,\Rd\setminus\thisset}}

Let $\thisset\subseteq\Rd$, let $\thispt\in\thisset$, and consider a $\normpair$-process started from $\thisconfig$. If $\cones$ is a convex pointed cone
such that $\cones\subseteq\thisset$ and $\thispt\in\cones$, then
\[
\conqset{1}{\normpair}{\thispt}{\thisset} \supseteq
\thispt + \conqset{1}{\normpair}{\vect 0}{\siversion{\cones}},
\]
where $\siversion{\cones}$ is the scale-invariant version of $\cones$.
%
\end{lem}

\begin{proof}
\newcommand{\thisset}{B}
\newcommand{\thispt}{\vect z}
\newcommand{\thisconfig}{\pair[1]{\thispt}{\,\Rd\setminus\thisset}}
\newcommand{\thissicone}{\siversion{\cones}}

Since $\thisset \supseteq\cones$, we have $\conqset{1}{\normpair}{\thispt}{\thisset} \supseteq \conqset{1}{\normpair}{\thispt}{\cones}$ by \eqref{conq_sets_monotone_eqn}.
Now, since $\thispt\in\cones$ and $\cones$ is convex and affine scale-invariant, Part~\ref{wedge_properties:convexity_part} of Lemma~\ref{wedge_properties_lem} implies that $\cones\supseteq \thispt+\thissicone$, and hence $\conqset{1}{\normpair}{\thispt}{\cones}\supseteq \conqset{1}{\normpair}{\thispt}{\thispt+\thissicone}$ by \eqref{conq_sets_monotone_eqn}. 
Finally, $\conqset{1}{\normpair}{\thispt}{\thispt+\thissicone} = \thispt + \conqset{1}{\normpair}{\vect 0}{\thissicone}$ by the translation-invariance of $\normpair$. Note that since $\cones$ is pointed by assumption, $\pair{\vect 0}{\thissicone}$ is a valid starting configuration and hence this conquered region is nonempty.
\end{proof}

Using Lemma~\ref{conqset_relation_lem}, we derive the next result as a special case of Part~\ref{cone_closer_set:infinite_part} of Lemma~\ref{cone_closer_set_lem}; this corresponds to
Proposition~\ref{wide_crit_survival_prop} in the special case where $\cones$ is a $\mu_2$-cone.

\begin{lem}[Conquering a $\mu_1$-subcone of a wide $\mu_2$-cone]
\label{determ_conq_mu1_subcone_lem}

\newcommand{\thisadv}{\alpha}
\newcommand{\bigthp}{\thpadv{2}{\dirvec}{\thisadv}}
\newcommand{\smallthp}{\symbolref{\thp}{determ_conq_mu1_subcone_lem} (\thisadv)}

\newcommand{\thisapex}{\vect a}
\newcommand{\conept}{\vect z}

\newcommand{\bigcone}{\cone{\mu_2}{\bigthp^-}{\thisapex}{\dirvec}}
\newcommand{\smallcone}{\cone{\mu_1}{\smallthp}{\conept}{\dirvec}}

Suppose $\esupnorm[2]{\frac{\mu_2}{\mu_1}}\ge \thisadv\ge 1$, and for $\dirsymb\in \Rd\setminus \setof{\vect 0}$, let $\thpadv{2}{\dirvec}{\thisadv} = \thisadv \frac{\mu_1(\dirsymb)}{\mu_2(\dirsymb)}$ as defined in \eqref{thpadv_def_eqn}. Then
\begin{enumerate}

\item \label{determ_conq_mu1_subcone:cone_part}
If $\frac{\mu_2(\dirsymb)}{\mu_1(\dirsymb)} \ge\thisadv$ and $\thisapex\in\Rd$, the set $\cones\definedas \bigcone$ is a small cone with $\coneadvantage{1}{\normpair}{\cones}\ge \thisadv$.

\item \label{determ_conq_mu1_subcone:conq_part}
There exists $\smallthp\ge 0$ such that if $\dirvec$ and $\cones$ are as in Part~\ref{determ_conq_mu1_subcone:cone_part}, then for any $\conept\in \cones$,
\[
\conqset[1]{1}{\normpair}{\conept}{\cones} \supseteq \smallcone,
\]
and $\smallthp>0$ if $\thisadv>1$.
\end{enumerate}
\end{lem}

\begin{proof}

\newcommand{\thisadv}{\alpha}
\newcommand{\bigthp}{\thpadv{2}{\dirvec}{\thisadv}}
\newcommand{\smallthp}{\symbolref{\thp}{determ_conq_mu1_subcone_lem} (\thisadv)}

\newcommand{\thisapex}{\vect a}
\newcommand{\conept}{\vect z}

\newcommand{\bigcone}{\cone{\mu_2}{\bigthp^-}{\thisapex}{\dirvec}}
\newcommand{\smallcone}{\cone{\mu_1}{\smallthp}{\conept}{\dirvec}}
\newcommand{\smallsicone}{\cone{\mu_1}{\smallthp}{\vect 0}{\dirvec}}

\proofpart{1}
If $\frac{\mu_2(\dirsymb)}{\mu_1(\dirsymb)} \ge\thisadv$, then $\bigthp \le \thisadv \cdot \thisadv^{-1} = 1$, so the set $\cones\definedas \bigcone$ is a small, convex, pointed cone by Lemma~\ref{mu_cone_description_lem}, and we have $\coneadvantage{1}{\normpair}{\cones}\ge \thisadv$ by Corollary~\ref{thp_adv_cor}.

\proofpart{2}
Let $\smallthp \definedas (\thisadv-1)\rmeet{1}{\thisadv}$, where $\rmeet{1}{\thisadv}>0$ is defined in \eqref{rmeet_def_eqn}. Then $\smallthp\ge0$, with $\smallthp>0$ if $\thisadv>1$, and by Part~\ref{cone_closer_set:infinite_part} of Lemma~\ref{cone_closer_set_lem} we have
\[
\smallsicone \subseteq
\cstarset[2]{\normpair}{1}{1}{\vect 0}{\Rd\setminus \siversion{\cones}}
\subseteq \conqset{1}{\normpair}{\vect 0}{\siversion{\cones}}.
\]
Then since $\cones$ is convex and pointed, Lemma~\ref{conqset_relation_lem} implies that for any $\conept\in\cones$,
\[
\conqset{1}{\normpair}{\conept}{\cones}
\supseteq \conept + \conqset{1}{\normpair}{\vect 0}{\siversion{\cones}}
\supseteq \conept + \smallsicone = \smallcone.
\qedhere
\]
\end{proof}

The next results follow from Proposition~\ref{conquering_from_point_prop} in the previous section.


\begin{prop}[Competition with norms that are scalar multiples of one another]
\label{mu_cone_criticality_prop}
\newcommand{\speedparam}{\lambda}
\newcommand{\thisnormpair}{\normpair_{\speedparam}}

\newcommand{\thisapex}{\vect a}
\newcommand{\thiscone}{\cone{\mu}{\thp}{\thisapex}{\dirvec}}
\newcommand{\thispt}{\vect z}
\newcommand{\thisconfig}{\pair[1]{\thispt}{\,\Rd\setminus \thiscone}}
\newcommand{\thisaxis}{\thisapex+\cdirvec}

\newcommand{\thisconqset}{\conqset[1]{1}{\thisnormpair}{\thispt}{\thiscone}}

\newcommand{\critspeed}{\speedparam_c(\cones)}
\newcommand{\critthp}{\thp_c}


Let $\mu$ be a norm on $\Rd$, and for each $\speedparam>0$ let $\thisnormpair = \pair[2]{\mu}{\,\speedparam^{-1}\mu}$. Consider a $\thisnormpair$-competition process started from $\thisconfig$, where $\thisapex\in\Rd$, $\dirvec\in\dirset{\Rd}$, $\thp\in (0,1]$, and $\thispt\in\thiscone$.
Then for any $\dirvec\in \dirset{\Rd}$,
\begin{enumerate}

\item $\displaystyle \coneadvantage[1]{1}{\thisnormpair}{\thiscone} = \frac{\thp}{\speedparam}$.

\item The cone $\thiscone$ is wide, critical, or narrow for species~1 according to whether $\thp>\speedparam$, $\thp=\speedparam$, or $\thp<\speedparam$, respectively.

\item For any $\thispt\in\thiscone$, species~1's conquered region, $\thisconqset$, contains a nondegenerate cone at $\thispt$ if $\thp>\speedparam$ and is bounded if $\thp<\speedparam$.

\item 
If $\thp=\speedparam$, then species~1 conquers at least the closed ray $\thispt+\cdirvec$. If additionally $\thispt=\thisapex$, then species~1 conquers precisely the closed ray $\thisaxis$. 



\end{enumerate}

\end{prop}

%
%
%
%
%
%
%

\begin{prop}[Existence of critical speed ratio in convex cones]
\label{mu_norms_critical_speed_prop}
\newcommand{\spone}{\lambda_1}
\newcommand{\sptwo}{\lambda_2}
\newcommand{\spratio}{\sptwo/\spone}

\newcommand{\newnormpair}{\normpair_{\spone,\sptwo}}

\newcommand{\thispt}{\vect z}
\newcommand{\thisconfig}{\pair[2]{\thispt}{\,\Rd\setminus\cones}}
\newcommand{\thisadv}{\coneadvantage{1}{\normpair}{\cones}}

\newcommand{\thisconqset}{\conqset{1}{\newnormpair}{\thispt}{\cones}}

\newcommand{\critratio}{\lambda_c(\cones)}
\newcommand{\bcritratio}{\Lambda_c(\cones)}

Let $\normpair = \pair{\mu_1}{\mu_2}$ be a pair of norms on $\Rd$.
For any $\spone,\sptwo>0$, consider the deterministic process using the pair of scaled norms $\newnormpair = \pair[2]{\spone^{-1}\mu_1}{\sptwo^{-1}\mu_2}$ and started from the configuration $\thisconfig$, where $\cones \subset \Rd$ is a convex cone, and $\vect z\in \cones$.
\begin{enumerate}
\item If $\sptwo/\spone<\thisadv$, then species~1 conquers a nondegenerate cone at $\thispt$.
\item If $\sptwo/\spone>\thisadv$, then species~1 conquers only a bounded region.
\end{enumerate}
Thus, for a given traversal norm pair $\normpair$, any nondegenerate convex cone $\cones$ has a critical speed ratio $\bcritratio \in \cosegment[1]{\esupnorm[2]{\frac{\mu_2}{\mu_1}}^{-1}}{\infty}$ such that in the $\newnormpair$-process started from $\thisconfig$, species~1 conquers a nondegenerate cone if $\lambda_1/\lambda_2 > \bcritratio$ and dies out if $\lambda_1/\lambda_2 < \bcritratio$.
%
\end{prop}

\part{Analyzing the Random Processes by Comparison with Deterministic Processes}


\chapter{Large Deviations Estimates for Growth in Cones}
\label{cone_growth_chap}

The aim of this chapter will be to obtain large deviations estimates for the growth of a restricted first-passage percolation process in $\Zd$, showing that for various restricting sets, the restricted process grows asymptotically at the same speed as an unrestricted process. It will be most convenient to treat the restricting sets as subsets of $\Rd$ rather than subgraphs of $\Zd$; to convert between the two, we use the lattice approximation and cube-expansion operations, as discussed in Section~\ref{fpc_basics:continuum_extension_sec}. We focus on restricting sets that are $\mu$-stars, i.e.\ unions of $\mu$-cone segments (cf.\ Sections~\ref{intro:ch4_sec} and \ref{deterministic:norm_geometry_sec}), because in these sets it is relatively easy to measure distances relative to the shape function $\mu$ of the random process, allowing for an easy comparison with the growth of the restricted deterministic $\mu$-process from Section~\ref{deterministic:one_type_sec}. In principle, one can use the estimates for growth in $\mu$-stars to obtain similar large deviations estimates for growth in other regions (for example, all convex bodies) by approximating these regions with sufficiently thin $\mu$-stars, though we do not explicitly state any such results. All the results in this chapter give good estimates only in ``large" restricting sets, getting exponentially better as the scale of the picture increases. For example, the results are useful for all sufficiently large times when the restricting set is scale-invariant (such as an infinite cone) or is a compact $\mu$-star with very large diameter.


The large deviations estimates we obtain can be seen as generalizing classical shape deviation estimates for the unrestricted process, such as Lemma~\ref{GM_large_dev_lem} below.
In fact, the large deviations estimate in Lemma~\ref{GM_large_dev_lem} for unrestricted growth will be the starting point for proving our main results about restricted growth.
Namely, we will bootstrap on the estimate in Lemma~\ref{GM_large_dev_lem}, extending it to increasingly general restricting sets --- first to $\mu$-balls (Lemma~\ref{ball_covering_lem}), then to $\mu$-cone segments (Lemma~\ref{thin_cone_covering_lem}), and finally to $\mu$-stars in Theorem~\ref{cone_segment_growth_thm}, which is the main general result of the chapter. The general strategy for proving each of these estimates (Lemmas~\ref{ball_covering_lem} and \ref{thin_cone_covering_lem} and Theorem~\ref{cone_segment_growth_thm}) will be to show that with high probability, the process first covers a ``permeating subset" of the restricting set within approximately the right amount of time, and then, with high probability, the process spreads from the permeating subset to the remainder of the restricting set in a relatively short amount of time. In Section~\ref{cone_growth:lattice_path_sec} below, we give a formal definition of ``$\epsilon$-permeating subset" (our definition coincides with the definition of ``$\epsilon$-net" from metric geometry), and we go into a bit more detail about how this concept fits into the structure of the later proofs.
%
This idea of employing ``permeating subsets" is essentially the same strategy that is used to prove the Shape Theorem (cf.\ \cite{Kesten:1986aa} or \cite{Howard:2004aa}); in the Shape Theorem, the set being covered is a $\mu$-ball, and the permeating subset is taken to be a collection of lattice points lying on a family of geometrically expanding spheres in the $\lone{d}$-metric. The main difference between the proofs of the estimates in this chapter and the proof of the Shape Theorem is the tool used to show that the permeating subset is covered in approximately the right amount of time. Namely, in the results below, the primary tool will be Lemma~\ref{GM_large_dev_lem} (a large deviations estimate for the Shape Theorem), while in the proof of the Shape Theorem itself, the primary tool is the subadditive ergodic theorem.


The organization of Chapter~\ref{cone_growth_chap} is as follows. In Section~\ref{cone_growth:cov_time_sec} we introduce the idea of ``covering times" and prove some elementary results about the restricted first-passage percolation process in $\Zd$, treated as a growth process in subsets of $\Rd$. In Section~\ref{cone_growth:lattice_path_sec}, we briefly discuss permeating subsets, and we prove a basic large deviations estimate for traversal of lattice paths that will be the basis of the large deviations estimates for ``short paths" that appear in all the main results of the chapter. In Section~\ref{cone_growth:ball_sec} we combine the main result of Section~\ref{cone_growth:lattice_path_sec} with the basic large deviations estimate in Lemma~\ref{GM_large_dev_lem} for unrestricted growth to obtain a large deviations estimate for growth restricted to $\mu$-balls, in Lemma~\ref{ball_covering_lem}. In Section~\ref{cone_growth:thin_cone_sec} we use Lemma~\ref{ball_covering_lem} to prove a large deviations estimate for covering a $\mu$-cone segment, in Lemma~\ref{thin_cone_covering_lem}; this will be the main technical result needed for the general results in the final two sections. In Section~\ref{cone_growth:cone_growth_sec} we use Lemma~\ref{thin_cone_covering_lem} to prove the main result of the chapter, Theorem~\ref{cone_segment_growth_thm}, which is a large deviations estimate for a first-passage percolation process restricted to any $\mu$-star. We then use this result to derive several corollaries, including a shape theorem in $\mu$-cones. Finally, in Section~\ref{cone_growth:logarithmic_tube_sec} we use Theorem~\ref{cone_segment_growth_thm} to prove a shape theorem for growth restricted to a $\mu$-tube whose width grows faster than logarithmically with its height.

\begin{remark}[Notational conventions for constants and events]
\label{constants_notation_rem}
Throughout Chapters~\ref{cone_growth_chap} and \ref{random_fpc_chap}, numerous constants will appear in the statements and proofs of the various results. I use the following conventions to help the reader keep the notation straight.
\begin{itemize}

\item Usually, time will be denoted by $t$ or $s$, and distances will be denoted by $r$ (for radius) or $h$ (for height). The thickness of cones will be denoted by $\thp$.

\item I will usually use $\spp$ for a small positive constant measuring  the ``relative speed error" in various events, similar to the $\spp$ appearing in the statement of the Shape Theorem. The variables $\alpha$ or $\beta$ may also play a similar role, corresponding roughly to $1-\spp$ or $1+\spp$, respectively (see below).

\item In large deviations estimates: $C =$ large positive constant, usually depending on $\spp$; $c=$ small positive constant, usually depending on $\spp$; and $K =$ (typically large) universal constant, not depending on $\spp$. In general, the constants $\bigcsymb$, $\smallcsymb$, and $\konstsymb$ also depend on the dimension $d$ and the distribution of the traversal measure $\tmeasure$, and sometimes on other parameters that will be identified in particular results. These constants will usually have subscripts identifying which lemma they are from.

\item Events: $\eventsymb_j(a;b)$ denotes an event in Lemma~$j$, whose probability we bound explicitly as a function of the parameter(s) $a$ and non-explicitly as a function of the parameter(s) $b$.

\item The variables $\alpha,\beta,\delta,\epsilon$ will be used to denote nonnegative constants. Usually, $\delta$ and $\epsilon$ should be thought of as positive constants that are close to 0, while $\alpha$ and $\beta$ should be thought of as constants that are close to 1, typically with $\alpha<1$ and $\beta>1$. In Chapters~\ref{random_fpc_chap} and \ref{coex_finite_chap}, $\alpha$ will also be used to denote the advantage of species~1, as defined in Section~\ref{deterministic:cone_competition_sec}.

%
%
%
%

\end{itemize}
\end{remark}

\section{Covering Times and Properties of the Continuum Process}
\label{cone_growth:cov_time_sec}


Let $\tmeasure$ be a (one-type) traversal measure on $\edges{\Zd}$, and let $\ptimefamily = \ptimefamily[\tmeasure]$ be the family of pseudometrics induced by $\tmeasure$. For $U,V,S\subseteq \Zd$, we define the \textdef{$S$-restricted covering time of $V$ from $U$} as
\begin{equation}\label{covering_time_def_eqn}
\ctime{}{S}{U}{V} \definedas \sup_{\vect v\in V} \rptime{}{S}{U}{\vect v}.
\end{equation}
To explain the terminology, observe that
$
\ctime{}{S}{U}{V} = 
\inf \setbuilder[1]{t}{\rreachedset{\tmeasure}{U}{S}{t} \supseteq V}.
$
That is, $\ctime{}{S}{U}{V}$ is the time at which the set $V$ is ``covered" by the $S$-restricted process $\rreachedfcn{\tmeasure}{U}{S}$ started from $U$.

We define covering times for continuum sets by taking lattice approximations, the same way we did for passage times: If $A,B,S\subseteq \Rd$, we define the \textdef{$S$-restricted covering time of $B$ from $A$} as
\begin{equation}\label{continuum_covering_time_def}
\ctime{}{S}{A}{B} \definedas \ctime[2]{}{\lat{S}}{\lat{A}}{\lat B}
= \sup_{\vect v\in \lat{B}} \rptime{}{S}{A}{\vect v}.
\end{equation}
If $S=B$, we call $\ctime{}{B}{A}{B}$ the \textdef{internal covering time} of $B$ from $A$, and we use the notation
\begin{equation}\label{int_cov_time_def_eqn}
\ictime{}{A}{B} \definedas \ctime{}{B}{A}{B},
\end{equation}
which will make some formulas more readable later on, particularly when $B$ is some subset of $\Rd$ with unwieldy notation. It follows from the definition \eqref{continuum_covering_time_def} that
\[
\ctime{}{S}{A}{B} = \sup_{\vect y\in B} \ctime{}{S}{A}{\vect y}
\quad\text{for any $B,S\subseteq \Rd$.}
\]
However, note that if $B$ is not a lattice set, then in general $\ctime{}{S}{A}{B} \ne \sup_{\vect y\in B} \rptime{}{S}{A}{\vect y}$, because if $\vect y \in \bd \cubify{\vect v}$ for some $\vect v\in\Zd$, then $\abs{\lat{\vect y}}>1$, and so typically $\ctime{}{S}{A}{\vect y} \ne \rptime{}{S}{A}{\vect y}$. On the other hand, note that if $\vect y\in \interior{\cubify{\vect v}}$ for some $\vect v\in\Zd$, then $\abs{\lat{\vect y}}=1$, and so  $\ctime{}{S}{A}{\vect y} = \rptime{}{S}{A}{\vect y}$.
We now enumerate some
elementary properties of covering times.
The proof is an easy exercise.

\begin{lem}[Properties of covering times]
\label{covering_time_properties_lem}
For any $A,B,C,S\subseteq \Rd$, the following properties hold.
\begin{enumerate}
\item  \label{covering_time_properties:monotone_S_part}
If $S'\subseteq S$, then $\ctime{}{S'}{A}{B} \ge \ctime{}{S}{A}{B}$ (monotonicity in $S$).

\item  \label{covering_time_properties:monotone_AB_part}
$\displaystyle \ctime{}{S}{A}{B} = \inf_{A'\subseteq A} \ctime{}{S}{A'}{B} = \sup_{B'\subseteq B} \ctime{}{S}{A}{B'}$ (monotonicity in $A$ and $B$).

\item  \label{covering_time_properties:tri_ineq_part}
$\ctime{}{S}{A}{C} \le \ctime{}{S}{A}{B} + \ctime{}{S}{B}{C}$ (triangle inequality).
%

\item \label{covering_time_properties:decomp_part}
If $B= \bigcup_{i\in I} B_i$, then $\displaystyle \ctime{}{S}{A}{B} = \sup_{i\in I} \ctime{}{S}{A}{B_i}$ (decomposition).

\item  \label{covering_time_properties:superdecomp_part}
If $B\subseteq \bigcup_{i\in I} B_i$, and if $A_i\subseteq A$ and $S_i\subseteq S$ for all $i\in I$, then
\[
\ctime{}{S}{A}{B} \le \sup_{i\in I} \ctime{}{S_i}{A_i}{B_i}
\quad\text{(superdecomposition).}
\]

\item  \label{covering_time_properties:tmeasure_bound_part}
If $A\subseteq \Rd$, and $B$ is a connected subset of $\Rd$ or $\Zd$ with $\lat B\subseteq \lat S$, then
\[
\ctime{}{S}{A}{B} \le \rptime{}{S}{A}{B} + \tmeasure \parens{B}
\quad\text{(total traversal measure bound),}
\]
where $\tmeasure(B) \definedas \tmeasure \parens[2]{\edges{\lat B}}$.
\end{enumerate}
\end{lem}

%

Note that, unlike the passage time $\rptime{}{S}{A}{B}$, the covering time $\ctime{}{S}{A}{B}$ is asymmetric with respect to $A$ and $B$. However, the usefulness of covering times comes from the fact that they satisfy the triangle inequality with arbitrary \emph{subsets} of $\Rd$ as arguments, rather than just \emph{points} in $\Zd$, whereas this property fails for passage times. We now prove a simple extension of the triangle inequality (as well as the monotonicity properties and finite decomposition property above)
that allows us to bound the covering time of a set $B$ by
taking a finite number of intermediate steps. This will be our main tool for analyzing the restricted process using covering times.
\begin{lem}[Chaining of covering times]
\label{chaining_lem}
Let $\setof{A_k}_{k=1}^n$, $\setof{B_k}_{k=0}^n$, and $\setof{S_k}_{k=1}^n$ be collections of subsets of $\Rd$, and suppose $A_k\subseteq \bigcup_{j=0}^{k-1} B_{j}$ for $1\le k\le n$. Then for any $B,S\subseteq\Rd$ with $B \subseteq \bigcup_{k=0}^n B_k$ and $\bigcup_{k=1}^n S_k\subseteq S$,
\[
\ctime{}{S}{B_0}{B} \le \sum_{k=1}^n \ctime{}{S_k}{A_k}{B_k}.
\]
\end{lem}

\begin{proof}
We prove this by induction on $n$. 
For $n=0$, we have $B\subseteq B_0$ and hence $\ctime{}{S}{B_0}{B}=0$, which equals the empty sum on the right, so the bound holds in this case. Now suppose $n\ge 1$, and for fixed sets $\setof{A_k}_{k=1}^n$, $\setof{B_k}_{k=0}^n$, $\setof{S_k}_{k=1}^n$, and $S \supseteq \bigcup_{k=1}^n S_k$, assume inductively that if $B'$ is any subset of $\Rd$ with $B'\subseteq  \bigcup_{k=0}^{n-1} B_k$, then
\[
\ctime{}{S}{B_0}{B'} \le \sum_{k=1}^{n-1} \ctime{}{S_k}{A_k}{B_k}.
\]
Set $B'= \bigcup_{k=0}^{n-1} B_k$. Then for any $B \subseteq \bigcup_{k=0}^n B_k$ we have $B\setminus B' \subseteq B_n$, and since $\ctime{}{S}{B'}{B} = \ctime{}{S}{B'}{B\setminus B'}$ by the decomposition property, we get
\begin{align*}
\ctime{}{S}{B_0}{B}
&\le \ctime{}{S}{B_0}{B'} + \ctime{}{S}{B'}{B}
	&& (\text{triangle inequality})\\
&\le \sum_{k=1}^{n-1} \ctime{}{S_k}{A_k}{B_k} + \ctime{}{S}{B'}{B\setminus B'}
	&& (\text{inductive hypothesis, decomposition})\\
&\le \sum_{k=1}^{n-1} \ctime{}{S_k}{A_k}{B_k} + \ctime{}{S_n}{A_n}{B_n}
	&& (\text{monotonicity})\\
&= \sum_{k=1}^{n} \ctime{}{S_k}{A_k}{B_k}.
	&& \qedhere
\end{align*}
%
%
\end{proof}

\begin{thmremark}
Note that Lemma~\ref{chaining_lem} also follows easily from the Markovesque property (Lemma~\ref{one_type_markovesque_lem}). However, the direct inductive proof of Lemma~\ref{chaining_lem} given above is much simpler than the rather subtle proof of Lemma~\ref{one_type_markovesque_lem}.
\end{thmremark}

Recall the definition of the one-type continuum process started from $A\subseteq\Rd$ and restricted to $S\subseteq\Rd$:
\[
\rreachedset{\tmeasure}{A}{S}{t} \definedas 
	\setbuilder[1]{\vect x\in S}{\rptime{\tmeasure}{S}{A}{\vect x}\le t},
\]
where $\rptime{\tmeasure}{S}{A}{\vect x} = \rptime[2]{\tmeasure}{\lat{S}}{\lat{A}}{\lat{\vect x}}$ by definition. We now define a similar continuum process using covering times instead of passage times, which we will call the continuum \textdef{covering process}:
\begin{equation}
\label{covering_process_def_eqn}
\crreachedset{\tmeasure}{A}{S}{t} \definedas 
	\setbuilder[1]{\vect x\in S}{\ctime{\tmeasure}{S}{A}{\vect x}\le t}.
\end{equation}
It is then natural to ask how the two processes $\reachedsymb$ and $\creachedsymb$ are related. Since for $\vect v\in\Zd$, the passage time $\rptime{\tmeasure}{S}{A}{\vect v}$ and the covering time $\ctime{\tmeasure}{S}{A}{\vect v}$ agree, both of the above processes can be viewed as continuum versions of the underlying lattice process on $\lat S$; that is,
\[
\rreachedset{\tmeasure}{A}{\lat S}{t} = \crreachedset{\tmeasure}{A}{\lat S}{t}
	= \setbuilder[1]{\vect v\in \lat S}{\rptime{\tmeasure}{S}{A}{\vect v}\le t}.
\]
In fact, $\rreachedset{\tmeasure}{A}{S}{t}$ is obtained by placing a closed unit cube at each vertex of $\rreachedset{\tmeasure}{A}{\lat S}{t}$, then taking the intersection with $S$, and $\crreachedset{\tmeasure}{A}{S}{t}$ is simply the relative interior of $\rreachedset{\tmeasure}{A}{S}{t}$ in $S$. That is, $\rreachedset{\tmeasure}{A}{S}{t} = S\cap \cubify[2]{\rreachedset{\tmeasure}{A}{\lat S}{t}}$, and $\crreachedset{\tmeasure}{A}{S}{t} = S\cap \interior{\cubify[2]{\rreachedset{\tmeasure}{A}{\lat S}{t}}}$.
The following proposition shows that both continuum versions of the process will be useful.

\begin{prop}[Hitting times and covering times for continuum processes]
\label{hitting_covering_prop}
Let $\tmeasure$ be a traversal measure on $\edges{\Zd}$, and let $A,B,S\subseteq \Rd$ with
$B\subseteq S$. Then
\begin{enumerate}
\item The hitting time of $B$ for the process $\rreachedset{\tmeasure}{A}{S}{t}$ and the hitting time of $\lat{B}$ for the process $\rreachedset{\tmeasure}{A}{\lat{S}}{t}$ are both equal to $\rptime{\tmeasure}{S}{A}{B}$; that is,
\[
\rptime{\tmeasure}{S}{A}{B}
= \inf \setbuilder[3]{t}{\rreachedset{\tmeasure}{A}{S}{t}\cap B \ne \emptyset}
= \inf \setbuilder[3]{t}{\rreachedset{\tmeasure}{A}{\lat S}{t}\cap \lat B \ne \emptyset}.
\]

\item The covering time of $B$ for the process $\crreachedset{\tmeasure}{A}{S}{t}$ and the covering time of $\lat{B}$ for the process $\rreachedset{\tmeasure}{A}{\lat{S}}{t}$ are both equal to $\ctime{\tmeasure}{S}{A}{B}$; that is,
\[
\ctime{\tmeasure}{S}{A}{B}
= \inf \setbuilder[3]{t}{\crreachedset{\tmeasure}{A}{S}{t}\supseteq B}
= \inf \setbuilder[3]{t}{\rreachedset{\tmeasure}{A}{\lat S}{t}\supseteq \lat B}.
\]
\end{enumerate}
\end{prop}

The proof of Proposition~\ref{hitting_covering_prop} is left as a trivial exercise in unravelling the definitions. In light of Proposition~\ref{hitting_covering_prop}, we refer to $\rreachedfcn{\tmeasure}{A}{S}$ as the \textdef{hitting process}, and $\crreachedfcn{\tmeasure}{A}{S}$ as the \textdef{covering process}, for the underlying lattice process $\rreachedset{\tmeasure}{A}{\lat S}{t}$. The hitting process is appropriate for proving upper bounds on the growth of the conquered set, while the covering process is appropriate for proving lower bounds on growth.


Here is another useful property of the continuum processes, which follows from Lemmas~\ref{restricted_proc_equiv_lem} and \ref{strong_restricted_proc_equiv_lem} in Section~\ref{fpc_basics:basic_properties_sec}.


\begin{lem}[Continuum process containment]
\label{continuum_containment_lem}
Let $\tmeasure$ be a traversal measure on $\edges{\Zd}$, and let $A,S\subseteq \Rd$ with $\lat A\subseteq \lat S$. If $\rreachedset{\tmeasure}{A}{\Rd}{t} \subseteq S$ for some $t\ge 0$, then
\[
\rreachedset{\tmeasure}{A}{S}{t'} = \rreachedset{\tmeasure}{A}{\Rd}{t'}
\quad\text{and}\quad
\crreachedset{\tmeasure}{A}{S}{t'} = \crreachedset{\tmeasure}{A}{\Rd}{t'}
\quad\text{for all } t'< t.
\]
If $\tmeasure$ satisfies \geodesicsexist, the above conclusions also hold for $t'=t$.
\end{lem}

\begin{proof}
After unwinding the definitions for the continuum processes, this follows directly from Lemmas~\ref{restricted_proc_equiv_lem} and \ref{strong_restricted_proc_equiv_lem} in Section~\ref{fpc_basics:basic_properties_sec}.
\end{proof}

The following formulas, based on a simple result in Appendix~\ref{asymptotic_chap}, will be useful for converting between statements about the process $\rreachedfcn{\tmeasure}{A}{S}$ (or $\crreachedfcn{\tmeasure}{A}{S}$) and statements about the underlying pseudometric $\rptmetric[\tmeasure]{S}$ (or ``quasipseudometric" $\ctmetric{\tmeasure}{S}$). I call these the ``inversion formulas" because they encode what Howard \cite{Howard:2004aa} refers to vaguely as an ``inversion argument" (but does not write down explicitly) in the proof of the Shape Theorem.

\begin{lem}[Inversion formulas for first-passage growth]
\label{inversion_formula_lem}
Let $A,A',S\subseteq \Rd$. 
\begin{enumerate}
\item For any $\alpha>0$ (in particular $\alpha<1$),
\begin{enumerate}
\item
\label{inversion_formulas:lbg} 
$\displaystyle
\eventthat[3]{\crreachedset{\tmeasure}{A}{S}{t}
	\supseteq \rreacheddset{\mu}{A'}{S}{\alpha t}
	\; \forall t\ge t_0}
= \eventthat[3]{\ctime{\tmeasure}{S}{A}{\vect x} \le t_0\vee \alpha^{-1} \idist{\mu}{S}{A'}{\vect x}
	\; \forall \vect x\in S}.
$

\item
\label{inversion_formulas:strong_lbg} 
$\displaystyle
\eventthat[3]{\crreachedset{\tmeasure}{A}{S}{t-}
	\supseteq \rreacheddset{\mu}{A'}{S}{\alpha t}
	\; \forall t\ge t_0}
= \eventthat[3]{\ctime{\tmeasure}{S}{A}{\vect x} < t_0\vee \alpha^{-1} \idist{\mu}{S}{A'}{\vect x}
	\; \forall \vect x\in S}.
$
\end{enumerate}

\item For any $\beta>0$ (in particular $\beta>1$),
\begin{enumerate}
\item
\label{inversion_formulas:ubg} 
$\displaystyle
\eventthat[3]{\rreachedset{\tmeasure}{A}{S}{t}
	\subseteq \rreacheddset{\mu}{A'}{S}{\beta t}
	\; \forall t\ge t_0}
= \eventthat[3]{t_0\vee \rptime{\tmeasure}{S}{A}{\vect x} \ge \beta^{-1} \idist{\mu}{S}{A'}{\vect x}
	\; \forall \vect x\in S}.
$

\item
\label{inversion_formulas:strong_ubg} 
$\displaystyle
\eventthat[3]{\rreachedset{\tmeasure}{A}{S}{t}
	\subseteq \rreacheddset{\mu}{A'}{S}{\beta t-}
	\; \forall t\ge t_0}
= \eventthat[3]{t_0\vee \rptime{\tmeasure}{S}{A}{\vect x} > \beta^{-1} \idist{\mu}{S}{A'}{\vect x}
	\; \forall \vect x\in S}.
$
\end{enumerate}
\end{enumerate}
\end{lem}

\begin{proof}
The formulas are special cases of the general formulas in Lemma~\ref{general_inversion_lem}; take the functions $f$ and $g$ to be $\ctime{\tmeasure}{S}{A}{\cdot}$, $\rptime{\tmeasure}{S}{A}{\cdot}$, or $\idist{\mu}{S}{A'}{\cdot}$, as appropriate.
\end{proof}

\section{Permeating Subsets and Large Deviations for Traversal of Lattice Paths}
\label{cone_growth:lattice_path_sec}

Here we introduce a simple concept that will be useful for understanding the structure of the proofs in this chapter. Let $(M,\dist{})$ be a metric space, and let $\epsilon\ge 0$. We say that $A$ is an \textdef{$\epsilon$-permeating subset} of $M$ (with respect to $\dist{}$) if $\dist{}(A,\vect y)\le \epsilon$ for every $\vect y\in M$. This definition of an $\epsilon$-permeating subset coincides precisely with the definition of an \textdef{$\epsilon$-net} of the metric space $(M,\dist{})$ (cf.\ \cite[p.~13]{Burago:2001aa}), but I opt for the term ``$\epsilon$-permeating" rather than ``$\epsilon$-net" because I feel it is more descriptive and avoids any potential confusion with the unrelated concept of nets in general topology.

In this chapter we are interested in showing that a restricted process doesn't take too long to cover some set $B\subseteq\Rd $. In this setting, we take $M$ to be the set $B$, typically endowed with the  ``natural distance" $\dist{\mu}$ associated with a first-passage growth process with shape function $\mu$. The general strategy for proving that the process covers $B$ within, say, time $(1+\epsilon)t$, is to break the the covering event into two steps: First, the process takes over a $\delta t$-permeating subset of $B$ within time $(1+\epsilon_1)t$ for some small $\delta>0$ and $\epsilon_1<\epsilon$; then, it spreads from the $\delta t$-permeating subset to cover the rest of $B$ in a short amount of time, say $\epsilon_2 t$, where $\epsilon_1+\epsilon_2 \le \epsilon$.
In order to apply this ``two-stage covering" technique, we need large deviations estimates putting a small upper bound on the probability that a ``bad" event happens in each of the two covering stages.

In the first stage, in order to show that the process is likely to cover a permeating subset of $B$ in about the right amount of time, we need some result showing that with high probability, the speed of the process on large scales doesn't deviate too much from the asymptotic speed identified in the Shape Theorem; Lemma~\ref{GM_large_dev_lem} in the next section provides the necessary starting point, and we will bootstrap our way from there to obtain similar, increasingly more general results. In the second stage, we need to show that with high probability, the process doesn't take too long to reach the rest of the set $B$ once it covers a permeating subset; for this we use elementary large deviations estimates for sums of \iid\ exponentially-tailed random variables to obtain a crude bound on the time it takes the process to travel a short distance. Lemma~\ref{tmeas_dev_lem} below provides the basic large deviations estimate in the appropriate context; we use Lemma~\ref{tmeas_dev_lem} to prove Lemma~\ref{segment_meas_lem}, which will be our primary tool for bounding the traversal times of ``short" paths in later proofs.

\begin{lem}[Upper large deviations bound for traversal measure]
\label{tmeas_dev_lem}

\newcommand{\thisrv}{\tmeasureof{\edge}}
\newcommand{\thislaw}{\law \parens[2]{\thisrv}}

\newcommand{\thisrvj}[1]{\tmeasureof{\edge_{#1}}} 

\newcommand{\thisC}{\bigcref{tmeas_dev_lem}(a)}
\newcommand{\thisc}{\smallcref{tmeas_dev_lem}(a)}

Let $\tmeasure$ be an \iid\ traversal measure on $\Zd$ such that $\thisrv$ satisfies \expmoment. For any $a>\E \thisrv$, there exist positive constants $\thisC$ and $\thisc$ (which depend also on $\thislaw$) such that for any finite set of edges $F\subseteq \edges{\Zd}$ with $\abs{F}\le n\in \N$,
\[
\Pr\eventthat[2]{ \tmeasure(F) \ge an}
\le \thisC e^{-\thisc n}.
\]
\end{lem}

\begin{proof}
\newcommand{\thisrv}{\tmeasureof{\edge}}
\newcommand{\thislaw}{\law \parens[2]{\thisrv}}

\newcommand{\thisrvj}[1]{\tmeasureof{\edge_{#1}}} 

\newcommand{\thisC}{\bigcref{tmeas_dev_lem}(a)}
\newcommand{\thisc}{\smallcref{tmeas_dev_lem}(a)}

By definition we have $\tmeasure(F)= \sum_{\edge\in F} \tmeasure(\edge)$. By assumption, the collection $\setof{\tmeasure(\edge)}_{\edge\in F}$ is \iid\ with $\tmeasure(\edge)\eqd \mttime_0$ for all $\edge$.  Thus, since $\abs{F}\le n$ we have
$
\tmeasure(F) \eqd \sum_{j=1}^{\abs{F}} \mttime_j
\le \sum_{j=1}^n \mttime_j
$,
and hence
\[
\Pr \eventthat[2]{\tmeasure(F)\ge an}  \le \Pr \eventthat[1]{\sum_{j=1}^{n} \mttime_j \ge an }
\le \bigcref{basic_large_dev_lem}(a) e^{-\smallcref{basic_large_dev_lem}(a) n},
\]
where the final inequality follows from Lemma~\ref{basic_large_dev_lem} since $\thisrv$ has a finite exponential moment and $a>\E\thisrv$.
\end{proof}

\begin{lem}[Traversal measure of lattice paths near a line segment]
\label{segment_meas_lem}
\newcommand{\thisnorm}{\mu}
\newcommand{\thisrv}{\tmeasureof{\edge}}
\newcommand{\thislaw}{\law \parens[2]{\thisrv}}

\newcommand{\thisC}{\bigcref{segment_meas_lem}}
\newcommand{\thisCp}{\thisC  \parens[2]{\phi, \thisnorm, \thislaw}}

\newcommand{\thisc}{\smallcref{segment_meas_lem}}
\newcommand{\thiscp}{\thisc \parens[2]{\phi, \thisnorm, \thislaw}}

\newcommand{\thiskonst}{\konstref{segment_meas_lem}^{\thisnorm}}
\newcommand{\thiskonstp}{\thiskonst \parens[2]{\E \thisrv}}

\newcommand{\thisevent}[1]{\eventref{segment_meas_lem} \parens{\vect u,\vect v, #1}}

Let $\tmeasure$ be an \iid\ traversal measure such that $\thisrv$ satisfies \expmoment, and for any $\vect u,\vect v\in\Zd$ and $t\ge 0$, define the event
\[
\thisevent{t} =  \eventthat[3]{\tmeasure \parens[2]{\segment{\vect x}{\vect y}} 
	\le t
	\text{\rm\ for all  $\vect x\in \cubify{\vect u}$ and $\vect y\in\cubify{\vect v}$} }.
\]
Then for any norm $\thisnorm$ and any function $\phi\colon \Rplus\to\R$ with $\phi(r)/r\to 0$ as $r\to\infty$, there exists some positive constant $\thiskonst = \thiskonstp$ and positive constants $\thisC = \thisCp$ and $\thisc = \thiscp$ such that for any $r\ge 0$ and any $\vect u, \vect v\in\Zd$ with
$
\idist{\thisnorm}{}{\vect u}{\vect v} \le r + \phi(r),
$
we have
\[
\Pr \parens[2]{\thisevent{\thiskonst r}}
	\ge 1- \thisC e^{-\thisc r}.
\]
\end{lem}

\begin{proof}
\newcommand{\thisnorm}{\mu}
\newcommand{\thisrv}{\tmeasureof{\edge}}
\newcommand{\thislaw}{\law \parens[2]{\thisrv}}

\newcommand{\thisC}{\bigcref{segment_meas_lem}}
\newcommand{\thisCp}{\thisC  \parens[2]{\phi, \thisnorm}}

\newcommand{\thisc}{\smallcref{segment_meas_lem}}
\newcommand{\thiscp}{\thisc \parens[2]{\phi, \thisnorm}}

\newcommand{\thiskonst}{\konstref{segment_meas_lem}^{\thisnorm}}

\newcommand{\thisevent}[1]{\eventref{segment_meas_lem} \parens{\vect u,\vect v, #1}}

\newcommand{\genkonst}{\konstsymb}

\newcommand{\thisx}{\vect x}
\newcommand{\thisy}{\vect y}
\newcommand{\xysegment}{\segment{\thisx}{\thisy}}

\newcommand{\thisu}{\vect u}
\newcommand{\thisv}{\vect v}
\newcommand{\uvsegment}{\segment{\thisu}{\thisv}}

\newcommand{\ucube}{\cubify{\thisu}}
\newcommand{\vcube}{\cubify{\thisv}}

\newcommand{\thisconvhull}{\conv \parens[2]{\ucube\cup \vcube}}

\newcommand{\chull}{{\rm co}_{\thisu,\thisv}}
\newcommand{\latchull}{\lat[2]{\chull}}
\newcommand{\edgeset}{F_{\thisu,\thisv}}

\newcommand{\thiscylconst}{\normcylconst{\linfty{d}}{\linfty{d}}}
\newcommand{\edgecardconst}{A}

\newcommand{\thisrad}{R}
\newcommand{\explicitrad}{\frac{3\sqrt{d}}{2}}
\newcommand{\thisheight}{\distnew{\ltwo{d}}{\thisu}{\thisv}}

\newcommand{\ballmeasure}[1]{B_{#1}}

\newcommand{\aconst}{a_d}
\newcommand{\bconst}{b_{d,\mu}}
\newcommand{\rconst}{r_0}



First, observe that the definition $\tmeasureof{A} \definedas \tmeasureof[1]{\edges[2]{\lat A}}$ for $A\subseteq\Rd$ preserves the monotonicity property of measures, i.e.\ if $A\subseteq B\subseteq \Rd$, then $\tmeasureof{A}\le \tmeasureof{B}$. Now note that the definition of convex hull implies that
\[
\bigcup_{\thisx\in \ucube, \thisy\in \vcube}
\segment{\thisx}{\thisy}
= \thisconvhull.
\]
In particular, $\segment{\thisx}{\thisy} \subseteq \thisconvhull$ for all $\thisx\in\ucube$ and $\thisy\in\vcube$, so
\[
\sup_{\thisx\in \ucube, \thisy\in \vcube}
\tmeasureof[2]{\xysegment}\le \tmeasureof[2]{\thisconvhull}.
\]
Therefore, for any $\genkonst>0$,
\begin{align}
\label{segment_meas:eventcomp_eqn}
\eventcomplement{\thisevent{\genkonst r}}
&= \eventthat[1]{ \tmeasureof[2]{\xysegment} > \genkonst r
	\text{ for some $\thisx\in\ucube$ and $\thisy\in\vcube$}} \notag\\
&\subseteq \eventthat[1]{\tmeasureof[2]{\thisconvhull} > \genkonst r}.
\end{align}
For the sake of more compact notation, let $\chull \definedas \thisconvhull$, and let $\edgeset \definedas \edges[1]{\latchull}$. We want to bound the cardinality of the edge set $\edgeset$ in order to choose $\genkonst$ large enough that Lemma~\ref{tmeas_dev_lem} can be applied to the event in \eqref{segment_meas:eventcomp_eqn}. First note that since $\Zd$ is $2d$-regular, we have $\card{\edges{V}}\le d\cdot \card{V}$ for any $V\subseteq \Zd$, so $\card{\edgeset}\le d\cdot \card[1]{\latchull}$.
By Lemma~\ref{lattice_cube_properties_lem}, we have $\card[1]{\latchull} = \lebmeasof[1]{d}{\cubify[2]{\latchull}}$, and $\cubify[2]{\latchull} \subseteq \chull + \unitball{\linfty{d}}$. Now Lemma~\ref{convex_hull_properties_lem} implies that $\chull = \uvsegment + \frac{1}{2}\unitball{\linfty{d}}$, and hence $\chull + \unitball{\linfty{d}} =  \uvsegment + \frac{3}{2}\unitball{\linfty{d}}$. Therefore, noting that $\dilnorm{\linfty{d}}{\ltwo{d}} = \sqrt d$, we get
\begin{equation}
\label{segment_meas:cyl_containment_eqn}
\cubify[2]{\latchull} \subseteq \uvsegment
+ \frac{3}{2} \sqrt{d} \cdot \unitball{\ltwo{d}}.
\end{equation}
The set on the right-hand side of \eqref{segment_meas:cyl_containment_eqn} is the union of a Euclidean cylinder of radius $\explicitrad $ and height $\thisheight$, with two solid hemispherical caps of radius $\explicitrad$ at its ends. Taking the Lebesgue measure of both sides of \eqref{segment_meas:cyl_containment_eqn}, we get
\begin{equation}
\label{segment_meas:latchull_card_eqn}
\card[1]{\latchull}= \lebmeasof[1]{d}{\cubify[2]{\latchull}} 
\le \ballmeasure{d} \cdot \parens[1]{\textstyle \explicitrad}^d
+\ballmeasure{d-1} \cdot \parens[1]{\textstyle \explicitrad}^{d-1} \cdot \thisheight,
\end{equation}
where $\ballmeasure{n}$ is the Lebesgue measure of $\unitball{\ltwo{n}}$, the Euclidean unit ball in $\R^n$. Setting $\aconst\definedas d \ballmeasure{d}  \parens[1]{ \explicitrad}^{d}$ and $\bconst\definedas d \ballmeasure{d-1}  \parens[1]{\textstyle \explicitrad}^{d-1} \dilnorm{\ltwo{d}}{\thisnorm}$, \eqref{segment_meas:latchull_card_eqn} implies that
\begin{equation}
\label{segment_meas:edge_card_eqn}
\card{\edgeset} \le d\cdot \card[2]{\latchull}
\le \aconst +\bconst \distnew{\thisnorm}{\thisu}{\thisv}
\le \aconst + \bconst \parens[2]{r+\phi(r)}.
\end{equation}
Now choose $\rconst \ge \aconst/\bconst$ large enough that $\phi(r)/r \le 1$ for all $r\ge \rconst$. Then for all $r\ge \rconst$ we have
\[
\aconst + \bconst \parens[2]{r+\phi(r)}
= \bconst r \parens[1]{\frac{\aconst}{\bconst r} + 1 + \frac{\phi(r)}{r}}
\le 3 \bconst r,
\]
so \eqref{segment_meas:edge_card_eqn} implies that
\begin{equation}
\label{segment_meas:edge_card2_eqn}
\forall r\ge \rconst,\quad
\card{\edgeset} \le 3 \bconst r.
\end{equation}
Now define $\thiskonst \definedas 4\bconst \E \thisrv$. Then for all $r\ge \rconst$ we have
\begin{align*}
\eventcomplement{\thisevent{\thiskonst r}}
&\subseteq \eventthat[3]{\tmeasureof[1]{\edges[1]{\latchull}} > \thiskonst r}
	&&\text{(by  \eqref{segment_meas:eventcomp_eqn})}\\
&\subseteq \eventthat[1]{\tmeasureof[2]{\edgeset}
	> \frac{\thiskonst}{3 \bconst} \cdot\floor{3 \bconst r}}
	&&\text{(by the definition of $\edgeset$)}\\
&= \eventthat[1]{\tmeasureof[2]{\edgeset}
	> \frac{4}{3} \E \thisrv \cdot\floor{3 \bconst r}}.
	&&\text{(by the definition of $\thiskonst$)}
\end{align*}
Therefore, by  \eqref{segment_meas:edge_card2_eqn} and Lemma \ref{tmeas_dev_lem}, for all $r\ge \rconst$ we have
\begin{align*}
\Pr \parens[1]{\eventcomplement{\thisevent{\thiskonst r}}}
&\le  \Pr \eventthat[1]{\tmeasureof[2]{\edgeset} > \frac{4}{3} \E \thisrv \cdot\floor{3 \bconst r}}\\
&\le \bigcref{tmeas_dev_lem} \parens[1]{ \tfrac{4}{3} \E \thisrv}
	e^{-\smallcref{tmeas_dev_lem} \parens[1]{ \frac{4}{3} \E \thisrv} \floor{3 \bconst r} }\\
&\le \thisC e^{-\thisc r},
\end{align*}
where, with $a=\tfrac{4}{3} \E \thisrv$,
\[
\thisc \definedas 3\bconst \smallcref{tmeas_dev_lem} \parens[1]{a}
\quad\text{and}\quad
\thisC \definedas
\bigcref{tmeas_dev_lem} \parens[1]{a}
e^{\smallcref{tmeas_dev_lem}  \parens[1]{a}}
\join e^{\thisc \rconst}. 
\]
Note that we took the maximum with $e^{\thisc \rconst}$ in the definition of $\thisC$ to get rid of the floor in the exponent. But then since $\thisC \ge e^{\thisc \rconst}$ by definition, it follows that if $r\le \rconst$, then $\thisC e^{-\thisc r}\ge 1$, so the probability bound trivially holds for $r\le \rconst$, hence for all $r\ge 0$.
\end{proof}

\section{Growth Restricted to a $\mu$-Ball}
\label{cone_growth:ball_sec}

%
%
%
%

Throughout the rest of the chapter, unless otherwise specified, $\tmeasure$ will be a one-type traversal measure on $\edges{\Zd}$ which is \iid\ with model traversal times $\setof{\mttime_j}_{j\in\N}$ satisfying \finitespeed{d} and \expmoment, and $\mu$ will denote the shape function for $\tmeasure$ as defined in Theorem~\ref{shape_thm}.

The goal of this section will be to prove Lemma~\ref{ball_covering_lem} below, which says that with high probability, the internal covering time of a large $\mu$-ball from its center is not greater than its radius plus a small linear factor. Lemma~\ref{ball_covering_lem} will be the key tool for analyzing growth restricted to cone segments in Section~\ref{cone_growth:thin_cone_sec}.
Before proving Lemma~\ref{ball_covering_lem} in Section~\ref{cone_growth:ball_proof_sec}, we first prove two basic large deviations estimates in Section~\ref{cone_growth:basic_dev_sec}. These two results are Lemmas~\ref{shifted_dev_bounds_lem} and \ref{distant_path_lem}, which correspond to steps one and two, respectively, in the ``two-stage covering" strategy from Section~\ref{cone_growth:lattice_path_sec}; that is, Lemma~\ref{shifted_dev_bounds_lem} gives a large deviations estimate for covering a permeating subset (a $\mu$-ball in this case), and Lemma~\ref{distant_path_lem} gives a large deviations estimate for traversing a particular collection of ``short" paths.
Lemma~\ref{shifted_dev_bounds_lem} follows easily from a standard large deviations estimate for the unrestricted process, stated in Lemma~\ref{GM_large_dev_lem} below.
Lemma~\ref{distant_path_lem} is proved using Lemma~\ref{segment_meas_lem} from the previous section. We will use both Lemmas~\ref{shifted_dev_bounds_lem} and \ref{distant_path_lem} in the main proof of Lemma~\ref{ball_covering_lem} as well as in later proofs in Chapters~\ref{cone_growth_chap} and \ref{random_fpc_chap}.
%
%

\subsection{The Basic Large Deviations Estimates}
\label{cone_growth:basic_dev_sec}

Garet and Marchand~\cite{\GMdensity} state the following result bounding the growth of the first-passage percolation process.
It is a generalization of a large deviations estimate proved by Grimmett and Kesten \cite{Grimmett:1984aa} for the passage times in a fixed direction. For a proof that works simultaneously in all directions, see Ahlberg \cite[Proposition~1.5]{Ahlberg:2011aa}.

\begin{lem}[Shape deviation bounds {\cite[Proposition 2.1]{\GMdensity}}]
\label{GM_large_dev_lem}
For any $\epsilon>0$, there exist positive constants $\bigcsymb$ and $\smallcsymb $ such that for all $t>0$,
\[
\Pr \eventthat[1]{ \reacheddset[2]{\mu}{\vect 0}{(1-\epsilon)t}
	\subseteq \rreachedset{\tmeasure}{\vect 0}{\Rd}{t}
	\subseteq \reacheddset[2]{\mu}{\vect 0}{(1+\epsilon)t}}
	\ge 1 - \bigcsymb e^{-\smallcsymb t}.
\]
\end{lem}


In Lemma~\ref{equiv_dev_lem} in Appendix~\ref{asymptotic_chap}, we prove a general result about first-passage growth processes that yields the following stronger version of Lemma~\ref{GM_large_dev_lem} as a corollary.


\begin{cor}[Stronger shape deviation bounds]
\label{strong_dev_cor}
For any $\epsilon>0$ and any $o(t)$ function $\phi\colon \R_+\to\R_+$, there exist positive constants $C_{\ref{strong_dev_cor}} = C_{\ref{strong_dev_cor}}(\epsilon,\phi)$ and $c_{\ref{strong_dev_cor}}= c_{\ref{strong_dev_cor}}(\epsilon,\phi)$ such that for all $t_0\ge 0$,
\begin{align*}
\Pr \eventthat[1]{
	\reacheddset[2]{\mu}{\vect 0}{(1-\epsilon)t+\phi(t)}
	\subseteq \crreachedset{\tmeasure}{\vect 0}{\Rd}{t-}
	\text{ and }
	\rreachedset{\tmeasure}{\vect 0}{\Rd}{t}
	\subseteq \reacheddset[2]{\mu}{\vect 0}{[(1+\epsilon)t-\phi(t)]-}
	\quad \forall t\ge t_0} \quad\\
	\ge 1- C_{\ref{strong_dev_cor}} e^{-c_{\ref{strong_dev_cor}} t_0}.
\end{align*}
\end{cor}

\begin{proof}
The above statement corresponds precisely to the second of the two equivalent statements in Lemma~\ref{equiv_dev_lem}, with $S=\Rd$ and $A=A' = \setof{\vect 0}$, and the first of the two equivalent statements in Lemma~\ref{equiv_dev_lem} is weaker than the statement in Lemma~\ref{GM_large_dev_lem}.
\end{proof}

Using Corollary~\ref{strong_dev_cor}, we can easily obtain the following simple result, bounding the growth of a one-type process started at $\vect v\in\Zd$ using $\mu$-balls centered at some nearby point $\vect z\in\Rd$, rather than requiring the $\mu$-balls to be centered at $\vect v$ itself. In fact, we can trivially get the bound to work \emph{simultaneously for all nearby points} $\vect z$, in particular for all $\vect z\in \cubify{\vect v}$. This formulation will obviate the need to find lattice points approximating some $\vect x\in\Rd$ later on, making many of the longer arguments cleaner than if we used Corollary~\ref{strong_dev_cor} directly. We will use Lemma~\ref{shifted_dev_bounds_lem} in the proofs of Lemma~\ref{ball_covering_lem} and Proposition~\ref{restricted_growth_ub_prop} below, as well as for some results in Chapter~\ref{random_fpc_chap}.

\begin{lem}[Bounding growth with slightly shifted $\mu$-balls]
\label{shifted_dev_bounds_lem}
\newcommand{\thisC}{\bigcref{shifted_dev_bounds_lem}}
\newcommand{\thisc}{\smallcref{shifted_dev_bounds_lem}}
\newcommand{\thisvertex}{\vect v}
\newcommand{\thiscenter}{\vect z}
\newcommand{\thiscube}{\cubify{\thisvertex}}
\newcommand{\thisevent}{\eventref{shifted_dev_bounds_lem}^{\thisvertex}(t_0;\spp)}

For $\thisvertex\in\Zd$, $t_0\ge 0$, and $\spp>0$, define the event
\[
\thisevent \definedas
\bigcap_{\thiscenter\in\thiscube} 
\eventthat[3]{\reacheddset[2]{\mu}{\thiscenter}{(1-\spp)t}
\subseteq \crreachedset{\tmeasure}{\thisvertex}{\Rd}{t-}
\text{ and }
\rreachedset{\tmeasure}{\thisvertex}{\Rd}{t}
\subseteq \reacheddset[2]{\mu}{\thiscenter}{(1+\spp)t-}
\;\; \forall t\ge t_0}.
\]
Then for any $\spp >0$, there exist positive constants $\thisC$ and $\thisc$ such that for any $\thisvertex\in\Zd$ and $t_0\ge 0$,
\[
\Pr \parens[2]{\thisevent} \ge 1-\thisC(\spp) e^{-\thisc(\spp) t_0}.
\]
\end{lem}

\begin{proof}
\newcommand{\thisC}{\bigcref{shifted_dev_bounds_lem}}
\newcommand{\thisc}{\smallcref{shifted_dev_bounds_lem}}

\newcommand{\oldC}{\bigcref{strong_dev_cor}(\spp,\thisdil)}
\newcommand{\oldc}{\smallcref{strong_dev_cor}(\spp,\thisdil)}
\newcommand{\genvertex}{\vect v}
\newcommand{\thisvertex}{\vect 0}
\newcommand{\thiscenter}{\vect z}
\newcommand{\thiscube}{\cubify{\thisvertex}}
\newcommand{\thisevent}{\eventref{shifted_dev_bounds_lem}^{\thisvertex}(t_0;\spp)}
\newcommand{\thisdil}{\kappa}



By translation invariance it suffices to assume $\genvertex = \thisvertex$. To get the desired bound simultaneously for all $\vect z\in \cubify{\vect 0}$, observe that if we set $\thisdil = \frac{1}{2} \dilnorm{\linfty{d}}{\mu}$, then by the triangle inequality for $\mu$, for any $r\ge 0$ we have
\[
\bigcup_{\vect z\in \cubify{\vect 0}} \reacheddset{\mu}{\vect z}{r}
\subseteq \reacheddset{\mu}{\vect 0}{r+\thisdil}
\quad\text{and}\quad
\reacheddset[2]{\mu}{\vect 0}{[r-\thisdil]-}
\subseteq \bigcap_{\vect z\in \cubify{\vect 0}} \reacheddset{\mu}{\vect z}{r-}.
\]
Applying Corollary~\ref{strong_dev_cor} with $\phi(t)=\thisdil$ (a constant function), this implies that for any $t_0\ge 0$ we have
\begin{align*}
1 &- \oldC e^{\oldc t_0} \\
&\le \Pr \eventthat[1]{
	\reacheddset[2]{\mu}{\vect 0}{(1-\spp)t+\thisdil}
	\subseteq \crreachedset{\tmeasure}{\vect 0}{\Rd}{t-}
	\text{ and }
	\rreachedset{\tmeasure}{\vect 0}{\Rd}{t}
	\subseteq \reacheddset[2]{\mu}{\vect 0}{[(1+\spp)t-\thisdil]-}
	\quad \forall t\ge t_0} \\
&\le \Pr \parens[1]{
	\bigcap_{t\ge t_0} \bigcap_{\vect z\in \cubify{\vect 0}}
	\eventthat[1]{
	\reacheddset[2]{\mu}{\vect z}{(1-\spp)t}
	\subseteq \crreachedset{\tmeasure}{\vect 0}{\Rd}{t-}
	\text{ and }
	\rreachedset{\tmeasure}{\vect 0}{\Rd}{t}
	\subseteq \reacheddset[2]{\mu}{\vect z}{(1+\spp)t-}}}\\
&=\Pr \parens[1]{\thisevent}.
\end{align*}
Thus we can take $\thisC(\spp) \definedas \oldC$ and $\thisc(\spp) \definedas \oldc$.
\end{proof}

%


Using Lemma~\ref{segment_meas_lem} from Section~\ref{cone_growth:lattice_path_sec} above, we prove the following result, which will be needed in the proofs of Lemmas~\ref{ball_covering_lem} and \ref{path_tail_bound_lem}.
Lemma~\ref{distant_path_lem} says that with high probability, the traversal measure of every radial lattice path connecting two concentric $\mu$-spheres is bounded above by a fixed constant times the difference in the spheres' radii. Similar to the statement of Lemma~\ref{shifted_dev_bounds_lem}, the event in Lemma~\ref{distant_path_lem} simultaneously treats all pairs of spheres centered in a unit cube around the origin.

\begin{lem}[Traversal of radial paths in a $\mu$-spherical shell]
\label{distant_path_lem}

\newcommand{\thisrv}{\tmeasureof{\edge}}
\newcommand{\thisevent}{\symbolref{\eventsymb}{distant_path_lem}(r;\epsilon)}

\newcommand{\thiskonst}{\konstref{distant_path_lem}}

\newcommand{\thisC}{\bigcref{distant_path_lem}}
\newcommand{\thisc}{\smallcref{distant_path_lem}}

Let $\konstref{distant_path_lem}= \konstref{segment_meas_lem}^{\mu} \parens[2]{\E \thisrv}$. For $\epsilon\in (0,1]$ and $r\ge 0$, define the event
\[
\thisevent =
	\bigcap_{\vect z\in \cubify{\vect 0}}
	\bigcap_{t\ge r}
	\eventthat[1]{
	\ttime \parens[3]{\segment[2]{\vect z+ \smash{(}1-\epsilon\smash{)}\vect y}{\vect z+\vect y}}
	\le \thiskonst \cdot \epsilon t
	\ \ \text{for all } \vect y\in \bd \reacheddset{\mu}{\vect 0}{t}}.
\]
Then for any $\epsilon\in (0,1]$, there exist positive constants $\thisC(\epsilon)$ and $\thisc(\epsilon)$ such that for any $r\ge 0$,
\[
\Pr \parens[2]{\thisevent }\ge 1-\thisC(\epsilon) e^{-\thisc(\epsilon)r}.
\]
\end{lem}

\begin{proof}

\newcommand{\thisevent}{\eventref{distant_path_lem}(r;\epsilon)}
\newcommand{\thiskonst}{\konstref{distant_path_lem}}
\newcommand{\thisC}{\bigcref{distant_path_lem}}
\newcommand{\thisc}{\smallcref{distant_path_lem}}

\newcommand{\oldC}{\bigcref{segment_meas_lem}}
\newcommand{\oldc}{\smallcref{segment_meas_lem}}

\newcommand{\oldCp}{\bigcref{segment_meas_lem}(1+\thisdil,\mu)}
\newcommand{\oldcp}{\smallcref{segment_meas_lem}(1+\thisdil,\mu)}

\newcommand{\thisdil}{\kappa}

\newcommand{\latticepoints}[3]{U^{#1}_{\segment{#2}{#3}}}
\newcommand{\thiscylconst}{Q_{d,\mu}}

\newcommand{\tevent}[1]{\eventsymb_{#1}}
\newcommand{\jevent}[1]{\widetilde{\eventsymb}_{#1}}


Fix $\epsilon\in (0,1]$. For any $t\ge 0$, define the event
\[
\tevent{t} \definedas 
\bigcap_{\vect z\in \cubify{\vect 0}}
	\eventthat[1]{
	\ttime \parens[1]{\segment[2]{\vect z+ \smash{(}1-\epsilon\smash{)}\vect y}{\vect z+\vect y}}
	\le \konstref{distant_path_lem}\cdot \epsilon t
	\ \ \text{for all } \vect y\in \bd \reacheddset{\mu}{\vect 0}{t}},
\]
and for $j\in\N$, define the event
$
\jevent{j} \definedas \bigcap_{j\le t<j+1} \tevent{t}.
$
Then for any $r\ge 0$,
\begin{equation}
\label{distant_path:event_containment_eqn}
\thisevent = \bigcap_{t\ge r} \tevent{t} \supseteq \bigcap_{j\ge \floor{r}}^\infty \jevent{j}.
\end{equation}
We will obtain a lower bound on the probability of $\thisevent$ by showing that each of the events $\jevent{j}$ occurs on a finite intersection of events of the form in Lemma~\ref{segment_meas_lem}. In order to do this, we need to count pairs of lattice points near the endpoints of the intervals in the definition of $\tevent{t}$.

Let $\thisdil = \dilnorm{\lspace{\infty}{d}}{\mu}$, and note that $\sup_{\vect z\in \cubify{\vect 0}} \mu(\vect z) = \thisdil/2$. Given $\vect z\in\cubify{\vect 0}$ and $\vect y\in  \bd \reacheddset{\mu}{\vect 0}{t}$, choose $\vect v,\vect u\in\Zd$ with $\vect z+\vect y\in \cubify{\vect v}$ and $\vect z+(1-\epsilon)\vect y\in \cubify{\vect u}$. Then $\dist{\mu} \parens{\vect v,\vect z+\vect y} \le {\thisdil}/{2}$ and $\dist{\mu} \parens[2]{\vect u,\vect z+(1-\epsilon)\vect y} \le {\thisdil}/{2}$, so the triangle inequality implies that
\begin{equation}\label{distant_path:uv_bounds}
t-\thisdil\le \mu(\vect v)\le t+\thisdil
\text{\quad and \quad}
\epsilon t - \thisdil \le \dist{\mu} \parens{\vect u,\vect v} \le \epsilon t+\thisdil. 
\end{equation}
Furthermore, since the point $\vect z+(1-\epsilon)\vect y$ is contained in the line segment $\segment{\vect z}{\vect z+\vect y}$, we have $\dist{\lspace{\infty}{d}} \parens[2]{\vect u, \segment{\vect z}{\vect z+\vect y}} \le 1/2$. Now note that since $\vect z\in \cubify{\vect 0}$ and $\vect z+\vect y\in \cubify{\vect v}$, the entire line segment $\segment{\vect z}{\vect z+\vect y}$ is contained in $\conv \parens[2]{\cubify{\vect 0}\cup \cubify{\vect v}} = \segment{\vect 0}{\vect v} +\cubify{\vect 0}$, which is the $\frac{1}{2}$-neighborhood of $\segment{\vect 0}{\vect v}$ in the $\lspace{\infty}{d}$ norm. Therefore,
\begin{equation}
\label{distant_path:segment_bound}
\dist{\lspace{\infty}{d}} \parens[2]{\vect u, \segment{\vect 0}{\vect v}}
\le \dist{\lspace{\infty}{d}} \parens[2]{\vect u, \segment{\vect z}{\vect z+\vect y}}
+\sup_{\vect x\in \segment{\vect z}{\vect z+\vect y}}
\dist{\lspace{\infty}{d}} \parens[2]{\vect x, \segment{\vect 0}{\vect v}}
\le \frac{1}{2}+\frac{1}{2} = 1.
\end{equation}
We wish to obtain an upper bound on the number of lattice points $\vect u$ which satisfy \eqref{distant_path:segment_bound} and are also within some specified distance of $\vect v$. Namely, for each $\vect v\in\Zd$ and $a,b\ge 0$, we wish to bound the cardinality of the set
\begin{equation}\label{distant_path:lattice_set_def}
\latticepoints{\vect v}{a}{b} \definedas \setbuilder[1]{\vect u\in\Zd}
{\dist{\lspace{\infty}{d}} \parens[2]{\vect u, \segment{\vect 0}{\vect v}} \le 1
\text{ and }  \dist{\mu}(\vect u,\vect v)\in \segment{a}{b}}.
\end{equation}
Using an argument similar to the one in the proof of Lemma~\ref{segment_meas_lem}, one can show that for any $\vect v\in\Zd$ and $0\le a \le b$, the set $\latticepoints{\vect v}{a}{b}$ is contained in a Euclidean cylinder of some radius depending only on $d$ and of height that is linear in $(b-a)$, so there is some constant $\thiscylconst<\infty$ such that
\begin{equation}\label{distant_path:lattice_set_bound}
\abs[1]{\latticepoints{\vect v}{a}{b}} \le \thiscylconst\cdot \brackets[2]{(b-a)\join 1}
\text{\quad for all $\vect v\in\Zd$ and $0\le a\le b$.}
\end{equation}
Now, returning to the definition of the event $\tevent{t}$, the inequalities \eqref{distant_path:uv_bounds}, \eqref{distant_path:segment_bound}, and the definition \eqref{distant_path:lattice_set_def} imply that
\begin{align*}
\tevent{t}
&\supseteq \bigcap_{\substack{\vect v\in\Zd \\ t-\thisdil \le\mu(\vect v)\le t+\thisdil}}
	\bigcap_{\vect u\in \latticepoints{\vect v}{\epsilon t-\thisdil}{\epsilon t+\thisdil}}
	\eventthat[1]{ \tmeasure \parens[2]{\segment{\vect x}{\vect x'}} \le\thiskonst \cdot \epsilon t
	\text{ for all } \vect x\in \cubify{\vect u} \text{ and } \vect x'\in\cubify{\vect v}}\\
&= \bigcap_{\substack{\vect v\in\Zd \\ t-\thisdil \le\mu(\vect v)\le t+\thisdil}}
	\bigcap_{\vect u\in \latticepoints{\vect v}{\epsilon t-\thisdil}{\epsilon t+\thisdil}}
	\eventref{segment_meas_lem} (\vect u,\vect v, \thiskonst \epsilon t).
\end{align*}
Recall that $\jevent{j} = \bigcap_{j\le t<j+1} \tevent{t}$, and observe that since the events $\eventref{segment_meas_lem} (\vect u,\vect v, \thiskonst \epsilon t)$ are increasing in $t$, it then follows that
\begin{align*}
\jevent{j}
&\supseteq  \bigcap_{j\le t<j+1} 
	\bigcap_{\substack{\vect v\in\Zd \\ t-\thisdil \le\mu(\vect v)\le t+\thisdil}}
	\bigcap_{\vect u\in \latticepoints{\vect v}{\epsilon t-\thisdil}{\epsilon t+\thisdil}}
	\eventref{segment_meas_lem} (\vect u,\vect v, \thiskonst \epsilon t)\\
&\supseteq \bigcap_{\substack{\vect v\in\Zd \\ j-\thisdil \le \mu(\vect v)\le j+1+\thisdil}}
	\bigcap_{\vect u\in \latticepoints{\vect v}{\epsilon j-\thisdil}{\epsilon (j+1)+\thisdil}}
	\eventref{segment_meas_lem} (\vect u,\vect v, \thiskonst \epsilon j).
\end{align*}
Now note that since $\epsilon\le 1$, if $\vect u\in  \latticepoints{\vect v}{\epsilon j-\thisdil}{\epsilon (j+1)+\thisdil}$, then $\dist{\mu}(\vect u,\vect v) \le \epsilon j+ (1+\thisdil)$. Thus, applying Lemma~\ref{segment_meas_lem} with $r=\epsilon j$ and $\phi(r) = 1+\thisdil$, and using Lemma~\ref{norm_shell_lem} and the bound \eqref{distant_path:lattice_set_bound} to estimate the number of lattice points in the sum, we get
\begin{align}
\label{distant_path:Ej_containment_eqn}
\Pr \parens[1]{\jevent{j}^\complement}
&\le \sum_{\substack{\vect v\in\Zd \\ j-\thisdil \le \mu(\vect v)\le j+1+\thisdil}}
	\sum_{\vect u\in \latticepoints{\vect v}{\epsilon j-\thisdil}{\epsilon (j+1)+\thisdil}}
	\Pr \parens[1]{\eventref{segment_meas_lem}
		(\vect u,\vect v, \thiskonst \epsilon j)^\complement} \notag\\
&\le \mushellbound (j+1+\thisdil)^{d-1} (1+2\thisdil)
	\cdot \thiscylconst (1+2\thisdil)
	\cdot \oldCp e^{-\oldcp \epsilon j}.
\end{align}
Set $\oldC = \oldCp$ and $\oldc = \oldcp$. Then combining \eqref{distant_path:event_containment_eqn} and \eqref{distant_path:Ej_containment_eqn}, and using Lemma~\ref{poly_geo_series_lem} to sum the resulting series, we get
\begin{align*}
\Pr \parens[1]{\eventcomplement{\thisevent}}
&\le \sum_{j= \floor{r}}^\infty \Pr \parens[1]{\eventcomplement{\jevent{j}}}\\
&\le \mushellbound \thiscylconst (1+2\thisdil)^2 \cdot \oldC
	\sum_{j= \floor{r}}^\infty (j+1+\thisdil)^{d-1} e^{-\oldc \epsilon j}\\
&\le \mushellbound \thiscylconst (1+2\thisdil)^2 \cdot \oldC
	\cdot \frac{4[(2d-2)\join (1+\thisdil)]^{d-1}}{[(\oldc\epsilon)\meet 1]^{d}}
	\cdot e^{-\frac{\oldc \epsilon}{2} \floor{r}}\\
&\le \thisC(\epsilon) e^{-\thisc(\epsilon) r},
\end{align*}
where $\thisC(\epsilon)$ and $\thisc(\epsilon)$ are some positive constants depending only on $\epsilon$, $d$, and the distribution of $\tmeasureof{\edge}$.
\end{proof}

\subsection{Proof of Large Deviations Estimate for Covering a $\mu$-Ball}
\label{cone_growth:ball_proof_sec}


Now we are ready to prove the main large deviations estimate for growth restricted to a $\mu$-ball, which will be the main ingredient for analyzing growth restricted to $\mu$-cone segments in the next section. As in Lemmas~\ref{shifted_dev_bounds_lem} and \ref{distant_path_lem}, we consider an event that simultaneously treats all $\mu$-balls centered near a given lattice point, in order to simplify later arguments. This convention will be propagated throughout most of the main results in Chapters~\ref{cone_growth_chap} and \ref{random_fpc_chap}.

\begin{lem}[Internal covering time for $\mu$-balls]
\label{ball_covering_lem}
For any $\epsilon>0$ there exist positive constants $C_{\ref{ball_covering_lem}}$ and $c_{\ref{ball_covering_lem}}$ such that for all $\vect v\in\Zd$ and $r_0\ge 0$,
\[
\Pr \eventthat[1]{
	\crreachedset[2]{\tmeasure}{\vect v}{ \smash[b]{ \reacheddset{\mu}{\vect z}{r} }}{(1+\epsilon)r}
	= \reacheddset{\mu}{\vect z}{r}
	\text{\rm\ for all $\vect z\in\cubify{\vect v}$ and } r\ge r_0}
	\ge 1- C_{\ref{ball_covering_lem}} e^{-c_{\ref{ball_covering_lem}} r_0}.
\]
\end{lem}

\begin{proof}
\newcommand{\radialsegment}[3]{L_{#1}^{#2}(#3)}

%
By translation invariance it suffices to assume $\vect v = \vect 0$. Given $\epsilon>0$, let $\epsilon_1 = \epsilon_1(\epsilon) \definedas \frac{\epsilon}{2 \konstref{distant_path_lem}}$, let $\alpha = \alpha(\epsilon) \definedas \frac{1-\epsilon_1}{1+\epsilon_1}$, and for each $r\ge r_0$ let $t_r \definedas (1+\epsilon_1)^{-1}r$. That is, $t_r$ is chosen so that $(1+\epsilon_1) t_r = r$, so the unrestricted random process is likely to be contained in the ball $\reacheddset{\mu}{\vect z}{r}$ at time $t_r$, and then $\alpha r = (1-\epsilon_1)t_r$ is the radius of a slightly smaller ball, which the random process is likely to have covered by time $t_r$.
Here, the smaller ball $\reacheddset{\mu}{\vect z}{\alpha r}$ plays the role of the permeating subset; namely, it is a $(1-\alpha)r$-permeating subset of $\reacheddset{\mu}{\vect z}{r}$. The value of $\epsilon_1$ was chosen so that $(1-\alpha)r \approx \frac{\epsilon}{\konstref{distant_path_lem}}\cdot r$, so that Lemma~\ref{distant_path_lem} implies that the remainder of the large ball is likely to be reached within an additional time of $\epsilon r$, via approximately straight paths.

More concretely, the outline of the proof is as follows. First we combine Lemma~\ref{shifted_dev_bounds_lem} with Lemma~\ref{continuum_containment_lem} to show that the restricted process covers the union $\bigcup_{\vect z\in \cubify{\vect 0}} \reacheddset{\mu}{\vect z}{\alpha r}$ by time $t_r$ with high probability. Then we use Lemma~\ref{distant_path_lem} to show that with high probability, the remaining spherical shells $\reacheddset{\mu}{\vect z}{r}\setminus \reacheddset{\mu}{\vect z}{\alpha r}$ are covered within an additional time of $\epsilon r$, and so the total covering time is at most $t_r+\epsilon r < (1+\epsilon)r$. More explicitly, the chaining property (Lemma~\ref{chaining_lem}) implies that on any realization of $\tmeasure$, for each $\vect z\in \cubify{\vect 0}$ we have
\begin{equation}
\label{ball_covering:chaining_eq}
\ctime[1]{}{\reacheddset{\mu}{\vect z}{r}}{\vect 0}{\reacheddset{\mu}{\vect z}{r}}
\le \ctime[1]{}{\reacheddset{\mu}{\vect z}{r}}{\vect 0}{\reacheddset{\mu}{\vect z}{\alpha r}}
	+ \ctime[1]{}{\reacheddset{\mu}{\vect z}{r}}
		{\bd \reacheddset{\mu}{\vect z}{\alpha r}}
		{\,\reacheddset{\mu}{\vect z}{r} \setminus \reacheddset{\mu}{\vect z}{\alpha r}}.
\end{equation}
Our goal, then, is to show that with exponentially high probability in $r_0$, for all $r\ge r_0$ and all $\vect z\in \cubify{\vect 0}$, the first term on the right in \eqref{ball_covering:chaining_eq} is bounded above by $r$, and the second term on the right is bounded above by $\epsilon r$.

We start with the first term on the right in \eqref{ball_covering:chaining_eq}. To get the desired bound simultaneously for all $\vect z\in \cubify{\vect 0}$, observe that if we set $\delta = \frac{1}{2} \dilnorm{\ell^\infty}{\mu}$, then (by the triangle inequality) for any $t\ge 0$,
\[
\bigcup_{\vect z\in \cubify{\vect 0}} \reacheddset{\mu}{\vect z}{t}
\subseteq \reacheddset{\mu}{\vect 0}{t+\delta}
\quad\text{and}\quad
\reacheddset{\mu}{\vect 0}{t-\delta}
\subseteq \bigcap_{\vect z\in \cubify{\vect 0}} \reacheddset{\mu}{\vect z}{t}.
\]
Applying Corollary~\ref{strong_dev_cor} with $\phi(t)=\delta$, this implies that for any $r_0\ge 0$ we have
\begin{align}\label{ball_covering:orig_eqn}
1 &- C_{\ref{strong_dev_cor}}(\epsilon_1,\delta) 
	e^{-c_{\ref{strong_dev_cor}}(\epsilon_1,\delta) t_{r_0}} \notag\\
&\le \Pr \eventthat[1]{
	\reacheddset[2]{\mu}{\vect 0}{(1-\epsilon_1)t_r+\delta}
	\subseteq \crreachedset{\tmeasure}{\vect 0}{\Rd}{t_r-}
	\text{ and }
	\rreachedset{\tmeasure}{\vect 0}{\Rd}{t_r}
	\subseteq \reacheddset[2]{\mu}{\vect 0}{(1+\epsilon_1)t_r-\delta}
	\quad \forall r\ge {r_0}}\notag \\
&\le \Pr \parens[1]{
	\bigcap_{r\ge r_0} \bigcap_{\vect z\in \cubify{\vect 0}}
	\eventthat[1]{
	\reacheddset[2]{\mu}{\vect z}{(1-\epsilon_1)t_r}
	\subseteq \crreachedset{\tmeasure}{\vect 0}{\Rd}{t_r-}
	\text{ and }
	\rreachedset{\tmeasure}{\vect 0}{\Rd}{t_r}
	\subseteq \reacheddset[2]{\mu}{\vect z}{(1+\epsilon_1)t_r}}}.
\end{align}
Note that $\reacheddset[2]{\mu}{\vect z}{(1+\epsilon_1)t_r} = \reacheddset{\mu}{\vect z}{r}$, and by Lemma~\ref{continuum_containment_lem} and Lemma~\ref{restricted_proc_equiv_lem},
\begin{align*}
\eventthat[1]{\rreachedset{\tmeasure}{\vect 0}{\Rd}{t_r}
	\subseteq \reacheddset{\mu}{\vect z}{r}}
&\subseteq
\eventthat[1]{\crreachedset{\tmeasure}{\vect 0}{\Rd}{t_r-}
	= \crreachedset{\tmeasure}{\vect 0}{\reacheddset{\mu}{\vect z}{r}}{t_r-}},
\end{align*}
so \eqref{ball_covering:orig_eqn} implies that
\begin{equation}\label{ball_covering:restrict_eqn}
\Pr \parens[1]{
	\bigcap_{r\ge r_0} \bigcap_{\vect z\in \cubify{\vect 0}}
	\eventthat[1]{
	\reacheddset[2]{\mu}{\vect z}{(1-\epsilon_1)t_r}
	\subseteq \crreachedset{\tmeasure}{\vect 0}{\reacheddset{\mu}{\vect z}{r}}{t_r-}}}
\ge 1 - \bigcref{strong_dev_cor}(\epsilon_1,\delta) 
	e^{-\smallcref{strong_dev_cor}(\epsilon_1,\delta) t_{r_0}}.
\end{equation}
Since $r>t_r$, the process at time $r$ contains the process at time $t_r-$, so \eqref{ball_covering:restrict_eqn} implies that
\begin{equation}\label{ball_covering:mostly_covered_eqn}
\Pr \parens[1]{
	\bigcap_{r\ge r_0} \bigcap_{\vect z\in \cubify{\vect 0}}
	\eventthat[1]{
	\reacheddset{\mu}{\vect z}{\alpha r}
	\subseteq \crreachedset{\tmeasure}{\vect 0}{\reacheddset{\mu}{\vect z}{r}}{r}}}
\ge 1 - C_{\ref{strong_dev_cor}}(\epsilon_1,\delta) 
	e^{-\smallc{1}(\epsilon) r_0},
\end{equation}
where we have set $\smallc{1}(\epsilon) = \parens[2]{1+\epsilon_1(\epsilon)}^{-1} \smallcref{strong_dev_cor}\parens[2]{\epsilon_1(\epsilon),\delta}$ and have recalled that $(1-\epsilon_1)t_r = \alpha r$.  The above event, call it  $\symboleqref{\eventsymb}{ball_covering:mostly_covered_eqn}$, says that $\ctime[1]{}{\reacheddset{\mu}{\vect z}{r}}{\vect 0}{\reacheddset{\mu}{\vect z}{\alpha r}} \le r$ for all $r\ge r_0$ and all $\vect z\in \cubify{\vect 0}$, so \eqref{ball_covering:mostly_covered_eqn} provides the desired bound for the first term in \eqref{ball_covering:chaining_eq}.

Now we turn to the final term in \eqref{ball_covering:chaining_eq}. The idea is that each point in the spherical shell $\reacheddset{\mu}{\vect z}{r} \setminus \reacheddset{\mu}{\vect z}{\alpha r}$ can be reached from the boundary of the inner ball $\reacheddset{\mu}{\vect z}{\alpha r}$ via a radial line segment of length $(1-\alpha)r \approx \frac{\epsilon}{\konstref{distant_path_lem}}\cdot r$, and  Lemma~\ref{distant_path_lem} implies that with high probability, the covering time of all such line segments is bounded above by $\epsilon r$. We now make this argument more precise.

First note that 
$\reacheddset{\mu}{\vect 0}{r}\setminus \interior{\reacheddset{\mu}{\vect 0}{\alpha r}}
= \bigcup_{\vect y\in \bd \reacheddset{\mu}{\vect 0}{r}}
	\segment[2]{\alpha \vect y}{\vect y}$.
If we set $\radialsegment{\alpha}{\vect z}{\vect y} = \vect z+ \segment[2]{\alpha \vect y}{\vect y}$ for $\vect z\in \cubify{\vect 0}$ and $\vect y\in \bd \reacheddset{\mu}{\vect 0}{r}$, then translating everything by $\vect z$ gives
\[
\reacheddset{\mu}{\vect z}{r}\setminus \interior{\reacheddset{\mu}{\vect z}{\alpha r}}
= \bigcup_{\vect y\in \bd \reacheddset{\mu}{\vect 0}{r}} \radialsegment{\alpha}{\vect z}{\vect y}.
\]
Using the fact that for any $\vect y\in \bd \reacheddset{\mu}{\vect 0}{r}$ we have $(\vect z+ \alpha \vect y) \in \bd \reacheddset{\mu}{\vect z}{\alpha r}\cap \radialsegment{\alpha}{\vect z}{\vect y}\subset \reacheddset{\mu}{\vect z}{r}$, the decomposition and monotonicity properties and the total traversal measure bound from Lemma~\ref{covering_time_properties_lem} then imply that
\begin{align}\label{ball_covering:shell_bound_eq}
\ctime[3]{}{\reacheddset{\mu}{\vect z}{r}}
	{\bd \reacheddset{\mu}{\vect z}{\alpha r}}
	{\: \reacheddset{\mu}{\vect z}{r}\setminus \reacheddset{\mu}{\vect z}{\alpha r}}
&\le \sup_{\vect y\in \bd \reacheddset{\mu}{\vect 0}{r}}
	\ctime[2]{}{\reacheddset{\mu}{\vect z}{r}}
	{\vect z+ \alpha \vect y}{\radialsegment{\alpha}{\vect z}{\vect y}} \notag\\
&\le \sup_{\vect y\in \bd \reacheddset{\mu}{\vect 0}{r}}
	\tmeasure \parens[2]{\radialsegment{\alpha}{\vect z}{\vect y}}.
\end{align}
Now define the event
\begin{equation}\label{ball_covering:shell_covered_eq}
\symboleqref{\eventsymb}{ball_covering:shell_covered_eq}
=\eventthat[1]{ 
	\ctime[3]{}{\reacheddset{\mu}{\vect z}{r}}
	{\bd \reacheddset{\mu}{\vect z}{\alpha r}}
	{\: \reacheddset{\mu}{\vect z}{r}\setminus \reacheddset{\mu}{\vect z}{\alpha r}}
\le \epsilon r
\text{ for all $\vect z\in \cubify{\vect 0}$ and } r\ge r_0}.
\end{equation}
The bound in \eqref{ball_covering:shell_bound_eq} implies that $\symboleqref{\eventsymb}{ball_covering:shell_covered_eq}$ occurs on the event $\symbolref{\eventsymb}{distant_path_lem}(r_0;1-\alpha)$, as follows:
\begin{align*}
\symboleqref{\eventsymb}{ball_covering:shell_covered_eq}
&= \bigcap_{\vect z\in \cubify{\vect 0}} \bigcap_{r\ge r_0}
	\eventthat[1]{ \ctime[3]{}{\reacheddset{\mu}{\vect z}{r}}
	{\bd \reacheddset{\mu}{\vect z}{\alpha r}}
	{\: \reacheddset{\mu}{\vect z}{r}\setminus \reacheddset{\mu}{\vect z}{\alpha r}}
	\le \epsilon r }\\
&\supseteq \bigcap_{\vect z\in \cubify{\vect 0}} \bigcap_{r\ge r_0}
	\eventthat[3]{ \tmeasure \parens[2]{\radialsegment{\alpha}{\vect z}{\vect y}}
	\le \epsilon r \text{ for all } \vect y\in \bd \reacheddset{\mu}{\vect 0}{r}}\\
&\supseteq \bigcap_{\vect z\in \cubify{\vect 0}}
	\bigcap_{r \ge r_0}
	\eventthat[3]{ \tmeasure \parens[2]{\radialsegment{\alpha}{\vect z}{\vect y}}
	\le \konstref{distant_path_lem}\cdot (1-\alpha)r
	\ \ \text{for all } \vect y\in \bd \reacheddset{\mu}{\vect 0}{r}}\\
&= \symbolref{\eventsymb}{distant_path_lem}(r_0;1-\alpha),
\end{align*}
where the inclusion in the second-to-last line holds because
\[
\konstref{distant_path_lem}\cdot (1-\alpha) =  \konstref{distant_path_lem}\cdot \frac{2\epsilon_1}{1+\epsilon_1} < \konstref{distant_path_lem} \cdot 2\epsilon_1 = \epsilon.
\]
Therefore, by Lemma~\ref{distant_path_lem}, we have
\begin{equation}\label{ball_covering:shell_event_prob_eq}
\Pr \parens[2]{\symboleqref{\eventsymb}{ball_covering:shell_covered_eq}}
\ge \Pr \parens[2]{ \symbolref{\eventsymb}{distant_path_lem}(r_0;1-\alpha)}
\ge 1-\bigcref{distant_path_lem}(1-\alpha) 
		e^{-\smallcref{distant_path_lem}(1-\alpha) r_0}.
\end{equation}
Finally, \eqref{ball_covering:chaining_eq}, \eqref{ball_covering:mostly_covered_eqn}, and \eqref{ball_covering:shell_event_prob_eq} imply that
\begin{align*}
\Pr \eventthat[1]{
	\crreachedset[2]{\tmeasure}{\vect 0}{ \smash[b]{ \reacheddset{\mu}{\vect z}{r} }}{(1+\epsilon)r}
	= \reacheddset{\mu}{\vect z}{r}
	\text{\rm\ for all $\vect z\in\cubify{\vect 0}$ and } r\ge r_0}
&\ge \Pr \parens[2]{\symboleqref{\eventsymb}{ball_covering:mostly_covered_eqn}
		\cap \symboleqref{\eventsymb}{ball_covering:shell_covered_eq}}\\
&\ge 1- \bigcref{ball_covering_lem}(\epsilon) e^{-\smallcref{ball_covering_lem}(\epsilon) r_0},
\end{align*}
where $\bigcref{ball_covering_lem}(\epsilon) = \bigcref{strong_dev_cor}\parens[2]{\epsilon_1(\epsilon),\delta}+ \bigcref{distant_path_lem}\parens[2]{1-\alpha(\epsilon)}$ and $\smallcref{ball_covering_lem}(\epsilon) =  \smallc{1}(\epsilon) \wedge \smallcref{distant_path_lem}\parens[2]{1-\alpha(\epsilon)}$.
\end{proof}

\section{Large Deviations Estimate for Covering a Thin $\mu$-Cone Segment}
\label{cone_growth:thin_cone_sec}

The goal of this section is to use Lemma~\ref{ball_covering_lem} from the previous section to prove Lemma~\ref{thin_cone_covering_lem} below, bounding the internal covering time of a ``thin" $\mu$-cone segment. Lemma~\ref{thin_cone_covering_lem} is the key technical result of the present chapter, and we will use it to prove our main theorems about growth in star sets, cones, and tubes in the next two sections. In order to prove Lemma~\ref{thin_cone_covering_lem}, we need four preliminary lemmas. The first two are purely geometric in nature, while the second two combine geometric arguments with the probability estimates of the previous section to bound the probability of events that will be important for the main proof.

\subsection{Preliminary Lemmas Needed for the Proof of Lemma~\ref{thin_cone_covering_lem}}

\newcommand{\centersymb}{\dirsymb}
\newcommand{\centerpt}[1]{\centersymb_{#1}} 

\newcommand{\gballsymb}{\mathfrak{B}}
\newcommand{\gball}[3]{\gballsymb^{#1}_{#2_{#3}}} 
\newcommand{\gballchain}[4]{\gballsymb^{#1,\direction{#2}}_{#3,#4}} 

\newcommand{\gballeventsymb}{F}
\newcommand{\gballevent}[5]{\gballeventsymb^{#1}_{{#2}_{#3}} \parens{#4; #5}}
\newcommand{\gballchainevent}[6]{\gballeventsymb^{#1,\direction{#2}}_{#3,#4} \parens{#5;#6}}


The first preliminary lemma defines a geometrically expanding family of $\mu$-balls and highlights some elementary properties of this construction that will be needed in the main proof.

\begin{lem}[A cone-permeating chain of $\mu$-balls]
\label{geometric_balls_lem}
Fix $\thp\in (0,1]$. For each $\dirsymb\in \sphere{d-1}{\mu}$ and $k\in \Z$, let $\centerpt{k} = \centerpt{k}(\thp) = (1+\thp)^k \dirsymb$, and for any $\vect z\in\Rd$, define $\gball{\vect z}{\dirsymb}{k} = \gball{\vect z}{\dirsymb}{k}(\thp)$ by
\[
\gball{\vect z}{\dirsymb}{k} \definedas 
\vect z + \reacheddset[2]{\mu}{\centerpt{k}}{\thp \mu(\centerpt{k})}
= \reacheddset[2]{\mu}{\vect z+ \centerpt{k}}{\thp (1+\thp)^k}.
\]
Further, for $m,n\in\Z$, define $\gballchain{\vect z}{\dirsymb}{m}{n} = \gballchain{\vect z}{\dirsymb}{m}{n}(\thp) \definedas \bigcup_{k=m}^n \gball{\vect z}{\dirsymb}{k}(\thp)$. Then the sequence $\setof{\gball{\vect z}{\dirsymb}{k}}_{k\in\Z}$ forms a chain of intersecting $\mu$-balls with geometrically increasing radii, such that for any $m,n\in\Z$, the union $\gballchain{\vect z}{\dirsymb}{m}{n}$ contains the line segment $\vect z+ \segment[2]{(1-\thp) \centerpt{m}}{\centerpt{n+1}}$ and has convex hull equal to the \conesegmenttext{\mu}{\thp} $\conesegment[2]{\mu}{\thp}{\vect z}{\dirvec}{(1+\thp)^m}{(1+\thp)^n}.$

\end{lem}

\begin{proof}
Observe that the family $\setof{\gball{\vect z}{\dirsymb}{k}}_{k\in\Z}$ is defined so that:
\begin{itemize}
\item All the balls are centered on the ray $\vect z+\dirvec$ (the center of $\gball{\vect z}{\dirsymb}{k}$ is $\vect z +\centerpt{k}$);
\item The center of $\gball{\vect z}{\dirsymb}{k+1}$ lies on the boundary of $\gball{\vect z}{\dirsymb}{k}$;
\item The radius of $\gball{\vect z}{\dirsymb}{0}$ is $\thp$, and for all $k\in\Z$, the radius of $\gball{\vect z}{\dirsymb}{k+1}$ is $1+\thp$ times as large as the radius of $\gball{\vect z}{\dirsymb}{k}$.
\end{itemize}
All the statements in the lemma follow easily from these properties.
\end{proof}

The union $\gballchain{\vect z}{\dirsymb}{m}{n}$ from Lemma~\ref{geometric_balls_lem} will play the role of the permeating subset in the proof of Lemma~\ref{thin_cone_covering_lem} below. The chain of $\mu$-balls $\gballchain{\vect z}{\dirsymb}{m}{n}$ is constructed so that (1) the chaining property (Lemma~\ref{chaining_lem}) implies that the the covering time of the entire union $\gballchain{\vect z}{\dirsymb}{m}{n}$ from the initial center $\vect z+\centerpt{m}$ is bounded by the sum of the covering times of all the balls from their centers, and (2) the sum of these covering times is likely to be on the order of $\distnew{\mu}{\centerpt{m}} {\centerpt{n+1}}$.
To make this statement precise, we first formally define the relevant covering events and indicate their relationship in the following lemma.

\begin{lem}[Covering the permeating chain one $\mu$-ball at a time]
\label{ball_chaining_lem}
Fix $\thp\in (0,1)$, and for any $\dirsymb\in \sphere{d-1}{\mu}$ and  $\vect z\in\Rd$, let the sequences $\setof{\centerpt{k}(\thp)}_{k\in\Z}$, $\setof{\gball{\vect z}{\dirsymb}{k}(\thp)}_{k\in\Z}$, and $\setof{\gballchain{\vect z}{\dirsymb}{m}{n}(\thp)}_{m,n\in\Z}$ be defined as in Lemma~\ref{geometric_balls_lem}. 
For each $k,m,n\in\Z$ and $\beta> 1$, define events
\begin{align*}
\gballevent{\vect z}{\dirsymb}{k}{\thp}{\beta}
&\definedas \eventthat[1]{
	\ictime[2]{}{\vect z+\centerpt{k}}{\gball{\vect z}{\dirsymb}{k}}
	\le \beta \thp \mu \parens{\centerpt{k}}},\\
\gballchainevent{\vect z}{\dirsymb}{m}{n}{\thp}{\beta}
&\definedas \eventthat[1]{
	\ictime[2]{}{\vect z+\centerpt{m}}{\gballchain{\vect z}{\dirsymb}{m}{n}}
	\le \beta \dist{\mu}(\centerpt{m}, \centerpt{n+1})}.
\end{align*}
Then $\gballchainevent{\vect z}{\dirsymb}{m}{n}{\thp}{\beta}$ occurs on the intersection $\bigcap_{k=m}^n \gballevent{\vect z}{\dirsymb}{k}{\thp}{\beta}$.
\end{lem}

\begin{proof}
Use Lemma~\ref{chaining_lem} (chaining) and Lemma~\ref{geometric_balls_lem}. Note that the statement remains true for any $\beta\ge 0$, but the only relevant case will be $\beta>1$.
\end{proof}

Observe that each of the events $\gballevent{\vect z}{\dirsymb}{k}{\thp}{\beta}$ in Lemma~\ref{ball_chaining_lem} is precisely the type of event we dealt with in Lemma~\ref{ball_covering_lem}. In particular, combining Lemma~\ref{ball_covering_lem} with Lemma~\ref{ball_chaining_lem} shows that for any $\beta>1$, if $m$ is sufficiently large and $n\ge m$, then the entire chain $\gballchain{\vect z}{\dirsymb}{m}{n}$ is likely to be covered from the point $\vect z+ \centerpt{m}$ by time $\beta \distnew{\mu}{\centerpt{m}} {\centerpt{n+1}}$. This corresponds to step~1 in the ``two-stage covering" strategy from Section~\ref{cone_growth:lattice_path_sec}.
However, rather than directly invoking Lemma~\ref{ball_covering_lem} in the proof of Lemma~\ref{thin_cone_covering_lem}, we will use Lemma~\ref{ball_covering_lem} to prove Lemma~\ref{space_covering_lem} below, which will be more convenient for the main proof.
First we prove the following lemma for dealing with the relevant ``short" paths in step~2 of the ``two-stage covering."

%

\begin{lem}[Traversal of radial lattice paths in a $\mu$-ball]
\label{path_tail_bound_lem}
\newcommand{\thisevent}{\symbolref{\eventsymb}{path_tail_bound_lem}(r)}
\newcommand{\thiskonst}{\konstref{path_tail_bound_lem}}

Let $\thiskonst = \konstref{segment_meas_lem}^{\mu} \parens[2]{\E \tmeasureof{\edge}}$, and for each $r\ge 0$, define the event
\[
\thisevent = \bigcap_{\vect z\in \cubify{\vect 0}} \bigcap_{t\ge r}
	\eventthat[1]{\tmeasure \parens[2]{\segment{\vect z}{\vect z+ \vect x}}
	\le \konstref{path_tail_bound_lem} t
	\text{ for all } \vect x\in \reacheddset{\mu}{\vect 0}{t} }.
\]
Then there exist (universal) positive constants $\bigcref{path_tail_bound_lem}$ and $\smallcref{path_tail_bound_lem}$ such that for all $r\ge 0$,
\[
\Pr \parens[2]{\thisevent} \ge 
1-\bigcref{path_tail_bound_lem}e^{-\smallcref{path_tail_bound_lem} r}.
\]
\end{lem}

\begin{proof}
\newcommand{\thisevent}{\eventref{path_tail_bound_lem}(r)}
\newcommand{\oldevent}[1]{\eventref{distant_path_lem} (r;#1)}

\newcommand{\thiscube}{\cubify{\vect 0}}
\newcommand{\centerpoint}{\vect z}
\newcommand{\ballpt}{\vect x}
\newcommand{\bdpt}{\vect y_{\ballpt}}

\newcommand{\thissegment}{\segment{\centerpoint}{\centerpoint+\ballpt}}
\newcommand{\longsegment}{\segment{\centerpoint}{\centerpoint +\bdpt}}


This can be proved directly, using an argument similar to but simpler than the one in the proof of Lemma~\ref{distant_path_lem}. Alternatively, we can simply observe that the event $\thisevent$ occurs on the event $\oldevent{1}$, as follows. Given any $\centerpoint\in\thiscube$ and $\ballpt\in \reacheddset{\mu}{\vect 0}{t} \setminus \setof{\vect 0}$, if we set $\bdpt \definedas \ballpt/\mu(\ballpt)$, then $\thissegment \subseteq \longsegment$, so
\[
\tmeasureof[2]{\thissegment} \le \tmeasureof{\longsegment}.
\]
Noting that $\bdpt\in \bd \reacheddset{\mu}{\vect 0}{t}$ and $\longsegment = \segment{\centerpoint+ (1-1) \bdpt}{\centerpoint +\bdpt}$, this shows that we get the desired bound for all $\ballpt\in \reacheddset{\mu}{\vect 0}{t}$ by setting $\epsilon=1$ in Lemma~\ref{distant_path_lem}, so we can take $\bigcref{path_tail_bound_lem} \definedas \bigcref{distant_path_lem}(1)$ and $\smallcref{path_tail_bound_lem} \definedas \smallcref{distant_path_lem}(1)$.
\end{proof}

\begin{lem}[Internal covering time of large $\mu$-balls far from the origin]
\label{space_covering_lem}
\newcommand{\thisevent}{\symbolref{\eventsymb}{space_covering_lem}(r, \thp; \beta)}
For any $\beta>1$, $\thp\in (0,1]$, and $r\ge 0$, define the event
\[
\thisevent = \bigcap_{\vect z\in \cubify{\vect 0}}
	\bigcap_{\substack{\vect x\in \Rd \\ \mu(\vect x)\ge r}}
	\eventthat[1]{\crreachedset{\tmeasure}{\vect z+\vect x}{\reacheddset{\mu}{\vect z+\vect x}{t}}{\beta t}
	= \reacheddset{\mu}{\vect z+\vect x}{t} \;\forall t\ge \thp \mu(\vect x)}.
\]
Then for any $\beta>1$, there exist positive constants $\bigcref{space_covering_lem}$ and $\smallcref{space_covering_lem}$ such that for any $\thp\in (0,1]$ and $r\ge 0$,
\[
\Pr \parens[2]{\thisevent} \ge
1- \frac{1}{\thp^d}\bigcref{space_covering_lem}(\beta)
e^{-\smallcref{space_covering_lem}(\beta) \thp r}.
\]
\end{lem}

\begin{proof}
\newcommand{\thisevent}{\eventref{space_covering_lem}(r, \thp; \beta)}
\newcommand{\thisC}{\bigcref{space_covering_lem}\parens{\beta}}
\newcommand{\thisc}{\smallcref{space_covering_lem}\parens{\beta}}
\newcommand{\ballevent}[4][0]{\eventref{ball_covering_lem}^{#2} \parens[#1]{#3;#4}}

\newcommand{\Cbeta}{\bigcsymb(\beta)}
\newcommand{\cbeta}{\smallcsymb(\beta)}

\newcommand{\thisdil}{\kappa}


Let $\thisdil = \dilnorm{\lspace{\infty}{d}}{\mu}$. Note for any point $\vect y\in\Rd$, there is some $\vect v\in\Zd$ such that $\vect y\in\cubify{\vect v}$, and since $\sup_{\vect z\in \cubify{\vect 0}} \mu(\vect z) = \thisdil/2$, it follows that $\dist{\mu}(\vect y,\vect v)\le \thisdil/2$. Therefore, the triangle inequality implies that if $\vect y\in\Rd$ and $\vect z\in \cubify{\vect 0}$, then
\begin{equation}\label{space_covering:yz_bounds}
\mu(\vect y -\vect z) \ge  \mu(\vect y) -\thisdil/2,
\quad\text{and}\quad
\mu(\vect y -\vect z)\ge r \implies \mu(\vect y)\ge r-\thisdil/2,
\end{equation}
and furthermore,
\begin{equation}\label{space_covering:yv_bounds}
\mu(\vect y)\ge r-\thisdil/2 \implies
\vect y\in \cubify{\vect v} \text{ for some $\vect v\in\Zd$ with }
	\mu(\vect v)\ge \mu(\vect y) -\thisdil/2 \ge r-\thisdil.
\end{equation}

Therefore, using \eqref{space_covering:yz_bounds} and \eqref{space_covering:yv_bounds} in the third and fourth lines, respectively, we have
\begin{align*}
\thisevent^\complement
&= \bigcup_{\vect z\in \cubify{\vect 0}}
	\bigcup_{\substack{\vect x\in \Rd \\ \mu(\vect x)\ge r}}
	\eventthat[1]{\ictime[2]{}{\vect z+\vect x}{\reacheddset{\mu}{\vect z+\vect x}{t}} >\beta t
	\text{ for some } t\ge \thp \mu(\vect x)}\\
&= \bigcup_{\vect z\in \cubify{\vect 0}}
	\bigcup_{\substack{\vect y\in \Rd \\ \mu(\vect y-\vect z)\ge r}}
	\eventthat[1]{\ictime[2]{}{\vect y}{\reacheddset{\mu}{\vect y}{t}} >\beta t
	\text{ for some } t\ge \thp \mu(\vect y-\vect z)}\\
&\subseteq \bigcup_{\substack{\vect y\in \Rd \\ \mu(\vect y)\ge r-\thisdil/2}}
	\eventthat[1]{\ictime[2]{}{\vect y}{\reacheddset{\mu}{\vect y}{t}} >\beta t
	\text{ for some } t\ge \thp \parens[2]{\mu(\vect y)-\thisdil/2}}\\
&\subseteq \bigcup_{\substack{\vect v\in \Zd \\ \mu(\vect v)\ge r-\thisdil}}
	\bigcup_{\vect y\in \cubify{\vect v}}
	\eventthat[1]{\ictime[2]{}{\vect v}{\reacheddset{\mu}{\vect y}{t}} >\beta t
	\text{ for some } t\ge \thp \parens[2]{\mu(\vect v)-\thisdil}}\\
&= \bigcup_{\substack{\vect v\in \Zd \\ \mu(\vect v)\ge r-\thisdil}}
	\ballevent[2]{\vect v}{\thp(\mu(\vect v)-\thisdil)}{\beta-1}^\complement,
\end{align*}
where $\eventref{ball_covering_lem}^{\vect v}(r_0;\epsilon)$ is the event in Lemma~\ref{ball_covering_lem}. Thus, setting $\Cbeta= \bigcref{ball_covering_lem}(\beta-1)$ and $\cbeta=\smallcref{ball_covering_lem}(\beta-1)$, Lemma~\ref{ball_covering_lem} and Lemmas~\ref{norm_shell_lem} and \ref{poly_geo_series_lem} imply that for all $r\ge 0$,
\begin{align*}
\Pr \parens[1]{\thisevent^\complement}
&\le \sum_{\substack{\vect v\in \Zd \\ \mu(\vect v)\ge r-\thisdil}}
	\Cbeta e^{-\cbeta \thp (\mu(\vect v)-\thisdil)}
	&& \text{(Lem.~\ref{ball_covering_lem})}\\
&\le \sum_{j= \floor{r-\thisdil}\join 0}^\infty
	\sum_{\substack{\vect v\in \Zd \\ j\le \mu(\vect v)< j+1}}
	\Cbeta e^{-\cbeta \thp (\mu(\vect v)-\thisdil)}
	\\
&\le \sum_{j= \floor{r-\thisdil}\join 0}^\infty
	\mushellbound \cdot (j+1)^{d-1}\cdot
	\Cbeta e^{-\cbeta \thp (j-\thisdil)}
	&& \text{(Lem.~\ref{norm_shell_lem})}\\
&\le \mushellbound \Cbeta e^{\cbeta \thp \thisdil}
	\cdot \frac{4(2d-2)^{d-1}}{\parens[2]{[\cbeta \thp] \meet 1}^d}
	\cdot e^{-\frac{\cbeta}{2} \thp \parens{\floor{r-\thisdil}\join 0}}
	&& \text{(Lem.~\ref{poly_geo_series_lem})}\\
&\le \mushellbound \Cbeta e^{\cbeta \thp \thisdil}
	\cdot \frac{4(2d-2)^{d-1}}{\parens[2]{\cbeta \meet 1}^d}
	\cdot \frac{1}{\thp^d}
	\cdot e^{-\frac{\cbeta}{2} \thp (r-\thisdil -1)}
	&& \text{($\because \thp\le 1$)}\\
&= \frac{4\mushellbound \Cbeta \cdot (2d-2)^{d-1}\cdot e^{\frac{\cbeta}{2} \thp (3\thisdil+1)}}
	{\parens[2]{\cbeta \meet 1}^d}
	\cdot \frac{1}{\thp^d} \cdot e^{-\frac{\cbeta}{2} \thp r}\\
&\le \frac{1}{\thp^d} \thisC e^{-\thisc \thp r},
	&& \text{($\because \thp\le 1$)}
\end{align*}
where (replacing $\thp$ with 1 in the exponent out front)
\[
\thisC \definedas
\frac{4\mushellbound \Cbeta \cdot (2d-2)^{d-1}\cdot e^{\frac{\cbeta}{2} (3\thisdil+1)}}
	{\parens[2]{\cbeta \meet 1}^d},
\quad\text{and}\quad
\thisc \definedas \frac{\cbeta}{2}.
\qedhere
\]
\end{proof}

\begin{thmremark}
\label{space_covering:fixed_dir_rem}
If the proof of Lemma~\ref{space_covering_lem} is carried out for a fixed direction rather than for all directions simultaneously (i.e.\ in the definition of $\eventref{space_covering_lem}(r, \thp; \beta)$, we replace the intersection over ``$\vect x\in \Rd$, $\mu(\vect x)\ge r$" by the intersection over ``$\vect x\in \dirvec$, $\mu(\vect x)\ge r$" for some fixed $\dirsymb\in \sphere{d-1}{\mu}$), then the resulting probability bound will be improved to
\[
\Pr\ge 1-\frac{1}{\thp} \bigcsymb(\beta) e^{-\smallcsymb(\beta) \thp r},
\]
for some positive constants $\bigcsymb(\beta)$ and $\smallcsymb(\beta)$. The factor $\thp^{d-1}$ will be gained in the estimate of the number of lattice points between levels $j$ and $j+1$ (constant vs.\ $(j+1)^{d-1}$), which affects the sum of the resulting series via Lemma~\ref{poly_geo_series_lem}.
\end{thmremark}

\subsection{Proof of Large Deviations Estimate for Covering a $\mu$-Cone Segment}

Now we are ready to prove the main result of the present section, bounding the internal covering time of a $\mu$-cone segment. In Lemma~\ref{thin_cone_covering_lem}, $\thp$ is the thickness of the $\mu$-cone segment, while $\spp$ is an independent ``speed parameter," specifying the accuracy with which we want to bound the covering time.
The quantity $\thp\join \spp$ appears in the covering time bound, which indicates that the bound in Lemma~\ref{thin_cone_covering_lem} only gives a good estimate for ``thin" cone segments. That is, if we want to use Lemma~\ref{thin_cone_covering_lem} to show that a cone segment of height $h$ is likely to be covered within time $\approx (1+\spp) h$, then the thickness $\thp$ of the cone segment needs to be on the same order as $\spp$ or smaller. 
However, if we want a good bound in a thicker cone segment, we can simply decompose the larger segment into a union of sufficiently thin segments using Lemma~\ref{mu_cone_decomp_lem}; this is what we will do in the proof of Theorem~\ref{cone_segment_growth_thm} in the next section.



The main tools for the proof of Lemma~\ref{thin_cone_covering_lem} are Lemmas~\ref{path_tail_bound_lem} and \ref{space_covering_lem}. Because of the uniform nature of the events in these lemmas, we obtain a similar uniform event in Lemma~\ref{thin_cone_covering_lem} essentially ``for free." Namely, the bound works simultaneously for all $\mu$-cone segments with apex in a fixed unit cube, axis in any direction, and sufficiently large thickness and height.
To simplify notation for events that are simultaneous over the various parameters, we will use the following shorthand for fixed $\vect v\in \Zd$, $\thp_0\in (0,1)$, and $h_0 >0$:
\begin{equation}
\label{thin_cone:param_restrictions_eqn}
\newcommand{\minth}{\thp_0}
\forall (\vect z, \dirsymb, \thp, h)
\quad\text{means}\quad
``\forall \vect z\in \cubify{\vect v},\ \forall \dirsymb\in \sphere{d-1}{\mu},\
\forall \thp \in [\minth,1],\ \forall h \ge h_0". 
\end{equation}

\begin{lem}[Internal covering time for a $\delta$-thin $\mu$-cone segment]
\label{thin_cone_covering_lem}
\newcommand{\minth}{\thp_0} 
\newcommand{\thisvertex}{\vect v}

\newcommand{\thisevent}{\eventref{thin_cone_covering_lem}^{\thisvertex}(h_0,\minth;\spp)}
\newcommand{\parameters}{\vect z, \dirsymb, \thp, h}
\newcommand{\thiskonst}{\konstref{thin_cone_covering_lem}}

Let $K_{\ref{thin_cone_covering_lem}}= 2\konstref{path_tail_bound_lem}+4$. For $\thisvertex\in \Zd$, $h_0 >0$, and $\minth,\spp\in (0,1)$, define the event
\[
\thisevent = \eventthat[1]{\forall (\parameters),\;
	\ictime[2]{}{\thisvertex}{\cone{\mu}{\thp}{\vect z}{\dirvec}(h)}
	\le \brackets[2]{1+\thiskonst (\thp\join\spp)} h},
\]
where the notation $\forall (\parameters)$ is defined in \eqref{thin_cone:param_restrictions_eqn}.
Then for any $\spp\in (0,1)$, there exist positive constants $\bigcref{thin_cone_covering_lem}$ and $\smallcref{thin_cone_covering_lem}$ such that for any $\thisvertex\in \Zd$, $\minth\in (0,1)$, and $h_0>0$,
\[
\Pr \parens[2]{\thisevent} 
\ge 1- \frac{1}{(\minth)^d} \bigcref{thin_cone_covering_lem}(\spp)
e^{-\smallcref{thin_cone_covering_lem}(\spp) \minth h_0}.
\]
\end{lem}

\begin{thmremark}
\label{thin_cone_covering:fixed_dir_rem}
The form of the probability bound in Lemma~\ref{thin_cone_covering_lem} is inherited directly from the bound in Lemma~\ref{space_covering_lem}. Thus, by Remark~\ref{space_covering:fixed_dir_rem}, if we bounded the covering time of $\mu$-cone segments in a fixed direction rather than simultaneously considering all directions, the probability bound in Lemma~\ref{thin_cone_covering_lem} would be improved by a factor of $(\thp_0)^{d-1}$.
\end{thmremark}

\begin{thmremark}
\label{thin_cone_covering:notation_rem}
Recall from \eqref{cone_fam_def_eqn} in Section~\ref{deterministic:mu_cones_sec} that we denoted the set of all initial $\mu$-cone segments with apex $\vect z\in\Rd$ and thickness $\thp\ge 0$ by
\begin{equation}
\label{cone_fam_def_reminder_eqn}
\conefamily{\mu}{\vect z,\thp} \definedas
\setbuilderbar[1]{\conetip{\mu}{\thp}{\vect z}{\dirvec}{h}}
	{\dirvec\in \dirset{\Rd},\, h\in [0,\infty]}.
\end{equation}
If we let $\conefamily{\mu}{\vect z,\thp, h}$ denote the set of $\mu$-cone segments in $\conefamily{\mu}{\vect z,\thp}$ with height $h$, i.e.
\begin{equation}
\label{cone_fam_h_def_eqn}
\conefamily{\mu}{\vect z,\thp, h} \definedas
\setbuilderbar[1]{\conetip{\mu}{\thp}{\vect z}{\dirvec}{h}}
	{\dirvec\in \dirset{\Rd}},
\end{equation}
then the event in Lemma~\ref{thin_cone_covering_lem} can also be written
\begin{multline*}
\eventref{thin_cone_covering_lem}^{\vect v} \parens{h_0,\thp_0;\spp}\\
= \bigcap_{\vect z\in\cubify{\vect v}} \bigcap_{\thp\in [\thp_0,1]} \bigcap_{h\ge h_0}
	\eventthat[1]{\ictime{}{\vect v}{\cones}
	\le \brackets{1+\konstref{thin_cone_covering_lem}(\thp\join\spp)} h
	\quad \forall \cones\in \conefamily{\mu}{\vect z,\thp,h}}.
\end{multline*}
We will use the notation \eqref{cone_fam_h_def_eqn} in the proof of Theorem~\ref{cone_segment_growth_thm} in the next section.
\end{thmremark}

\begin{proof}[Proof of Lemma~\ref{thin_cone_covering_lem}]

\newcommand{\minth}{\thp_0} 
\newcommand{\thisvertex}{\vect 0}
\newcommand{\genvertex}{\vect v}

\newcommand{\thisevent}{\eventref{thin_cone_covering_lem}^{\thisvertex}(h_0,\minth;\spp)}
\newcommand{\parameters}{\vect z, \dirsymb, \thp, h}

\newcommand{\mh}{M_h}
\newcommand{\nh}{N_h}

\newcommand{\thiscone}{\conetip{\mu}{\thp}{\vect z}{\dirvec}{h}}
\newcommand{\thisconetip}{\conetip{\mu}{\thp}{\vect z}{\dirvec}{\spp h}}
\newcommand{\zeroconetip}{\conetip{\mu}{\thp}{\vect 0}{\dirvec}{\spp h}}
\newcommand{\thisconesegment}{\conesegment{\mu}{\thp}{\vect z}{\dirvec}{\spp h}{h}}
\newcommand{\thisballchain}{\gballchain{\vect z}{\dirsymb}{\mh}{\nh}}
\newcommand{\thissegment}{\vect z+\segment{\spp h}{h}\dirsymb}
\newcommand{\zerosegment}{\segment{\spp h}{h}\dirsymb}

\newcommand{\conetipevent}{E(h_0,\minth; \spp)}
\newcommand{\ballchainevent}{F(h_0,\minth; \spp)}
\newcommand{\finalevent}{G(h_0,\minth; \spp)}
\newcommand{\intersectionevent}{H(h_0,\minth; \spp)}

\newcommand{\lineevent}[1][0]{\eventref{path_tail_bound_lem}\parens[#1]{\spp h_0}}
\newcommand{\ballevent}[1][0]{\eventref{space_covering_lem}%
	\parens[#1]{\spp h_0,\, \minth\, ;\, 1+\spp}}

\newcommand{\oldkonst}{\konstref{path_tail_bound_lem}}
\newcommand{\thiskonst}{\konstref{thin_cone_covering_lem}}


By translation invariance it suffices to assume $\genvertex = \thisvertex$. Fix $\spp,\minth\in (0,1)$ and $h_0>0$. For any $\vect z\in \cubify{\vect 0}$, $\dirsymb\in \sphere{d-1}{\mu}$, and $\thp\in [\minth,1]$, let the sequences $\setof{\centerpt{k}(\thp)}_{k\in\Z}$, $\setof{\gball{\vect z}{\dirsymb}{k}(\thp)}_{k\in\Z}$, and $\setof{\gballchain{\vect z}{\dirsymb}{m}{n}(\thp)}_{m,n\in\Z}$ be defined as in Lemma~\ref{geometric_balls_lem}. For the rest of the proof, we will omit the $\thp$ from this notation; that is, whenever we write $\centerpt{k}$, $\gball{\vect z}{\dirsymb}{k}$, or $\gballchain{\vect z}{\dirsymb}{m}{n}$, the argument is always implicitly taken to be $\thp$. For any $\thp\in [\minth,1]$ and $h \ge h_0$, let $\mh = \mh(\thp)$ and $\nh = \nh(\thp)$ be the unique integers satisfying
\begin{equation}\label{thin_cone:mh_nh_def_eqn}
(1+\thp)^{\mh-1}<\spp h \le (1+\thp)^{\mh}
\quad\text{and}\quad
(1+\thp)^{\nh} \le h < (1+\thp)^{\nh+1}.
\end{equation}
Again, we will omit the $\thp$ from the notation; any subsequent appearances of $\mh$ and $\nh$ are interpreted to have argument $\thp$. Observe that with these definitions, Lemma~\ref{geometric_balls_lem} implies that the chain of $\mu$-balls $\thisballchain$ is a $\thp h$-permeating subset of the cone sub-segment $\thisconesegment$, hence a $(\thp+\spp) h$-permeating subset of the whole cone segment $\thiscone$.

When writing events below, we will use the notation ``$\forall(\parameters)$" defined in \eqref{thin_cone:param_restrictions_eqn}, with $\genvertex=\thisvertex$.
We may omit some parameters when using this notation; for example, $\forall(\vect z,\thp)$ means $\forall \vect z\in \cubify{\vect 0}$ and $\forall \thp \in [\minth,1]$. Our goal is to bound the internal covering time of $\thiscone$ from $\vect 0$, simultaneously over the parameters $\vect z, \dirsymb, \thp, h$. It will be easiest to think about this covering event for fixed parameters --- for example, fix some $\vect z$ and $\dirsymb$, and assume $h=h_0$ and $\spp = \thp = \minth$. The uniformity over all parameter values will be a byproduct
of the uniform nature of the events in Lemmas~\ref{path_tail_bound_lem} and \ref{space_covering_lem}, which will be our primary technical tools.

For fixed parameter values, we bound the internal covering time of $\thiscone$ by  using the chaining lemma to break the covering event into three steps: First the process takes over the initial cone segment $\thisconetip$ up to height $\spp h$; second, the process takes over a chain of $\mu$-balls, namely $\thisballchain$, starting at height $\approx \spp h$ and ending near the end of the cone at height $\approx h$; finally, the process spreads from the $\thp h$-permeating subset $\thisballchain$ to the remainder of the cone, $\thisconesegment$. The precise statement is the following.

\begin{claim}\label{thin_cone:chaining_claim}
Lemma~\ref{chaining_lem} (chaining) implies that for any realization of $\tmeasure$, $\forall (\parameters)$,
\begin{align*}
\ictime[2]{}{\vect 0}{\,\thiscone}
&\le \ictime[2]{}{\vect 0}{\,\thisconetip} \notag\\
&+ \ictime[2]{}{\vect z+ \centerpt{\mh}}{\,\thisballchain} \notag\\
&+\ictime[2]{}{\thissegment}{\,\thisconesegment}.
\end{align*}
\end{claim}

\begin{proof}[Proof of Claim~\ref{thin_cone:chaining_claim}]
It is left to the reader to verify the hypotheses of Lemma~\ref{chaining_lem}, using Lemma~\ref{geometric_balls_lem} together with the definitions of $\mh$ and $\nh$ in \eqref{thin_cone:mh_nh_def_eqn}.
\end{proof}

We now define three events corresponding to the three terms in the bound from Claim~\ref{thin_cone:chaining_claim}. 
Let
\begin{align*}
\conetipevent &= \eventthat[3]{\forall (\parameters),\;
	\ictime[2]{}{\vect 0}{\,\thisconetip}
	\le \oldkonst (1+\thp)\spp h},\\
\ballchainevent &= \eventthat[3]{\forall (\parameters),\;
	\ictime[2]{}{\vect z+ \centerpt{\mh}}{\,\thisballchain}
	\le (1+\spp) \dist{\mu} \parens[2]{\centerpt{\mh},\centerpt{\nh+1}} },\\
\finalevent &= \eventthat[3]{\forall (\parameters),\;
	\ictime[2]{}{\thissegment}{\,\thisconesegment}
	\le (1+\spp) \thp h},
\end{align*}
where the dependence on $h_0$ and $\minth$ comes from the restrictions on the values of the parameters $h$ and $\thp$ in \eqref{thin_cone:param_restrictions_eqn}. Consider the intersection of these three events,
\[
\intersectionevent \definedas \conetipevent \cap \ballchainevent \cap \finalevent.
\]

\begin{claim}\label{thin_cone:total_event_claim}
The event $\thisevent$ occurs on the intersection $\intersectionevent$.
\end{claim}

\begin{proof}[Proof of Claim~\ref{thin_cone:total_event_claim}]
First note that the upper inequality for $\mh$ and the lower inequality for $\nh$ in \eqref{thin_cone:mh_nh_def_eqn}, together with the definition of $\setof{\centerpt{k}}_{k\in\Z}$, imply that $\forall (\dirsymb,\thp,h)$,
\[
\dist{\mu} \parens[2]{\centerpt{\mh},\centerpt{\nh+1}}
= (1+\thp)^{\nh+1} - (1+\thp)^{\mh}
\le (1+\thp)h - \spp h.
\]
Using this bound, the chaining inequality in Claim~\ref{thin_cone:chaining_claim} implies that on the intersection event $\intersectionevent$, it holds that $\forall (\parameters)$,
\begin{align}\label{thin_cone:total_covering_bound_eqn}
\ictime[2]{}{\vect 0}{\thiscone}
&\le \oldkonst (1+\thp)\spp h
	+ (1+\spp)\cdot \dist{\mu} \parens[2]{\centerpt{\mh},\centerpt{\nh+1}}
	+(1+\spp) \thp h \notag\\
&\le \oldkonst (1+\thp)\spp h
	+ (1+\spp) (1+\thp-\spp) h
	+(1+\spp) \thp h \notag\\
&= \brackets[1]{1+ 2\thp
	+ \parens[2]{\oldkonst (1+\thp)+ 2\thp -\spp}\spp } h \notag\\
&\le \brackets[2]{1 + 2\thp + \parens{2\oldkonst +2} \spp} h,
\end{align}
where the final line in \eqref{thin_cone:total_covering_bound_eqn} follows from the assumption that $\thp\le 1$. Now,
\[
2\thp + \parens{2\oldkonst +2} \spp
\le
\parens{2\oldkonst +4} (\thp\join\spp)
= \thiskonst (\thp\join\spp),
\]
so \eqref{thin_cone:total_covering_bound_eqn} implies that
\begin{align*}
\intersectionevent
&\subseteq \eventthat[1]{\forall (\parameters),\;
	\ictime[2]{}{\vect 0}{\thiscone}
	\le \brackets[2]{1+\thiskonst (\thp\join\spp)} h}\\
&= \thisevent. \qedhere
\end{align*}
\end{proof}

To finish the proof, we need a lower bound on the probability of the event $\intersectionevent$. This will be provided by Lemmas~\ref{path_tail_bound_lem} and \ref{space_covering_lem}, together with the following claim.


\begin{claim}\label{thin_cone:event_containment_claim}
The event $\conetipevent$ occurs on the event $\lineevent[2]$, and the events $\ballchainevent$ and $\finalevent$ both occur on the event $\ballevent[2]$.
\end{claim}

Assume for the moment that Claim~\ref{thin_cone:event_containment_claim} is true. Then Claims~\ref{thin_cone:total_event_claim} and \ref{thin_cone:event_containment_claim}, together with Lemmas~\ref{path_tail_bound_lem} and \ref{space_covering_lem}, imply that
\begin{align*}
\Pr \parens[2]{\thisevent}
&\ge \Pr \parens[2]{\intersectionevent}\\
&\ge \Pr \brackets[1]{\lineevent[1] \cap \ballevent[1]}\\
&\ge 1- \bigcref{path_tail_bound_lem}
	e^{-\smallcref{path_tail_bound_lem}\cdot \spp h_0}
	- \frac{\bigcref{space_covering_lem}(1+\spp)}{(\minth)^{d}}
	e^{-\smallcref{space_covering_lem}(1+\spp)
		\cdot \minth \cdot \spp h_0}\\
&\ge 1- (\minth)^{-d}\bigcref{thin_cone_covering_lem}(\spp)
	e^{-\smallcref{thin_cone_covering_lem}(\spp) \minth h_0},
\end{align*}
where
\[
\bigcref{thin_cone_covering_lem}(\spp)
=\bigcref{path_tail_bound_lem} + \bigcref{space_covering_lem}(1+\spp)
\quad\text{and}\quad
\smallcref{thin_cone_covering_lem}(\spp)
= \parens[2]{\smallcref{path_tail_bound_lem}
	\meet \smallcref{space_covering_lem}(1+\spp)}
	\cdot \spp.
\]
To complete the proof, it remains only to prove Claim~\ref{thin_cone:event_containment_claim}.

\begin{proof}[Proof of Claim~\ref{thin_cone:event_containment_claim}]

\renewcommand{\proofpart}[1]{\vskip 2mm \noindent \fbox{{Proof that #1:}} \vskip 2mm}

\newcommand{\centralball}[1][0]{\reacheddset[#1]{\mu}{\vect 0}{(1+\thp) \spp h}}
\newcommand{\hmin}{h_0}
\newcommand{\kball}[1][0]{\reacheddset[#1]{\mu}{\vect z+\centerpt{k}}{\thp \mu(\centerpt{k})}}
\newcommand{\xball}[1][0]{\reacheddset[#1]{\mu}{\vect z+\vect x}{\thp \mu(\vect x)}}
\newcommand{\xballt}[1][0]{\reacheddset[#1]{\mu}{\vect z+\vect x}{t}}

We divide the proof into three parts corresponding to the three claimed containments.

\proofpart{$\conetipevent \supseteq \lineevent[2]$} First note that $\forall (\parameters)$,
\begin{equation}\label{thin_cone:cone_tip_eqn}
\thisconetip = \bigcup_{\vect x\in\zeroconetip} \segment{\vect z}{\vect z+\vect x},
\quad\text{and}\quad
\zeroconetip \subseteq \centralball[2].
\end{equation}
Therefore, the decomposition property and total traversal measure bound in Lemma~\ref{covering_time_properties_lem} imply that for any realization of $\tmeasure$, $\forall (\parameters)$,
\begin{align*}
\ictime[1]{}{\vect 0}{\,\thisconetip}
&= \sup_{\vect x\in\zeroconetip} 
	\ctime[2]{}{\thisconetip}{\vect 0}{\segment{\vect z}{\vect z+\vect x}}
	&&\text{(decomposition, \eqref{thin_cone:cone_tip_eqn})}\\
&\le \sup_{\vect x\in\zeroconetip}
	\brackets[1]{\rptime[0]{}{\thisconetip}{\vect 0}{\vect z}
	+ \tmeasure \parens[2]{\segment{\vect z}{\vect z+\vect x}}}
	&&\text{(total $\tmeasure$ bound)}\\
&\le \sup_{\vect x\in\centralball}  \tmeasure \parens[2]{\segment{\vect z}{\vect z+\vect x}},
\end{align*}
where the last line follows from \eqref{thin_cone:cone_tip_eqn} and the fact that $\rptime[0]{}{\thisconetip}{\vect 0}{\vect z} = 0$ for any $\vect z\in \cubify{\vect 0}$. Note that the final bound doesn't depend on $\dirsymb$. Therefore,
\begin{align*}
\conetipevent &= \eventthat[1]{\forall (\parameters),\;
	\ictime[2]{}{\vect 0}{\thisconetip}
	\le \oldkonst (1+\thp)\spp h}\\
&\supseteq \eventthat[1]{\forall (\vect z,\thp, h),\;
	\sup_{\vect x\in\centralball}  \tmeasure \parens[2]{\segment{\vect z}{\vect z+\vect x}}
	\le \oldkonst (1+\thp)\spp h}\\
&\supseteq \eventthat[1]{\forall (\vect z,\thp, h),\;
	\sup_{\vect x\in \reacheddset{\mu}{\vect 0}{t}} 
	\tmeasure \parens[2]{\segment{\vect z}{\vect z+\vect x}} \le  \oldkonst t
	\text{ for all } t\ge (1+\thp)\spp h}\\
&\supseteq \eventthat[1]{\forall \vect z\in \cubify{\vect 0},\;
	\sup_{\vect x\in \reacheddset{\mu}{\vect 0}{t}} 
	\tmeasure \parens[2]{\segment{\vect z}{\vect z+\vect x}} \le \oldkonst t
	\text{ for all } t\ge \spp\hmin }\\
&= \eventref{path_tail_bound_lem} \parens[2]{\spp\hmin}.
\end{align*}

\proofpart{$\ballchainevent\supseteq \ballevent[2]$} Let $\setof[2]{\gballevent{\vect z}{\dirsymb}{k}{\thp}{\beta}}_{\vect z,\dirsymb, k, \thp, \beta}$ and $\setof[2]{\gballchainevent{\vect z}{\dirsymb}{m}{n}{\thp}{\beta}}_{\vect z,\dirsymb, m,n, \thp, \beta}$ be the collections of events defined in Lemma~\ref{ball_chaining_lem}. Then by Lemma~\ref{ball_chaining_lem} and the definition of $\ballchainevent$, we have
\begin{equation}\label{thin_cone:ballchain_event_eqn}
\ballchainevent = \bigcap_{\vect z,\dirsymb, \thp,h}
	\gballchainevent{\vect z}{\dirsymb}{\mh}{\nh}{\thp}{1+\spp}
\supseteq \bigcap_{\parameters}\: \bigcap_{k=\mh}^{\nh}
	\gballevent{\vect z}{\dirsymb}{k}{\thp}{1+\spp},
\end{equation}
where, \textit{\`a la} \eqref{thin_cone:param_restrictions_eqn}, the intersection over $(\parameters)$ means the intersection over $\vect z\in \cubify{\vect 0}$, $\dirsymb\in \sphere{d-1}{\mu}$, $\thp \in [\minth,1]$, and $h \ge \hmin$. Recall that, since $\gball{\vect z}{\dirsymb}{k} = \kball[2]$,
\begin{equation}\label{thin_cone:gball_event_eqn}
\gballevent{\vect z}{\dirsymb}{k}{\thp}{1+\spp}
= \eventthat[1]{\ictime[3]{}{\vect z+\centerpt{k}}{\kball[2]}
	\le (1+\spp)\thp \mu(\centerpt{k}) }.
\end{equation}
Observe that for any fixed $\thp\in (0,1]$,
\begin{align}\label{thin_cone:index_union_eqn}
\bigcup_{\substack{\dirsymb\in \sphere{d-1}{\mu} \\ h\ge \hmin }}
	\bigcup_{k = \mh(\thp)}^\infty
	\setof[2]{\centerpt{k}(\thp)}
&\subseteq \bigcup_{\substack{\dirsymb\in \sphere{d-1}{\mu} \\ h\ge \hmin}}
	\setbuilder[3]{\vect x\in\Rd}
	{\mu(\vect x)\ge \mu \parens[2]{\centerpt{\mh}(\thp)}} \notag\\
&\subseteq \setbuilder[3]{\vect x\in\Rd}{\mu(\vect x)\ge \spp\hmin },
\end{align}
where the final inclusion follows from the fact that
\[
\mu \parens[2]{\centerpt{\mh}(\thp)} = (1+\thp)^{\mh(\thp)} \ge \spp h \ge \spp\hmin 
\quad \forall \dirsymb\in \sphere{d-1}{\mu},\; \forall \thp\in (0,1],\; \forall h\ge \hmin
\]
by the definitions of the sequence $\setof{\centerpt{k}}$ in Lemma~\ref{geometric_balls_lem} and of $\mh(\thp)$ in \eqref{thin_cone:mh_nh_def_eqn}. Therefore, using \eqref{thin_cone:ballchain_event_eqn}, \eqref{thin_cone:gball_event_eqn}, and \eqref{thin_cone:index_union_eqn} in lines one, two, and three, respectively, we get
\begin{align*}
\ballchainevent
&\supseteq \bigcap_{\parameters}\: \bigcap_{k=\mh}^{\nh}
	\gballevent{\vect z}{\dirsymb}{k}{\thp}{1+\spp}\\
&\supseteq \bigcap_{\parameters}\: \bigcap_{k=\mh}^{\infty}
	\eventthat[1]{\ictime[3]{}{\vect z+\centerpt{k}}{\kball[2]}
	\le (1+\spp)\thp \mu(\centerpt{k}) }\\
&\supseteq \bigcap_{\vect z,\thp}
	\bigcap_{\substack{\vect x\in\Rd \\ \mu(\vect x)\ge \spp \hmin}}
	\eventthat[1]{\ictime[3]{}{\vect z+\vect x}{\xball[2]}
	\le (1+\spp)\thp \mu(\vect x) }\\
&\supseteq \bigcap_{\vect z,\thp}
	\bigcap_{\substack{\vect x\in\Rd \\ \mu(\vect x)\ge \spp \hmin}}
	\eventthat[1]{\ictime[1]{}{\vect z+\vect x}{\xballt}
	\le (1+\spp) t \;\; \forall t\ge \thp \mu(\vect x) }\\
&\supseteq \bigcap_{\vect z\in\cubify{\vect 0}}
	\bigcap_{\substack{\vect x\in\Rd \\ \mu(\vect x)\ge \spp \hmin}}
	\eventthat[1]{ \ictime[1]{}{\vect z+\vect x}{\xballt}
	\le (1+\spp) t \;\; \forall t\ge \minth \mu(\vect x) }\\
&= \eventref{space_covering_lem}\parens[2]{\spp \hmin, \minth; 1+\spp}.
\end{align*}

\proofpart{$\finalevent\supseteq\ballevent[2]$} First note that $\forall (\parameters)$,
\[
\thisconesegment = \bigcup_{\vect x\in \zerosegment} \xball[2].
\]
Thus, the superdecomposition property in Lemma~\ref{covering_time_properties_lem} implies that for any realization of $\tmeasure$, $\forall (\parameters)$,
\begin{equation}\label{thin_cone:cone_segment_covering_bound}
\ictime[1]{}{\thissegment}{\,\thisconesegment}
\le \sup_{\vect x\in\zerosegment} \ictime[3]{}{\vect z+\vect x}{\,\xball[2]}.
\end{equation}
Now note that for any $\dirsymb\in\sphere{d-1}{\mu}$ and $h\ge \hmin$, 
\begin{equation}\label{thin_cone:segment_end_bounds}
\vect x\in \zerosegment \implies
\mu(\vect x)\le h \quad\text{and}\quad
\mu(\vect x)\ge \spp h\ge \spp\hmin.
\end{equation}
Therefore, using \eqref{thin_cone:cone_segment_covering_bound} and \eqref{thin_cone:segment_end_bounds} in lines two and three, respectively,
\begin{align*}
\finalevent
&= \eventthat[3]{\forall (\parameters),\;
	\ictime[2]{}{\thissegment}{\,\thisconesegment}
	\le (1+\spp) \thp h}\\
&\supseteq \eventthat[1]{\forall (\parameters),\;
	\sup_{\vect x\in\zerosegment} \ictime[3]{}{\vect z+\vect x}{\,\xball[2]}
	\le (1+\spp) \thp h}\\
&\supseteq \eventthat[1]{\forall (\vect z,\thp),\;
	\sup_{\vect x : \mu(\vect x)\ge \spp\hmin} 
	\ictime[3]{}{\vect z+\vect x}{\,\xball[2]}
	\le (1+\spp) \thp \mu(\vect x)}\\
&\supseteq \bigcap_{\substack{\vect z\in\cubify{\vect 0} \\ \thp \in [\minth,1]}} 
	\bigcap_{\substack{\vect x\in\Rd \\ \mu(\vect x)\ge \spp\hmin}}
	\eventthat[1]{
	\ictime[1]{}{\vect z+\vect x}{\,\xballt}
	\le (1+\spp)t \;\; \forall t\ge \thp \mu(\vect x)}\\
&\supseteq \bigcap_{\vect z\in\cubify{\vect 0}} 
	\bigcap_{\substack{\vect x\in\Rd \\ \mu(\vect x)\ge \spp\hmin}}
	\eventthat[1]{\ictime[1]{}{\vect z+\vect x}{\,\xballt}
	\le (1+\spp)t \;\; \forall t\ge \minth \mu(\vect x)}\\
&= \eventref{space_covering_lem} \parens[2]{\spp\hmin, \minth; 1+\spp}.\qedhere
\end{align*}
\end{proof}
This completes the proof of Lemma~\ref{thin_cone_covering_lem}.
\end{proof}

\section{Growth in Star Sets and Cones}
\label{cone_growth:cone_growth_sec}
In this section we use Lemmas~\ref{path_tail_bound_lem} and \ref{thin_cone_covering_lem} to prove Theorem~\ref{cone_segment_growth_thm} below, which provides a large deviations estimate for first-passage growth restricted to $\mu$-stars. Theorem~\ref{cone_segment_growth_thm} is the main result of Chapter~\ref{cone_growth_chap} and, together with Proposition~\ref{restricted_growth_ub_prop} below, gives a stochastic analogue of Proposition~\ref{determ_star_growth_prop} (which described the growth of the deterministic $\mu$-process in star sets).
Recall from \eqref{fat_stars_def_eqn} that for a norm $\mu$ on $\Rd$ and $\thp\ge 0$, we defined the set of \textdef{$(\mu,\thp)$-stars at $\vect z\in\Rd$} to be
\begin{equation}
\label{fat_stars_def_reminder_eqn}
\fatstars{\mu}{\thp}{\vect z}
\definedas \setof[2]{\text{Unions of elements in $\conefamily{\mu}{\vect z,\thp}$}}
= \setbuilderbar[2]{{\textstyle \bigcup}\, {\conesubfam}}
	{\conesubfam \subseteq \conefamily{\mu}{\vect z,\thp}},
\end{equation}
where $\conefamily{\mu}{\vect z,\thp}$ is the set of $\thp$-thick initial $\mu$-cone segments at $\vect z$ (cf.\ \eqref{cone_fam_def_reminder_eqn} in Remark~\ref{thin_cone_covering:notation_rem} above). Also recall, from Section~\ref{deterministic:convex_geom_defs_sec}, that the collection of all \textdef{star sets at $\vect z$} is
\begin{equation}
\label{star_set_def_reminder_eqn}
\stars{\vect z}\definedas
\setbuilderbar[1]{S\subseteq\Rd}{\segment{\vect z}{\vect x}\subseteq S \text{ for all } \vect x\in S}.
\end{equation}
We observed in Section~\ref{deterministic:norm_geometry_sec} that $\fatstars{\mu}{0}{\vect z}= \stars{\vect z}$, and if $\thp>0$, then $\fatstars{\mu}{\thp}{\vect z}\subsetneq \stars{\vect z}$.

Section~\ref{cone_growth:cone_growth_sec} is divided into two subsections. In Section~\ref{cone_growth:mu_stars_dev_sec} we prove the main large deviations estimates, and in Section~\ref{cone_growth:mu_stars_shape_thm_sec} we use these estimates to derive Shape Theorems for restricted growth.

\subsection{Large Deviations Estimates for Growth in $\mu$-Stars and Cones}
\label{cone_growth:mu_stars_dev_sec}



Before proving Theorem~\ref{cone_segment_growth_thm}, which provides a lower bound on growth in $\mu$-stars, we first derive the following upper bound on growth restricted to any subset of $\Rd$. This result follows trivially from Lemma~\ref{shifted_dev_bounds_lem}; the bound it gives for an arbitrary restricting set may not be very accurate, but it gives a useful bound for growth restricted to star sets, using Proposition~\ref{determ_star_growth_prop} for the deterministic process.

\begin{prop}[Upper bound on restricted growth]
\label{restricted_growth_ub_prop}
\newcommand{\thisset}{S}
\newcommand{\thiscenter}{\vect z}
\newcommand{\thisvertex}{\vect v}
\newcommand{\thiscube}{\cubify{\thisvertex}}
\newcommand{\thist}{t_0}

\newcommand{\theseparams}{^{\thisvertex} \parens{\thist; \spp}}
\newcommand{\thisC}{\bigcref{shifted_dev_bounds_lem}(\spp)}
\newcommand{\thisc}{\smallcref{shifted_dev_bounds_lem}(\spp)}

\newcommand{\setsevent}{\eventref{restricted_growth_ub_prop} \theseparams}
\newcommand{\starsevent}{\eventsymb_{\ref{restricted_growth_ub_prop}}^{\thisvertex ;\bigstar} \parens{\thist; \spp}}
\newcommand{\oldevent}{\eventref{shifted_dev_bounds_lem} \theseparams}


%

Suppose $\tmeasure$ is an \iid\ traversal measure on $\Zd$ satisfying \finitespeed{d} and \expmoment, and let $\mu$ be the shape function for $\tmeasure$ from Theorem~\ref{shape_thm}.
For any $\thisvertex\in\Zd$, $\thist\ge 0$, and $\spp>0$, the following two events occur on $\oldevent$ and hence occur with probability at least $1-\thisC e^{-\thisc \thist}$:
\begin{enumerate}
\item
$\displaystyle
\setsevent \definedas
\eventthat[3]{ \forall \thiscenter\in \thiscube \text{ and } \thisset \subseteq \Rd,\;
\rreachedset{\tmeasure}{\thisvertex}{\thisset}{t}
\subseteq \thisset\cap \reacheddset[2]{\mu}{\thiscenter}{(1+\spp) t-}
\text{ for all } t\ge \thist};
$

\item
$\displaystyle
\starsevent \definedas
\eventthat[3]{ \forall\thiscenter\in \thiscube \text{ and } \thisset\in \stars{\thiscenter},\;
\rreachedset{\tmeasure}{\thisvertex}{\thisset}{t}
\subseteq \rreacheddset[2]{\mu}{\thiscenter}{\thisset}{(1+\spp) t-}
\text{ for all } t\ge \thist}.
$

\end{enumerate}
\end{prop}

\begin{proof}
\newcommand{\thisset}{S}
\newcommand{\thiscenter}{\vect z}
\newcommand{\thisvertex}{\vect v}
\newcommand{\thiscube}{\cubify{\thisvertex}}
\newcommand{\thisC}{\bigcref{shifted_dev_bounds_lem}(\spp)}
\newcommand{\thisc}{\smallcref{shifted_dev_bounds_lem}(\spp)}

The event in Part 1 occurs on $\eventref{shifted_dev_bounds_lem}^{\thisvertex}(t_0;\spp)$ because for any $\thisset\subseteq\Rd$, we have $\rreachedset{\tmeasure}{\thisvertex}{\thisset}{t} \subseteq \thisset\cap \rreachedset{\tmeasure}{\thisvertex}{\Rd}{t}$ for all $t\ge 0$. Now observe that the event in Part 2 occurs on the event in Part 1, because if $S$ is star-shaped at $\thiscenter$, then $\rreacheddset[2]{\mu}{\thiscenter}{\thisset}{(1+\spp) t-} = \thisset\cap \reacheddset[2]{\mu}{\thiscenter}{(1+\spp) t-}$ by Proposition~\ref{determ_star_growth_prop}.
\end{proof}

We are now ready to prove the main result of Chapter~\ref{cone_growth_chap}.

\begin{thm}[Lower bound on growth restricted to $\thp$-thick $\mu$-stars]
\label{cone_segment_growth_thm}
\newcommand{\minth}{\thp}
\newcommand{\tmin}{t_0} 

\newcommand{\thisvertex}{\vect v}
\newcommand{\thisapex}{\vect z}
\newcommand{\thisunions}{\fatstars{\mu}{\minth}{\thisapex}}
\newcommand{\thisset}{S}

\newcommand{\thisevent}{\eventref{cone_segment_growth_thm}^{\thisvertex} (\tmin,\minth; \spp)}
\newcommand{\thisbigc}{\bigcref{cone_segment_growth_thm}}
\newcommand{\thissmallc}{\smallcref{cone_segment_growth_thm}}

Suppose $\tmeasure$ is an \iid\ traversal measure on $\Zd$ satisfying \finitespeed{d} and \expmoment, and let $\mu$ be the shape function for $\tmeasure$ from Theorem~\ref{shape_thm}. For any $\thisvertex\in \Zd$, $\tmin\ge 0$, $\minth\in (0,1]$, and $\spp\in (0,1)$, define the event
\[
\thisevent \definedas
\bigcap_{\thisapex\in \cubify{\thisvertex}}
\eventthat[3]{\forall \thisset\in \thisunions,\;
	\crreachedset{\tmeasure}{\thisvertex}{\thisset}{t}
	\supseteq \rreacheddset[2]{\mu}{\thisapex}{\thisset}{(1-\spp) t}
	\text{ for all } t\ge \tmin}.
\]
Then given $\spp\in (0,1)$, there exist positive constants $\thisbigc$ and $\thissmallc$ such that for any $\thisvertex\in \Zd$, $\tmin\ge 0$, and $\minth\in (0,1]$,
\[
\Pr \parens[2]{\thisevent} \ge
1- \frac{1}{\minth^d} \thisbigc(\spp) e^{-\thissmallc(\spp) \minth \tmin}.
\]
\end{thm}

\begin{proof}
\newcommand{\minth}{\thp}
\newcommand{\tmin}{t_0} 

\newcommand{\genvertex}{\vect v} 
\newcommand{\thisvertex}{\vect 0}
\newcommand{\thisapex}{\vect z}
\newcommand{\thiscube}{\cubify{\thisvertex}}
\newcommand{\thisset}{S}

\newcommand{\thisevent}{\eventref{cone_segment_growth_thm}^{\thisvertex} (\tmin,\minth; \spp)}
\newcommand{\thisbigc}{\bigcref{cone_segment_growth_thm}}
\newcommand{\thissmallc}{\smallcref{cone_segment_growth_thm}}

\newcommand{\linekonst}{\konstref{path_tail_bound_lem}}
\newcommand{\oldkonst}{\konstref{thin_cone_covering_lem}}
\newcommand{\thiskonst}{\konstsymb}

\newcommand{\sppk}{\spp_0}
\newcommand{\thpk}{\thp_0}
\newcommand{\hk}{h_0}

\newcommand{\conept}{\vect y} 
\newcommand{\centeredpt}{\vect x} 

\newcommand{\cutradius}{\frac{\tmin}{\linekonst}}

\newcommand{\nearevent}{E(\tmin,\minth;\spp)}
\newcommand{\farevent}{F(\tmin,\minth;\spp)}

\newcommand{\lineevent}{\eventref{path_tail_bound_lem} \parens[1]{\cutradius}}
\newcommand{\conesevent}{\eventref{thin_cone_covering_lem}^{\thisvertex} \parens{\hk, \thpk; \sppk}}

\newcommand{\thisunions}{\fatstars{\mu}{\minth}{\thisapex}}

\newcommand{\family}{\conesubfam}
\newcommand{\famcone}{\coneunion{\family}}

\newcommand{\ptsetcone}{\conesymb^{\conept}_{\thisset,\thpk}}
\newcommand{\ptseth}{h^{\conept}_{\thisset,\thpk}}


By translation invariance it suffices to assume $\genvertex = \thisvertex$. The idea of the proof is to divide the points in $\thisset$ into two groups, those ``near the origin," and those ``far away," and then to bound the covering times of the ``near points" using Lemma~\ref{path_tail_bound_lem} and to bound the covering times of the ``far points" using Lemma~\ref{thin_cone_covering_lem}. We get a simultaneous bound over all $\thisapex$ and $\thisset$ essentially ``for free" because of the uniform nature of these two lemmas. First we define two events
\begin{align*}
\nearevent &\definedas \bigcap_{\thisapex,\thisset}
	\eventthat[1]{
	\ctime{}{\thisset}{\thisvertex}{\conept} \le \tmin 
	\text{ for all $\conept\in \thisset$ with $\distnew{\mu}{\thisapex}{\conept} \le \cutradius$}},\\
\farevent &\definedas  \bigcap_{\thisapex,\thisset}
	\eventthat[1]{
	\ctime{}{\thisset}{\thisvertex}{\conept} \le \frac{\distnew{\mu}{\thisapex}{\conept}}{1-\spp}
	\text{ for all $\conept\in \thisset$ with $\distnew{\mu}{\thisapex}{\conept} > \cutradius$}},
\end{align*}
where the intersection over $(\thisapex,\thisset)$ means the intersection over $\thisapex\in\thiscube$ and $\thisset\in \thisunions$, and we have abbreviated $\ctmetric{\tmeasure}{\thisset}$ to $\ctmetric{}{\thisset}$. These events correspond to the ``near points" and ``far points" mentioned above, and the following claim shows why these definitions are relevant.

\begin{claim}
\label{cone_segment_growth:inversion_claim}
The event $\thisevent$ occurs on the intersection $\nearevent\cap \farevent$.
\end{claim}

\begin{proof}[Proof of Claim~\ref{cone_segment_growth:inversion_claim}]
Using the inversion formula \eqref{inversion_formulas:lbg} from Lemma~\ref{inversion_formula_lem}, with $A=\thisvertex$, $A'=\thisapex$, and $\alpha = 1-\spp$, for any $\thisapex\in\thiscube$ and $\thisset\in \thisunions$, we have
\begin{equation}
\label{cone_segment_growth:inversion_eqn}
\eventthat[3]{ \crreachedset{\tmeasure}{\thisvertex}{\thisset}{t}
	\supseteq \rreacheddset[2]{\mu}{\thisapex}{\thisset}{(1-\spp) t}
	\; \forall t\ge \tmin}
=\eventthat[1]{ \ctime{\tmeasure}{\thisset}{\thisvertex}{\conept} \le 
	\tmin \join \parens[1]{\frac{\distnew{\mu}{\thisapex}{\conept}}{1-\spp}}
	\; \forall \conept\in \thisset}.
\end{equation}
Note that we used the fact that $\thisset$ is star-shaped at $\thisapex$ to replace $\imetric{\mu}{S}$ with $\normmetric{\mu}$ by Lemma~\ref{segment_idist_lem}. Now observe that for any $M\ge 0$, the event on the right-hand side of \eqref{cone_segment_growth:inversion_eqn} occurs on the intersection
\begin{multline*}
\eventthat[1]{\ctime{\tmeasure}{\thisset}{\thisvertex}{\conept} \le \tmin\quad
	\forall\conept\in \thisset \text{ with } \distnew{\mu}{\thisapex}{\conept}\le M}\\
\hspace{3cm} \cap \eventthat[1]{\ctime{\tmeasure}{\thisset}{\thisvertex}{\conept}
	\le \frac{\distnew{\mu}{\thisapex}{\conept}}{1-\spp}\quad
	\forall\conept\in \thisset \text{ with } \distnew{\mu}{\thisapex}{\conept}> M}.
\end{multline*}
Thus, taking $M=\cutradius$ and intersecting over all $\thisapex\in\thiscube$ and $\thisset\in \thisunions$ proves the claim.
%
\end{proof}

The trick to bounding the covering time of points in $\thisset$, as well as to getting the bounds to hold simultaneously for all $\thisset$, rests on the following claim.

\begin{claim}
\label{cone_segment_growth:set_reductions_claim}
\renewcommand{\ptsetcone}{\conesymb^{\conept}_{\thisset,\thp'}}
\renewcommand{\ptseth}{h^{\conept}_{\thisset,\thp'}}

Let $\thisset\in \thisunions$.
\begin{enumerate}
\item \label{cone_segment_growth:tmeasure_bound}
For any $\conept\in \thisset$ and any realization of $\tmeasure$, we have $\ctime{\tmeasure}{\thisset}{\thisvertex}{\conept} \le \tmeasureof[2]{\segment{\thisapex}{\conept}}$. 

\item \label{cone_segment_growth:cone_exists}
For any $\thp'\in [0,\minth]$ and $\conept\in \thisset$, there exists $\ptseth\ge 0$ and $\ptsetcone\in \conefamily[2]{\mu}{\thisapex,\thp',\ptseth}$ such that
\[
\conept\in \ptsetcone\subseteq \thisset
\quad\text{and}\quad
(1-\thp')\ptseth \le \distnew{\mu}{\thisapex}{\conept} \le (1+\thp')\ptseth,
\]
where $\conefamily[2]{\mu}{\thisapex,\thp',\ptseth}$ is the family of initial $\mu$-cone segments at $\vect z$ defined by \eqref{cone_fam_h_def_eqn} in Remark~\ref{thin_cone_covering:notation_rem}.
\end{enumerate}
\end{claim}

\begin{proof}[Proof of Claim~\ref{cone_segment_growth:set_reductions_claim}]
\renewcommand{\ptsetcone}{\conesymb^{\conept}_{\thisset,\thp'}}
\renewcommand{\ptseth}{h^{\conept}_{\thisset,\thp'}}

\newcommand{\pthbig}{h_{\conept}}
\newcommand{\pth}{h_{\conept}'}
\newcommand{\ptdirsymb}{\dirsymb_{\conept}}
\newcommand{\ptdirvec}{\dirvec_{\conept}}

For Part~\ref{cone_segment_growth:tmeasure_bound}, note that since any $\thisset\in \thisunions$ is star-shaped at $\thisapex$, we have $\segment{\thisapex}{\conept}\subseteq \thisset$ for all $\conept\in\thisset$. Thus, using the total traversal measure bound from Lemma~\ref{covering_time_properties_lem}, we have
\[
\ctime{}{\thisset}{\thisvertex}{\conept}
\le \ctime[2]{}{\thisset}{\thisvertex}{\segment{\thisapex}{\conept}}
\le \rptime[2]{}{\thisset}{\thisvertex}{\segment{\thisapex}{\conept}}
	+ \tmeasureof[2]{\segment{\thisapex}{\conept}}
=\tmeasureof[2]{\segment{\thisapex}{\conept}}.
\]
Note that the last step follows because $\thisapex\in\thiscube$, so $\thisvertex \in \lat \thisapex \subseteq \lat \thisset$, and hence $\rptime[2]{}{\thisset}{\thisvertex}{\segment{\thisapex}{\conept}} \le \rptime[2]{}{\thisset}{\thisvertex}{\thisapex} = 0$.

For Part~\ref{cone_segment_growth:cone_exists}, first recall that if $\thp'\le \minth$, then $\thisunions\subseteq \fatstars{\mu}{\thp'}{\thisapex}$ by \eqref{fatstars_monotone_eqn}, so $\thisset\in \fatstars{\mu}{\thp'}{\thisapex}$. Thus, $\thisset$ is a union of $(\mu,\thp')$-cone segments at $\thisapex$, so for any $\conept\in \thisset$ there is some $\ptdirsymb\in \sphere{d-1}{\mu}$ and $\pthbig\ge 0$ such that $\conept\in \conetip{\mu}{\thp'}{\thisapex}{\ptdirvec}{\pthbig} \subseteq S$. Since
\[
\conetip{\mu}{\thp'}{\thisapex}{\ptdirvec}{\pthbig}
= \bigcup_{h\in [0,\pthbig]} \reacheddset{\mu}{\thisapex + h\ptdirsymb}{\thp' h},
\]
there is some $\pth \in [0,\pthbig]$ such that $\conept\in \reacheddset{\mu}{\thisapex + \pth \ptdirsymb}{\thp' \pth}$. It then follows that
$
\conept \in \conetip{\mu}{\thp'}{\thisapex}{\ptdirvec}{\pth} 
\subseteq \conetip{\mu}{\thp'}{\thisapex}{\ptdirvec}{\pthbig}
\subseteq \thisset,
$
and the triangle inequality for $\mu$ implies that
\[
(1-\thp')\pth \le \distnew{\mu}{\thisapex}{\conept} \le (1+\thp')\pth.
\]
Thus we can take $\ptseth \definedas \pth$ and $\ptsetcone \definedas \conetip{\mu}{\thp'}{\thisapex}{\ptdirvec}{\pth} \in  \conefamily[2]{\mu}{\thisapex,\thp',\ptseth}$.
\end{proof}

Parts~\ref{cone_segment_growth:tmeasure_bound} and \ref{cone_segment_growth:cone_exists} of Claim~\ref{cone_segment_growth:set_reductions_claim} will be used to bound the covering times of the ``near points" and ``far points" in $\thisset$, respectively. The ``near points" are dispatched in the following claim.

\begin{claim}
\label{cone_segment_growth:lineevent_claim}
The event $\nearevent$ occurs on the event $\lineevent$.
\end{claim}

\begin{proof}[Proof of Claim~\ref{cone_segment_growth:lineevent_claim}]
Note that the bound in Part~\ref{cone_segment_growth:tmeasure_bound} of Claim~\ref{cone_segment_growth:set_reductions_claim} does not depend on $\thisset$. Therefore,
\begin{align*}
\nearevent
&= \bigcap_{\thisapex,\thisset}
	\eventthat[1]{ \ctime{}{\thisset}{\thisvertex}{\conept} \le \tmin 
	\text{ for all $\conept\in \thisset$ with $\distnew{\mu}{\thisapex}{\conept} \le \cutradius$}}\\
&\supseteq \bigcap_{\thisapex \in \thiscube}
	\eventthat[1]{ \tmeasureof[2]{\segment{\thisapex}{\conept}} \le \tmin 
	\text{ for all $\conept\in \Rd$ with $\distnew{\mu}{\thisapex}{\conept} \le \cutradius$}}\\
&= \bigcap_{\thisapex \in \thiscube}
	\eventthat[1]{ \tmeasureof[2]{\segment{\thisapex}{\thisapex+\centeredpt}} \le \tmin 
	\text{ for all $\centeredpt\in \Rd$ with
	$\distnew{\mu}{\thisvertex}{\centeredpt} \le \cutradius$}}\\
&\supseteq \lineevent. \qedhere
\end{align*}
\end{proof}

In order to tackle the ``far points," we first define
\[
\sppk(\spp) \definedas \frac{\spp}{1+(1-\spp)\oldkonst},
\quad \thpk(\thp, \spp) \definedas \minth \meet \sppk,
\quad\text{and}\quad
\hk(\tmin,\spp) \definedas \frac{\tmin}{(1+\sppk) \linekonst}.
\]
These definitions were specially formulated to obtain the following relations, which will be needed in the subsequent claim.
\begin{equation}
\label{cone_segment_growth:param_relations_eqn}
\frac{1-\sppk}{1-\spp} = 1+\oldkonst \sppk,
\qquad\thpk \le \sppk,
\quad\text{and}\qquad
(1+\thpk)h\ge \cutradius \implies h\ge \hk.
\end{equation}

\begin{claim}
\label{cone_segment_growth:conesevent_claim}
The event $\farevent$ occurs on the event $\conesevent$.
\end{claim}

\begin{proof}[Proof of Claim~\ref{cone_segment_growth:conesevent_claim}]
Using Part~\ref{cone_segment_growth:cone_exists} of Claim~\ref{cone_segment_growth:set_reductions_claim} (with $\thp'=\thpk\le \minth$) in the second and third lines, and using the relations in \eqref{cone_segment_growth:param_relations_eqn} in the fourth and fifth lines, we have
\begin{align*}
\farevent &=  \bigcap_{\thisapex,\thisset}
	\eventthat[1]{ \ctime{}{\thisset}{\thisvertex}{\conept}
	\le \frac{\distnew{\mu}{\thisapex}{\conept}}{1-\spp}\quad
	\text{$\forall\conept\in \thisset$ with $\distnew{\mu}{\thisapex}{\conept} > \cutradius$}}\\
&\supseteq  \bigcap_{\thisapex,\thisset}
	\eventthat[1]{ \ictime[2]{}{\thisvertex}{\ptsetcone}
	\le \frac{\distnew{\mu}{\thisapex}{\conept}}{1-\spp}
	\quad\forall \conept\in \thisset
	\text{ with } \distnew{\mu}{\thisapex}{\conept} > \cutradius}\\
&\supseteq  \bigcap_{\thisapex,\thisset}
	\eventthat[1]{ \ictime[2]{}{\thisvertex}{\ptsetcone}
	\le \frac{1-\thpk}{1-\spp} \ptseth
	\quad \text{$\forall\conept\in \thisset$ with $(1+\thpk)\ptseth > \cutradius$}}\\
&\supseteq  \bigcap_{\thisapex\in\thiscube} \bigcap_{h\ge \hk}
	\eventthat[1]{ \ictime[2]{}{\thisvertex}{\cones}
	\le \frac{1-\sppk}{1-\spp} h
	\quad \forall \cones\in \conefamily{\mu}{\thisapex,\thpk,h}}\\
&=  \bigcap_{\thisapex\in\thiscube} \bigcap_{h\ge \hk}
	\bigcap_{\dirvec \in \dirset{\Rd}}
	\eventthat[1]{ \ictime[2]{}{\thisvertex}{\conetip{\mu}{\thpk}{\thisapex}{\dirvec}{h}}
	\le \parens[2]{1+\oldkonst \sppk} h}\\
&\supseteq \conesevent. \qedhere
\end{align*}
\end{proof}

Finally, combining Claims~\ref{cone_segment_growth:inversion_claim}, \ref{cone_segment_growth:lineevent_claim}, and \ref{cone_segment_growth:conesevent_claim}, applying Lemmas~\ref{path_tail_bound_lem} and \ref{thin_cone_covering_lem}, and noting that $\thpk \ge \minth \sppk$ since $\minth,\sppk \in (0,1]$, we get
\begin{align*}
\Pr \parens[2]{\thisevent}
&\ge \Pr \parens[2]{\nearevent \cap \farevent}\\
&\ge \Pr \parens[1]{\textstyle \lineevent \cap \conesevent}\\
&\ge 1-\bigcref{path_tail_bound_lem} e^{-\smallcref{path_tail_bound_lem} \cutradius}
	-\frac{1}{\thpk^d} \bigcref{thin_cone_covering_lem}(\sppk)
	e^{-\smallcref{thin_cone_covering_lem}(\sppk)\thpk \hk}\\
&\ge 1-\bigcref{path_tail_bound_lem} e^{-\smallcref{path_tail_bound_lem} \cutradius}
	-\frac{1}{\minth^d}\cdot \frac{1}{\sppk^d}
		\bigcref{thin_cone_covering_lem}(\sppk)
	e^{-\smallcref{thin_cone_covering_lem}(\sppk)
		\frac{\sppk}{(1+\sppk) \linekonst} \minth \tmin}\\
&\ge 1-\frac{1}{\minth^d} \thisbigc(\spp) e^{-\thissmallc(\spp) \minth \tmin},
\end{align*}
where
\[
\thisbigc(\spp) \definedas \bigcref{path_tail_bound_lem}
	+ \frac{1}{\sppk(\spp)^d} \bigcref{thin_cone_covering_lem}\parens[2]{\sppk(\spp)}
\quad\text{and}\quad
\thissmallc(\spp) \definedas \frac{\smallcref{path_tail_bound_lem}}{\linekonst} \meet
	\frac{\sppk(\spp)\smallcref{thin_cone_covering_lem} \parens[2]{\sppk(\spp)}}
	{\parens[2]{1+\sppk(\spp)} \linekonst}.
\]
This completes the proof of Theorem~\ref{cone_segment_growth_thm}.
\end{proof}

The following corollary gives a large deviations estimate for first-passage growth restricted to a $\mu$-cone (or more generally any cone that is a $\mu$-star), showing that first-passage percolation restricted to a cone grows asymptotically as fast as unrestricted first-passage percolation.

\begin{cor}[Large deviations estimate for growth in $\mu$-cones]
\label{mu_cone_large_dev_cor}
\newcommand{\thiscone}{\cones}
\newcommand{\thisC}{\bigcref{mu_cone_large_dev_cor}(\spp)}
\newcommand{\thisc}{\smallcref{mu_cone_large_dev_cor}(\spp)}

Let $\tmeasure$ be an \iid\ traversal measure on $\edges{\Zd}$ satisfying \finitespeed{d} and \expmoment, and let $\mu$ be the shape function for $\tmeasure$ from Theorem~\ref{shape_thm}.
Then for any $\epsilon\in (0,1)$, there exist constants $\thisC$ and $\thisc$ such that if $\cones$ is any cone in $\Rd$ that is also a $(\mu,\thp)$-star at $\vect 0$ for some $\thp \in (0,1]$, then
\begin{multline*}
\Pr \eventthat[3]{ {(1-\epsilon)t}\unitball{\mu} \cap \thiscone
\subseteq \crreachedset{\tmeasure}{\vect 0}{\thiscone}{t} \subseteq
{(1+\epsilon)t \unitball{\mu}} \cap \thiscone
\text{ for all $t\ge t_0$}} \ge 1- \frac{1}{\thp^d} \thisC e^{-\thisc \thp t_0}.
\end{multline*}
In particular, this holds if $\thiscone = \cone{\mu}{\thp}{\vect 0}{\dirvec}$ for some $\dirvec\in \dirset{\Rd}$.
%
\end{cor}

\begin{proof}
\newcommand{\thiscone}{\cones}
\newcommand{\thisC}{\bigcref{mu_cone_large_dev_cor}(\spp)}
\newcommand{\thisc}{\smallcref{mu_cone_large_dev_cor}(\spp)}

\newcommand{\upperC}{\bigcref{shifted_dev_bounds_lem}(\spp)}
\newcommand{\upperc}{\smallcref{shifted_dev_bounds_lem}(\spp)}

\newcommand{\lowerC}{\bigcref{cone_segment_growth_thm}(\spp)}
\newcommand{\lowerc}{\smallcref{cone_segment_growth_thm}(\spp)}

This follows directly from Proposition~\ref{restricted_growth_ub_prop} and Theorem~\ref{cone_segment_growth_thm}: Take
\[
\thisC \definedas \upperC+\lowerC
\quad\text{and}\quad
\thisc \definedas \upperc\meet\lowerc,
\]
and note that $\rreacheddset{\mu}{\vect 0}{\thiscone}{r} = \parens{r\unitball{\mu}} \cap \thiscone$ for all $r\ge 0$ by Proposition~\ref{determ_star_growth_prop} since $\thiscone\in \stars{\vect 0}$.
\end{proof}

\begin{thmremark}
\newcommand{\thisC}{\bigcsymb(\spp)}
\newcommand{\thisc}{\smallcsymb(\spp)}
Note that for a fixed $(\mu,\thp)$-cone $\cones = \cone{\mu}{\thp}{\vect 0}{\dirvec}$, Remark~\ref{thin_cone_covering:fixed_dir_rem} implies that the bound in Corollary~\ref{mu_cone_large_dev_cor} can be improved to $\Pr \ge 1- \frac{1}{\thp} \thisC e^{-\thisc \thp t_0}$ for some positive constants $\thisC$ and $\thisc$.
\end{thmremark}

\subsection{Shape Theorems in $\mu$-Stars and Cones}
\label{cone_growth:mu_stars_shape_thm_sec}


Together, Proposition~\ref{restricted_growth_ub_prop} and Theorem~\ref{cone_segment_growth_thm} immediately yield a Shape Theorem for growth restricted to a $\mu$-star. In fact, we trivially get various formulations of the Shape Theorem that are simultaneous over various collections of $\mu$-stars; the following theorem states several such results.
Recall from \eqref{more_fat_stars_defs_eqn} that $\normstars{\mu}{\vect z} \definedas \bigcup_{\thp>0} \fatstars{\mu}{\thp}{\vect z}$ is the set of all $\mu$-stars at $\vect z\in\Rd$.

\begin{thm}[Shape Theorems for growth in $\mu$-stars]
\label{mu_stars_shape_thm}

\newcommand{\thisset}{S} 
\newcommand{\thisvertex}{\vect v}
\newcommand{\thiscenter}{\vect z}
\newcommand{\thiscube}{\cubify{\thisvertex}}
\newcommand{\thisstars}{\fatstars{\mu}{\thp}{\thiscenter}}

\newcommand{\thisevent}[4][0]{\eventref{mu_stars_shape_thm}^{#2} \parens[#1]{#3;#4}}

Suppose $\tmeasure$ is an \iid\ traversal measure on $\Zd$ satisfying \finitespeed{d} and \expmoment, and let $\mu$ be the shape function for $\tmeasure$ from Theorem~\ref{shape_thm}.
Let $\thisvertex\in\Zd$ and $\thp,\spp\in (0,1]$. Then each of the following events occurs almost surely:
\begin{enumerate}
\item \label{mu_stars_shape:all_fixed_part}
$\displaystyle
\thisevent{\thiscube}{\thp}{\spp}
\definedas \eventthat[1]{ \mathstack{\exists t_0\ge 0 \text{ such that for any }
	\thiscenter\in\thiscube \text{ and } \thisset\in \thisstars,\, \forall t\ge t_0,}
{\rreacheddset[2]{\mu}{\vect z}{\thisset}{(1-\spp)t}
\subseteq \crreachedset{\tmeasure}{\thisvertex}{\thisset}{t}
\subseteq \rreacheddset[2]{\mu}{\vect z}{\thisset}{(1+\spp)t}}};
$

\item \label{mu_stars_shape:fixed_d_e_part}
$\displaystyle
\thisevent{\Rd}{\thp}{\spp} 
\definedas \eventthat[1]{
	\mathstack{\text{For all } \thiscenter\in\Rd,\,
	\exists t_0\ge 0 \text{ such that if }
	\thisset\in \thisstars, \text{ then } \forall t\ge t_0,}
{\rreacheddset[2]{\mu}{\vect z}{\thisset}{(1-\spp)t}
\subseteq \crreachedset{\tmeasure}{\thiscenter}{\thisset}{t}
\subseteq \rreacheddset[2]{\mu}{\vect z}{\thisset}{(1+\spp)t}}};
$

\item \label{mu_stars_shape:fixed_v_e_part}
$\displaystyle
\thisevent[2]{\thiscube}{0^+}{\spp} 
\definedas \eventthat[1]{
	\mathstack{\text{For all } \thiscenter\in\thiscube
	\text{ and } \thisset\in \normstars{\mu}{\thiscenter},\; 
	\exists t_0\ge 0 \text{ such that } \forall t\ge t_0,}
{\rreacheddset[2]{\mu}{\vect z}{\thisset}{(1-\spp)t}
\subseteq \crreachedset{\tmeasure}{\thisvertex}{\thisset}{t}
\subseteq \rreacheddset[2]{\mu}{\vect z}{\thisset}{(1+\spp)t}}};
$

\item \label{mu_stars_shape:fixed_e_part}
$\displaystyle
\thisevent[2]{\Rd}{0^+}{\spp} 
\definedas \eventthat[1]{
	\mathstack{\text{For all } \thiscenter\in\Rd
	\text{ and } \thisset\in \normstars{\mu}{\thiscenter},\; 
	\exists t_0\ge 0 \text{ such that } \forall t\ge t_0,}
{\rreacheddset[2]{\mu}{\vect z}{\thisset}{(1-\spp)t}
\subseteq \crreachedset{\tmeasure}{\thiscenter}{\thisset}{t}
\subseteq \rreacheddset[2]{\mu}{\vect z}{\thisset}{(1+\spp)t}}};
$

\end{enumerate}
\end{thm}

\begin{proof}
\newcommand{\thisset}{S} 
\newcommand{\thisvertex}{\vect v}
\newcommand{\thiscenter}{\vect z}
\newcommand{\thiscube}{\cubify{\thisvertex}}
\newcommand{\thisstars}{\fatstars{\mu}{\thp}{\thiscenter}}

\newcommand{\thisevent}[4][0]{\eventref{mu_stars_shape_thm}^{#2} \parens[#1]{#3;#4}}
\newcommand{\starsdevevent}[1]{\eventref{restricted_growth_ub_prop}^{\thisvertex;\bigstar} \parens{#1;\spp}}
\newcommand{\mustarsdevevent}[1]{\eventref{cone_segment_growth_thm}^{\thisvertex} \parens{#1,\thp;\spp}}

For Part~\ref{mu_stars_shape:all_fixed_part}, simply observe that
\begin{equation*}
\thisevent{\thiscube}{\thp}{\spp}
\supseteq \bigcup_{t_0\ge 0} \starsdevevent{t_0} \cap \mustarsdevevent{t_0},
\end{equation*}
so Proposition~\ref{restricted_growth_ub_prop} and Theorem~\ref{cone_segment_growth_thm} imply that $\Pr \parens[2]{\thisevent{\thiscube}{\thp}{\spp}} = 1$ by monotnicity of measure.
Part~\ref{mu_stars_shape:fixed_d_e_part} then follows immediately since $\thisevent{\Rd}{\thp}{\spp} = \bigcap_{\vect v\in\Zd} \thisevent{\thiscube}{\thp}{\spp}$. Then Part~\ref{mu_stars_shape:fixed_v_e_part} follows from Part~\ref{mu_stars_shape:all_fixed_part} because $\thisevent[2]{\thiscube}{0^+}{\spp} = \bigcap_{n\in \naturalpositive} \thisevent{\thiscube}{1/n}{\spp}$, and Part~\ref{mu_stars_shape:fixed_e_part} follows from Part~\ref{mu_stars_shape:fixed_d_e_part} because $\thisevent[2]{\Rd}{0^+}{\spp} = \bigcap_{n\in \naturalpositive} \thisevent{\Rd}{1/n}{\spp}$.
\end{proof}

We now state two corollaries of Theorem~\ref{mu_stars_shape_thm} for the case when the restricting set $S$ is a cone.


\begin{cor}[A Shape Theorem in cones that are $\mu$-stars]
\label{mu_cone_shape_thm_cor}
\newcommand{\thiscone}{\cones}

Let $\tmeasure$ be an \iid\ traversal measure on $\edges{\Zd}$ satisfying \finitespeed{d} and \expmoment.
If $\cones$ is any cone in $\Rd$ that is also a $\mu$-star at $\vect 0$, then for any $\epsilon\in (0,1)$,
\[
\Pr \eventthat[3]{ {(1-\epsilon)t}\unitball{\mu} \cap \thiscone
\subseteq \crreachedset{\tmeasure}{\vect 0}{\thiscone}{t} \subseteq
{(1+\epsilon)t \unitball{\mu}} \cap \thiscone
\quad \text{for all large $t$}} = 1,
\]
where $\mu$ is the shape function for $\tmeasure$ from Theorem~\ref{shape_thm}. In particular, this holds if $\thiscone$ is any $\mu$-cone at $\vect 0$.
%
\end{cor}

\begin{proof}
\newcommand{\thiscone}{\cones}
This follows directly from Part~\ref{mu_stars_shape:fixed_v_e_part} of Theorem~\ref{mu_stars_shape_thm} since $\thiscone\in \normstars{\mu}{\vect 0}$ by assumption, and $\rreacheddset{\mu}{\vect 0}{\thiscone}{r} = \parens{r \unitball{\mu}} \cap \thiscone$ for any $r\ge 0$ by Proposition~\ref{determ_star_growth_prop} since $\thiscone\in \stars{\vect 0}$.
\end{proof}

The following result will be needed in Chapter~\ref{coex_finite_chap}.
Using the notation \eqref{cone_fam_h_def_eqn} in Remark~\ref{thin_cone_covering:notation_rem}, denote the set of all $\mu$-cones at $\vect z\in\Rd$ by
\begin{equation}
\label{cone_fam_infty_def_eqn}
\conefamily[2]{\mu}{\vect z,0^+, \infty}
\definedas \bigcup_{\thp>0} \conefamily{\mu}{\vect z,\thp, \infty}
= \setbuilderbar[1]{\cone{\mu}{\thp}{\vect z}{\dirvec}}{\dirvec\in\dirset{\Rd},\, \thp>0}.
\end{equation}

\begin{cor}[A simultaneous Shape Theorem in $\mu$-cones]
\label{simultaneous_cone_shape_thm_cor}
\newcommand{\thiscone}{\cones}

\newcommand{\thispt}{\vect z}
\newcommand{\thisevent}{\eventref{simultaneous_cone_shape_thm_cor} (\spp)}
\newcommand{\thisconefam}{\conefamily[2]{\mu}{\thispt,0^+, \infty}}

Let $\tmeasure$ and $\mu$ be as in Theorem~\ref{mu_stars_shape_thm}. For $\spp\in (0,1)$, let
\[
\thisevent
\definedas \bigcap_{\thispt\in\Rd}
\eventthat[1]{
	\mathstack{
	\forall \thiscone\in \thisconefam,\; 
	\exists t_0\ge 0 \text{ such that } \forall t\ge t_0,}
{\rreacheddset[2]{\mu}{\vect z}{\thiscone}{(1-\spp)t}
\subseteq \crreachedset{\tmeasure}{\thispt}{\thiscone}{t}
\subseteq \rreacheddset[2]{\mu}{\vect z}{\thiscone}{(1+\spp)t}}},
\]
where $\thisconefam$ is defined in \eqref{cone_fam_infty_def_eqn}.
Then $\Pr \parens[2]{\thisevent} = 1$ for any $\spp\in (0,1)$.

%
%
\end{cor}

\begin{proof}
\newcommand{\thiscone}{\cones}
\newcommand{\thispt}{\vect z}
\newcommand{\thisevent}{\eventref{simultaneous_cone_shape_thm_cor} (\spp)}
\newcommand{\thisconefam}{\conefamily[2]{\mu}{\thispt,0^+, \infty}}

This follows directly from Part~\ref{mu_stars_shape:fixed_e_part} of Theorem~\ref{mu_stars_shape_thm} since $\thisconefam \subseteq \normstars{\mu}{\thispt}$.
\end{proof}

\section{Growth in a Super-Logarithmically Expanding Tube}
\label{cone_growth:logarithmic_tube_sec}

{
\newcommand{\hmax}{h_{\max}}


In this section we apply the large deviations estimate from Theorem~\ref{cone_segment_growth_thm} to obtain a Shape Theorem for the growth of a first-passage percolation process restricted to a tube-shaped region. Our definition of ``tube" is a natural generalization of cylinders and cones; specifically, the result will be stated for ``$\mu$-tubes" as defined in Section~\ref{deterministic:norm_geometry_sec}.

It is well known (see, e.g.\ \cite{\DHunbounded} or \cite{Deijfen:2004aa}) that for growth in a cylinder or half-cylinder---that is, a tube of constant width---the restricted process has an asymptotic speed, and this speed can be made arbitrarily close to the asymptotic speed of an unrestricted process by choosing a cylinder of sufficiently large width.
On the other hand, the results of Section~\ref{cone_growth:cone_growth_sec} above show that for growth in a cone---that is, a tube whose width grows linearly with its height---the restricted process has an asymptotic speed in every direction of the cone, and this speed \emph{equals} the asymptotic speed of an unrestricted process. It is natural to ask what happens in intermediate tubes---those in which the width increases at a rate that is $\omega(1)$ but $o(h)$, where $h$ is the ``height," or distance from the origin. We show in Theorem~\ref{tube_shape_thm} below that an expansion rate of $\omega(\log h)$ in the tube is sufficient for the restricted process to grow asymptotically as fast as an unrestricted process. Before giving a formal proof of this result in Section~\ref{cone_growth:logarithmic_tube_proof_sec}, we first give a heuristic argument  in Section~\ref{cone_growth:heuristic_tube_sec} showing how the expansion rate of $\omega(\log h)$ is obtained.

\begin{remark}
\label{logarithmic_tube_suboptimal_rem}
In fact, recent papers by Chatterjee and Dey \cite{Chatterjee:2009aa} and Ahlberg \cite{Ahlberg:2011aa}, \cite{Ahlberg:2011ab} have shown that restricted first-passage percolation grows asymptotically as fast as an unrestricted process in tubes with expansion rate $\omega(1)$, i.e.\ any tube whose width goes to $\infty$, and moreover that the asymptotic speed of the restricted process is strictly slower in tubes of bounded width. Thus, the hypothesis of expansion rate $\omega(\log h)$ in Theorem~\ref{tube_shape_thm} is suboptimal, and as noted in Section~\ref{cone_growth:heuristic_tube_sec} below, it cannot be improved using the method presented here.

On the other hand, since the basis of the proof of Theorem~\ref{tube_shape_thm} is Lemma~\ref{thin_cone_covering_lem}, we obtain a Shape Theorem for tubes in all directions simultaneously.
As noted in Remark~\ref{thin_cone_covering:fixed_dir_rem}, the probability bound in Lemma~\ref{thin_cone_covering_lem} improves by a factor of $\thp^{d-1}$ if we instead consider a fixed direction. This improvement is almost but not quite enough to make the argument in Section~\ref{cone_growth:heuristic_tube_sec} go through for a single fixed tube of expansion rate $\omega(1)$; to make it work, we would need to remove the final factor of $\frac{1}{\thp}$ from the front of the bound in Lemma~\ref{thin_cone_covering_lem}, so that $\thp$ only appears in the exponent. However, it does not seem to be possible to make this improvement in the proof of Lemma~\ref{thin_cone_covering_lem} (or rather, Lemma~\ref{space_covering_lem}, which is where the bound originally came from). This raises the question of whether the bound in (the fixed direction version of) Lemmas~\ref{space_covering_lem} and \ref{thin_cone_covering_lem} is in fact optimal, since this factor of $\frac{1}{\thp}$ is the only obstruction to attaining the optimal expansion rate of $\omega(1)$ in the proof below.
\end{remark}

\subsection{Heuristic Argument Underlying the Proof of Theorem~\ref{tube_shape_thm}}
\label{cone_growth:heuristic_tube_sec}

The idea behind the proof of Theorem~\ref{tube_shape_thm} below is to exploit the dependence on $\thp$ in the probability estimate for the internal covering time of a $\mu$-cone segment $\cones$ of thickness $\thp$ and height $h$. Lemma~\ref{thin_cone_covering_lem} shows that for fixed $\spp>0$ and $\thp \in (0,\spp]$, the probability that the internal covering time of $\cones$ and the unrestricted covering time of $\cones$ differ by more than a constant times $\spp h$ is bounded above by
\begin{equation}
\label{yakyakyak}
\frac{1}{\thp^d} C(\spp) e^{-c(\spp) \thp h}
= C(\spp) \exp \brackets[2]{-c(\spp) \thp h + d \log (\thp^{-1}) },
\end{equation}
where $C(\spp)$ and $c(\spp)$ are positive constants depending on $\spp$. Now consider a tube
\[
\tubesymb = \bigcup_{h\ge 0} \conetip{\mu}{\thp_h}{\vect z}{\dirvec}{h},
\]
consisting of a union of $\mu$-cone segments with common axis $\vect z+\dirvec$, such that the thickness $\thp_h$ of each cone segment is allowed to depend on its height $h$.
The tube $\tubesymb$ is a $\mu$-tube as defined by \eqref{mu_tube_segment_def_eqn} and \eqref{mu_tube_def_eqn} in Section~\ref{deterministic:norm_geometry_sec}.
If $\thp_h\to 0$ as $h\to \infty$ (i.e.\ the width of  $\tubesymb$ expands at a strictly sublinear rate), then $\thp_h \le \spp$ for all sufficiently large $h$, so the probability bound in \eqref{yakyakyak} applies (with $\thp = \thp_h$) to the internal covering time of $\conetip{\mu}{\thp_h}{\vect z}{\dirvec}{h}$ for all large $h$. We can now ask: How does this probability bound depend on the choice of function $\thp_h$? 

The key idea is to choose the function $\thp_h$ so that the resulting probability bound in \eqref{yakyakyak} is summable with respect to $h$ for any fixed $\spp$, in order to obtain a shape theorem in $\tubesymb$ using the Borel-Cantelli Lemma. For example, if the probability bound in \eqref{yakyakyak} is summable, then Borel-Cantelli implies that we almost surely get an $\epsilon$-good bound on the internal covering time of, say, $\conetip{\mu}{\thp_n}{\vect z}{\dirvec}{n}$ for all sufficiently large $n\in\N$. This will imply that the ratio between the growth rates of the $\tubesymb$-restricted process and the unrestricted process is on the order of $1+\spp$.

How quickly can we allow $\thp_h$ to approach 0 if we want the expression in \eqref{yakyakyak} to be summable in $h$ for any $\spp$? To start with, the bound must approach 0 as $h\to\infty$, which means that we need
\[
-c(\spp) \thp_h h + d \log (\thp_h^{-1}) \to -\infty
\quad\text{as } h\to\infty.
\]
Since $c(\spp)$ becomes arbitrarily small as $\spp\to 0$, the only way this can hold for all $\epsilon>0$ is if $\thp_h h$ grows faster than $\log (\thp_h^{-1})$, i.e.\ the ratio $\frac{\thp_h h}{\log (\thp_h^{-1})} \to\infty$ as $h\to\infty$. Setting $\rfcn(h) = \thp_h h$, this necessary condition becomes
\begin{equation}
\label{yaksnork}
\frac{\rfcn(h)}{\log \brackets[2]{h/\rfcn(h)}} \to \infty
\quad\text{as } h\to\infty.
\end{equation}
Clearly, $\rfcn(h)$ satisfies \eqref{yaksnork} iff $\frac{\rfcn(h)}{\log h}\to \infty$, i.e.\ $\rfcn(h) = \omega(\log h)$. Moreover, if $\rfcn(h) = \omega(\log h)$, then it is easily seen that for any fixed $\epsilon$, the function
\[
f(h) = \exp \brackets[1]{-c(\spp) \rfcn(h) + d \log \parens[1]{\tfrac{h}{\rfcn(h)}}}
\]
decreases faster than any power of $h$ and hence is summable. Thus, the necessary condition $\rfcn(h) = \omega(\log h)$ is also sufficient to obtain a summable probability bound. Now observe that the function $\rfcn(h) = \thp_h h$ simply gives the radius of the final $\mu$-ball in the cone segment $\conetip{\mu}{\thp_h}{\vect z}{\dirvec}{h}$; if $\thp_h$ is a decreasing function, then this is equal to the ``inner $\mu$-radius" of the tube $\tubesymb$ at height $h$, i.e.\ the radius of the largest $\mu$-ball contained in $\tubesymb$ and centered at $\vect z+h\munitvec$. Thus, taking a tube of expansion rate $\omega(\log h)$ is sufficient to make the above argument go through, and the argument will not work for any expansion rate that is $O(\log h)$.


\subsection{Proof of Shape Theorem in a $\mu$-Tube with Expansion Rate $\omega(\log h)$}
\label{cone_growth:logarithmic_tube_proof_sec}
 

Recall from \eqref{mu_tube_segment_def_eqn} and \eqref{mu_tube_def_eqn} in Section~\ref{deterministic:norm_geometry_sec} that for a norm $\mu$ on $\Rd$, the \textdef{initial $\mu$-tube segment} with parameters $\vect z\in\Rd$, $\dirvec\in\dirset{\Rd}$, $h \in [0,\infty]$, and $\rfcn\colon \Rplus\to\Rplus$ is defined to be
\begin{equation*}
\tubetip{\mu}{\rfcn}{\vect z}{\dirvec}{h}
	\definedas \bigcup_{s\in [0,h]}
	\reacheddset[2]{\mu}{\vect z+h \munitvec}{\rfcn(s)},
\end{equation*}
and the corresponding infinite \textdef{$\mu$-tube} is
\begin{equation*}
\tube{\mu}{\rfcn}{\vect z}{\dirvec}
	\definedas \tubetip{\mu}{\rfcn}{\vect z}{\dirvec}{\infty}
	\definedas \bigcup_{h\ge 0} \tubetip{\mu}{\rfcn}{\vect z}{\dirvec}{h}.
\end{equation*}
Theorem~\ref{tube_shape_thm} below will use the following result about the geometry of $\mu$-tubes.

\begin{lem}[Geometric properties of $\mu$-tubes with sublinear expansion]
\label{sublinear_tube_geometry_lem}
\newcommand{\thisz}{\vect z}
\newcommand{\thisdir}{\dirvec}
\newcommand{\thisu}{\dirsymb}
\newcommand{\thistube}{\tubesymb}
\newcommand{\explicittube}{\tube{\mu}{\rfcn}{\thisz}{\thisdir}}

\newcommand{\thistubeh}[1]{\tubesymb_{#1}}
\newcommand{\explicittubeh}[2][0]{\tubetip[#1]{\mu}{\rfcn}{\thisz}{\thisdir}{#2}}

\newcommand{\thph}[1]{\thp_{#1}}
\newcommand{\thisconeh}[2][0]{\conetip[#1]{\mu}{\thph{#2}}{\thisz}{\thisdir}{#2}}

\newcommand{\initialball}{\reacheddsymb_0}

Fix $\thisz\in\Rd$, $\thisu\in \sphere{d-1}{\mu}$, a norm $\mu$ on $\Rd$, and a nonnegative function $\rfcn$ on $\Rplus$. Let $\thistube = \explicittube$, and for each $h\ge 0$ let $\thistubeh{h} = \explicittubeh{h}$, so $\setof{\thistubeh{h}}_{h\ge 0}$ forms an increasing family of sets in $\Rd$ whose union is $\thistube$.
If the function $\thph{h}\definedas \rfcn(h)/h$ is nonincreasing on $(0,\infty)$, then for all $h\ge 0$,
\begin{enumerate}
\item \label{sublinear_tube_geometry:union_part}
$\displaystyle
\thistubeh{h} =  \initialball \cup \bigcup_{h'\in (0,h]} \thisconeh{h'},
$
where $\initialball \definedas \reacheddset[2]{\mu}{\thisz}{\rfcn(0)}$.

\item \label{sublinear_tube_geometry:star_part}
$\thistubeh{h}\in \fatstars{\mu}{\thph{h}}{\thisz}$, and $\thistube\in \stars{\thisz}$.

\item \label{sublinear_tube_geometry:growth_part}
If $s\le (1-\thph{h})h$, then $\rreacheddset{\mu}{\thisz}{\thistube}{s} \subseteq \thistubeh{h}$.

\end{enumerate}
\end{lem}

\begin{proof}
\newcommand{\thisz}{\vect z}
\newcommand{\thisdir}{\dirvec}
\newcommand{\thisu}{\dirsymb}

\newcommand{\thistube}{\tubesymb}
\newcommand{\explicittube}{\tube{\mu}{\rfcn}{\thisz}{\thisdir}}
\newcommand{\thistubeh}[1]{\tubesymb_{#1}}
\newcommand{\explicittubeh}[2][0]{\tubetip[#1]{\mu}{\rfcn}{\thisz}{\thisdir}{#2}}

\newcommand{\thph}[1]{\thp_{#1}}
\newcommand{\thisconeh}[2][0]{\conetip[#1]{\mu}{\thph{#2}}{\thisz}{\thisdir}{#2}}

\newcommand{\initialball}{\reacheddsymb_0}

\newcommand{\thisx}{\vect x}
\newcommand{\thish}{h_{\vect x}}
\newcommand{\thiscenter}{\thisz+ \thish \thisu}
\newcommand{\thisxball}{\reacheddset[2]{\mu}{\thiscenter}{\rfcn(\thish)}}

For each $h>0$, let $\thph{h} \definedas \rfcn(h)/h$, and let $\thph{0} \definedas 0$.

\proofpart{\ref{sublinear_tube_geometry:union_part}}
By assumption, $\thph{h}$ is a nonincreasing function of $h$ on $(0,\infty)$, so for any $h\ge 0$ we have 
\begin{align*}
\thisconeh{h} &= 
	\bigcup_{h' \in [0,h]} \reacheddset[2]{\mu}{\thisz+h'\thisu}{\thph{h} \cdot h'}
\subseteq \reacheddset[2]{\mu}{\thisz}{\rfcn(0)}
	\cup \bigcup_{h' \in (0,h]} \reacheddset[2]{\mu}{\thisz+h'\thisu}{\thph{h'} \cdot h'}
=\thistubeh{h}.
\end{align*}
Since the collection $\setof{\thistubeh{h}}_{h\ge 0}$ is increasing in $h$, this shows that for any $h\ge 0$ we have
\[
\thistubeh{h} \supseteq  \reacheddset[2]{\mu}{\thisz}{\rfcn(0)} \cup \bigcup_{h'\in (0,h]} \thisconeh{h'},
\]
and since $\reacheddset[2]{\mu}{\thisz+h'\thisu}{\rfcn(h')} \subseteq \thisconeh{h'}$, the reverse inclusion is trivial.

\proofpart{\ref{sublinear_tube_geometry:star_part}} Since $\thph{h'}\ge \thph{h}$ for all $0<h'\le h$, it follows from Part~\ref{sublinear_tube_geometry:union_part} and Lemma~\ref{mu_cone_decomp_lem} that every point in $\thistubeh{h}$ is contained in a $\mu$-cone with thickness $\thph{h}$ and apex $\thisz$, so $\thistubeh{h}\in \fatstars{\mu}{\thph{h}}{\thisz} \subseteq \normstars{\mu}{\thisz}$ for all $h\ge 0$. Since any union of elements of $\normstars{\mu}{\thisz}$ is star-shaped at $\thisz$, we then have $\thistube = \bigcup_{h\ge 0} \thistubeh{h}\in \stars{\thisz}$.

\proofpart{\ref{sublinear_tube_geometry:growth_part}} We will show that if $\thisx\in \thistube\setminus \thistubeh{h}$, then $\distnew{\mu}{\thisz}{\thisx} > (1-\thph{h})h$. Suppose $\thisx\in \thistube\setminus \thistubeh{h}$. Since $\thisx\in\thistube$, there is some $\thish \in [0,\infty)$ such that $\thisx\in \thisxball$, and since $\thisx\not\in \thistubeh{h}$, we must have $\thish >h$ for any such $\thish$. Therefore, using the reverse triangle inequality for $\mu$ and the fact that $\thph{h}$ is nonincreasing,
\begin{align*}
\distnew{\mu}{\thisz}{\thisx}
&\ge \distnew{\mu}{\thisz}{\thiscenter} -\distnew{\mu}{\thiscenter}{\thisx}
	&& \parens[2]{\text{reverse triangle inequality}}\\
&\ge \thish - \rfcn(\thish)
	&& \parens[2]{\text{since } \thisx\in\thisxball\, }\\
&= \parens{1-\thph{\thish}} \thish
	&& \parens[2]{\text{definition of $\thph{\thish}$}} \\
&\ge \parens{1-\thph{h}} \thish
	&& \parens[2]{\thish>h \implies \thph{\thish}\le \thph{h}}\\
&> \parens{1-\thph{h}} h.
	&& \parens[2]{\text{since } \thish>h}
\qedhere
\end{align*}
\end{proof}

We are now ready to prove the main result of the section. Rather than directly applying Lemma~\ref{thin_cone_covering_lem} as described in Section~\ref{cone_growth:heuristic_tube_sec}, the proof of Theorem~\ref{tube_shape_thm} will use Theorem~\ref{cone_segment_growth_thm}, since much of the work for translating Lemma~\ref{thin_cone_covering_lem} into a Shape Theorem was already done there.

\begin{thm}[Shape theorem in a super-logarithmically expanding $\mu$-tube]
\label{tube_shape_thm}
\newcommand{\maxratio}{\thp}
\newcommand{\thisz}{\vect z}
\newcommand{\thisdir}{\dirvec}
\newcommand{\thisdirsymb}{\dirsymb}
\newcommand{\thistube}{\tubesymb}
\newcommand{\explicittube}{\tube{\mu}{\rfcn}{\thisz}{\thisdir}}

Let $\tmeasure$ be an \iid\ traversal measure on $\Zd$ satisfying \finitespeed{d} and \expmoment, and suppose $\rfcn$ is a nonnegative function on $\Rplus$ satisfying:
\begin{enumerate}
\item $\dfrac{\rfcn(h)}{h}$ is nonincreasing;
\item $\displaystyle \lim_{h\to\infty} \frac{\rfcn(h)}{h} = \maxratio<1$;
\item $\displaystyle \lim_{h\to\infty}\frac{\rfcn(h)}{\log h} = \infty$.
\end{enumerate}
If $\mu$ is the shape function for $\tmeasure$ from Theorem~\ref{shape_thm}, and $\thistube = \explicittube$ for some $\thisz\in\Rd$ and $\thisdirsymb\in \sphere{d-1}{\mu}$, then for any $\spp\in (0,1)$,
\begin{equation}
\label{tube_shape_thm_eqn}
\Pr \eventthat[3]{
\rreacheddset[2]{\mu}{\thisz}{\thistube}{(1-\spp)t}
\subseteq \crreachedset{\tmeasure}{\thisz}{\thistube}{t}
\subseteq \rreacheddset[2]{\mu}{\thisz}{\thistube}{(1+\spp)t}
\text{ for all large $t$}} =1.
\end{equation}
\end{thm}

\begin{proof}
\newcommand{\thplim}{\thp}

\newcommand{\thisv}{\vect v}
\newcommand{\thisz}{\vect z}
\newcommand{\thisdir}{\dirvec}
\newcommand{\thisdirsymb}{\dirsymb}
\newcommand{\thistube}{\tubesymb}
\newcommand{\explicittube}{\tube{\mu}{\rfcn}{\thisz}{\thisdir}}

\newcommand{\thistubeh}[1]{\tubesymb_{#1}}
\newcommand{\explicittubeh}[2][0]{\tubetip[#1]{\mu}{\rfcn}{\thisz}{\thisdir}{#2}}

\newcommand{\thph}[1]{\thp_{#1}}
\newcommand{\thisconeh}[2][0]{\conetip[#1]{\mu}{\thph{#2}}{\thisz}{\thisdir}{#2}}

\newcommand{\badevent}{F(\spp)}
\newcommand{\bighevent}[1]{G_{#1}(\spp)}
\newcommand{\badhevent}[1]{F_{#1}(\spp)}

\newcommand{\thisC}{\bigcsymb(\spp)}
\newcommand{\thisc}{\smallcsymb(\spp)}

\newcommand{\firstn}{N}
\newcommand{\secondn}{N'}
\newcommand{\maxn}{N''}

Fix $\spp\in (0,1)$ and choose any $\thisv\in \lat{\thisz}$. Observe that since $\thistube\in \stars{\thisz}$ by Part~\ref{sublinear_tube_geometry:star_part} of Lemma~\ref{sublinear_tube_geometry_lem}, the upper containment in \eqref{tube_shape_thm_eqn} follows immediately from Proposition~\ref{restricted_growth_ub_prop} and the Borel-Cantelli Lemma. Thus it remains to prove the lower containment, or equivalently, to show that the complementary event
\[
\badevent \definedas \eventthat[3]{
\exists \setof{t_n}_{n\in\N} \text{ with $t_n\to\infty$ such that } 
\rreacheddset[2]{\mu}{\thisz}{\thistube}{(1-\spp)t_n}
\not \subseteq \crreachedset{\tmeasure}{\thisv}{\thistube}{t_n}}
\]
has probability zero. We will do this by finding a suitable sequence of events $\setof{\badhevent{n}}_{n\in\N}$ such that $\badevent = \setof{\badhevent{n} \text{ i.o.}}$, and then applying Borel-Cantelli.

Let $\alpha \definedas \frac{1-\thplim}{2} \in (0,1-\thplim)$. For each $h>0$, let $\thistubeh{h} \definedas \explicittubeh{h}$ and $\thph{h} \definedas \rfcn(h)/h$ as in Lemma~\ref{sublinear_tube_geometry_lem}, and define the event
\[
\bighevent{h} \definedas \eventthat[3]{
\crreachedset{\tmeasure}{\thisv}{\thistubeh{h}}{t}
\supseteq \rreacheddset[2]{\mu}{\thisz}{\thistubeh{h}}{(1-\spp)t}
\text{ for all } t\ge \alpha h}.
\]
Then since $\thistubeh{h}\in \fatstars{\mu}{\thph{h}}{\thisz}$ by Part~\ref{sublinear_tube_geometry:star_part} of Lemma~\ref{sublinear_tube_geometry_lem}, Theorem~\ref{cone_segment_growth_thm} implies that there are positive constants $\thisC$ and $\thisc$ such that
\begin{equation}
\label{tube_shape:badevent_prob_eqn}
\Pr \parens[2]{\bighevent{h}} \ge 1- \thph{h}^{-d} \thisC e^{-\thisc \thph{h} \alpha h}
\quad\forall h>0.
\end{equation}
The following claim illuminates a connection between $\bighevent{h}$ and $\eventcomplement{\badevent}$.

\begin{claim}
\label{tube_shape:event_containment_claim}
There exists $\firstn \in\N$ such that for any $h> \firstn$,
\[
\bighevent{h} \subseteq \eventthat[3]{
\crreachedset{\tmeasure}{\thisv}{\thistube}{t}
\supseteq \rreacheddset[2]{\mu}{\thisz}{\thistube}{(1-\spp)t}
\text{ for all } t \in \segment[2]{\alpha h}{\alpha (h+1)}}.
\]
\end{claim}

\begin{proof}[Proof of Claim~\ref{tube_shape:event_containment_claim}]

\newcommand{\thisx}{\vect x}
\newcommand{\thish}{h_{\vect x}}
\newcommand{\thiscenter}{\thisz+ \thish \thisdirsymb}
\newcommand{\thisxball}{\reacheddset[2]{\mu}{\thiscenter}{\rfcn(\thish)}}

For any $h>0$, we have $\thistube \supseteq \thistubeh{h}$ and hence $\crreachedset{\tmeasure}{\thisv}{\thistube}{t} \supseteq \crreachedset{\tmeasure}{\thisv}{\thistubeh{h}}{t}$; and Part~\ref{sublinear_tube_geometry:growth_part} of Lemma~\ref{sublinear_tube_geometry_lem} implies that if $(1-\spp)t \le (1-\thph{h})h$, then $\rreacheddset[2]{\mu}{\thisz}{\thistube}{(1-\spp)t} = \rreacheddset[2]{\mu}{\thisz}{\thistubeh{h}}{(1-\spp)t}$. Therefore, for any $h>0$,
\begin{align*}
\bighevent{h}
&= \eventthat[3]{ \crreachedset{\tmeasure}{\thisv}{\thistubeh{h}}{t}
\supseteq \rreacheddset[2]{\mu}{\thisz}{\thistubeh{h}}{(1-\spp)t}
\text{ for all } t\ge \alpha h}\\
&\subseteq \eventthat[3]{ \crreachedset{\tmeasure}{\thisv}{\thistube}{t}
\supseteq \rreacheddset[2]{\mu}{\thisz}{\thistube}{(1-\spp)t}
\text{ for all } t\in \segment[1]{\alpha h}{\tfrac{1-\thph{h}}{1-\spp}h}}.
\end{align*}
Thus, to complete the proof of Claim~\ref{tube_shape:event_containment_claim}, we need only show that
\begin{equation*}
\alpha (h+1) \le \frac{1-\thph{h}}{1-\spp}h
\quad\text{for all sufficiently large $h$}.
\end{equation*}
To see this, first note that $h+1\le \frac{1}{1-\spp} h$ for all $h\ge \spp^{-1} -1$, and next, since $\thplim<1$ and $\thph{h}\searrow \thplim$ as $h\to\infty$, there is some $h_0<\infty$ such that $\thph{h} \le \frac{1+\thplim}{2}$ for all $h\ge h_0$. Then $1-\thph{h} \ge \frac{1-\thplim}{2} =\alpha$ for all $h\ge h_0$, and hence
\[
h\ge h_0\join \parens{\spp^{-1} -1} \implies
\alpha (h+1)
\le \alpha \cdot \frac{1}{1-\spp}h
\le \frac{1-\thph{h}}{1-\spp}h,
\]
so the claim holds with $N \definedas \ceil[1]{h_0\join \parens{\spp^{-1} -1}}$.
%
%
\end{proof}

Now for $h\ge 0$ define
\[
\badhevent{h} \definedas \eventthat[3]{
\exists t \in \segment[2]{\alpha h}{\alpha (h+1)} \text{ such that } 
\rreacheddset[2]{\mu}{\thisz}{\thistube}{(1-\spp)t}
\not\subseteq \crreachedset{\tmeasure}{\thisv}{\thistube}{t}},
\]
and note that since $\alpha>0$, if we consider the subcollection $\setof{\badhevent{n}}_{n\in\N}$, then
\begin{equation}
\label{tube_shape:badevent_io_eqn}
\badevent = \setof{\badhevent{n} \text{ i.o.}}.
\end{equation}
By Claim~\ref{tube_shape:event_containment_claim} and \eqref{tube_shape:badevent_prob_eqn}, for any $h> \firstn$ we have $\badhevent{h} \subseteq \eventcomplement{\bighevent{h}}$ and hence
\begin{align*}
\Pr \parens[2]{\badhevent{h}}
&\le \thph{h}^{-d} \thisC e^{-\thisc \thph{h} \alpha h}
=  \thisC \exp \brackets[1]{-\thisc \alpha\, \rfcn(h) + d \log \parens[1]{\tfrac{h}{\rfcn(h)}}}.
\end{align*}
Since $\frac{\rfcn(h)}{\log h} \to\infty$ as $h\to\infty$ by assumption, there is some $\secondn \in \N$ such that $\rfcn(h)\ge 1\join \parens[1]{\frac{d+2}{\thisc \alpha} \log h}$ for all $h> \secondn$. Thus, for $h>\firstn\join \secondn$ we have
\begin{align*}
\Pr \parens[2]{\badhevent{h}}
&\le \thisC \exp \brackets[1]{-(d+2)\log h + d \log \parens[1]{\tfrac{h}{1}}}
= \thisC h^{-2}.
\end{align*}
Therefore, setting $\maxn =\firstn\join \secondn$,
\[
\sum_{n=0}^\infty \Pr \parens[2]{\badhevent{n}}
\le \sum_{n=0}^{N''} \Pr \parens[2]{\badhevent{n}}
	+ \sum_{n=N''+1}^\infty \thisC n^{-2} <\infty,
\]
so by the First Borel-Cantelli Lemma, we have $\Pr \setof{\badhevent{n} \text{ i.o.}} = 0$ and hence $\Pr \parens[2]{\badevent} =  0$ by \eqref{tube_shape:badevent_io_eqn}.
\end{proof}

}

%


\chapter{Stochastic Analogues of Deterministic Competition}
\label{random_fpc_chap}

In this chapter we combine the results of Chapters~\ref{fpc_basics_chap}, \ref{deterministic_chap}, and \ref{cone_growth_chap} to study the random two-type competition process in $\Zd$. Our primary goal will be to prove various results showing that on large scales, the random $\Zd$-process behaves similarly to the deterministic $\Rd$-process run with a traversal norm pair consisting of the two species' shape functions. In particular, after developing some general tools for analyzing the random two-type process, we will focus in Section~\ref{random_fpc:cone_competition_sec} on the same ``point vs.\ conical shell" configurations studied for the deterministic process in Section~\ref{deterministic:cone_competition_sec}. Our main results can be seen as extending Deijfen and \Haggstrom's Theorem~\ref{dh_unbounded_survival_thm} \cite[Theorem~1.1]{Deijfen:2007aa}, and we will use the results of Chapter~\ref{random_fpc_chap} to analyze competition from finite starting configurations in Chapter~\ref{coex_finite_chap}.

Chapter~\ref{random_fpc_chap} is organized as follows. In Section~\ref{random_fpc:continuum_process_sec} we introduce notation to facilitate treating the random two-type process as a subset of $\Rd$, similar to the notation introduced in Chapter~\ref{cone_growth_chap} for the one-type process. In Section~\ref{random_fpc:conq_upper_bds_sec} we prove a simple conquering property for the two-type process, as a corollary of Proposition~\ref{conquering_property_prop}, and then we prove some probabilistic upper bounds on the growth of a one-type process, to complement the lower bounds from Chapter~\ref{cone_growth_chap}. In Section~\ref{random_fpc:segment_conquering_sec} we develop the main technical tools needed for analyzing the two-type process. The main general result is Theorem~\ref{cone_segment_conquering_thm}, which is a stochastic analogue of Proposition~\ref{determ_star_conquering_prop}, showing that if species~1 starts at a point and species~2's initial set is far away, then with high probability species~1 conquers most of its $\beta$-Voronoi star-cell for any fixed $\beta>1$. In Section~\ref{random_fpc:cone_competition_sec} we study the random two-type process when species~1 starts at a point and species~2 starts on the exterior of a cone. The results of Section~\ref{random_fpc:cone_competition_sec} are stochastic analogues of the results from Section~\ref{deterministic:cones:conquered_regions_sec}, showing that species~1 can survive with high probability in wide cones but loses almost surely in narrow cones. In Section~\ref{random_fpc:additional_results_sec} we prove some additional results about the two-type competition process as corollaries of the results in Section~\ref{random_fpc:cone_competition_sec}, and we state a result of \Haggstrom and Pemantle, Proposition~\ref{hp_point_ball_prop} \cite[Proposition~2.2]{Haggstrom:2000aa}, to help put our results in context.

\section{Extending the Random Two-Type Process to $\Rd$}
\label{random_fpc:continuum_process_sec}

In this chapter we study a random two-type process evolving according to some (sufficiently nice) random two-type traversal measure $\tmeasurepair$ on $\edges{\Zd}$ whose components $\tmeasure_1$ and $\tmeasure_2$ have corresponding shape functions $\mu_1$ and $\mu_2$. It is most convenient to describe the process as a subset of $\Rd$ rather than $\Zd$. In order to do this, we need to extend our definitions for the random two-type process analogously to the definitions for the continuum one-type process in Section~\ref{cone_growth:cov_time_sec}.

For the random two-type process, an \textdef{initial configuration $\initconfig$ in $\Rd$} means a pair of subsets of $\Rd$ such that $\lat{\initialset{1}}\cap \lat{\initialset{2}} = \emptyset$, and we define the two-type process with such an initial configuration using lattice approximation as usual:
\begin{equation}
\label{continuum_entangled_proc_def}
\reachedset{\tmeasurepair}{\initconfig}{t} \definedas
\reachedset{\tmeasurepair}{ \pair[1]{\lat{\initialset{1}}}{\lat{\initialset{2}}} }{t}
\quad\text{for all } t\in [0,\infty].
\end{equation}
We also use lattice approximation to extend the definitions of the entangled passage times and the occupied sets: For any $B\subseteq \Rd$ and $i\in \setof{1,2}$, we define, respectively, the \textdef{entangled passage time} to $B$ and the \textdef{entangled covering time} of $B$ for species $i$ as
\begin{equation}
\label{continuum_entangled_times_def_eqn}
\entangledtime{i}{\initconfig}{B} \definedas
	\inf_{\vect v\in \lat B} \entangledtime{i}{ \pair[1]{\lat{\initialset{1}}}{\lat{\initialset{2}}} }{\vect v}
\quad\text{and}\quad
\entangledctime{i}{\initconfig}{B} \definedas
	\sup_{\vect v\in \lat B} \entangledtime{i}{ \pair[1]{\lat{\initialset{1}}}{\lat{\initialset{2}}} }{\vect v},
\end{equation}
where the definition of the entangled passage times appearing on the right is given in Section~\ref{fpc_basics:additional_notation_sec}.
The continuum \textdef{hitting process} and \textdef{covering process} for each species are then defined, respectively, as
\begin{align}
\label{entangled_hitting_proc_def_eqn}
\rreachedset{i}{\initconfig}{\Rd}{t} &\definedas
	\setbuilder[1]{\vect x\in\Rd}{\entangledtime{i}{\initconfig}{\vect x}\le t}
	=\cubify[1]{\reachedset{i}{\initconfig}{t}},\\
\intertext{and}
\label{entangled_covering_proc_def_eqn}
\crreachedset{i}{\initconfig}{\Rd}{t} &\definedas
	\setbuilder[1]{\vect x\in\Rd}{\entangledctime{i}{\initconfig}{\vect x}\le t}
	=\interior{\cubify[1]{\reachedset{i}{\initconfig}{t}}}.
\end{align}
One can verify from the definitions (and Lemma~\ref{lattice_cube_properties_lem}) that
\begin{equation}
\label{entangled_covering_proc_equiv_eqn}
\crreachedset{i}{\initconfig}{\Rd}{t}
= \crreachedset{\tmeasure_i}{\initialset{i}}{\intcubify{\finalset{i}}}{t}, 
\end{equation}
where $\finalset{i} = \reachedset{i}{\initconfig}{\infty}$ is species~$i$'s finally conquered set. The geometric interpretation of the definitions \eqref{entangled_hitting_proc_def_eqn} and \eqref{entangled_covering_proc_def_eqn} is the same as for the corresponding one-type processes, as described in Chapter~\ref{cone_growth_chap}.

\section{A Simple Conquering Property and Upper Bounds on Growth}
\label{random_fpc:conq_upper_bds_sec}

%

Our first result follows directly from the basic conquering property in Proposition~\ref{conquering_property_prop} and generalizes one of the main ideas of Deijfen and \Haggstrom in the proof of \cite[Proposition~3.1]{\DHunbounded}, which considered the case where $S$ is a half-cylinder in the lattice.\footnote{%
In fact, it was the proof of \cite[Proposition~3.1]{Deijfen:2007aa} that inspired me to attempt to find a generalized version of the statement in Lemma~\ref{sufficient_conquering_lem} and eventually led me to the conquering property in Proposition~\ref{conquering_property_prop} and the subsequent development of Definition~\ref{two_type_proc_def} for the two-type process.
} 
Lemma~\ref{sufficient_conquering_lem} will be our primary tool for showing that one species conquers a given set in the two-type process.

\begin{lem}[Sufficient condition for conquering a subset of $\Rd$ in the two-type process]
\label{sufficient_conquering_lem}
Let $\initconfig$ be an initial configuration in $\Rd$, and let $S\subseteq\Rd$. If on some realization of $\tmeasurepair$ it holds that
\[
\ctime{1}{S}{\initialset{1}}{\vect x} < \ptimenew{2}{\initialset{2}}{\vect x}
\quad\text{for all } \vect x\in S,
\]
then species 1 conquers $S$ in the two-type covering process $\crreachedset{\tmeasurepair}{\initconfig}{\Rd}{t}$,
and moreover, 
\[
\crreachedset{1}{\initconfig}{\Rd}{t} \supseteq \crreachedset{\tmeasure_1}{\initialset{1}}{S}{t}
\quad\text{for all $t\ge 0$.}
\]
A symmetric statement holds with the roles of species 1 and 2 reversed.
\end{lem}

\begin{proof}
%
%
Going back to the definitions, the stated hypothesis means
\[
\sup_{\vect v\in \lat{\vect x}} \rptime[2]{1}{\lat S}{\lat{\initialset{1}}}{\vect v}
< \inf_{\vect v\in \lat{\vect x}} \ptimenew[2]{2}{\lat{\initialset{2}}}{\vect v}
\quad\text{for all } \vect x\in S.
\]
Now for any $\vect v\in \lat{S}$ there is some $\vect x\in S$ such that $\vect v\in \lat{\vect x}$, so this implies that
\begin{equation}\label{sufficient_conquering:lattice_times_eqn}
\rptime[2]{1}{\lat S}{\lat{\initialset{1}}}{\vect v} < \ptimenew[2]{2}{\lat{\initialset{2}}}{\vect v}
\quad\text{for all } \vect v\in \lat S.
\end{equation}
Notice that \eqref{sufficient_conquering:lattice_times_eqn} implies that $\rptime[2]{1}{\lat S}{\lat{\initialset{1}}}{\vect v}<\infty$ and $\ptimenew[2]{2}{\lat{\initialset{2}}}{\vect v}>0$ for all $\vect v\in \lat S$, so $\lat{\initialset{1}}\cap \lat S \ne\emptyset$ and $\lat S\cap \lat{\initialset{2}} = \emptyset$. Now, since we trivially have $\ptimenew[2]{2}{\lat{\initialset{2}}}{\vect v} \le \starptimenew[2]{2}{(\Zd\setminus \lat S)}{\lat{\initialset{2}}}{\vect v}$ for any $\vect v$, \eqref{sufficient_conquering:lattice_times_eqn} implies that 
\[
\rptime[2]{1}{\lat S}{\lat{\initialset{1}}}{\vect v} <
\starptimenew[2]{2}{(\Zd\setminus \lat S)}{\lat{\initialset{2}}}{\vect v}
\quad\text{for all } \vect v\in \bd \lat{S} \setminus \lat{\initialset{1}},
\]
and since $\lat S\cap \lat{\initialset{2}} = \emptyset$, Proposition~\ref{conquering_property_prop} implies that species 1 conquers $\lat S$ in the lattice process, i.e.\ $\lat S \subseteq \reachedset{1}{\initconfig}{\infty}$. 
Lemma~\ref{lattice_cube_properties_lem} then implies that
\[
S\subseteq
\intcubify[1]{\reachedset{1}{\initconfig}{\infty}}=
\crreachedset{1}{\initconfig}{\Rd}{\infty},
\]
so species 1 conquers $S$ in the covering process. Moreover, Lemma~\ref{conquered_set_growth_lem} implies that for all $t\in [0,\infty]$, we have $\reachedset{1}{\initconfig}{t} \supseteq \rreachedset{\tmeasure_1}{\lat{\initialset{1}}}{\lat S}{t}$ and hence 
\[
\crreachedset{1}{\initconfig}{\Rd}{t} =
\interior{\cubify[1]{\reachedset{1}{\initconfig}{t}}} \supseteq
\interior{\cubify[1]{\rreachedset{\tmeasure_1}{\lat{\initialset{1}}}{\lat S}{t}}} \supseteq
S\cap \interior{\cubify[1]{\rreachedset{\tmeasure_1}{\lat{\initialset{1}}}{\lat S}{t}}}
= \crreachedset{\tmeasure_1}{\initialset{1}}{S}{t}.
\qedhere
\]
%
%
%
%
\end{proof}

Lemma~\ref{sufficient_conquering_lem} is useful because it reduces a statement about the entangled process $\crreachedfcn{\tmeasurepair}{\initconfig}{\Rd}$ to a statement about the two disentangled processes $\crreachedfcn{\tmeasure_1}{\initialset{1}}{S}$ and $\rreachedfcn{\tmeasure_2}{\initialset{2}}{\Rd}$. Chapter~\ref{cone_growth_chap} provides the necessary lower bounds on the growth of the restricted process $\crreachedfcn{\tmeasure_1}{\initialset{1}}{S}$. The next three lemmas provide the relevant upper bounds on the growth of an unrestricted one-type process $\rreachedfcn{\tmeasure}{A}{\Rd}$; in the next section, we will apply these results with $\tmeasure = \tmeasure_2$ and $A = \initialset{2}$. 


If we want to show that species 1 conquers some set $S\subseteq\Rd$ using Lemma~\ref{sufficient_conquering_lem}, we need a lower bound on $\ptimenew{2}{\initialset{2}}{\vect x}$ for $\vect x\in S$. Instead of considering the disentangled type-2 process $\rreachedfcn{\tmeasure_2}{\initialset{2}}{\Rd}$ spreading out from $\initialset{2}$, it is convenient to make use of the symmetry of $\ptmetric[2]$ and consider a dual%
\footnote{
Two stochastic processes $X_t,Y_t$ defined on some interval $[0,T]$ and taking values in $\mathcal{X},\mathcal{Y}$, respectively, are said to be \textdef{dual} with respect to the function $H\colon \mathcal{X}\times\mathcal{Y}\to \R$ if $\E H(X_t,Y_0) = \E H(X_0,Y_t)$ for all $t\in [0,T]$ (cf.\ \cite[p.~120]{Mueller:2009aa} or \cite[p.~115]{Liggett:2010aa}). Thus, for any $A,B\subseteq\Zd$, the processes $\reachedset{\tmeasure}{A}{t}$ and $\reachedset{\tmeasure}{B}{t}$ on $[0,\infty]$ are dual with respect to the indicator function $\ind{\eventthat{A'\cap B'\ne\emptyset}} \colon \powerset{\Zd}\times \powerset{\Zd}\to \Rplus$ since $\mbox{\ensuremath{\Pr \eventthat[2]{\reachedset{\tmeasure}{A}{t} \cap B\ne \emptyset}}} = \Pr \eventthat[2]{\reachedset{\tmeasure}{B}{t} \cap A\ne \emptyset}$ for all $t\in [0,\infty]$.
}
process $\rreachedfcn{\tmeasure_2}{S}{\Rd}$ starting at points in $S$. This dual perspective, formalized in Lemma~\ref{dual_times_lem} below, is useful because for a given $\vect x\in S$, the dual process $\rreachedfcn{\tmeasure_2}{\vect x}{\Rd}$ originating from the single point $\vect x$ is easier to deal with than the more complicated process $\rreachedfcn{\tmeasure_2}{\initialset{2}}{\Rd}$ starting from multiple points.
This idea for ``dualizing" the process appears in Deijfen and \Haggstrom's proof of \cite[Proposition~3.1]{Deijfen:2007aa}, and a similar idea is explicitly referred to as a ``duality principle" by Kordzakhia and Lalley \cite[p.~8]{Kordzakhia:2005aa}.
Lemma~\ref{dual_times_lem} can be viewed as a one-sided analogue or variant of Lemma~\ref{closer_char_lem}.


\begin{lem}[Dual growth principle for passage time from a set to a point]
\label{dual_times_lem}
Fix a realization $\tmeasure$ of a one-type process and a norm $\mu$ on $\Rd$, and let $A\subseteq \Rd$, $\vect x\in\Rd$, and $\alpha >0$. Then
\[
``\rreachedset{\tmeasure}{\vect x}{\Rd}{\alpha t}\subseteq
\reacheddset{\mu}{\vect x}{t-} \text{ for all } t\ge \distnew{\mu}{A}{\vect x}"
\implies
\ptimenew{\tmeasure}{A}{\vect x} > \alpha \distnew{\mu}{A}{\vect x}.
\]
\end{lem}

\begin{proof}
We will deduce this result from the inversion formulas. First note that we must have $\distnew{\mu}{A}{\vect x} = \distnew{\mu}{\vect y}{\vect x}$ for some $\vect y\in \closure{A}$. Thus,
\begin{align*}
\eventthat[2]{\ptimenew{\tmeasure}{A}{\vect x} > \alpha \distnew{\mu}{A}{\vect x}}
&\supseteq \eventthat[2]{\,\ptimenew{\tmeasure}{A}{\vect x} 
	> \alpha \distnew{\mu}{\vect y}{\vect x}
	\; \forall \vect y\in \closure{A}\,}\\
&\supseteq \eventthat[1]{\textstack{$\ptimenew{\tmeasure}{A}{\vect x} 
	> \alpha \distnew{\mu}{\vect y}{\vect x}$ for all $\vect y\in\Rd$}
	{with $\distnew{\mu}{\vect y}{\vect x} \ge \distnew{\mu}{A}{\vect x}$}}\\
&= \eventthat[1]{\alpha \distnew{\mu}{A}{\vect x} \join \ptimenew{\tmeasure}{\vect y}{\vect x}
	> \alpha \distnew{\mu}{\vect y}{\vect x} \;\forall \vect y\in\Rd}\\
&= \eventthat[1]{\rreachedset{\tmeasure}{\vect x}{\Rd}{\alpha t}
	\subseteq \reacheddset{\mu}{\vect x}{t-} \;\forall t\ge \distnew{\mu}{A}{\vect x} },
\end{align*}
where the last line follows from the inversion formula \eqref{inversion_formulas:strong_ubg} in Lemma~\ref{inversion_formula_lem}.
\end{proof}

Note that Lemma~\ref{dual_times_lem} is really a deterministic statement, holding pointwise on the canonical sample space $\onetypespace$; in order to obtain probability estimates from it, we will use Lemma~\ref{shifted_dev_bounds_lem} to get a lower bound on the probability of the event ``$\rreachedset{\tmeasure}{\vect x}{\Rd}{\alpha t}\subseteq \reacheddset{\mu}{\vect x}{t-}$ for all $ t\ge \distnew{\mu}{A}{\vect x}$."
For example, in the following simple lemma we combine Lemma~\ref{dual_times_lem} with Lemma~\ref{shifted_dev_bounds_lem} to get a probabilistic lower bound on the passage time between two subsets of $\Rd$ when one of the sets is bounded.

\begin{lem}[Passage time to a bounded set]
\label{time_to_finite_set_lem}
Let $\tmeasure$ be a random one-type traversal measure on $\Zd$ satisfying \finitespeed{d} and \expmoment, with shape function $\mu$ (a norm on $\Rd$). Given $\alpha\in (0,1)$, there exist positive constants $\bigcref{time_to_finite_set_lem}$ and $\smallcref{time_to_finite_set_lem}$ such that if $A$ and $B$ are nonempty subsets of $\Rd$, and $B$ is bounded, then
\begin{align*}
\Pr \eventthat[2]{\ptimenew{\tmeasure}{A}{B} > \alpha \distnew{\mu}{A}{B}} 
&\ge 1- \card[1]{\lat{B}}\cdot \bigcref{time_to_finite_set_lem} e^{-\smallcref{time_to_finite_set_lem} \distnew{\mu}{A}{B}}.
\end{align*}
\end{lem}

\begin{proof}
\newcommand{\thisv}{\vect v}
\newcommand{\thisx}{\vect x}
\newcommand{\thiscube}{\cubify{\thisv}}

\newcommand{\bound}{b}

\newcommand{\thisdil}{\smallkonstsymb}

\newcommand{\distab}{\distnew{\mu}{A}{B}}
\newcommand{\distav}{\distnew{\mu}{A}{\thisv}}
\newcommand{\distax}{\distnew{\mu}{A}{\thisx}}

\newcommand{\ptab}{\ptimenew{\tmeasure}{A}{B}}
\newcommand{\ptav}{\ptimenew{\tmeasure}{A}{\thisv}}
\newcommand{\ptax}{\ptimenew{\tmeasure}{A}{\thisx}}

\newcommand{\thisxreached}{\rreachedset{\tmeasure}{\thisx}{\Rd}{\alpha t}}
\newcommand{\thisvreached}{\rreachedset{\tmeasure}{\thisv}{\Rd}{\alpha t}}

\newcommand{\thisxreachedd}{\reacheddset{\mu}{\thisx}{t-}}

\newcommand{\ballevent}{\eventref{shifted_dev_bounds_lem}^{\thisv} \parens[2]{\alpha[\distab -\thisdil];\alpha^{-1}-1}}

\newcommand{\oldC}{\bigcref{shifted_dev_bounds_lem}\parens{\alpha^{-1}-1}}
\newcommand{\oldc}{\smallcref{shifted_dev_bounds_lem}\parens{\alpha^{-1}-1}}

\newcommand{\thisC}{\bigcref{time_to_finite_set_lem}(\alpha)}
\newcommand{\thisc}{\smallcref{time_to_finite_set_lem}(\alpha)}

Recall that $\ptab = \ptimenew{\tmeasure}{A}{\lat B}$ by definition. Since $B$ is bounded, $\lat B$ is finite, and thus
\begin{align}
\label{time_to_finite_set:strict_inequality_eqn}
\eventthat[2]{\ptab > \alpha \distab}
&=\eventthat[1]{\min_{\thisv\in \lat B} \ptav > \alpha \distab} \notag\\
&=\eventthat[1]{\inf_{\thisx\in B} \ptimenew{\tmeasure}{A}{\lat \thisx} > \alpha \distab}.
\end{align}
Next, note that if $\thisdil \definedas \dilnorm{\linfty{d}}{\mu}$, then Lemma~\ref{lattice_cube_properties_lem} implies that
\begin{equation}
\label{time_to_finite_set:close_points_eqn}
\distnew{\mu}{B}{\thisx}\le \thisdil
\text{\quad for all } \thisx\in \cubify[2]{\lat B}.
\end{equation}
Now, using \eqref{time_to_finite_set:strict_inequality_eqn} in the first line, Lemma~\ref{dual_times_lem} in the third line, and \eqref{time_to_finite_set:close_points_eqn} in the fifth line, we have
\begin{align*}
\eventthat[2]{\ptab > \alpha \distab}
&= \bigcap_{\thisx\in B} \eventthat[2]{\ptax > \alpha \distab}\\
&\supseteq \bigcap_{\thisx\in B} \eventthat[2]{\ptax > \alpha \distax}\\
&\supseteq  \bigcap_{\thisx\in B}
	\eventthat[1]{\thisxreached \subseteq \thisxreachedd \; \forall t\ge \distax}\\
&\supseteq  \bigcap_{\thisv\in \lat B} \bigcap_{\thisx\in \thiscube}
	\eventthat[1]{\thisvreached \subseteq \thisxreachedd  \; \forall t\ge \distax}\\
&\supseteq  \bigcap_{\thisv\in \lat B} \bigcap_{\thisx\in \thiscube}
	\eventthat[1]{\thisvreached \subseteq \thisxreachedd  \; \forall t\ge \distab -\thisdil}\\
&\supseteq \bigcap_{\thisv\in \lat B} \ballevent.
\end{align*}
Therefore, since $\alpha<1$, Lemma~\ref{shifted_dev_bounds_lem} implies that
\begin{align*}
\Pr \eventthat[2]{\ptab \le \alpha \distab}
&\le \sum_{\thisv\in\lat B} \Pr \parens[1]{\eventcomplement{\ballevent}}\\
&\le \sum_{\thisv\in\lat B} \oldC e^{-\oldc \alpha[\distab -\thisdil]}\\
&\le \card{\lat B}\cdot \thisC e^{-\thisc \distab},
\end{align*}
where
\[
\thisc \definedas \alpha\cdot \oldc
\quad\text{and}\quad
\thisC \definedas \oldC e^{\thisc \thisdil}.
\qedhere
\]
%
%
\end{proof}

Using Lemma~\ref{time_to_finite_set_lem}, we obtain the following lower bound on the passage time between two concentric $\mu$-spheres. In the style Chapter~\ref{cone_growth_chap}'s events, we prove the bound simultaneously for all pairs of spheres with common center in a uint cube centered at $\vect v\in\Zd$, in order to simplify later proofs. We will use Lemma~\ref{time_between_spheres_lem} in the proof of Lemma~\ref{ball_conquering_lem} below, which will be our first nontrivial result about the two-type process. Recall from Section~\ref{intro:euc_lat_geo_sec} that for $A\subseteq\Rd$ we define $\shellof{A} \definedas \Rd\setminus \interior{A}$.


\begin{lem}[Passage time between nested $\mu$-spheres]
\label{time_between_spheres_lem}

\newcommand{\thisz}{\vect z}
\newcommand{\thisvertex}{\vect v}
\newcommand{\thiscube}{\cubify{\thisvertex}}

\newcommand{\bigrad}{R}
\newcommand{\smallrad}{r}

\newcommand{\bigball}{\reacheddset{\mu}{\thisz}{\bigrad}}
\newcommand{\smallball}{\reacheddset{\mu}{\thisz}{\smallrad}}

\newcommand{\thisA}{\shellof{\bigball}}

\newcommand{\thisC}{\bigcref{time_between_spheres_lem}}
\newcommand{\thisc}{\smallcref{time_between_spheres_lem}}

Let $\tmeasure$ and $\mu$ be as in Lemma~\ref{time_to_finite_set_lem}. Given $\alpha<1$, there exist positive constants $\thisC$ and $\thisc$ such that for any $\thisvertex\in\Zd$ and any $0\le \smallrad \le \bigrad$,
\[
\Pr \eventthat[3]{ \forall \thisz\in\thiscube,\;
\ptimenew[2]{\tmeasure}{\thisA}{\,\smallball} > \alpha (\bigrad-\smallrad)}
\ge 1-\thisC \cdot (\smallrad\join 1)^d e^{-\thisc \cdot (\bigrad-\smallrad)},
\]
where $\shellof{B} \definedas \Rd\setminus \interior{B}$ for $B\subseteq\Rd$.
\end{lem}

\begin{proof}

\newcommand{\thisz}{\vect z}
\newcommand{\thisvertex}{\vect v}
\newcommand{\thiscube}{\cubify{\thisvertex}}

\newcommand{\bigrad}{R}
\newcommand{\smallrad}{r}
\newcommand{\rmin}{M_{\alpha}}

\newcommand{\bigball}{\reacheddset{\mu}{\thisz}{\bigrad}}
\newcommand{\smallball}{\reacheddset{\mu}{\thisz}{\smallrad}}
\newcommand{\vball}[2][0]{\reacheddset[#1]{\mu}{\thisvertex}{#2}}

\newcommand{\thisA}{\shellof{\bigball}}
\newcommand{\newA}{\shellof{\vball{\bigrad-\thisdil}}}
\newcommand{\newB}{\vball{\smallrad+\thisdil}}

\newcommand{\thisC}{\bigcref{time_between_spheres_lem}}
\newcommand{\thisc}{\smallcref{time_between_spheres_lem}}
\newcommand{\thisevent}{\eventref{time_between_spheres_lem}^{\thisvertex} \parens{r,R;\alpha}}

\newcommand{\oldC}{\bigcref{time_to_finite_set_lem}(\alpha')}
\newcommand{\oldc}{\smallcref{time_to_finite_set_lem}(\alpha')}

\newcommand{\thisdil}{\smallkonstsymb}

The proof idea is similar to that of Lemma~\ref{shifted_dev_bounds_lem}. Let $\thisevent$ denote the event in the statement of the lemma, and set $\thisdil \definedas \frac{1}{2} \dilnorm{\linfty{d}}{\mu}$. Then
\[
\bigcup_{\thisz\in\thiscube} \shellof{\bigball}
\subseteq \shellof{\vball{\bigrad-\thisdil}}
\quad\text{and}\quad
\bigcup_{\thisz\in\thiscube} \smallball \subseteq \vball{\smallrad+\thisdil},
\]
and therefore,
\begin{align*}
\thisevent
&= 
	\eventthat[3]{\forall \thisz\in\thiscube,\;
	\ptimenew[2]{\tmeasure}{\thisA}{\,\smallball} > \alpha (\bigrad-\smallrad)}\\
&\supseteq \eventthat[3]{\ptimenew[2]{\tmeasure}{\newA}{\,\newB}
	> \alpha (\bigrad-\smallrad)}.
\end{align*}
Now note that
\begin{equation}
\label{time_between_spheres:sphere_dist_eqn}
\distnew[2]{\mu}{\newA}{\newB} = \bigrad-\smallrad - 2\thisdil.
\end{equation}
Set $\alpha' \definedas \frac{1+\alpha}{2}<1$ and $\rmin \definedas \frac{2(1+\alpha)\thisdil}{1-\alpha}<\infty$. Then
\[
\bigrad-\smallrad \ge \rmin \implies
\alpha' (\bigrad-\smallrad - 2\thisdil)\ge \alpha (\bigrad-\smallrad),
\]
so if $\bigrad-\smallrad\ge \rmin$, then Lemma~\ref{time_to_finite_set_lem} and \eqref{time_between_spheres:sphere_dist_eqn} imply that
\begin{align*}
\Pr \parens[2]{\thisevent}
&\ge \Pr\eventthat[3]{\ptimenew[2]{\tmeasure}{\newA}{\,\newB}
	> \alpha (\bigrad-\smallrad)}\\
&\ge \Pr \eventthat[3]{\ptimenew[2]{\tmeasure}{\newA}{\,\newB
	 > \alpha' (\bigrad-\smallrad-2\thisdil)}}\\
&\ge 1-\card[1]{\lat[2]{\newB}} \cdot \oldC e^{-\oldc (\bigrad-\smallrad-2\thisdil)}.
\end{align*}
Since $\lat[2]{\newB}=\Zd\cap \cubify[2]{\newB} \subseteq \Zd\cap \vball{\smallrad+2\thisdil}$, Lemma~\ref{norm_shell_lem} implies that there is some constant $\mushellbound$ such that
\[
\card[1]{\lat[2]{\newB}} \le \mushellbound \cdot [(\smallrad+2\thisdil)\join1]^d
\le \mushellbound (1+2\thisdil)^d (\smallrad \join1)^d.
\]
Thus, setting
\[
\thisC(\alpha) \definedas
\brackets[1]{\mushellbound (1+2\thisdil)^d e^{\oldc \cdot 2\thisdil} \join e^{\oldc \cdot \rmin}}
\cdot \oldC
\quad\text{and}\quad
\thisc(\alpha) \definedas \oldc,
\]
we get $\Pr \parens[2]{\thisevent} \ge 1-\thisC(\alpha) (\smallrad\join 1)^d e^{-\thisc(\alpha) (\bigrad-\smallrad)}$ for all $0\le \smallrad\le \bigrad$.
\end{proof}

\section{Conquering a $\mu_1$-Star in the Two-Type Process}
\label{random_fpc:segment_conquering_sec}

Throughout this section we assume that $\tmeasurepair$ satisfies \finitespeed{d}, \expmoment, and \notiesprob, meaning that the component measures $\tmeasure_1$ and $\tmeasure_2$ each satisfy \finitespeed{d} and \expmoment, and the pair satisfies \notiesprob.

\subsection{General Results for Conquering Star-Shaped Sets}
\label{random_fpc:main_gen_thm_sec}

Recall from \eqref{rmeet_def_eqn} that for $\beta>0$, species 1's $\beta$-meeting radius for the traversal norm pair $\pair{\beta \mu_1}{\mu_2}$ was defined to be
\begin{equation}
\label{rmeet_reminder_eqn}
\rmeet{1}{\beta} = \frac{1}{\beta + \esupnorm{\mu_2/\mu_1}},
\end{equation}
and that $\rmeet{1}{\beta}$ is the radius of the largest $\mu_1$-ball conquered by species 1 starting from $\vect z\in\Rd$ when species 2 starts on the boundary of a unit $\mu_2$-ball centered at $\vect z$.
Scaling this starting configuration by $r>0$, Lemma~\ref{determ_bullseye_lem} showed that for the traversal norm pair
$\adjustednormpair{2}{\beta^{-1}} = \pair{\mu_1}{\beta^{-1}\mu_2}$,
\begin{equation}
\label{bullseye_config_eqn}
\reacheddset[2]{1}{\pair{\vect z}{\,\shellof{\reacheddset{\mu_2}{\vect z}{r}}}}{\rmeet{1}{\beta}r}
= \reacheddset[2]{\mu_1}{\vect z}{\rmeet{1}{\beta} r},
\end{equation}
where $\shellof{A} \definedas \Rd\setminus \interior{A}$ for $A\subseteq\Rd$. The following lemma gives a weakened version of \eqref{bullseye_config_eqn} for the random process, showing that with high probability in $r$, species 1 will conquer the ball that a $\beta$-slowed version of species 1 would conquer in the corresponding deterministic process.
Lemma~\ref{ball_conquering_lem} will be the starting point for proving the much more general Theorem~\ref{cone_segment_conquering_thm} below.



\begin{lem}
[Conquering a $\mu_1$-ball inside a $\mu_2$-spherical shell]
\label{ball_conquering_lem}


\newcommand{\thisz}{\vect z}
\newcommand{\thisvertex}{\vect v}
\newcommand{\thiscube}{\cubify{\thisvertex}}

\newcommand{\thisAone}{\thisvertex}
\newcommand{\thisAtwo}{\shellof{\reacheddset{\mu_2}{\thisz}{r}}}
\newcommand{\thisconfig}{\pair[1]{\thisAone}{\,\thisAtwo}}

\newcommand{\thisevent}{\eventref{ball_conquering_lem}^{\thisvertex} \parens{r;\beta}}

\newcommand{\thisC}{\bigcref{ball_conquering_lem}(\beta)}
\newcommand{\thisc}{\smallcref{ball_conquering_lem}(\beta)}

\newcommand{\thisdil}{\dilnorm{\linfty{d}}{\mu_2}}
\newcommand{\thisdilsymb}{\smallkonstsymb}

Suppose $\tmeasurepair$ satisfies \finitespeed{d}, \expmoment, and \notiesprob, and let $\mu_1$ and $\mu_2$ be the shape functions for the component measures $\tmeasure_1$ and $\tmeasure_1$.
Let $\thisdilsymb \definedas \thisdil$, and for any $\thisvertex\in\Zd$, $r> \thisdilsymb$, and $\beta >1$, define the event
\[
\thisevent \definedas
\eventthat[1]{\text{For all } \thisz\in\thiscube,\;
\crreachedset[2]{1}{\thisconfig}{\Rd}{\rmeet{1}{1}r}
\supseteq \reacheddset[2]{\mu_1}{\thisz}{\rmeet{1}{\beta}r}},
\]
where $\rmeet{1}{1}$ and $\rmeet{1}{\beta}$ are defined by \eqref{rmeet_reminder_eqn}. Then given $\beta >1$, there are positive constants $\thisC$ and $\thisc$ such that for all $\thisvertex\in\Zd$ and $r> \thisdilsymb$,
\[
\Pr \parens[2]{\thisevent}
\ge 1-\thisC e^{-\thisc r}.
\]
\end{lem}

\begin{proof}


\newcommand{\thisx}{\vect x}
\newcommand{\thisz}{\vect z}
\newcommand{\thisvertex}{\vect v}
\newcommand{\thiscube}{\cubify{\thisvertex}}

\newcommand{\thisAone}{\thisvertex}
\newcommand{\thisAtwo}{\shellof{\reacheddset{\mu_2}{\thisz}{r}}}
\newcommand{\thisconfig}{\pair[1]{\thisAone}{\,\thisAtwo}}

\newcommand{\thisevent}{\eventref{ball_conquering_lem}^{\thisvertex} \parens{r;\beta}}
\newcommand{\strongballevent}{\widetilde{\eventsymb}_{\ref{ball_conquering_lem}}^{\thisvertex} \parens{r;\beta}}

\newcommand{\thisC}{\bigcref{ball_conquering_lem}(\beta)}
\newcommand{\thisc}{\smallcref{ball_conquering_lem}(\beta)}

\newcommand{\thisdil}{\dilnorm{\linfty{d}}{\mu_2}}
\newcommand{\thisdilsymb}{\smallkonstsymb}

\newcommand{\eone}{\spp}
\newcommand{\etwo}{\spp_2}

\newcommand{\ballC}{\bigcref{ball_covering_lem}(\eone)}
\newcommand{\ballc}{\smallcref{ball_covering_lem}(\eone)}
\newcommand{\bddC}{\bigcref{time_to_finite_set_lem} (\alpha)}
\newcommand{\bddc}{\smallcref{time_to_finite_set_lem} (\alpha)}
\newcommand{\sphereC}{\bigcref{time_between_spheres_lem} (\alpha)}
\newcommand{\spherec}{\smallcref{time_between_spheres_lem} (\alpha)}

\newcommand{\bigrad}{\rmeet{1}{1} r}
\newcommand{\smallrad}{\rmeet{1}{\beta} r}
\newcommand{\smalltworad}[1][0]{%
\ifcase#1
	\smallrad \esupnorm[2]{\tfrac{\mu_2}{\mu_1}}		
	\or \smallrad \esupnorm[1]{\frac{\mu_2}{\mu_1}}		
	\else \dilnorm{\mu_1}{\mu_2}\cdot \smallrad		
\fi}

\newcommand{\smallball}[1][0]{\reacheddset[#1]{\mu_1}{\thisz}{\rmeet{1}{\beta} r}}
\newcommand{\smalltwoball}[1][2]{\reacheddset[#1]{\mu_2}{\thisz}{\smalltworad}}

Note that the requirement $r>\thisdilsymb = \thisdil$ guarantees that the starting configuration $\thisconfig$ is valid, because if $\thisz\in\thiscube$ and $\distnew{\mu_2}{\vect y}{\thisz} \ge r>\thisdilsymb$, then $\distnew{\linfty{d}}{\vect y}{\thisvertex} > \frac{1}{2}$, showing that $\thisvertex \notin \lat[2]{\thisAtwo}$.

Our strategy is to combine Lemma~\ref{sufficient_conquering_lem} with Lemmas~\ref{ball_covering_lem} and \ref{time_between_spheres_lem} to show that species 1 conquers each of the balls $\smallball$ by time $\bigrad$ with high probability. First, we want an upper bound on the time it takes species 1 to internally cover the ball $\smallball$. Choose $\eone$ so that $(1+\eone) \smallrad = \bigrad$, i.e.
\[
\eone(\beta) \definedas \frac{\rmeet{1}{1}}{\rmeet{1}{\beta}} - 1
=\frac{\beta-1}{1+ \esupnorm{\mu_2/\mu_1}} >0.
\]
Then by Lemma~\ref{ball_covering_lem} (with $\tmeasure = \tmeasure_1$, $\mu=\mu_1$, and $r_0 = \smallrad$),
\begin{align}
\label{ball_conquering:cover_eqn}
1-\ballC e^{-\ballc \smallrad}
&\le \Pr \eventthat[1]{\forall \thisz\in\thiscube,\;
	\crreachedset[2]{\tmeasure_1}{\thisvertex}{\smallball}{\bigrad} = \smallball[2]} \notag\\
&= \Pr \eventthat[1]{\forall \thisz\in\thiscube,\;
	\ctime{1}{\smallball}{\thisvertex}{\thisx} \le \bigrad
	\text{ for all } \thisx\in \smallball[2]}.
\end{align}
Next we want a lower bound on the time it takes species 2 reach the ball $\smallball$, starting from $\thisAtwo$. Observe that 
\begin{equation}
\label{ball_conquering:ball_radii_eqn}
\smallball \subseteq \smalltwoball,
\quad\text{and}\quad
r- \smalltworad[1]
= \frac{\beta r}{\beta+\esupnorm{\mu_2/\mu_1}} = \beta \smallrad.
\end{equation}
Choose $\alpha$ so that $\alpha\cdot \beta \smallrad = \bigrad$, i.e.
\[
\alpha(\beta) \definedas \frac{\rmeet{1}{1}}{\beta \rmeet{1}{\beta}}
=\frac{1 + \beta^{-1} \esupnorm{\mu_2/\mu_1}}{1 + \esupnorm{\mu_2/\mu_1}}<1.
\]
Then by Lemma~\ref{time_between_spheres_lem} (with $\tmeasure = \tmeasure_2$ and $\mu=\mu_2$) and \eqref{ball_conquering:ball_radii_eqn},
\begin{align}
\label{ball_conquering:no_hit_eqn}
1- \sphereC &\parens[2]{\smalltworad[0] \join 1}^d e^{-\spherec \beta\smallrad} \notag\\
&\le \Pr \eventthat[3]{\forall \thisz\in\thiscube,\;
	\ptimenew[2]{2}{\thisAtwo}{\, \smalltwoball} > \alpha\cdot \beta \smallrad} \notag\\
&\le \Pr \eventthat[3]{\forall \thisz\in\thiscube,\;
	\ptimenew[2]{2}{\thisAtwo}{\, \smallball} > \bigrad}.
\end{align}
Now define the event
\begin{equation}
\label{ball_conquering:strong_event_eqn}
\strongballevent \definedas \bigcap_{\thisz\in\thiscube}
\eventthat[3]{\ctime{1}{\smallball}{\thisvertex}{\thisx} \le \bigrad
	< \ptimenew[2]{2}{\thisAtwo}{\thisx} \text{ for all } \thisx\in \smallball[2]}.
\end{equation}
Then $\strongballevent$ is the intersection of the two events in \eqref{ball_conquering:cover_eqn} and \eqref{ball_conquering:no_hit_eqn}, so we have
\[
\Pr \parens[2]{\strongballevent} \ge 1-\thisC e^{-\thisc r},
\]
where, using the bound \eqref{poly_geo_series:sup_bound_eqn} from the proof of Lemma~\ref{poly_geo_series_lem} for the polynomial$\times$exponential term in \eqref{ball_conquering:no_hit_eqn},
\begin{align*}
\thisC &\definedas
\ballC +   \brackets[1]{\frac{d\esupnorm{\mu_2/\mu_1}}{\spherec}}^d\cdot \sphereC,
\quad\text{and}\\
\thisc &\definedas \brackets[2]{\ballc \meet (\beta-1)\spherec}\cdot \rmeet{1}{\beta}.
\end{align*}
Finally, observe that $\thisevent$ occurs on $\strongballevent$ by Lemma~\ref{sufficient_conquering_lem}.
%
\end{proof}


The main result of this section will be Theorem~\ref{cone_segment_conquering_thm} below, which shows that if species 1 starts at the single vertex $\vect z$, then with high probability it conquers any union of sufficiently thick cone segments radiating from $\vect z$ and contained within the set of points that are at least $\beta$ times closer to $\vect z$ than to species 2's starting set, for any fixed $\beta>1$. Theorem~\ref{cone_segment_conquering_thm} is a very general result that extends Lemma~\ref{ball_conquering_lem} above and is a stochastic analogue of Proposition~\ref{determ_star_conquering_prop} for the deterministic process. We will use Theorem~\ref{cone_segment_conquering_thm} to prove Lemmas~\ref{bowling_pin_lem} and \ref{bowling_pin2_lem} below, which will be our main tools for analyzing the random process in the next section, where we consider the case when one species' starting set is the boundary of a cone. 

In order to prove Theorem~\ref{cone_segment_conquering_thm}, we will use the following lemma bounding the growth of species 2.  Lemma~\ref{growth_containment_lem} is similar in structure and proof to Lemma~\ref{space_covering_lem} in the previous chapter, so we have relegated its proof to Appendix~\ref{leftover_chap}. The difference between the two lemmas is that that Lemma~\ref{growth_containment_lem} provides an upper bound on growth whereas Lemma~\ref{space_covering_lem} provided a lower bound. The main probabilistic tool for proving Lemma~\ref{growth_containment_lem} is Lemma~\ref{shifted_dev_bounds_lem}, with the rest of the proof being geometric in nature.

\begin{lem}[Containment of growth within large $\mu_2$-balls far from the origin]
\label{growth_containment_lem}

\newcommand{\thisevent}{\eventref{growth_containment_lem}\parens{r,\beta;\alpha}}
\newcommand{\thisnormkonst}{k_{\mu_1}}

For any constants $\alpha\in (0,1)$, $\beta\ge 1$, and $r>0$, define the event
\[
\thisevent =
 \bigcap_{\vect z\in \cubify{\vect 0}} \bigcap_{\substack{\vect x\in \Rd \\ \mu_1(\vect x) \ge r}}
\eventthat[1]{ \rreachedset{\tmeasure_2}{\vect z+\vect x}{\Rd}{\alpha t} \subseteq
\reacheddset{\mu_2}{\vect z+\vect x}{t-} \;\forall t\ge \beta \mu_1(\vect x) },
\]
and let $\thisnormkonst = 3\dilnorm{\lspace{\infty}{d}}{\mu_1}+1$. Then given $\alpha \in (0,1)$, there exist positive constants $C_{\ref{growth_containment_lem}}$ and $c_{\ref{growth_containment_lem}}$ such that for any $\beta\ge 1$ and $r \ge 0$,
\[
\Pr \parens[2]{\thisevent} \ge 1- \bigcref{growth_containment_lem}(\alpha)
e^{-\smallcref{growth_containment_lem}(\alpha)\cdot \beta \cdot (r - \thisnormkonst)}.
\]
\end{lem}


To state Theorem~\ref{cone_segment_conquering_thm}, we make the following definition.
For $\beta,\thp>0$, $\vect z\in \Rd$, and $\initialset{2}\subseteq\Rd$, define \textdef{species 1's $\thp$-thick $\beta$-Voronoi star-cell} for the norm pair $\normpair = \pair{\mu_1}{\mu_2}$ and the initial configuration $\pair{\vect z}{\initialset{2}}$ to be
\begin{equation}
\label{fatcstarset_def_eqn}
\begin{split}
\fatcstarset{\normpair}{\beta}{\thp}{1}{\vect z}{\initialset{2}}
\definedas \bigcup &\setbuilder[1]{\cones\in \conefamily{\mu_1}{\vect z,\thp}}
	{\cones \subseteq \closerset{\normpair}{\beta}{1}{\vect z}{\initialset{2}}}\\
= \bigcup &\setbuilder[1]{S\in \fatstars{\mu_1}{\thp}{\vect z}}
	{S \subseteq \closerset{\normpair}{\beta}{1}{\vect z}{\initialset{2}}}.
\end{split}
\end{equation}
We will also use the following alternate notation to make some formulas more readable:
\begin{equation}
\label{fatcstar_def_eqn}
\fatcstar{\normpair}{\beta}{\thp}{1}{\vect z}{\initialset{2}} \definedas
\fatcstarset{\normpair}{\beta}{\thp}{1}{\vect z}{\initialset{2}}.
\end{equation}
Both forms of the notation are intended to be suggestive of previously defined objects. Namely, the star-shaped set $\fatcstar{\normpair}{\beta}{\thp}{1}{\vect z}{\initialset{2}} \equiv \fatcstarset{\normpair}{\beta}{\thp}{1}{\vect z}{\initialset{2}}$ is the largest element of $\fatstars{\mu_1}{\thp}{\vect z}$ that is contained in $\closerset{\normpair}{\beta}{1}{\vect z}{\initialset{2}}$; it is thus a subset and analogue of the $\beta$-Voronoi star-cell $\cstarset{\normpair}{\beta}{1}{\vect z}{\initialset{2}}$, defined using  $(\mu_1,\thp)$-cone segments originating at $\vect z$ instead of line segments originating at $\vect z$. It follows directly from the definition \eqref{fatcstarset_def_eqn} that for $S\subseteq\Rd$,
\begin{equation*}
S\in \fatstars{\mu_1}{\thp}{\vect z}\text{ and }
S\subseteq \closerset{\normpair}{\beta}{1}{\vect z}{\initialset{2}}
\implies S\subseteq \fatcstar{\normpair}{\beta}{\thp}{1}{\vect z}{\initialset{2}}.
\end{equation*}
Note that it follows from Lemma~\ref{mu_cone_decomp_lem} that the sets $\fatcstar{\normpair}{\beta}{\thp}{1}{\vect z}{\initialset{2}}$ are decreasing with $\thp$. This fact will be used in Lemma~\ref{bowling_pin2_lem} below.

\begin{thm}[Conquering a fat star-shaped set]
\label{cone_segment_conquering_thm}


\newcommand{\thisz}{\vect z}
\newcommand{\thisvertex}{\vect v}
\newcommand{\thiscube}{\cubify{\thisvertex}}

\newcommand{\thisAone}{\thisvertex}
\newcommand{\thisAtwo}{\initialset{2}}
\newcommand{\thisconfig}{\pair[1]{\thisAone}{\thisAtwo}}
\newcommand{\thisshell}{\shellof{\reacheddset{\mu_2}{\thisz}{r}}}

\newcommand{\thisevent}{\eventref{cone_segment_conquering_thm}^{\thisvertex} \parens{r,\thp;\beta}}

\newcommand{\thisC}{\bigcref{cone_segment_conquering_thm}(\beta)}
\newcommand{\thisc}{\smallcref{cone_segment_conquering_thm}(\beta)}

\newcommand{\thisdil}{\dilnorm{\linfty{d}}{\mu_2}}
\newcommand{\thisdilsymb}{\smallkonstsymb}

\newcommand{\thisfatstar}{\fatcstar{\normpair}{\beta}{\thp}{1}{\thisz}{\thisAtwo}}

Suppose $\tmeasurepair$ satisfies \finitespeed{d}, \expmoment, and \notiesprob, and let $\mu_1$ and $\mu_2$ be the shape functions for the component measures $\tmeasure_1$ and $\tmeasure_1$.
Let $\thisdilsymb \definedas \thisdil$, and for any $\thisvertex\in\Zd$, $r> \thisdilsymb$, $\thp\in (0,1]$, and $\beta >1$, define the event
\begin{multline*}
\thisevent \definedas \\
\bigcap_{\thisz\in\thiscube}
\eventthat[1]{\forall \thisAtwo \subseteq \thisshell,\;\;
\crreachedset[2]{1}{\thisconfig}{\Rd}{\rmeet{1}{1} s}
\supseteq \rreacheddset[2]{\mu_1}{\thisz}{\thisfatstar}{\rmeet{1}{\beta} s}
\text{ for all } s\ge r}.
\end{multline*}
Then given any $\beta>1$, there exist positive constants $\thisC$ and $\thisc$ such that for any $\thisvertex\in \Zd$, $r>\thisdilsymb$, and $\thp\in (0,1]$,
\[
\Pr \parens[2]{\thisevent}
\ge 1- \frac{1}{\thp^d} \thisC e^{-\thisc \thp r}.
\]
A symmetric statement holds with the roles of species 1 and 2 reversed.
\end{thm}

\begin{proof}


\newcommand{\thisz}{\vect z}
\newcommand{\thisy}{\vect y}
\newcommand{\genvertex}{\vect v}
\newcommand{\thisvertex}{\vect 0}
\newcommand{\thiscube}{\cubify{\thisvertex}}

\newcommand{\Aone}{\thisvertex}
\newcommand{\Atwo}{\initialset{2}}
\newcommand{\thisconfig}{\pair[1]{\Aone}{\Atwo}}
\newcommand{\thisshell}{\shellof{\reacheddset{\mu_2}{\thisz}{r}}}

\newcommand{\thisevent}{\eventref{cone_segment_conquering_thm}^{\thisvertex} \parens{r,\thp;\beta}}

\newcommand{\thisC}{\bigcref{cone_segment_conquering_thm}(\beta)}
\newcommand{\thisc}{\smallcref{cone_segment_conquering_thm}(\beta)}

\newcommand{\thisdil}{\dilnorm{\linfty{d}}{\mu_2}}
\newcommand{\thisdilsymb}{\smallkonstsymb}
\newcommand{\thisnormkonst}{k_{\mu_1}} 

\newcommand{\thisfatstar}{\fatcstar{\normpair}{\beta}{\thp}{1}{\thisz}{\Atwo}}
\newcommand{\shortstar}{\bigstar_{\Atwo}^{\thisz}}
\newcommand{\thiscloserset}{\closerset{\normpair}{\beta}{1}{\vect z}{\Atwo}}

\newcommand{\smallrad}{\rmeet{1}{\beta}r}
\newcommand{\bigrad}{\rmeet{1}{1}r}
\newcommand{\aone}{\alpha_1}
\newcommand{\atwo}{\alpha_2}

\newcommand{\smallball}{\reacheddset{\mu_1}{\thisz}{\smallrad}}

\newcommand{\nearevent}{E \parens{r,\thp;\beta}}
\newcommand{\farevent}{F \parens{r,\thp;\beta}}

\newcommand{\strongballevent}{\widetilde{\eventsymb}_{\ref{ball_conquering_lem}}^{\thisvertex} \parens{r;\beta}}
\newcommand{\conesevent}{\eventref{cone_segment_growth_thm}^{\thisvertex} \parens[2]{\bigrad,\thp;\spp(\beta)}}
\newcommand{\containmentevent}{ \eventref{growth_containment_lem} \parens[2]{\smallrad,\beta;\atwo(\beta)}}

\newcommand{\ballC}{\bigcref{ball_conquering_lem} (\beta)}
\newcommand{\ballc}{\smallcref{ball_conquering_lem} (\beta)}
\newcommand{\conesC}{\bigcref{cone_segment_growth_thm} \parens[2]{\spp(\beta)}}
\newcommand{\conesc}{\smallcref{cone_segment_growth_thm} \parens[1]{\spp(\beta)}}
\newcommand{\contC}{\bigcref{growth_containment_lem} \parens[2]{\atwo(\beta)}}
\newcommand{\contc}{\smallcref{growth_containment_lem} \parens[1]{\atwo(\beta)}}

By translation invariance it suffices to assume $\genvertex = \thisvertex$. As in Lemma~\ref{ball_conquering_lem}, the requirement $r>\thisdilsymb$ guarantees that $\thisconfig$ is a valid starting configuration for any $\Atwo\subseteq\thisshell$ if $\thisz\in\thiscube$.
The structure of the proof will be similar to that of Theorem~\ref{cone_segment_growth_thm}, in that we will define two events to deal separately with the ``near points" and ``far away points" in $\thisfatstar$.
Fix $\beta>1$ and $\thp\in (0,1]$, and for notational convenience, let $\shortstar \definedas \thisfatstar$. Define two events, corresponding to the ``near points" and ``far points", respectively, by
\begin{align*}
\nearevent &\definedas \bigcap_{\thisz,\Atwo}
\eventthat[3]{\ctime{1}{\shortstar}{\Aone}{\thisy} < \ptimenew{2}{\Atwo}{\thisy}
\text{ for all $\thisy\in \shortstar$ with $\distnew{\mu_1}{\thisz}{\thisy}\le\smallrad$}},\\
\farevent &\definedas \bigcap_{\thisz,\Atwo}
\eventthat[3]{\ctime{1}{\shortstar}{\Aone}{\thisy} < \ptimenew{2}{\Atwo}{\thisy}
\text{ for all $\thisy\in \shortstar$ with $\distnew{\mu_1}{\thisz}{\thisy}>\smallrad$}},
\end{align*}
where the intersection over $(\thisz,\Atwo)$ means the intersection over $\thisz\in\thiscube$ and $\Atwo\subseteq \thisshell$. Clearly, Lemma~\ref{sufficient_conquering_lem} implies that on the intersection $\nearevent\cap \farevent$, species 1 conquers the set $\shortstar$ for any $\thisz$ and $\Atwo$, so we seek a lower bound on the probability of this intersection. The following claim deals with the near points.

\begin{claim}
\label{cone_segment_conquering:nearevent_claim}
The event $\nearevent$ occurs on the event $\strongballevent$ defined in \eqref{ball_conquering:strong_event_eqn}.
\end{claim}

\begin{proof}[Proof of Claim~\ref{cone_segment_conquering:nearevent_claim}]
By Lemma~\ref{mu_cone_decomp_lem}, every $\mu_1$-ball with center $\thisz$ is a union of $(\mu_1,\thp)$-cone segments at $\thisz$ (for any $\thp\ge 0$), so Lemma~\ref{determ_bullseye_lem} and the monotonicity of the $\beta$-Voronoi cells imply that for any $\thisz\in\Rd$ and any $\Atwo\subseteq \thisshell$, we have
\[
\smallball \subseteq \fatcstar[1]{\normpair}{\beta}{\thp}{1}{\thisz}{\,\thisshell}
\subseteq \thisfatstar.
\]
This implies that for any $\thisz,\thisy\in\Rd$ and any $\Atwo\subseteq \thisshell$, we have $\ctime{1}{\thisfatstar}{\Aone}{\thisy} \le \ctime{1}{\smallball}{\Aone}{\thisy}$, and that
\[
\setbuilderbar[1]{\thisy\in \thisfatstar}{\distnew{\mu_1}{\thisz}{\thisy}\le\smallrad}
= \smallball.
\]
Moreover, we also have $\ptimenew{2}{\Atwo}{\thisy} \ge \ptimenew[2]{2}{\thisshell}{\thisy}$ for any such $\thisz$, $\thisy$, and $\Atwo$, and therefore,
\begin{align*}
\nearevent &= \bigcap_{\thisz,\Atwo}
\eventthat[4]{\ctime{1}{\thisfatstar}{\Aone}{\thisy} < \ptimenew{2}{\Atwo}{\thisy}
\text{ for all $\thisy\in \smallball$}}\\
&\supseteq \bigcap_{\thisz\in\thiscube}
\eventthat[1]{\ctime{1}{\smallball}{\Aone}{\thisy} \le \bigrad 
< \ptimenew{2}{\thisshell}{\thisy} \text{ for all $\thisy\in \smallball$}}\\
&=\strongballevent. \qedhere
\end{align*}
\end{proof}

Now it remains to deal with the far points. First define three constants in terms of $\beta$,
\[
\spp(\beta) \definedas 1-\frac{\rmeet{1}{\beta}}{\rmeet{1}{1}},
\qquad
\aone(\beta) \definedas \frac{\rmeet{1}{\beta}}{\rmeet{1}{1}},
\qquad
\atwo(\beta) \definedas \frac{\rmeet{1}{1}}{\beta \rmeet{1}{\beta}},
\]
and note that $\spp,\aone,\atwo\in (0,1)$, with $\spp\searrow 0$ and $\aone,\atwo\nearrow 1$ as $\beta\searrow 1$. 
These constants were chosen to satisfy the following relations, which will be needed in the proof of the subsequent two claims:
\[
1-\spp = \aone, \qquad
\aone \bigrad = \smallrad, \qquad
\aone^{-1} = \atwo\beta.
\]

\begin{claim}
\label{cone_segment_conquering:farevent_claim}
The event $\farevent$ occurs on the intersection
\[
\conesevent \cap \containmentevent,
\]
with $\tmeasure = \tmeasure_1$ and $\mu=\mu_1$ in the event $\conesevent$.
\end{claim}

\begin{proof}[Proof of \ref{cone_segment_conquering:farevent_claim}]
\newcommand{\genset}{S}
\newcommand{\thisunions}{\fatstars{\mu_1}{\thp}{\thisz}}

\newcommand{\speciesoneevent}{G_1(r,\thp;\beta)}
\newcommand{\speciestwoevent}{G_2(r,\thp;\beta)}

Since $\shortstar\in \thisunions$ for any $\thisz$ and $\Atwo$, we have
\begin{align}
\label{cone_segment_conquering:sp1_growth_eqn}
\conesevent
&= \bigcap_{\thisz\in\thiscube}
	\eventthat[3]{\forall \genset\in \thisunions,\;
	\crreachedset{\tmeasure_1}{\thisvertex}{\genset}{t}
	\supseteq \rreacheddset[2]{\mu_1}{\thisz}{\genset}{(1-\spp) t}
	\text{ for all } t\ge \bigrad} \notag\\
&\subseteq \bigcap_{\thisz,\Atwo} \eventthat[3]{
	\crreachedset{\tmeasure_1}{\thisvertex}{\shortstar}{t}
	\supseteq \rreacheddset[2]{\mu_1}{\thisz}{\shortstar}{(1-\spp) t}
	\text{ for all } t\ge \bigrad} \notag\\
&= \bigcap_{\thisz,\Atwo} \eventthat[4]{\ctime{1}{\shortstar}{\Aone}{\thisy}
	\le \bigrad\join \aone^{-1} \distnew{\mu_1}{\thisz}{\thisy}
	\text{ for all $\thisy\in \shortstar$}},
\end{align}
where the last line follows from inversion formula~\eqref{inversion_formulas:lbg} in Lemma~\ref{inversion_formula_lem}, since \mbox{$1-\spp = \aone$}. Note that on the final event in \eqref{cone_segment_conquering:sp1_growth_eqn}, for any fixed $\thisz$ and $\Atwo$, if $\thisy\in \shortstar$ and $\distnew{\mu_1}{\thisz}{\thisy}\ge \smallrad = \aone \bigrad$, then $\aone^{-1} \distnew{\mu_1}{\thisz}{\thisy} \ge \bigrad$, so it must hold that $\ctime{1}{\shortstar}{\Aone}{\thisy} \le \aone^{-1} \distnew{\mu_1}{\thisz}{\thisy}$. Therefore, \eqref{cone_segment_conquering:sp1_growth_eqn} implies that 
the event
\[
\speciesoneevent \definedas
\bigcap_{\thisz,\Atwo} \eventthat[3]{\ctime{1}{\shortstar}{\Aone}{\thisy}
	\le  \aone^{-1} \distnew{\mu_1}{\thisz}{\thisy}
	\text{ for all $\thisy\in \shortstar\setminus \interior{\smallball}$}} 
\]
occurs on the event $\conesevent$.
%
Thus, to show that $\farevent$ occurs on the intersection in Claim~\ref{cone_segment_conquering:farevent_claim}, it will suffice to prove that the event
\[
\speciestwoevent \definedas \bigcap_{\thisz,\Atwo}
\eventthat[3]{\ptimenew{2}{\Atwo}{\thisy} > \aone^{-1} \distnew{\mu_1}{\thisz}{\thisy}
\text{ for all $\thisy\in \shortstar \setminus \interior{\smallball}$}} 
\]
occurs on the event $\containmentevent$.

Note that since $\shortstar\subseteq \thiscloserset$ by definition, for any $\thisz$ and $\Atwo$ we have
\begin{equation}
\label{cone_segment_conquering:closer_points_eqn}
\forall \thisy\in\shortstar,\quad
\beta \distnew{\mu_1}{\thisz}{\thisy}
\le \distnew{\mu_2}{\Atwo}{\thisy}.
\end{equation}
Now recall that $\aone^{-1} = \atwo\beta$. Therefore, using \eqref{cone_segment_conquering:closer_points_eqn} in lines two and four, and using Lemma~\ref{dual_times_lem} (with $\tmeasure = \tmeasure_2$ and $\mu=\mu_2$) in line three,
\begin{align*}
\speciestwoevent
&= \bigcap_{\thisz,\Atwo}
\bigcap_{\substack{\thisy\in\shortstar\\ \distnew{\mu_1}{\thisz}{\thisy}\ge \smallrad}}
\eventthat[3]{ \ptimenew{2}{\Atwo}{\thisy} >  \atwo \beta \distnew{\mu_1}{\thisz}{\thisy}}\\
&\supseteq \bigcap_{\thisz,\Atwo}
\bigcap_{\substack{\thisy\in\shortstar\\ \distnew{\mu_1}{\thisz}{\thisy}\ge \smallrad}}
\eventthat[3]{ \ptimenew{2}{\Atwo}{\thisy} >  \atwo \distnew{\mu_2}{\Atwo}{\thisy}}\\
&\supseteq \bigcap_{\thisz,\Atwo}
\bigcap_{\substack{\thisy\in\shortstar\\ \distnew{\mu_1}{\thisz}{\thisy}\ge \smallrad}}
\eventthat[3]{\rreachedset{\tmeasure_2}{\thisy}{\Rd}{\atwo t} \subseteq \reacheddset{\mu_2}{\thisy}{t-} \; \forall t\ge \distnew{\mu_2}{\Atwo}{\thisy}}\\
&\supseteq \bigcap_{\thisz\in\thiscube}
\bigcap_{\substack{\thisy\in\Rd\\ \distnew{\mu_1}{\thisz}{\thisy}\ge \smallrad}}
\eventthat[3]{\rreachedset{\tmeasure_2}{\thisy}{\Rd}{\atwo t} \subseteq \reacheddset{\mu_2}{\thisy}{t-} \; \forall t\ge \beta \distnew{\mu_1}{\thisz}{\thisy}}\\
&=\containmentevent.
\end{align*}
Thus we have
\[
\conesevent\cap \containmentevent
\subseteq \speciesoneevent\cap \speciestwoevent
\subseteq \farevent. \qedhere
\]
\end{proof}

\newcommand{\intersectionevent}{H \parens{r,\thp;\beta}}

\begin{claim}
\label{cone_segment_conquering:finalevent_claim}
The event $\thisevent$ occurs on the intersection 
\[
\intersectionevent \definedas
\strongballevent\cap \conesevent \cap \containmentevent,
\]
with $\tmeasure = \tmeasure_1$ and $\mu=\mu_1$ in the event $\conesevent$.
\end{claim}

\begin{proof}[Proof of Claim~\ref{cone_segment_conquering:finalevent_claim}]
Combining Claims~\ref{cone_segment_conquering:nearevent_claim} and \ref{cone_segment_conquering:farevent_claim}, then applying Lemma~\ref{sufficient_conquering_lem}, we have
\begin{align}
\label{cone_segment_conquering:sp1_comparison_eqn}
\intersectionevent 
\subseteq \nearevent\cap \farevent
&\subseteq \bigcap_{\thisz,\Atwo}
	\eventthat[3]{\ctime{1}{\shortstar}{\Aone}{\thisy} < \ptimenew{2}{\Atwo}{\thisy}
	\text{ for all $\thisy\in \shortstar$}}, \notag\\
&\subseteq \bigcap_{\thisz,\Atwo}
	\eventthat[3]{\crreachedset{1}{\thisconfig}{\Rd}{t} \supseteq
	\crreachedset{\tmeasure_1}{\Aone}{\shortstar}{t}
	\text{ for all $t\ge 0$}}.
\end{align}
Additionally, since $\intersectionevent \subseteq \conesevent$, by \eqref{cone_segment_conquering:sp1_growth_eqn} we have
\begin{equation}
\label{cone_segment_conquering:sp1_growth2_eqn}
\intersectionevent \subseteq \bigcap_{\thisz,\Atwo} \eventthat[3]{
	\crreachedset{\tmeasure_1}{\thisvertex}{\shortstar}{t}
	\supseteq \rreacheddset[2]{\mu_1}{\thisz}{\shortstar}{(1-\spp) t}
	\text{ for all } t\ge \bigrad}.
\end{equation}
Claim~\ref{cone_segment_conquering:finalevent_claim} now follows directly from \eqref{cone_segment_conquering:sp1_comparison_eqn} and \eqref{cone_segment_conquering:sp1_growth2_eqn} by setting $t=\rmeet{1}{1} s$ and noting that \mbox{$(1-\spp) \rmeet{1}{1} = \rmeet{1}{\beta}$}.
\end{proof}

Finally, combining Claim~\ref{cone_segment_conquering:finalevent_claim} with Lemma~\ref{ball_conquering_lem}, Theorem~\ref{cone_segment_growth_thm}, and Lemma~\ref{growth_containment_lem}, we get
\begin{align*}
\Pr \parens[2]{\thisevent} \ge \Pr \parens[2]{\intersectionevent}
\ge 1	&- \ballC e^{-\ballc r}\\
	& -\frac{1}{\thp^d} \conesC e^{-\conesc \thp \bigrad}\\
	& -\contC e^{-\contc \cdot \beta \cdot [\smallrad - 3\dilnorm{\linfty{d}}{\mu_1} -1]}\\
\ge 1	&- \frac{1}{\thp^d} \thisC e^{-\thisc \thp r},
\end{align*}
where $\thisC$ and $\thisc$ are some positive constants depending only on $\beta$ and the distributions of $\tmeasure_1(\edge)$ and $\tmeasure_2(\edge)$.
\end{proof}

The following corollary restates the conclusion of Theorem~\ref{cone_segment_conquering_thm} with the starting configuration fixed ahead of time and without explicitly bounding the growth of species 1.

\begin{cor}[Conquering species 1's fat star set for fixed $\vect z$ and $\initialset{2}$]
\label{star_point_conquering_cor}

\newcommand{\thisv}{\vect v}
\newcommand{\thisz}{\vect z}
\newcommand{\thiscube}{\cubify{\thisv}}

\newcommand{\Aone}{\thisv}
\newcommand{\Atwo}{\initialset{2}}
\newcommand{\thisconfig}{\pair{\Aone}{\Atwo}}

\newcommand{\thisfatcstar}{\fatcstar{\normpair}{\beta}{\thp}{1}{\thisz}{\Atwo}}

\newcommand{\thisC}{\bigcref{cone_segment_conquering_thm} (\beta)}
\newcommand{\thisc}{\smallcref{cone_segment_conquering_thm} (\beta)}

Let $\thisv\in\Zd$ and $\Atwo\subseteq\Rd$, and suppose $\thisv\notin\lat{\Atwo}$. Then for any $\beta>1$, $\thp\in (0,1]$, and $\thisz\in\thiscube$,
\[
\Pr \eventthat[1]{\crreachedset{1}{\thisconfig}{\Rd}{\infty} \supseteq \thisfatcstar}
\ge 1- \frac{1}{\thp^d} \thisC e^{-\thisc \thp \cdot \distnew{\mu_2}{\thisz}{\Atwo}}.
\]
A symmetric statement holds with the roles of species 1 and 2 switched.
\end{cor}

\begin{proof}
\newcommand{\thisv}{\vect v}
\newcommand{\thisz}{\vect z}
\newcommand{\thiscube}{\cubify{\thisv}}

\newcommand{\Aone}{\thisv}
\newcommand{\Atwo}{\initialset{2}}
\newcommand{\thisconfig}{\pair{\Aone}{\Atwo}}

\newcommand{\thisfatcstar}{\fatcstar{\normpair}{\beta}{\thp}{1}{\thisz}{\Atwo}}
\newcommand{\thisshell}{\shellof{\reacheddset{\mu_2}{\thisz}{r}}}

\newcommand{\thisC}{\bigcref{cone_segment_conquering_thm} (\beta)}
\newcommand{\thisc}{\smallcref{cone_segment_conquering_thm} (\beta)}

\newcommand{\oldevent}{\eventref{cone_segment_conquering_thm}^{\thisv} \parens{r,\thp;\beta}}

\newcommand{\thisdil}{\dilnorm{\linfty{d}}{\mu_2}}

Set $r= \distnew{\mu_2}{\thisz}{\Atwo}$. Then $\Atwo \subseteq \thisshell$, so if $r>\thisdil$, then the result follows directly from Theorem~\ref{cone_segment_conquering_thm} because on the event $\oldevent$, species 1 conquers the set $\thisfatcstar$ from the initial configuration $\thisconfig$.
If $r\le\thisdil$, then the proof of Theorem~\ref{cone_segment_conquering_thm} goes through for any fixed $\thisz\in\thiscube$ and $\Atwo\subseteq \thisshell$ as long as the configuration $\thisconfig$ is valid, which we assumed explicitly. Alternatively, we could increase $\thisC$ if necessary to make the probability bound hold trivially for small $r$.
\end{proof}

The next result generalizes Corollary~\ref{star_point_conquering_cor} to the case where $A_1$ is any bounded set, and provides a stochastic analogue of Proposition~\ref{determ_star_conquering_prop} in which the set $\cstarset{\normpair}{1}{1}{\initialset{1}}{\initialset{2}}$ is replaced with $\fatcstarset{\normpair}{\beta}{\thp}{1}{\initialset{1}}{\initialset{2}}$ for some $\beta>1$ and $\thp>0$. We leave the proof to the reader.

\begin{cor}[Conquering a union of thick star-shaped sets when $\initialset{1}$ is bounded]
\label{star_conquering_cor}

\newcommand{\thisC}{\bigcref{star_conquering_cor}(\beta,\thp)}
\newcommand{\thisc}{\smallcref{star_conquering_cor}(\beta,\thp)}


Let $\initconfig$ be an initial configuration in $\Rd$ with $\initialset{1}$ bounded and $\distnew{\mu_2}{\initialset{1}}{\initialset{2}}\ge r$. Then for any $\beta>1$ and $\thp>0$, there exist positive constants $\thisC$ and $\thisc$ such that in the $\tmeasurepair$-competition process started from $\initconfig$, \speciesname{species 1} conquers all of $\fatcstarset{\normpair}{\beta}{\thp}{1}{\initialset{1}}{\initialset{2}}$ with probability at least $1-\card[2]{\lat{\initialset{1}}} \cdot \thisC e^{-\thisc r}$.
\end{cor}

\subsection{Bowling Pin Lemmas}

In this subsection we combine Corollary~\ref{star_point_conquering_cor} with the results in Chapter~\ref{deterministic_chap} analyzing the specific geometry of $\mu$-balls and $\mu$-cones, obtaining lower bounds on the probability that species $i$ will conquer a particular $\mu_i$-bowling pin. The following two ``bowling pin lemmas" are the key tools we will use in the next section to analyze the the two-type process when species 2 starts on the shell of a cone and species 1 starts at some point inside the cone.

The first bowling pin lemma, which will be used in the proof of Lemma~\ref{conq_mu1_subcone_lem} below, is a stochastic analogue of Lemma~\ref{cone_closer_set_lem} in which cones are replaced by bowling pins, and the parameter $\beta_0$ is required to be greater than 1.
Recall the definitions
\[
\thpadv{2}{\dirvec}{\beta} \definedas \beta \frac{\mu_1(\dirsymb)}{\mu_2(\dirsymb)}
\quad\text{and}\quad
\rmeet{1}{\beta} \definedas \frac{1}{\beta + \esupnorm{\mu_2/\mu_1}}
\]
from \eqref{thpadv_def_eqn} and \eqref{rmeet_def_eqn}.

\begin{lem}[Conquering a $\mu_1$-bowling pin inside the shell of a $\mu_2$-bowling pin]
\label{bowling_pin_lem}
\newcommand{\bigbeta}{\beta_0}
\newcommand{\smallbeta}{\beta}
\newcommand{\bigthp}{\thpadv{2}{\dirvec}{\bigbeta}}
\newcommand{\smallthp}{\thp}
\newcommand{\smallthpexp}{(\bigbeta-\smallbeta)\rmeet{1}{\smallbeta}}

\newcommand{\vertex}{\vect v}
\newcommand{\origin}{\vect z}

\newcommand{\bigrad}{r_0}
\newcommand{\smallrad}{r}

\newcommand{\bigpin}{\bpin[2]{\mu_2}{\bigrad}{\bigthp}{\origin}{\dirvec}{\mu_2(\dirsymb)}}
\newcommand{\smallpin}{\bpin[2]{\mu_1}{\smallrad}{\smallthp}{\origin}{\dirvec}{\mu_1(\dirsymb)}}

\newcommand{\biginfpin}{\bpin{\mu_2}{\bigrad}{\bigthp}{\origin}{\dirvec}{\infty}}
\newcommand{\smallinfpin}{\bpin{\mu_1}{\smallrad}{\smallthp}{\origin}{\dirvec}{\infty}}

\newcommand{\Atwo}{\initialset{2}}
\newcommand{\thisconfig}{\pair{\vertex}{\Atwo}}

\newcommand{\thisC}{\bigcref{cone_segment_conquering_thm}(\smallbeta)}
\newcommand{\thisc}{\smallcref{cone_segment_conquering_thm}(\smallbeta)}

Let $\tmeasurepair$ be a two-type traversal measure with limiting traversal norm pair $\normpair = \pair{\mu_1}{\mu_2}$ as described above. Fix $\vertex\in\Zd$, $\origin\in \cubify{\vertex}$, $\dirsymb\in \Rd\setminus \setof{\vect 0}$, $\bigrad >0$, and $\bigbeta > 1$. Suppose $1<\smallbeta<\bigbeta$, and let
\[
\smallrad \definedas \bigrad \rmeet{1}{\smallbeta}
\quad\text{and}\quad
\smallthp \definedas \smallthpexp \meet 1.
\]
Then in the $\tmeasurepair$-competition process,
\begin{enumerate}
\item \label{bowling_pin:finite_part}
If $\Atwo\subseteq \shellof{\bigpin}$ and $\vertex\notin \lat{\Atwo}$, then \speciesname{species 1} conquers $\smallpin$ from the starting configuration $\thisconfig$ with probability at least $1- \smallthp^{-d}\thisC e^{-\thisc \smallthp  \bigrad}.$

\item \label{bowling_pin:infinite_part}
If $\Atwo\subseteq \shellof{\biginfpin}$ and $\vertex\notin \lat{\Atwo}$, then \speciesname{species 1} conquers $\smallinfpin$ from the starting configuration $\thisconfig$ with probability at least $1- \smallthp^{-d}\thisC e^{-\thisc \smallthp  \bigrad}$.
\end{enumerate}
A symmetric statement holds with \speciesname{species 1} and \speciesname{species 2} switched.
\end{lem}

\begin{proof}

\newcommand{\bigbeta}{\beta_0}
\newcommand{\smallbeta}{\beta}
\newcommand{\bigthp}{\thpadv{2}{\dirvec}{\bigbeta}}
\newcommand{\smallthp}{\thp}
\newcommand{\smallthpexp}{(\bigbeta-\smallbeta)\rmeet{1}{\smallbeta}}

\newcommand{\vertex}{\vect v}
\newcommand{\origin}{\vect z}

\newcommand{\bigrad}{r_0}
\newcommand{\smallrad}{r}

\newcommand{\twoht}{h_2}
\newcommand{\oneht}{h_1}

\newcommand{\bigpin}{\bpin{\mu_2}{\bigrad}{\bigthp}{\origin}{\dirvec}{\twoht}}
\newcommand{\intbigpin}{\interior{\bigpin}}
\newcommand{\smallpin}{\bpin{\mu_1}{\smallrad}{\smallthp}{\origin}{\dirvec}{\oneht}}

\newcommand{\bigball}{\reacheddset{\mu_2}{\origin}{\bigrad}}
\newcommand{\intbigball}{\interior{\bigball}}
\newcommand{\smallball}{\reacheddset{\mu_1}{\origin}{\smallrad}}

\newcommand{\bigcone}{\conetip{\mu_2}{\bigthp}{\origin}{\dirvec}{\twoht}}
\newcommand{\intbigcone}{\interior{\bigcone}}
\newcommand{\smallcone}{\conetip{\mu_1}{\smallthp}{\origin}{\dirvec}{\oneht}}

\newcommand{\Atwo}{\initialset{2}}
\newcommand{\thisconfig}{\pair{\vertex}{\Atwo}}

\newcommand{\thisC}{\bigcref{cone_segment_conquering_thm}(\smallbeta)}
\newcommand{\thisc}{\smallcref{cone_segment_conquering_thm}(\smallbeta)}

\newcommand{\thisfatstars}{\fatstars{\mu_1}{\smallthp}{\origin}}
\newcommand{\thisfatcstarset}{\fatcstarset{\normpair}{\smallbeta}{\smallthp}{1}{\origin}{\Atwo}}
\newcommand{\thiscloserset}{\closerset{\normpair}{\smallbeta}{1}{\origin}{\Atwo}}

For Part~\ref{bowling_pin:finite_part}, set $\oneht = \mu_1(\dirsymb)$ and $\twoht = \mu_2(\dirsymb)$, and for Part~\ref{bowling_pin:infinite_part}, set $\oneht = \twoht = \infty$. 
We will show that $\smallpin\subseteq \thisfatcstarset$ and then apply Corollary~\ref{star_point_conquering_cor}. Note that $\smallpin\in \thisfatstars$ (cf.\ \eqref{bpins_are_fatstars_eqn}), so it remains to check that $\smallpin \subseteq \thiscloserset$.
Recall that
\[
\bigpin = \bigball\cup \bigcone,
\]
so the hypothesis $\Atwo \subseteq \shellof{\bigpin}$ implies that
\[
\Atwo \subseteq \shellof{\bigball}
\quad\text{and}\quad
\Atwo \subseteq \shellof{\bigcone}.
\]
Therefore, since $\smallrad = \bigrad \rmeet{1}{\smallbeta}$, Lemma~\ref{determ_bullseye_lem} implies that $\smallball \subseteq \thiscloserset$, and since $\smallthp \le \smallthpexp$, Lemma~\ref{cone_closer_set_lem} implies that $\smallcone \subseteq \thiscloserset$. Thus, since
\[
\smallpin = \smallball\cup\smallcone,
\]
we have $\smallpin\subseteq \thiscloserset$ and hence $\smallpin\subseteq \thisfatcstarset$.
Finally, by hypothesis we have $\smallthp\in (0,1]$ and $\distnew{\mu_2}{\origin}{\Atwo} \ge \bigrad$, so the result now follows directly from Corollary~\ref{star_point_conquering_cor}.
%
%
\end{proof}

\begin{thmremark}
As a function of $\beta$, the probability bound in Lemma~\ref{bowling_pin_lem} becomes worse both as $\beta\searrow 1$ (meaning that species 1 has a smaller margin of error with which to outrun species 2) and as $\beta\nearrow \beta_0$ (equivalently $\thp\searrow 0$, meaning that species 1 has a narrower cone segment in which to move freely).
\end{thmremark}


The second bowling pin lemma will be combined with Lemma~\ref{cross_section_covering_lem} to prove Theorem~\ref{narrow_survival_time_thm} below.

\begin{lem}[Conquering a $\mu_2$-bowling pin in species 2's $\beta$-Voronoi star-cell]
\label{bowling_pin2_lem}

\newcommand{\rad}{r}
\newcommand{\origin}{\vect y}

\newcommand{\thispin}{\bpin{\mu_2}{\rad}{\thp}{\origin}{\dirvec}{h}}
\newcommand{\thiscloser}{\closerset{\normpair}{\beta}{2}{\initialset{1}}{\origin}}
\newcommand{\thisconfig}{\pair{\initialset{1}}{\origin}}

\newcommand{\thisC}{\bigcref{bowling_pin2_lem}(\beta,\thp)}
\newcommand{\thisc}{\smallcref{bowling_pin2_lem}(\beta,\thp)}

Given $\beta>1$ and $\thp>0$, there exist positive constants $\thisC$ and $\thisc$ such that if $\thisconfig$ is an initial configuration in $\Rd$ such that $\thispin \subseteq \thiscloser$ for some $\rad>0$, $\dirvec\in \dirset{\Rd}$, and $h\in [0,\infty]$, then
\[
\Pr \eventthat[1]{\crreachedset{2}{\thisconfig}{\Rd}{\infty} \supseteq \thispin}
\ge 1- \thisC e^{-\thisc \rad}.
\]
%
%
\end{lem}

\begin{proof}


\newcommand{\rad}{r}
\newcommand{\origin}{\vect y}

\newcommand{\Aone}{\initialset{1}}
\newcommand{\Atwo}{\origin}
\newcommand{\thisconfig}{\pair{\Aone}{\Atwo}}

\newcommand{\thispin}{\bpin{\mu_2}{\rad}{\thp}{\origin}{\dirvec}{h}}
\newcommand{\thiscloser}{\closerset{\normpair}{\beta}{2}{\initialset{1}}{\origin}}
\newcommand{\thisball}{\reacheddset{\mu_2}{\origin}{\rad}}

\newcommand{\oldC}{\bigcref{cone_segment_conquering_thm}(\beta)}
\newcommand{\oldc}{\smallcref{cone_segment_conquering_thm}(\beta)}
\newcommand{\thisC}{\bigcref{bowling_pin2_lem}(\beta,\thp)}
\newcommand{\thisc}{\smallcref{bowling_pin2_lem}(\beta,\thp)}

\newcommand{\thisfatstars}{\fatstars{\mu_2}{\thp}{\origin}}
\newcommand{\thisfatcstar}{\fatcstar{\normpair}{\beta}{\thp}{2}{\Aone}{\Atwo}}
\newcommand{\thisfatcstarone}{\fatcstar{\normpair}{\beta}{\thp\meet 1}{2}{\Aone}{\Atwo}}

We will apply Corollary~\ref{star_point_conquering_cor} with the roles of species 1 and 2 switched.
Since $\thispin \in \thisfatstars$, the assumption $\thispin\subseteq \thiscloser$ implies that
\[
\thispin\subseteq \thisfatcstar\subseteq \thisfatcstarone.
\]
Therefore, noting that $\lat{\Aone} \cap \lat{\Atwo} = \emptyset$ by hypothesis, Corollary~\ref{star_point_conquering_cor} implies that
\[
\Pr \eventthat[1]{\crreachedset{2}{\thisconfig}{\Rd}{\infty} \supseteq \thispin}
\ge 1- \frac{1}{(\thp\meet 1)^d}\cdot \oldC
	e^{-\oldc (\thp\meet 1) \cdot \distnew{\mu_1}{\Aone}{\origin}}.
\]
Now note that by definition we have $\thisball \subseteq \thispin$ and $\thiscloser\cap \Aone = \emptyset$, so $\thisball \cap \Aone = \emptyset$. Therefore, $\distnew{\mu_2}{\Aone}{\origin} \ge r$, and hence
\[
\distnew{\mu_1}{\Aone}{\origin}
\ge \esupnorm{\mu_2/\mu_1}^{-1}\cdot \distnew{\mu_2}{\Aone}{\origin}
\ge \esupnorm{\mu_2/\mu_1}^{-1}\cdot r.
\]
Thus the conclusion of the lemma holds with
\[
\thisC \definedas \frac{1}{(\thp\meet 1)^d}\cdot \oldC
\quad\text{and}\quad
\thisc \definedas \oldc (\thp\meet 1) \cdot \esupnorm{\mu_2/\mu_1}^{-1}.
\qedhere
\]
\end{proof}

\section{Competition When One Species Starts Outside a Cone}
\label{random_fpc:cone_competition_sec}

The results in this section run parallel to those in Section~\ref{deterministic:cones:conquered_regions_sec}, with the deterministic $\normpair$-process replaced by the random $\tmeasurepair$-process. We consider a $\tmeasurepair$-competition process in which species 1 starts at a single point $\vect z$ inside a cone $\cones$, and species 2 starts on the exterior of $\cones$. Our main result is the following theorem, which is a stochastic analogue of Proposition~\ref{conquering_from_point_prop} for non-critical cones. Namely, in the random two-type process started from $\pair[2]{\lat{\vect z}}{\Zd\setminus \lat \cones}$ as just described, the probability of survival for species 1 is positive if $\cones$ is wide, and zero if $\cones$ is narrow and contained in a half-space. As we did for Proposition~\ref{conquering_from_point_prop}, we will break the proof of Theorem~\ref{wide_narrow_survival_thm} into smaller pieces that will be treated individually below.


\begin{thm}[The probability of survival in wide and narrow cones]
\label{wide_narrow_survival_thm}
Suppose $\tmeasurepair$ satisfies the assumptions of Theorem~\ref{irrelevance_thm} and also \expmoment. Let $\cones$ be a cone, let $\vect z\in \cones$, and consider a $\tmeasurepair$-process with starting configuration $\pair[2]{\vect z}{\Zd\setminus \lat \cones}$.
\begin{enumerate}
\item \label{wide_narrow_survival:wide_part}
If $\cones$ is wide for species 1, then the probability that species 1 survives and conquers a nondegenerate subcone of $\cones$ is positive. Moreover, we can assume that the thickness of this conquered subcone is bounded below by a constant that depends only on the advantage of species 1 in $\cones$.

\item \label{wide_narrow_survival:narr_part}
If $\cones$ is narrow for species 1 and additionally is contained in a half-space, then the probability that species 1 survives is zero.
\end{enumerate}
\end{thm}

\begin{proof}
\newcommand{\conept}{\vect z}

Part~\ref{wide_narrow_survival:wide_part} follows from Corollary~\ref{wide_cone_survival_possible_cor} below, and Part~\ref{wide_narrow_survival:narr_part} follows from Corollary~\ref{narrow_extinction_cor} below.
\end{proof}

\begin{thmremark}
Note that the deterministic analogue of Part~\ref{wide_narrow_survival:wide_part} of Theorem~\ref{wide_narrow_survival_thm} is false in general for nonconvex cones: A nonconvex wide cone may have regions which look locally like very narrow cones, and if $\vect z$ is in such a region, then species 1 cannot survive from $\vect z$ in the deterministic process.
Correspondingly, the survival probability from such a $\vect z$ in the random process will be very low, but nevertheless positive as long as $\tmeasurepair$ satisfies the hypotheses of Theorem~\ref{irrelevance_thm}.
\end{thmremark}

We will deduce Part~\ref{wide_narrow_survival:wide_part} of Theorem~\ref{wide_narrow_survival_thm} as a consequence of the following large deviations estimate, which follows easily from the Bowling Pin Lemma~\ref{bowling_pin_lem} above.

\begin{thm}[Survival in wide cones with high probability]
\label{wide_cone_survival_thm}
\newcommand{\thisadv}{\alpha}
\newcommand{\conept}{\vect z}
\newcommand{\bddist}{\distnew{\mu_2}{\bd \cones}{\conept}}

\newcommand{\thisconfig}{\pair[2]{\conept}{\, \Zd\setminus \lat{\cones}}}
\newcommand{\thisconfigadv}{\configadvantage{1}{\normpair}{\conept}{\cones}}

\newcommand{\thisthp}{\symbolref{\thp}{wide_cone_survival_thm}(\thisadv)}
\newcommand{\thisC}{\bigcref{wide_cone_survival_thm}(\thisadv)}
\newcommand{\thisc}{\smallcref{wide_cone_survival_thm}(\thisadv)}

%

Suppose $\tmeasurepair$ satisfies \finitespeed{d}, \expmoment, and \notiesprob. Let $\cones\subseteq\Rd$ be a cone that is wide for species 1, and let $\vect z\in \cones$ with $\thisconfigadv = \thisadv >1$. Then there exist positive constants $\thisthp$, $\thisC$, and $\thisc$ such that if $\bddist \ge r$, then the probability that species 1 conquers some nondegenerate $\mu_1$-cone of thickness $\thisthp$ at $\conept$ in the $\tmeasurepair$-process started from $\thisconfig$ is at least $1-\thisC e^{-\thisc r}$. If $\cones$ is convex and has advantage $\thisadv>1$, then this bound holds for all $\conept\in\cones$ with $\bddist\ge r$.
%
\end{thm}

\begin{proof}
\newcommand{\conept}{\vect z}
\newcommand{\thisadv}{\alpha}
\newcommand{\smalleradv}{\thisadv'}
\newcommand{\smallestadv}{\thisadv''}

\newcommand{\bddist}{\distnew{\mu_2}{\bd \cones}{\conept}}
\newcommand{\thisdist}{r}
\newcommand{\smalldist}{\thisdist_{\thisadv}}

\newcommand{\thisconfigadv}{\configadvantage{1}{\normpair}{\conept}{\cones}}

\newcommand{\thisthp}{\symbolref{\thp}{wide_cone_survival_thm}(\thisadv)}
\newcommand{\thisC}{\bigcref{wide_cone_survival_thm}(\thisadv)}
\newcommand{\thisc}{\smallcref{wide_cone_survival_thm}(\thisadv)}

\newcommand{\bigthp}{\thpadv{2}{\dirvec}{\smalleradv}}
\newcommand{\twocone}{\cone{\mu_2}{\bigthp}{\conept}{\dirvec}}
\newcommand{\onecone}{\cone{\mu_1}{\thisthp}{\conept}{\dirvec}}

\newcommand{\bigball}{\reacheddset{\mu_2}{\conept}{\thisdist-}}

\newcommand{\bigpin}{\bpin{\mu_2}{\thisdist}{\bigthp}{\conept}{\dirvec}{\infty}}
\newcommand{\bigpinint}{\interior{\bigpin}}
\newcommand{\bigpinshell}{\shellof{\bigpin}}
\newcommand{\smallpin}{\bpin{\mu_1}{\smalldist}{\thisthp}{\conept}{\dirvec}{\infty}}

\newcommand{\oldC}{\bigcref{cone_segment_conquering_thm}(\smallestadv)}
\newcommand{\oldc}{\smallcref{cone_segment_conquering_thm}(\smallestadv)}

Let $\smalleradv \definedas \frac{1}{3} +\frac{2}{3}\thisadv$ and $\smallestadv \definedas \frac{2}{3} + \frac{1}{3}\thisadv$, so $1<\smallestadv<\smalleradv<\thisadv$.
Since $\smalleradv <\thisadv = \thisconfigadv$, there is by definition some $\dirvec\in \dirset{\Rd}$ such that $\twocone \subseteq \cones$. Since $\bddist \ge \thisdist$, we also have $\bigball\subseteq \cones$. Therefore, $\bigpinint\subseteq \cones$, and if $\initialset{2} = \Zd\setminus \lat \cones$, then Lemma~\ref{lattice_cube_properties_lem} implies that
\[
\initialset{2} \subseteq \cubify{\initialset{2}} \subseteq \Rd\setminus \cones
\subseteq \Rd\setminus \bigpinint =\bigpinshell.
\]
Now let $\thisthp \definedas (\smalleradv - \smallestadv) \rmeet{1}{\smallestadv} \meet 1 >0$ and $\smalldist \definedas \thisdist \rmeet{1}{\smallestadv}$. Then since $\smallestadv>1$, Part~\ref{bowling_pin:infinite_part} of Lemma~\ref{bowling_pin_lem} implies that species 1 conquers $\smallpin \supset \onecone$ from the configuration $\pair{\conept}{\initialset{2}}$ with probability at least $1- \thisC e^{-\thisc \thisdist}$, where
\[
\thisC \definedas \thisthp^{-d} \oldC,
\quad\text{and}\quad
\thisc \definedas \thisthp \oldc.
\]
If $\cones$ is convex, then $\thisconfigadv = \thisadv$ for all $\conept\in\cones$ by Lemma~\ref{convex_advantage_lem}, so in this case the bound holds for all $\conept\in\cones$ with $\bddist\ge \thisdist$.
\end{proof}

\begin{cor}[Survival in wide cones with positive probability]
\label{wide_cone_survival_possible_cor}
\newcommand{\thisadv}{\alpha}
\newcommand{\conept}{\vect z}
\newcommand{\thisthp}{\symbolref{\thp}{wide_cone_survival_thm}(\thisadv)}
\newcommand{\thisconfig}{\pair[2]{\conept}{\,\Zd\setminus \lat{\cones}\,}}

Suppose $\tmeasurepair$ satisfies the assumptions of Theorem~\ref{wide_narrow_survival_thm}. Let $\cones\subseteq \Rd$ be a cone with $\coneadvantage{1}{\normpair}{\cones} = \thisadv >1$, and let $\conept\in\cones$. Then in the $\tmeasurepair$-process started from $\thisconfig$, the probability that species 1 conquers some $\mu_1$-cone of thickness $\thisthp$ is positive.
\end{cor}

\begin{proof}
Combine Theorems~\ref{wide_cone_survival_thm} and \ref{irrelevance_thm}.
\end{proof}


Part~\ref{wide_narrow_survival:narr_part} of Theorem~\ref{wide_narrow_survival_thm} will be a consequence of
Theorem~\ref{narrow_survival_time_thm} below, which gives a
large deviations estimate for the size of the set conquered by species~1 in a small narrow cone. We refer to this estimate as a ``tail bound for the survival time of species~1" 
because with high probability, the duration of species 1's survival is comparable to the $\mu_1$-radius of the set it conquers.
In order to state the result, we introduce notation for the conquered region in the random process, analogous to the definition \eqref{conq_set_def_eqn} for the deterministic process.  For $\cones \subseteq\Rd$ and $\initialset{1}\subseteq \cones$, we write $\conqset[1]{1}{\tmeasurepair}{\initialset{1}}{\cones}$ for species 1's conquered region in the $\tmeasurepair$-process started from $\pair[2]{\lat{\initialset{1}}}{\,\Zd\setminus \lat \cones}$, i.e.
\begin{equation}
\label{random_conq_set_def_eqn}
\conqset[1]{1}{\tmeasurepair}{\initialset{1}}{\cones}
\definedas \conqset[1]{1}{\tmeasurepair}{\lat{\initialset{1}}}{\lat \cones}
\definedas
\crreachedset{1}{\pair{\lat{\initialset{1}}}{\,\Zd\setminus \lat \cones}}{\Rd}
	{\infty}_{\tmeasurepair}.
\end{equation}
Note that Lemma~\ref{monotonicity_lem} implies that for any $\cones,\cones'\subseteq\Rd$ and any realization of $\tmeasurepair$,
\begin{equation}
\label{rand_conq_sets_monotone_eqn}
\lat{\cones}\subseteq \lat{\cones'} \implies
\conqset[1]{1}{\tmeasurepair}{\initialset{1}}{\cones}
\subseteq \conqset[1]{1}{\tmeasurepair}{\initialset{1}}{\cones'}.
\end{equation}
The proof of Theorem~\ref{narrow_survival_time_thm} is essentially a probabilistic, ``lattice-ized" version of the proof of Proposition~\ref{narrow_extinction_prop} for the deterministic process, relying on the ``cross-section covering" Lemma~\ref{cross_section_covering_lem} and the probability estimate in Bowling Pin Lemma~\ref{bowling_pin2_lem} above.

\begin{thm}[Tail bound for survival time in small narrow cones]
\label{narrow_survival_time_thm}

\newcommand{\smallrad}{\alpha}

Suppose $\tmeasurepair$ satisfies \finitespeed{d}, \expmoment, and \notiesprob. Let $\cones\subset\Rd$ be a closed cone with apex $\vect a\in\Rd$ such that $\cones$ is contained in some half-space and $\coneadvantage{1}{\normpair}{\cones} <1$. Then there exist positive constants $R_0$, $\smallrad$, $\bigcref{narrow_survival_time_thm}$, and $\smallcref{narrow_survival_time_thm}$ such that for any $\vect z\in \cones$ and  $R\ge R_0\join \smallrad^{-1}\distnew{\mu_1}{\vect a}{\vect z}$,
\[
\Pr
\eventthat[3]{ \conqset[1]{1}{\tmeasurepair}{\vect z}{\cones}
	\subseteq \reacheddset{\mu_1}{\vect a}{R}}
\ge 1- \bigcref{narrow_survival_time_thm}
e^{-\smallcref{narrow_survival_time_thm} R}.
\]
\end{thm}

\begin{proof}
\newcommand{\thisadv}{\coneadvantage{1}{\normpair}{\cones}}
\newcommand{\thiscone}{\cones}
\newcommand{\thispt}{\vect z}

\newcommand{\thisapex}{\vect a}
\newcommand{\crosssection}{K}
\newcommand{\pincollection}{\boldsymbol{P}}
\newcommand{\interiorcollection}{\interior{\pincollection}}

\newcommand{\pinsymb}{P} 

\newcommand{\pin}{\pinsymb}
\newcommand{\pininterior}{\interior{P}}
\newcommand{\pinorigin}{\vect y_{\pin}}

\newcommand{\nbhd}{U_{\thisapex}}


\newcommand{\conecomp}{\Rd\setminus \thiscone}
\newcommand{\startset}{\Rd \setminus \cubify[2]{\lat{\thiscone}}}
\newcommand{\latstartset}{\Zd\setminus \lat{\thiscone}}
\newcommand{\cubelatcone}{\cubify[2]{\lat{\thiscone}}}

\newcommand{\bluntcone}{\cones \setminus \setof{\thisapex}}
\newcommand{\thisconfig}{\pair[2]{\thispt}{\startset}}
\newcommand{\thislatconfig}[1][0]{\pair[#1]{\thispt}{\latstartset}}

\newcommand{\thisconqset}{\conqset[1]{1}{\tmeasurepair}{\thispt}{\thiscone}}

\newcommand{\genpt}{\vect x}

\newcommand{\dilconst}{\kappa}

\newcommand{\coverrad}{\varepsilon}
\newcommand{\smallrad}{\alpha}

\newcommand{\pinfamrad}{R_0}
\newcommand{\pinfactor}{R_{\pin}}

\newcommand{\zrad}{R_{\min}(\thispt)}
\newcommand{\bigrad}{R}

\newcommand{\thishom}{\homothe{\bigrad}{\thisapex}}

\newcommand{\thisthick}{\thp_\beta}
\newcommand{\pinheadrad}{r_{\pin}}
\newcommand{\minheadrad}{r_0}

\newcommand{\thisunitball}{\reacheddset{\mu_1}{\thisapex}{1}}

\newcommand{\smallball}{\reacheddset{\mu_1}{\thisapex}{\smallrad}} 
\newcommand{\bigball}{\reacheddset{\mu_1}{\thisapex}{\bigrad \smallrad}}
\newcommand{\biggerball}{\reacheddset{\mu_1}{\thisapex}{\bigrad}}

\newcommand{\thisevent}{\eventref{narrow_survival_time_thm}(\bigrad)}
\newcommand{\blockingevent}{\eventref{narrow_survival_time:blocking_event_claim}(\bigrad)}


Referring to \eqref{random_conq_set_def_eqn} and Lemma~\ref{lattice_cube_properties_lem}, note that $\thisconqset$ corresponds to the starting configuration
\begin{equation}
\label{narrow_survival_time:initconfig_eqn}
\initconfig = \pair[1]{\lat\thispt}{\,\latstartset} = \pair[1]{\lat\thispt}{\,\lat[2]{\startset}}.
\end{equation}
Also note that since we are working in $\Zd$ rather than $\Rd$, it suffices to assume $\anapex\in \cubify{\vect 0}$, but it does not suffice to assume $\anapex=\vect 0$. Thus, we will carry out the proof for a general apex $\anapex\in\Rd$.
Set $\beta \definedas \frac{1}{2}\parens[2]{1+\thisadv^{-1}}$, 
so $1<\beta< \thisadv^{-1}$.
Let $\crosssection = \bd \thisunitball \cap \thiscone$, and let $\pincollection = \pincollection(K,\beta)$ be the finite collection of $\mu_2$-bowling pins from Lemma~\ref{cross_section_covering_lem}. Then $\pincollection$ covers $\crosssection$, and there is some open neighborhood $\nbhd$ of $\thisapex$ such that
\begin{equation}
\label{narrow_survival_time:pin_closer_eqn}
\forall \pin\in \pincollection, \quad
\pin \subseteq \cstarset{\normpair}{\beta}{2}{\nbhd}{\pinorigin},
\end{equation}
where $\pinorigin \in \conecomp$ is the origin of $\pin$.
Recall that for $R\ge 0$, $\homothe{R}{\thisapex} \colon \Rd\to \Rd$ denotes the homothety with scale factor $R$ and center $\thisapex$, i.e.
\[
\homothe{R}{\thisapex}\vect x \definedas \thisapex + R (\vect x - \thisapex)
\quad\text{for $\vect x\in \Rd$}.
\]
Since $\pincollection$ covers $\crosssection$ and is finite, there is some $\coverrad>0$ such that $\crosssection+\coverrad \unitball{\linfty{d}} \subseteq \bigcup \pincollection$. Assume without loss of generality that $\coverrad\le 1$.

\begin{claim}
\label{narrow_survival_time:K_containment_claim}
There is some $\dilconst<\infty$ such that if $\bigrad\ge \dilconst \coverrad^{-1}$, then $\bd \biggerball \cap \cubelatcone \subseteq \thishom \parens[2]{\crosssection + \coverrad \unitball{\linfty{d}}}$.
\end{claim}

\begin{proof}[Proof of Claim~\ref{narrow_survival_time:K_containment_claim}]
\newcommand{\spherept}{\vect x}
\newcommand{\conept}{\vect y_{\spherept}}
\newcommand{\projpt}{\vect z_{\spherept}}

\newcommand{\dilone}{\dilnorm{\linfty{d}}{\mu_1}}
\newcommand{\dilinfty}{\dilnorm{\mu_1}{\linfty{d}}}

Let $\dilconst \definedas 1+ \parens[2]{1\join \dilinfty}\cdot \dilone$, and let $\bigrad \ge \dilconst \coverrad^{-1}$.
By Lemma~\ref{lattice_cube_properties_lem}, we have 
$
\bd \biggerball \cap \cubelatcone
\subseteq \bd \biggerball \cap \parens[2]{\thiscone+\unitball{\linfty{d}}}.
$
Thus, given any $\spherept\in \bd \biggerball \cap \cubelatcone$, there exists some $\conept\in \thiscone$ with $\distnew{\linfty{d}}{\spherept}{\conept}\le 1$.
Since $\coverrad\le 1$, we have $\bigrad\ge \dilconst > \dilone$, and hence
\[
\distnew{\mu_1}{\,\thisapex}{\conept}
\ge \distnew{\mu_1}{\thisapex}{\spherept} - \distnew{\mu_1}{\spherept}{\conept}
\ge \bigrad - 1\cdot \dilone>0.
\]
Thus, for any $\spherept\in \bd \biggerball \cap \cubelatcone$, we have $\homothe{\mu_1(\conept-\thisapex)^{-1}}{\thisapex} (\conept) \in \bd \thisunitball \cap \thiscone = \crosssection$, and
we can define
\[
\projpt \definedas \thishom \circ \homothe{\mu_1(\conept-\thisapex)^{-1}}{\thisapex} (\conept)
= \thisapex +\bigrad \cdot \frac{\conept-\thisapex}{\mu_1(\conept-\thisapex)}
\in 
\thishom\crosssection.
\]
That is, $\projpt$ is the radial projection (from $\thisapex$) of $\conept\in \thiscone$ onto $\thishom \crosssection$. Then $\projpt$ is a $\mu_1$-closest point in $\bd \biggerball$ to $\conept$, so since $\spherept \in \bd \biggerball$ and $\distnew{\linfty{d}}{\conept}{\spherept}\le 1$, we have
\[
\distnew{\mu_1}{\conept}{\projpt} \le \distnew{\mu_1}{\conept}{\spherept}
\le 1 \cdot \dilone.
\]
Therefore, for any $\spherept\in \bd \biggerball \cap \cubelatcone$, we have
\begin{align*}
\distnew{\linfty{d}}{\spherept}{\thishom \crosssection}
\le \distnew{\linfty{d}}{\spherept}{\projpt}
&\le \distnew{\linfty{d}}{\spherept}{\conept} + \distnew{\linfty{d}}{\conept}{\projpt}\\
&\le 1 + \distnew{\mu_1}{\conept}{\projpt}\cdot \dilinfty\\
&\le 1+ \dilone\cdot \dilinfty\\
&\le \dilconst.
\end{align*}
Therefore, since $\dilconst \le \bigrad\coverrad$, we have
\begin{align*}
\bd \biggerball \cap \cubelatcone
&\subseteq \setbuilder[1]{\spherept \in \Rd}{\distnew{\linfty{d}}
	{\spherept}{\thishom \crosssection}\le \dilconst}\\
&\subseteq \setbuilder[1]{\spherept \in \Rd}{\distnew{\linfty{d}}
	{\spherept}{\thishom \crosssection}\le \bigrad\coverrad}
= \thishom \crosssection+\bigrad \coverrad \unitball{\linfty{d}}.
\qedhere
\end{align*}
\end{proof}
%
%
Since $\thiscone$ is closed and $\pinorigin\in \conecomp$, for each $\pin\in \pincollection$ there exists some $\pinfactor<\infty$ such that
$
\thishom \pinorigin \in \startset
$
for all $\bigrad\ge \pinfactor$.
Let $\dilconst$ be the constant from Claim~\ref{narrow_survival_time:K_containment_claim}, and set
\[
\pinfamrad \definedas  \frac{\dilconst}{\coverrad} \join \max_{\pin\in\pincollection} \pinfactor <\infty.
\]
Let $\smallrad>0$ be small enough that $\smallball \subset \nbhd$. Then it follows as in the proof of Proposition~\ref{narrow_extinction_prop} that $\smallrad<1$.
Now fix $\thispt \in \thiscone$ and set
\[
\zrad \definedas \pinfamrad \join \frac{1}{\smallrad}\distnew{\mu_1}{\thisapex}{\thispt}.
\]
Then if $\bigrad\ge \zrad$, we have $\distnew{\mu_1}{\thisapex}{\thispt} \le \bigrad \smallrad$, and hence
\begin{equation}
\label{narrow_survival_time:z_containment_eqn}
\thispt \in \bigball =\thishom \smallball 
\quad \forall \bigrad\ge \zrad.
\end{equation}
Moreover, the following hold for any $\bigrad\ge \zrad$.
\begin{enumerate}
\item \label{narrow_survival_time:scaled_pins_cover}
The scaled family of $\mu_2$-bowling pins $\thishom \pincollection \definedas \setof{\thishom \pin}_{\pin \in \pincollection}$ covers the thickened, scaled cross section $\thishom \parens[2]{\crosssection + \coverrad \unitball{\linfty{d}}}$, since $\pincollection$ covers $\crosssection+\coverrad \unitball{\linfty{d}}$.

\item \label{narrow_survival_time:K_containment}
$\bd \biggerball \cap \cubelatcone \subseteq \thishom \parens[2]{\crosssection + \coverrad \unitball{\linfty{d}}}$ by Claim~\ref{narrow_survival_time:K_containment_claim}, since $\bigrad\ge \dilconst \coverrad^{-1}$.

\item \label{narrow_survival_time:z_in_biggerball}
$\thispt \in \biggerball$ by \eqref{narrow_survival_time:z_containment_eqn}, since $\smallrad<1$.

\item \label{narrow_survival_time:origins_outside}
For all $\pin\in\pincollection$, we have $\thishom \pinorigin \in \startset$, since $\bigrad\ge  \max_{\pin\in\pincollection} \pinfactor$.

\item \label{narrow_survival_time:z_in_nbhd}
$\thispt \in \thishom \nbhd$, by \eqref{narrow_survival_time:z_containment_eqn} and the choice of $\smallrad$.

\item \label{narrow_survival_time:scaled_pins_closer}
For all $\pin\in \pincollection$, we have $\thishom\pin \subseteq \cstarset{\normpair}{\beta}{2}{\thispt}{\thishom \pinorigin}$, by \eqref{narrow_survival_time:pin_closer_eqn} and Statement~\ref{narrow_survival_time:z_in_nbhd} above.
\end{enumerate}
Now let $\thisevent$ be the event in the statement of the theorem, i.e.
\[
\thisevent \definedas
\eventthat[3]{\conqset[1]{1}{\tmeasurepair}{\thispt}{\thiscone} \subseteq \biggerball}.
\]
Then

\begin{claim}
\label{narrow_survival_time:blocking_event_claim}
For any $\bigrad\ge \zrad$,  $\thisevent$ occurs almost surely on the event
\begin{equation*}
\blockingevent \definedas
\eventthat[1]{
\textstack{Species 2 conquers each $\mu_2$-bowling pin in $\thishom \pincollection$}
	{from the initial configuration $\thisconfig$}}.
\end{equation*}
\end{claim}

\begin{proof}[Proof of Claim \ref{narrow_survival_time:blocking_event_claim}]
Referring to \eqref{narrow_survival_time:initconfig_eqn}, the above event can be written
\[
\blockingevent = \bigcap_{\pin\in \pincollection}
\eventthat[1]{\crreachedset{2}{\thislatconfig}{\Rd}{\infty}_{\tmeasurepair} \supseteq \thishom \pin}.
\]
Since the family $\thishom \pincollection$ covers $\thishom \parens[2]{\crosssection + \coverrad \unitball{\linfty{d}}}$ by Statement~\ref{narrow_survival_time:scaled_pins_cover} above, on the event $\blockingevent$ it holds that
$\crreachedset{2}{\thislatconfig}{\Rd}{\infty}_{\tmeasurepair} \supseteq \thishom \parens[2]{\crosssection + \coverrad \unitball{\linfty{d}}}$.
By Statement~\ref{narrow_survival_time:K_containment} above, we have $\bd \biggerball \cap \cubelatcone \subseteq \thishom \parens[2]{\crosssection + \coverrad \unitball{\linfty{d}}}$, and since $\biggerball \setminus \cubelatcone \subseteq \startset = \intcubify{\initialset{2}}$, it follows that
\begin{equation}
\label{narrow_survival_time:sp2_sphere_eqn}
\crreachedset{2}{\thislatconfig}{\Rd}{\infty}_{\tmeasurepair} \supseteq \bd \biggerball
\text{ on the event $\blockingevent$.}
\end{equation}
Let
\[
\finalset{1} = \reachedset{1}{\thislatconfig}{\infty}_{\tmeasurepair}
\quad\text{and}\quad
\finalset{2} = \reachedset{2}{\thislatconfig}{\infty}_{\tmeasurepair}.
\]
By the assumptions on $\tmeasurepair$, we have $\finalset{1}\cap \finalset{2} =\emptyset$ almost surely. Thus, \eqref{narrow_survival_time:sp2_sphere_eqn} implies that $\finalset{1}\cap \lat[2]{\bd \biggerball} = \emptyset$ almost surely on the event $\blockingevent$.
Now, by Statement~\ref{narrow_survival_time:z_in_biggerball} above, we have $\thispt \in \biggerball$ for all $\bigrad\ge \zrad$, and hence $\lat \thispt \in \lat[2]{\biggerball}$. Lemma~\ref{lat_euc_bd_separation_lem} then implies that any lattice path from $\lat\thispt$ to $\lat[2]{\Rd\setminus\biggerball}$ must intersect $\lat[2]{\bd\biggerball}$. Therefore, if $\finalset{1}\cap \lat[2]{\bd \biggerball} = \emptyset$, then species~1 cannot conquer any points in $\lat[2]{\Rd\setminus \biggerball}$ because any such point would have to be connected to $\lat\thispt$ by some path in $\finalset{1}$, by Part~\ref{final_sets:components_part} of Lemma~\ref{final_sets_properties_lem}. This shows that almost surely on the event $\blockingevent$, $\lat[2]{\thisconqset}$ is contained in the set $\Zd \setminus \parens[2]{\lat[2]{\Rd\setminus \biggerball}}$. By Lemma~\ref{lattice_cube_properties_lem}, this implies that
\[
\thisconqset \subseteq \cubify[1]{\Zd \setminus \parens[3]{\lat[2]{\Rd\setminus \biggerball}}}
\subseteq \biggerball,
\]
and hence $\thisevent$ occurs almost surely on $\blockingevent$.
\end{proof}

To complete the proof, we need a lower bound on the probability of the event $\blockingevent$ when $\bigrad \ge \zrad$. By Statement~\ref{narrow_survival_time:origins_outside} above we have  $\thishom \pinorigin \in \startset$ and hence $\lat[2]{\thishom \pinorigin} \subseteq \latstartset$ for each $\pin\in \pincollection$, and therefore
\begin{equation}
\label{narrow_survival_time:pin_conq_eqn}
\blockingevent = \bigcap_{\pin\in \pincollection}
	\eventthat[1]{\crreachedset{2}{\thislatconfig}{\Rd}{\infty} \supseteq \thishom \pin}
\supseteq \bigcap_{\pin\in \pincollection}
	\eventthat[1]{\crreachedset{2}{(\thispt,\thishom \pinorigin)}{\Rd}{\infty} \supseteq \thishom \pin}.
\end{equation}
Let $\thisthick$ be the common thickness of the $\mu_2$-bowling pins in $\pincollection$, and for each $\pin\in \pincollection$, let $\pinheadrad>0$ be the head radius of $\pin$. Then the head radius of the scaled bowling pin $\thishom \pin$ is $\bigrad \pinheadrad$, so Statement~\ref{narrow_survival_time:scaled_pins_closer} above and Lemma~\ref{bowling_pin2_lem} imply that
\begin{equation}
\label{narrow_survival_time:pin_conq_prob_eqn}
\Pr \eventthat[1]{\crreachedset{2}{(\thispt,\thishom \pinorigin)}{\Rd}{\infty} \supseteq \thishom \pin}
\ge 1-\bigcref{bowling_pin2_lem}(\beta,\thisthick)
e^{-\smallcref{bowling_pin2_lem}(\beta,\thisthick) \bigrad \pinheadrad}.
\end{equation}
Thus if we set $\minheadrad \definedas \min_{\pin\in \pincollection} \pinheadrad>0$, then \eqref{narrow_survival_time:pin_conq_eqn} and \eqref{narrow_survival_time:pin_conq_prob_eqn} imply that
\begin{equation*}
\Pr \parens[1]{\eventcomplement{\blockingevent }}
\le \sum_{\pin\in \pincollection}
	\Pr \parens[10]{\eventthat[1]{\crreachedset{2}{(\thispt,\thishom \pinorigin)}{\Rd}{\infty}
	\not\supseteq \thishom \pin}}
\le \card{\pincollection} \bigcref{bowling_pin2_lem}(\beta,\thisthick)
e^{-\smallcref{bowling_pin2_lem}(\beta,\thisthick) \minheadrad \bigrad}.
\end{equation*} 
Therefore, Claim~\ref{narrow_survival_time:blocking_event_claim} implies that the statement in Theorem~\ref{narrow_survival_time_thm} holds with $\bigcref{narrow_survival_time_thm} \definedas \card{\pincollection} \bigcref{bowling_pin2_lem}(\beta,\thisthick)$ and $\smallcref{narrow_survival_time_thm} \definedas \smallcref{bowling_pin2_lem}(\beta,\thisthick) \minheadrad$.
\end{proof}

\begin{thmremark}
\newcommand{\thisapex}{\vect a}
\newcommand{\crosssection}{K}
\newcommand{\pincollection}{\boldsymbol{P}}
\newcommand{\interiorcollection}{\interior{\pincollection}}

\newcommand{\pin}{P}
\newcommand{\pininterior}{\interior{P}}
\newcommand{\pinorigin}{\vect y_{\pin}}

\newcommand{\nbhd}{U_{\thisapex}}
\newcommand{\smallrad}{\alpha}
\newcommand{\coverrad}{\varepsilon}

A priori, the four constants in Theorem~\ref{narrow_survival_time_thm} can depend on the specific geometry of the cone $\cones$. However, with a little more work, it should be possible to show that all four constants only depend on the single parameter $\coneadvantage{1}{\normpair}{\cones}$. In particular, this would follow if we could strengthen the statement of Lemma~\ref{cross_section_covering_lem} by proving the following additional statements about the collection of $\mu_2$-bowling pins $\pincollection$ and the neighborhood $\nbhd$ for a fixed spherical section $K$:
\begin{enumerate}

\item $\distnew{\mu_2}{\pinorigin}{\cones}$ is bounded below by some constant depending only on $\beta$ and $\coneadvantage{1}{\normpair}{\cones}$.

\item There is some positive constant $\smallrad$ depending only on $\beta$ and $\coneadvantage{1}{\normpair}{\cones}$ such that the neighborhood $\nbhd$ contains the $\mu_1$-ball $\reacheddset{\mu_1}{\thisapex}{\smallrad}$.

\item The head radius $r_{\pin}$ of each $\pin\in \pincollection$ is bounded below by some positive constant $r$ depending only on $\beta$ and $\coneadvantage{1}{\normpair}{\cones}$.

\item There exists some $\coverrad>0$ depending only on $\beta$ and $\coneadvantage{1}{\normpair}{\cones}$ such that for each $\vect x\in\crosssection$, the cube $\reacheddset{\linfty{d}}{\vect x}{\coverrad}$ is contained in some $\pin\in \pincollection$.

\item $\card{\pincollection}$ is bounded above by some constant that depends only on $\beta$ and $\coneadvantage{1}{\normpair}{\cones}$.

\end{enumerate}
I believe this improved version of Lemma~\ref{cross_section_covering_lem} should be true, but at present I have not worked out the geometric details needed for the proof.
\end{thmremark}

The following corollary shows that Part~\ref{wide_narrow_survival:narr_part} of Theorem~\ref{wide_narrow_survival_thm} follows from Theorem~\ref{narrow_survival_time_thm}, even without the additional assumption that $\tmeasurepair$ satisfies the hypotheses of Theorem~\ref{irrelevance_thm}.

\begin{cor}[Extinction in small narrow cones]
\label{narrow_extinction_cor}
Suppose $\tmeasurepair$ satisfies the hypotheses of Theorem~\ref{narrow_survival_time_thm}. Let $\cones$ be a cone that is narrow for species 1 and contained in some half-space, and let $\vect z\in \cones$. Then the probability that species 1 survives from the starting configuration $\pair[2]{\vect z}{\Zd\setminus \lat \cones}$ is zero.
\end{cor}

\begin{proof}
Note that
\[
\eventthat[3]{\text{\speciesname{species~1} survives from } \pair[2]{\vect z}{\Zd\setminus \lat \cones}}
= \bigcap_{N\in\N}
\eventthat[3]{ \conqset[1]{1}{\tmeasurepair}{\vect z}{\cones}
	\not\subseteq \reacheddset{\mu_1}{\vect a}{N}},
\]
so the result follows from Theorem~\ref{narrow_survival_time_thm} and continuity of measure.
\end{proof}

\section{Additional Results for Random First-Passage Competition in $\Zd$}
\label{random_fpc:additional_results_sec}


Here we briefly discuss the result of \Haggstrom and Pemantle \cite[Proposition~2.2]{Haggstrom:2000aa} which is the key step in proving Theorem~\ref{hp_slow_growth_thm} \cite[Lemma~5.2]{\HPb}. We state a slightly more general version of this result (valid for infinite starting configurations), which appears in \cite[Proposition~5.2]{Deijfen:2007aa}. Recall that $\statespace = \statespace[1]$ is the state space for the two-type process, and $\processlaw{\lambda_1,\lambda_2}{X}$ denotes the law of the two-type Richardson model with exponential rates $\lambda_1,\lambda_2>0$ and initial configuration $X\in\statespace$.

\begin{prop}[Low coexistence probability from ``point vs.\ ball" configurations {\cite[Proposition~2.2]{Haggstrom:2000aa}}]
\label{hp_point_ball_prop}

Consider the two-type Richardson model on $\Zd$ with rates $\lambda_1=1$ and $\lambda_2 = \lambda<1$, and let $\mu$ be the shape function for \speciesname{species~1}. For each $r> 0$ and $\beta>1$, define the set of initial configurations
\[
\pvbset{\beta}{r} \definedas
\setbuilder[2]{X\in \statespace}{\stateset{2}{X}\subseteq r\unitball{\mu}
	\text{ and } \stateset{1}{X} \not\subseteq \beta r\unitball{\mu}}.
\]
Then for any $\beta>1$,
\[
\lim_{r\to\infty} \sup_{X\in \pvbset{\beta}{r}}
\processlaw{1,\lambda}{X} \eventthat[2]{\text{\speciesname{species~2} survives}}
=0.
\]
\end{prop}

We describe Proposition~\ref{hp_point_ball_prop} as a result about low coexistence probability from ``point vs.\ ball" configurations because these are the minimal configurations in the set $\pvbset{\beta}{r}$ (where, of course, $\pvbsymb$ stands for ``point vs.\ ball") with respect to the ordering from Section~\ref{fpc_basics:monotonicity_sec}. That is, every configuration $X\in \pvbset{\beta}{r}$ dominates a configuration $X'$ in which 
\speciesname{species~2} occupies an entire $\mu$-ball of radius $r$ while \speciesname{species~1} occupies a single vertex outside a larger $\mu$-ball of radius $\beta r$,
so these ``point vs.\ ball" configurations are the worst possible elements of $\pvbset{\beta}{r}$ for \speciesname{species~1}. 
Proposition~\ref{hp_point_ball_prop} shows that if \speciesname{species~1} is the faster species, it still has a high probability of winning from such a minimal configuration when $r$ is large.

%
%


The following result is a special case of Theorem~\ref{wide_cone_survival_thm} in the case where $\cones$ is a $\mu_2$-cone; it is a stochastic analogue of Lemma~\ref{determ_conq_mu1_subcone_lem} that applies when species 1's advantage in the $\mu_2$-cone is \emph{strictly} greater than 1.

\begin{lem}[Conquering a $\mu_1$-subcone of a wide $\mu_2$-cone with high probability]
\label{conq_mu1_subcone_lem}

\newcommand{\thisadv}{\alpha}
\newcommand{\zdist}{r}

\newcommand{\bigthp}{\thpadv{2}{\dirvec}{\thisadv}}
\newcommand{\smallthp}{\symbolref{\thp}{conq_mu1_subcone_lem} (\thisadv)}

\newcommand{\thisapex}{\vect a}
\newcommand{\conept}{\vect z}

\newcommand{\bigcone}{\cone{\mu_2}{\bigthp^-}{\thisapex}{\dirvec}}
\newcommand{\smallcone}{\cone{\mu_1}{\smallthp}{\conept}{\dirvec}}

\newcommand{\thisC}{\bigcref{conq_mu1_subcone_lem}(\thisadv)}
\newcommand{\thisc}{\smallcref{conq_mu1_subcone_lem}(\thisadv)}

Let $\tmeasurepair$ be a two-type traversal measure on $\edges{\Zd}$ satisfying \finitespeed{d}, \expmoment, and \notiesprob, and let $\normpair = \pair{\mu_1}{\mu_2}$ be the pair of shape functions (norms) corresponding to $\tmeasurepair$.
Suppose $\esupnorm[2]{\frac{\mu_2}{\mu_1}}\ge \thisadv>1$, and for $\dirsymb\in \Rd\setminus \setof{\vect 0}$, let $\thpadv{2}{\dirvec}{\thisadv} = \thisadv \frac{\mu_1(\dirsymb)}{\mu_2(\dirsymb)}$ as defined in \eqref{thpadv_def_eqn}. Then

\begin{enumerate}

\item \label{conq_mu1_subcone:cone_part}
If $\frac{\mu_2(\dirsymb)}{\mu_1(\dirsymb)} \ge\thisadv$ and $\thisapex\in\Rd$, the set $\cones\definedas \bigcone$ is a small cone with $\coneadvantage{1}{\normpair}{\cones}\ge \thisadv$.

\item \label{conq_mu1_subcone:conq_part}
There exist positive constants $\smallthp$, $\thisC$, and $\thisc$ such that if $\dirvec$ and $\cones$ are as in Part~\ref{conq_mu1_subcone:cone_part}, then for any $\conept\in \cones$ with $\distnew{\mu_2}{\conept}{\bd \cones} \ge \zdist$,
\[
\Pr \eventthat[1]{
	\conqset[1]{1}{\tmeasurepair}{\conept}{\cones} \supseteq \smallcone}
\ge 1 - \thisC e^{-\thisc \zdist}.
\]
\end{enumerate}
\end{lem}

\begin{proof}
\newcommand{\thisapex}{\vect a}
\newcommand{\conept}{\vect z}

\newcommand{\thisadv}{\alpha}
\newcommand{\smalladv}{\thisadv'}

\newcommand{\bigthp}{\thpadv{2}{\dirvec}{\thisadv}}
\newcommand{\smallthp}{\symbolref{\thp}{conq_mu1_subcone_lem} (\thisadv)}

\newcommand{\zdist}{r}
\newcommand{\smalldist}{r'}

\newcommand{\bigconea}{\cone{\mu_2}{\bigthp^-}{\thisapex}{\dirvec}}
\newcommand{\bigconez}{\cone{\mu_2}{\bigthp^-}{\conept}{\dirvec}}

\newcommand{\smallcone}{\cone{\mu_1}{\smallthp}{\conept}{\dirvec}}

\newcommand{\bigball}{\reacheddset{\mu_2}{\conept}{\zdist-}}

\newcommand{\bigpin}{\bpin{\mu_2}{\zdist}{\bigthp}{\conept}{\dirvec}{\infty}}
\newcommand{\bigpinint}{\interior{\bigpin}}
\newcommand{\bigpinshell}{\shellof{\bigpin}}
\newcommand{\smallpin}{\bpin{\mu_1}{\smalldist}{\smallthp}{\conept}{\dirvec}{\infty}}

\newcommand{\thisC}{\bigcref{conq_mu1_subcone_lem}(\thisadv)}
\newcommand{\thisc}{\smallcref{conq_mu1_subcone_lem}(\thisadv)}

\newcommand{\oldC}{\bigcref{cone_segment_conquering_thm}(\smalladv)}
\newcommand{\oldc}{\smallcref{cone_segment_conquering_thm}(\smalladv)}

Part~\ref{conq_mu1_subcone:cone_part} is the same as in Lemma~\ref{determ_conq_mu1_subcone_lem}. For Part~\ref{conq_mu1_subcone:conq_part}, let $\smalladv \definedas \frac{1}{2}(1+\thisadv)>1$, let $\smallthp \definedas (\thisadv - \smalladv) \rmeet{1}{\smalladv} \meet 1 >0$, and let $\smalldist \definedas \zdist \rmeet{1}{\smalladv}$.
Since $\cones = \bigconea$ is convex and $\conept\in\cones$, Part~\ref{wedge_properties:convexity_part} of Lemma~\ref{wedge_properties_lem} implies that $\bigconez \subseteq \bigconea$, and since $\distnew{\mu_2}{\vect z}{\bd \cones} \ge \zdist$, we have $\bigball \subseteq \cones = \bigconea$. Therefore, $\bigpinint \subseteq \bigconea$, and if $\initialset{2} = \Zd\setminus \lat \cones$, then Lemma~\ref{lattice_cube_properties_lem} implies that
\[
\initialset{2} \subseteq \cubify{\initialset{2}} \subseteq \Rd\setminus \cones
= \Rd\setminus \bigconea \subseteq \Rd\setminus \bigpinint
=\bigpinshell.
\]
Thus, by Part~\ref{bowling_pin:infinite_part} of Lemma~\ref{bowling_pin_lem}, species 1 conquers $\smallpin \supset \smallcone$ from the configuration $\pair{\conept}{\initialset{2}}$ with probability at least $1- \thisC e^{-\thisc \zdist}$, where
\[
\thisC \definedas \smallthp^{-d} \oldC,
\quad\text{and}\quad
\thisc \definedas \smallthp \oldc.
\qedhere
\]
\end{proof}

Using Lemma~\ref{conq_mu1_subcone_lem}, we prove the following result, which will be needed in Chapter~\ref{coex_finite_chap}.

\begin{lem}[Uniform lower bound for survival probability in cones with nearby apex]
\label{unif_cone_survival_lem}
\newcommand{\thisadv}{\alpha}
\newcommand{\thisvertex}{\vect 0}

\newcommand{\bigthp}{\thpadv{2}{\dirvec}{\thisadv}}
\newcommand{\smallthp}{\thp_{\ref{conq_mu1_subcone_lem}}(\thisadv)}

\newcommand{\bigapex}{\vect a}
\newcommand{\smallapex}{\vect z}

\newcommand{\bigcone}{\cone{\mu_2}{\bigthp}{\bigapex}{\dirvec}}
\newcommand{\smallcone}{\cone{\mu_1}{\smallthp}{\smallapex}{\dirvec}}


Suppose $\tmeasurepair$ satisfies the hypotheses of Theorem~\ref{irrelevance_thm} and also \expmoment, and suppose $\frac{\mu_2(\dirsymb)}{\mu_1(\dirsymb)}\ge \thisadv>1$, where $\normpair = \pair{\mu_1}{\mu_2}$ is the pair of shape functions for $\tmeasurepair$, and $\dirsymb\in \Rd\setminus \setof{\vect 0}$. Let $\bigthp = \thisadv \frac{\mu_1(\dirsymb)}{\mu_2(\dirsymb)} \in (0,1]$, and let $\smallthp$ be the constant from Lemma~\ref{conq_mu1_subcone_lem}. Then there exists $\vect z\in \dirvec$ such that
\[
\inf_{\bigapex \in \cubify{\thisvertex}}
\Pr \eventthat[1]{\conqset[1]{1}{\tmeasurepair}{\thisvertex}{\bigcone} \supseteq \smallcone}
>0.
\]
\end{lem}

\begin{proof}
\newcommand{\thisadv}{\alpha}
\newcommand{\thisvertex}{\vect 0}

\newcommand{\bigthp}{\thpadv{2}{\dirvec}{\thisadv}}
\newcommand{\thisthp}{\thp} 
\newcommand{\smallthp}{\thp_{\ref{conq_mu1_subcone_lem}}(\thisadv)}

\newcommand{\bigapex}{\vect a}
\newcommand{\smallapex}{\vect z}

\newcommand{\bigcone}{\cone{\mu_2}{\thisthp}{\bigapex}{\dirvec}}
\newcommand{\smallcone}{\cone{\mu_1}{\smallthp}{\smallapex}{\dirvec}}

\newcommand{\thisdil}{\kappa}
\newcommand{\thiscube}{\cubify{\thisvertex}}
\newcommand{\thisball}{\reacheddset{\mu_2}{\thisvertex}{\thisdil}}

\newcommand{\midapex}{\vect z_0}
\newcommand{\midcone}{\cone{\mu_2}{\thisthp}{\midapex}{\dirvec}}

\newcommand{\thisdist}{R}
\newcommand{\axispt}{\vect z_{\thisdist}}

\newcommand{\apath}{\gamma_{\bigapex}}

\newcommand{\bigball}{\reacheddset[2]{\mu_2}{\thisvertex}{\mu_2(\smallapex)}}

\newcommand{\pathcone}{\cones_{\bigapex}}


Set $\thisthp \definedas \bigthp$, set $\thisdil \definedas \frac{1}{2} \dilnorm{\linfty{d}}{\mu_2}$, and define $\midapex \definedas \frac{\thisdil}{\thisthp} \munitvec$. Then $\thiscube \subseteq\thisball$, and by Part~\ref{special_mu:cone_intersection_part} of Lemma~\ref{special_mu_properties_lem}, we have
\begin{equation}
\label{unif_cone_survival:intersection_eqn}
\bigcap_{\bigapex\in \thiscube} \bigcone
\;\;\supseteq \bigcap_{\bigapex\in\thisball} \bigcone
\;\;=\;\; \midcone.
\end{equation}
For each $\thisdist\ge 0$, set $\axispt \definedas \midapex+\thisdist \munitvec \in\midcone$.
By Lemma~\ref{mu_cone_bd_lem} (Part~\ref{mu_cone_bd:axis_part}), $\distnew[2]{\mu_2}{\bd \midcone}{\,\axispt} \to\infty$ as $\thisdist\to\infty$, so Lemma~\ref{conq_mu1_subcone_lem} implies that for all sufficiently large $\thisdist$,
\begin{equation}
\label{unif_cone_survival:high_prob_eqn}
\Pr \eventthat[1]{\conqset[1]{1}{\tmeasurepair}{\axispt}{\,\midcone}
\supseteq \cone{\mu_1}{\smallthp}{\axispt}{\dirvec}}>0.
\end{equation}
Choose $\thisdist$ large enough that \eqref{unif_cone_survival:high_prob_eqn} holds, and set $\smallapex\definedas \axispt$. Now for each $\bigapex\in\thiscube$, let $\apath$ be a lattice path in $\lat[2]{\bigcone\cap \bigball}$ from $\thisvertex$ to $\lat{\smallapex}$, and set $\pathcone \definedas \apath \cup \midcone$. Note that the path $\apath$ exists for each $\bigapex\in\thiscube$ because the lattice set $\lat[2]{\bigcone\cap \bigball}$ is connected and contains both $\thisvertex$ and $\lat{\smallapex}$. Moreover, since each of the paths $\apath$ is contained in the ball $\lat[2]{\bigball}$, there are only finitely many such paths, and therefore the collection
$
\setbuilder{\pathcone}{\bigapex\in\thiscube}
$
is finite.
Furthermore, since $\lat[2]{\bigcone}\supseteq \apath$ by definition, \eqref{unif_cone_survival:intersection_eqn} implies that $\lat[2]{\bigcone}\supseteq \lat{\pathcone}$ for all $\bigapex\in\thiscube$.

Now, since $\pathcone\supseteq\midcone$ by construction, we have $\conqset[1]{1}{\tmeasurepair}{\smallapex}{\pathcone}\supseteq \conqset[2]{1}{\tmeasurepair}{\smallapex}{\,\midcone}$ by \eqref{rand_conq_sets_monotone_eqn}, so \eqref{unif_cone_survival:high_prob_eqn} implies that
\begin{equation*}
\Pr \eventthat[1]{\conqset[1]{1}{\tmeasurepair}{\smallapex}{\pathcone}
\supseteq \smallcone}>0
\quad\text{for all $\bigapex\in\thiscube$,}
\end{equation*}
and since $\thisvertex$ and $\lat{\smallapex}$ are in the same component of $\Zd\setminus \lat{\pathcone}$, Theorem~\ref{irrelevance_thm} then implies that
\begin{equation}
\label{unif_cone_survival:pos_prob_eqn}
\Pr \eventthat[1]{\conqset[1]{1}{\tmeasurepair}{\thisvertex}{\pathcone}
\supseteq \smallcone}>0
\quad\text{for all $\bigapex\in\thiscube$.}
\end{equation}
Finally, since $\lat[2]{\bigcone}\supseteq \lat{\pathcone}$, we have $\conqset[2]{1}{\tmeasurepair}{\thisvertex}{\bigcone}\supseteq \conqset{1}{\tmeasurepair}{\thisvertex}{\pathcone}$, so \eqref{unif_cone_survival:pos_prob_eqn} plus the fact that the collection $\setbuilder{\pathcone}{\bigapex\in\thiscube}$ is finite imply that
\begin{align*}
\inf_{\bigapex \in \thiscube}
\Pr \eventthat[1]{\conqset[1]{1}{\tmeasurepair}{\thisvertex}{\bigcone} \supseteq \smallcone}
&\ge \inf_{\bigapex \in \thiscube}
\Pr \eventthat[1]{\conqset[1]{1}{\tmeasurepair}{\thisvertex}{\pathcone}
\supseteq \smallcone}\\
&=\; \min_{\pathcone}\,
\Pr \eventthat[1]{\conqset[1]{1}{\tmeasurepair}{\thisvertex}{\pathcone}
\supseteq \smallcone}
>0.
\qedhere
\end{align*}
\end{proof}



The next two results follow from Theorem~\ref{wide_narrow_survival_thm} and are analogues of Propositions~\ref{mu_cone_criticality_prop} and \ref{mu_norms_critical_speed_prop} for the deterministic process.

%
%

\begin{thm}[Existence of critical thickness for $\mu$-cones]
\label{mu_cone_criticality_thm}
\newcommand{\speedparam}{\lambda}
\newcommand{\thisnormpair}{\normpair_{\speedparam}}

\newcommand{\thisttimepair}{\tmeasurepair_{\speedparam}}

\newcommand{\thisapex}{\vect a}
\newcommand{\thiscone}{\cone{\mu}{\thp}{\thisapex}{\dirvec}}
\newcommand{\thispt}{\vect z}
\newcommand{\thisconfig}{\pair[1]{\thispt}{\,\Zd\setminus \lat[2]{\thiscone}}}
\newcommand{\thisaxis}{\thisapex+\cdirvec}

\newcommand{\thisconqset}{\conqset[1]{1}{\thisttimepair}{\thispt}{\thiscone}}

\newcommand{\critspeed}{\speedparam_c(\cones)}
\newcommand{\critthp}{\thp_c}

Let $\tmeasure$ be a one-type \iid\ traversal measure on $\edges{\Zd}$ such that $\ttime(\edge)$ has a continuous distribution function with finite exponential moment, and suppose $0\in \supp \ttime(\edge)$.
Let $\mu$ be the shape function for $\tmeasure$, and for each $\speedparam>0$ define $\thisttimepair \definedas \pair[2]{\tmeasure}{\, \speedparam^{-1} \tmeasure}$ and $\thisnormpair \definedas \pair[2]{\mu}{\,\speedparam^{-1}\mu}$. 
Consider a $\thisttimepair$-competition process started from $\thisconfig$, where $\thisapex\in\Rd$, $\dirvec\in\dirset{\Rd}$, $\thp\in (0,1]$, and $\thispt\in\thiscone$.
\begin{enumerate}
\item $\displaystyle \coneadvantage[1]{1}{\thisnormpair}{\thiscone} = {\thp}/{\speedparam}$, and hence the cone $\thiscone$ is wide, critical, or narrow for species~1 according to whether $\thp>\speedparam$, $\thp=\speedparam$, or $\thp<\speedparam$, respectively.

\item If $\thp > \speedparam$, then $\displaystyle \Pr \eventthat[1]{\thisconqset \text{\rm contains a nondegenerate cone}}>0$.

\item If $\thp < \speedparam$, then $\displaystyle \Pr \eventthat[1]{\thisconqset \text{\rm is bounded}}=1$.
\end{enumerate}
%
%
%
%
%
%
\end{thm}

\begin{thm}[Existence of critical speed ratio in small cones]
\label{mu_norms_critical_speed_thm}
\newcommand{\spone}{\lambda_1}
\newcommand{\sptwo}{\lambda_2}
\newcommand{\spratio}{\sptwo/\spone}

\newcommand{\newnormpair}{\normpair_{\spone,\sptwo}}
\newcommand{\newttimepair}{\tmeasurepair_{\spone,\sptwo}}

\newcommand{\thispt}{\vect z}
\newcommand{\thisconfig}{\pair[2]{\thispt}{\,\Rd\setminus \lat{\cones}}}

\newcommand{\thisadv}{\coneadvantage{1}{\normpair}{\cones}}

\newcommand{\thisconqset}{\conqset{1}{\newnormpair}{\thispt}{\cones}}

\newcommand{\critratio}{\lambda_c(\cones)}
\newcommand{\bcritratio}{\Lambda_c(\cones)}


Let $\tmeasurepair = \pair{\tmeasure_1}{\tmeasure_2}$ be a two-type \iid\ traversal measure on $\edges{\Zd}$ such that $0\in \supp \ttime_1(\edge)$, and for $i\in \setof{1,2}$, $\ttime_i(\edge)$ is a continuous nonnegative random variable with finite exponential moment,
and let $\normpair = \pair{\mu_1}{\mu_2}$ be the pair of shape functions corresponding to $\tmeasurepair$.  For any $\spone,\sptwo>0$, consider the random two-type process run with the pair of scaled traversal times $\newttimepair = \pair[2]{\spone^{-1}\tmeasure_1}{\sptwo^{-1}\tmeasure_2}$ and started from the configuration $\thisconfig$, where $\cones\subset\Rd$ is any cone contained in a half-space, and $\vect z\in \cones$.
%
%
%
\begin{enumerate}
\item If $\sptwo/\spone<\thisadv$, then species~1 conquers a nondegenerate subcone of $\cones$ with positive probability.
\item If $\sptwo/\spone>\thisadv$, then species~1 almost surely conquers only a bounded region.
\end{enumerate}
Thus, any nondegenerate small cone $\cones$ has a critical speed ratio $\bcritratio \in \cosegment[1]{\esupnorm[2]{\frac{\mu_2}{\mu_1}}^{-1}}{\infty}$ such that in the $\newttimepair$-process started from $\thisconfig$, species~1 conquers a nondegenerate cone with positive probability if $\lambda_1/\lambda_2 > \bcritratio$ and dies out with probability one if $\lambda_1/\lambda_2 < \bcritratio$.
\end{thm}

\chapter{Coexistence from Finite Starting Configurations}
\label{coex_finite_chap}

In this chapter, our goal will be to strengthen \Haggstrom and Pemantle's Theorem~\ref{hp_slow_growth_thm} \cite[Lemma~5.2]{\HPb}, which showed that in the two-type Richardson model started from finite initial configurations, if coexistence occurs, then the total infected region grows asymptotically at the same rate as the slow species. Namely, we will prove the following theorem:

\begin{thm}[Coexistence implies fast species can't touch support hyperplanes infinitely often]
\label{conv_hull_thm}
\newcommand{\thisstate}{X}
\newcommand{\thisrate}{\lambda}
\newcommand{\rateone}{\lambda_1}
\newcommand{\ratetwo}{\lambda_2}
\newcommand{\thislaw}{\processlaw{\rateone,\ratetwo}{\thisstate}}
\newcommand{\chullevent}{H}
Consider a two-type process in which $\cmttime{1}$ and $\cmttime{2}$ are exponentially distributed with rates $\rateone$ and $\ratetwo$, respectively, and assume $\rateone>\ratetwo$.
Let $\thisstate \equiv \initconfig$ be a finite initial configuration, and let $\chullevent(\thisstate)$ be the event that in the process started from $\thisstate$, there exists a sequence of times $t_1,t_2,\dotsc$ with $t_n\to\infty$ such that at each time $t_n$, the total occupied region $\bothreachedset{\thisstate}{t_n}$ has a support point in $\reachedset{1}{\thisstate}{t_n}$, the set of vertices occupied by \speciesname{species~1}.
Then $\Pr \parens[2]{\coexstate{\tmeasurepair}{\thisstate} \cap \chullevent(\thisstate)} = 0$.
\end{thm}

The fact that Theorem~\ref{conv_hull_thm} implies Theorem~\ref{hp_slow_growth_thm} follows from the Shape Theorem and the observation that the growth of the slow species in the two-type process is dominated by its disentangled version.
More explicitly, suppose Theorem~\ref{conv_hull_thm} holds, and let $X \equiv\initconfig$ be a finite initial configuration. Then almost surely on the event $\coexstate{\tmeasurepair}{X}$, there exists $t_0<\infty$ such that for all $t\ge t_0$ we have $\conv \bothreachedset{X}{t} = \conv \reachedset{2}{X}{t}$ (cf.\ Lemma~\ref{convex_hull_properties_lem}). By Part~\ref{entangled_disentangled:unrestricted_part} of Lemma~\ref{entangled_disentangled_lem}, we also have $\reachedset{2}{X}{t} \subseteq \reachedset{\tmeasure_2}{\initialset{2}}{t}$ for all $t\ge 0$. Now since $\initialset{2}$ is finite, the Shape Theorem for species 2 implies that almost surely, for any $\epsilon>0$ there is some $t_\epsilon<\infty$ such that $\reachedset{\tmeasure_2}{\initialset{2}}{t} \subseteq (1+\epsilon) t \unitball{\mu_2}$ for all $t\ge t_\epsilon$, where $\mu_2$ is the shape function for $\tmeasure_2$. Thus, since $\unitball{\mu_2}$ is convex, almost surely on the event $\coexstate{\tmeasurepair}{X}$, for any $\epsilon>0$ we have $\bothreachedset{X}{t} \subseteq \conv \reachedset{2}{X}{t} \subseteq (1+\epsilon) t \unitball{\mu_2}$ for all $t\ge t_0\join t_\epsilon$, so the conclusion of Theorem~\ref{hp_slow_growth_thm} holds. On the other hand, it is clear that Theorem~\ref{hp_slow_growth_thm} does not immediately imply Theorem~\ref{conv_hull_thm}, as Theorem~\ref{conv_hull_thm} gives more precise information than Theorem~\ref{hp_slow_growth_thm} about the growth of the process on the event of coexistence.

As noted in the introduction, the proof of Theorem~\ref{conv_hull_thm} is actually a bootstrap argument that relies on Theorem~\ref{hp_slow_growth_thm} in an essential way, combining it with the results of Chapter~\ref{random_fpc_chap} to obtain a stronger result. Theorem~\ref{conv_hull_thm} will be a corollary of Theorem~\ref{advantage_sup_thm} below, in which we show that on the event of coexistence, the advantage of the fast species in the process configuration at time $t$ is bounded above by 1 in the limit as $t\to\infty$. 

Chapter~\ref{coex_finite_chap} is organized as follows. In Section~\ref{coex_finite:setup_sec} we first describe the general setup and notational conventions relating specifically to the two-type Richardson model, then we define the advantage for discrete configurations in order to state the main result of Chapter~\ref{coex_finite_chap}, Theorem~\ref{advantage_sup_thm}. In Section~\ref{coex_finite:proof_sec} we prove Theorem~\ref{advantage_sup_thm} through a sequence of lemmas, culminating in the proof of the main result via the strong Markov property, and then we show that Theorem~\ref{advantage_sup_thm} implies Theorem~\ref{conv_hull_thm}.


\section{Setup and Statement of Main Result}
\label{coex_finite:setup_sec}

Throughout this chapter, we assume that $\tmeasurepair = \pair{\tmeasure_1}{\tmeasure_2}$ is an \iid\ traversal measure on $\edges{\Zd}$ with exponentially distributed model traversal times, say $\cmttime{i} \sim \explaw{\exprate_i}$ for some $\exprate_1,\exprate_2>0$. As noted in Section~\ref{fpc_basics:ips_def_sec}, this implies that the two-type first-passage competition process $\process{}$ is a Markov process on the state space $\statespace=\statespace[1]$, called the \emph{two-type Richardson model}.
Note that the process $\process{}$ is time-homogeneous and a.s.\
right-continuous, and for finite starting configurations, it is
also a.s.\ constant except for isolated jumps, hence of \emph{pure
jump type} (e.g.\ as defined in \cite{Kallenberg:2002aa}).
By time scaling, we can without loss of generality fix the parameter $\lambda_1$ for species 1, and allow $\lambda_2$ to vary. Thus, we will assume throughout that $\lambda_1=1$ and $\lambda_2 = \lambda<1$.
We let $\mu_1 = \mu$ denote the corresponding shape function for species 1, in which case it follows that species 2's shape function is $\mu_2 = \expmean[1] \mu$. We let $\expnormpair = \pair{\mu_1}{\mu_2} = \pair{\mu}{\expmean[1] \mu}$ denote the traversal norm pair for the corresponding deterministic process. Observe that the limit shapes for the two species are then $\unitball{\mu_1} = \unitball{\mu}$ and $\unitball{\mu_2} = \exprate \unitball{\mu}$.

Our main goal for the rest of the chapter will be to prove Theorem~\ref{advantage_sup_thm} below. First we need a suitable definition for the advantage of species~1 in a discrete configuration $X\in\statespace$ as opposed to a continuous starting configuration in $\Rd$. Following the definition \eqref{config_advantage_def2_eqn} for the deterministic process, but making adjustments for lattice approximation, we define the \textdef{advantage of species 1} in a state $X\in\statespace$ with respect to the norm pair $\normpair = \pair{\mu_1}{\mu_2}$ by
\begin{align}
\label{state_advantage_def_eqn}
\stateadvantage{1}{\normpair}{X}
&\definedas \sup \setbuilderbar[1]{ \alpha\ge 0}
{\textstack{$\exists \dirsymb\in \Rd\setminus \setof{\vect 0}$ and
	 $\vect z\in \cubify{\stateset{1}{X}}$ with}
	 {$\cone{\mu_2}{\thpadv{2}{\dirvec}{\alpha}}{\vect z}{\dirvec} 
	 \cap \cubify{\stateset{2}{X}} = \emptyset$, where 
	 $\thpadv{2}{\dirvec}{\alpha} = \alpha \frac{\mu_1(\dirsymb)}{\mu_2(\dirsymb)}$}}\\
\notag\\
&= \sup \setbuilderbar[1]{\coneadvantage{1}{\normpair}{\cones}}
{\textstack{$\cones$ is a cone with $\apexset{\cones}\cap \cubify{\stateset{1}{X}} \ne \emptyset$}{and
	$\cones \cap \cubify{\stateset{2}{X}} = \emptyset$}}.\notag
\end{align}

\begin{remark}
\label{exp_norms_adv_rem}
For the norm pair $\expnormpair = (\mu,\expmean[1] \mu)$,  Corollary~\ref{exact_adv_cor} says that $\coneadvantage[2]{1}{\expnormpair}{\cone{\mu}{\thp}{\vect z}{\dirvec}} = \thp/\exprate$ for any $\dirvec\in \dirset{\Rd}$ and $\thp\in [0,1]$. Thus, if $\alpha\in [0,\expmean[1]]$, the advantage of species~1 in the $\mu$-cone $\cone{\mu}{\exprate\alpha}{\vect z}{\dirvec}$ is $\alpha$, and if $\exprate \in (0,1)$, the $\mu$-cone $\cone{\mu}{\thp}{\vect z}{\dirvec}$ is wide, critical, or narrow for species~1 according to whether $\thp>\exprate$, $\thp=\exprate$, or $\thp <\exprate$, respectively.
%
\end{remark}


Now we are ready to state the main result of Chapter~\ref{coex_finite_chap}.
Recall that for the two-type Richardson model with rates $\exprate_1$ and $\exprate_2$, we write $\coexstate{\exprate_1,\exprate_2}{X}$ for the event that coexistence occurs when the starting configuration is $X$.

\begin{thm}[Coexistence implies advantage is bounded above by 1 in the limit]
\label{advantage_sup_thm}
\newcommand{\thisstate}{X}
For any fixed $\exprate\in (0,1)$, if $\thisstate\in \statespace$ is finite, then almost surely on the event $\coexstate{1,\exprate}{\thisstate}$, we have
\[
\limsup_{t\to\infty} \stateadvantage[1]{1}{\expnormpair}{\processt{X}{t}} \le 1,
\]
where $\expnormpair = \pair{\mu}{\expmean[1] \mu}$ is the pair of shape functions for the two species, and the advantage of species~1 in the configuration $\processt{X}{t}$ is defined by \eqref{state_advantage_def_eqn}.
\end{thm}

The idea for the proof of Theorem~\ref{advantage_sup_thm} is that if at some random time $\sigma$ the advantage of species~1 in the configuration $\processt{X}{\sigma}$ is greater than 1, then by Part~\ref{wide_narrow_survival:wide_part} of Theorem~\ref{wide_narrow_survival_thm}, species~1 has a positive probability of conquering a cone from the configuration $\processt{X}{\sigma}$. Since the growth of species~1 in a cone is asymptotically as fast as its unrestricted growth,
this implies that species~1 will eventually occupy points
outside the ball $\reacheddset{\mu_2}{\vect 0}{(1+\epsilon)t}$,
precluding the possibility of coexistence by Theorem~\ref{hp_slow_growth_thm}. Thus, if there is some unbounded sequence of times $\sigma_1,\sigma_2,\dotsc$ such that at each time $\sigma_n$ the advantage of species~1 is at least $\alpha>1$, we can get a uniform positive lower bound on the probability that species~1 eventually wins from the configuration $\processt{X}{\sigma_n}$; applying the strong Markov property at the stopping times $\sigma_n$ will show that species~1 must in fact win if such a sequence exists.
We prove Theorem~\ref{advantage_sup_thm} in the following section by formalizing this argument via a sequence of lemmas, and we then show that Theorem~\ref{conv_hull_thm} follows as a corollary of Theorem~\ref{advantage_sup_thm}.


\section{Proof of Theorems~\ref{advantage_sup_thm} and \ref{conv_hull_thm}}
\label{coex_finite:proof_sec}

We start by defining several increasing subsets of the state space $\statespace = \statespace[1]$. Recall that $H\subseteq \statespace$ is \textdef{increasing} if for any state $X\in H$, the result of changing some of the sites in $X$ from 2's or 0's to 0's or 1's yields a configuration that is still in $H$.
Also recall that a configuration $X=\initconfig$ is \textdef{finite} if both $\initialset{1}$ and $\initialset{2}$ are finite, \textdef{fertile} if neither of the sets $\initialset{1}$ or $\initialset{2}$ surrounds the other, and \textdef{final} if $\initialset{1}\cup\initialset{2}=\Zd$.

\begin{definition}[Some increasing subsets of $\statespace$]

\mbox{}

\begin{description}
\item[$\setaasymb$] For each $\exprate\in(0,1)$ and $\alpha\in [0,\expmean[1])$, define
\[
\advset{\exprate}{\alpha} \definedas
\setof[1]{\text{finite, fertile $X\in\statespace$
such that $\stateadvantage{1}{\expnormpair}{X}> \alpha$}},
\]
where $\expnormpair = \pair{\mu}{\expmean[1] \mu}$, and for each $r>0$ let
\[
\advsetr{\exprate}{\alpha}{r} \definedas
\setbuilder[2]{X\in \advset{\exprate}{\alpha}}
{{\stateset{2}{X}}\subseteq r {\unitball{\mu}}}.
\]
Note that the map $r\mapsto \advsetr{\exprate}{\alpha}{r}$ is increasing as a function from $\Rplus$ to $2^{\statespace}$.

\item[$\setbbsymb$] For each $D\subseteq \dirset{\Rd}$, let
\[
\coneset{D} \definedas
\setof[1]{\textstack{final $X\in\statespace$ such that there exists $\dirvec\in D$,}
{$\vect z\in\Rd$, and $\thp>0$, with $\cone{\mu}{\thp}{\vect z}{\dirvec} \subseteq \intcubify{\stateset{1}{X}}$}},
\]
i.e.\ final (not finite) configurations in which species 1 has conquered an entire $\mu$-cone
with some axis in $D$.

\item[$\setccsymb$] For each $r>0$, let
\[
\escset{r} \definedas
\setof[2]{\text{finite $X\in\statespace$ such that $\cubify{\stateset{1}{X}}\not\subseteq r\unitball{\mu}$}},
\]
i.e.\ finite configurations in which species 1 has ``escaped" the ball $r\unitball{\mu}$, in that it occupies some site outside this ball.

\item[$\setddsymb$] Finally, let
\[
\winset \definedas
\setof[2]{\text{$X\in\statespace$ such that species 1 surrounds species 2}},
\]
i.e.\ the set of configurations (finite or not) in which species 1 has won. Observe that $\eventthat[1]{\processt{X}{\infty}\in \winset}$ is the event that species 1 wins when the process starts with the initial configuration $X$. Observe also that $\winset$ is absorbing, i.e.\ if $\processt{X}{t_0}\in \winset$, then $\processt{X}{t}\in \winset$ for all $t\ge t_0$.
\end{description}
\end{definition}

Very loosely speaking, the idea of the proof will be to show that for a given starting configuration $X\in\statespace$, we have ``$\setaasymb \implies \setbbsymb\implies \setccsymb \implies \setddsymb$" with probability bounded away from 0, and then to combine this result with the strong Markov property to show that a configuration with large advantage can't happen infinitely often. The following lemma is the first step in this chain of implications, i.e.\ ``$\setaasymb \implies \setbbsymb$."

\begin{lem}[Large advantage implies species 1 can conquer a cone]
\label{adv_implies_cone_lem}
Fix $\exprate<1$ and $\alpha\in (1,\exprate^{-1})$. There exists a finite set $D\subseteq \dirset{\Rd}$ (depending on $\exprate$ and $\alpha$) and a positive number $p_0 = p_0(\exprate,\alpha)$ such that
\[
\inf_{X\in\advset{\exprate}{\alpha}}
\Pr \eventthat[1]{\processt{X}{\infty} \in \coneset{D}}
\ge p_0.
\]
\end{lem}

\begin{proof}
\newcommand{\coversetsymb}{U}
\newcommand{\coverset}[1]{\coversetsymb_{#1}}
\newcommand{\opencoverset}[1]{\interior{\coverset{#1}}}
\newcommand{\finitedirset}{D}
\newcommand{\conestate}[3]{X_{#1}^{(#2,#3)}} 

Let $\alpha' \definedas (1+\alpha)/2$, so $1<\alpha'<\alpha$. For each $\dirvec\in \dirset[2]{\Rd}$, define
\begin{align*}
\coverset{\dirvec} &\definedas
	\setbuilder[1]{\dirvec[1]\in \dirset[2]{\Rd}}
	{\cone{\mu}{\exprate \alpha'}{\vect 0}{\dirvec[0]} \subseteq
		\cone{\mu}{\exprate \alpha}{\vect 0}{\dirvec[1]}}.
\end{align*}

\begin{claim}
\label{adv_implies_cone:open_cover_claim}
The collection $\setof[1]{\opencoverset{\dirvec}}_{\dirvec\in \dirset{\Rd}}$ is an open cover of $\dirset[2]{\Rd}$.
\end{claim}

\begin{proof}[Proof of Claim~\ref{adv_implies_cone:open_cover_claim}]
\newcommand{\thisuvec}[1]{\unitvecnorm{\dirsymb[#1]}{\mu}}
\newcommand{\scoverset}[1]{\widehat{\coversetsymb}_{#1}} 
\newcommand{\socoverset}[1]{\interior{\scoverset{#1}}} 

We have to show that $\setof[1]{\opencoverset{\dirvec}}_{\dirvec\in \dirset{\Rd}}$ is a cover by showing that $\dirvec\in \opencoverset{\dirvec}$. 
It follows from the definition of $\mu$-cones that
\begin{equation}
\label{adv_implies_cone:coverset_equiv_eqn}
\coverset{\dirvec} \supseteq \setbuilder[1]{\dirvec[1]\in \dirset[2]{\Rd}}
	{\reacheddset{\mu}{\thisuvec{1}}{\exprate\alpha'}
	\subseteq \reacheddset{\mu}{\thisuvec{0}}{\exprate\alpha}}.
\end{equation}
Let
$
\scoverset{\dirvec} \definedas
\setbuilder[2]{\thisuvec{1}\in \sphere{d-1}{\mu}}{\dirvec[1]\in \coverset{\dirvec[0]}}.
$
Since $\dirset[2]{\Rd}$ is homeomorphic to $\sphere{d-1}{\mu}$ under the unit vector projection map, it will suffice to show that $\thisuvec{0}\in \socoverset{\dirvec}$. We claim that in fact
\begin{equation}
\label{adv_implies_cone:small_ball_eqn}
\sphere{d-1}{\mu}\cap \reacheddset[2]{\mu}{\thisuvec{0}}{\exprate(\alpha-\alpha')}
\subseteq \scoverset{\dirvec}.
\end{equation}
To see that \eqref{adv_implies_cone:small_ball_eqn} holds, let $\thisuvec{1}\in \sphere{d-1}{\mu}$ with $\distnew{\mu}{\thisuvec{0}}{\thisuvec{1}}\le \exprate (\alpha-\alpha')$. Then by the triangle inequality, if $\vect x\in\Rd$ and $\distnew{\mu}{\vect x}{\thisuvec{0}}\le \exprate\alpha'$, we have
\begin{align*}
\distnew{\mu}{\vect x}{\thisuvec{1}}
\le \distnew{\mu}{\vect x}{\thisuvec{0}} + \distnew{\mu}{\thisuvec{0}}{\thisuvec{1}}
\le \exprate\alpha' + \exprate (\alpha-\alpha')
= \exprate\alpha.
\end{align*}
Therefore, $\reacheddset{\mu}{\thisuvec{1}}{\exprate\alpha'} \subseteq \reacheddset{\mu}{\thisuvec{0}}{\exprate\alpha}$, so $\thisuvec{1}\in \scoverset{\dirvec}$ by the characterization \eqref{adv_implies_cone:coverset_equiv_eqn}. Since $\thisuvec{1}\in \sphere{d-1}{\mu}\cap \reacheddset[2]{\mu}{\thisuvec{0}}{\exprate(\alpha-\alpha')}$ was arbitrary,
the containment \eqref{adv_implies_cone:small_ball_eqn} holds as claimed.
Finally, note that since $\alpha'<\alpha$, the open ball $\interior{\reacheddset[2]{\mu}{\thisuvec{0}}{\exprate(\alpha-\alpha')}}$ is nonempty and contains $\thisuvec{0}$, so \eqref{adv_implies_cone:small_ball_eqn} implies that $\thisuvec{0}\in \socoverset{\dirvec}$. Therefore, $\dirvec \in \opencoverset{\dirvec}$, so $\setof[1]{\opencoverset{\dirvec}}_{\dirvec\in \dirset{\Rd}}$ is an open cover of $\dirset{\Rd}$.
\end{proof}

Since $\dirset[2]{\Rd} \cong \sphere{d-1}{\lspace{2}{d}}$ is compact, Claim~\ref{adv_implies_cone:open_cover_claim} implies that there is a finite subcover $\coverset{\dirvec_1},\dotsc,\coverset{\dirvec_n}$; let $\finitedirset \definedas \setof{\dirvec_1,\dotsc,\dirvec_n}$. For each $j\in \setof{1,\dotsc,n}$ and each pair $(\vect v,\vect z)\in\Zd\times \Rd$ with $\vect z\in\cubify{\vect v}$, define the state $\conestate{j}{\vect v}{\vect z}\in \statespace$ by
\[
\conestate{j}{\vect v}{\vect z} \equiv
\pair[1]{\vect v}
{\,\Zd\setminus \lat[2]{\cone{\mu}{\exprate \alpha'}{\vect z}{\dirvec_j}}}.
\]
By Corollary~\ref{exact_adv_cor}, for any $j$ and $\vect z$, the advantage of species~1 in $\cone{\mu}{\exprate \alpha'}{\vect z}{\dirvec_j}$ is $\alpha'>1$, and by Lemma~\ref{unif_cone_survival_lem}, we have
\[
p_j(\lambda,\alpha) \definedas \inf_{\vect z\in \cubify{\vect 0}}
\Pr \eventthat[1]{\processt{\conestate{j}{\vect 0}{\vect z}}{\infty} \in \coneset{\dirvec_j}}>0.
\]
By translation invariance, for any index $j$ and any pair $(\vect v,\vect z)$ with $\vect z\in \cubify{\vect v}$ we have
\[
\Pr \eventthat[1]{\processt{\conestate{j}{\vect v}{\vect z}}{\infty} \in \coneset{\dirvec_j}}
= \Pr \eventthat[1]{\processt{\conestate{j}{\vect 0}{\vect z-\vect v}}{\infty} \in \coneset{\dirvec_j}}
\ge p_j(\lambda,\alpha).
\]
Thus, if we define $p_0$ by
\[
p_0(\exprate, \alpha) \definedas
\min_{j\in \setof{1,\dotsc,n}} p_j(\lambda,\alpha)>0,
\]
then
\[
\Pr \eventthat[1]{\processt{\conestate{j}{\vect v}{\vect z}}{\infty} \in \coneset{\finitedirset}}
\ge p_0
\]
for all $j\in \setof{1,\dotsc,n}$ and all pairs $(\vect v,\vect z)\in\Zd\times \Rd$ with $\vect z\in \cubify{\vect v}$ (equivalently $\vect v\in \lat{\vect z}$).

Now suppose $X\in \advset{\exprate}{\alpha}$. Then $\stateadvantage{1}{\expnormpair}{X} > \alpha$, so by definition \eqref{state_advantage_def_eqn} and Lemma~\ref{lattice_cube_properties_lem} there exists $\vect z\in \cubify{\stateset{1}{X}}$ and $\dirvec[1]\in \dirset{\Rd}$ such that $\stateset{2}{X} \cap \lat[2]{ \cone{\mu}{\exprate \alpha}{\vect z}{\dirvec[1]}} = \emptyset$.
Since $\vect z\in \cubify{\stateset{1}{X}}$, there is some $\vect v\in \lat{\vect z}\cap \stateset{1}{X}$. Since the sets $\coverset{\dirvec_1},\dotsc,\coverset{\dirvec_n}$ cover $\dirset{\Rd}$, we have $\dirvec[1]\in \coverset{\dirvec_j}$ for some $j$. Therefore, $\cone{\mu}{\exprate \alpha'}{\vect z}{\dirvec_j} \subseteq \cone{\mu}{\exprate \alpha}{\vect z}{\dirvec[1]}$ by the definition of $\coverset{\dirvec_j}$, so $\stateset{2}{X} \cap \lat[2]{ \cone{\mu}{\exprate \alpha'}{\vect z}{\dirvec_j}} = \emptyset$. Thus,
\[
\stateset{1}{X}\supseteq \setof{\vect v}
\quad\text{and}\quad
\stateset{2}{X} \subseteq \Zd\setminus \lat[2]{ \cone{\mu}{\exprate \alpha'}{\vect z}{\dirvec_j}}.
\]
That is, $X\ge \conestate{j}{\vect v}{\vect z}$, so since $\coneset{\finitedirset}$ is increasing, Lemma~\ref{increasing_sets_lem} implies that
\[
\Pr \eventthat[1]{\processt{X}{\infty} \in \coneset{\finitedirset}}
\ge \Pr \eventthat[1]{\processt{\conestate{j}{\vect v}{\vect z}}{\infty} \in \coneset{\finitedirset}}
\ge p_0,
\]
which proves the lemma.
\end{proof}

\begin{thmremark}
The assumption that $X$ is finite in the definition of $\advset{\exprate}{\alpha}$ was made to simplify the statement of Lemma~\ref{sp1_wins_lem} and the proof of Theorem~\ref{advantage_sup_thm} below, but this assumption was not used in the proof of Lemma~\ref{adv_implies_cone_lem}. Thus, Lemma~\ref{adv_implies_cone_lem} remains true if we drop the finiteness assumption in the definition of $\advset{\exprate}{\alpha}$, meaning that we get the same uniform lower bound on the probability of species~1 conquering a cone with axis in $D$ from any starting configuration with advantage $>\alpha$, even if species~2 initially occupies an infinite set.
\end{thmremark}

The next lemma is the second implication in the chain, ``$\setbbsymb\implies \setccsymb$."

\begin{lem}[If species 1 conquers a cone, then it escapes from species 2]
\label{cone_implies_headstart_lem}
\newcommand{\speedfactor}{a}

Let $D$ be any subset of $\dirset{\Rd}$. Then for any finite $X\in\statespace$ and any positive number $\speedfactor<1$,
\[
\eventthat[1]{\processt{X}{\infty}\in \coneset{D}}
\assubseteq
\eventthat[1]{\processt{X}{t} \in \escset{\speedfactor t}
	\text{ for all large $t$}}.
\]
\end{lem}

\begin{proof}
\newcommand{\thiscone}{\cones}
\newcommand{\thist}{t_0}
\newcommand{\speedfactor}{a}
\newcommand{\newspf}{\speedfactor'}
\newcommand{\bigt}{\thist'}
\newcommand{\shapet}{s_\spp}
\newcommand{\thisx}{\vect x_t}

Fix a finite initial configuration $X$, and suppose that $\eventthat[1]{\processt{X}{\infty}\in \coneset{D}}$ occurs. Then there is some $\mu$-cone $\thiscone = \cone{\mu}{\thp}{\vect z}{\dirvec}$ such that $\crreachedset{1}{X}{\Rd}{\infty} \supseteq \thiscone$, where $\vect z\in\Rd$, $\dirvec\in D$, and $\thp>0$ are random. Let $\thist = \entangledctime{1}{X}{\vect z}<\infty$, i.e.\ $\thist$ is the (random) entangled covering time of $\vect z$ by species~1 (defined in \eqref{continuum_entangled_times_def_eqn}).

\begin{claim}
\label{cone_implies_headstart:domination_claim}
$\crreachedset{1}{X}{\Rd}{t} \supseteq \crreachedset{\tmeasure_1}{\vect z}{\thiscone}{t-\thist}$ for all $t\ge \thist$.
\end{claim}

\begin{proof}[Proof of Claim~\ref{cone_implies_headstart:domination_claim}]
\newcommand{\thisz}{\vect z}
\newcommand{\latz}{\lat{\thisz}}
\newcommand{\avertex}{\vect v}
\newcommand{\conept}{\vect x}
\newcommand{\thispath}{\gamma}
\newcommand{\allpaths}{{\gamma}}
\newcommand{\geov}{\gamma_{\avertex}}

\newcommand{\latcone}{\lat{\thiscone}}
\newcommand{\uset}{S}
\newcommand{\iset}{\initialset{1}}
\newcommand{\fset}{\finalset{1}}
\newcommand{\cubefset}{\intcubify{\fset}}

Let $\fset = \reachedset{1}{X}{\infty}$ be species~1's finally conquered set, and let $\allpaths$ be the union of all $\rptmetric[1]{\fset}$-geodesics from $\iset$ to some vertex in $\latz$. That is, $\allpaths \definedas \bigcup_{\avertex\in\latz} \geov$, where $\geov$ is the $\rptmetric[1]{\fset}$-geodesic from $\iset$ to $\avertex$. Then $\rptime{1}{\geov}{\iset}{\avertex} = \rptime{1}{\fset}{\iset}{\avertex}$ for each $\avertex\in \latz$, so
\begin{equation}
\label{cone_implies_headstart:ctime_ineq_eqn}
\ctime{1}{\allpaths}{\iset}{\thisz}
= \sup_{\avertex\in\latz} \rptime{1}{\allpaths}{\iset}{\avertex}
\le \sup_{\avertex\in\latz} \rptime{1}{\geov}{\iset}{\avertex}
= \sup_{\avertex\in\latz} \rptime{1}{\fset}{\iset}{\avertex}
= \ctime{1}{\fset}{\iset}{\thisz}
=\thist.
\end{equation}
Now set $\uset \definedas \allpaths \cup \thiscone$, and note that $\uset \subseteq \cubefset$ by assumption. By Lemma~\ref{chaining_lem} (chaining) and \eqref{cone_implies_headstart:ctime_ineq_eqn}, for any $\conept\in \thiscone$ we have
\begin{equation}
\label{cone_implies_headstart:chaining_eqn}
\ctime{1}{\uset}{\iset}{\conept} \le \ctime{1}{\allpaths}{\iset}{\thisz}
	+ \ctime{1}{\thiscone}{\thisz}{\conept}
\le \thist + \ctime{1}{\thiscone}{\thisz}{\conept}.
\end{equation}
Since $\uset\supseteq \thiscone$, \eqref{cone_implies_headstart:chaining_eqn} implies that for any $t\ge \thist$,
\begin{align}
\label{cone_implies_headstart:rset_containment_eqn}
\crreachedset{\tmeasure_1}{\iset}{\uset}{t}
= \setbuilder[1]{\conept\in \uset}{\ctime{1}{\uset}{\iset}{\conept} \le t}
&\supseteq \setbuilder[1]{\conept\in \uset}{\thist +\ctime{1}{\thiscone}{\thisz}{\conept} \le t}
	\notag\\
&\supseteq \setbuilder[1]{\conept\in \thiscone}{\ctime{1}{\thiscone}{\thisz}{\conept} \le t-\thist}
= \crreachedset{\tmeasure_1}{\thisz}{\thiscone}{t-\thist}.
\end{align}
Finally, since $\cubefset \supseteq \uset$, it follows from \eqref{cone_implies_headstart:rset_containment_eqn} and the definition of $\crreachedset{1}{X}{\Rd}{t}$ (cf.\ \eqref{entangled_covering_proc_equiv_eqn}) that
\[
\crreachedset{1}{X}{\Rd}{t}
=\crreachedset{\tmeasure_1}{\iset}{\cubefset}{t}
\supseteq \crreachedset{\tmeasure_1}{\iset}{\uset}{t}
\supseteq \crreachedset{\tmeasure_1}{\thisz}{\thiscone}{t-\thist}
\quad\text{for all $t\ge\thist$,}
\]
which proves Claim~\ref{cone_implies_headstart:domination_claim}.
\end{proof}

Now by Corollary~\ref{simultaneous_cone_shape_thm_cor}, the Shape Theorem for $\tmeasure_1$ holds simultaneously in all $\mu$-cones and hence, in particular, in the cone $\thiscone$. That is, for any $\spp\in (0,1)$, there is almost surely some (random) $\shapet<\infty$ such that
\begin{equation}
\label{cone_implies_headstart:shape_thm_eqn}
\crreachedset{\tmeasure_1}{\vect z}{\thiscone}{s}
\supseteq
\rreacheddset[2]{\mu}{\vect z}{\thiscone}{(1-\spp)s}
\text{ for all $s\ge \shapet$.}
\end{equation}
We will take $s=t-\thist$ in \eqref{cone_implies_headstart:shape_thm_eqn} and combine this with Claim~\ref{cone_implies_headstart:domination_claim} to finish the proof.
Choose some $\newspf \in (\speedfactor,1)$, and choose some $\spp>0$ small enough that $1-\spp>\newspf$.
Now define the random time
\[
\bigt\definedas
\frac{\mu(\vect z) + (1-\spp)\thist}{1-\spp-\newspf} \in [\thist,\infty),
\]
and for each $t\ge 0$ define the random point
\[
\thisx \definedas \vect z + \parens[2]{\mu(\vect z) + \newspf t} \munitvec \in\thiscone.
\]
Then by the triangle inequality and the definitions of $\thisx$ and $\bigt$, 
\begin{equation}
\label{cone_implies_headstart:x_dist_eqn}
\mu(\thisx)\ge \newspf t
\text{ for all $t\ge 0$,\quad and\quad}
\distnew{\mu}{\vect z}{\thisx} \le (1-\spp) (t-\thist)
\text{ for all $t\ge \bigt$.}
\end{equation}
Taking $s=t-\thist$ in \eqref{cone_implies_headstart:shape_thm_eqn}, it follows from Claim~\ref{cone_implies_headstart:domination_claim} and the second inequality in \eqref{cone_implies_headstart:x_dist_eqn} that 
almost surely on the event $\eventthat[1]{\processt{X}{\infty}\in \coneset{D}}$, we have
\begin{equation}
\label{cone_implies_headstart:esc_event_eqn}
\crreachedset{1}{X}{\Rd}{t}
\supseteq \crreachedset{\tmeasure_1}{\vect z}{\thiscone}{t-\thist}
\supseteq \rreacheddset[2]{\mu}{\vect z}{\thiscone}{(1-\spp)(t-\thist)}
\ni \thisx
\text{ for all $t\ge \shapet+\bigt$.}
\end{equation}
Finally, since $\newspf>\speedfactor$, the first inequality in \eqref{cone_implies_headstart:x_dist_eqn} shows that on the event \eqref{cone_implies_headstart:esc_event_eqn} we have $\processt{X}{t} \in\escset{\speedfactor t}$ for all $t\ge \shapet+\bigt$. 
\end{proof}

The final step in the implication chain is ``$\setccsymb\implies \setddsymb$"; this is precisely the content of \Haggstrom and Pemantle's Theorem~\ref{hp_slow_growth_thm} \cite[Lemma~5.2]{Haggstrom:2000aa}. Given this, it then follows trivially from the preceding lemma that ``$\setbbsymb \iff\setccsymb \iff \setddsymb$". More precisely, we have the following result.

%

\begin{lem}[Species 1 conquers a cone if and only if it wins]
\label{equivalent_events_lem}
For any finite $X\in\statespace$ and any subset $D\subseteq \dirset{\Rd}$,
\[
\eventthat[1]{\processt{X}{\infty}\in \coneset{D}}
\asequal
\eventthat[1]{\processt{X}{\infty}\in \winset}.
\]
\end{lem}

\begin{proof}
\newcommand{\speedfactor}{a}

Recall our assumption that the process has rates $\lambda_1=1$ and $\lambda_2=\lambda$ for some $\lambda\in (0,1)$. Choose any $\speedfactor\in (\lambda,1)$, and set $\spp \definedas \lambda^{-1}\speedfactor -1$. Since $\speedfactor >\lambda$, we have $\spp>0$, and Theorem~\ref{hp_slow_growth_thm} implies that if there exist arbitrarily large $t$ for which $\cubify{\reachedset{1}{X}{t}}$ is not contained in the ball $\speedfactor t\unitball{\mu} = \reacheddset[2]{\mu_2}{\vect 0}{\lambda^{-1}\speedfactor t} = \reacheddset[2]{\mu_2}{\vect 0}{(1+\spp) t}$, then species~1 almost surely wins. That is,
\[
\eventthat[1]{\processt{X}{t} \in \escset{\speedfactor t}
	\text{ for all large $t$}}
\assubseteq
\eventthat[1]{\processt{X}{\infty}\in \winset}.
\]
On the other hand, if species~1 wins, then it conquers everything outside some finite region, so we trivially have
\[
\eventthat[1]{\processt{X}{\infty}\in \winset}
\assubseteq
\eventthat[1]{\processt{X}{\infty}\in \coneset{D}}
\text{ for any $D\subseteq \dirset{\Rd}$.}
\]
Since $\speedfactor<1$, it then follows trivially from Lemma~\ref{cone_implies_headstart_lem} that
\[
\eventthat[1]{\processt{X}{\infty}\in \coneset{D}}
\asequal
\eventthat[1]{\processt{X}{t} \in \escset{\speedfactor t}
	\text{ for all large $t$}}
\asequal
\eventthat[1]{\processt{X}{\infty}\in \winset}.
\qedhere
\]
\end{proof}

%
%

Combining Lemmas~\ref{adv_implies_cone_lem} and \ref{equivalent_events_lem} immediately yields the following.

\begin{lem}[Large advantage implies species 1 can win]
\label{sp1_wins_lem}
For any $\lambda\in (0,1)$ and $\alpha\in (1,\lambda^{-1})$, we have
\[
\inf_{X\in\advset{\exprate}{\alpha}}
\Pr \eventthat[1]{\processt{X}{\infty} \in \winset}
\ge p_0,
\]
where $p_0  = p_0(\lambda,\alpha) > 0$ is the constant from Lemma~\ref{adv_implies_cone_lem}.
\end{lem}

Thus, if $X$ is any finite configuration in which species~1 has advantage greater than
$\alpha$, the probability that species~1 eventually wins from the
initial configuration $X$ is at least $p_0$. We now show that if we
consider only configurations $X\in\advset{\exprate}{\alpha}$ such that
species~2 is contained in a ball of some fixed radius $r$, then
there is some fixed \emph{finite} time (depending on $r$) such that
the probability that species~1 wins by this \emph{fixed} time is
also uniformly bounded below.

\begin{lem}[Positive probability of winning within a fixed time]
\label{sp1_wins_soon_lem}
Let $p_0  = p_0(\lambda,\alpha)$ be the constant from Lemma~\ref{adv_implies_cone_lem}. For any positive $p<p_0$, there is an increasing function
$N_p\colon (0,\infty)\to (0,\infty)$ such that
\[
    \inf_{X\in \advsetr{\exprate}{\alpha}{r}}
    \Pr \eventthat[1]{\processt{X}{N_p(r)} \in \winset}
    \ge p.
\]
\end{lem}

\begin{proof}
\newcommand{\wintime}[2]{n_{#1}(#2)} 
\newcommand{\Wintime}[2]{N_{#1}(#2)} 
\newcommand{\Wintimefcn}[1]{N_{#1}} 
\newcommand{\advsetrr}[3]{{\setaasymb}_{#1}^{#2} \parens{#3}^*}
\newcommand{\finitize}[1]{{#1}^*} 

First note that
\[
\eventthat[1]{\processt{X}{\infty} \in \winset}
= \bigcup_{n\in\N} \eventthat[1]{\processt{X}{n} \in \winset},
\]
so by Lemma~\ref{sp1_wins_lem} and continuity of measure, for any $X\in \advset{\exprate}{\alpha}$ and $p<p_0$ we can find $\wintime{p}{X}\in\N$ so that
\[
\Pr \eventthat[1]{\processt{X}{\wintime{p}{X}}\in \winset}\ge p.
\]
Recall that
\[
\advsetr{\exprate}{\alpha}{r} =
\setbuilder[2]{X\in \advset{\exprate}{\alpha}}
{{\stateset{2}{X}}\subseteq r{\unitball{\mu}}}.
\]
Now let $\advsetrr{\exprate}{\alpha}{r}$ be the subset of $\advsetr{\exprate}{\alpha}{r}$ consisting of configurations in which \emph{both} species are contained in a ball of radius $r$, i.e.
\[
\advsetrr{\exprate}{\alpha}{r} \definedas
\setbuilder[2]{X\in \advset{\exprate}{\alpha}}
{{\nonzeroset{X}}\subseteq r{\unitball{\mu}}}.
\]
Then $\advsetrr{\exprate}{\alpha}{r}$ is finite, and if we set
\[
\Wintime{p}{r} \definedas 
\max_{X\in \advsetrr{\exprate}{\alpha}{r}} \wintime{p}{X},
\]
it follows that
\[
    \inf_{X\in \advsetrr{\exprate}{\alpha}{r}}
    \Pr \eventthat[1]{\processt{X}{\Wintime{p}{r}} \in \winset}
    \ge p.
\]
That is, the desired property holds with the finite set $\advsetrr{\exprate}{\alpha}{r}$ in place of the infinite set $\advsetr{\exprate}{\alpha}{r}$.
Now, for arbitrary $X\in \advsetr{\exprate}{\alpha}{r}$, define $\finitize{X}\in \advsetrr{\exprate}{\alpha}{r}$ by
\[
\finitize{X}(\vect v)=
\begin{cases}
X(\vect v) & \text{if } \vect v\in r{\unitball{\mu}}\\
    0 & \text{otherwise}.
\end{cases}
\]
Then since $X$ contains no 2's outside the ball $r{\unitball{\mu}}$ by assumption, the modified state $\finitize{X}$ is obtained from $X$ by (possibly) changing some 1's to 0's, and doing nothing to the set of 2's, so we have $\stateset{1}{\finitize{X}} \subseteq \stateset{1}{X}$ and $\stateset{2}{\finitize{X}} = \stateset{2}{X}$. That is, $\finitize{X}\le X$, and since $\winset$ is increasing, Lemma~\ref{increasing_sets_lem} then implies that
\[
    \Pr \eventthat[1]{\processt{X}{\Wintime{p}{r}} \in \winset}
    \ge  
    \Pr \eventthat[1]{\processt{\finitize{X}}{\Wintime{p}{r}} \in \winset}
    \ge p.
\]
Since $X\in \advsetr{\exprate}{\alpha}{r}$ was arbitrary, this provides the desired lower bound. Finally, observe that the function $\Wintimefcn{p}$ is increasing because the family of sets $\setof[2]{\advsetrr{\exprate}{\alpha}{r}}_{r\ge 0}$ is increasing in $r$.
\end{proof}

Now we are ready to prove Theorem~\ref{advantage_sup_thm}.
The idea is to combine Lemma~\ref{sp1_wins_soon_lem} with the strong
Markov property and the Shape Theorem for the slow species, in a
manner similar to the proof of Theorem~\ref{hp_slow_growth_thm} \cite[Lemma~5.2]{Haggstrom:2000aa} from Proposition~\ref{hp_point_ball_prop} \cite[Proposition~2.2]{Haggstrom:2000aa}.
However, a bit more care is needed in our case
because whereas 
Proposition~\ref{hp_point_ball_prop} 
states that the probability that
species~1 wins approaches 1 as a particular starting configuration
is scaled up, our Lemma~\ref{sp1_wins_soon_lem} merely provides a
nonzero lower bound on this probability for a certain set of
starting configurations. This is why we need Lemma~\ref{sp1_wins_soon_lem} rather than proving Theorem~\ref{advantage_sup_thm} directly from Lemma~\ref{sp1_wins_lem}.

\begin{proof}[Proof of Theorem~\ref{advantage_sup_thm}]
\newcommand{\mainevent}{G(\alpha)}
\newcommand{\maineventn}{G_n(\alpha)}
\newcommand{\neweventn}{G_n'(\alpha)}
\newcommand{\thisstate}{X} 
\newcommand{\genstate}{X'} 
\newcommand{\posprob}{q_0}
\newcommand{\smallprob}{q}

\newcommand{\Wintime}[2]{N_{#1}(#2)} 
\newcommand{\Wintimefcn}[1]{N_{#1}} 

\newcommand{\stoptime}[1]{\sigma_{#1}} 
\newcommand{\killtime}{\kappa} 

\newcommand{\filtration}[1]{\mathcal{F}_{#1}} 

\newcommand{\projop}[1]{\pi_{#1}} 

It will suffice to prove that for every $\alpha\in (1,\expmean[1])$, the
event
\[
\mainevent \definedas
\eventthat[1]{\text{$\exists \setof{t_n}_{n\in\N}$ with $t_n\to\infty$
such that $\processt{\thisstate}{t_n}$ is fertile and
$\stateadvantage[1]{1}{\expnormpair}{\processt{\thisstate}{t_n}} > \alpha$}}
\]
has probability 0. Observe that $\mainevent = \eventthat{\maineventn \text{ i.o.}}$, where
\[
\maineventn \definedas
\eventthat[3]{\exists t\in [n,n+1) \text{ such that }
\processt{\thisstate}{t}\in \advset{\exprate}{\alpha}}.
\]
Recall that the limit shape for $\tmeasure_2$ is $\exprate \unitball{\mu}$. Since species~2's entangled version $\reachedset{2}{\thisstate}{t}$ is dominated by its disentangled version $\reachedset{\tmeasure_2}{\stateset{2}{\thisstate}}{t}$, and $\lambda<1$, the Shape Theorem for $\tmeasure_2$ implies that
\[
\Pr \eventthat[3]{\reachedset{2}{\thisstate}{t} \subseteq t\unitball{\mu}
\text{ for all large $t$}} = 1.
\]
Thus, if we define
\[
\neweventn \definedas\maineventn \cap
\eventthat[3]{\reachedset{2}{\thisstate}{t} \subseteq t\unitball{\mu}\ \forall t\ge n},
\]
it follows that $\eventthat{\maineventn \text{ i.o.}} \asequal \eventthat{\neweventn \text{ i.o.}}$. Note that by the definitions of $\maineventn$ and $\advsetr{\exprate}{\alpha}{r}$ we have
\[
\neweventn \subseteq
\eventthat[1]{\exists t\in [n,n+1) \text{ such that }
\processt{\thisstate}{t} \in \advsetr{\exprate}{\alpha}{t}}.
\]

We proceed by contradiction. Suppose $\Pr \parens[2]{\mainevent} = \posprob>0$ and hence $\Pr \eventthat{\neweventn \text{ i.o.}} = \posprob$. Fix $\smallprob\in (0,\posprob)$. Then it follows from continuity of measure that given any $\ell\ge 0$ there exists
$r(\ell)>\ell$ such that
\[
\Pr \eventthat[1]{\text{$\exists n\in [\ell, r(\ell)]$ such that $\neweventn$ occurs}}
\ge \smallprob.
\]
Now fix $p\in (0,p_0)$, where $p_0 = p_0(\exprate,\alpha)$ is the constant from Lemma~\ref{adv_implies_cone_lem}, and let $\Wintimefcn{p}$ be the increasing function defined in Lemma~\ref{sp1_wins_soon_lem}. We recursively define a sequence of disjoint intervals $I_k=[\ell_k,r_k)$, $k\in\N$, as follows:
\begin{itemize}
\item Set $\ell_0=0$;

\item For each $k\ge 0$, let $r_k\definedas 1+r(\ell_k)$ and $\ell_{k+1}\definedas r_k+\Wintime{p}{r_k}$.
\end{itemize}
Now for each $k\in\N$ define
\[
\stoptime{k} \definedas
\inf \setbuilder[1]{t\in [\ell_k,r_k)}
{\processt{\thisstate}{t} \in \advsetr{\exprate}{\alpha}{t}},
\]
with the convention that $\inf\emptyset =\infty$. Since the process
$\processt{}{t}$ is a.s.\ right-continuous, the time $\stoptime{k}$ is
optional with respect to the induced right-continuous filtration
$\filtration{} = \setof{\filtration{t}}_{t\ge 0}$. Moreover, it is clear that
\[
\eventthat{\stoptime{k}<\infty} =
\eventthat[3]{\exists t\in [\ell_k,r_k) \text{ such that }
\processt{\thisstate}{t}\in  \advsetr{\exprate}{\alpha}{t}},
\]
and it follows from the above definitions that for each $k\in\N$ we
have
\[
\Pr\eventthat{\stoptime{k}<\infty}\ge \smallprob.
\]
We introduce one further stopping time, the time that species~2 first gets surrounded by species~1:
\[
\killtime \definedas \inf \setbuilder[1]{t>0}{\processt{\thisstate}{t} \in \winset}.
\]
Since $\stoptime{k}\ge \ell_k$ by definition, and species~2 must still be alive at time $\stoptime{k}$ whenever $\stoptime{k}$ is finite (this follows by construction, from the definition of $\advset{\exprate}{\alpha}$), for each $k\in\N$ we have
\[
\eventthat[2]{\processt{\thisstate}{\ell_{k+1}}\in \winset}
\cap \eventthat{\stoptime{k}<\infty}
\subseteq \eventthat{\killtime\in (\ell_k,\ell_{k+1}]},
\]
and hence
\begin{equation}
\label{advantage_sup:winning_time_eqn}
\sum_{k=0}^\infty \Pr \parens[1]{
	\eventthat[2]{\processt{\thisstate}{\ell_{k+1}}\in \winset}
	\cap \eventthat{\stoptime{k}<\infty}}
\le \Pr \eventthat{\killtime<\infty}\le 1.
\end{equation}
On the other hand, we now use the strong Markov property to show
that Lemma~\ref{sp1_wins_soon_lem} implies that each summand on the
left-hand side of \eqref{advantage_sup:winning_time_eqn} is bounded below by
$p \smallprob>0$, which provides the desired contradiction.

For each $t\ge 0$, let $\projop{t}\colon \processspace \to \statespace$ be the natural projection operator defined by $\projop{t}\zeta_\cdot = \zeta_t$, and let $\shiftop{t}\colon \processspace\to \processspace$ be the natural shift operator defined by $(\shiftop{t} \zeta_\cdot)_s = \zeta_{t+s}$. Recall that $\ell_{k+1}=r_k+\Wintime{p}{r_k}$. Since $\winset$ is absorbing, and $\stoptime{k}<r_k$ on the event $\eventthat{\stoptime{k}<\infty}$, we have
\begin{align*}
\eventthat[1]{\processt{\thisstate}{r_k+\Wintime{p}{r_k}}\in \winset}
	\cap \eventthat{\stoptime{k}<\infty}
&\supseteq \eventthat[1]{\processt{\thisstate}{\stoptime{k}+\Wintime{p}{r_k}}\in \winset}
	\cap \eventthat{\stoptime{k}<\infty}\\
&= \eventthat[1]{\shiftop{\stoptime{k}} \process{\thisstate}
	\in \projop{\Wintime{p}{r_k}}^{-1} \winset}
	\cap \eventthat{\stoptime{k}<\infty}.
\end{align*}
Since the event $\eventthat{\stoptime{k}<\infty}$ is in the $\sigma$-field $\filtration{\stoptime{k}} = \setbuilder[2]{A\in \sigmafield}{A\cap \eventthat{\stoptime{k}\le t} \in\filtration{t}, t\ge 0}$, we have
\begin{multline*}
\Pr \parens[1]{
	\eventthat[1]{\shiftop{\stoptime{k}} \process{\thisstate}
	\in \projop{\Wintime{p}{r_k}}^{-1} \winset}
	\cap \eventthat{\stoptime{k}<\infty}}\\
= \int\displaylimits_{\eventthat{\stoptime{k}<\infty}}
	\cprobbrack[1]{\shiftop{\stoptime{k}} \process{\thisstate}
	\in \projop{\Wintime{p}{r_k}}^{-1} \winset}{\filtration{\stoptime{k}}}
	\, d(\Pr).
\end{multline*}
Now, since the process $\process{}$ is of pure jump type, the
strong Markov property is valid at any stopping time, so we can
apply it at $\stoptime{k}$ to compute the integrand:
\[
\cprobbrack[1]{\shiftop{\stoptime{k}} \process{\thisstate}
	\in \projop{\Wintime{p}{r_k}}^{-1} \winset}{\filtration{\stoptime{k}}}
= \processlaw{1,\exprate}{\processt{\thisstate}{\stoptime{k}}}
	\parens[1]{\projop{\Wintime{p}{r_k}}^{-1} \winset}
\quad\text{on } \eventthat{\stoptime{k}<\infty}.
\]
Since $\processt{\thisstate}{\stoptime{k}}\in \advsetr{\exprate}{\alpha}{\stoptime{k}} \subseteq \advsetr{\exprate}{\alpha}{r_k}$ on the event $\eventthat{\stoptime{k}<\infty}$, on this event we have
\[
\processlaw{1,\exprate}{\processt{\thisstate}{\stoptime{k}}}
	\parens[1]{\projop{\Wintime{p}{r_k}}^{-1} \winset}
\ge \inf_{\genstate\in \advsetr{\exprate}{\alpha}{r_k}}
	\Pr \eventthat[1]{\processt{\genstate}{\Wintime{p}{r_k}}\in \winset}
\ge p
\]
by Lemma~\ref{sp1_wins_soon_lem}. Putting everything together we see
that
\[
\Pr \parens[1]{
	\eventthat[2]{\processt{\thisstate}{\ell_{k+1}}\in \winset}
	\cap \eventthat{\stoptime{k}<\infty}}
\ge \int\displaylimits_{\eventthat{\stoptime{k}<\infty}} p\;d(\Pr)
    =p\cdot \Pr \eventthat{\stoptime{k}<\infty} \ge pq
\]
for all $k\in\N$, contradicting the inequality in \eqref{advantage_sup:winning_time_eqn}.
\end{proof}

Finally, we show that Theorem~\ref{conv_hull_thm} follows from Theorem~\ref{advantage_sup_thm}.

\begin{proof}[Proof of Theorem~\ref{conv_hull_thm}]
\newcommand{\thisstate}{X}
\newcommand{\thisrate}{\lambda}
\newcommand{\rateone}{\lambda_1}
\newcommand{\ratetwo}{\lambda_2}
\newcommand{\thislaw}{\processlaw{\rateone,\ratetwo}{\thisstate}}
\newcommand{\chullevent}{H(\thisstate)}

\newcommand{\occset}{\bothreachedset{\thisstate}{t_n}}
\newcommand{\oneset}{\reachedset{1}{\thisstate}{t_n}}
\newcommand{\twoset}{\reachedset{2}{\thisstate}{t_n}}

\newcommand{\supppoint}{\vect v_n}
\newcommand{\oppspace}{\Hspace_n}
\newcommand{\newoppspace}{\oppspace'}
\newcommand{\thisapex}{\vect z_n}
\newcommand{\thisdir}{\dirvec_n}
\newcommand{\onecone}{\cone{\mu}{1}{\thisapex}{\thisdir}}
\newcommand{\thpcone}{\cone{\mu}{\thp}{\thisapex}{\thisdir}}

\newcommand{\thisproc}{\processt{\thisstate}{t_n}}
\newcommand{\thisadv}{\alpha}

As noted above, we can assume that $\rateone = 1$ and $\ratetwo = \exprate \in (0,1)$, and we set $\expnormpair = \pair{\mu_1}{\mu_2} = \pair{\mu}{\expmean[1] \mu}$, where $\mu$ is species~1's shape function.
Choose any $\thp \in (\exprate,1)$ and let $\thisadv \definedas \thp/\exprate$, so $1<\thisadv<\exprate^{-1}$.
Let $\thisstate\in \statespace$ be a finite initial configuration, and suppose that $\coexstate{1,\exprate}{\thisstate} \cap \chullevent$ occurs. Then there is some sequence of times $t_1,t_2,\dotsc$ with $t_n\to\infty$ such that at each time $t_n$, the total occupied region $\occset$ has a support point $\supppoint\in \oneset$. For each $n$, let $\oppspace$ be an opposing half-space of $\occset$ at $\supppoint$.
Then there is some extremal point $\thisapex\in \cubify{\supppoint}\setminus \cubify{\twoset}$ such that $\newoppspace \definedas (\thisapex-\supppoint) + \oppspace$ is an opposing half-space of $\cubify[1]{\occset}$ at $\thisapex$, so $\thisapex\in \newoppspace$ and $\interior{(\newoppspace)}\cap \cubify{\occset} = \emptyset$. Since any $\mu$-cone of thickness 1 is contained in a half-space by Part~\ref{mu_cone_description:th=1_part} of Lemma~\ref{mu_cone_description_lem}, for each $n$ there is some $\thisdir\in \dirset{\Rd}$ such that $\onecone\subseteq \newoppspace$. Since $\thp<1$, Part~\ref{mu_cone_description:th-_part} of Lemma~\ref{mu_cone_description_lem} implies that
\[
\thpcone \subseteq \cone{\mu}{1-}{\thisapex}{\thisdir}
= \setof{\thisapex} \cup \interior{\parens[1]{\onecone}}
\subseteq \setof{\thisapex} \cup \interior{(\newoppspace)}.
\]
Therefore, since $\parens[2]{\setof{\thisapex} \cup \interior{(\newoppspace)}} \cap \cubify{\twoset} =\emptyset$, we have
\[
\thisapex\in \cubify{\supppoint}\subseteq \cubify{\oneset}
\quad\text{and}\quad
\thpcone \cap \cubify{\twoset} =\emptyset,
\]
so by definition \eqref{state_advantage_def_eqn} we have
\begin{equation}
\label{conv_hull:adv_eqn}
\stateadvantage[1]{1}{\expnormpair}{\thisproc}
\ge \coneadvantage[1]{1}{\expnormpair}{\thpcone}
=\frac{\thp}{\exprate} = \thisadv >1.
\end{equation}
Since $t_n\to\infty$, \eqref{conv_hull:adv_eqn} shows that
\begin{equation}
\label{conv_hull:event_eqn}
\coexstate{1,\exprate}{\thisstate} \cap \chullevent
\subseteq \coexstate{1,\exprate}{\thisstate}\cap
\eventthat[1]{
	\limsup_{t\to\infty} \stateadvantage[1]{1}{\expnormpair}{\processt{\thisstate}{t}}
	\ge \thisadv}.
\end{equation}
Since $\thisadv>1$, Theorem~\ref{advantage_sup_thm} implies that the event on the right-hand side of \eqref{conv_hull:event_eqn} has probability zero, so we also have $\Pr \parens[2]{\coexstate{1,\exprate}{\thisstate} \cap \chullevent} = 0$.
%
%
\end{proof}

\bibliographystyle{halpha}
\bibliography{thesis}

\begin{thebibliography}{{Ahl}11b}

\bibitem[AD11a]{Auffinger:2011aa}
A.~{Auffinger} and M.~{Damron}.
\newblock {A simplified proof of the relation between scaling exponents in
  first-passage percolation}.
\newblock {\em ArXiv e-prints}, September 2011, 1109.0523.

\bibitem[AD11b]{Auffinger:2011ab}
A.~{Auffinger} and M.~{Damron}.
\newblock {Differentiability at the edge of the percolation cone and related
  results in first-passage percolation}.
\newblock {\em ArXiv e-prints}, May 2011, 1105.4172.

\bibitem[ADMP11]{Antunovic:2011aa}
T.~{Antunovi{\'c}}, Y.~{Dekel}, E.~{Mossel}, and Y.~{Peres}.
\newblock {Competing first passage percolation on random regular graphs}.
\newblock {\em ArXiv e-prints}, September 2011, 1109.2575.

\bibitem[{Ahl}11a]{Ahlberg:2011ab}
D.~{Ahlberg}.
\newblock {Asymptotics of first-passage percolation on 1-dimensional graphs}.
\newblock {\em ArXiv e-prints}, July 2011, 1107.2276.

\bibitem[{Ahl}11b]{Ahlberg:2011aa}
D.~{Ahlberg}.
\newblock {The asymptotic shape, large deviations and dynamical stability in
  first-passage percolation on cones}.
\newblock {\em ArXiv e-prints}, July 2011, 1107.2280.

\bibitem[Ahl11c]{Ahlberg:2011ac}
Daniel Ahlberg.
\newblock {\em Asymptotics and dynamics in first-passage and continuum
  percolation}.
\newblock {Ph.D.} dissertation, G{\"o}teborgs universitet. Naturvetenskapliga
  fakulteten, Gothenburg, Sweden, sep 2011, http://hdl.handle.net/2077/26666.

\bibitem[AT07]{Aliprantis:2007aa}
Charalambos~D. Aliprantis and Rabee Tourky.
\newblock {\em Cones and duality}, volume~84 of {\em Graduate Studies in
  Mathematics}.
\newblock American Mathematical Society, Providence, RI, 2007.

\bibitem[BBI01]{Burago:2001aa}
Dmitri Burago, Yuri Burago, and Sergei Ivanov.
\newblock {\em A course in metric geometry}, volume~33 of {\em Graduate Studies
  in Mathematics}.
\newblock American Mathematical Society, Providence, RI, 2001.

\bibitem[{Bla}10]{Blair-Stahn:2010aa}
N.~D. {Blair-Stahn}.
\newblock {First passage percolation and competition models}.
\newblock {\em ArXiv e-prints}, May 2010, 1005.0649.

\bibitem[Boi90]{Boivin:1990aa}
Daniel Boivin.
\newblock First passage percolation: the stationary case.
\newblock {\em Probab. Theory Related Fields}, 86(4):491--499, 1990.

\bibitem[CD81]{Cox:1981aa}
J.~Theodore Cox and Richard Durrett.
\newblock Some limit theorems for percolation processes with necessary and
  sufficient conditions.
\newblock {\em Ann. Probab.}, 9(4):583--603, 1981.

\bibitem[CD09]{Chatterjee:2009aa}
S.~{Chatterjee} and P.~S. {Dey}.
\newblock {Central limit theorem for first-passage percolation time across thin
  cylinders}.
\newblock {\em ArXiv e-prints}, November 2009, 0911.5702.

\bibitem[{Cha}11]{Chatterjee:2011aa}
S.~{Chatterjee}.
\newblock {The universal relation between scaling exponents in first-passage
  percolation}.
\newblock {\em ArXiv e-prints}, May 2011, 1105.4566.

\bibitem[Dei03]{Deijfen:2003aa}
Maria Deijfen.
\newblock Asymptotic shape in a continuum growth model.
\newblock {\em Adv. in Appl. Probab.}, 35(2):303--318, 2003.

\bibitem[DH04]{Deijfen:2004ab}
Maria Deijfen and Olle H{\"a}ggstr{\"o}m.
\newblock Coexistence in a two-type continuum growth model.
\newblock {\em Adv. in Appl. Probab.}, 36(4):973--980, 2004.

\bibitem[DH06a]{Deijfen:2006aa}
Maria Deijfen and Olle H{\"a}ggstr{\"o}m.
\newblock The initial configuration is irrelevant for the possibility of mutual
  unbounded growth in the two-type {R}ichardson model.
\newblock {\em Combin. Probab. Comput.}, 15(3):345--353, 2006.

\bibitem[DH06b]{Deijfen:2006ab}
Maria Deijfen and Olle H{\"a}ggstr{\"o}m.
\newblock Nonmonotonic coexistence regions for the two-type {R}ichardson model
  on graphs.
\newblock {\em Electron. J. Probab.}, 11:no. 13, 331--344 (electronic), 2006.

\bibitem[DH07]{Deijfen:2007aa}
Maria Deijfen and Olle H{\"a}ggstr{\"o}m.
\newblock The two-type {R}ichardson model with unbounded initial
  configurations.
\newblock {\em Ann. Appl. Probab.}, 17(5-6):1639--1656, 2007.

\bibitem[DH08]{Deijfen:2008aa}
Maria Deijfen and Olle H{\"a}ggstr{\"o}m.
\newblock The pleasures and pains of studying the two-type {R}ichardson model.
\newblock In {\em Analysis and stochastics of growth processes and interface
  models}, pages 39--54. Oxford Univ. Press, Oxford, 2008.

\bibitem[DH10]{Damron:2010aa}
M.~{Damron} and M.~{Hochman}.
\newblock {Examples of non-polygonal limit shapes in i.i.d. first-passage
  percolation and infinite coexistence in spatial growth models}.
\newblock {\em ArXiv e-prints}, September 2010, 1009.2523.

\bibitem[DHB04]{Deijfen:2004aa}
M.~Deijfen, O.~H{\"a}ggstr{\"o}m, and J.~Bagley.
\newblock A stochastic model for competing growth on {$\bold R^d$}.
\newblock {\em Markov Process. Related Fields}, 10(2):217--248, 2004.

\bibitem[Die10]{Diestel:2010aa}
Reinhard Diestel.
\newblock {\em Graph theory}, volume 173 of {\em Graduate Texts in
  Mathematics}.
\newblock Springer, Heidelberg, fourth edition, 2010.

\bibitem[DL81]{Durrett:1981aa}
Richard Durrett and Thomas~M. Liggett.
\newblock The shape of the limit set in {R}ichardson's growth model.
\newblock {\em Ann. Probab.}, 9(2):186--193, 1981.

\bibitem[Fol99]{Folland:1999aa}
Gerald~B. Folland.
\newblock {\em Real analysis}.
\newblock Pure and Applied Mathematics (New York). John Wiley \& Sons Inc., New
  York, second edition, 1999.
\newblock Modern techniques and their applications, A Wiley-Interscience
  Publication.

\bibitem[GK84]{Grimmett:1984aa}
Geoffrey Grimmett and Harry Kesten.
\newblock First-passage percolation, network flows and electrical resistances.
\newblock {\em Z. Wahrsch. Verw. Gebiete}, 66(3):335--366, 1984.

\bibitem[GK12]{Grimmett:2012aa}
G.~R. {Grimmett} and H.~{Kesten}.
\newblock {Percolation since Saint-Flour}.
\newblock {\em ArXiv e-prints}, July 2012, 1207.0373.

\bibitem[GM05]{Garet:2005aa}
Olivier Garet and R{\'e}gine Marchand.
\newblock Coexistence in two-type first-passage percolation models.
\newblock {\em Ann. Appl. Probab.}, 15(1A):298--330, 2005.

\bibitem[GM08]{Garet:2008aa}
Olivier Garet and R{\'e}gine Marchand.
\newblock First-passage competition with different speeds: positive density for
  both species is impossible.
\newblock {\em Electron. J. Probab.}, 13:no. 70, 2118--2159, 2008.

\bibitem[Gou07]{Gouere:2007aa}
Jean-Baptiste Gou{{\'e}}r{{\'e}}.
\newblock Shape of territories in some competing growth models.
\newblock {\em Ann. Appl. Probab.}, 17(4):1273--1305, 2007.

\bibitem[Gro99]{Gromov:1999aa}
Misha Gromov.
\newblock {\em Metric structures for {R}iemannian and non-{R}iemannian spaces},
  volume 152 of {\em Progress in Mathematics}.
\newblock Birkh{\"a}user Boston Inc., Boston, MA, 1999.
\newblock Based on the 1981 French original, With appendices by M. Katz, P.
  Pansu and S. Semmes, Translated from the French by Sean Michael Bates.

\bibitem[HM95]{Haggstrom:1995aa}
Olle H{\"a}ggstr{\"o}m and Ronald Meester.
\newblock Asymptotic shapes for stationary first passage percolation.
\newblock {\em Ann. Probab.}, 23(4):1511--1522, 1995.

\bibitem[Hof05]{Hoffman:2005aa}
Christopher Hoffman.
\newblock Coexistence for {R}ichardson type competing spatial growth models.
\newblock {\em Ann. Appl. Probab.}, 15(1B):739--747, 2005.

\bibitem[Hof08]{Hoffman:2008aa}
Christopher Hoffman.
\newblock Geodesics in first passage percolation.
\newblock {\em Ann. Appl. Probab.}, 18(5):1944--1969, 2008.

\bibitem[How04]{Howard:2004aa}
C.~Douglas Howard.
\newblock Models of first-passage percolation.
\newblock In {\em Probability on discrete structures}, volume 110 of {\em
  Encyclopaedia Math. Sci.}, pages 125--173. Springer, Berlin, 2004.

\bibitem[HP98]{Haggstrom:1998aa}
Olle H{\"a}ggstr{\"o}m and Robin Pemantle.
\newblock First passage percolation and a model for competing spatial growth.
\newblock {\em J. Appl. Probab.}, 35(3):683--692, 1998.

\bibitem[HP00]{Haggstrom:2000aa}
Olle H{\"a}ggstr{\"o}m and Robin Pemantle.
\newblock Absence of mutual unbounded growth for almost all parameter values in
  the two-type {R}ichardson model.
\newblock {\em Stochastic Process. Appl.}, 90(2):207--222, 2000.

\bibitem[HW65]{Hammersley:1965aa}
J.~M. Hammersley and D.~J.~A. Welsh.
\newblock First-passage percolation, subadditive processes, stochastic
  networks, and generalized renewal theory.
\newblock In {\em Proc. Internat. Res. Semin., Statist. Lab., Univ. California,
  Berkeley, Calif.}, pages 61--110. Springer-Verlag, New York, 1965.

\bibitem[Kal02]{Kallenberg:2002aa}
Olav Kallenberg.
\newblock {\em Foundations of modern probability}.
\newblock Probability and its Applications (New York). Springer-Verlag, New
  York, second edition, 2002.

\bibitem[Kes86]{Kesten:1986aa}
Harry Kesten.
\newblock Aspects of first passage percolation.
\newblock In {\em \'Ecole d'\'et\'e de probabilit\'es de Saint-Flour,
  XIV---1984}, volume 1180 of {\em Lecture Notes in Math.}, pages 125--264.
  Springer, Berlin, 1986.

\bibitem[Kes87]{Kesten:1987aa}
Harry Kesten.
\newblock Percolation theory and first-passage percolation.
\newblock {\em Ann. Probab.}, 15(4):1231--1271, 1987.

\bibitem[KL05]{Kordzakhia:2005aa}
George Kordzakhia and Steven~P. Lalley.
\newblock A two-species competition model on {$\Bbb Z^d$}.
\newblock {\em Stochastic Process. Appl.}, 115(5):781--796, 2005.

\bibitem[KS02]{Kobayashi:2002aa}
Kei Kobayashi and Kokichi Sugihara.
\newblock Crystal {V}oronoi diagram and its applications.
\newblock {\em Future Generation Computer Systems}, 18(5):681 -- 692, 2002.
\newblock ICCS2001.

\bibitem[{LaG}11]{LaGatta:2011aa}
T.~{LaGatta}.
\newblock {Dissertation: Geodesics of Random Riemannian Metrics}.
\newblock {\em ArXiv e-prints}, July 2011, 1108.0098.

\bibitem[Lee00]{Lee:2000aa}
John~M. Lee.
\newblock {\em Introduction to topological manifolds}, volume 202 of {\em
  Graduate Texts in Mathematics}.
\newblock Springer-Verlag, New York, 2000.

\bibitem[Lig85]{Liggett:1985aa}
Thomas~M. Liggett.
\newblock {\em Interacting particle systems}, volume 276 of {\em Grundlehren
  der Mathematischen Wissenschaften [Fundamental Principles of Mathematical
  Sciences]}.
\newblock Springer-Verlag, New York, 1985.

\bibitem[Lig10]{Liggett:2010aa}
Thomas~M. Liggett.
\newblock {\em Continuous time {M}arkov processes}, volume 113 of {\em Graduate
  Studies in Mathematics}.
\newblock American Mathematical Society, Providence, RI, 2010.
\newblock An introduction.

\bibitem[LW10]{LaGatta:2010aa}
T.~{LaGatta} and J.~Wehr.
\newblock A shape theorem for {R}iemannian first-passage percolation.
\newblock {\em J. Math. Phys.}, 51(5):053502, 14, 2010.

\bibitem[Mos06]{Moszynska:2006aa}
Maria Moszy{{\'n}}ska.
\newblock {\em Selected topics in convex geometry}.
\newblock Birkh{\"a}user Boston Inc., Boston, MA, 2006.
\newblock Translated and revised from the 2001 Polish original.

\bibitem[Mue09]{Mueller:2009aa}
Carl Mueller.
\newblock Some tools and results for parabolic stochastic partial differential
  equations.
\newblock In {\em A minicourse on stochastic partial differential equations},
  volume 1962 of {\em Lecture Notes in Math.}, pages 111--144. Springer,
  Berlin, 2009.

\bibitem[Pim07]{Pimentel:2007aa}
Leandro P.~R. Pimentel.
\newblock Multitype shape theorems for first passage percolation models.
\newblock {\em Adv. in Appl. Probab.}, 39(1):53--76, 2007.

\bibitem[Pim11]{Pimentel:2011aa}
Leandro P.~R. Pimentel.
\newblock Asymptotics for first-passage times on {D}elaunay triangulations.
\newblock {\em Combin. Probab. Comput.}, 20(3):435--453, 2011.

\bibitem[Ric73]{Richardson:1973aa}
Daniel Richardson.
\newblock Random growth in a tessellation.
\newblock {\em Proc. Cambridge Philos. Soc.}, 74:515--528, 1973.

\bibitem[RJ88]{Rice:1888aa}
J.M. Rice and W.W. Johnson.
\newblock {\em An elementary treatise on the differential calculus founded on
  the method of rates or fluxions}.
\newblock J. Wiley and sons, 1888.

\bibitem[Sch92]{Schaudt:1992aa}
Barry~F. Schaudt.
\newblock {Multiplicatively Weighted Crystal Growth Voronoi Diagrams (Thesis)}.
\newblock Technical Report PCS-TR92-177, Dartmouth College, Computer Science,
  Hanover, NH, 1992.

\bibitem[SD91]{Schaudt:1991aa}
Barry~F. Schaudt and R.~L.~Scot Drysdale.
\newblock Multiplicatively weighted crystal growth {V}oronoi diagrams (extended
  abstract).
\newblock In {\em Proceedings of the seventh annual symposium on Computational
  geometry}, SCG '91, pages 214--223, New York, NY, USA, 1991. ACM.

\bibitem[Zha10]{Zhang:2010aa}
Yu~Zhang.
\newblock On the concentration and the convergence rate with a moment condition
  in first passage percolation.
\newblock {\em Stochastic Process. Appl.}, 120(7):1317--1341, 2010.

\end{thebibliography}


\appendix


\chapter{Proofs Omitted from Main Text}
\label{leftover_chap}


\section{Proofs from Chapter~\ref{fpc_basics_chap}}
\label{leftover:fpc_basics_sec}

The following elementary lemma is used in the proofs of Lemma~\ref{final_sets_properties_lem} and Proposition~\ref{entangled_full_prop} in Section~\ref{fpc_basics:two_type_construction_sec}; Lemma~\ref{final_sets_properties_lem} is proved below.

\begin{lem}[Equivalent characterizations of conquering sets]
\label{conq_set_equiv_conditions_lem}
Let $\initconfig$ be an initial configuration in $\Zd$, let $\tmeasurepair$ be a two-type traversal measure on $\edges{\Zd}$, and let $S\subseteq\Zd\setminus \initialset{2}$.
\begin{enumerate}
\item The following are equivalent.
\begin{enumerate}
\item $\rptime{1}{S}{\initialset{1}}{\vect v} < \starptimenew{2}{(\Zd\setminus S)}{\initialset{2}}{\vect v}$ for all $\vect v\in \bd S\setminus \initialset{1}$.
\item $\rptime{1}{S}{\initialset{1}}{\vect v} < \starptimenew{2}{(\Zd\setminus S)}{\initialset{2}}{\vect v}$ for all $\vect v\in S\setminus \initialset{1}$.
\end{enumerate}

\item The following are equivalent.
\begin{enumerate}
\item For all $\vect v\in \bd S$, $\rptime{1}{S}{\initialset{1}}{\vect v} <\infty$ and $\rptime{1}{S}{\initialset{1}}{\vect v} \le \starptimenew{2}{(\Zd\setminus S)}{\initialset{2}}{\vect v}$.
\item For all $\vect v\in S$, $\rptime{1}{S}{\initialset{1}}{\vect v} <\infty$ and $\rptime{1}{S}{\initialset{1}}{\vect v} \le \starptimenew{2}{(\Zd\setminus S)}{\initialset{2}}{\vect v}$.
\end{enumerate}

\item The equivalent statements in part 1 imply those in part 2. If $\tmeasurepair$ satisfies \notiesdet, and its component measures each satisfy \stargeodesicsexist, then the equivalent statements in part 2 imply those in part 1.

\end{enumerate}
By symmetry, the corresponding statements hold with the labels 1 and 2 switched.
\end{lem}

\begin{proof}

\proofpart{1}

(1a $\Rightarrow$ 1b) If (1a) holds, then the inequality in (1b) is trivial for all $\vect v\in \bd S$. On the other hand, for $\vect v$ in the interior of $S$, it is trivial that $ \starptimenew{2}{(\Zd\setminus S)}{\initialset{2}}{\vect v} = \infty$ (since $\vect v\not\in \nbrhood[2]{\Zd\setminus S}$ and $\vect v\not\in \initialset{2}$), and (1a) implies that $\rptime{1}{S}{\initialset{1}}{\vect v}<\infty$ (since the property ``$\rptime{1}{S}{\initialset{1}}{\vect v}<\infty$" is constant on connected components of $S$, and any connected component of $S$ has a nonempty intersection with $\bd S$). Thus $\rptime{1}{S}{\initialset{1}}{\vect v} < \starptimenew{2}{(\Zd\setminus S)}{\initialset{2}}{\vect v}$ in either case, so (1b) holds.


(1b $\Rightarrow$ 1a) Trivial.

\proofpart{2}

(2a $\Rightarrow$ 2b) Same argument as (1a) $\Rightarrow$ (1b).

(2b $\Rightarrow$ 2a) Trivial.

\proofpart{3} 

First we show that (1b) $\Rightarrow$ (2b). If $\vect v\in \initialset{1}$, then $\rptime{1}{S}{\initialset{1}}{\vect v} \le 0\le \starptimenew{2}{(\Zd\setminus S)}{\initialset{2}}{\vect v}$. Thus, if (1b) holds, then we get $\rptime{1}{S}{\initialset{1}}{\vect v} \le  \starptimenew{2}{(\Zd\setminus S)}{\initialset{2}}{\vect v}$ for all $\vect v\in S$, and moreover, $\rptime{1}{S}{\initialset{1}}{\vect v} <\infty$ for all $\vect v\in S$, so (2b) holds.

Finally, suppose that $\tmeasurepair$ satisfies \notiesdet, and its component measures each satisfy \stargeodesicsexist. We will show that (2a) $\Rightarrow$ (1a).
Suppose (2a) holds, so for all $\vect v\in \bd S$ we have $\rptime{1}{S}{\initialset{1}}{\vect v} <\infty$ and $\rptime{1}{S}{\initialset{1}}{\vect v} \le \starptimenew{2}{(\Zd\setminus S)}{\initialset{2}}{\vect v}$. Let $\vect v\in \bd S\setminus \initialset{1}$. If $ \starptimenew{2}{(\Zd\setminus S)}{\initialset{2}}{\vect v}=\infty$, we're done. Otherwise, let $\gamma_1$ be a $\ptime_1^S$-geodesic from $\initialset{1}$ to $\vect v$, and let $\gamma_2$ be a $ \starptmetric[2]{(\Zd\setminus S)}$-geodesic from $\initialset{2}$ to $\vect v$ (both of which exist by the hypothesis \stargeodesicsexist\ and the assumption that the respective passage times from $\initialset{i}$ to $\vect v$ are finite). Then we have
\[
\tau_1(\gamma_1) = \ptime_1^{\gamma_1}(\initialset{1},\vect v)
=\rptime{1}{S}{\initialset{1}}{\vect v} \le \starptimenew{2}{(\Zd\setminus S)}{\initialset{2}}{\vect v}
=\ptime_2^{\gamma_2}(\initialset{2},\vect v) = \tau_2(\gamma_2).
\]
The paths $\gamma_1$ and $\gamma_2$ are edge-disjoint by definition, and since $\vect v$ is not an element of either $\initialset{1}$ or $\initialset{2}$, both $\gamma_1$ and $\gamma_2$ must contain at least one edge. Therefore, since $\tmeasurepair$ satisfies \notiesdet, the left-hand and right-hand sides of the above inequality cannot be equal, so the inequality must be strict. Thus, $\rptime{1}{S}{\initialset{1}}{\vect v} < \starptimenew{2}{(\Zd\setminus S)}{\initialset{2}}{\vect v}$ for all $\vect v\in \bd S\setminus \initialset{1}$, so (1a) holds.
\end{proof}

\begin{proof}[Proof of Lemma~\ref{final_sets_properties_lem} (Properties of the finally conquered sets)]
The proof for each part is numbered accordingly.
\begin{enumerate}
\item (trivial relationships with initial sets) These properties are either explicit in the definition of $\finalset{i}$ or are trivial.
 
\item (inequality for vertices in final sets) We will prove the statement with $i=1$ for concreteness; the case $i=2$ follows by symmetry.

($\implies$) Let $\vect v\in \finalset{1}$. If $\vect v\in \initialset{1}$, we're done, so assume $\vect v\in \finalset{1} \setminus \initialset{1}$. By the definition of $\finalset{1}$, $\vect v$ is contained in some $\finalset{1}$-set $S$; that is, $S\subseteq \finalset{1}\subseteq \Zd\setminus \initialset{2}$, and using Lemma~\ref{conq_set_equiv_conditions_lem},
\[
\rptime{1}{S}{\initialset{1}}{\vect v'} < \starptimenew{2}{(\Zd\setminus S)}{\initialset{2}}{\vect v'}
	\text{ for all } \vect v'\in S\setminus \initialset{1}.
\]
Moreover, since $S\subseteq \finalset{1}$, by the monotonicity of $\ptimefamily$ we have
\[
\rptime{1}{\finalset{1}}{\initialset{1}}{\vect v'}\le \rptime{1}{S}{\initialset{1}}{\vect v'}
< \starptimenew{2}{(\Zd\setminus S)}{\initialset{2}}{\vect v'}\le \starptimenew{2}{(\Zd\setminus \finalset{1})}{\initialset{2}}{\vect v'}
	\quad\forall \vect v'\in S\setminus \initialset{1}.
\]
Setting $\vect v' = \vect v\in S\setminus \initialset{1}$, we get $\rptime{1}{\finalset{1}}{\initialset{1}}{\vect v} < \starptimenew{2}{(\Zd\setminus \finalset{1})}{\initialset{2}}{\vect v}$.
 
($\followsfrom$) Let $S$ be the set of all vertices $\vect v\not\in \initialset{2}$ satisfying $\vect v\in \initialset{1}$ or $\rptime{1}{\finalset{1}}{\initialset{1}}{\vect v} < \starptimenew{2}{(\Zd\setminus \finalset{1})}{\initialset{2}}{\vect v}$.
By the above argument for the converse, we have $\finalset{1}\subseteq S$.
Thus we have (by the monotonicity of $\ptime^{\cdot}$ and the definition of $S$)
\[
\rptime{1}{S}{\initialset{1}}{\vect v}\le \rptime{1}{\finalset{1}}{\initialset{1}}{\vect v} 
< \starptimenew{2}{(\Zd\setminus \finalset{1})}{\initialset{2}}{\vect v}\le \rptime{1}{(\Zd\setminus S)}{\initialset{2}}{\vect v)}
\quad\forall \vect v\in S\setminus \initialset{1}.
\]
Therefore, $S$ is a $\finalset{1}$-set, so $S\subseteq \finalset{1}$.

\item (connected components) Let $V$ be a component of $\finalset{i}$, and let $\vect v\in V$. Since $V\subseteq \finalset{i}$, Part~\ref{final_sets:characterization_part} implies that $\rptime{i}{\finalset{i}}{\initialset{i}}{\vect v} < \infty$, so there exists some $\finalset{i}$-path $\gamma$ from $\initialset{i}$ to $\vect v$. Since $\gamma$ is connected and contains $\vect v\in V$, we must have $\gamma\subseteq V$ since $V$ is a component. Thus, if $\vect u$ is the first vertex of $\gamma$, we have $\vect u\in V$. Since $\vect u\in \initialset{i}$, this implies that $V$ contains the component of $\initialset{i}$ that contains $\vect u$.
%
%
%
%
%

\item (final sets are union) We have $\reachedset{i}{\initconfig}{t}\subseteq \finalset{i}$ for all $t\ge 0$ by definition. On the other hand, Part~\ref{final_sets:characterization_part} implies that for any $\vect v\in \finalset{i}$ we have $\rptime{i}{\finalset{i}}{\initialset{i}}{\vect v}<\infty$, so there is some $t<\infty$ such that  $\rptime{i}{\finalset{i}}{\initialset{i}}{\vect v}\le t$, and hence $\vect v\in \rreachedset{i}{\initialset{i}}{\finalset{i}}{t} = \reachedset{i}{\initconfig}{t}$.
%
%
\end{enumerate}
\end{proof}

\section{Proofs from Chapter~\ref{deterministic_chap}}
\label{leftover:deterministic_sec}

The following lemma is needed to prove Lemma~\ref{wedge_properties_lem}; Lemma~\ref{wedge_properties_lem} is proved below, following the proof of Lemma~\ref{affine_si_span_lem}.

\begin{lem}[Affine spans of apexes]
\label{affine_si_span_lem}
\newcommand{\apexspan}{A}
\newcommand{\intpoint}{\vect z}

Let $\asiset\subseteq\Rd$ be affine scale-invariant, let $\vect a_1,\dotsc, \vect a_n$ be apexes of $\asiset$, and let $\apexspan = \aff \setof{\vect a_1,\dotsc, \vect a_n}$. If $\intpoint\in \asiset$, then $\asiset$ contains the set
\[
\Hspace = \setbuilderbar[1]{c_0 \intpoint + \sum_{i=1}^n c_i \vect a_i}
	{c_0 >0,\: c_1,\dotsc,c_n\in\R,\: \sum_{i=0}^n c_i = 1}.
\]
Moreover, if $\intpoint\in \apexspan$, then $\Hspace = \apexspan$, and if $\intpoint\not\in \apexspan$ and $n>0$, then $\Hspace$ is an open affine half-space of dimension $(\dim \apexspan)+1$ with $\affrelbd \Hspace = \apexspan$ and $\intpoint\in J$.
\end{lem}

\begin{proof}
\newcommand{\apexspan}{A}
\newcommand{\intpoint}{\vect z}
\newcommand{\spanpt}{\vect y}
\newcommand{\smallspanpt}{\vect x}

\newcommand{\icoef}{\beta}
\newcommand{\scoef}{\alpha}

\newcommand{\projpt}{\intpoint_{\apexspan}}

We prove the first statement by induction on $n$. If $n=0$, then $\Hspace=\setof{\intpoint}$, so the statement is a tautology. Now let $n\ge 1$, and suppose the statement holds for any collection of $n-1$ apexes. Let $\spanpt\in \Hspace$, so $\spanpt = c_0 \intpoint + c_1\anapex_1 +\dotsb + c_n \anapex_n$ for some $c_0,c_1,\dotsc,c_n\in \R$ with $c_0>0$ and $\sum_{i=0}^n c_i =1$. Since $0<c_0 = 1- \sum_{i=1}^n c_i$, there must exist some $i\in \setof{1,\dotsc,n}$ such that $c_i<\frac{1}{n}\le 1$. By relabeling the apexes if necessary, assume without loss of generality that $c_n<1$. Now define
\[
\smallspanpt \definedas \frac{c_0}{1-c_n} \intpoint
+ \frac{c_1}{1-c_n} \anapex_1 + \dotsb + \frac{c_{n-1}}{1-c_n} \anapex_{n-1}.
\]
Then
\[
\frac{c_0}{1-c_n}>0, 
\quad\text{and}\quad
\sum_{i=0}^{n-1} \frac{c_i}{1-c_n} = \frac{1}{1-c_n} \sum_{i=0}^{n-1} c_i
=\frac{1-c_n}{1-c_n} = 1,
\]
so $\smallspanpt \in\asiset$ by the inductive hypothesis. Now note that since $1-c_n>0$, we have
\begin{align*}
\spanpt = (1-c_n) \smallspanpt + c_n \anapex_n 
&\in \setbuilder{r \smallspanpt + (1-r)\anapex_n}{r>0}
= \setbuilder{\anapex_n + r(\smallspanpt - \anapex_n)}{r>0}
=\ray{\anapex_n}{\smallspanpt},
\end{align*}
and the ray $\ray{\anapex_n}{\smallspanpt}$ is contained in $\asiset$ because $\smallspanpt\in\asiset$ and $\asiset$ is affine scale-invariant at $\anapex_n$. Therefore $\spanpt\in \asiset$, which proves the first statement.

If $n=0$, the second statement is vacuous, so let $n\ge 1$. Suppose first that $\intpoint\in \apexspan$, so $\intpoint = \sum_{i=1}^n \icoef_i \anapex_i$ for some $\icoef_1,\dotsc,\icoef_n\in \R$ with $\sum_{i=1}^n \icoef_i = 1$. If $\spanpt\in J$, then there exist $c_0,c_1,\dotsc,c_n\in \R$ with $c_0>0$ and $\sum_{i=0}^n c_i =1$ such that
\[
\spanpt = c_0\intpoint + \sum_{i=1}^n c_i\anapex_i
= \sum_{i=1}^n (c_0\icoef_i+c_i)\anapex_i.
\]
Since $\sum_{i=1}^n \icoef_i = 1$, we have
$
\sum_{i=1}^n (c_0\icoef_i+c_i)
=c_0+ \sum_{i=1}^n c_i
=1,
$
so $\spanpt\in \aff \setof{\anapex_1,\dotsc,\anapex_n} = \apexspan$. Conversely, if $\spanpt\in \apexspan$, then $\spanpt = \sum_{i=1}^n \scoef_i \anapex_i$ for some $\scoef_1,\dotsc,\scoef_n\in\R$ with $\sum_{i=1}^n \scoef_i = 1$. Let $c_0=1$ and $c_i = \scoef_i-\icoef_i$ for $1\le i\le n$. Then $c_0>0$ and $c_0+\sum_{i=1}^n c_i = 1+\sum_{i=1}^n\scoef_i - \sum_{i=1}^n \icoef_i = 1$, and
\[
c_0\intpoint + \sum_{i=1}^n c_i\anapex_i
= 1\cdot \sum_{i=1}^n \icoef_i\anapex_i + \sum_{i=1}^n (\scoef_i-\icoef_i)\anapex_i
=  \sum_{i=1}^n \scoef_i\anapex_i
=\spanpt,
\]
so $\spanpt\in \Hspace$. Thus $\Hspace = \apexspan$.

Now suppose $\intpoint\notin \apexspan$.
Let $\projpt$ be the orthogonal projection of $\intpoint$ onto $\apexspan$ with respect to $\innerprod{\cdot}{\cdot}$. Then the vector $\dirsymb \definedas \intpoint -\projpt$ is orthogonal to $\apexspan$, and we have $\innerprod{\anapex_i-\projpt}{\dirsymb} = 0$ for $1\le i\le n$ since $\anapex_i-\projpt$ is parallel to $\apexspan$ for each $i$. Moreover, $\dirsymb\ne\vect 0$ since $\intpoint\notin \apexspan$. Now let $\spanpt$ be any point in $\aff \setof{\intpoint,\anapex_1,\dotsc,\anapex_n}$, so $\spanpt = c_0\intpoint+\sum_{i=1}^n c_i \anapex_i$ for some $c_0,\dotsc,c_n\in \R$ with $\sum_{i=0}^n c_i =1$. Then
\begin{align}
\label{affine_si_span:innerprod_eqn}
\innerprod{\spanpt - \projpt}{\dirsymb}
&= \innerprod[2]{\textstyle c_0\intpoint+\sum_{i=1}^n c_i \anapex_i - \projpt}{\dirsymb} \notag\\
&= \innerprod[2]{\textstyle c_0\intpoint+\sum_{i=1}^n c_i \anapex_i - \projpt
	+\parens[2]{\sum_{i=1}^n c_i \projpt - \sum_{i=1}^n c_i \projpt}}
	{\dirsymb} \notag\\
&= \innerprod[2]{\textstyle c_0\intpoint - \parens[1]{1-\sum_{i=1}^n c_i} \projpt}{\dirsymb}
	+\sum_{i=1}^n c_i \innerprod{ \anapex_i -\projpt }{\dirsymb} \notag\\
&= c_0 \innerprod{\intpoint - \projpt}{\dirsymb} + 0 \notag\\
&= c_0 \norm{\dirsymb}_{\ltwo{d}}^2.
\end{align}
The equation \eqref{affine_si_span:innerprod_eqn} shows that $\innerprod{\spanpt - \projpt}{\dirsymb} >0$ if and only if $c_0>0$, and therefore
\[
\Hspace =
\setbuilder[2]{\spanpt \in \aff \setof{\intpoint,\anapex_1,\dotsc,\anapex_n}}
	{\innerprod{\spanpt - \projpt}{\dirsymb} >0}.
\]
Since $\dirsymb$ is a nonzero vector parallel to $\affspan{\intpoint,\apexspan}$, the set on the right is an open half-space of dimension 
$\dim \affspan{\intpoint,\apexspan}$,
which equals $1+ (\dim \apexspan)$ since $\intpoint\notin \apexspan$. To see that $\affrelbd \Hspace = \apexspan$, first note that
\[
\affrelbd \Hspace = \setbuilder[2]{\spanpt \in \affspan{\intpoint,\apexspan} }
	{\innerprod{\spanpt - \projpt}{\dirsymb} =0}.
\]
Thus, $\affrelbd \Hspace \supseteq \apexspan$ since $\innerprod{\spanpt - \projpt}{\dirsymb} =0$ for any $\spanpt\in \apexspan$ by the definition of $\dirsymb$. On the other hand, \eqref{affine_si_span:innerprod_eqn} shows that if $\spanpt \in \affspan{\intpoint,\apexspan}$ and $\innerprod{\spanpt - \projpt}{\dirsymb} =0$, then $\spanpt\in \apexspan$, so we also have $\apexspan \subseteq \affrelbd \Hspace$.
\end{proof}

\begin{proof}[Proof of Lemma~\ref{wedge_properties_lem}]
\newcommand{\mysiset}{\siversion{\asiset}}
\newcommand{\myapexset}{\apexset{\asiset}}

\newcommand{\apexspan}{A}

\newcommand{\thisz}{\vect z}
\newcommand{\thisy}{\vect y}
\newcommand{\thisx}{\vect x}

\newcommand{\stab}[3][0]{\operatorname{Stab}_{#2} \parens[#1]{#3}}

Let $\asiset$ be a wedge in $\Rd$.
\begin{enumerate}
\item 
The scale-invariant version $\mysiset$ of $\asiset$ is unique.

\begin{proof}
Suppose $\asiset = \anapex + \mysiset$ and $\asiset = \anapex' + \mysiset'$ for some $\anapex,\anapex'\in\Rd$ and some scale-invariant sets $\mysiset,\mysiset'\subseteq\Rd$. We will show that $\mysiset\subseteq\mysiset'$, whence the reverse inclusion follows by symmetry. Let $\thisy \in \mysiset$, so $\thisy = \thisz -\anapex$ for some $\thisz\in \asiset$. Let $\thisz' \definedas \thisz -\anapex+\anapex'$.
Since $\mysiset$ and $\mysiset'$ are both scale-invariant by assumption, the points $\anapex$ and $\anapex'$ are both apexes of $\asiset$ by definition. Thus, since $\thisz\in \asiset$, Lemma~\ref{affine_si_span_lem} implies that $\asiset$ contains the point $\thisz' = 1\cdot \thisz -1\cdot \anapex + 1\cdot \anapex'$. Therefore, $\thisz' - \anapex' \in \asiset - \anapex' = \mysiset'$. But $\thisz' - \anapex' = \thisz-\anapex = \thisy$, so $\thisy\in \mysiset'$. Since $\thisy$ was an arbitrary point in $\mysiset$, we have $\mysiset \subseteq \mysiset'$. By symmetry, we also have $\mysiset' \subseteq \mysiset$ and hence $\mysiset = \mysiset'$.
\end{proof}

\item 
$\myapexset$ is an affine subspace of $\Rd$ which can be identified with the vector space of translations that fix $\mysiset$. 

\begin{proof}
Let $\mysiset$ be the scale-invariant version of $\asiset$, and let 
\[
\stab{\Rd}{\mysiset} \definedas \setbuilder[2]{\anapex\in\Rd}{\anapex+\mysiset = \mysiset}.
\]
That is, $\stab{\Rd}{\mysiset}$ is the stabilizer of $\mysiset$ in the additive group $\Rd$, acting on $\powerset{\Rd}$ via translation. Since $\mysiset$ is scale-invariant, we have
\[
\parens[2]{r>0 \text{ and } \anapex+\mysiset = \mysiset} \implies
r\anapex + \mysiset = r\anapex + r\mysiset
=r(\anapex + \mysiset) = r\mysiset = \mysiset.
\]
That is, if $\anapex \in \stab{\Rd}{\mysiset}$, then $r\anapex\in \stab{\Rd}{\mysiset}$ for all $r>0$. Since $\stab{\Rd}{\mysiset}$ is a subgroup of $\Rd$ under addition, we also have ${-\anapex}\in \stab{\Rd}{\mysiset}$ and ${\vect 0}\in \stab{\Rd}{\mysiset}$, so it follows that
\[
{\anapex}\in \stab{\Rd}{\mysiset} \implies {r\anapex}\in \stab{\Rd}{\mysiset}
\text{ for all $r\in\R$}.
\]
Therefore, $\stab{\Rd}{\mysiset}$ is in fact a vector subspace of $\R^d$. We claim that $\stab{\Rd}{\mysiset} = \apexset{\mysiset}$. First suppose $\anapex\in \stab{\Rd}{\mysiset}$. Then for any $r>0$,
\[
\homothe{r}{\anapex} \mysiset
= \anapex+ r(\mysiset -\anapex)
= (1-r)\anapex + r\mysiset
= (1-r)\anapex + \mysiset
=\mysiset,
\]
since $\mysiset$ is scale-invariant and $(1-r)\anapex \in \stab{\Rd}{\mysiset}$, so we have $\anapex\in \apexset{\mysiset}$. On the other hand, if $\anapex\in \apexset{\mysiset}$, then we have $\mysiset - \anapex = \mysiset$ by Part~\ref{wedge_properties:siversion_unique_part}, since $\mysiset - \anapex$ is scale-invariant by the definition of an apex, and $\mysiset$ is its own unique scale-invariant version. Therefore, $-\anapex\in \stab{\Rd}{\mysiset}$, and hence $\anapex\in \stab{\Rd}{\mysiset}$ since $\stab{\Rd}{\mysiset}$ is closed under additive inverses. Thus we have $\stab{\Rd}{\mysiset} = \apexset{\mysiset}$ as claimed. Finally, note that if $\anapex'$ is any apex of $\asiset$, we have $\myapexset = \anapex' + \apexset{\mysiset}$, so $\myapexset$ is an affine subspace of $\Rd$ parallel to the linear subspace $\stab{\Rd}{\mysiset}$.
\end{proof}

\item 
If $\asiset$ is nontrivial (i.e.\ $\asiset\notin \setof[2]{\emptyset,\Rd}$),
then $\myapexset\subseteq \bd \asiset$. 

\begin{proof}
If $\asiset\ne\emptyset$, then Lemma~\ref{affine_si_span_lem} implies that $\myapexset \subseteq \closure{\asiset}$. Replacing $\asiset$ with its complement, it follows that if $\Rd\setminus \asiset \ne\emptyset$, then $\myapexset \subseteq \shellof[0]{\asiset}$, since $\apexset[2]{\Rd\setminus \asiset} = \myapexset$. Thus, if both $\asiset$ and $\Rd\setminus \asiset$ are nonempty, we have $\myapexset \subseteq \closure{\asiset} \cap \shellof[0]{\asiset} = \bd \asiset$.
\end{proof}

\item 
The sets $\interior{\asiset}$, $\closure{\asiset}$, $\bd \asiset$, $\exterior{\asiset}$, and $\shellof{\asiset}$ and  are all wedges whose apex sets contain $\apexset{\asiset}$. 

\begin{proof}
Let $\anapex$ be an apex of $\asiset$. For any $r>0$, the homothety $\homothe{r}{\anapex}\colon \Rd\to\Rd$ is a homeomorphism. Thus, since $\homothe{r}{\anapex} \asiset = \asiset$ for all $r>0$ by assumption, we also have $\homothe{r}{\anapex} \interior{\asiset} = \interior{\asiset}$ and  $\homothe{r}{\anapex} \closure{\asiset} = \closure{\asiset}$ for all $r>0$, so $\interior{\asiset}$ and $\closure{\asiset}$ are scale-invariant at $\anapex$. Replacing $\asiset$ with its complement, it then follows that $\exterior{\asiset}$ and $\shellof{\asiset}$ are scale-invariant at $\anapex$, and then so is $\bd \asiset = \closure{\asiset} \cap \shellof{\asiset}$.
\end{proof}

%
%
%

\item 
If $\asiset$ is pointed, then it contains all of its apexes, and $\Rd\setminus \asiset$ is blunt.

\begin{proof}
Suppose $\asiset$ is pointed, so it contains one of its apexes $\anapex$. Let $\anapex'$ be an arbitrary apex of $\asiset$, and let $\apexspan \definedas \aff \setof{\anapex,\anapex'}$. By Lemma~\ref{affine_si_span_lem}, $\asiset$ contains the set $\apexspan$ since $\anapex\in \asiset \cap \apexspan$, so we have $\anapex'\in\asiset$. Thus $\asiset$ contains all of its apexes since $\anapex'$ was arbitrary. Since $\apexset{\Rd\setminus \asiset} = \apexset{\asiset}$, this implies that the wedge $\Rd\setminus \asiset$ contains none of its apexes, hence is blunt.
\end{proof}

\item 
If $\asiset$ is pointed and contained in some half-space, then every apex of $\asiset$ is a support point of $\asiset$. 

\begin{proof}
\newcommand{\suppspace}{\Hspace_{\anapex}}
\newcommand{\bdpt}{\vect z}
\newcommand{\thishom}{\homothe{r}{\anapex}}

Suppose $\asiset$ is pointed and contained in the half-space $\Hspace$. Let $\anapex\in \myapexset$, and let $\bdpt\in \bd \Hspace$.  We claim that the half-space $\suppspace \definedas \Hspace + \anapex-\bdpt$ is a support half-space of $\asiset$ at $\anapex$. Let $\dirsymb$ be the outer normal vector of $\Hspace$. Since $\asiset$ is pointed and contained in $\Hspace$, we have $\anapex\in \asiset\subseteq \Hspace$, which implies that $\innerprod{\anapex-\bdpt}{\dirsymb}\le 0$. Now suppose for a contradiction that there is some $\thisx\in \asiset\setminus \suppspace$. Then since $\dirsymb$ is also the outer normal to $\suppspace$, and $\thisx\notin \suppspace$, we have $\innerprod{\thisx-\anapex}{\dirsymb}>0$. Since $\thisx\in \asiset$, and $\asiset$ is scale-invariant at $\anapex$, we have $\thishom \thisx \in \asiset$ for all $r>0$. Since $\asiset\subseteq\Hspace$, this implies that
\begin{equation}
\label{wedge_properties:innerprod_neg_eqn}
\innerprod{\thishom\thisx-\bdpt}{\dirsymb} \le 0
\quad\text{for all $r>0$.}
\end{equation}
On the other hand,
\begin{align*}
\innerprod{\thishom\thisx-\bdpt}{\dirsymb}
=\innerprod[2]{\anapex + r(\thisx-\anapex)-\bdpt}{\dirsymb}
=\innerprod{\anapex-\bdpt}{\dirsymb} + r\innerprod{\thisx-\anapex}{\dirsymb},
\end{align*}
so
\[
\innerprod{\thishom\thisx-\bdpt}{\dirsymb}>0
\quad\text{for all}\quad
r> -\frac{\innerprod{\anapex-\bdpt}{\dirsymb}}{\innerprod{\thisx-\anapex}{\dirsymb}}.
\]
This contradicts \eqref{wedge_properties:innerprod_neg_eqn}, so we conclude that such an $\thisx$ cannot exist. Therefore, $\asiset\subseteq\suppspace$, and $\anapex\in \asiset\cap \bd\suppspace$ by construction, so $\suppspace$ is a support half-space of $\asiset$ at $\anapex$.
\end{proof}

\item 
$\asiset \subseteq \asiset + \mysiset$, with equality if and only if $\asiset$ is convex.

\begin{proof}
The inclusion $\subseteq$ is trivial for pointed wedges and almost trivial for blunt wedges. For the reverse inclusion when $\asiset$ is convex, see Lemma~\ref{convex_si_sets_lem}.
\end{proof}
\end{enumerate}
This concludes the proof of Lemma~\ref{wedge_properties_lem}.
\end{proof}

\begin{lem}[Convex scale-invariant sets]
\label{convex_si_sets_lem}
A scale-invariant set $A$ is convex if and only if $A+A\subseteq A$.
\end{lem}

\begin{proof}
We need to show that if $A$ is scale-invariant, then $A$ is convex if and only if $\vect x+\vect y\in A$ for all $\vect x,\vect y\in A$. If $A$ is convex and $\vect x,\vect y\in A$, then $\frac{\vect x+\vect y}{2} = \frac{1}{2}\vect x + (1-\frac{1}{2})\vect y\in A$. Thus since $A$ is scale-invariant, $\vect x+\vect y = 2\cdot \frac{\vect x+\vect y}{2} \in A$. On the other hand, suppose that $\vect x+\vect y\in A$ for all $\vect x,\vect y\in A$. Choose any $\vect x,\vect y\in A$ and $\alpha \in [0,1]$. Then $\vect x' = \alpha \vect x\in A$ and $\vect y' = (1-\alpha)\vect y \in A$ since $A$ is scale-invariant, so $\alpha \vect x + (1-\alpha)\vect y = \vect x' + \vect y'\in A$ by assumption. Thus $A$ is convex.
\end{proof}

The following two lemmas are the main results needed to prove Lemma~\ref{unobtrusive_star_set_lem}, which is the basis for Proposition~\ref{narrow_extinction_prop} and Part~\ref{conquering_from_apex:narr_part} of Proposition~\ref{conquering_from_apex_prop}. The proof uses several lemmas from Appendix~\ref{notmyjob_chap}.

\begin{lem}[Existence of blocking segments in small cones]
\label{blocking_segments_lem}
\newcommand{\bluntcone}{C\setminus \setof{\vect 0}}
\newcommand{\conept}{\vect x} 
\newcommand{\bdpt}{\vect y_{\conept}} 
\newcommand{\segpt}{\vect y} 

Let $C$ be a closed cone at $\vect 0$ which is contained in some half-space. For every $\conept\in \bluntcone$ and $\alpha >\coneadvantage{1}{\normpair}{C}$, there exists $\bdpt\in \bd C \setminus \setof{\vect 0}$ such that $\distnew{\mu_2}{\bdpt}{\segpt} < \alpha \distnew{\mu_1}{\vect 0}{\segpt}$ for all $\segpt\in \segment{\bdpt}{\conept}$.
\end{lem}

\begin{proof}
\newcommand{\bluntcone}{C\setminus \setof{\vect 0}}
\newcommand{\conept}{\vect x} 
\newcommand{\bdpt}{\vect y_{\conept}} 
\newcommand{\segpt}{\vect y} 


\newcommand{\abdpt}{\vect y'} 
\newcommand{\ibdpt}{\vect y^{\conept}_0} 
\newcommand{\bdpts}[1]{\vect y^{\conept}_{#1}} 
\newcommand{\hpt}{\widetilde{\segpt}} 

\newcommand{\ffcn}{f_{\conept}} 
\newcommand{\gfcn}{g_{\conept}} 

\newcommand{\nbhd}{U_{(\conept)}} 
\newcommand{\delnbhd}{\nbhd^*} 

\newcommand{\proofcase}[2]{\vskip 2mm \noindent \fbox{\underline{Case #1:}\ #2} \vskip 2mm}


First, let $H$ be a support half-space of $C$ at $\vect 0$, which exists by Part~\ref{wedge_properties:apex_support_part} of Lemma~\ref{wedge_properties_lem} because $C$ is contained in some half-space by assumption. For each $\conept\in \bluntcone$, let $\ibdpt$ be a $\mu_2$-closest point in $\bd C$ to $\conept$. That is, we choose some $\ibdpt\in \bd C$ satisfying
\[
\distnew{\mu_2}{\ibdpt}{\conept} = \distnew{\mu_2}{\bd C}{\conept},
\]
which is possible because $\bd C$ is closed and $\setof{\conept}$ is compact (cf.\ \cite[p.~15]{Burago:2001aa}). We divide the proof into three cases depending on the location of $\ibdpt$ in relation to $\bd H$.

\proofcase{1}{$\ibdpt\ne \vect 0$.}

In this case, we simply take $\bdpt = \ibdpt$.
Then we claim that $\segment{\bdpt}{\conept}\subseteq \bluntcone$. To see this, choose $\segpt_0 \in \segment{\bdpt}{\conept} \setminus \interior{C}$ with $\distnew{\mu_2}{\conept}{\segpt_0} = \distnew[2]{\mu_2}{\conept}{\segment{\bdpt}{\conept} \setminus \interior{C}}$, which is possible because $\segment{\bdpt}{\conept} \setminus \interior{C}$ is compact and nonempty (it contains $\bdpt$). It follows from this choice of $\segpt_0$ that $\ocsegment{\segpt_0}{\vect x}\subseteq \interior{C}$ and $\segpt_0\in \bd C$ (see Lemma~\ref{closest_bd_point_lem}). But then we must have $\segpt_0 = \bdpt$, because otherwise $\segpt_0$ would be a point in $\bd C$ which is strictly closer to $\conept$ than $\bdpt$ in the $\mu_2$-metric, contradicting our choice of $\ibdpt$. Thus we have $\segment{\bdpt}{\conept}\subseteq \bluntcone$ since $\bdpt\in \bluntcone$ and $\ocsegment{\segpt_0}{\vect x}\subseteq \interior{C}$.

Moreover, we claim that for all $\segpt\in \segment{\bdpt}{\conept}$, we have $\distnew{\mu_2}{\bd C}{\segpt} = \distnew{\mu_2}{\bdpt}{\segpt}$, because if there were some point $\abdpt\in \bd C$ with $\distnew{\mu_2}{\abdpt}{\segpt} < \distnew{\mu_2}{\bdpt}{\segpt}$, then since $\segpt\in \segment{\bdpt}{\conept}$, we would also have $\distnew{\mu_2}{\abdpt}{\conept} < \distnew{\mu_2}{\bdpt}{\conept}$ by the triangle inequality, contradicting our choice of $\bdpt = \ibdpt$ as a $\mu_2$-closest point in $\bd C$ to $\conept$.
Therefore, since $\alpha>\coneadvantage{1}{\normpair}{C}$ and $\segment{\bdpt}{\conept}\subseteq \bluntcone$, we have
\[
\forall \segpt\in \segment{\bdpt}{\conept},\quad
\distnew{\mu_2}{\bdpt}{\segpt} = \distnew{\mu_2}{\bd C}{\segpt}
 < \alpha \distnew{\mu_1}{\vect 0}{\segpt},
\]
which proves Lemma~\ref{blocking_segments_lem} in the case $\ibdpt \ne \vect 0$.

\proofcase{2}{$\ibdpt =\vect 0$ and $\distnew{\mu_1}{\conept}{\bd H} = \distnew{\mu_1}{\conept}{\vect 0}$.}

The hypothesis in this case means that the apex $\vect 0$ is both a $\mu_2$-closest point in $\bd C$ to $\conept$ and a $\mu_1$-closest point in $\bd H$ to $\conept$. Note that since $\conept \ne \vect 0 =\ibdpt$, we must have $\conept\in \interior{C}$ (because $\conept\in \bd C$ implies that $\ibdpt= \conept$). Since $\conept\ne \vect 0$, the function
\[
\ffcn (\vect z) \definedas \frac{\distnew{\mu_2}{\vect z}{\conept}}{\distnew{\mu_1}{\vect 0}{\conept}},
\qquad \vect z\in \Rd,
\]
is continuous on $\Rd$, and since $\ibdpt =\vect 0$, we have
\[
\ffcn (\vect 0)
=  \frac{\distnew{\mu_2}{\vect 0}{\conept}}{\distnew{\mu_1}{\vect 0}{\conept}}
=  \frac{\distnew{\mu_2}{\bd C}{\conept}}{\distnew{\mu_1}{\vect 0}{\conept}}
\le \coneadvantage{1}{\normpair}{C}.
\]
Since $\alpha>\coneadvantage{1}{\normpair}{C}$, there is some open neighborhood $\nbhd$ of $\vect 0$ such that $\ffcn(\vect z) < \alpha$ for all $\vect z\in \nbhd$. Let $\delnbhd = \nbhd\setminus \setof{\vect 0}$, and choose any $\bdpts{1}\in \delnbhd \cap \bd H$ (which is nonempty since $\nbhd \cap \bd H$ is a relatively open subset of $\bd H$ containing $\vect 0$). Then $\bdpts{1}\in \bd H\setminus \setof{\vect 0}$, and $\distnew{\mu_2}{\bdpts{1}}{\conept} < \alpha \distnew{\mu_1}{\vect 0}{\conept}$ since $\bdpts{1}\in \nbhd$. Now for each $\segpt\in \segment{\bdpts{1}}{\conept}$, choose $\hpt\in \segment{\bdpts{1}}{\vect 0}\subset \bd H$ such that the triangles $\trianglev{\bdpts{1}}{\hpt}{\segpt}$ and $\trianglev{\bdpts{1}}{\vect 0}{\conept}$ are similar (i.e.\ if $\segpt = \lambda \bdpts{1} + (1-\lambda) \conept$ with $\lambda\in [0,1]$, then $\hpt = \lambda \bdpts{1}$). Then since $\vect 0$ is a $\mu_1$-closest point in $\bd H$ to $\conept$ by assumption, the similarity of the triangles (or the explicit formula for $\hpt$) implies that $\hpt$ is a $\mu_1$-closest point in $\bd H$ to $\segpt$.
In particular, since $\vect 0\in \bd H$, this implies that
\begin{equation}
\label{blocking_segments:hpt_closer_eqn}
\forall \segpt\in \segment{\bdpts{1}}{\conept},\quad
\distnew{\mu_1}{\hpt}{\segpt} \le \distnew{\mu_1}{\vect 0}{\segpt}.
\end{equation}
Moreover, again using the similarity of the triangles (or the formula for $\hpt$), we have
\begin{equation}
\label{blocking_segments:hpt_ratio_equal_eqn}
\forall \segpt\in \ocsegment{\bdpts{1}}{\conept},\quad
\frac{\distnew{\mu_2}{\bdpts{1}}{\segpt}}{\distnew{\mu_1}{\hpt}{\segpt}}
= \frac{\distnew{\mu_2}{\bdpts{1}}{\conept}}{\distnew{\mu_1}{\vect 0}{\conept}}
=\ffcn \parens[2]{\bdpts{1}} <\alpha.
\end{equation}
Combining \eqref{blocking_segments:hpt_ratio_equal_eqn} and \eqref{blocking_segments:hpt_closer_eqn}, we get
\[
\forall \segpt\in \ocsegment{\bdpts{1}}{\conept},\quad
\distnew{\mu_2}{\bdpts{1}}{\segpt}
< \alpha \distnew{\mu_1}{\hpt}{\segpt}
\le \alpha \distnew{\mu_1}{\vect 0}{\segpt},
\]
and since $\distnew{\mu_2}{\bdpts{1}}{\bdpts{1}} = 0 < \alpha \distnew{\mu_1}{\vect 0}{\bdpts{1}}$, we therefore have
\begin{equation}
\label{blocking_segments:hseg1_closer_eqn}
\forall \segpt\in \segment{\bdpts{1}}{\conept},\quad
\distnew{\mu_2}{\bdpts{1}}{\segpt}
< \alpha \distnew{\mu_1}{\vect 0}{\segpt}.
\end{equation}
Finally, let $\bdpt$ be the $\mu_2$-closest point in $\segment{\bdpts{1}}{\vect x}\setminus \interior{C}$ to $\conept$. Then $\bdpt\in \bd C$ by Lemma~\ref{closest_bd_point_lem}, and we claim that $\bdpt\ne \vect 0$. To see this, observe that $\vect 0\not\in \segment{\bdpts{1}}{\conept}$, because $\bdpts{1}\ne \vect 0$ by definition, and since $\conept \in \interior{C}\subseteq\interior{H}$ and $H$ is convex, we must have  $\ocsegment{\bdpts{1}}{\conept} \subset \interior{H} \subseteq H\setminus \setof{\vect 0}$ by Lemma~\ref{convex_ray_bd_point_lem}. Thus, $\bdpt\in \bd C\setminus \setof{\vect 0}$, and using \eqref{blocking_segments:hseg1_closer_eqn}, we have
\[
\forall \segpt\in \segment{\bdpt}{\conept},\quad
\distnew{\mu_2}{\bdpt}{\segpt}
\le \distnew{\mu_2}{\bdpts{1}}{\segpt}
< \alpha \distnew{\mu_1}{\vect 0}{\segpt},
\]
which proves Lemma~\ref{blocking_segments_lem} in Case 2.

\proofcase{3}{$\ibdpt =\vect 0$ and $\distnew{\mu_1}{\conept}{\bd H} < \distnew{\mu_1}{\conept}{\vect 0}$.}

In this case, let $\bdpts{1}$ be a $\mu_1$-closest point in $\bd H$ to $\conept$, so $\bdpts{1}\ne \vect 0$ by hypothesis. Now observe that since $\conept\ne \vect 0$, the function
\[
\gfcn (\vect z)
\definedas \frac{\distnew{\mu_2}{\vect z}{\conept}}{\distnew{\mu_1}{\vect z}{\conept}},
\qquad \vect z\in \Rd\setminus \setof{\conept},
\]
is continuous at $\vect 0$, and since $\ibdpt = \vect 0$,
\[
\gfcn(\vect 0)
= \frac{\distnew{\mu_2}{\vect 0}{\conept}}{\distnew{\mu_1}{\vect 0}{\conept}}
= \frac{\distnew{\mu_2}{\bd C}{\conept}}{\distnew{\mu_1}{\vect 0}{\conept}}
\le \coneadvantage{1}{\normpair}{C}.
\]
Since $\alpha>\coneadvantage{1}{\normpair}{C}$, there is some open neighborhood $\nbhd$ of $\vect 0$ such that $\gfcn(\vect z) < \alpha$ for all $\vect z\in \nbhd$.
Choose any $\bdpts{2}\in \nbhd \cap \ocsegment{\vect 0}{\bdpts{1}}$, which is nonempty since $\bdpts{1}\ne \vect 0$. Then since for any $\segpt\in \ocsegment{\bdpts{2}}{\conept}$, the vectors $\segpt-\bdpts{2}$ and $\conept-\bdpts{2}$ point in the same direction, the scale-equivariance of the norm metrics implies that
\begin{equation}
\label{blocking_segments:hpt2_ratio_equal_eqn}
\forall \segpt\in \ocsegment{\bdpts{2}}{\conept},\qquad
\frac{\distnew{\mu_2}{\bdpts{2}}{\segpt}}{\distnew{\mu_1}{\bdpts{2}}{\segpt}}
= \frac{\distnew{\mu_2}{\bdpts{2}}{\conept}}{\distnew{\mu_1}{\bdpts{2}}{\conept}}
= \gfcn(\bdpts{2}) <\alpha.
\end{equation}
Now, for each $\segpt\in \segment{\bdpts{2}}{\conept}$, choose $\hpt\in \segment{\bdpts{2}}{\bdpts{1}} \subset \bd H$ so that the triangles $\trianglev{\bdpts{2}}{\hpt}{\segpt}$ and $\trianglev{\bdpts{2}}{\bdpts{1}}{\conept}$ are similar (i.e.\ if $\segpt = \lambda \bdpts{2} + (1-\lambda) \conept$ with $\lambda\in [0,1]$, then $\hpt = \lambda \bdpts{2} + (1-\lambda)\bdpts{1}$). Then since $\bdpts{1}$ is a $\mu_1$-closest point in $\bd H$ to $\conept$, the similarity of the triangles (or the formula for $\hpt$) implies that $\hpt$ is a $\mu_1$-closest point in $\bd H$ to $\segpt$. Therefore, Lemma~\ref{support_point_closest_lem} implies that $\bd H$ is a support hyperplane of $\reacheddset{\mu_1}{\conept}{r_{\segpt}}$ at $\hpt$, where $r_{\segpt} = \distnew{\mu_1}{\segpt}{\bd H}$. We claim that Lemma~\ref{dist_support_point_lem} then implies that
\begin{equation}
\label{blocking_segments:hpt2_closer_eqn}
\forall \segpt\in \segment{\bdpts{2}}{\conept},\qquad
\distnew{\mu_1}{\bdpts{2}}{\segpt} \le \distnew{\mu_1}{\vect 0}{\segpt}.
\end{equation}
To see this, let $\dirsymb = -\bdpts{1}\ne \vect 0$. Then $\dirvec$ is a direction in $\bd H$, 
and for any $\segpt\in \segment{\bdpts{2}}{\conept}$, the four points $\bdpts{1}$, $\hpt$, $\bdpts{2}$, $\vect 0$ are collinear along the ray $\bdpts{1}+\dirvec$, appearing in the listed order as we travel in the direction $\dirvec$ starting from $\bdpts{1}$. Therefore, for fixed $\segpt\in \segment{\bdpts{2}}{\conept}$, we have $\bdpts{2} = \hpt + s \dirsymb$ and $\vect 0 = \hpt + s' \dirsymb$ for some $0\le s< s'$, and since $\bd H$ is a support hyperplane of $\reacheddset{\mu_1}{\conept}{r_{\segpt}}$ at $\hpt$, Lemma~\ref{dist_support_point_lem} implies that
\[
\distnew{\mu_1}{\bdpts{2}}{\segpt}
= \distnew{\mu_1}{\hpt + s \dirsymb}{\segpt}
\le \distnew{\mu_1}{\hpt + s' \dirsymb}{\segpt}
=\distnew{\mu_1}{\vect 0}{\segpt},
\]
proving \eqref{blocking_segments:hpt2_closer_eqn}. Combining \eqref{blocking_segments:hpt2_ratio_equal_eqn} and \eqref{blocking_segments:hpt2_closer_eqn}, we get
\[
\forall \segpt\in \ocsegment{\bdpts{2}}{\conept},\qquad
\distnew{\mu_2}{\bdpts{2}}{\segpt}
< \alpha \distnew{\mu_1}{\bdpts{2}}{\segpt}
\le \alpha \distnew{\mu_1}{\vect 0}{\segpt},
\]
and since $\distnew{\mu_2}{\bdpts{2}}{\bdpts{2}} = 0 < \alpha \distnew{\mu_1}{\vect 0}{\bdpts{2}}$, we therefore have
\begin{equation}
\label{blocking_segments:hseg2_closer_eqn}
\forall \segpt\in \segment{\bdpts{2}}{\conept},\qquad
\distnew{\mu_2}{\bdpts{2}}{\segpt}
< \alpha \distnew{\mu_1}{\vect 0}{\segpt}.
\end{equation}
Finally, let $\bdpt$ be the $\mu_2$-closest point in $\segment{\bdpts{2}}{\vect x}\setminus \interior{C}$ to $\conept$. Then it follows as in Case 2 that $\bdpt\in \bd C\setminus \setof{\vect 0}$, and using \eqref{blocking_segments:hseg2_closer_eqn} we get
\[
\forall \segpt\in \segment{\bdpt}{\conept},\qquad
\distnew{\mu_2}{\bdpt}{\segpt}
\le \distnew{\mu_2}{\bdpts{2}}{\segpt}
< \alpha \distnew{\mu_1}{\vect 0}{\segpt},
\]
which proves Lemma~\ref{blocking_segments_lem} in Case 3.
\end{proof}

We now prove a strengthened version of Lemma~\ref{blocking_segments_lem}, showing that we can move each of the $\vect y_{\vect x}$'s within some small open set without affecting the conclusion of the lemma.

\begin{lem}[Wiggle room for blocking segments]
\label{strong_blocking_segments_lem}
\newcommand{\bdpoint}[1]{\vect y_{#1}}
\newcommand{\ballpoint}[1]{\bdpoint{#1}'}
\newcommand{\segpoint}{\vect y}
\newcommand{\radius}[1]{\epsilon_{#1}}
Let $C$ be a closed cone at $\vect 0$ which is contained in some half-space, and let $\alpha> \coneadvantage{1}{\normpair}{C}$. Then for any $\vect x\in C\setminus \setof{\vect 0}$, there exists $\vect y_{\vect x}\in \bd C\setminus \setof{\vect 0}$ and $\radius{\vect x}>0$ 
such that for all $\vect y_{\vect x}'\in \reacheddset{\mu_2}{\vect y_{\vect x}}{\radius{\vect x}}$,
\[
\sup_{\segpoint \in \segment{\ballpoint{\vect x}}{\vect x}}
\frac{\distnew{\mu_2}{\ballpoint{\vect x}}{\segpoint}}
{\distnew{\mu_1}{\vect 0}{\segpoint}}
<\alpha.
\]
\end{lem}

\begin{proof}
\newcommand{\bluntcone}{C\setminus \setof{\vect 0}}
\newcommand{\conept}{\vect x} 
\newcommand{\bdpt}{\vect y_{\conept}} 
\newcommand{\segpt}{\vect y} 


\newcommand{\smallalpha}{\alpha_0} 
\newcommand{\ballpt}{\bdpt'} 

\newcommand{\sepdist}{\delta_{\conept}} 
\newcommand{\radius}{\epsilon_{\conept}} 

\newcommand{\param}{\lambda} 
\newcommand{\psegpt}{\vect y_\param} 
\newcommand{\bpsegpt}{\psegpt'} 


Let $\smallalpha = \frac{1}{2} \parens[2]{\alpha + \coneadvantage{1}{\normpair}{C}}$, and fix $\conept\in \bluntcone$. Then $\smallalpha >  \coneadvantage{1}{\normpair}{C}$, so Lemma~\ref{blocking_segments_lem} implies that there exists $\bdpt\in \bd C\setminus \setof{\vect 0}$ such that
\begin{equation}
\label{strong_blocking:start_ineq_eqn}
\distnew{\mu_2}{\bdpt}{\segpt} < \smallalpha \distnew{\mu_1}{\vect 0}{\segpt}
\quad \text{for all $\segpt\in \segment{\bdpt}{\conept}$.}
\end{equation}
Note that \eqref{strong_blocking:start_ineq_eqn} implies $\distnew{\mu_1}{\vect 0}{\segpt}> 0$ for all $\segpt\in \segment{\bdpt}{\conept}$, so since the line segment is compact we have
\[
\sepdist \definedas \inf_{\segpt\in \segment{\bdpt}{\conept}} \mu_1(\segpt) >0.
\]
Choose any $\radius$ with
 \[
 0< \radius < \frac{(\alpha-\smallalpha)\sepdist}{2+ \alpha \esupnorm{\mu_1/\mu_2}},
 \]
and let $\ballpt \in \reacheddset{\mu_2}{\bdpt}{\radius}$. For each $\param\in [0,1]$, let
\[
\psegpt \definedas \param\bdpt + (1-\param) \conept
\quad\text{and}\quad
\bpsegpt \definedas \param\ballpt + (1-\param) \conept,
\]
so $\segment{\bdpt}{\conept} = \setbuilder[0]{\psegpt}{\param\in [0,1]}$ and $\segment{\ballpt}{\conept} = \setbuilder[0]{\bpsegpt}{\param\in [0,1]}$, and the triangles $\trianglev{\conept}{\bdpt}{\ballpt}$ and $\trianglev{\conept}{\psegpt}{\bpsegpt}$ are similar for any $\param$. Then for all $\param\in [0,1]$ we have
\begin{equation}
\label{strong_blocking:segments_close_eqn}
\distnew{\mu_i}{\psegpt}{\bpsegpt}
= \param \distnew{\mu_i}{\bdpt}{\ballpt}
\le \distnew{\mu_i}{\bdpt}{\ballpt}
\quad\text{for } i\in\setof{1,2}.
\end{equation}
Using the reverse triangle inequality, \eqref{strong_blocking:segments_close_eqn}, the definition of $\sepdist$, and the assumption that $\distnew{\mu_2}{\bdpt}{\ballpt} \le \radius$, we have
\begin{align}
\label{strong_blocking:reverse_triangle_eqn}
\distnew{\mu_1}{\vect 0}{\bpsegpt}
&\ge \distnew{\mu_1}{\vect 0}{\psegpt} - \distnew{\mu_1}{\psegpt}{\bpsegpt}
	\notag\\
&\ge \distnew{\mu_1}{\vect 0}{\psegpt} - \distnew{\mu_1}{\bdpt}{\ballpt}
	\notag\\
&\ge \sepdist - \esupnorm[1]{\mu_1/\mu_2} \radius.
\end{align}
Now, for all $\param\in [0,1]$,
\begin{align*}
\distnew{\mu_2}{\ballpt}{\bpsegpt}
&\le \distnew{\mu_2}{\ballpt}{\bdpt} + \distnew{\mu_2}{\bdpt}{\psegpt}
	+ \distnew{\mu_2}{\psegpt}{\bpsegpt}
	&& \text{(triangle inequality)}\\
&\le \radius + \smallalpha \distnew{\mu_1}{\vect 0}{\psegpt}
	+ \distnew{\mu_2}{\bdpt}{\ballpt}
	&& \text{(using \eqref{strong_blocking:start_ineq_eqn},
		\eqref{strong_blocking:segments_close_eqn})}\\
&\le 2\radius + \smallalpha \distnew{\mu_1}{\vect 0}{\psegpt}
	&& \text{(since $\ballpt\in \reacheddset{\mu_2}{\bdpt}{\radius}$)}\\
&\le 2\radius + \smallalpha \distnew{\mu_1}{\vect 0}{\bpsegpt}
	+  \smallalpha \distnew{\mu_1}{\bpsegpt}{\psegpt}
	&& \text{(triangle inequality)}\\
&\le 2\radius + \smallalpha \distnew{\mu_1}{\vect 0}{\bpsegpt}
	+  \smallalpha \distnew{\mu_1}{\ballpt}{\bdpt}
	&& \text{(by \eqref{strong_blocking:segments_close_eqn})}\\
&\le 2\radius + \smallalpha \distnew{\mu_1}{\vect 0}{\bpsegpt}
	+\smallalpha \esupnorm[1]{\mu_1/\mu_2} \radius
	&& \text{(since $\ballpt\in \reacheddset{\mu_2}{\bdpt}{\radius}$)}\\
&= \smallalpha \distnew{\mu_1}{\vect 0}{\bpsegpt}
	+ \parens[2]{2 +\alpha 
	\esupnorm[1]{\mu_1/\mu_2}
	} \radius\\
	&\hspace{3.1cm} -(\alpha -\smallalpha ) 
	 \esupnorm[1]{\mu_1/\mu_2} \radius\\
&< \smallalpha \distnew{\mu_1}{\vect 0}{\bpsegpt}
	+(\alpha - \smallalpha) \parens[2]{\sepdist - \esupnorm[1]{\mu_1/\mu_2}\radius}
	&& \text{(choice of $\radius$)}\\
&\le \smallalpha \distnew{\mu_1}{\vect 0}{\bpsegpt}
	+ (\alpha - \smallalpha) \distnew{\mu_1}{\vect 0}{\bpsegpt}
	&& \text{(by \eqref{strong_blocking:reverse_triangle_eqn})}\\
&= \alpha \distnew{\mu_1}{\vect 0}{\bpsegpt}.
\end{align*}
Since $\param\in [0,1]$ was arbitrary, this shows that
\[
\forall \segpt \in \segment{\ballpt}{\conept},\qquad
\distnew{\mu_2}{\ballpt}{\segpt}
<\alpha \distnew{\mu_1}{\vect 0}{\segpt},
\]
and since $\ballpt\in \reacheddset{\mu_2}{\bdpt}{\radius}$ was arbitrary, this proves Lemma~\ref{strong_blocking_segments_lem}.
\end{proof}


\section{Proofs from Chapter~\ref{random_fpc_chap}}

The following lemma is needed for the proof of Theorem~\ref{cone_segment_conquering_thm}.

\begin{lem}[Lemma~\ref{growth_containment_lem}: Containment of growth within large $\mu_2$-balls far from the origin]
\newcommand{\thisevent}{\eventref{growth_containment_lem}\parens{r,\beta;\alpha}}
\newcommand{\thiskonst}{k_{\mu_1}}

For any constants $\alpha\in (0,1)$, $\beta\ge 1$, and $r>0$, define the event
\[
\thisevent =
 \bigcap_{\vect z\in \cubify{\vect 0}} \bigcap_{\substack{\vect x\in \Rd \\ \mu_1(\vect x) \ge r}}
\eventthat[1]{ \rreachedset{\tmeasure_2}{\vect z+\vect x}{\Rd}{\alpha t} \subseteq
\reacheddset{\mu_2}{\vect z+\vect x}{t-} \;\forall t\ge \beta \mu_1(\vect x) },
\]
and let $\thiskonst = 3\dilnorm{\lspace{\infty}{d}}{\mu_1}+1$. Then given $\alpha \in (0,1)$, there exist positive constants $C_{\ref{growth_containment_lem}}$ and $c_{\ref{growth_containment_lem}}$ such that for any $\beta\ge 1$ and $r \ge 0$,
\[
\Pr \parens[2]{\thisevent} \ge 1- \bigcref{growth_containment_lem}(\alpha)
e^{-\smallcref{growth_containment_lem}(\alpha)\cdot \beta \cdot (r - \thiskonst)}.
\]
\end{lem}

\begin{proof}
\newcommand{\thisevent}{\eventref{growth_containment_lem}\parens{r,\beta;\alpha}}
\newcommand{\thisC}{\bigcref{growth_containment_lem}(\alpha)}
\newcommand{\thisc}{\smallcref{growth_containment_lem}(\alpha)}

\newcommand{\thiskonst}{k_{\mu_1}}

\newcommand{\ballevent}[4][0]{\eventsymb^{#2} \parens[#1]{#3;#4}}
\newcommand{\Calpha}{\bigcsymb(\alpha)}
\newcommand{\calpha}{\smallcsymb(\alpha)}

\newcommand{\thisz}{\vect z}
\newcommand{\thisy}{\vect y}
\newcommand{\thisx}{\vect x}
\newcommand{\thisv}{\vect v}

\newcommand{\thiscube}{\cubify{\vect 0}}

\newcommand{\thisdil}{\kappa}
\newcommand{\thisshellbound}{\normshellbound{\mu_1}}


Fix $\alpha\in (0,1)$, and for each $\vect v\in\Zd$ and $r\ge 0$, define the event
\[
\ballevent{\thisv}{r}{\alpha} \definedas \bigcap_{\vect y\in \cubify{\vect v}}
\eventthat[1]{\rreachedset{\tmeasure_2}{\thisv}{\Rd}{\alpha t} \subseteq \reacheddset{\mu_2}{\vect y}{t-}
\;\forall t\ge r}.
\]
Then $\ballevent{\thisv}{r}{\alpha}\supseteq \eventref{shifted_dev_bounds_lem}^{\thisv} \parens[2]{\alpha r; \alpha^{-1}-1}$, with $\tmeasure = \tmeasure_2$ and $\mu=\mu_2$. Thus, since $\alpha<1$, Lemma~\ref{shifted_dev_bounds_lem} implies there are positive constants $\Calpha$ and $\calpha$ (depending on $\law \parens[2]{\tmeasure_2 (\edge)}$) such that for any $\vect v\in\Zd$ and $r\ge 0$,
\begin{equation}
\label{growth_containment:ballevent_prob_eqn}
\Pr \parens[2]{\ballevent{\thisv}{r}{\alpha}} \ge 1- \Calpha e^{-\calpha r}.
\end{equation}
Let $\thisdil = \dilnorm{\lspace{\infty}{d}}{\mu_1}$. Then using \eqref{space_covering:yz_bounds} and \eqref{space_covering:yv_bounds} (from the proof of Lemma~\ref{space_covering_lem}, with $\mu_1$ in place of $\mu$) in lines three and four, respectively, we have
\begin{align*}
\eventcomplement{\thisevent}
&= \bigcup_{\thisz\in\thiscube} \bigcup_{\substack{\thisx\in\Rd \\ \mu_1(\thisx)\ge r}}
	\eventthat[3]{ \rreachedset{\tmeasure_2}{\thisz+\thisx}{\Rd}{\alpha t} \not\subseteq
	\reacheddset{\mu_2}{\thisz+\thisx}{t-} \text{ for some } t\ge \beta \mu_1(\thisx) }\\
&= \bigcup_{\thisz\in\thiscube} \bigcup_{\substack{\thisy\in\Rd \\ \mu_1(\thisy-\thisz)\ge r}}
	\eventthat[3]{ \rreachedset{\tmeasure_2}{\thisy}{\Rd}{\alpha t} \not\subseteq
	\reacheddset{\mu_2}{\thisy}{t-} \text{ for some } t\ge \beta \mu_1(\thisy-\thisz) }\\
&\subseteq \bigcup_{\substack{\thisy\in\Rd \\ \mu_1(\thisy)\ge r - \thisdil/2}}
	\eventthat[3]{ \rreachedset{\tmeasure_2}{\thisy}{\Rd}{\alpha t} \not\subseteq
	\reacheddset{\mu_2}{\thisy}{t-} \text{ for some }
	t\ge \beta \parens[2]{\mu_1(\thisy) - \thisdil/2} }\\
&\subseteq \bigcup_{\substack{\thisv\in\Zd \\ \mu_1(\thisv)\ge r - \thisdil}}
	\bigcup_{\thisy\in \cubify{\thisv}}
	\eventthat[3]{ \rreachedset{\tmeasure_2}{\thisv}{\Rd}{\alpha t} \not\subseteq
	\reacheddset{\mu_2}{\thisy}{t-} \text{ for some }
	t\ge \beta \parens[2]{\mu_1(\thisv) - \thisdil} }\\
&=\bigcup_{\substack{\thisv\in\Zd \\ \mu_1(\thisv)\ge r - \thisdil}}
	\eventcomplement{\ballevent[2]{\thisv}{\beta \parens[0]{\mu_1(\thisv) - \thisdil}}{\alpha}}.
\end{align*}
Therefore, \eqref{growth_containment:ballevent_prob_eqn} and Lemmas~\ref{norm_shell_lem} and \ref{poly_geo_series_lem} imply that for all $r\ge 0$,
\begin{align*}
\Pr \parens[1]{\eventcomplement{\thisevent}}
&\le \sum_{\substack{\thisv\in\Zd \\ \mu_1(\thisv)\ge r - \thisdil}}
	\Calpha e^{-\calpha \beta \parens{\mu_1(\thisv) - \thisdil}}
	&& \text{(by \eqref{growth_containment:ballevent_prob_eqn})}\\
&\le \sum_{j = \floor{r-\thisdil}\join 0}^\infty
	\sum_{\substack{\thisv\in\Zd \\ j\le \mu_1(\thisv)<j+1}}
	\Calpha e^{-\calpha \beta \parens{\mu_1(\thisv) - \thisdil}}\\
&\le \sum_{j = \floor{r-\thisdil}\join 0}^\infty
	\thisshellbound (j+1)^{d-1}
	\Calpha e^{-\calpha \beta (j - \thisdil)}
	&& \text{(Lem.~\ref{norm_shell_lem})}\\
&\le \thisshellbound \Calpha e^{\calpha \beta \thisdil}
	\cdot \frac{4(2d-2)^{d-1}}{\parens[2]{[\calpha \beta]\meet 1}^d}
	\cdot e^{-\frac{\calpha}{2} \beta \parens{\floor{r-\thisdil}\join 0}}
	&& \text{(Lem.~\ref{poly_geo_series_lem})}\\
&\le \thisshellbound \Calpha e^{\calpha \beta \thisdil}
	\cdot \frac{4(2d-2)^{d-1}}{\parens[2]{\calpha \meet 1}^d}
	\cdot e^{-\frac{\calpha}{2} \beta \parens{r-\thisdil-1}}
	&& \text{($\because \beta\ge 1$)}\\
&= \frac{4\thisshellbound \Calpha(2d-2)^{d-1}}{\parens[2]{\calpha \meet 1}^d}
	\cdot e^{-\frac{\calpha}{2} \beta \parens{r-3\thisdil-1}}\\
&= \thisC e^{-\thisc \beta (r-\thiskonst)},
\end{align*}
where
\[
\thisC \definedas
\frac{4\thisshellbound \Calpha(2d-2)^{d-1}}{\parens[2]{\calpha \meet 1}^d},
\quad\text{and}\quad
\thisc \definedas \frac{\calpha}{2}.
\qedhere
\]
\end{proof}


\chapter{Supplementary Results}
\label{notmyjob_chap}

\section{Geometry of $\Rd$ and $\Zd$}


The following observation is invoked in several proofs.

\begin{lem}[Closest points are boundary points]
\label{closest_bd_point_lem}
\newcommand{\thisnorm}{\norm{\cdot}}
\newcommand{\closedset}{A}
\newcommand{\otherset}{B}
\newcommand{\thispt}{\vect z}
\newcommand{\closept}{\vect y}
\newcommand{\thissegment}{\ocsegment{\closept}{\thispt}}

Let $\thisnorm$ be a norm on $\Rd$, let $\otherset$ be a nonempty subset of $\Rd$, and suppose that $\thispt\in \Rd\setminus \otherset$. Then for any closed star set $\closedset\in \stars{\thispt}$ with $\closedset\cap \closure{\otherset}\ne\emptyset$, there exists $\closept\in \closedset\cap \closure{\otherset}$ which is a $\thisnorm$-closest point in $\closedset\cap \closure{\otherset}$ to $\thispt$, and for any such $\vect y$ we have $\closept\in \bd \otherset$ and $\thissegment \subseteq \exterior{\otherset}$.
\end{lem}

\begin{proof}
\newcommand{\thisnorm}{\norm{\cdot}}
\newcommand{\closedset}{A}
\newcommand{\otherset}{B}
\newcommand{\cotherset}{\closure{\otherset}}
\newcommand{\thisext}{\exterior{\otherset}}

\newcommand{\thispt}{\vect z}
\newcommand{\closept}{\vect y}

\newcommand{\thissegment}{\ocsegment{\closept}{\thispt}}
\newcommand{\csegment}{\segment{\closept}{\thispt}}
\newcommand{\param}{\alpha}
\newcommand{\segpt}{\vect y_{\param}}

The point $\closept$ exists because $\closedset\cap \cotherset\ne\emptyset$ is closed and $\setof{\thispt}$ is compact; see \cite[Exercise~1.6.10, p.~15]{Burago:2001aa}.
If $\thispt\in \cotherset$, then the unique $\thisnorm$-closest point in $\cotherset$ to $\thispt$ is $\thispt$ itself, in which case we trivially have $\closept=\thispt\in \cotherset \cap \parens[2]{\Rd\setminus \otherset} \subseteq \bd\otherset$, and $\thissegment = \emptyset \subseteq \thisext$.
Thus, suppose $\thispt\notin \cotherset$, so $\thispt\in \thisext$.
Then for any point $\closept$ which is a $\thisnorm$-closest point in $\closedset\cap \cotherset$ to $\thispt$, we have $\closept\ne\thispt$, so $\thissegment\ne\emptyset$, and we claim that $\thissegment\subseteq\exterior{\otherset}$.
To see this, first note that since $\closept\in \closedset$, we must have $\csegment\subseteq \closedset$ because $\closedset$ is star-shaped at $\thispt$ by assumption.
Now, for $\param\in [0,1]$ let $\segpt \definedas \param \closept + (1-\param)\thispt$, so $\csegment = \setbuilder{\segpt}{\param\in [0,1]}$. Then
\begin{equation*}
\distnew{\thisnorm}{\segpt}{\thispt} = \param \distnew{\thisnorm}{\closept}{\thispt}
\quad\text{for all $\param\in [0,1]$,}
\end{equation*}
so we have
\begin{equation}
\label{closest_bd_point:smaller_dist_eqn}
\distnew{\thisnorm}{\segpt}{\thispt} < \distnew{\thisnorm}{\closept}{\thispt}
\quad\text{ for all $\param\in [0,1)$.}
\end{equation}
Thus, the set $\thissegment = \setbuilder{\segpt}{\param\in [0,1)} \subseteq\closedset$ cannot contain any points in $\cotherset$, because if it did, then \eqref{closest_bd_point:smaller_dist_eqn} would contradict the fact that $\closept$ is a $\thisnorm$-closest point in $\closedset\cap \cotherset$ to $\thispt$.
Therefore, $\thissegment\subseteq\thisext$ as claimed. Since $\thissegment\ne\emptyset$ by assumption, this implies that there are points in $\Rd\setminus \otherset$ that are arbitrarily close to $\closept$, and since $\closept\in \cotherset$ by definition, we must have $\closept\in \bd \otherset$.
\end{proof}

Recall that for any $A\subseteq\Rd$, we define $\cubify{A} \definedas A+\segment[1]{-\frac{1}{2}}{\frac{1}{2}}^d$ (the \textdef{cube expansion} of $A$) and $\lat A \definedas \Zd\cap \cubify{A}$ (the \textdef{lattice approximation} of $A$).
The following lemma enumerates some elementary properties of these operations; the proof is straightforward and is left to the reader.

\begin{lem}[Properties of cube expansion and lattice approximation]
\label{lattice_cube_properties_lem}
Let $\vect v\in V\subseteq\Zd$ and $\vect z\in A\subseteq\Rd$. Then
\begin{enumerate}
\item $\vect v\in \lat{\vect z}$ if and only if $\vect z\in\cubify{\vect v}$.
\item $\lat V= V$ and $\lat A \supseteq A\cap \Zd$.
\item $\lat[2]{\intcubify{V}} = V$ and $\intcubify[2]{\,\lat{A}\,} \supseteq A$.
\item If $A$ is a connected subset of $\Rd$, then $\lat A$ is a connected subgraph of $\Zd$.
\item $V$ is a connected subgraph of $\Zd$ if and only if $\intcubify{V}$ is a (path) connected subset of $\Rd$.
\item $V\cap \lat{A} = \emptyset$ if and only if $\cubify{V}\cap A = \emptyset$.
\item $\Rd\setminus \intcubify{V} = \cubify[1]{\Zd\setminus V}$ and $\Rd\setminus \cubify{V} = \intcubify[1]{\Zd\setminus V}$.
\item $\lat A\subseteq V$ if and only if $A\subseteq \intcubify{V}$.
\item $\lat[2]{\Rd\setminus \cubify{V}} = \Zd\setminus V$ and $\cubify[2]{\Zd\setminus \lat{A}\,} \subseteq \Rd\setminus A$.
\item $\card{V} = \lebmeasof[2]{d}{\cubify{V}} = \lebmeasof[2]{d}{\intcubify{V}}$, where $\lebmeas{d}$ is $d$-dimensional Lebesgue measure.
\item $\lat A \subseteq \cubify{A} = A+\frac{1}{2} \unitball{\linfty{d}}$, and $\cubify[2]{\lat A}\subseteq \cubify{\cubify{A}} =  A+\unitball{\linfty{d}}$.
\item If $A=\bigcup_{i\in I} A_i$, where $I$ is an arbitrary index set and $A_i\subseteq \Rd$, then $\cubify{A} = \bigcup_{i\in I} \cubify{A_i}$ and $\lat{A} = \bigcup_{i\in I} \lat{A_i}$.
\end{enumerate}
\end{lem}


The following estimate is used in Chapters~\ref{cone_growth_chap} and \ref{random_fpc_chap}.

\begin{lem}[The number of lattice points in a spherical shell]
\label{norm_shell_lem}
For any norm $\norm{\cdot}$ on $\Rd$, there is some constant $\normshellbound{\norm{\cdot}}<\infty$ such that for any $\vect x\in\Rd$ and $0\le a\le b$,
\[
\card[1]{\setbuilder[1]{\vect v\in\Zd}{a\le \norm{\vect v-\vect x}\le b}}
\le \normshellbound{\norm{\cdot}} (b\join 1)^{d-1} \brackets[2]{(b-a)\join 1}.
\]
\end{lem}

\begin{proof}
Let $V$ be the set of lattice points in the statement. By Lemma~\ref{lattice_cube_properties_lem}, $\card{V} = \lebmeasof[2]{d}{\cubify{V}}$, where $\lebmeas{d}$ denotes Lebesgue measure in $\Rd$. Since $\cubify{V}$ is contained in some $\norm{\cdot}$-ball of radius $\approx b$ and contains some $\norm{\cdot}$-ball of radius $\approx a$, the result follows from the fact that $\lebmeasof[2]{d}{\reacheddset{\norm{\cdot}}{\vect x}{r}} = \Theta(r^d)$.
%
%
%
%
\end{proof}

The following lemma is needed in the proof of Theorem~\ref{narrow_survival_time_thm} (specifically Claim~\ref{narrow_survival_time:blocking_event_claim}).


\begin{lem}[Separation of $\Zd$ via boundaries in $\Rd$]
\label{lat_euc_bd_separation_lem}
Let $B\subseteq \Rd$, and suppose that $B$ and $\Rd\setminus B$ are both nonempty. Then $\lat[1]{\bd B}$ separates $\lat{B}$ from $\lat[2]{\Rd\setminus B}$ in $\Zd$. That is, any $\Zd$-path from $\lat{B}$ to $\lat[2]{\Rd\setminus B}$ must intersect $\lat[1]{\bd B}$.
\end{lem}

\begin{proof}
\newcommand{\thisset}{B}
\newcommand{\bdset}{\bd \thisset}
\newcommand{\setcomp}{\Rd\setminus \thisset}

\newcommand{\latset}{\lat{\thisset}}
\newcommand{\latcomp}{\lat[2]{\Rd\setminus \thisset}}
\newcommand{\latbd}{\lat[1]{\bd \thisset}}

\newcommand{\invert}{\vect u}
\newcommand{\outvert}{\vect v}
\newcommand{\midvert}{\vect w}

\newcommand{\bdpt}{\vect x}
\newcommand{\bdptt}{\vect y}
\newcommand{\segbdpt}{\vect z}

\newcommand{\thispath}{\gamma}

Let $\thispath$ be a lattice path from $\latset$ to $\latcomp$, and let $\invert$ be the last vertex of $\thispath$ lying in $\latset$. We will show that $\invert\in \latbd$.
%
Recall that $\latbd = \Zd\cap \parens[2]{\bdset + \frac{1}{2} \unitball{\linfty{d}}}$; since $\bd B$ is closed, it follows that $\invert\in \latbd$ if and only if $\distnew{\linfty{d}}{\invert}{\bdset}\le 1/2$.
First suppose $\invert\notin \thisset$. Since $\invert\in \latset$ by assumption, we have $\distnew{\linfty{d}}{\invert}{\thisset} \le 1/2$. Thus, if $\bdpt$ is any $\linfty{d}$-closest point in $\closure{\thisset}$ to $\invert$, then $\bdpt\in \bdset$ by Lemma~\ref{closest_bd_point_lem}, so we have
\[
\distnew{\linfty{d}}{\invert}{\bdset}
= \distnew{\linfty{d}}{\invert}{\bdpt} 
=\distnew[2]{\linfty{d}}{\invert}{\closure{\thisset}\,}\le 1/2,
\]
and hence $\invert\in \latbd$.

Now suppose $\invert\in \thisset$. If $\invert$ is the final vertex of $\thispath$, then we must have $\invert\in \latcomp$, and hence $\distnew[2]{\linfty{d}}{\invert}{\,\setcomp} \le 1/2$. Thus, if $\bdptt$ is any $\linfty{d}$-closest point in $\shellof[0]{\thisset}$ to $\invert$, then $\bdptt\in \bdset$ by Lemma~\ref{closest_bd_point_lem}, so we have
\[
\distnew{\linfty{d}}{\invert}{\bdset}
= \distnew{\linfty{d}}{\invert}{\bdptt} 
=\distnew[2]{\linfty{d}}{\invert}{\,\shellof[0]{\thisset}}\le 1/2,
\]
and hence $\invert\in \latbd$ as before.
If $\invert$ is not the final vertex of $\thispath$, let $\outvert$ be the vertex in $\thispath$ that follows $\invert$, so $\outvert$ is adjacent to $\invert$ in $\Zd$, and $\outvert$ necessarily lies in $\Zd\setminus \latset$ since $\invert$ was the \emph{last} vertex of $\thispath$ in $\latset$. Since $\Zd\cap \thisset\subseteq \latset$, this implies that $\outvert\notin\thisset$. Since $\invert\in \thisset$ and $\outvert\notin\thisset$, the line segment $\segment{\invert}{\outvert}$ intersects $\bdset$ by Lemma~\ref{closest_bd_point_lem}; let $\segbdpt\in \segment{\invert}{\outvert} \cap \bdset$. Since $\invert$, $\segbdpt$, and $\outvert$ are collinear (occurring in that order), we have
\[
\distnew{\linfty{d}}{\invert}{\segbdpt} + \distnew{\linfty{d}}{\segbdpt}{\outvert}
= \distnew{\linfty{d}}{\invert}{\outvert} =1,
\]
where the final equality holds because $\invert$ and $\outvert$ are adjacent in $\Zd$. Since $\segbdpt\in \closure{\thisset}$, we have
\[
\distnew{\linfty{d}}{\segbdpt}{\outvert}
\ge \distnew{\linfty{d}}{\thisset}{\outvert}\ge1/2,
\]
where the final inequality holds because $\outvert\notin \latset$. Therefore,
\[
\distnew{\linfty{d}}{\invert}{\bdset}
\le \distnew{\linfty{d}}{\invert}{\segbdpt}
=1- \distnew{\linfty{d}}{\segbdpt}{\outvert}
\le 1/2,
\]
so again, $\invert\in \latbd$.
%
%
%
%
%
%
%
\end{proof}

\section{Graph Theory}
\label{graph_theory_supp_sec}

\begin{lem}[Edge sets]
\label{edge_sets_lem}
Let $U,V\subseteq \Zd$.
\begin{enumerate}
\item $\edges{U\cup V} = \edges{U} \cup \edges{V} \cup \jedges{U}{V}$.
\item $\bdedges{U\cup V} = \bdedges{U} \cup \bdedges{V} \setminus \jedges{U}{V}$.
\item $\compedges{U\cup V} = \compedges{U} \cap \compedges{V} \setminus \jedges{U}{V}$.
\item $\staredges{U\cup V} = \staredges{U}\cup \staredges{V}$.
\item $\edges{U\cap V} = \edges{U} \cap \edges{V}$, and $\compedges{U\cap V} = \compedges{U}\cup \compedges{V}$.
\item $\edges{U\cap V} \subseteq \jedges{U}{V}$.
\item $\staredges{V}\subseteq \edges[2]{\nbrhood{V}}$.
\end{enumerate}
\end{lem}

%
%
%

\begin{lem}[Boundaries and neighbor sets]
\label{bd_nbr_lem}
For any $U,V\subseteq\Zd$,
\[
\bd (U\cup V) = \setbuilder[2]{\vect v\in \bd U\cup \bd V}{\nbrset{\vect v} \not\subseteq U\cup V}.
\]
\end{lem}

\begin{proof}
By definition, $\bd (U\cup V) = \setbuilder[2]{\vect v\in U\cup V}{\nbrset{\vect v} \not\subseteq U\cup V}$. Thus, the inclusion $\supseteq$ is trivial, and it remains to show that $\bd (U\cup V) \subseteq \bd U\cup \bd V$. Suppose $\vect v\in \bd (U\cup V)$, i.e.\ $\vect v\in U\cup V$ and $\nbrset{\vect v} \not\subseteq U\cup V$. If $\vect v\in U$, then $\vect v\in \bd U$ since $\nbrset{\vect v}\not\subseteq U$, and if $\vect v\in V$, then $\vect v\in \bd V$ since $\nbrset{\vect v}\not\subseteq V$. Thus $\vect v\in \bd U\cup \bd V$.
\end{proof}

\begin{lem}[Minimal paths]
\label{minimal_path_lem}
\newcommand{\pclass}[2]{\pathclass_{#1\to #2}}
Let  $\pathclass$ be a class of lattice paths that is closed under taking subpaths, and for $U,V,S\subseteq \Zd$, let $\pclass{U}{V}$ be the set of $\pathclass$-paths from $U$ to $V$.
\begin{enumerate}
\item A path $\gamma\in \pclass{U}{V}$ is minimal if and only if $\gamma$ is simple and the intersections $\gamma \cap U$ and $\gamma\cap V$ both consist of exactly one vertex and no edges.

\item Any path in $\pclass{U}{V}$ contains a $\pclass{U}{V}$-minimal subpath. That is, if $\gamma$ is a $\pathclass$-path from $U$ to $V$, then $\gamma$ contains a minimal $\pathclass$-path from $U$ to $V$ as a subpath.
\end{enumerate}
\end{lem}

%
%
%
%
%
%

\section{Convex Geometry in $\Rd$}
\label{convex_geometry_sec}


The following four lemmas are needed for the proof of Lemma~\ref{blocking_segments_lem}, which is the basis for the proof of Proposition~\ref{narrow_extinction_prop}.

\begin{lem}[Distance increases from a support point]
\label{dist_support_point_lem}
Let $\mu$ be a norm, and let $H$ be a supporting hyperplane of the unit ball $\unitball{\mu}$ at the point $\vect x\in \bd \unitball{\mu}$. Then for any $\dirsymb\in H -\vect x$, the function $f_{\dirsymb}(s) = \distnew{\mu}{\vect 0}{\vect x + s\dirsymb}$ is increasing on $[0,\infty)$.
\end{lem}

\begin{proof}
Observe that
\[
f_{\dirsymb}(s) = \inf \setbuilder[2]{t\ge 1}{\vect x + s\dirsymb \in t\unitball{\mu}},
\]
where the infimum is over $t\ge 1$ rather than $t\ge 0$ because the fact that $H$ is a supporting hyperplane of $\unitball{\mu}$ implies that the set on the right would be empty for $t\in [0,1)$.
Since $H$ is a supporting hyperplane at $\vect x$, we have $\vect x\in t\unitball{\mu}$ for all $t\ge 1$. Thus, since the ball $t\unitball{\mu}$ is convex, if $\vect x + s\dirsymb \in t\unitball{\mu}$ for some $t\ge 1$, then we have $\segment{\vect x}{\vect x + s\dirsymb} \subseteq t\unitball{\mu}$. Now observe that
\[
\segment{\vect x}{\vect x + s\dirsymb}
= \setbuilder[2]{\vect x + s' \dirsymb}{0\le s'\le s},
\]
so
\[
\vect x + s\dirsymb \in t\unitball{\mu}
\implies
\vect x + s'\dirsymb \in t\unitball{\mu}
\quad \forall s'\in [0, s],
\]
which shows that if $s\ge s'\ge 0$, then $f_{\dirsymb}(s)\ge f_{\dirsymb}(s')$.
\end{proof}

\begin{lem}[Support hyperplanes of bodies with connected interior]
\label{support_hplane_body_lem}
Let $B$ be a body in $\Rd$ such that $\interior{B}$ is connected. Then a hyperplane $\Hplane$ in $\Rd$ is a support hyperplane of $B$ if and only if $H\cap B\ne \emptyset$ and $H\cap \interior{B} = \emptyset$.
\end{lem}

\begin{proof}

\newcommand{\thisx}{\vect x}

\newcommand{\intB}{\interior{B}}
\newcommand{\suppspace}{\Hspace}
\newcommand{\oppspace}{\widehat{\Hspace}}

\newcommand{\intsuppspace}{\interior{\suppspace}}
\newcommand{\intoppspace}{\interior{\oppspace}}

First suppose $\Hplane$ is a support hyperplane of $B$. Then by definition we have $\Hplane\cap B\ne\emptyset$ and $\Hplane = \bd \Hspace$ for some support half-space $\Hspace$ of $B$. Then $B\subseteq\Hspace$, so $\intB\subseteq \interior{\Hspace}$, and hence
$
\Hplane\cap \intB \subseteq \Hplane\cap \interior{\Hspace} = \emptyset.
$
Thus the ``only if" direction is true for any $B\subseteq\Rd$, not just bodies with connected interior.

For the ``if" direction, let $\Hplane\subset\Rd$ be a hyperplane such that $H\cap B\ne \emptyset$ and $H\cap \interior{B} = \emptyset$. Choose any $\thisx\in \interior{B}$ (such an $\thisx$ exists because $B$ is a body). Then since $\thisx\notin\Hplane$, exactly one of the two half-spaces with boundary $\Hspace$ contains the point $\thisx$; let $\suppspace$ be the unique closed half-space with $\bd \suppspace = \Hplane$ and $\thisx\in\suppspace$, and let $\oppspace$ be the unique closed half-space with $\bd \oppspace = \Hplane$ and $\thisx\notin \oppspace$.
Then the sets $\intsuppspace$, $\Hplane$, and $\intoppspace$ are mutually disjoint, and $\Rd = \intsuppspace\cup \Hplane\cup \intoppspace$. Thus, since $\interior{B}\cap \Hplane = \emptyset$, we have $\interior{B} \subseteq \intsuppspace\cup \intoppspace$, and since $\intB$ is connected, this implies that either $\intB\subseteq \intsuppspace$ or $\intB\subseteq \intoppspace$. Thus we must have $\intB\subseteq\intsuppspace$ since $\intB\cap \intsuppspace \ne\emptyset$ by construction. Therefore, since $B$ and $\suppspace$ are both bodies, we have
\[
B\subseteq \closure{B} = \closure{\intB}
\subseteq \closure{\intsuppspace} = \suppspace.
\]
Since $\bd \suppspace = \Hplane$ has nonempty intersection with $B$ by assumption, this shows that $\suppspace$ is a support half-space of $B$ and hence $\Hplane$ is a support hyperplane of $B$.
\end{proof}

\begin{lem}[Support hyperplanes and closest points]
\label{support_point_closest_lem}
Suppose $\distnew{\norm{\cdot}}{\vect x}{\vect y} = r$. Then $H$ is a support hyperplane of $\reacheddset{\norm{\cdot}}{\vect x}{r}$ at $\vect y$ if and only if $\vect y$ is a $\norm{\cdot}$-closest point in $H$ to $\vect x$.
\end{lem}

\begin{proof}
The implication $\Longrightarrow$ follows from Lemma~\ref{dist_support_point_lem}. For the other direction, suppose that $\vect y$ is a $\norm{\cdot}$-closest point in $H$ to $\vect x$. Then $\vect y\in H$, and for any $\vect y'\in H$ we have $\distnew{\norm{\cdot}}{\vect y'}{\vect x} \ge \distnew{\norm{\cdot}}{\vect y}{\vect x} = r$, so $H\cap \interior{\reacheddset{\norm{\cdot}}{\vect x}{r}} = \emptyset$. Thus, $H$ is a support hyperplane of $\reacheddset{\norm{\cdot}}{\vect x}{r}$ by Lemma~\ref{support_hplane_body_lem}, since the ball is a convex body.
\end{proof}

\begin{lem}[Radial projection in convex sets]
\label{convex_ray_bd_point_lem}
If $C\subseteq\Rd$ is convex and $\vect x\in \interior{C}$, then any ray originating at $\vect x$ intersects $\bd C$ in at most one point.
If $C$ is compact, the intersection contains exactly one point. 
\end{lem}

\begin{proof}
The version of this for compact $C$ is the main claim in the proof of Proposition~4.26 on p.\ 80 of \cite{Lee:2000aa} (see Lemma~\ref{convex_euclidean_lem} below). The case when $C$ is noncompact can be handled similarly.
\end{proof}

\begin{lem}[Convex sets in Euclidean space]
\label{convex_euclidean_lem}
Suppose $C$ is a closed convex subset of $\Rd$, and let $k = \dim (\aff C)$. Then $C$ is homeomorphic to either the closed ball $\unitball{\lspace{2}{k}}$, the closed half-space $\closure{\halfspace[k]}$, or the Euclidean space $\R^k$.
\end{lem}

\begin{proof}
The proof when $C$ is compact and $k=d$ is given in Proposition~4.26 of \cite[pp.~80--81]{Lee:2000aa}.
The idea of the proof is to use Lemma~\ref{convex_ray_bd_point_lem} above to define a homeomorphism from $\bd C$ to the sphere (using the closed map lemma \cite[Lemma~4.25, p.~79]{Lee:2000aa}), and then extend the inverse of this map to a homeomorphism from the closed ball to $C$, again using the closed map lemma.
The more general version stated here can be proved similarly.
%
%
%
\end{proof}

The following lemma enumerates some simple properties of the convex hull. The proof is left to the reader.

\begin{lem}[Properties of the convex hull]
\mbox{} 
\label{convex_hull_properties_lem}
\begin{enumerate}
\item If $\mu$ is a norm on $\Rd$, then for any $\vect x,\vect y\in\Rd$ and $r\ge 0$,
\[
\conv \parens[2]{\reacheddset{\mu}{\vect x}{r} \cup \reacheddset{\mu}{\vect y}{r}}
= \segment{\vect x}{\vect y} + r\unitball{\mu}.
\]

\item For any $A\subseteq\Rd$, $\conv A = \conv \setof{\text{support points of $A$}}$.
\end{enumerate}
\end{lem}

\section{Metric Geometry}
\label{metric_geom_sec}

Let $(X,\distOp)$ be a metric space. A \textdef{path} in $X$ is a continuous function $\gamma\colon I\to X$ from an interval $I\subseteq \R$ to $X$. The \textdef{length} of a path $\gamma\colon [a,b]\to X$ with respect to $\distOp$ is
\begin{equation}
\label{metric_length_def_eqn}
\length{\distOp}{\gamma} \definedas \sup \sum_{i=1}^N
\distnew[2]{}{\gamma(t_{i-1})}{\gamma(t_i)},
\end{equation}
where the supremum is over all partitions $a=t_0<t_1<\dotsb<t_N =b$ of $I=[a,b]$. The curve $\gamma$ is called \textdef{rectifiable} with respect to the metric $\distOp$ if $\length{\distOp}{\gamma}<\infty$ (see \cite[p.~34]{Burago:2001aa}).
The length operator defined by \eqref{metric_length_def_eqn} satisfies the \textdef{generalized triangle inequality} (cf.\ \cite[p.~35]{Burago:2001aa}):
\begin{equation}
\label{gen_triangle_ineq_eqn}
\text{If $\gamma\colon [a,b]\to X$ is a path, then }
\length{\distOp}{\gamma} \ge \distnew[2]{}{\gamma(a)}{\gamma(b)}.
\end{equation}
The length operator $\lengthOp_{\distOp}$ can be used to define a new metric on $X$:
\begin{equation}
\label{gen_intrinsic_metric_def_eqn}
\dist{\lengthOp} (\vect x,\vect y) \definedas 
\inf \setbuilder[3]{\length{\distOp}{\gamma}}{\gamma\colon [a,b] \to X
	\text{ is a path with $\gamma(a) = \vect x$, $\gamma(b) = \vect y$}}.
\end{equation}
The metric defined in \eqref{gen_intrinsic_metric_def_eqn} is called the \textdef{intrinsic metric induced by $\distOp$}. If it happens that $\dist{\lengthOp} = \distOp$, then the metric $\distOp$ is called \textdef{intrinsic}, and the metric space $(X,\distOp)$ is called a \textdef{length space}.
It can be shown that the metric $\dist{\lengthOp}$ itself is intrinsic according to this definition (see \cite[pp.~37--39]{Burago:2001aa}).
The metric $\dist{\lengthOp}$ is called \textdef{strictly intrinsic} if any two points in $X$ can be joined by a path that achieves the infimum in \eqref{gen_intrinsic_metric_def_eqn}, i.e.\ if for any $\vect x,\vect y\in X$ there exists a path $\gamma$ from $\vect x$ to $\vect y$ such that $\distnew{\lengthOp}{\vect x}{\vect y} = \length{\distOp}{\gamma}$; in this case $\gamma$ is called a \textdef{shortest path} or \textdef{minimizing geodesic} from $\vect x$ to $\vect y$.

\begin{lem}[Line segments minimize distance in any norm metric]
\label{norm_strictly_intrinsic_lem}
If $\mu$ is any norm on $\Rd$, then $\normmetric{\mu}$ is strictly intrinsic, and any line segment is a shortest path.
\end{lem}

\begin{proof}
Write $\lengthOp_{\mu}$ for the length operator induced by $\normmetric{\mu}$ via the definition \eqref{metric_length_def_eqn}. Let $\gamma\colon [a,b]\to\Rd$ be a simple path whose image is the line segment from $\vect x$ to $\vect y$. Then for any partition $a=t_0<t_1<\dotsb<t_N =b$ of $[a,b]$, the points $\vect x = \gamma(t_0), \gamma(t_1),\dotsc,\gamma(t_N) =\vect y$ are collinear and appear in the listed order, so by the definition \eqref{metric_length_def_eqn} we have
\[
\length{\mu}{\gamma}  =
\sup \sum_{i=1}^N \distnew[2]{\mu}{\gamma(t_{i-1})}{\gamma(t_i)}
= \sup \sum_{i=1}^N \mu\parens[2]{\gamma(t_i) - \gamma(t_{i-1})}
= \mu(\vect y-\vect x),
\]
where the final equality follows from the collinearity of the points $\vect x = \gamma(t_0), \gamma(t_1),\dotsc,\gamma(t_N) =\vect y$. Thus, $\gamma$ is a shortest path from $\vect x$ to $\vect y$. Since $\vect x$ and $\vect y$ were arbitrary, this shows that $\normmetric{\mu}$ is strictly intrinsic.
\end{proof}

%

\section{Probability and Analysis}

\begin{lem}
\label{sufficient_moment_cond_lem}
\Lspace{1/2} implies \minimalmoment{d} for all $d\ge 1$.
\end{lem}

\begin{proof}

Suppose $Y$ is a nonnegative random variable satisfying \Lspace{1/2}, and let $M_{2d} \definedas \min \setof{Y_1,\dotsc,Y_{2d}}$, where $d$ is a positive integer, and $Y_1,\dotsc,Y_{2d}$ are \iid\ random variables with $Y_j \eqd Y$ for $1\le j\le 2d$. By definition, $Y$ satisfies \minimalmoment{d} if and only if $M_{2d}\in L^{d}$.
By Exercise 38 on p.~199 of \cite[Section~6.4]{Folland:1999aa}, for any $p>0$ and any nonnegative random variable $X$ we have
\begin{equation}
\label{sufficient_moment_cond:Lp_sum_eqn}
X\in L^p \iff \sum_{k=0}^\infty 2^{kp}\, \Pr\eventthat[2]{X>2^{k}}< \infty.
\end{equation}
By \eqref{sufficient_moment_cond:Lp_sum_eqn}, since $Y\in L^{1/2}$ we have
\begin{equation}
\label{sufficient_moment_cond:Y_sum_eqn}
\sum_{k=0}^\infty 2^{k/2}\, \Pr\eventthat[2]{Y>2^k}< \infty,
\end{equation}
and \eqref{sufficient_moment_cond:Y_sum_eqn} implies that we can choose $k_0$ large enough that
\[
2^{k/2} \Pr \eventthat[2]{Y > 2^{k}}<1
\quad \text{for all $k\ge k_0$.}
\]
Then
\[
2^{k d} \Pr \eventthat[2]{Y > 2^k}^{2d}
= \parens[1]{2^{k/2} \Pr \eventthat[2]{Y > 2^k}}^{2d}
<2^{k/2} \Pr \eventthat[2]{Y > 2^k}
\quad \text{for all $k\ge k_0$.}
\]
Therefore we have
\begin{align*}
\sum_{k=k_0}^\infty 2^{k d}\, \Pr\eventthat[2]{M_{2d}>2^k}
&= \sum_{k=k_0}^\infty 2^{k d} \Pr \eventthat[2]{Y > 2^k}^{2d}\\
&< \sum_{k=k_0}^\infty 2^{k/2}\, \Pr\eventthat[2]{Y>2^k}
<\infty,
\end{align*}
where the convergence of the sum follows from \eqref{sufficient_moment_cond:Y_sum_eqn}. Thus we have $M_{2d}\in L^{d}$ by \eqref{sufficient_moment_cond:Lp_sum_eqn}, so $Y$ satisfies \minimalmoment{d}.
\end{proof}

The following lemma is a basic large deviations estimate which is the basis of many of the results in Chapter~\ref{cone_growth_chap}, starting with Lemma~\ref{tmeas_dev_lem}.

\begin{lem}[Large deviations estimate for \iid\ sums]
\label{basic_large_dev_lem}
\newcommand{\thisrv}{X}
\newcommand{\thislaw}{\law (\thisrv)}
\newcommand{\expconst}{b}

\newcommand{\thisC}{\bigcref{basic_large_dev_lem}}
\newcommand{\thisCp}[2][0]{\thisC \parens[#1]{#2}}
\newcommand{\thisc}{\smallcref{basic_large_dev_lem}}
\newcommand{\thiscp}[2][0]{\thisc \parens[#1]{#2}}

Let $\thisrv,\thisrv_1,\thisrv_2,\dotsc$ be \iid\ random variables such that $\E e^{b \thisrv} <\infty$ for some $b>0$. Then $\E\thisrv \in [-\infty,+\infty)$, and for any $x> \E\thisrv$, there are positive constants $\thisCp{x}$ and $\thiscp{x}$ (which depend on $\thislaw$) such that for any $n\in\N$,
\[
\Pr \eventthat[1]{\sum_{j=1}^n \thisrv_j \ge nx}
\le \thisCp{x} e^{-\thiscp{x} n}.
\]
\end{lem}

\begin{proof}
This follows from basic large deviations theory. See for example Corollary~27.4 on p.~541 of \cite{Kallenberg:2002aa}.
\end{proof}

The following lemma is cited in the proof of Theorem~\ref{irrelevance_thm}.

{
\newcommand{\sigfield}{\mathcal{F}}
\begin{lem}
\label{zero_cond_prob_lem}
Let $\probabilityspace$ be a probability space, and let $A\in \sigmafield$ be an event. Then the following are equivalent.
\begin{enumerate}
\item $\Pr(A)=0$.
\item $\cprob{A}{\sigfield} =0$ $\Pr$-a.s.\ for some $\sigma$-field $\sigfield\subseteq \sigmafield$.
\item $\cprob{A}{\sigfield} =0$ $\Pr$-a.s.\ for all $\sigma$-fields $\sigfield\subseteq \sigmafield$.
\end{enumerate}
\end{lem}

\begin{proof}
For any $A\in \sigmafield$ and any $\sigma$-field $\sigfield$, $\Pr(A) = \E \cprob{A}{\sigfield} = 0$ if and only if $\cprob{A}{\sigfield} = 0$ a.s.
\end{proof}
}

The following bound is used in several proofs in Chapters~\ref{cone_growth_chap} and \ref{random_fpc_chap}.

\begin{lem}[A bound for ``polynomial$\times$geometric" series]
\label{poly_geo_series_lem}
For any $n\in\N$, $c>0$, and $a,m\ge 0$,
\[
\sum_{j=n}^\infty (j+a)^m e^{-cj}
\le \frac{4}{(c\meet 1)^{m+1}}\cdot  \parens{2m \join a}^m \cdot e^{-\frac{c}{2} n},
\]
with the convention $0^0=1$ when $m=0$.
\end{lem}

\begin{proof}
\newcommand{\xmax}{x_{\max}}



For any $a,m\ge 0$ and $b>0$, consider the function $f\colon [-a,\infty)\to \Rplus$ defined by $f(x) = (x+a)^m e^{-bx}$, with the convention $0^0=1$ when $m=0$. Computing $f'(x)$ shows that $f$ attains a global maximum of $f(\xmax) = \parens[2]{\frac{m}{be}}^m e^{ba}$ at $\xmax = \frac{m}{b} -a$, and $f$ is decreasing on $[\xmax,\infty)$. Therefore, if $\xmax\le 0$, then $\restrict{f}{[0,\infty)}$ attains a maximum of $f(0) = a^m$ at $x=0$, so we have
\begin{equation}
\label{poly_geo_series:sup_cases_eqn}
\sup_{x\ge 0}\; (x+a)^m e^{-bx}=
\begin{cases}
\parens[1]{\frac{m}{be}}^m e^{ba} &\text{if } ab\le m,\\
a^m &\text{if } ab\ge m.
\end{cases}
\end{equation}
Note that if $ab\le m$, then $\parens[2]{\frac{m}{be}}^m e^{ba} \le \parens[2]{\frac{m}{be}}^m e^{m} = \parens[1]{\frac{m}{b}}^m$, so \eqref{poly_geo_series:sup_cases_eqn} implies that
\begin{equation}
\label{poly_geo_series:sup_bound_eqn}
\sup_{x\ge 0}\; (x+a)^m e^{-bx} \le
\parens[1]{\frac{m}{b} \join a}^m.
\end{equation}
Using the bound \eqref{poly_geo_series:sup_bound_eqn} with $b=c/2$ and $x=j$, for any $n\in\N$ we get
\begin{align*}
\sum_{j=n}^\infty (j+a)^m e^{-cj}
&= \sum_{j=n}^\infty (j+a)^m e^{-\frac{c}{2}j} e^{-\frac{c}{2}j}\\
&\le \parens[1]{\frac{2m}{c} \join a}^m \cdot \sum_{j=n}^\infty  e^{-\frac{c}{2}j}\\
&= \parens[1]{\frac{2m}{c} \join a}^m \cdot \frac{e^{-\frac{c}{2}n}}{1-e^{-\frac{c}{2}}}\\
&\le \parens[1]{\frac{2m}{c\meet 1} \join a}^m
	\cdot \frac{4}{c\meet 1}\cdot e^{-\frac{c}{2}n}\\
&\le \frac{4}{(c\meet 1)^{m+1}}\cdot \parens{2m \join a}^m
	\cdot e^{-\frac{c}{2}n},
\end{align*}
where in the second-to-last line we have used the bound
\[
\frac{1}{1-e^{-\frac{c}{2}}} = \frac{e^{\frac{c}{2}}}{e^{\frac{c}{2}} - 1}
\le \frac{2}{\parens[2]{e^{\frac{c}{2}} - 1}\meet 1} \le \frac{2}{\parens[1]{\frac{c}{2}}\meet 1}
\le \frac{4}{c\meet 1}. \qedhere
\]
\end{proof}


\chapter{Asymptotic Equivalence and Growth of Territory}
\label{asymptotic_chap}

In this appendix we prove several elementary results that are applied to first-passage percolation in Chapter~\ref{cone_growth_chap}, but in a more abstract setting that reduces the notational distractions that arise when considering the restricted first-passage percolation process.


\section{Sublevel Sets and Inversion Formulas}

Let $X$ be a nonempty set. For a function $f:X\to [0,\infty)$ and $t\ge 0$, we let $\reachedset{f}{}{t}$ denote the \textdef{sublevel set} of $f$ up to the value $t$. That is,
\[
\reached_f(t) \definedas f^{-1}[0,t] = \{x\in X: f(x)\le t\}.
\]
Clearly the function $\reachedfcn{f}{}: [0,\infty)\to 2^X$ is increasing with respect to set inclusion. Thus, if we interpret $t$ as a time parameter, we can think of $\reachedfcn{f}{}$ as a ``pure growth process" in $X$, where $\reachedset{f}{}{t}$ represents the set of points in $X$ that the process has reached by time $t$. Observe that the function $\reachedfcn{f}{}$ is also right-continuous,
in that $\reached_f(t) = \bigcap_{s>t} \reached_f(s)$. Moreover, the left-hand limit $\reachedset{f}{}{t-}$ exists at each $t$, and the function $t\mapsto \reachedset{f}{}{t-}$ gives a left-continuous version of the growth process:
\[
\reached_f(t-) \definedas \bigcup_{s<t} \reached_f(s) =  \{x\in X: f(x)< t\} = f^{-1}[0,t).
\]
We now give some general formulas for converting a particular statement about two functions $f$ and $g$ into a statement relating their sublevel sets. We refer to these as ``inversion formulas" because they allow us to easily go back and forth between the functions $f$ and $g$ and their inverse images $\reachedfcn{f}{}$ and $\reachedfcn{g}{}$.

\begin{lem}[Inversion formulas]
\label{general_inversion_lem}
Let $X$ be a nonempty set, and let $f,g:X\to [0,\infty)$. Then for any $\alpha>0$ and $t_0\ge 0$,
\begin{enumerate}
\item
$\displaystyle
\braces[3]{\forall x\in X,\; f(x) < t_0 \text{ or } f(x)<\alpha^{-1} g(x)}
\Longleftrightarrow
\braces[3]{\forall t\ge t_0,\; \reachedset{g}{}{\alpha t} \subseteq \reachedset{f}{}{t-}}.
$

\item
$\displaystyle
\braces[3]{ \forall x\in X,\; f(x) < t_0 \text{ or } f(x)\le \alpha^{-1} g(x)}
\Longleftrightarrow
\braces[3]{ \forall t\ge t_0,\; \reachedset{g}{}{\alpha t-} \subseteq \reachedset{f}{}{t-} }.
$

\item
$\displaystyle
\braces[3]{\forall x\in X,\; f(x) \le t_0 \text{ or } f(x)<\alpha^{-1} g(x)}
\Longleftrightarrow
\braces[3]{\forall t> t_0,\; \reachedset{g}{}{\alpha t} \subseteq \reachedset{f}{}{t-}}.
$

\item
$\displaystyle
\braces[3]{\forall x\in X,\; f(x) \le t_0 \text{ or } f(x)\le \alpha^{-1} g(x)} 
\Longleftrightarrow
\begin{cases}
\forall t\ge t_0,\; \reachedset{g}{}{\alpha t} \subseteq \reachedset{f}{}{t} & or\\
\forall t\ge t_0,\; \reachedset{g}{}{\alpha t-} \subseteq \reachedset{f}{}{t} & or\\
\forall t> t_0,\; \reachedset{g}{}{\alpha t} \subseteq \reachedset{f}{}{t} & or\\
\forall t> t_0,\; \reachedset{g}{}{\alpha t-} \subseteq \reachedset{f}{}{t-} & or\\
\forall t> t_0,\; \reachedset{g}{}{\alpha t-} \subseteq \reachedset{f}{}{t}.
\end{cases}
$
\end{enumerate}
\end{lem}

\begin{proof}
The proof is an elementary exercise and is left to the reader.
\end{proof}

\begin{thmremark}
In Lemma~\ref{general_inversion_lem}, the four statements on the left can be viewed as points in the Boolean lattice $\setof{\le,<}^2\cong \setof{0,1}^2$ if we partially order the statements by their strength (i.e.\ statement 1 is the strongest, statement 4 is the weakest, and statements 2 and 3 are intermediate). Similarly, the eight statements on the right correspond to points in the lattice $\setof{>,\ge}^2 \times \setof{\le,<} \cong \setof{0,1}^3$, and the lemma posits a lattice homomorphism $\setof{0,1}^3 \to \setof{0,1}^2$ in which each statement is logically equivalent to its image.
\end{thmremark}

\begin{lem}[More equivalent statements]
\label{equiv_statements_lem}
Let $X$ be a nonempty set, and let $f,g:X\to [0,\infty)$. Then for any $\alpha>0$ and $t_0\ge 0$, the following equivalencies hold, where each statement is implicitly prefaced by ``$\forall x\in X$," and the braces separating the statements on the right indicate disjunction.
\begin{enumerate}
\item
$\displaystyle
{f(x) < t_0 \text{ or } f(x)<\alpha^{-1} g(x)}
\Longleftrightarrow
\begin{cases}
\text{If $f(x)\ge t_0$ or $g(x)\ge \alpha t_0$, then } f(x)< \alpha^{-1}g(x).\\
\text{If $f(x)\ge t_0$ or $g(x)> \alpha t_0$, then } f(x)< \alpha^{-1}g(x).
\end{cases}
$

\item
$\displaystyle
{f(x) < t_0 \text{ or } f(x)\le \alpha^{-1} g(x)}
\Longleftrightarrow
\begin{cases}
\text{If $f(x)\ge t_0$ or $g(x)\ge \alpha t_0$, then } f(x)\le \alpha^{-1}g(x).\\
\text{If $f(x)\ge t_0$ or $g(x)> \alpha t_0$, then } f(x)\le \alpha^{-1}g(x).
\end{cases}
$

\item
$\displaystyle
{f(x) \le t_0 \text{ or } f(x)< \alpha^{-1} g(x)}
\Longleftrightarrow
\begin{cases}
\displaystyle
\parens[1]{\textstack{If $f(x)> t_0$ or $g(x)\ge \alpha t_0$, then either\qquad}
	{\quad $f(x)< \alpha^{-1}g(x)$ or $f(x)=t_0 = \alpha^{-1}g(x)$.}}\\
\text{If $f(x)> t_0$ or $g(x)> \alpha t_0$, then } f(x)< \alpha^{-1}g(x).
\end{cases}
$

\item
$\displaystyle
{f(x) \le t_0 \text{ or } f(x)\le \alpha^{-1} g(x)}
\Longleftrightarrow
\begin{cases}
\text{If $f(x)> t_0$ or $g(x)\ge \alpha t_0$, then } f(x)\le \alpha^{-1}g(x).\\
\text{If $f(x)> t_0$ or $g(x)> \alpha t_0$, then } f(x)\le \alpha^{-1}g(x).
\end{cases}
$
\end{enumerate}
\end{lem}

\begin{proof}
The implications ``$\Longleftarrow$" are all trivial, and it is elementary to prove ``$\Longrightarrow $" in each case.
\end{proof}

\begin{thmremark}
As in Lemma~\ref{general_inversion_lem}, the equivalencies in Lemma~\ref{equiv_statements_lem} can be interpreted as a map between the Boolean lattices $\setof{0,1}^3$ and $\setof{0,1}^2$ in which the image of a statement in $\setof{0,1}^3$ is an equivalent statement in $\setof{0,1}^2$. The anomalous statement in Part~3 is strictly weaker than the more obvious ``\textit{If $f(x)> t_0$ or $g(x)\ge \alpha t_0$, then $f(x)< \alpha^{-1}g(x)$}," which is intermediate between Statements 1 and 3.
\end{thmremark}

\section{Limits and Asymptotic Equivalence}

We introduce a simple definition of ``co-limiting behavior" for functions.

\begin{prop}[Characterizations of co-limiting behavior]
\label{colimit_char_prop}
Let $X$ be a nonempty set and $Y$ a topological space, with $y\in Y$. For a given pair of functions $g:X\to Y$ and $h:X\to [0,\infty)$, the following are equivalent.
\begin{enumerate}
\item For every neighborhood $U$ of $y$, there exists $M<\infty$ such that $g(x)\in U$ whenever $h(x)>M$ (equivalently, $\forall$ nbhd $U$ of $y$, $\exists M$ such that $g(x)\in Y\setminus U \Rightarrow h(x)\le M$).
\item If $\{x_n\}_{n\in\N}$ is any sequence of points in $X$ such that $h(x_n)\to\infty$, then $g(x_n)\to y$.
\end{enumerate}
\end{prop}

\begin{proof}

\underline{1 implies 2:} Let  $\{x_n\}_{n\in\N}\subseteq X$ with $h(x_n)\to\infty$, and let $U\subseteq Y$ be a neighborhood of $y$. By 1, there exists $M$ so that $g(x)\in U$ whenever $h(x)>M$. Since $h(x_n)\to\infty$, we have $h(x_n)>M$ eventually, so $g(x_n)\in U$ eventually. Thus $g(x_n)\to y$.

\underline{2 implies 1:} We will prove this by contradiction. Suppose 2 holds but 1 does not. Then there exists a neighborhood $U$ of $y$ for which there is no $M<\infty$ such that $h(x)\le M$ for all $x\in g^{-1} (Y\setminus U) = X\setminus g^{-1}U$. Therefore, for every $n\in\N$, there exists $x_n\in g^{-1}(Y\setminus U)$ (i.e.\ $x_n\in X$ with $g(x_n)\not\in U$) such that $h(x_n)>n$. Thus $h(x_n)\to\infty$, so $g(x_n)\to y$ by assumption.  But then $g(x_n)\in U$ eventually, contradicting our choice of $x_n$.
\end{proof}

If $g$ and $h$ satisfy either of the equivalent conditions in Proposition~\ref{colimit_char_prop}, we say that ``$g(x)\to y$ as $h(x)\to \infty$", ``$g(x)\to y$ if $h(x)\to \infty$", or ``$h(x)\to \infty$ implies $g(x)\to y$", and we write
\begin{equation}
\label{lim_g_as_h_def_eqn}
\lim_{h(x)\to \infty} g(x) = y.
\end{equation}
We can also say that ``$g(x)\to y$ \emph{uniformly in $x$ with respect to $h$} as $h(x)\to \infty$," since the speed of convergence to the limit doesn't depend directly on $x$ but only on the value of $h(x)$. This construction provides a convenient way of describing asymptotics of functions defined on an arbitrary space $X$ by choosing $h$ to be some function that approaches $\infty$ as $x\in X$ approaches the relevant values of the domain.

\begin{remark}
The definition \eqref{lim_g_as_h_def_eqn} can be generalized to allow $h: X\to Z$, where $Z$ is an arbitrary topological space instead of $[0,\infty]$. In this case we would write $\lim_{h(x)\to z} g(x) = y$ to mean $g(x)\to y\in Y$ as $h(x)\to z\in Z$. If $Z$ is not first-countable, then the second condition in Proposition~\ref{colimit_char_prop} should be replaced with the analogous statement for nets.
\end{remark}

Let $X$ be a nonempty set and let $h:X\to [0,\infty)$ be some ``reference function" on $X$.
We say that two complex-valued functions $f$ and $g$ on $X$ are \textdef{asymptotically equivalent} as $h\to \infty$, and write $f\sim g$ as $h\to \infty$ if
\begin{equation}
\label{asymp_equiv_def_eqn}
\frac{f(x)}{g(x)}\to 1 \text{ as } h(x)\to \infty,\ x\in X.
\end{equation}
Note that for a given $h$, $\sim$ is an equivalence relation. Now suppose that $f,g:X\to [0,\infty)$, and take $h = f\vee g:X\to [0,\infty)$. We say that $f$ and $g$ are \textdef{asymptotic at infinity}, and write $f\sim g$ at $\infty$, if
\begin{equation}
\label{asymp_at_infty_def_eqn}
\frac{f(x)}{g(x)}\to 1 \text{ as } (f\vee g)(x)\to \infty.
\end{equation}
Note that ``at $\infty$" refers to the range of $f$ or $g$, not the domain. Note also that $\sim$ defines an equivalence relation on $\R_+$-valued functions on $X$ (and this time there's no need to fix a reference function $h$ first). Using the inversion formulas (Lemmas~\ref{general_inversion_lem} and \ref{equiv_statements_lem}), we now give several equivalent characterizations for two functions $f$ and $g$ being asymptotic at $\infty$, in terms of their sublevel sets $\reachedfcn{f}{}$ and $\reachedfcn{g}{}$.

\begin{prop}[Asymptotic equivalence at $\infty$]
\label{asymptotic_equivalence_prop}
Let $X$ be a nonempty set, and let $f,g:X\to [0,\infty)$. The following are equivalent.
\begin{enumerate}
\item $f\sim g$ at $\infty$, as defined in \eqref{asymp_at_infty_def_eqn}.


\item For all $\epsilon>0$, there exists $t_0<\infty$ such that
\[
\reachedset[2]{g}{}{(1-\epsilon)t}
\subseteq \reachedset{f}{}{t}\subseteq
\reachedset[2]{g}{}{(1+\epsilon)t}
\quad\text{for all } t\ge t_0.
\]

\item There exists a nonnegative function $k$ on $[0,\infty)$ with $k(t) = o(t)$ such that
\[
\reachedset[2]{g}{}{t-k(t)}
\subseteq \reachedset{f}{}{t}\subseteq
\reachedset[2]{g}{}{t+k(t)}
\quad\text{for all } t\ge 0.
\]

\item There is some nonnegative $o(t)$ function $h(t)$ on $[0,\infty)$ such that for all $\epsilon>0$, there exists $t_0<\infty$ such that for all $t > t_0$,
\[
\reachedset[2]{g}{}{[(1-\epsilon)t - h(t)]-}
\subseteq \reachedset{f}{}{t}
\quad\text{and}\quad
\reachedset{f}{}{t-}\subseteq
\reachedset[2]{g}{}{(1+\epsilon)t+h(t)}.
\]

\item For any $\epsilon>0$ and any $o(t)$ function $h(t)$ on $[0,\infty)$, there exists $t_0<\infty$ such that for all $t\ge t_0$,
\[
\reachedset[2]{g}{}{(1-\epsilon)t + h(t)}
\subseteq \reachedset{f}{}{t-}
\quad\text{and}\quad
\reachedset{f}{}{t}\subseteq
\reachedset[2]{g}{}{[(1+\epsilon)t-h(t)]-}.
\]

\end{enumerate}
Because of the symmetry of Statement 1, the roles of $f$ and $g$ can be switched in any of the other statements.
\end{prop}

\begin{proof}

The equivalence $1\iff 2$ follows easily from Lemmas~\ref{general_inversion_lem} and \ref{equiv_statements_lem}, and the implications $5 \implies 2\implies 4$ are trivial. To complete the proof, we will show that $4\implies 3\implies 5$.

($4\implies 3$) Suppose Statement 4 holds. By the equivalence of the various processes in Part~4 of Lemma~\ref{general_inversion_lem}, we can replace the strict inequality on $t_0(\epsilon)$ with a weak one, and we can replace the left-continuous processes with right-continuous ones. Thus, suppose there exists a nonnegative function $h(t) = o(t)$ such that for every $\epsilon>0$, there exists $t(\epsilon)<\infty$ such that if $t\ge t(\epsilon)$ then
\[
\reached_g \bigl((1-\epsilon)t-h(t)\bigr)
\subseteq \reached_f(t)\subseteq
\reached_g \bigl((1+\epsilon)t+h(t)\bigr).
\]
For each positive integer $n$, let $\epsilon_n = 1/n$. Set $t_0=0$ and for each $n\ge 1$ let $t_n = t({\epsilon_n})\vee t_{n-1}$. Let $t_\infty = \lim_{n\to\infty} t_n = \sup_{n\ge 0} t_n \le \infty$, and let $N = \sup \{n: t_n = t_{n-1}\} \le \infty$. Now let $k(t)=(t+1)\vee \sup \{ g(x)-t : f(x)\le t\}$ on the interval $[0,t_1)$ (we can prove by contradiction that this $\sup$ is finite for each $t<t_1 = t(1)$), let $k(t) = \epsilon_n t +h(t)$ on the interval $[t_{n}, t_{n+1})$ for $1\le n< N$, and let $k(t)= h(t)$ on $[t_N,\infty)$.

($3\implies 5$) Let $k(t)$ be as in Statement~3. Given $\epsilon>0$ and $h(t) = o(t)$, choose $t_0$ large enough that $k(t)\le \epsilon t - h(t)$ for all $t\ge t_0$ (which is possible since $k(t)$ and $h(t)$ are $o(t)$). Then for $t\ge t_0$ we get $t-\epsilon t +h(t)\le t-k(t)$ and $t+k(t)\le t+\epsilon t -h(t)$, so
\[
\reached_g \bigl((1-\epsilon)t+h(t)\bigr)\subseteq
\reached_g \bigl(t-k(t)\bigr)
\subseteq \reached_f(t)\subseteq
\reached_g \bigl(t+k(t)\bigr)
\subseteq \reached_g \bigl((1+\epsilon)t-h(t)\bigr).
\]
This gives the corresponding statement with all the processes right-continuous. We can get the stated version involving left-continuous processes by replacing $h(t)$ with $1+h(t)\vee h(t+1)$, which is still $o(t)$, and then replacing the $t_0$ we obtain with $t_0+1$.
\end{proof}

\begin{thmremark}
The equivalent statements in Proposition~\ref{asymptotic_equivalence_prop} generalize the event that appears in the Shape Theorem (Theorem~\ref{shape_thm} or \ref{shape_thm1}). That is, the Shape Theorem says precisely that if $\tmeasure = \setof{\ttime(\edge)}_{\edge\in\edges{\Zd}}$ is an \iid\ collection of traversal times satisfying \finitespeed{d} and \minimalmoment{d}, then there is a norm $\mu$ on $\Rd$ such that the functions $f= \ptimenew{\tmeasure}{\vect 0}{\cdot}$ and $g = \distnew{\mu}{\vect 0}{\cdot}$ on $\Rd$ are almost surely asymptotic at $\infty$.
\end{thmremark}


The following proposition shows that if the functions $f$ and $g$ define random growth processes in $X$ such that we have an exponentially good large deviations estimate for the corresponding Shape Theorem holding at a fixed time $t$,
then we automatically get an exponentially good large deviations estimate for the event that the Shape Theorem takes hold for good, i.e.\ that the containment in Part 2 of Proposition~\ref{asymptotic_equivalence_prop} holds for all $t\ge t_0$.

\begin{prop}
\label{exponential_trapping_time_prop}
Let $X$ be a nonempty set, and let $f,g$ be random elements of $[0,\infty)^X$ in some probability space $\probabilityspace$. The following are equivalent.
\begin{enumerate}
\item For all $\epsilon>0$, there exist positive constants $\bigc{}$ and $\smallc{}$ such that for all $t\ge 0$,
\[
\Pr \eventthat[3]{
\reachedset[2]{g}{}{(1-\epsilon)t}
\subseteq \reachedset{f}{}{t}\subseteq
\reachedset[2]{g}{}{(1+\epsilon)t} }
\ge 1- \bigc{}e^{-\smallc{} t}.
\]

\item For all $\epsilon>0$, there exist positive constants $\bigc{}'$ and $\smallc{}'$ such that for any $t_0\ge 0$,
\[
\Pr \eventthat[3]{
\reachedset[2]{g}{}{(1-\epsilon)t}
\subseteq \reachedset{f}{}{t}\subseteq
\reachedset[2]{g}{}{(1+\epsilon)t}
\text{ for all } t\ge t_0}
\ge 1- \bigc{}'e^{-\smallc{}' t_0}.
\]
\end{enumerate}
\end{prop}

\begin{proof}
This follows from an easy adaptation of the proof of Lemma~\ref{trapping_time_lem} below, which is a special case.
\end{proof}

The following proposition gives various equivalent formulations of the probability bound appearing in Part~2 of Proposition~\ref{exponential_trapping_time_prop}, similar to the equivalent statements in Proposition~\ref{asymptotic_equivalence_prop}.

\begin{prop}
\label{equivalent_trapping_time_prop}
Let $X$ be a nonempty set, and let $f,g$ be random elements of $[0,\infty)^X$ in some probability space $\probabilityspace$. The following are equivalent. 
\begin{enumerate}


\item For all $\epsilon>0$, there exist positive constants $\bigc{}$ and $\smallc{}$ such that for all $t\ge t_0$,
\[
\Pr \eventthat[1]{
\abs[1]{\frac{f(x)}{g(x)}-1} < \epsilon
\text{ for all $x\in X$ with } (f\vee g)(x)\ge t_0}
\ge 1- \bigc{}e^{-\smallc{} t_0}.
\]

\item For all $\epsilon>0$, there exist positive constants $\bigc{}$ and $\smallc{}$ such that for any $t_0\ge 0$,
\[
\Pr \eventthat[3]{
\reachedset[2]{g}{}{(1-\epsilon)t}
\subseteq \reachedset{f}{}{t}\subseteq
\reachedset[2]{g}{}{(1+\epsilon)t}
\text{ for all } t\ge t_0}
\ge 1- \bigc{}e^{-\smallc{} t_0}.
\]


\item There is some $o(t)$ function $h:\R_+\to\R_+$ such that for all $\epsilon>0$, there exist positive constants $\bigc{}$ and $\smallc{}$ such that for all $t\ge 0$,
\[
\Pr \eventthat[3]{
\reachedset[2]{g}{}{[(1-\epsilon)t - h(t)]-}
\subseteq \reachedset{f}{}{t}
\text{ and }
\reachedset{f}{}{t-}\subseteq
\reachedset[2]{g}{}{(1+\epsilon)t+h(t)} }
\ge 1- \bigc{}e^{-\smallc{} t}.
\]

\item For any $\epsilon>0$ and any $o(t)$ function $h:\R_+\to\R_+$, there exist positive constants $\bigc{}$ and $\smallc{}$ such that for all $t_0\ge 0$,
\begin{align*}
\Pr \eventthat[3]{
\reachedset[2]{g}{}{(1-\epsilon)t + h(t)}
\subseteq \reachedset{f}{}{t-}
\text{ and }
\reachedset{f}{}{t}\subseteq
\reachedset[2]{g}{}{[(1+\epsilon)t-h(t)]-}
\text{ for all } t\ge t_0 }\quad \\
\ge 1- \bigc{}e^{-\smallc{} t_0}.
\end{align*}

\end{enumerate}
\end{prop}

\begin{proof}
Use Lemmas~\ref{general_inversion_lem} and \ref{equiv_statements_lem}, and an argument similar to the proof of Proposition~\ref{asymptotic_equivalence_prop}.
\end{proof}

\section{Applications to Restricted First-Passage Percolation}

%

The next two results strengthen the large deviations estimate in Lemma~\ref{GM_large_dev_lem}.

\begin{lem}\label{stronger_dev_lem}
Let $\tmeasure$ be an \iid\ traversal measure on $\edges{\Zd}$ satisfying \finitespeed{d} and \expmoment, and let $\mu$ be the shape function for $\tmeasure$ from Theorem~\ref{shape_thm}.
For any $\epsilon>0$, there exist two positive constants $C_1$ and $c_1$ such that for any $t>0$,
\begin{enumerate}
\item
$\displaystyle
	\Pr
	\left\{\reachedsymb_{\tmeasure}^{\vect 0;\Rd}(t) \subseteq (1+\epsilon)t \interior{\limitshape} \right\}
	\ge 1 - C_1 e^{-c_1 t},
$
and

\item
$\displaystyle
	\Pr
	\left\{(1-\epsilon)t \limitshape \subseteq \reachedsymb_{\tmeasure}^{\vect 0;\Rd}(t-) \right\}
	\ge 1 - C_1 e^{-c_1 t}.
$
\end{enumerate}
\end{lem}

\begin{proof}
Assume without loss of generality that $\epsilon<1$, and let $C_1(\epsilon)=C(\epsilon/2)$ and $c_1(\epsilon) = c(\epsilon/2) \cdot (1-\epsilon/2)$, where $C$ and $c$ are the $\epsilon$-dependent constants from Lemma~\ref{GM_large_dev_lem}.
For part 1, since $(1+\epsilon/2)t\limitshape \subseteq (1+\epsilon)t \interior{\limitshape}$, we have
\begin{align*}
	\Pr \left\{\reachedsymb_{\tmeasure}^{\vect 0;\Rd}(t) \not\subseteq (1+\epsilon)t \interior{\limitshape} \right\}
	&\le \Pr \left\{\reachedsymb_{\tmeasure}^{\vect 0;\Rd}(t) \not\subseteq (1+\epsilon/2)t \limitshape \right\} \\
	&\le C(\epsilon/2)\cdot \exp \left[-c(\epsilon/2)t\right]
	< C_1(\epsilon)\cdot \exp[-c_1 (\epsilon)t].
\end{align*}

For part 2, note that since the set of infected sites increases with time, if some vertex hasn't been reached at a particular time $s< t$, then it hasn't been reached at any earlier time, so
\[
	\left\{(1-\epsilon)t \limitshape \not\subseteq \reachedsymb_{\tmeasure}^{\vect 0;\Rd}(t-) \right\}
	= \bigcap_{s<t} \left\{(1-\epsilon)t \limitshape \not\subseteq \reachedsymb_{\tmeasure}^{\vect 0;\Rd}(s) \right\}.
\]
Now for any $s \ge \frac{1-\epsilon}{1-\epsilon/2} t$, we have $(1-\epsilon)t\limitshape \subseteq (1-\epsilon/2)s \limitshape$ and hence
\[
	\left\{(1-\epsilon)t \limitshape \not\subseteq \reachedsymb_{\tmeasure}^{\vect 0;\Rd}(s) \right\}
	\subseteq \left\{(1-\epsilon/2)s \limitshape \not\subseteq \reachedsymb_{\tmeasure}^{\vect 0;\Rd}(s) \right\}.
\]
Take $s=(1-\epsilon/2)t$. Then $\frac{1-\epsilon}{1-\epsilon/2} t<s< t$, so the above inclusions imply
\begin{align*}
	\Pr \left\{(1-\epsilon)t \limitshape \not\subseteq \reachedsymb_{\tmeasure}^{\vect 0;\Rd}(t-) \right\}
	&\le \Pr \left\{(1-\epsilon/2)s \limitshape \not\subseteq \reachedsymb_{\tmeasure}^{\vect 0;\Rd}(s) \right\} \\
	&\le C(\epsilon/2)\cdot \exp \left[-c(\epsilon/2)s\right]
	= C_1(\epsilon)\cdot \exp\left[-c_1(\epsilon)t\right].
\qedhere
\end{align*}
\end{proof}

In fact, while Lemma~\ref{GM_large_dev_lem} tells us that the set of sites reached at a specific time $t$ lies within $\mu$-distance $\epsilon t$ of the deterministic process with high probability as $t\to\infty$, we can strengthen this result even further to get an exponential bound on the probability that $\reachedsymb_{\tmeasure}^{\vect 0;\Rd}(t)$ is in the right range for \emph{all} times larger than a given time. In other words, we can get exponential tail bounds on the amount of time we have to wait for the Shape Theorem to kick in for good.

\begin{lem}\label{trapping_time_lem}
Let $\tmeasure$ be an \iid\ traversal measure on $\edges{\Zd}$ satisfying \finitespeed{d} and \expmoment, and let $\mu$ be the shape function for $\tmeasure$ from Theorem~\ref{shape_thm}.
For any $\epsilon>0$, there exist two positive constants $C_1'$ and $c_1'$ such that for any $t>0$,
\begin{enumerate}
\item
$\displaystyle
	\Pr
	\left\{\reachedsymb_{\tmeasure}^{\vect 0;\Rd}(s) \subseteq (1+\epsilon)s \interior{\limitshape}
	 \text{ for all } s\ge t\right\}
	\ge 1 - C_1' e^{-c_1' t},
$
and

\item
$\displaystyle
	\Pr
	\left\{(1-\epsilon)s \limitshape \subseteq \reachedsymb_{\tmeasure}^{\vect 0;\Rd}(s-)
	 \text{ for all } s\ge t\right\}
	\ge 1 - C_1' e^{-c_1' t}.
$
\end{enumerate}

\end{lem}

\begin{proof}
First we bound
\[
	\Pr
	\left\{(1-\epsilon)s \limitshape \subseteq \reachedsymb_{\tmeasure}^{\vect 0;\Rd}(s-) \subseteq
	\reachedsymb_{\tmeasure}^{\vect 0;\Rd}(s) \subseteq (1+\epsilon)s \interior{\limitshape}
	\text{ for all } s\in [t,t+1)\right\}.
\]
Let $C_1''(\epsilon)=C_1(\epsilon/2)\vee \exp [c_1(\epsilon/2) \cdot (1+2/\epsilon)]$ and $c_1'(\epsilon)=c_1(\epsilon/2)$. For $t\in [0,\infty)$, define
\[
E_t =  \left\{\reachedsymb_{\tmeasure}^{\vect 0;\Rd}(s) \not\subseteq (1+\epsilon)s \interior{\limitshape}
	 \text{ for some } s\in [t, t+1)\right\}.
\]
Using the monotonicity of $\reachedsymb_{\tmeasure}^{\vect 0;\Rd}(t)$ and $t \interior{\limitshape}$ with $t$, for any $t\ge 0$ we have
\begin{align*}
E_t &= \bigcup_{t\le s< t+1} \left\{\reachedsymb_{\tmeasure}^{\vect 0;\Rd}(s)
	 \not\subseteq (1+\epsilon)s \interior{\limitshape}\right\}\\
&\subseteq \bigcup_{t\le s\le t+1} \left\{\reachedsymb_{\tmeasure}^{\vect 0;\Rd}(s)
	 \not\subseteq (1+\epsilon)t \interior{\limitshape}\right\}
=\left\{\reachedsymb_{\tmeasure}^{\vect 0;\Rd}(t+1) \not\subseteq (1+\epsilon)t \interior{\limitshape}\right\}.
\end{align*}
Now, if $t\ge \frac{2}{\epsilon} +1$, then $(1+\frac{\epsilon}{2})(t+1)\interior{\limitshape} \subseteq (1+\epsilon)t \interior{\limitshape}$, so
\[
\left\{\reachedsymb_{\tmeasure}^{\vect 0;\Rd}(t+1) \not\subseteq (1+\epsilon)t \interior{\limitshape}\right\}
\subseteq \left\{\reachedsymb_{\tmeasure}^{\vect 0;\Rd}(t+1) \not\subseteq (1+\epsilon/2)(t+1) \interior{\limitshape}\right\}
\]
Therefore, using the above inclusions and applying Part 1 of Lemma~\ref{stronger_dev_lem}, for all $t\ge \frac{2}{\epsilon} +1$ we have
\begin{align*}
\Pr(E_t)
\le \Pr\left\{\reachedsymb_{\tmeasure}^{\vect 0;\Rd}(t+1) \not\subseteq (1+\epsilon/2)(t+1) \interior{\limitshape}\right\}
&\le C_1(\epsilon/2) \cdot e^{-c_1(\epsilon/2) (t+1)}\\
&\le  C_1''(\epsilon) \cdot e^{-c_1'(\epsilon) t}.
\end{align*}
On the other hand, if $t\le \frac{2}{\epsilon}+1$, then $C_1''(\epsilon) \cdot e^{-c_1'(\epsilon) t}\ge 1$ by our choice of $C_1''$ and $c_1'$, so the above inequality is in fact valid for all $t>0$.


Now, for $t\in [0,\infty)$, define
\[
F_t = \left\{(1-\epsilon)s \limitshape \not\subseteq \reachedsymb_{\tmeasure}^{\vect 0;\Rd}(s-)
	 \text{ for some } s\in [t,t+1)\right\}.
\]
Using the monotonicity of $\reachedsymb_{\tmeasure}^{\vect 0;\Rd}(t-)$ and $t\limitshape$ with $t$, for any $t\ge 0$ we have
\begin{align*}
F_t &= \bigcup_{t\le s< t+1} 
	\left\{(1-\epsilon)s \limitshape \not\subseteq \reachedsymb_{\tmeasure}^{\vect 0;\Rd}(s-)\right\}\\
&\subseteq \bigcup_{t\le s\le t+1} 
	\left\{(1-\epsilon)s \limitshape \not\subseteq \reachedsymb_{\tmeasure}^{\vect 0;\Rd}(t-)\right\}
=\left\{(1-\epsilon)(t+1) \limitshape \not\subseteq \reachedsymb_{\tmeasure}^{\vect 0;\Rd}(t-)\right\}.
\end{align*}
Now if $t\ge \frac{2}{\epsilon}-2$, then $(1-\epsilon)(t+1) \limitshape\subseteq (1-\frac{\epsilon}{2})t \limitshape$, so
\[
\left\{(1-\epsilon)(t+1) \limitshape \not\subseteq \reachedsymb_{\tmeasure}^{\vect 0;\Rd}(t-)\right\}
\subseteq
\left\{(1-\epsilon/2)t \limitshape \not\subseteq \reachedsymb_{\tmeasure}^{\vect 0;\Rd}(t-)\right\}.
\]
Therefore, using the above inclusions and applying Part 2 of Lemma~\ref{stronger_dev_lem}, for all $t\ge \frac{2}{\epsilon} -2$ we have
\begin{align*}
\Pr(F_t)\le
\Pr\left\{(1-\epsilon/2)t \limitshape \not\subseteq \reachedsymb_{\tmeasure}^{\vect 0;\Rd}(t-)\right\}
\le C_1(\epsilon/2) \cdot e^{-c_1(\epsilon/2) t}
\le  C_1''(\epsilon) \cdot e^{-c_1'(\epsilon) t}.
\end{align*}
On the other hand, if $t\le \frac{2}{\epsilon}-2<  \frac{2}{\epsilon}+1$, then $C_1''(\epsilon) \cdot e^{-c_1'(\epsilon) t}\ge 1$ by our choice of $C_1''$ and $c_1'$, so the above inequality is in fact valid for all $t>0$.

Finally, since
\[
\bigcup_{\substack{s\ge t\\ s\in \R}} (E_t \cup F_t) \subseteq 
\bigcup_{\substack{n=\lfloor t \rfloor\\ n\in\N}}^\infty (E_n\cup F_n),
\]
we have
\begin{align*}
\Pr \left(\bigcup_{s\ge t} (E_t \cup F_t)\right)
&\le \sum_{n=\lfloor t \rfloor}^\infty \Pr(E_n \cup F_n)\\
&\le \sum_{n=\lfloor t \rfloor}^\infty 2 C_1'' e^{-c_1' n}\\
&= \frac{2C_1''}{1-e^{-c_1'}}\cdot e^{-c_1' \lfloor t \rfloor}\\
&\le \frac{2C_1''}{e^{c_1'}-1}\cdot e^{-c_1' t}
\defines C_1' e^{-c_1' t}.
\qedhere
\end{align*}
\end{proof}


The following lemma generalizes the above results and is used in Chapter~\ref{cone_growth_chap}.

\begin{lem}[Equivalence of deviation bounds in first-passage growth]
\label{equiv_dev_lem}
Let $\tmeasure$ be a random traversal measure on $\edges{\Zd}$ in some probability space $\probabilityspace$, and let $\mu$ be a norm on $\Rd$. For any $A,A',S\subseteq\Rd$, the following are equivalent.
\begin{enumerate}
\item There exists an $o(t)$ function $\phi:\R_+\to\R_+$ such that for any $\epsilon>0$, there exist positive constants $C$ and $c$ such that for all $t>0$,
\begin{align*}
\Pr \eventthat[1]{
	\rreacheddset[2]{\mu}{A'}{S}{[(1-\epsilon)t-\phi(t)]-}
	\subseteq \rreachedset{\tmeasure}{A}{S}{t}
	\text{ and }
	\crreachedset{\tmeasure}{A}{S}{t-}
	\subseteq \rreacheddset[2]{\mu}{A'}{S}{(1+\epsilon)t+\phi(t)}}\quad\\
	\ge 1- Ce^{-ct}.
\end{align*}

\item For any $\epsilon>0$ and any $o(t)$ function $\phi:\R_+\to \R_+$, there exist positive constants $C'$ and $c'$ such that for all $t_0\ge 0$, 
\begin{align*}
\Pr \eventthat[1]{\mathstack{\forall t\ge t_0,\;
	\rreacheddset[2]{\mu}{A'}{S}{(1-\epsilon)t +\phi(t)}
	\subseteq \crreachedset{\tmeasure}{A}{S}{t-}
	\hfill}
	{ \hskip 2.9cm \text{and}\quad \rreachedset{\tmeasure}{A}{S}{t}
	\subseteq \rreacheddset[2]{\mu}{A'}{S}{[(1+\epsilon)t-\phi(t)]-}
	}} \ge 1- C'e^{-c't_0}.
\end{align*}

\end{enumerate}
\end{lem}

\begin{proof}
This is a generalization of Lemma~\ref{trapping_time_lem}; it follows from Propositions~\ref{exponential_trapping_time_prop} and \ref{equivalent_trapping_time_prop} by choosing the functions $f$ and $g$ to be $\ctime{\tmeasure}{S}{A}{\cdot}$, $\rptime{\tmeasure}{S}{A}{\cdot}$, or $\idist{\mu}{S}{A'}{\cdot}$, as appropriate.
\end{proof}





\vita{
%
%
Nathaniel Blair-Stahn was born in San Fransisco, California, in 1978, and has lived in several places in the western United States, including Hawaii, California, South Dakota, Arizona, and Washington State. He attended the University of Arizona in Tucson, earning a Bachelor of Science in Mathematics in 2002. In 2012 he earned a Doctor of Philosophy in Mathematics at the University of Washington in Seattle.
}

\end{document}